\documentclass[12pt]{article}
\usepackage{epsfig}
\usepackage{rotating}
\title{Outer Billiards on Kites}
\author{Richard Evan Schwartz \thanks{\hskip 5 pt 
This research is supported by 
N.S.F. Grant DMS-0604426}}

\newtheorem{theorem}{Theorem}[section]

\newtheorem{lemma}[theorem]{Lemma}

\newtheorem{corollary}[theorem]{Corollary}
\newtheorem{conjecture}[theorem]{Conjecture}

\def\startproof{{\bf {\medskip}{\noindent}Proof: }}

\def\endproof{$\spadesuit$  \newline}

\def\C{\mbox{\boldmath{$C$}}}%
\def\H{\mbox{\boldmath{$H$}}}%
\def\J{\mbox{\boldmath{$J$}}}%
\def\N{\mbox{\boldmath{$N$}}}%
\def\P{\mbox{\boldmath{$P$}}}%
\def\Q{\mbox{\boldmath{$Q$}}}%
\def\R{\mbox{\boldmath{$R$}}}%
\def\Z{\mbox{\boldmath{$Z$}}}%

\begin{document}

\nonumber
\noindent
{\huge {\bf Outer Billiards on Kites\/}
\newline
\newline
\newline
\indent 
\indent
\indent
\indent
\indent
\indent
\indent
{\bf by\/}
\newline
\newline
\newline
{\bf Richard Evan Schwartz\/}
\/}

\newpage

\noindent
{\huge {\bf Preface\/}\/}
\newline
\newline
Outer billiards is a basic dynamical system defined
relative to a convex shape in the plane.  B.H. Neumann
introduced outer billiards in the 1950s, and J. Moser
popularized outer billiards in the 1970s as a toy
model for celestial mechanics.  Outer billiards is
an appealing dynamical system because of its simplicity
and also because of its connection to such topics as
interval
exchange maps, piecewise isometric actions, and
area-preserving actions.
There is a lot left to
learn about these kinds of dynamical systems, and
a good understanding of outer billiards might
shed light on the more general situation.

The {\it Moser-Neumann question\/}, one of the central
problems in this subject, asks
{\it Does there exist an outer billiards system with
an unbounded orbit\/}?  Until recently, all the work
on this subject has been devoted to proving that
all the orbits are bounded for various classes of
shapes.  We will detail these results in the
introduction.

Recently we answered the Moser-Neumann question in the
affirmative by showing that outer billiards has an
unbounded orbit when defined relative to the Penrose
kite, the convex quadrilateral that arises in the
famous Penrose tiling.  Our proof involves special
properties of the Penrose kite, and naturally
raises questions about generalizations.

In this monograph we will give a more general and robust answer
to the Moser-Neumann question.  We will prove that
outer billiards has unbounded orbits when defined
relative to any irrational kite.  A {\it kite\/} is
probably best defined as a ``kite-shaped'' quadrilateral.
(See the top of \S \ref{mainresult} for a non-circular definition.)
The kite is irrational if
it is not affinely equivalent to a quadrilateral with
rational vertices.
Our analysis uncovers some of
the deep structure underlying outer billiards on
kites, including connections to
self-similar tilings, higher dimensional
polytope exchange maps,
Diophantine approximation, the modular group, 
and the universal odometer.

I discovered 
every result in this monograph by experimenting
with my computer program, Billiard King, a Java-based graphical
user interface.  For the most part,
the material here is logically independent from
Billiard King, but I encourage the serious reader of
this monograph to download
Billiard King from my website \footnote{www.math.brown.edu/$\sim$res} and play with it. 
My website also has an interactive guide to this monograph,
in which many of the basic ideas and constructions
are illustrated with interactive Java applets.
\newline

There are a number of people I would like to thank.
I especially thank Sergei
Tabachnikov, whose great book {\it Geometry and Billiards\/}
first taught me about outer billiards. Sergei has 
constantly encouraged me
as I have investigated this topic, and he
has provided much mathematical insight along the way.

I thank Yair Minsky for his work on the
punctured-torus case of the Ending Lamination
Conjecture. It might seem strange to relate outer
billiards to punctured-torus bundles, but there
seems to me to be a common theme.  In both cases,
one studies the  limit of geometric objects indexed by rational
numbers and controlled in some sense by the
Farey triangulation.

I thank Eugene Gutkin for the explanations
he has given me about his work on outer billiards.
The work of Gutkin-Simanyi and others on the periodicity
of the orbits for rational polygons provided
the theoretical underpinnings for some of my initial computer
investigations.

I thank Jeff Brock, Peter Doyle, Dmitry Dolgopyat,
David Dumas, Giovanni Forni, 
Richard Kent, Howie Masur, Curt McMullen,
John Smillie, and Ben Wieland, for 
various helpful conversations about this work.

I thank the National Science Foundation for their
continued support, currently
in the form of the grant DMS-0604426.

I thank my home institution, Brown University, for
providing an excellent research environment during
the genesis of most of this work.  I also
thank the Institut des Hautes Etudes Scientifiques
for providing a similarly excellent
research enviromnent, during the summer of 2008.

I dedicate this monograph to my parents, Karen and Uri.

\newpage

\noindent
{\bf {\huge Table of Contents\/}\/}

\noindent
\begin{tabular}{lr}
\\
1. Introduction $\hskip 179 pt$    & 6
\end{tabular}

\noindent
\begin{tabular}{lr}
\\
{\bf Part I\/} $\hskip 217 pt$ & 22 \\
2. The Arithmetic Graph & 23 \\
3. The Hexagrid Theorem & 34 \\
4. Period Copying & 42 \\
5. Proofs of the Basic Results & 48
\end{tabular}

\vskip 20 pt

\noindent
\begin{tabular}{lr}
\\
{\bf Part II\/} $\hskip 210 pt$ & 56 \\
6. The Master Picture Theorem &  57\\
7. The Pinwheel Lemma & 66\\
8. The Torus Lemma & 79\\
9. The Strip Functions & 88\\
10. Proof of the Master Picture Theorem &  95 \\
11. Some Formulas &  100
\end{tabular}

\vskip 20 pt

\noindent
\begin{tabular}{lr}
\\
{\bf Part III\/}  $\hskip 200 pt$   & 108 \\
12. Proof of the Embedding Theorem  & 109\\
13. Extension and Symmetry & 116 \\
14. The Structure of the Doors & 123 \\
15. Proof of the Hexagrid Theorem I & 128\\
16. Proof of the Hexagrid Theorem II & 135 \\
17. The Barrier Theorem & 147 \\
\end{tabular}

\vskip 20 pt

\noindent
\begin{tabular}{lr}
\\
{\bf Part IV\/}  $\hskip 200 pt$   & 156 \\
18. Proof of the Superior Sequence Lemma  & 157\\
19. The Diophantine Lemma & 164 \\
20. Existence of Strong Sequences & 174 \\
21. Proof of the Decomposition Theorem & 179 \\
\end{tabular}

\newpage

\noindent
\begin{tabular}{lr}
\\
{\bf Part V\/} $\hskip 217 pt$ & 188\\
22. Odd Approximation Results &  189 \\
23. The Fundamental Orbit &  194 \\
24. Most of the Comet Theorem & 208 \\
25. Dynamical Consequences & 220 \\
26. Geometric Consequences & 230\\
\end{tabular}

\noindent
\begin{tabular}{lr}
\\
{\bf Part VI\/} $\hskip 217 pt$ & 238\\
27. Proof of the Copy Theorem &239 \\
28. Pivot Arcs in the Even Case & 247 \\
29. Proof of the Pivot Theorem & 260\\
30. Proof of the Period Theorem & 273 \\
31. The End of the Comet Theorem & 280 \\
\end{tabular}

\noindent
\begin{tabular}{lr}
\\
32. References & 295
\end{tabular}

\newpage

\newpage
\noindent

\section{Introduction}

\subsection{History of the Problem}

B.H. Neumann [{\bf N\/}] introduced {\it outer billiards\/} in
the late 1950s.  In the 1970s, J. Moser [{\bf M1\/}]
popularized outer billiards as a toy model
for celestial mechanics.  One appealing
feature of polygonal outer billiards is that it 
gives rise to a piecewise isometric mapping
of the plane. Such maps have close connections
to interval exchange transformations and
more generally to polygon exchange maps.
See [{\bf T1\/}] and $[{\bf DT\/}]$ for an
exposition of outer billiards and many references.

To define an outer billiards system,
one starts with a bounded convex set $K \subset \R^2$
and considers a point $x_0 \in \R^2-K$.  
One defines $x_1$ to be the point such that
the segment $\overline{x_0x_1}$ is tangent to $K$ at its
midpoint and $K$ lies to the right of the ray
$\overrightarrow{x_0x_1}$.  (See Figure 1.1 below.)
The iteration $x_0 \to x_1 \to x_2...$
is called the {\it forwards outer billiards orbit\/} of $x_0$.
It is defined for almost every point of $\R^2-K$.
The backwards orbit is defined similarly.

\begin{center}
\psfig{file=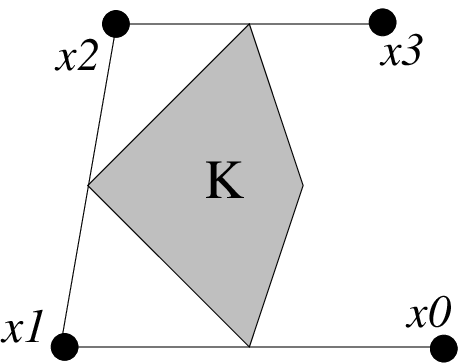}
\newline
{\bf Figure 1.1:\/} Outer Billiards
\end{center}

Moser [{\bf M2\/}, p. 11]
attributes the
following question to Neumann {\it circa\/} 1960, though it is sometimes
called {\it Moser's Question\/}.
\newline
\newline
{\bf Question:\/} {\it Is there an outer billiards system with an unbounded
orbit? \/}
\newline

This question is an idealized version of the question about the
stability of the solar system.  The Moser-Neumann question has
been considered by various authors.  Here is a list of the main
results on the question.

\begin{itemize}
\item J. Moser [{\bf M\/}] sketches a proof, inspired by
K.A.M. theory, that outer
billiards on $K$ has all bounded orbits provided
that $\partial K$ is at least $C^6$ smooth and positively curved.
R. Douady gives a complete proof in his thesis, [{\bf D\/}].
\item P. Boyland [{\bf B\/}] gives examples of $C^1$ smooth convex
domains for which an orbit can contain the domain boundary in
its $\omega$-limit set.
\item In [{\bf VS\/}], [{\bf Ko\/}], and (later, but with different methods)
[{\bf GS\/}], it is proved
that outer billiards on a {\it quasirational polygon\/} has
all orbits bounded.  This class of polygons
includes rational polygons and also regular polygons.
In the rational case,
all defined orbits are periodic.   
\item  S. Tabachnikov analyzes
the outer billiards system for the regular pentagon
and shows that there are some non-periodic (but bounded)
orbits.  See [{\bf T1\/}, p 158] and the references there. 
\item D. Genin [{\bf G\/}] shows that all orbits are bounded
for the outer billiards systems associated to trapezoids.
He also makes a brief numerical study of
a particular irrational kite based on the
square root of $2$, observes possibly
unbounded orbits, and indeed conjectures that this is
the case.
\item Recently, in [{\bf S\/}] we proved that outer billiards on
the Penrose kite has unbounded orbits, thereby answering
the Moser-Neumann question in the affirmative.
The Penrose kite is the convex quadrilateral that arises in the
Penrose tiling. 
\item Very recently, D. Dolgopyat and B. Fayad [{\bf DF\/}]
show that outer billiards
around a semicircle has some unbounded orbits.   Their proof
also works for ``circular caps'' sufficiently close to
the semicircle.  This is a second affirmative answer
to the Moser-Neumann question.
\end{itemize}

The result in [{\bf S\/}] naturally raises questions about
generalizations.
The purpose of this monograph
is to develop the theory of outer billiards on kites and
show that the phenomenon of unbounded orbits for polygonal
outer billiards is (at least for kites) quite robust.
We think that the theory we develop here will work,
to some extent, for polygonal outer
billiards in general, though right now a general theory
is beyond us.

We mention again that we discovered all the results in the
monograph through computer experimentation.  The interested
reader can download my program, Billiard King, from my
website \footnote{www.math.brown.edu/$\sim$res}.

\subsection{The Basic Results}
\label{mainresult}

For us, a
{\it kite\/} is a quadrilateral of the form $K(A)$, with vertices
\begin{equation}
\label{parameter}
(-1,0); \hskip 20 pt (0,1) \hskip 20 pt (0,-1) \hskip 20 pt (A,0); \hskip 30 pt
A \in (0,1).
\end{equation}
Figure 1.1 shows an example.
We call $K(A)$ {\it (ir)rational\/}
iff $A$ is (ir)rational.
Outer billiards is an affinely invariant system, and
any quadrilateral that is traditionally called a kite
is affinely equivalent to some $K(A)$.  

Let $\Z_{\rm odd\/}$ denote the set of odd integers.  Reflection
in each vertex of $K(A)$ preserves $\R \times \Z_{\rm odd\/}$.
Hence, outer billiards on $K(A)$ 
preserves $\R \times \Z_{\rm odd\/}$.  We say that
a {\it special orbit\/} on $K(A)$ is an orbit contained
in $\R \times \Z_{\rm odd\/}$.  This monograph \footnote{Some of our
theory works for the general orbit, and there seems to be quite a
lot to say, but we will not say it here.  The special orbits are
hard enough for us already.} only
discusses special orbits.

We call an orbit {\it forwards erratic\/} if the forwards orbit
is unbounded and also returns to every neighborhood of a kite
vertex.   We make the same definition for the backwards
direction.  We call an orbit {\it erratic\/} if it is both
forwards and backwards erratic.  Say that a
{\it trimmed Cantor set\/} is a set of the form
$C-C'$, where $C$ is a Cantor set and $C'$ is countable.
Note that a trimmed Cantor set is an uncountable set.
Here are our $3$ basic results.

\begin{theorem}[Erratic Orbits]
\label{erratic}
On any irrational kite, the union of
special erratic orbits contains a
trimmed Cantor set.
\end{theorem}

\begin{theorem}[Dichotomy]
\label{dichotomy0}
On any irrational kite, every special
orbit is either periodic or else unbounded in
both directions.  
\end{theorem}

\begin{theorem}[Density]
\label{density}
On any irrational kite, the union of
periodic special orbits is open dense in
$\R \times \Z_{\rm odd\/}$.
\end{theorem}

Thanks to the work mentioned above, we already know
that all orbits are bounded on rational kites.
The Erratic Orbits Theorem therefore has the
following simple corollary.

\begin{corollary}
\label{unbounded}
 Outer billiards on a kite has an unbounded
orbit if and only if the kite is irrational.
\end{corollary}

Our monograph comes in 6 parts.  Parts I-IV constitute
a self-contained subset of the monograph designed to
prove the results listed above. Parts V-VI go deeper
into the subject, and establish the Comet Theorem,
a fairly complete description of the set of unbounded
special orbits.

\subsection{The Comet Theorem}
\label{dmr}

The Comet Theorem has a number of corollaries that are
much easier to state than the result itself.
We state some of these first.  Let $U_A$ denote
the set of unbounded special orbits relative to
an irrational $A \in (0,1)$.

\begin{itemize}
\item $U_A$ is minimal:  Every orbit in $U_A$ is
dense in $U_A$. 
\item $U_A$ is locally homogeneous:
Every two points in $U_A$ have arbitrarily small neighborhoods
that are isometric to each other.
\item $U_A$ has length $0$.  Hence, almost every point in
$\R \times \Z_{\rm odd\/}$ is periodic.
\item Let $u(A)$ denote the Hausdorff dimension of $U_A$.
The function $u$ maps every 
open subset of $[0,1]-\Q$ onto $[0,1]$ and yet is
almost everywhere constant.  (We don't
know ``the constant''.)
\item Let $\Gamma_2 \subset SL_2(\Z)$ denote the subgroup
consisting of matrices congruent to the identity mod $2$.
As usual, $\Gamma_2$ acts on $\R \cup \infty$ by linear
fractional transformations.  The function $u$ is constant
on $\Gamma_2$-orbits.
\end{itemize}

We will deduce these results, and many others,
from the Comet Theorem in \S \ref{comet2} and
\S \ref{lengthdim}.
We turn now to the statement of the Comet Theorem.
Consider the regions
\begin{equation}
\label{atmosphere}
I=[0,2] \times \{-1\}; \hskip 30 pt 
J=[-2,2] \times \{-1,1\}=\bigcup_{k=0}^3 (\psi')^k(I).
\end{equation}
Here $\psi'$ is the outer billiards map.
The domains $I$ and $J$ turn out to be canonical
domains for outer billiards on kites, as we have
normalized them.   The Comet Theorem provides
a model for the way the unbounded orbits return to $I$.

\begin{center}
\resizebox{!}{1.8in}{\includegraphics{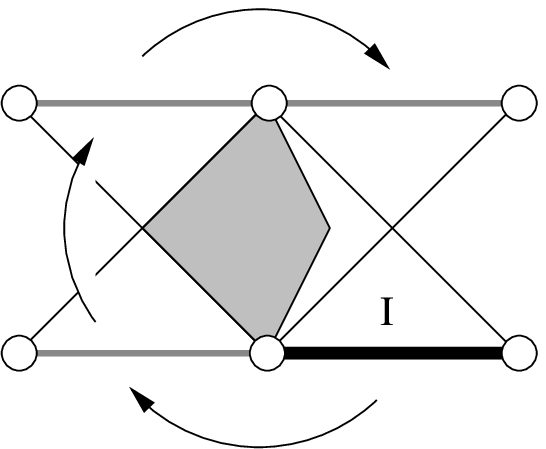}}
\newline
{\bf Figure 1.2:\/} $I$ is black and $J$ is grey.
\end{center}

Say that $p/q$ is {\it odd\/} or
{\it even\/} according to whether $pq$ is odd or even.
There is a unique sequence
$\{p_n/q_n\}$ of distinct odd rationals, converging to
$A$, such that $p_0/q_0=1/1$ and
$|p_nq_{n+1}-q_np_{n+1}|=2$ for all $n$.
We call this sequence the {\it inferior sequence\/}. See
\S \ref{diop}.
This sequence is closely related to continued fractions.

We define
\begin{equation}
d_n={\rm floor\/}\bigg(\frac{q_{n+1}}{2q_n}\bigg); \hskip 30 pt n=0,1,2...
\end{equation}
Say that a {\it superior term\/} is a term $p_n/q_n$ such that
$d_n \geq 1$.   We will show that there are
infinitely many superior terms.
Say that the {\it superior sequence\/} 
is the subsequence of superior terms.
Say that the
{\it renormalization sequence\/} is the corresponding
subsequence of $\{d_n\}$.  We re-index so that the
superior and renormalization sequences are indexed by
$0,1,2...$.   The definitions that follow work
entirely with the superior sequence.

We define ${\cal Z\/}_A$ to be the inverse limit of the
system
\begin{equation}
\ldots \to \Z/D_3 \to \Z/D_2 \to \Z/D_1; \hskip 30 pt
D_n=\prod_{i=0}^{n-1}(d_i+1)
\end{equation}
We equip ${\cal Z\/}$ with a metric, defining
$d(x,y)=q_{n-1}^{-1}$,
where $n$ is the smallest index
such that $[x]$ and $[y]$ disagree in $\Z/D_n$.
In case $d_i+1=p$ for all $i$,
the metric $d$ is the $p$-adic metric.
In general,
${\cal Z\/}_A$ is a metric abelian group.
The map $x \to x+1$ is a canonical self-homeomorphism
called the {\it odometer map\/}.

We can identify the points of ${\cal Z\/}_A$ with the
sequence space
\begin{equation}
\label{productspace}
\Pi_A=\prod_{i=0}^{\infty} \{0,...,d_i\}
\end{equation}
Our identification works like this
\begin{equation}
\label{phi1}
\phi_1: \hskip 15 pt \sum_{j=0}^{\infty} \widetilde k_j D_j \in {\cal Z\/}_A \hskip 20 pt
\longrightarrow \hskip 20 pt
\{k_j\} \in \Pi_A.
\end{equation}
The elements on the left hand side are formal series, and
\begin{equation}
\label{tildek}
\widetilde k_j= \Bigg{\{} \matrix{k_j & {\rm if\/} & p_j/q_j<A \cr \cr
d_j-k_j & {\rm if\/}& p_j/q_j>A} 
\end{equation}
Our identification is a bit nonstandard, in that it
uses $\widetilde k_j$ in place of the more obvious
choice of $k_j$.

There is a map $\phi_2: \Pi_A \to \R \times \{-1\}$, defined
as follows.
\begin{equation}
\phi_2: \hskip 15 pt \{k_j\} \hskip 20 pt
\longrightarrow \hskip 20pt \bigg(\sum_{j=0}^{\infty} 2k_j \lambda_j,-1\bigg); \hskip 30 pt
\lambda_j=|Aq_j-p_j|.
\end{equation}
We define
$C_A=\phi_2(\Pi_A)$.
Equivalently,
\begin{equation}
C_A=\phi({\cal Z\/}_A); \hskip 50 pt
\phi=\phi_2 \circ \phi_1.
\end{equation}
It turns out that $\phi: {\cal Z\/}_A \to C_A$ is
a homeomorphism and
$C_A$ is a Cantor set whose convex hull is exactly $I$.
Let $C_A^{\#}$ denote the set obtained from $C_A$ by deleting
the endpoints of the complementary intervals in $I-C_A$.

The points $\phi(-1), \phi(0) \in C_A$ are
important points for us.
It turns out that these points have
well-defined orbits iff they lie in $C_A^{\#}$, which
happens iff
the superior sequence for $A$ is not eventually monotone.
Reflection in the midpoint of $I$ preserves
$C_A$ and swaps $\phi(-1)$ and $\phi(0)$.

Define
\begin{equation}
\Z[A]=\Z \oplus \Z A=\{mA+n|\ m,n \in \Z\}.
\end{equation}

Say that the {\it excursion distance\/} of a portion
of an outer billiards orbit is the maximum distance
from a point on this orbit-portion to the origin.

\begin{theorem}[Comet]
\label{cantor}
Let $U_A$ denote the set of unbounded special orbits
relative to an irrational $A \in (0,1)$.
\begin{enumerate}
\item  For any $N$ there is an $N'$ with the following
property.  If $\zeta \in U_A$ satisfies $\|\zeta\|<N$ then 
the $k$th outer billiards
iterate of $\zeta$ lies in $I$ for some $|k|<N'$.
Here $N'$ only depends on $N$ and $A$.
\item
$U_A \cap I=C_A^{\#}$.
The first return map $\rho_A: C_A^{\#} \to C_A^{\#}$ is defined precisely
on $C_A^{\#}-\phi(-1)$. The map $\phi^{-1} \circ \rho_A \circ \phi$,
wherever defined on ${\cal Z\/}_A$, equals the odometer.
\item For any $\zeta \in C_A^{\#}-\phi(-1)$, the
orbit-portion
between $\zeta$ and $\rho_A(\zeta)$ has excursion
distance in $[c_1^{-1} d^{-1},c_1 d^{-1}]$ and length in
$[c_2^{-1}d^{-2},c_2d^{-3}]$.
 Here $d=d(-1,\phi^{-1}(\zeta))$, and $c_1,c_2$ are universal positive
constants.
\item $C_A^{\#}=C_A-(2\Z[A] \times \{-1\})$. 
Two points in $U_A$ lie on the same orbit if and
only if the
difference of their first coordinates lies in
$2\Z[A]$.
\end{enumerate}
\end{theorem}

\noindent
{\bf Remarks:\/} \newline
(i)
To use a celestial analogy, we think of $I$ as
the visible sky, and the
special unbounded orbits as the comets. 
Item 1 says (in particular) that any comet visits $I$.
Item 2 describes the exact locations and
combinatorial structure of the visits.
Item 3 gives a coarse model for the excursion
distances and return times between visits.
Item 4 gives an algebraic view.  In short,
the Comet Theorem tells us where and (approximately)
when to point the telescope.
\newline
\newline
(ii)
Item 1 of the Comet Theorem is somewhat loose, in
that we don't know how $N'$ depends on $N$ and $A$.
In principle, one could extract estimates from
our proof, but we didn't try to do this.
\newline
\newline
(iii)
Lemma \ref{precomet} replaces our
bounds in Item 3 with explicit estimates.
The orders on all our bounds in Item 3 are
sharp except perhaps for the length upper-bound.
See the remarks following
Lemma \ref{precomet} for a discussion.
\newline
\newline
(iv)
The Comet Theorem has an analogue for the
backwards orbits.  The statement is
the same except that the point $\phi(0)$
replaces the point $\phi(-1)$ and the map
$x \to x-1$ replaces the odometer.
\newline
\newline
(v)
Our analysis will show that 
$\phi(0)$ and $\phi(-1)$ have well-defined
orbits iff they lie in $C_A^{\#}$.  Again,
this happens iff the superior sequence for
$A$ is not eventually monotone.   This
happens for a full-measure set of 
parameters.
Item 1 implies that the forwards orbit
of $\phi(-1)$, when defined, only
accumulates at $\infty$.  The same goes
for the backwards orbit of $\phi(0)$.
We think of $\phi(-1)$ as the ``cosmic ejector''.
When a comet comes close to this point, it gets
ejected way out into space.  Likewise, we think
of $\phi(0)$ as the ``cosmic attractor''.
\newline
\newline
(vi)
By symmetry $U_A \cap J$ consists of $4$ copies of
$C_A^{\#}$ arranged in a symmetric pattern about
the two kite vertices $(0,\pm 1)$.  Thus, the
Comet Theorem completely controls how the unbounded
special orbits return near the origin. 
\newline
\newline
(vii)
In \S \ref{model0} we will formalize the idea of
constructing a model for the dynamics of the
outer billiards map on $U_A$. 
The {\it solenoid\/} ${\cal S\/}_A$ is
the mapping cylinder for the odometer on
${\cal Z\/}_A$.  We delete a point from
${\cal S\/}_A$ and alter the metric so that this
deleted point lies infinitely far away.
That is, we create a cusp.
There is a fairly canonical way to do this, and
we call the result ${\cal C\/}_A$ the
{\it cusped solenoid\/}.
We will see that, in some sense
the time-one map for the geodesic flow
on ${\cal C\/}_A$ serves as a good model for
the dynamics on $U_A$.  This result is
really just a re-packaging of the Comet Theorem.
\newline
\newline
(viii)
For almost all choices of $A$, the object
${\cal Z\/}_A$ and its odometer coincide
with the {\it universal odometer\/}. This
is the profinite completion of $\Z$ -- i.e.,
the inverse limit over all finite cyclic
groups.   We call the corresponding object
${\cal C\/}_A$ the {\it universal cusped
solenoid\/}.  As we formalize in
\S \ref{model0} and \S \ref{UNIVERSAL}, 
the time-one map of
the geodesic flow on the universal
cusped solenoid serves as a good model,
in some sense,
for the dynamics on $U_A$ for almost all $A$.
In particular, for almost all $A$,
the return map to $C_A^{\#}$
is conjugate (modulo a countable set)
to the universal odometer.
\newline
\newline
(ix)
By the Dichotomy Theorem, all well-defined orbits
in $I-C_A$ are periodic. Conjecture \ref{gapmap3} describes
the dynamics of these points.  In brief, we can
identify $C_A$ with the ends of a certain directed
tree.  The first return map to $C_A$ is induced
by a certain automorphism of the directed tree.
The complementary intervals in $I-C_A$ are naturally
in bijection with the forward cones of the directed
tree.  Conjecture \ref{gapmap3} says that the
first return map to $I-C_A$ permutes these
intervals just as the tree automorphism
permutes the forward cones.
\newline
\newline
(x)
The $\Gamma_2$-invariance of the dimension function
$\dim(U_A)$ is a small reflection of the beautiful structure
of the sets $C_A$.   This monograph only scratches the surface.
Here is a structural result is outside the scope of the monograph.
Letting $C_A'$ denote the scaled-in-half version of $C_A$
that lives in the unit interval, it seems that
\begin{equation}
\label{butterfly}
C=\bigcup_{A \in [0,1]} \Big(C_A' \times \{A\}\Big) \subset [0,1]^2 \subset \R\P^2
\end{equation}
is the limit set of a semigroup $S \subset SL_3(\Z)$ that acts
by projective transformations. ($C_A$ can be defined even for
rational $A$.) The group closure of $S$ has finite index in
a maximul cusp of $SL_3(\Z)$. 
 The projective geometry underlying
the set $C$ emerges almost immediately from a good plot.
We might have included a plot here, but 
we don't know how to draw a good picture without
producing a huge picture file.  
 My website \footnote{www.math.brown.edu/$\sim$res/BilliardKing/Butterfly0}
 has a picture of $C$.   We produced this picture using 
the formula in Theorem \ref{discrete} below.

\subsection{Rational Kites}
\label{rationalresults}

We find it convenient to work with the square of the
outer billiards map.
Let $O_2(x)$ denote the square outer billiards orbit of $x$.
Let $I$ be as above, and let
\begin{equation}
\Xi=\R_+ \times \{-1,1\}.
\end{equation}
When $\epsilon \in (0,2/q)$, the orbit
$O_2(\epsilon,-1)$ has a combinatorial structure independent
of $\epsilon$. See Lemma \ref{irrat}.
Thus, $O_2(1/q,-1)$ is a natural
 representative of this orbit.
This orbit plays a crucial role in our proofs. Reflection
The following result is our basic mechanism for producing
unbounded orbits.

\begin{theorem}
\label{ex0}
\label{ex}
Let $p/q \in (0,1)$ be any rational. Relative to $A=p/q$ the
following is true. \begin{itemize}
\item If $p/q$ is odd then $O_2(1/q,-1)$ has diameter
between $(p+q)/2$ and $p+q$.
\item If $p/q$ is even then $O_2(1/q,-1)$ has diameter
between $(p+q)$ and $2(p+q)$.
\end{itemize}
\end{theorem}

Here is an amplification of the upper bound in Theorem \ref{ex0}.

\begin{theorem}
\label{orbit structure1}
If $p/q$ is odd let $\lambda=1$. If $p/q$ is even let
$\lambda=2$.
Each special orbit intersects $\Xi$ in exactly one set
of the form $I_k \times \{-1,1\}$, where
$$
I_k=(\lambda k(p+q),\lambda (k+1)(p+q))  \hskip 30 pt k=0,1,2,3...
$$
Hence, any special orbit intersects $\Xi$ in a set of
diameter at most $\lambda \cdot (p+q)$.
\end{theorem}

Theorem \ref{orbit structure1} is similar in spirit to a
result in [{\bf K\/}].  See \S \ref{discussK} for
a discussion.  

An outer billiards orbit on $K(A)$ is called {\it stable\/} if
there are nearby and combinatorially identical orbits
on $K(A')$ for all $A'$ sufficiently close to $A$.  
Otherwise, the orbit is called
{\it unstable\/}. 
In the odd case, 
$O_2(1/q,\pm 1)$ is unstable.  This fact is
of crucial importance to all our proofs.
Here is a classification of special
orbits in terms of stability.

\begin{theorem}
\label{orbit structure2}
In the even rational case, all special orbits are stable.
In the odd case,
the set $I_{k} \times \{-1,1\}$ contains exactly two unstable orbits,
$U_k^+$ and $U_k^-$, and these are conjugate
by reflection in the $x$-axis. In particular, we have
$U_0^{\pm}=O_2(1/q,\pm 1)$.
\end{theorem}

The preceding results give effective but somewhat
coarse global pictures of the special orbits. 
Here we describe a very precise picture of our
fundamental orbit $O_2(1/q,-1)$ near the origin.
Any odd rational $p/q$ appears as a term in a
superior sequence, and the terms before $p/q$ are
uniquely determined by $p/q$.  This is similar to what
happens for continued fractions.

\begin{theorem}
\label{discrete}
Let $\mu_i=|p_nq_i-q_np_i|$.  
$$O_2(1/q_n,-1) \cap I = \bigcup_{\kappa \in \Pi_n} \Big(X_n(\kappa),-1\Big); \hskip 20 pt
X_n(\kappa)=\frac{1}{q_n} \bigg(1+\sum_{i=0}^{n-1} 2k_i\mu_i\bigg).$$
\end{theorem}
To prove the Comet Theorem, we
will combine Theorem \ref{ex0} and Theorem \ref{discrete} and
then take a geometric limit.

Here we show Theorem \ref{discrete}
in action.  Relative to $19/49$, the intersection
$O_2(1/49) \cap \Xi$ has diameter between $34$ and $68$.
The odd rational $19/49$ determines the inferior sequence 
$$\frac{p_0}{q_0}=\frac{1}{1}, \frac{1}{3}, \frac{5}{13}, \frac{19}{49}=\frac{p_3}{q_3}.$$
All terms are superior, so this is also the superior sequence.
$n=3$ in our example, and the renormalization sequence is $1,2,1$.
The $\mu$ sequence is $30,8,2$.   The first coordinates of the
 $12$ points of 
$O_2(1/49) \cap I$ are given by 
$$\bigcup_{k_0=0}^1\  \bigcup_{k_1=0}^2\  \bigcup_{k_2=0}^1 \frac{2(30k_0+8k_1+2k_2)+1}{49}.$$
Writing these numbers in a suggestive way, the union above works out to
$$\frac{1}{49} \times \big(1\hskip 8pt 5\hskip 20 pt    17\hskip 8pt 21\hskip 20pt   33\hskip 8pt 37 \hskip 60 pt   
61\hskip 8pt 65\hskip 20 pt    77\hskip 8pt 81\hskip 20 pt    93\hskip 8pt 97\big).$$
The reader can check this example, and many others, using Billiard King.
\newline
\newline
{\bf Remarks:\/} \newline
(i)
A version of Theorem \ref{discrete} holds in the even case as well.
We will discuss the even case of Theorem \ref{discrete} in
\S \ref{extra}.  \newline
(ii) 
Theorem \ref{discrete} has a nice conjectural extension, which 
describes the entire return map to $I$.  See \S \ref{gapmap}.

\subsection{The Arithmetic Graph}
\label{preag}

All our results about special orbits derive from our analysis of
a fundamental object, which we call the {\it arithmetic graph\/}.
  One should think of the first return map to
$\Xi$, for rational parameters, as an essentially combinatorial
object.  The idea behind the arithmetic graph is to give a
$2$ dimensional pictorial representation of this combinatorial
object.  

The principle guiding our construction is that
sometimes it is better to understand the abelian group
$\Z[A]:=\Z \oplus\Z A$ as a module over $\Z$ rather than as a subset
of $\R$.  
 Our arithmetic graph is similar in spirit to the
lattice vector fields studied by Vivaldi et. al.
in connection with interval exchange transformations.  See e.g.
[{\bf VL\/}].
In this section we will explain the idea behind the
arithmetic graph.  In \S \ref{basic0} we will give a precise
construction.  

The arithmetic graph is most easily explained in the rational case.
Let $\psi$ be the square of the outer billiards map.
Let $\Xi=\R_+ \times \{-1,1\}$ as above.  It turns out that
every orbit starting on $\Xi$ eventually returns to
$\Xi$.  See Lemma \ref{rlemma}.   Thus we can define
the return map
\begin{equation}
\Psi: \Xi \to \Xi.
\end{equation}
We define the map $T: \Z^2 \to 2\Z[A] \times \{-1,1\}$ by the formula
\begin{equation}
T(m,n)=\Big(2Am+2n+\frac{1}{q},(-1)^{p+q+1}\Big).
\end{equation}
Here $A=p/q$.  

Up to the reversal of the direction of the dynamics, every point of
$\Xi$ has the same orbit as a point of the form $T(m,n)$, where
$(m,n) \in \Z^2$.   For instance, the orbit of $T(0,0)=(1/q,-1)$
is the same as the orbit of our favorite point $(1/q,-1)$ up to
reversing the dynamics.  The point here is that reflection in the
$x$-axis conjugates the outer billiards map to its inverse.

We form the graph $\widehat \Gamma(p/q)$ by joining the points
$(m_1,m_2)$ to $(m_2,n_2)$ when these points are sufficiently
close together and also
$T(m_1,n_1)=\Psi^{\pm 1}(m_2,n_2)$.
(The map $T$ is not injective, so we have choices to make.  That
is the purpose of the {\it sufficiently close\/} condition.)

We let $\Gamma(p/q)$ denote the component of $\widehat \Gamma(p/q)$
that contains $(0,0)$. This component tracks the orbit 
$O_2(1/q,-1)$, the main orbit of interest to us. 
When $p/q$ is odd, $\Gamma(p/q)$ is an infinite periodic
polygonal arc, invariant under translation by the vector
$(q,-p)$.  Note that $T(q,-p)=T(0,0)$.  When
$p/q$ is even, $\Gamma(p/q)$ is an embedded polygon.
This difference causes us to prefer the odd case.

We prove many structural theorems about the arithmetic
graph.  Here we mention $3$ central ones. 
We state these results vaguely here, and refer the reader to
the chapters where the precise statements are given.

\begin{itemize}
\item {\bf The Embedding Theorem\/} (\S 2):  $\widehat \Gamma(p/q)$ is a
disjoint union of embedded polygons and infinite embedded polygonal
arcs.  Every edge of $\widehat \Gamma(p/q)$ has length at most $\sqrt 2$.
The stable orbits correspond to closed polygons, and the unstable
orbits correspond to infinite (but periodic) polygonal arcs.
\item {\bf The Hexagrid Theorem\/} (\S \ref{hexagrid0}): The structure of
$\widehat \Gamma(p/q)$ is controlled by $6$ infinite families of
parallel lines.  The {\it quasiperiodic\/} 
structure is similar to what one sees in
DeBruijn's famous pentagrid construction of the Penrose tilings.
See [{\bf DeB\/}].
\item {\bf The Copy Theorem:\/} (\S \ref{copyproof}; also
Lemmas \ref{strongcopy} and Lemma \ref{weakcopy}.)
If $A_1$ and $A_2$ are two rationals
that are close in the sense of Diophantine approximation
then the corresponding arithmetic graphs
$\Gamma_1$ and $\Gamma_2$ have substantial agreement.
\end{itemize}

The Hexagrid Theorem causes $\Gamma(p/q)$ to have
an oscillation (relative to the line of slope $-p/q$ through the origin)
on the order of $p+q$.  The Hexagrid Theorem is responsible for
Theorems \ref{ex0}, \ref{orbit structure1}, and
\ref{orbit structure2}.  Referring to the superior sequence,
the Copy Theorem guarantees that
the structure the graph $\Gamma(p_n/q_n)$ is
copied by the graph $\Gamma(p_{n+1}/q_{n+1})$.  The
Copy Theorem 
is responsible for Theorem \ref{discrete}.
Thus, the Hexagrid Theorem and the Copy Theorem serve as
a kind of a team, with one result forcing large oscillations
in certain orbits, and the other result guaranteeing that
the oscillations are coherently organized in the family
of arithmetic graphs corresponding to the superior sequence.

We illustrate these ideas with some pictures.
Each picture shows $\Gamma(p/q)$ in reference to
the line of slope $-p/q$ through the origin.
The rationals 
$$\frac{7}{13},\ \frac{19}{35},\ \frac{45}{83}$$
form $3$ terms in a superior sequence.
Figure 1.3 shows a bit more than one period of $\Gamma(7/13)$.

\begin{center}
\resizebox{!}{1.9in}{\includegraphics{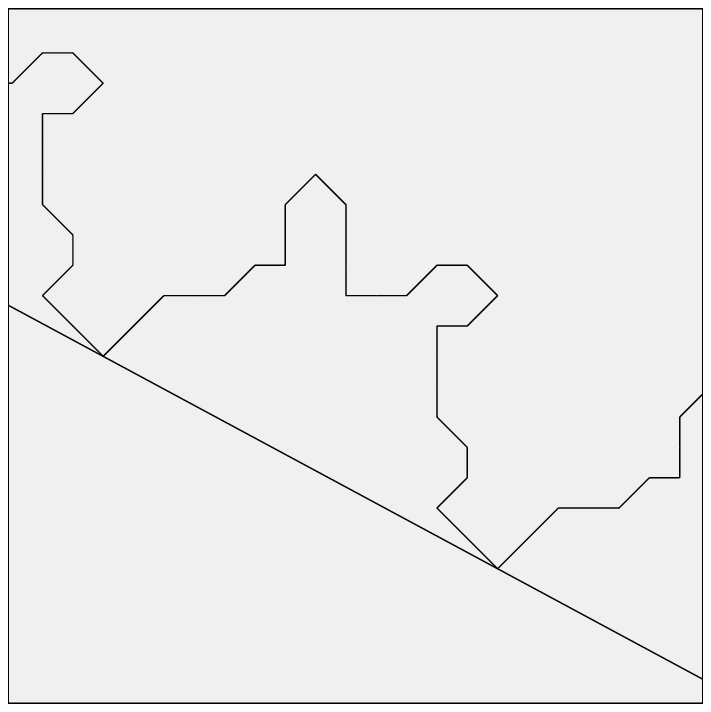}}
\newline
{\bf Figure 1.3:\/} The graph $\Gamma(7/13)$.
\end{center}

Figure 1.4 shows a picture of $\Gamma(19/35)$.
Notice that $\Gamma(19/35)$ has a much
wider oscillation, but also manages to copy a bit more
than one period of
$\Gamma(7/13)$. 
The reader can see many
more pictures like this using either
Billiard King or our interactive
guide to the monograph

\begin{center}
\resizebox{!}{4.1in}{\includegraphics{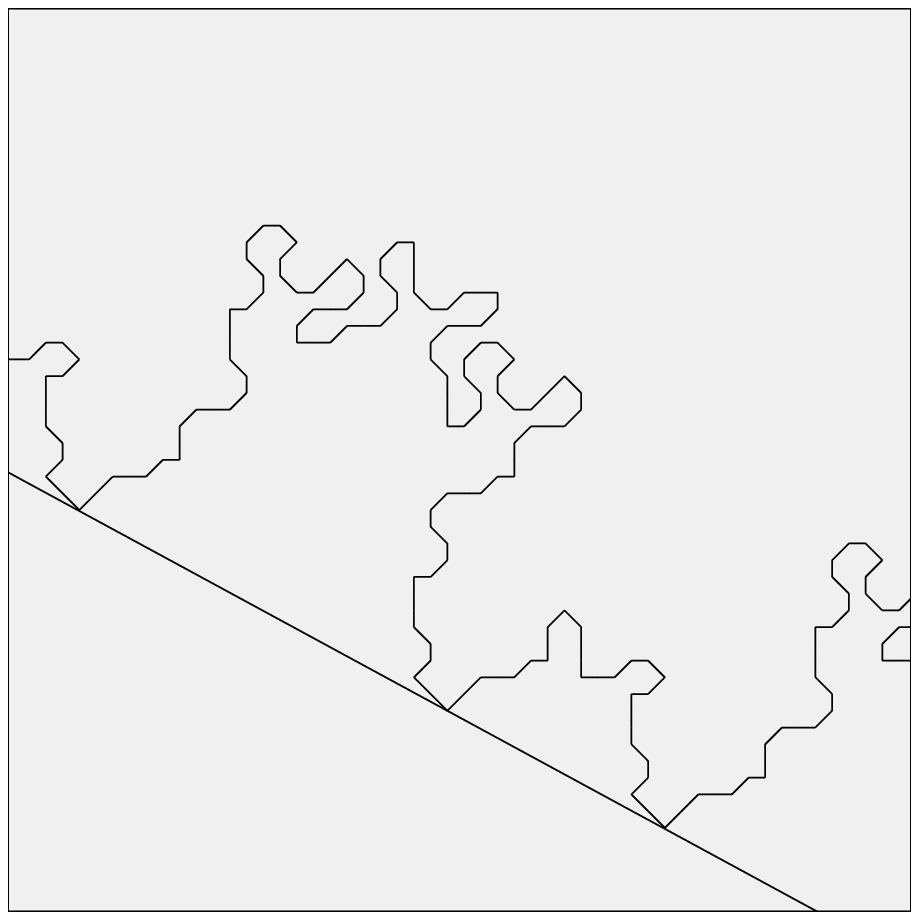}}
\newline
{\bf Figure 1.4:\/} The graph $\Gamma(19/35)$.
\end{center}

Figure 1.5 shows the same phenomenon for
$\Gamma(45/83)$.   This graph oscillates on
a large scale but still manages to copy a bit more than one
period of $\Gamma(19/35)$.  Hence $\Gamma(45/83)$
copies a period of $\Gamma(7/13)$ and a period
of $\Gamma(19/35)$.  That is, $\Gamma(45/83)$
oscillates on $3$ scales. 

\begin{center}
\resizebox{!}{5.2in}{\includegraphics{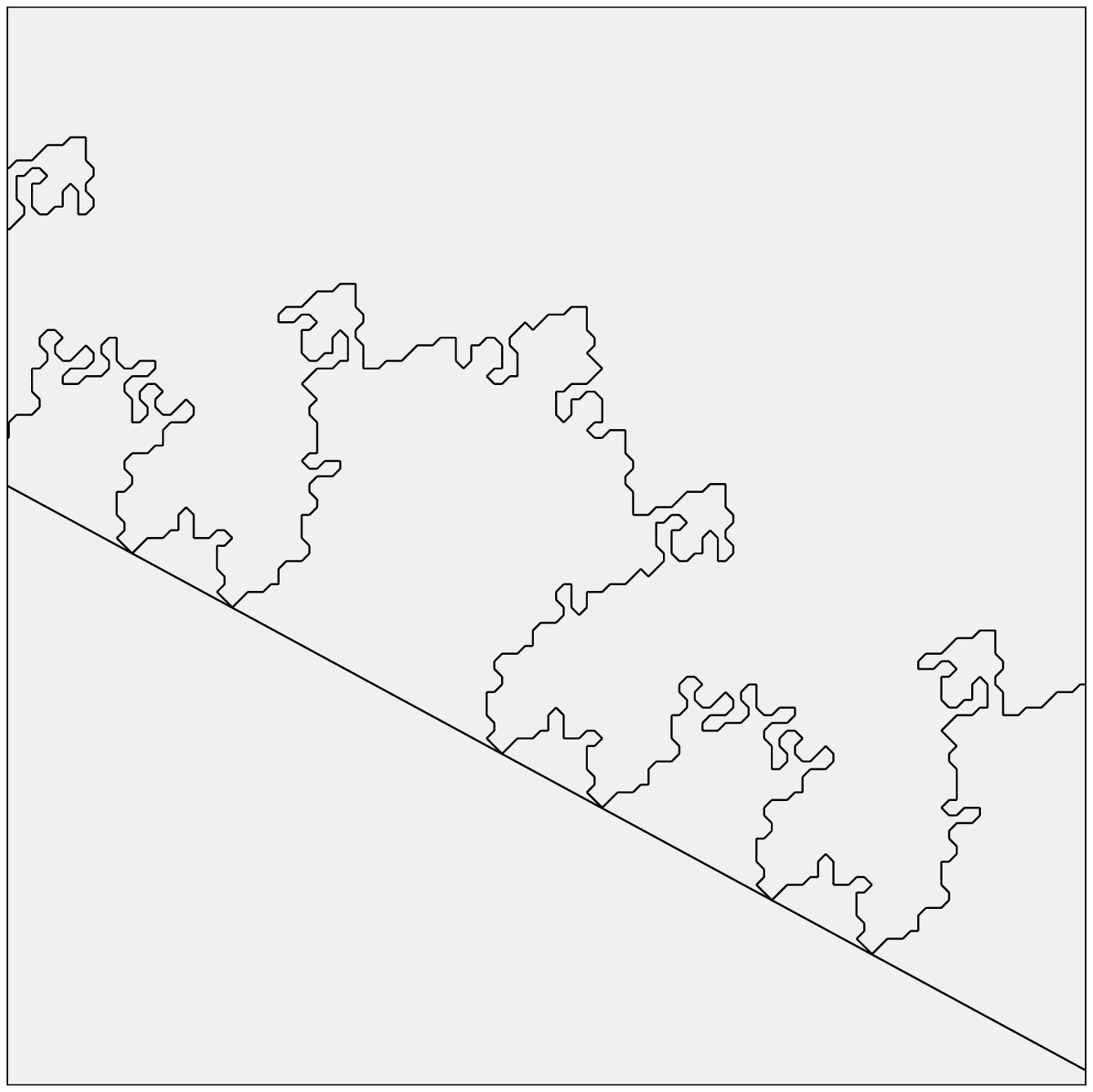}}
\newline
{\bf Figure 1.5:\/} The graph $\Gamma(45/83)$.
\end{center}

What emerges in these pictures is both the wide
excursions predicted by Theorem \ref{ex0} an
also the Cantor-like structure predicted by
Theorem \ref{discrete}.  Taking a limit of
this process, we produce a graph that oscillates
on all scales.  This limiting graph tracks the
sort of unbounded orbit discussed  above.

\subsection{The Master Picture Theorem}

Essentially all of our results have a common source,
the Master Picture Theorem.
The Master Picture Theorem 
is, in some sense, a closed form expression
for the arithmetic graph. 
We formulate and prove the Master Picture Theorem
 in Part II of the monograph.
Here we will give the reader
a feel for the result.

Recall that $\Xi=\R_+ \times \{-1,1\}$.  The
arithmetic graph encodes the dynamics of the
first return map $\Psi: \Xi \to \Xi$.
It turns out that $\Psi$ is an infinite
interval exchange map.  The Master Picture
Theorem reveals the following structure.
\begin{enumerate}
\item There is a locally affine map $\mu$ from $\Xi$ into
a union $\widehat \Xi$ of two $3$-dimensional
tori.  
\item There is a polyhedron exchange map $\widehat \Psi:
\widehat \Xi \to \widehat \Xi$, defined relative to
a partition of $\widehat \Xi$ into $28$ polyhedra.
\item The map $\mu$ is a semi-conjugacy between
$\Psi$ and $\widehat \Psi$.
\end{enumerate}
In other words, the return dynamics of $\widehat \Psi$
has a kind of compactification into a $3$ dimensional
polyhedron exchange map.  All the objects above
depend on the parameter $A$, but we have suppressed them
from our notation.

There is one master picture, a union of
two $4$-dimensional convex lattice polytopes partitioned into
$28$ smaller convex lattice polytopes, that controls
everything. For each parameter, one obtains the
$3$-dimensional picture by taking a suitable slice.

The fact that nearby slices give almost the same
picture is the source of our Copy Theorem.
The interaction between the map $\mu$ and the
walls of our convex polytope partitions is the
source of the Hexagrid Theorem.  The Embedding
Theorem follows from basic geometric properties
of the polytope exchange map in an elementary
way that is hard to summarize here.

My investigation of the Master Picture Theorem is
really just starting, and this monograph only has
the beginnings of a theory. First, I believe that a version of the Master Theorem
should hold much more generally. (This is something
that John Smillie and I hope to work out together.)
Second, some recent experiments convince me that
there is a renormalization theory for this object,
grounded in real projective geometry. 
All of this is perhaps the subject of a future
work.

\subsection{Computational Issues}

I discovered all the structure
of outer billiards by experimenting with Billiard King.
Ultimately, I am trying to verify the structure I noticed
on the computer, and so one might expect there to be
some computation in the proof.
The proof here uses considerably less computation
than the proof in [{\bf S\/}], but I still use
a computer-aided proof in several places.  For
example, I use the computer to check that
various $4$ dimensional convex integral polytopes
have disjoint interiors.

To the reader who does not like computer-aided proofs
(however mild) I would like to remark that the
experimental method here has the advantage that I
checked all the results with massive
and visually-based computation.  The reader can
make the same checks, by downloading Billiard
King or else playing with the interactive online
guide to the monograph.

\subsection{Organization of the Monograph}
\label{org}

As we mentioned above, our monograph comes
in 6 parts.  
Parts I-IV comprise the core of the monograph.
In part I we define the arithmetic graph and
state its basic properties, such as the
Embedding Theorem and the Hexagrid Theorem.
Modulo these structural results,
Part I proves the results listed in
\S \ref{mainresult}, and all the results
in \S \ref{rationalresults} except
Theorem \ref{discrete}.

Part II proves the Master Picture Theorem, our main
structural result.  Part III deduces the Embedding
Theorem and the Hexagrid Theorem from the
Master Picture Theorem.  Part IV establishes various
period copying results needed for the results in
\S \ref{mainresult}.

Parts V and VI describes the close connection between the
arithmetic graph and the modular group.  In Part V,
we prove the Comet Theorem modulo technical details
that we resolve in Part VI.

Before each part of the monograph, we include an
overview of that part.

\newpage
\noindent
{\bf {\Huge Part I\/}\/}
\newline

Here is an overview of this part of the monograph.

\begin{itemize}

\item In \S 2 we establish some
basic results that allow for the definition of the arithmetic graph.
The arithmetic graph is our main object of study. We also state the
Embedding Theorem, a basic structural result about the arithmetic
graph that we prove in Part III of the monograph.

\item In \S 3 we state our main structural
result, the Hexagrid Theorem.  We then
deduce Theorems \ref{ex}, \ref{orbit structure1},
and \ref{orbit structure2} from the Hexagrid Theorem.
We prove the Hexagrid Theorem in Part III of
the monograph.

\item In \S 4 we discuss the period
copying results needed to prove
the Erratic Orbits Theorem. 
Along the way, we introduce the
inferior and superior sequences,
two basic ingredients in our overall theory.
To illustrate the connection between outer billiards and
these sequences, we state the Decomposition Theorem,
a basic structural result that helps with
the period copying. We prove the period copying
results and the Decomposition Theorem in Part IV.

\item In \S 5 we assemble the ingredients from
previous chapters and prove the Erratic Orbits Theorem.
At the end of the chapter, we prove 
Theorems \ref{dichotomy0} and \ref{density}.

\end{itemize}

We mention several conventions that we use repeatedly
throughout the monograph. Recall that 
$p/q$ is an odd rational if $pq$ is odd.  When we say
{\it odd rational\/} we mean that the odd rational
lies in $(0,1)$. On very rare occasions,
we also consider the odd rational $1/1$.
However, we never consider negative odd
rationals, or odd rationals $>1$.
Also, $A$ always stands for a kite
parameter, and we write $A=p/q$. Similarly,
$A_n$ stands for $p_n/q_n$, and
$A_+$ stand for $p_+/q_+$, etc.  Sometimes
we will mention these conventions explicitly,
and sometimes we will forget to mention them.

\newpage

\section{The Arithmetic Graph}
\label{graph_definition}

\subsection{Polygonal Outer Billiards}

Let $P$ be a polygon.
We denote the outer billiards map by $\psi'$, and the square
of the outer billiards map by $\psi=(\psi')^2$.  Our
convention is that a person walking from $p$ to $\psi'(p)$
sees the $P$ on the right side.   These maps are
defined away from a countable set of line segments in
$\R^2-P$.  This countable set of line segments is sometimes
called the {\it limit set\/}.

\begin{center}
\resizebox{!}{3.2in}{\includegraphics{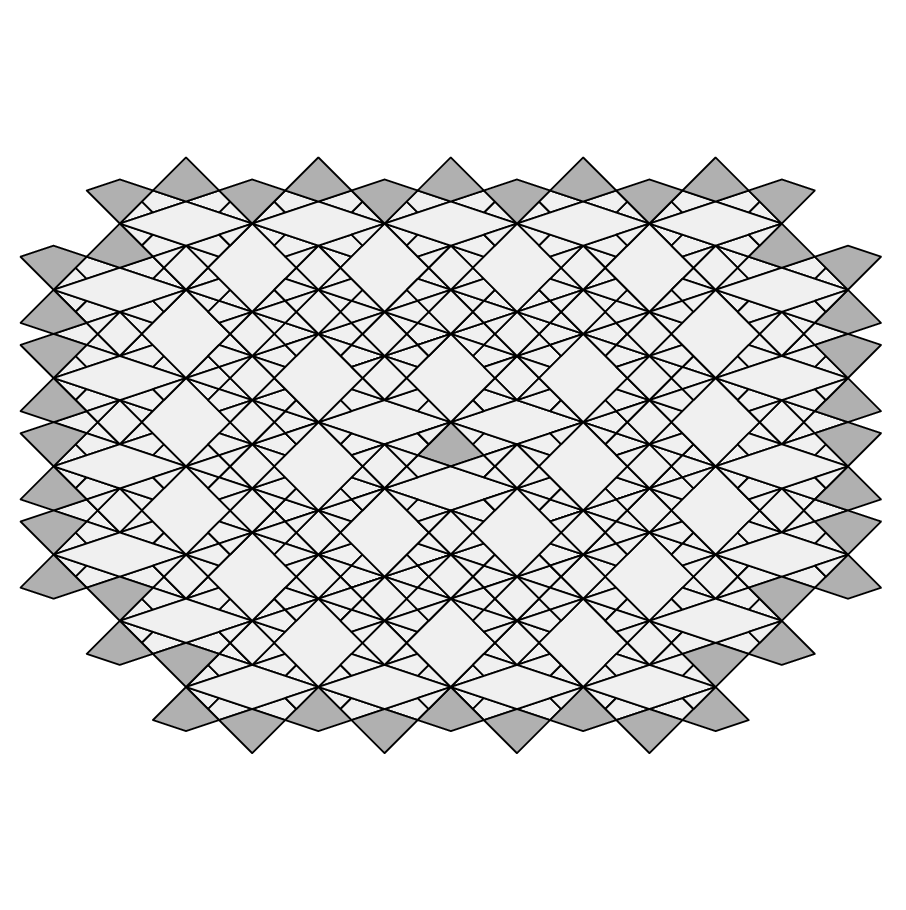}}
\newline
{\bf Figure 2.1:\/} Part of the Tiling for $K(1/3)$.
\end{center}

The result in [{\bf VS\/}], [{\bf K\/}] and
[{\bf GS\/}] states, in particular, that the orbits
for rational polygons are all periodic.  In this
case, the complement of the limit set is tiled
dy dynamically invariant convex polygons.
Figure 2.1 shows the picture for the kite $K(1/3)$.

This is the simplest tiling \footnote{Note that the picture is rotated by
$90$ degrees from our usual normalization.}
we see amongst all the kites.
We have only drawn part of the tiling.  The reader
can draw more of these pictures, and in color,
using Billiard King.
The existence of these tilings was what motivated me to
study outer billiards.  I wanted to understand how the
tiling changed with the rational parameter and saw that
the kites gave rise to highly nontrivial pictures.

\subsection{Special Orbits}
\label{definebasic}

Until the last result in this section the parameter
$A=p/q$ is rational.
Say that a {\it special interval\/} is an open
horizontal interval of length
$2/q$ centered at a point 
of the form $(a/q,b)$, with $a$ odd.  Here $a/q$ need not
be in lowest terms.

\begin{lemma}
\label{basic}
The outer billiards map is entirely defined on any special
interval, and indeed permutes the special intervals.
\end{lemma}

\startproof
We note first that the order $2$ rotations about the vertices
of $K(A)$ 
send the point $(x,y)$ to the point:
\begin{equation}
\label{trans}
(-2-x,-y); \hskip 20 pt
(-x,2-y); \hskip 20 pt
(-x,-2-y); \hskip 20 pt
(2A-x,-y).
\end{equation}

Let $\psi'$ denote the outer billiards map on $K(A)$.  The
map $\psi'$ is built out of the $4$ transformations
from Equation \ref{trans}.  
The set 
$\R \times \Z_{\rm odd\/}$ is a countable collection of lines.
Let $\Lambda \subset \R \times \Z_{\rm odd\/}$ 
denote the set of points of the form
$(2a+2bA,2c+1)$, with $a,b,c \in \Z$.
The complementary set $\Lambda^c=\R \times \Z_{\rm odd\/}-\Lambda$ 
is the union of the
special intervals.

Looking at
Equation \ref{trans}, we see that $\psi'(x) \in \Lambda^c$
provided that $x \in \Lambda^c$ and $\psi'$ is defined on
$x$.  To prove this lemma,
it suffices to show that
$\psi'$ is defined on any point of $\Lambda^c$.

To find the points of $\R \times \Z_{\rm odd\/}$ 
where $\psi'$
is not defined, we extend the sides of
$K(A)$ and intersect them with $\R \times \Z_{\rm odd\/}$.
We get $4$ families of points. 
$$(2n,2n+1); \hskip 30 pt
(2n,-2n-1); \hskip 30 pt
(2An,2n-1); \hskip 30 pt
(2An,-2n+1).
$$
Here $n \in \Z$.
Notice that all these points lie in $\Lambda$.  
\endproof

Let $\Z[A]=\Z\oplus\Z A$.
More generally, the same proof gives:

\begin{lemma}
\label{irrat}
Suppose that $A \in (0,1)$ is any number. Relative to
$K(A)$, the entire outer billiards orbit of any point
$(\alpha,n)$ is defined provided that
$\alpha \not \in 2\Z[A]$ and $n \in \Z_{\rm odd\/}$.
\end{lemma}

When $A$ is irrational, the set
$2\Z[A]$ is dense in $\R$.  However, it is always
a countable set.

\subsection{Structure of the Square Map}
\label{squarestruct}

As we mentioned in  \S \ref{preag}, we have
$\psi(p)-p=V$, where $V$ is twice a vector that points
from one vertex of $K(A)$ to another.  See Figure 1.2. There are
$12$ possilities for $V$, namely
\begin{equation}
\pm(0,4); \hskip 10 pt
\pm(2,2); \hskip 10 pt 
\pm(-2,2); \hskip 10 pt
\pm(2,2A); \hskip 10 pt
\pm(-2,2A); \hskip 10 pt
\pm(2+2A,0).
\end{equation}
These vectors are drawn, for the parameter $A=1/3$, in Figure 2.2.
The grey lines are present to guide the reader's eye.

\begin{center}
\resizebox{!}{4.1in}{\includegraphics{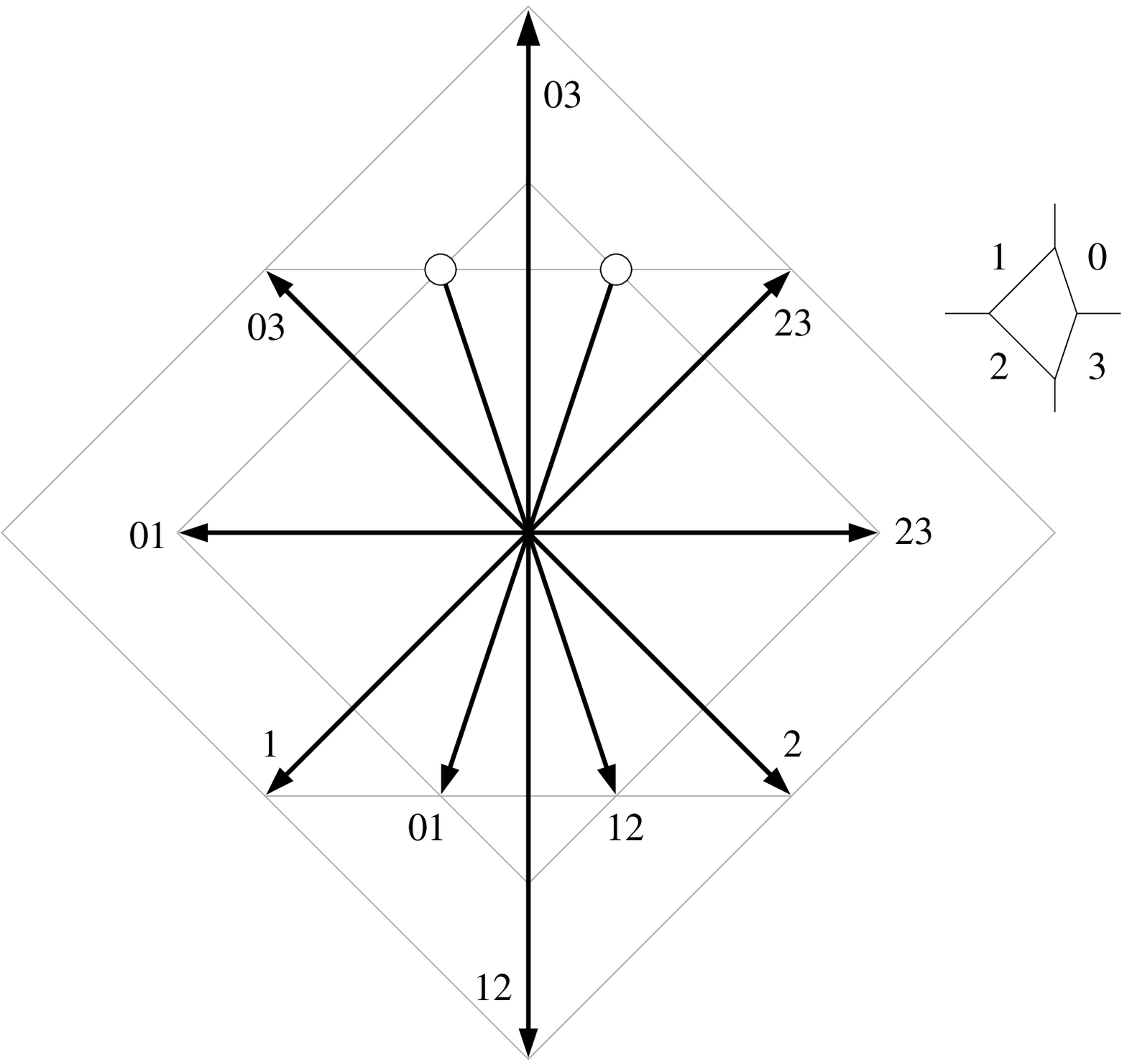}}
\newline
{\bf Figure 2.2:\/} The $12$ direction vectors
\end{center}

The labelling of the vectors works as follows.  We divide the plane
into its quadrants, according to the numbering scheme shown in
Figure 2.2.  A vector $V$ gets the label $k$ if there exists a
parameter $A$ and a point $p \in Q_k$ such that $\psi(v)-v=V$.
Here $Q_k$ is the $k$th quadrant. For instance,
$(0,-4)$ gets the labels $1$ and $2$.
 The two vectors with dots
never occur. In \S \ref{partition} we
will give a much more precise version of Figure 2.2.
For now, Figure 2.2 is sufficient for our purposes.

\subsection{The Return Lemma}
\label{return2proof}

As in the introduction let $\Xi=\R_+ \times \{-1,1\}$.

\begin{lemma}[Return]
\label{rlemma}
Let $p \in \R \times \Z_{\rm odd\/}$
be a point with a well-defined outer billiards
orbit. Then there is some $a>0$ such that
$\psi^a(p) \in \Xi$.  Likewise, there is some
$b<0$ such that $\psi^b(p) \in \Xi$.
\end{lemma}

Consider the sequence $\{\psi^k(p)\}$ for $k=1,2,3...$.
We order the quadrants of $\R^2-\{0\}$ cyclically.
Let $Q_0$ be the $(++)$ quadrant.
We include the positive $x$-axis in $Q$.  Let
$Q_{n+1}$ be the quadrant obtained by rotating
$Q_n$ clockwise by $\pi/2$.  We take indices mod $4$.

\begin{lemma}
The sequence $\{\psi^k(p)\}$ cannot
remain in a single quadrant.
\end{lemma}

\startproof
We prove this for $Q_0$.  The other cases are similar.
Let $q=\psi^k(p)$ and $r=\psi(q)$.  We write $q=(q_1,q_2)$ and
$r=(r_1,r_2)$.  Looking at Figure 2.2, we see that either
\begin{enumerate}
\item $r_2 \geq q_2+2$ and $r_1 \leq q_1$.
\item $r_1 \leq q_1-2A$.
\end{enumerate}
Moreover, Option $1$ cannot happen 
if the angle between $\overrightarrow{0r}$
and the $x$-axis is sufficiently close to $\pi/2$. Hence,
as we iterate, Option 2 occurs every so often until
the first coordinate is negative and our sequence leaves $Q_0$.
\endproof

Call $p \in \R^2-K$ a {\it bad point\/} if
$p \in Q_k$ and $\psi(p) \not \in Q_k \cup Q_{k+1}$.

\begin{lemma}
If $q$ is bad, then either $q$ or $\psi(q)$ lies in $\Xi$.
\end{lemma}

\startproof
Let $r=\psi(q)$.   If
$q$ is bad then $q_2$ and $r_2$ have opposite signs.
if $q_2-r_2=\pm 2$ then $q_2=\pm 1$ and $r_2=\pm 1$.
A short case-by-case analysis shows that this forces
$q_1$ and $r_1$ to have opposite signs.  
The other possibility is that $q_2-r_2=4$. But then
$r-q = \pm (0,4)$.  A routine 
case-by-case analysis shows that
$r-q=(0,4)$ only if $q_1>0$ and
$r-q=(0,-4)$ only if $q_1<0$. But $q$ is not
bad in these cases.
\endproof

If the Return Lemma is false, then our sequence is entirely good.
But then we must have some $k$ such that $\psi^k(p) \in Q_3$ and
$\psi^{k+1}(p) \in Q_0$.  Since the second coordinates differ
by at most $4$, we must have either $\psi^k(p) \in \Xi$
or $\psi^{k+1}(p) \in \Xi$. This proves the first
statement.  The second statement
follows from the first statement and symmetry.

\subsection{The Return Map}
\label{basic0}

The Return Lemma implies that the
{\it first return map\/} $\Psi: \Xi \to \Xi$ is
well defined on any point with an outer billiards
orbit.  This includes the set
$$(\R_+-2\Z[A]) \times \{-1,1\},$$
as we saw in Lemma \ref{irrat}.

Given the nature of the
maps in Equation \ref{trans} comprising $\psi$, we see that
$$
\Psi(p)-(p) \in 2\Z[A] \times \{-2,0,2\}.
$$

In Part II, we will prove our main structural result about
the first return map, namely the Master Picture Theorem.
We will also prove the Pinwheel Lemma, in Part II.  Combining
these two results, we have a much stronger
result about the nature of the first return map:

\begin{equation}
\label{shortreturn}
\Psi(p)-(p)= 2(A\epsilon_1+ \epsilon_2,\epsilon_3); \hskip 10 pt
\epsilon_j \in \{-1,0,1\}; \hskip 10 pt \sum_{j=1}^3 \epsilon_j \equiv 0 \hskip 5 pt 
{\rm mod\/} \hskip 5 pt 2.
\end{equation}

\noindent
{\bf Remarks:\/} \newline
(i) Some notion of the return map is also used in
[{\bf K\/}] and [{\bf GS\/}].  This is quite a natural
object to study. \newline
(ii)
We can at least roughly
explain the first statement of Equation \ref{shortreturn}
in an elementary way.  At least far from the origin,
the square outer billiards orbit circulates around the kite in
such a way as to nearly make an octagon with $4$-fold
symmetry.  Compare Figure 11.3.
  The {\it return pair\/} $(\epsilon_1(p),\epsilon_2(p))$
essentially measures the {\it approximation error\/} between
the true orbit and the closed octagon. \newline
(iii)
On a nuts-and-bolts level, this monograph concerns how to determine
$(\epsilon_1(p),\epsilon_2(p))$ as a function
of $p \in \Xi$.   (The pair
$(\epsilon_1,\epsilon_2)$ and the parity condition determine $\epsilon_3$.)
I like to tell people
that this book is really about the infinite accumulation
of small errors.
\newline
(iv) Reflection in the $x$-axis conjugates the map $\psi$ to
the map $\psi^{-1}$.  Thus, once we understand the orbit of
the point $(x,1)$ we automatically understand the orbit
of the point $(x,-1)$.  Put another way, the unordered pair of
return points $\{\Psi(p),\Psi^{-1}(p)\}$ for $p=(x,\pm 1)$
only depends on $x$.
 \newline

\subsection{The Arithmetic Graph}
\label{ag}

Recall that $\Xi=\R_+ \times \{-1,1\}$.
Define $M=M_{A,\alpha}: \R \times \{-1,1\}$
by
\begin{equation}
\label{map}
M_{A,\alpha}(m,n)=\Big(2Am + 2m + 2\alpha,(-1)^{m+n+1}\Big)
\end{equation}
The second coordinate of $M$ is either $1$ or $-1$ depending
on the parity of $m+n$.  This definition is adapted to the parity
condition in Equation \ref{shortreturn}.
We call $M$ a {\it fundamental map\/}.  Each choice of
$\alpha$ gives a different map.

When $A$ is irrational, $M$ is injective.
In the rational case, $M$ is injective on any
disk of radius $q$.
Given $p_1,p_2 \in \Z^2$, we write
$p_1 \to p_2$ iff the following holds.
\begin{itemize}
\item $\zeta_j = M(p_j) \in \Xi$.
\item $\Psi(\zeta_1)=\zeta_2$.
\item $\|p_1-p_2\| \leq \sqrt 2$.
\end{itemize}
The third condition is only relevant in the rational case.
See Equation \ref{shortreturn}.
Our construction gives a directed graph with vertices in $\Z^2$.
We call this graph the {\it arithmetic graph\/} and denote it by
$\widehat \Gamma_{\alpha}(A)$.

When $A=p/q$, any choice of
 $\alpha \in (0,2/q)$
gives the same result. This is a consequence of
Lemma \ref{basic}.  
To simplify our formulas,  we choose $\alpha=0_+$, where $0_+$ is
an infinitesimally small positive number.  The reader who does
not like infinitesimally small positive numbers can
take
$$\alpha=\exp(-(\exp(q))).$$

  When we write our
formulas, we usually take $\alpha=0$, but we always use the
convention that the lattice point $(m,n)$ tracks the orbits
just to the right of the points $(2Am+2n,\pm 1)$.   
With this convention, we have
\begin{equation}
\label{funm}
\widehat \Gamma(p/q)=\widehat \Gamma_{0_+)}(p/q); \hskip 30 pt
M(m,n)=\Big(2(p/q)m+2n,(-1)^{m+n+1}\Big)
\end{equation}
We say that the
{\it baseline\/} of $\widehat \Gamma(A)$ is the line
$M^{-1}(0)$.  We think of the baseline
essentially as the line $L$ of slope $-A$ through the
origin.  However, we really want to think of the
baseline as lying infinitesimally beaneath $L$, so that
the entire arithmetic graph lies above the baseline.

we will prove the following result.

\begin{theorem}[Embedding]
Any well-defined arithmetic graph is the disjoint union of
embedded polygons and bi-infinite embedded polygonal curves.
\end{theorem}

\begin{center}
\resizebox{!}{4.2in}{\includegraphics{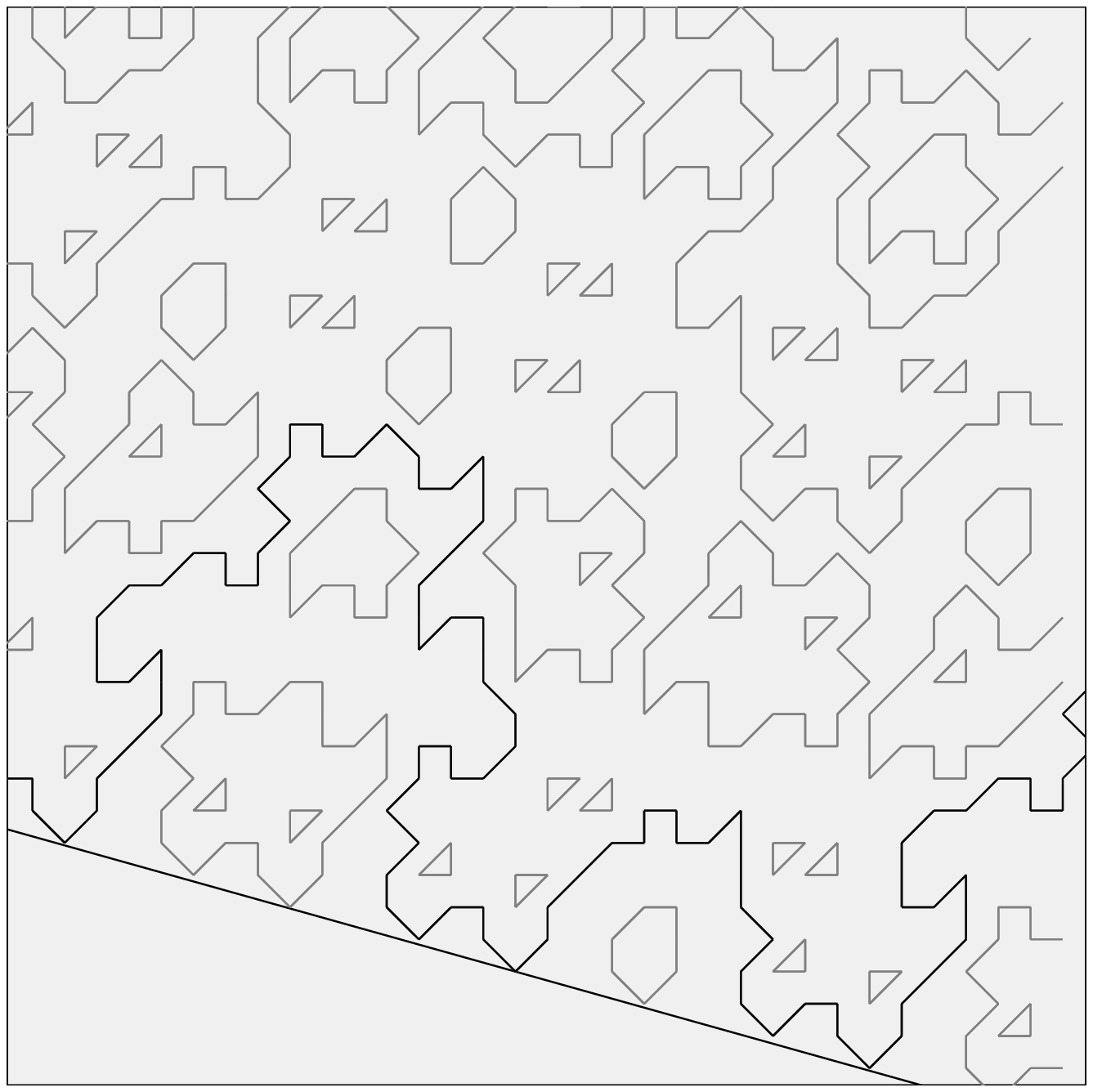}}
\newline
{\bf Figure 2.3:\/} Some of $\widehat \Gamma(7/25)$, with
$\Gamma(7/25)$ in black
\end{center}

\noindent
{\bf Remark:\/}
In the arithmetic graph, there are some lattice points
having no edges emanating from them.  These isolated points
correspond to points where the return map is the identity
and hence the orbit is periodic in the simplest possible
way.  We usually ignore these trivial components.
\newline

We are mainly interested in the component
of $\widehat \Gamma$ that
contains $(0,0)$.  We denote this component by $\Gamma$.
In the rational case, $\Gamma(p/q)$ encodes
the structure of the orbit $O_2(1/q,-1)$.
The orbit $O_2(1/q,1)$, the subject of
Theorems \ref{ex} and \ref{discrete}, 
is conjugate to $O_2(1/q,-1)$ {\it via\/}
reflection in the $x$-axis.
As we have said several times above, $\Gamma(p/q)$
really tracks the orbit of the special interval
bounded by $(0,-1)$ and $(2/q,-1)$.

\subsection{The Continuity Principle}
\label{contp}

Given two compact subsets $K_1,K_2 \subset \R^2$, we define
$d(K_1,K_2)$ to be the infimal $\epsilon>0$ such that
$K_1$ is contained in the $\epsilon$-tubular neighborhood
of $K_2$, and {\it vice versa\/}.  The function
$d(K_1,K_2)$ is known as the {\it Hausdorff metric\/}.
A sequence $\{C_n\}$ of closed subsets of $\R^2$ is
said to {\it Hausdorff converge\/} to $C \subset \R^2$ if
$d(C_n \cap K,C \cap K) \to 0$ for every compact subset
$K \subset \R^2$.

In the cases of interest to us, $C_n$ will always be an
arc of an arithmetic graph that contains $(0,0)$.
In this case, the Hausdorff convergence
has a simple meaning.  $\{C_n\}$ converges to $C$ if
and only if the following property holds true.
For any $N$ there some $N'$ such that
$n>N'$ implies that the first $N$ steps of
$C_n$ away from $(0,0)$ in either direction agree with the
corresponding steps of $C$.

Given a parameter $A \in (0,1)$ and a point
$\zeta \in \Xi$, we say that
a pair $(A,\zeta)$ is $N$-{\it defined\/}
if the first $N$ iterates of the outer billiards
map of $\zeta$ are defined relative to $A$,
in both directions.  We let
$\Gamma(A,\zeta)$ be as much of the
arithmetic graph as is defined.
We call $\Gamma(A,\zeta)$ a {\it partial arithmetic
graph\/}.

\begin{lemma}[Continuity Principle]
Let $\{\zeta_n\} \in \Xi$ converge 
to $\zeta \in \Xi$.  Let $\{A_n\}$ converge to $A$.
Suppose the orbit of $\zeta$ is defined relative to
$A$.  Then for any $N$ there is some $N'$ such that
$n>N'$ implies that $(\zeta_n,A_n)$ is $N$-defined.
The corresponding sequence $\{\Gamma(A_n,\zeta_n)\}$ of partially defined
arithmetic graphs Hausdorff converges to $\Gamma(A,\zeta)$.
\end{lemma}

\startproof
Let $\psi'_n$ be the outer billiards map relative
to $A_n$.  Let $\psi'$ be the outer billiards map
defined relative to $A$. If
$p_n \to p$ and
$\psi'$ is defined at $p$
then $\psi'_n$ is defined at $p_n$ for $n$ sufficiently large and
$\psi'_n(p_n) \to \psi(p)$.
This follows from the fact that $K(A_n) \to K(A)$ and
from the fact that a piecewise isometric map is
defined and continuous in open sets.  Our
continuity principle now follows from induction.
\endproof

In case the orbit of $\zeta_n$ relative
to $A_n$ is already well defined, the partial
arithmetic graph is the same as one component
of the ordinary arithmetic graph.  In this case,
we can state the Continuity Principle more simply.

\begin{corollary}
Let $\{\zeta_n\} \in \Xi$ converge 
to $\zeta \in \Xi$.  Let $\{A_n\}$ converge to $A$.
Suppose the orbit of $\zeta$ is defined relative to
$A$ and the orbit of $\zeta_n$ is defined relative to $A_n$
for all $n$.  Then
$\{\Gamma(A_n,\zeta_n)\}$ Hausdorff converges to $\Gamma(A,\zeta)$.
\end{corollary}

We will have occasion to use both versions in our arguments.

\subsection{Low Vertices and Parity}
\label{low vertex}

Let $A$ be any kite parameter. 
We define the {\it parity\/} of a low vertex $(m,n)$ to
be the parity of $m+n$.  
Here we explain the structure of the arithmetic
graph at low vertices. Our answer will be given
in terms of a kind of phase portrait.  Given a
point $(x,A) \in (0,2) \times (0,1)$, we have
\begin{equation}
\Psi^{\pm 1}(x,-1)=(x,-1)+2(\epsilon_1^{\pm}A+\epsilon_2^{\pm},\epsilon_3^{\pm}).
\end{equation}
For the point $(x,t)$ we associate the directed graph
$$(\epsilon_1^-,\epsilon_2^-) \to (0,0) \to
(\epsilon_1^+,\epsilon_2^+).$$
This gives a local picture of the arithmetic
at the low vertex $(m,n)$ such that
$M_A(m,n)=(x,-1)$. If
$M_A(m,n)=(x,1)$ then we get the local
picture by reversing the edges.
Figure 2.2 shows the final result.
 The grey edges in the picture,
present for reference,
connect $(0,0)$ to $(0,-1)$.  The grey
triangle represents the places where the return
map is the identity.

\begin{center}
\resizebox{!}{2.1in}{\includegraphics{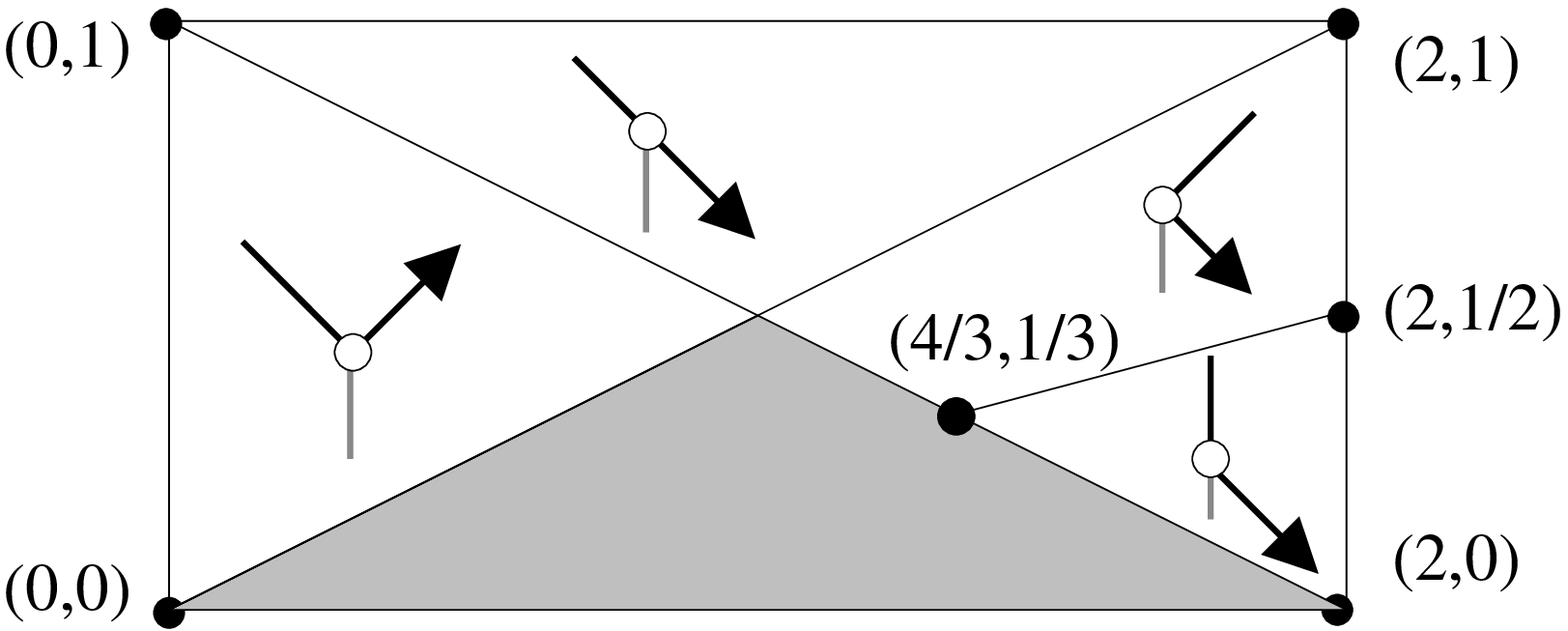}}
\newline
{\bf Figure 2.4:\/} Low vertex phase portrait
\end{center}

\noindent
{\bf Example:\/}
Relative to
$A=1/3$, the vertex $-7,3$ is a low vertex.
We compute that $$M_{1/3}(A)=(4/3+\alpha,-1).$$ 
Here $\alpha$ is an infinitesimally small positive
number.  To see the
local picture of the arithmetic graph $\Gamma(1/3)$ 
at $(-7,3)$ we observe that the point
$p=(4/3+\alpha,1/3)$ lies
infinitesimally to the right of the point
$(4/3,1/3)$.   Hence
$(\epsilon_1^-,\epsilon_2^-)=(0,1)$ and
$(\epsilon_1^+,\epsilon_2^+)=(1,-1)$.
\newline

In principle, one can derive Figure 2.4 by hand.
We will explain how to derive it in \S \ref{computations},
as a corollary of the Master Picture Theorem.

\begin{lemma}
\label{parity}
No component 
of $\widehat \Gamma(A)$ contains low vertices of
both parities.
\end{lemma}

\startproof
Recall that $\widehat \Gamma$ is an oriented
graph.  
If $v$ is a nontrivial low vertex of $\widehat \Gamma$
we can say whether $\widehat \Gamma$ is {\it left travelling\/}
at $v$ or {\it right travelling\/}.  The definition is this:
As we travel along the orientation and pass through $v$,
the line segment connecting $v$ to $v-(0,1)$ either lies
on our left or our right.  This gives the name to our
definition.  Figure 2.5 shows examples in each case.
Our convention is that $\widehat \Gamma$ is right
oriented at $(0,0)$.

\begin{center}
\resizebox{!}{1.4in}{\includegraphics{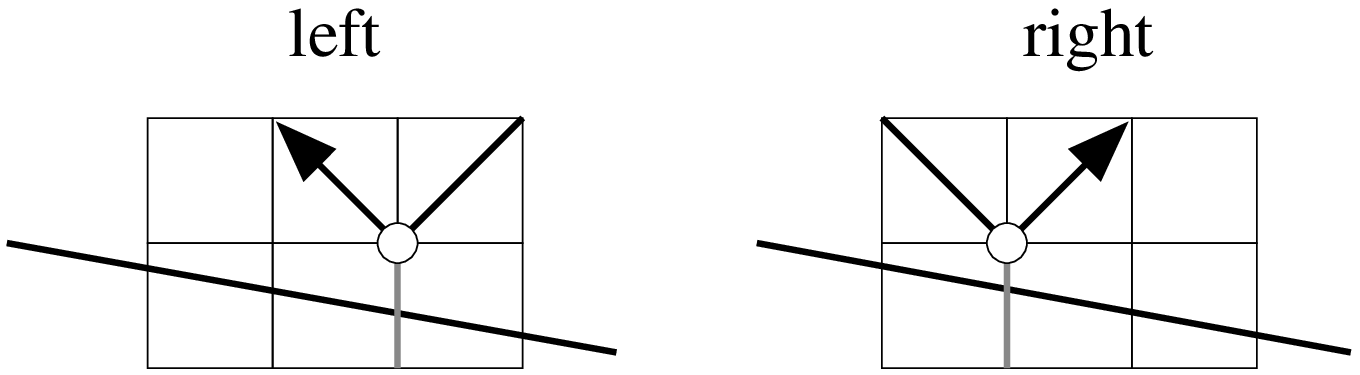}}
\newline
{\bf Figure 2.5:\/} Left travelling and right travelling.
\end{center} 

A component of $\Gamma$ cannot right-travel
at one low vertex and left travel at another.
Figure 2.6 shows 
the problem.  The curve $\gamma$ would
create a pocket for itself.  $\gamma$ could not
escape from this pocket because it must stay
above the baseline.  The low vertices of
$\gamma$ serve as barriers.
  Travelling into the pocket,
$\gamma$ would have only a finite number of
steps before it would have to cross itself.
(Recall that $\gamma$ is either a closed
polygon or an infinite periodic arc.)  
But then we contradict the Embedding Theorem.

\begin{center}
\resizebox{!}{1.7in}{\includegraphics{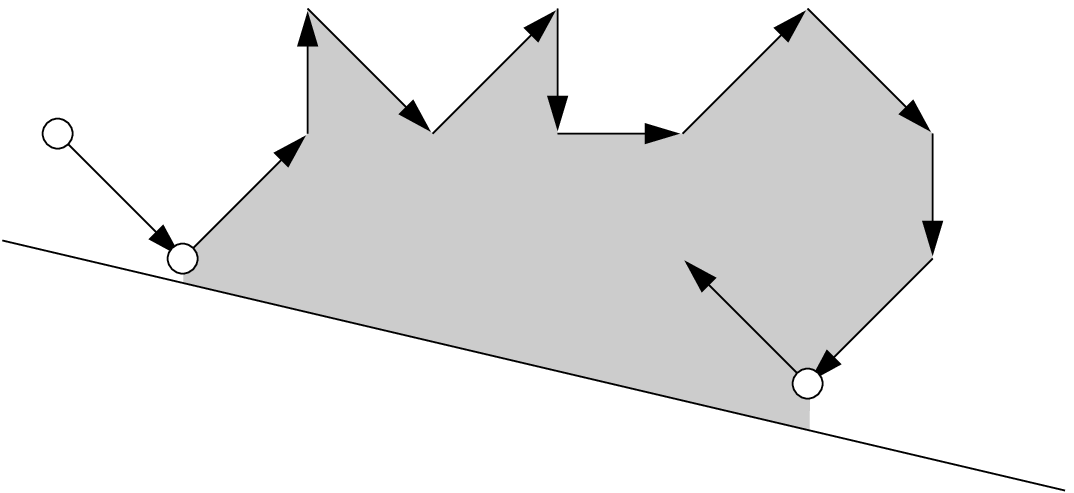}}
\newline
{\bf Figure 2.6:\/} $\gamma$ travels into a pocket.
\end{center}

To finish our proof, we just have to show that
a component of $\widehat \Gamma$ 
right-travels at a low vertex  $v$ if and only if $v$ has even
parity.   We will show that a component of
$\widehat \Gamma$ always right-travels at
low vertices of even parity.  Let us explain why
this suffices.  Recall that the
fundamental map $M$ maps vertices of
even parity to $\R_+ \times \{-1\}$ and
vertices of odd parity to $\R \times \{1\}$.
Also, recall that reflection in the
$x$-axis conjugates the return map
$\Psi$ to $\Psi^{-1}$.  It follows from
this symmetry that $\widehat \Gamma$
left-travels at all low vertices of odd
parity if and only if $\widehat \Gamma$
right-travels at all vertices of even parity.
But a glance at Figure 2.4 shows that
$\widehat \Gamma$ right travels at all
vertices of even parity.  
The grey line segment always lies on the right.
\endproof

\begin{corollary}
\label{exclude}
Let $A$ be any rational parameter.
Let $\xi_{\pm}$ be any point
in $(0,2) \times \{\pm 1\}$ with
a well defined orbit relative to $A$. 
 Then the two
orbits $O_2(\xi_+)$ and
$O_2(\xi_-)$ are disjoint.
\end{corollary}

\startproof
Our orbits are either disjoint or identical.
By perturbing $\xi_{\pm}$ slightly we arrange that
$\xi_{\pm}=M(m_{\pm},n_{\pm})$.   Here
$(m_+,n_+)$ has odd parity and
$(m_-,n_-)$ has even parity.
If we have $O_2(\xi_+)=O_2(\xi_-)$ then one and
the same component of
$\widehat \Gamma(A)$ contains both
$(m_+,n_+)$ and $(m_-,n_-)$.  But
this contradicts our previous result.
\endproof

\newpage

\section{The Hexagrid Theorem}
\label{hexagrid0}

\subsection{The Arithmetic Kite}
\label{arithkite}

In this section we describe a certain quadrilateral, which we
call the {\it arithmetic kite\/}.  This object is meant to
``live'' in the same plane as the arithmetic graph.
  The diagonals and sides of
this quadrilateral define
$6$ special directions.
 In the next section we describe
a grid made from $6$ infinite families of parallel lines,
based on these $6$ directions.

\begin{center}
\resizebox{!}{2in}{\includegraphics{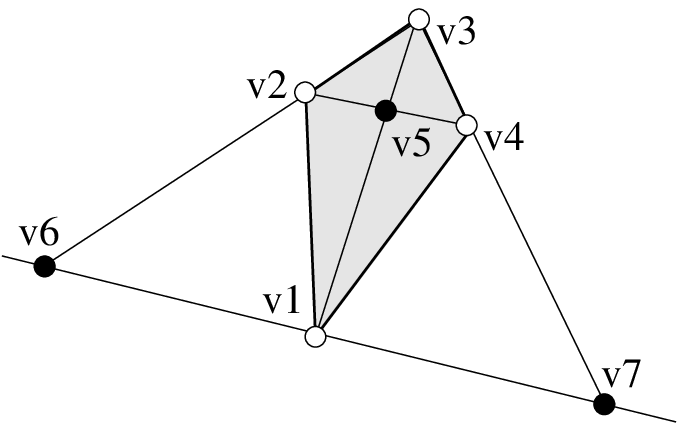}}
\newline
{\bf Figure 3.1:\/} The arithmetic kite
\end{center}

Let $A=p/q$.
Figure 3.1 shows a schematic picture of ${\cal K\/}(A)$.  The vertices
are given by the equations.

\begin{enumerate}
\item $v_1=(0,0)$.
\item $v_2=\frac{1}{2}(0,p+q)$.
\item $v_3=\frac{1}{2q}(2pq,(p+q)^2-2p^2)$.
\item $v_4=\frac{1}{2(p+q)}(4pq,(p+q)^2-4p^2)$.
\item $v_5=\frac{1}{2(p+q)}(2pq,(p+q)^2-2p^2)$.
\item $-v_6=v_7=(q,-p)$.
\end{enumerate}

A short calculation, which we omit, shows that $K(A)$
and ${\cal K\/}(A)$ are actually affinely equivalent. ${\cal K\/}(A)$
does not have Euclidean bilateral symmetry, but it does
have affine bilateral symmetry. 
We especially single out the vectors $V=v_7$ and $W-=v_5$. That is,
\begin{equation}
\label{boxvector}
\label{boxvectors}
V=(q,-p); \hskip 30 pt 
W=\bigg(\frac{pq}{p+q},\frac{pq}{p+q}+\frac{q-p}{2}\bigg).
\end{equation}

The {\it hexagrid\/} $G(A)$ consists of two interacting grids,
which we call the {\it room grid\/} $RG(A)$ and the {\it door grid\/} $DG(A)$.
\newline
\newline
{\bf Room Grid:\/}
When $A$ is an odd rational, $RG(A)$ consists of the lines obtained by extending
the diagonals
of ${\cal K\/}(A)$ and then taking the orbit
under the lattice $\Z[V/2,W]$.  These are the black lines
in Figure 3.2.
In case $A$ is an even rational, we would make the same
definition, but use the lattice $\Z[V,2W]$ instead.
\newline
\newline
{\bf Door Grid:\/}
The {\it door grid\/} $DG(A)$ is the same for both
even and odd rationals. It is obtained by extending
the sides of ${\cal K\/}(A)$ and then taking their orbit under the
one dimensional lattice $\Z[V]$.  These are the grey lines in Figure 3.2.

\begin{center}
\resizebox{!}{3.2in}{\includegraphics{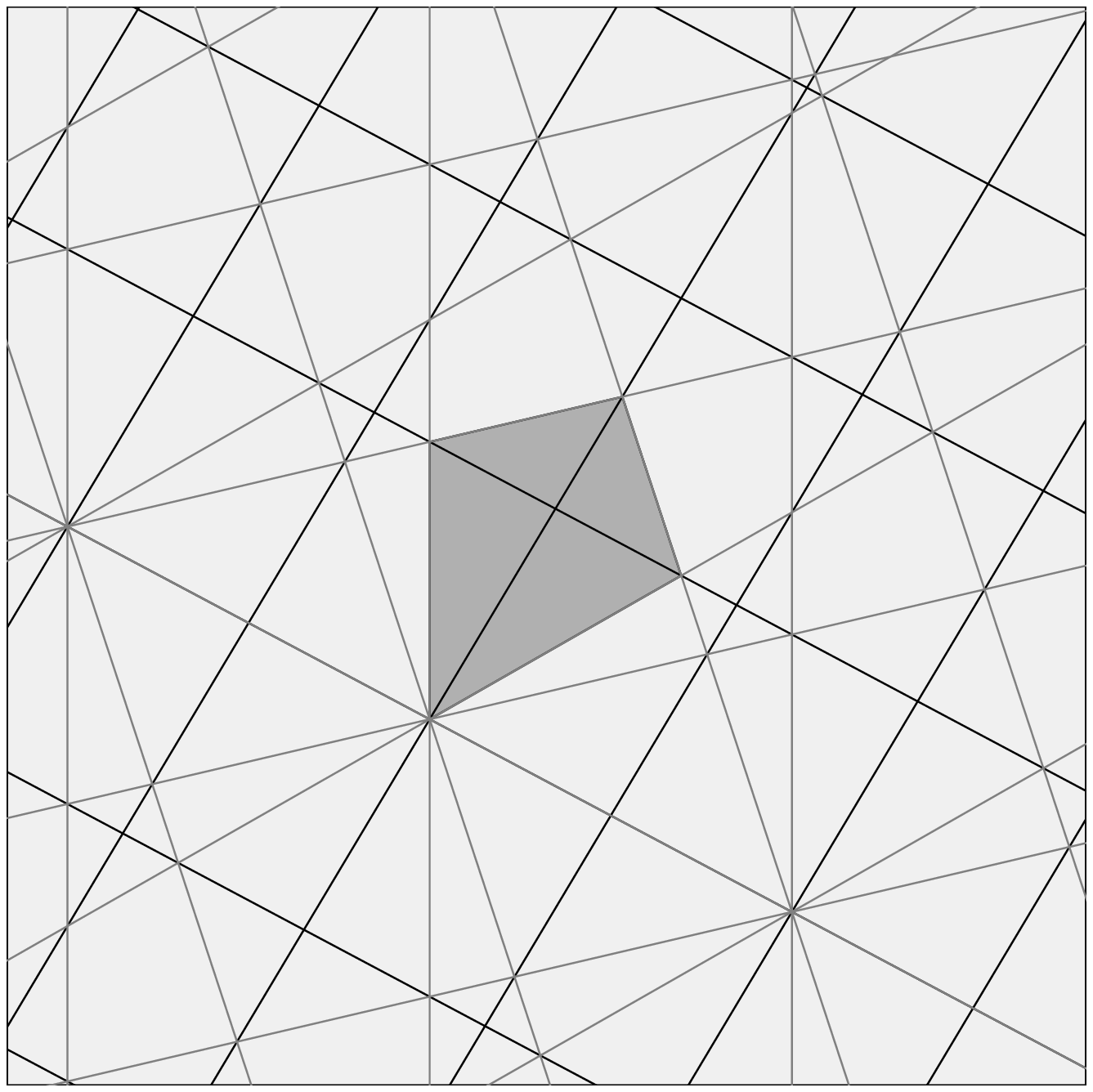}}
\newline
{\bf Figure 3.2:\/} $G(25/47)$. and ${\cal K\/}(25/47)$.
\end{center}

\subsection{The Hexagrid Theorem}
\label{roomproof}
\label{hex}

The Hexagrid Theorem relates two kinds of objects,
{\it wall crossings\/} and {\it doors\/}.  Informally,
the Hexagrid Theorem says that the arithmetic
graph only crosses a wall at a door. Here are
formal definitions.
\newline
\newline
{\bf Rooms and Walls:\/}
$RG(A)$ divides $\R^2$ into different
connected components which we call {\it rooms\/}.
Say that a {\it wall\/} is the line segment of positive
slope that divides two adjacent rooms.
\newline
\newline
{\bf Doors:\/} 
When $p/q$ is odd,
we say that a {\it door\/} is a point of intersection between a
wall of $RG(A)$ and a line of $DG(A)$.  When $p/q$ is
even, we make the same definition, except that we
exclude crossing points of the form $(x,y)$, where
$y$ is a half-integer.  
Every door is a triple point, and every wall has one door.
The first coordinate of a door is always an integer.
(See Lemma \ref{doorint}.)
In exceptional cases -- when the second coordinate
is also an integer -- the door lies in the
corner of the room.  In this case, we associate
the door to both walls containing it.
The door $(0,0)$ has this
property. 
\newline
\newline
{\bf Crossing Cells:\/}
Say that an edge $e$ of $\widehat \Gamma$ {\it crosses a wall\/}
if $e$ intersects a wall at an interior point.
Say that a union of two incident edges of $\Gamma$ {\it crosses a wall\/}
if the common vertex lies on a wall, and the two edges  point to opposite sides of the wall.
The point $(0,0)$ has this property.
We say that a {\it crossing cell\/} is either an edge or a union of
two edges that crosses a wall in the manner just described.
For instance $(-1,1)\to(0,0) \to (1,1)$ is a crossing cell
for any $A \in (0,1)$.
\newline

In Part III of the monograph we will prove the following
result. Let $\underline y$ denote the greatest integer less than $y$.

\begin{theorem}[Hexagrid]
Let $A \in (0,1)$ be rational.
\begin{enumerate} 
\item $\widehat \Gamma(A)$ never crosses a floor of $RG(A)$.
Any edges of $\widehat \Gamma(A)$ incident to a vertex contained
on a floor rise above that floor (rather than below it.)
\item There is a bijection between the set of doors
and the set of crossing cells.  If $y$ is not an integer, then
the crossing cell corresponding to the
door $(m,y)$ contains 
$(m,\underline y) \in \Z^2$.
If $y$ is an integer, then $(x,y)$ corresponds to $2$ doors.
One of the corresponding crossing cells contains $(x,y)$ and
the other one contains $(x,y-1)$.
\end{enumerate}
\end{theorem}

\noindent
{\bf Remark:\/} We really only care about the odd case of
the hexagrid theorem.  We include the even case for the
sake of completeness.
\newline

Figure 3.3 illustrates the Hexagrid Theorem
for $p/q=25/47$.  
We will explain the shaded parallelogram $R(25/47)$ in the
next section.  We have only drawn the
unstable components in Figure 3.3. 
The reader can see much better pictures
of the Hexagrid Theorem using either
Billiard King or our interactive
guide to the monograph.  (The interactive
guide only shows the odd case, but Billiard
King also shows the even case.)

\begin{center}
\resizebox{!}{5.3in}{\includegraphics{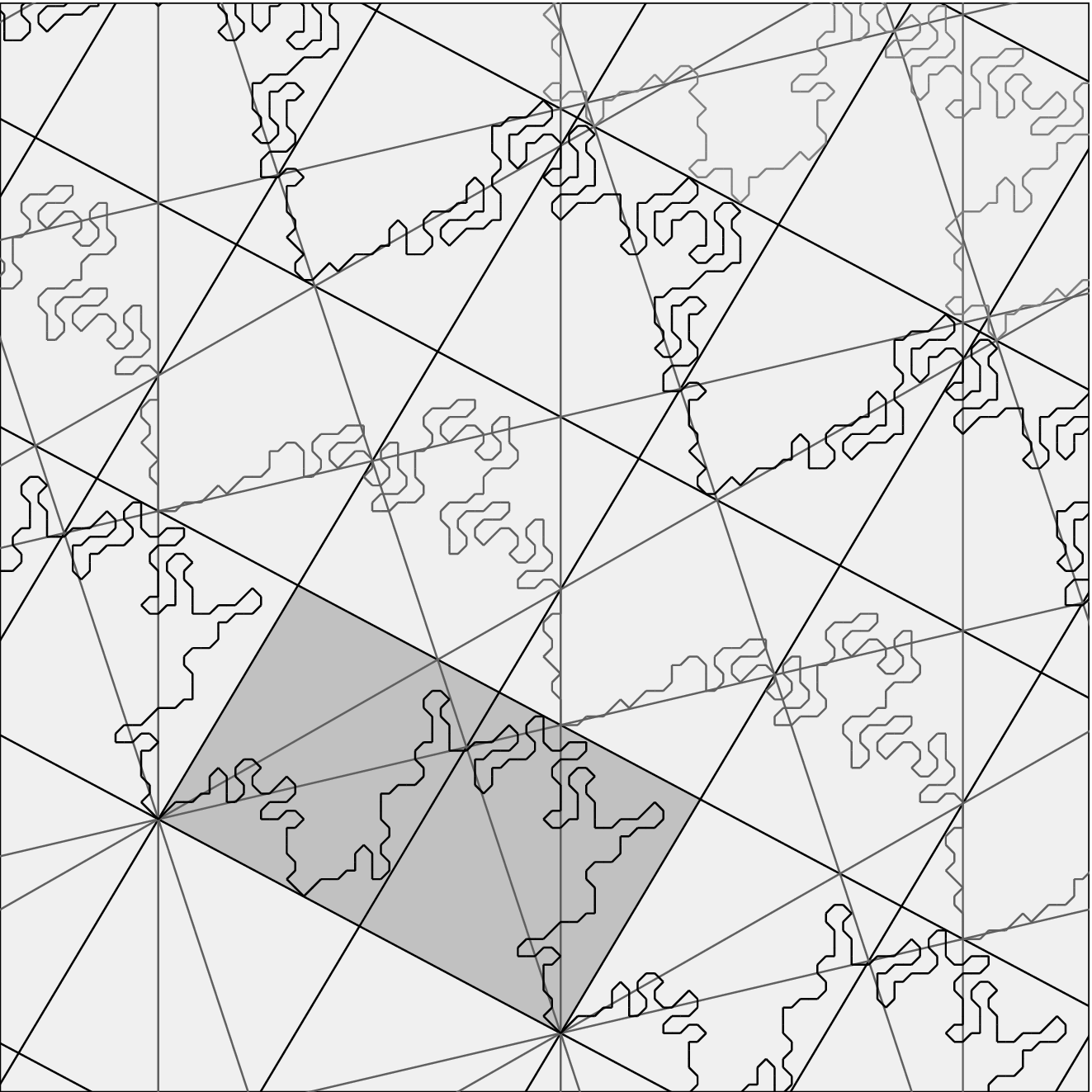}}
\newline
{\bf Figure 3.3:\/} $G(25/47)$, $R(25/47)$, and some of $\widehat \Gamma(25/47)$.
\end{center}

\subsection{The Room Lemma}
\label{roomend}

Let $R(p/q)$ denote the parallelogram whose vertices are
\begin{equation}
\label{rectdef}
(0,0); \hskip 30 pt V; \hskip 30 pt W; \hskip 30 pt V+W.
\end{equation}
Here $V$ and $W$ are as in Equation \ref{boxvectors}.
See Figure 3.3.
We also define
\begin{equation}
\label{maindoor}
d_0=(x,\underline y); \hskip 30 pt
x=\frac{p+q}{2}; \hskip 30 pt
y=\frac{q^2-p^2}{4q}
\end{equation}
$d_0$ lies within $1$ vertical unit of the centerline of $R(p/q)$,
above the center. $d_0$ is just below the door contained inside the shaded
parallelogram in Figure 3.3.  Figure 3.1 is an enlargement of
this parallelogram. 

\begin{lemma}[Room]
\label{room}
$\Gamma(p/q)$ is an open polygonal curve.  One period of
$\Gamma(p/q)$ connects $(0,0)$ to $d_0$ to $(q,-p)$. This
period is contained in $R(p/q)$.
\end{lemma}

\startproof
First of all, for any value of $A$, it is easy to check that
$\Gamma(A)$ contains the arc $(-1,1) \to (0,0) \to (1,1)$.
This is to say that $\Gamma(p/q)$ enters $R(p/q)$ from
the left at $(0,0)$.
Now, $R(p/q)$ is the union of two adjacent rooms, $R_1$ and $R_2$.
Note that $(0,0)$ is the only door on the left wall of $R_1$ and
$(x,y)$ is the only door on the wall separating $R_1$ and $R_2$,
and $(q,-p)$ is the only door on the right wall of $R_2$.
Here $(x,y)$ is as in Equation \ref{maindoor}.  From the
Hexagrid Theorem and the Embedding Theorem,
$\Gamma(p/q)$ must connect
$(0,0)$ to $d_0$ to $(q,-p)$. 
The arithmetic graph $\widehat \Gamma(p/q)$ is invariant
under translation by $(q,-p)$, and so the whole picture
repeats endlessly to the left and the right of $R(p/q)$.
Hence $\Gamma(p/q)$ is an open polygonal curve.
\endproof

We remark that we did not really need the Embedding
Theorem in our proof above. 
All we require is that
$\Gamma(p/q)$ cannot backtrack as we travel from
one corner of $R(p/q)$ to the other. Lemma \ref{val2}
below gives a self-contained proof of what we need.

\begin{lemma}
\label{val2}
$\Gamma(p/q)$ has valence $2$ at every vertex.
\end{lemma}

\startproof
As in our proof of the Room Lemma, $\Gamma(p/q)$ has
valence $2$ at $(0,0)$.  But $\Gamma(p/q)$
describes the forward orbit of $p=(1/q,1)$ under $\Psi$.
If some vertex of $\Gamma(p/q)$ has valence $1$ then
$\Psi$ has order $2$ when evaluated at the corresponding point.
But then $\Psi$ has order $2$ when evaluated at $v$.  But
then $\Gamma(p/q)$ has valence $1$ at $(0,0)$. This
is a contradiction.
\endproof

\subsection{Proof of Theorems \ref{ex} and \ref{orbit structure1}}
\label{discussK}

The bounds in Theorem \ref{orbit structure1} imply
the upper bound in Theorem \ref{ex}.  First we
establish the lower bound in Theorem \ref{ex}.
Suppose that $p/q$ is an odd rational.
Let $M_1$ be the first coordinate for the
fundamental map associated to $p/q$.
We compute that $M_1(d_0)>(p+q)/2$, at least when $p>1$. 
Technically, $\Gamma(p/q)$ describes 
$O_2(1/q,-1)$, but the two orbits
$O_2(1/q,1)$ and $O_2(1/q,-1)$ are conjugate by reflection
in the $x$-axis.

Now suppose that $p/q$ is even.
Referring to the plane containing the arithmetic graph,
let $S_0$ be the line segment connecting the origin to
$v_3$, the very tip of the arithmetic kite.
Then $S_0$ is bounded by two consecutive doors on $L_0$.
The bottom endpoint of $S_0$ is $(0,0)$, one of the
vertices of $\Gamma(0,0)$.  We know already that
$\Gamma(p/q)$ is a closed polygon.  By the hexagrid
Theorem, $\Gamma(p/q)$, cannot cross $S_0$ except within
$1$ unit of the door $v_3$.  Hence, $\Gamma(p/q)$ must
engulf all but the top $1$ unit of $S_0$.

Essentially the same calculation as in the odd
case now shows that $\Gamma(p/q)$ rises up
at least $(p+q)$ units from the baseline
when $p>1$.  When $p=1$ the same result holds,
but the calculation is a bit harder.  The
reason why we get an extra factor of $2$ in the
even case is that $v_3$ is twice as far
from the baseline as is the door near $d_0$. See
Equation \ref{maindoor}.

First suppose that $p/q$ is odd.
Let $M_1$ be the first coordinate of the fundamental
map associated to $p/q$.
Since $p$ and $q$ are
relatively prime, we can realize any integer
as an integer combination of $p$ and $q$.  From this we
see that every point of the form $s/q$, with $s$ odd,
lies in the image of $M_1$.  Hence, some point of
$\Z^2$, above the baseline of $\widehat \Gamma(p/q)$, corresponds to
the orbit of either $(s/q,1)$ or $(s/q,-1)$.

Let the {\it floor grid\/} denote the lines of negative
slope in the room grid.  These lines all have slope $-p/q$.
The $k$th line $L_k$ of the floor grid contains the point
$$\zeta_k=\bigg(0,\frac{k(p+q)}{2}\bigg).$$
Modulo translation by $V$, the
point $\zeta_k$ is the only lattice point on
$L_k$.  Statement 1 of the Hexagrid Theorem
contains that statement that
the edges of $\Gamma$ incident to $\zeta_k$
lie between $L_k$ and $L_{k+1}$ (rather than between
$L_{k-1}$ and $L_k$). 

We compute that 
$$M_1(\zeta_k)=k(p+q)+\frac{1}{q}.$$
For all lattice points $(m,n)$ between
$L_k$ and $L_{k+1}$ we therefore have
\begin{equation}
\label{confine1}
M_1(m,n) \in I_k,
\end{equation}
the interval from Theorem \ref{orbit structure1}.
Theorem \ref{orbit structure1} now follows from
Equation \ref{confine1}, Statement 1 of the
Hexagrid Theorem, and our remarks about $\zeta_k$.

The proof of Theorem \ref{orbit structure1} in the even
case is exactly the same, except that we get a
factor of $2$ due to the different definition of
the room grid. 
\newline
\newline
\noindent
{\bf Remark:\/} 
We compare Theorem \ref{orbit structure1} to a
result in [{\bf K\/}].  
The result in [{\bf K\/}] is quite
general, and so we will specialize 
it to kites.
In this case, a kite is quasi-rational iff it is rational.
The (special case of the) result in [{\bf K\/}], interpreted in
our language, says that every special orbit
is contained in one of the intervals $J_0,J_1,J_2,...$, where
$$J_a=\bigcup_{i=0}^{p+q-1} I_{ak+i}.$$
The endpoints of the $J$ intervals correspond to
{\it necklace orbits\/}.
A necklace
orbit (in our case) is an outer billiards orbit consisting
of copies of the kite, touching vertex to vertex. Compare
Figure 2.1. 

\subsection{Proof of Theorem \ref{orbit structure2}}

Let $p/q$ be some rational and let
$\widehat \Gamma$ be the
corresponding arithmetic graph.  Let
$O_2(m,n)$ denote the orbit corresponding to
the component $\widehat \Gamma(m,n)$.

\begin{lemma}
\label{stable0}
A periodic orbit 
$O_2(m,n)$ is stable iff
$\widehat \Gamma(m,n)$ is a polygon.
\end{lemma}

\startproof
Let $K$ be the period of $\Psi$ on $p_0$.
Tracing out $\widehat \Gamma(m,n)$, we get integers
$(m_k,n_k)$ such that 
\begin{equation}
\label{stability}
\Psi^k(p_0)-p_0=(2m_kA+2n_k,2\epsilon_k); \hskip 30 pt k=1,...,K.
\end{equation}
Here $\epsilon_k \in \{0,1\}$, and
$\epsilon_k=0$ iff $m_k+n_k$ is even.
The integers $(m_k,n_k)$ are determined
by the combinatorics of a finite portion
of the orbit. Hence, 
Equation \ref{stability} holds true for
all nearby parameters $A$.

If $\widehat \Gamma(m,n)$ is a closed polygon, then 
$(m_K,n_K)=0$.  But then
$\Psi^k(p_0)=p_0$ for all parameters near $A$.
If $O_2(m,n)$ is stable then
$(m_K,n_K)=(0,0)$.  Otherwise, the equation
$m_KA+n_K=0$ would force $A=-n_K/m_K$.
\endproof

\noindent
{\bf Odd Case:\/}
Assume that $A=p/q$ is an odd rational.
Say that a {\it suite\/} is
the region between two floors of the room grid.   Each suite
is partitioned into rooms.  Each room has two walls, and each
wall has a door in it.  From the Hexagrid Theorem, we see that
there is an infinite polygonal arc of $\widehat \Gamma(p/q)$
that lives in each suite.  Let
$\Gamma_k(p/q)$ denote the infinite polygonal arc that
lies in the $k$th suite.   Here
$\Gamma_0(p/q)=\Gamma(p/q)$.

We have just described the infinite family of unstable components
listed in Theorem \ref{orbit structure2}.  
All the other components of $\widehat \Gamma(p/q)$ are closed polygons
and must be confined to single rooms. The corresponding
orbits are stable, by Lemma \ref{stable0}.
The already-described polygonal arcs use up
all the doors.

Each vertex $(m,n)$ in the
arithmetic graph corresponds to the two points
$(M_1(m,n),\pm 1)$.  Thus, each component of $\widehat \Gamma$ tracks
either $1$ or $2$ orbits.  By the parity result in
Equation \ref{shortreturn}, these two points lie on
different $\psi$-orbits.
Therefore, each component of $\widehat \Gamma$ 
tracks two special orbits.
In particular, 
there are exactly two unstable orbits $U_k^+$ and $U_k^-$
contained in the interval $I_k$, and these correspond to
$\Gamma_k(p/q)$.  This completes the proof in the odd case.
\newline
\newline
{\bf Even Case:\/}
Now let $p/q$ be even.
By Lemma \ref{stable0}, it suffices to show
that all nontrivial components of
$\widehat \Gamma$ are polygons.
Suppose $\widehat \Gamma(m,n)$ is not a polygon.
Let $R$ denote reflection in the $x$-axis.
We have
\begin{equation}
\label{reflect}
R \Psi R^{-1}=\Psi^{-1}; \hskip 30 pt
R(M(m,n))=M(m+q,n-p).
\end{equation}
From this equation we see that translation by
$(q,-p)$ preserves $\widehat \Gamma$ but reverses the
orientation of all components. But then
$(m,n)+(q,-p) \not \in \widehat \Gamma(m,n)$.

\begin{center}
\resizebox{!}{1.5 in}{\includegraphics{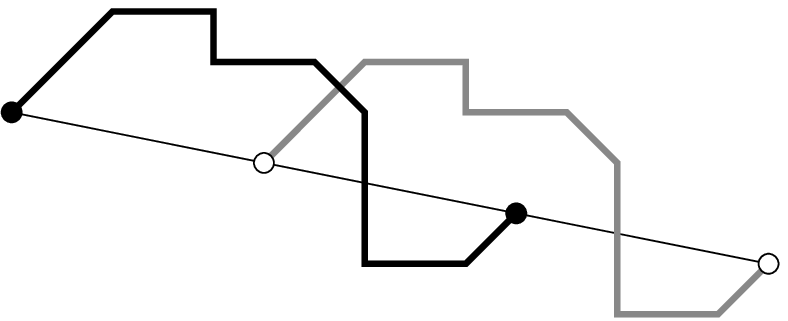}}
\newline
{\bf Figure 3.4:\/} $\gamma$ and $\gamma+(q,-p)$.
\end{center}

Since all orbits are periodic,
$(m,n)+k(p,-q) \in \widehat \Gamma(m,n)$ for some integer $k \geq 2$.
Let $\gamma$ be the arc of $\widehat \Gamma(m,n)$
connecting $(m,n)$ to $(m,n)+k(q,-p)$.
By the Embedding Theorem,  $\gamma$ and
$\gamma'=\gamma+(q,-p)$ are disjoint.   But this 
situation violates the Jordan Curve Theorem.
See Figure 3.4.

\newpage

\section{Period Copying}
\label{pc}

\subsection{Inferior and Superior Predecessors}

Let $p/q \in (0,1)$ be any odd rational.
There are unique 
rationals $p_-/q_-$ and $p_+/q_+$ such that
\begin{equation}
\label{inferior0}
\label{induct0}
\frac{p_-}{q_-}<\frac{p}{q}<\frac{p_+}{q_+}; \hskip 30 pt
\max(q_-,q_+)<q; \hskip 30 pt
qp_{\pm}-pq_{\pm}=\pm 1.
\end{equation}
See \S \ref{ssp} for more details.

We define the odd rational.
\begin{equation}
\label{inferior}
\frac{p'}{q'}=\frac{|p_+-p_-|}{|q_+-q_-|},
\end{equation}
$p'/q'$ is the unique odd rational satisfying the equation
\begin{equation}
\label{oddfarey}
q'<q; \hskip 30 pt |pq'-qp'|=2.
\end{equation}
We call $p'/q'$ the {\it inferior predecessor\/} of $p/q$, and we write
$p'/q' \leftarrow p/q$ or $p/q \to p'/q'$.
We can iterate this procedure.  Any $p/q$ belongs to a finite
chain
\begin{equation}
\frac{1}{1} \leftarrow \frac{p_1}{q_1} \leftarrow ... \leftarrow \frac{p_n}{q_n}=
\frac{p}{q}.
\end{equation}
Corresponding to this sequence we define
\begin{equation}
\label{renorm}
d_k={\rm floor\/}\bigg(\frac{q_{k+1}}{2q_k}\bigg).
\end{equation}
We define the {\it superior predecessor\/} of $p/q$ to be 
$p_k/q_k$, where $k$ is the largest index such that
$d_k \geq 1$.   It might happen that the inferior and
superior predecessors coincide, and it might not.

Here is an example, where the terms are highlighted in a suggestive way.
$$ {\bf \frac{1}{1}\/} \leftarrow \frac{1}{3} {\bf \leftarrow \frac{1}{5}\/} 
\leftarrow \frac{3}{13} \leftarrow {\bf \frac{5}{21}\/} \leftarrow \frac{13}{55}
 {\bf \leftarrow \frac{21}{89}\/} \leftarrow  \frac{55}{233}\leftarrow {\bf \frac{89}{377}\/} \ldots$$
$3/13$ has $1/5$ as both a superior and an inferior predecessor.
$5/21$ has $3/13$ as an inferior predecessor and $1/5$ as a superior prececessor.
The implied limit of this sequence is $\sqrt 5-2$, the Penrose kite parameter.

\subsection{Inferior and Superior Sequences}
\label{diop}

The inferior predecessor construction organizes all the odd
rationals into a directed tree of infinite valence.  
The rational $1/1$ is the terminal node of this tree.
The nodes incident to $1/1$ are $1/3$, $3/5$, $5/7$, etc.
  Figure 4.1 shows
part of this tree.  The edges are labelled with the
$d$ values from Equation \ref{renorm}.

\begin{center}
\resizebox{!}{2in}{\includegraphics{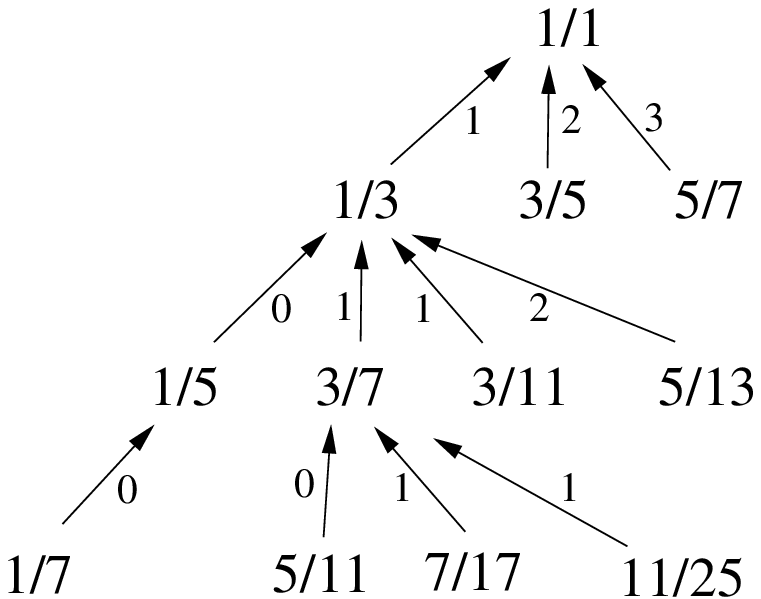}}
\newline
{\bf Figure 4.1:\/} The odd tree.
\end{center} 

The next result identifies certain of the ends of this
tree with the irrationals in $(0,1)$.  In Part IV,
we prove the following result.

\begin{lemma}[Superior Sequence]
Let $A \in (0,1)$ be irrational.
There is a unique sequence $\{p_n/q_n\}$ of odd rationals such that
\begin{equation}
\frac{p_0}{q_0}=\frac{1}{1}; \hskip 30 pt
\frac{p_{n+1}}{q_{n+1}} \to \frac{p_n}{q_n} \hskip 8 pt \forall n; \hskip 30 pt
A=\lim_{n \to \infty} \frac{p_n}{q_n}.
\end{equation}
There are infinitely many indices $n$ such that
$2q_n<q_{n+1}$. 
\end{lemma}

We call the sequence $\{p_n/q_n\}$ the {\it inferior sequence\/}.
We call $n$ a {\it superior index\/} if $2q_n<q_{n+1}$.
In terms of Equation \ref{renorm}, the index $n$
is superior if and only if $d_n \geq 1$.
We define the {\it superior sequence\/} to
be the subsequence that is indexed
by the superior indices.
Though there are many inferior and superior sequences
containing $p_n/q_n$, the initial parts of these sequences
are determined by $p_n/q_n$.  This comes from the
directed tree structure we have already mentioned.
\newline
\newline
{\bf Remark:\/}
The converse result is also true.
Any inferior sequence with infinitely
many superior terms as an irrational limit.
 This is a consequence of 
Lemma \ref{superior}.

\subsection{Strong Sequences}
\label{strongpair}

Let $A_1$ and $A_2$ be two odd rationals.  Let $\Gamma_1$ and $\Gamma_2$
be the corresponding arithmetic graphs.  We fix 
\begin{equation}
\epsilon=\frac{1}{8}.
\end{equation}
This is an arbitrary but convenient choice.

Let $V_1=(q_1,-p_1)$.  Let $\Gamma_1^1$ denote the period of
$\Gamma_1$ connecting $(0,0)$ to $V_1$.  Let
$\Gamma_1^{-1}$ denote the period of $\Gamma_1$ connecting
$(0,0)$ to $-V_1$.
We define
\begin{equation}
\Gamma_1^{1+\epsilon}=\Gamma_1^1 \cup \bigg(\Gamma_1 \cap B_{\epsilon q_1}(V_1)\bigg);
\hskip 40 pt
\Gamma_1^{-1-\epsilon}=\Gamma_1^{-1} \cup  \bigg(\Gamma_1 \cap B_{\epsilon q_1}(-V_1)\bigg).
\end{equation}
We are extending one period of $\Gamma_1$ slightly beyond one of its endpoints.
Call a monotone convergent 
sequence of odd rationals $\{p_n/q_n\}$ {\it strong\/} if it
has the following properties.
\begin{enumerate}
\item $|A-A_n|<Cq_n^{-2}$ for some universal constant $C$.
\item If $A_n<A_{n+1}$ then $\Gamma_n^{1+\epsilon} \subset \Gamma_{n+1}^1$.
\item If $A_n>A_{n+1}$ than $\Gamma_n^{-1-\epsilon} \subset \Gamma_{n+1}^{-1}$.
\end{enumerate}
In other words, $\Gamma_{n+1}$ copies about $1+\epsilon$ periods of $\Gamma_n$
for every $n$.  As usual, we have set $A_n=p_n/q_n$.

In Part IV we will prove the following result.
\begin{lemma}
\label{strongcopy}
Any superior sequence has a strong subsequence.  In
particular, any irrational in $(0,1)$ is the limit
of a strong subsequence.
\end{lemma}

In the next chapter we will prove that any
limit of a strong sequence satisfies the
conclusions of the Erratic Orbits Theorem.
Thus, Lemma \ref{strongcopy} is one of the
key ingredients in the
proof of the Erratic Orbits Theorem.  The proof of
Lemma \ref{strongcopy}, however, is rather involved.
We can prove a result nearly as strong as
the Erratic Orbits Theorem based on a slightly
weaker result that is much easier to prove.
We now describe this alternate result.

Let $\Delta_k \subset (0,1)$ denote the set of irrationals $A$ such that
the equation
\begin{equation}
\label{delta3}
0<\bigg|A-\frac{p}{q}\bigg|<\frac{1}{kq^2}; \hskip 30 pt
p,q \in \Z_{\rm odd\/}
\end{equation}
holds infinitely often.

In Part IV we prove the following result.
\begin{lemma}
\label{weakcopy}
Let $A_j=p_j/q_j$ be odd rationals such that
$|A_1-A_2|<1/(2q_1^2)$. 
\begin{itemize}
\item If $A_1<A_2$ then
$\Gamma_1^{1+\epsilon} \subset \Gamma_2^1$.
\item If $A_1>A_2$ then
$\Gamma_1^{-1-\epsilon} \subset \Gamma_2^{-1}$.
\end{itemize}
\end{lemma}

\begin{corollary}
Every $A \in \Delta_2$ is the limit of a strong sequence.
\end{corollary}

\startproof
If $A \in \Delta_2$, then there exists a monotone
sequence of solutions to Equation \ref{delta3} for
$k=2$.  This sequence is strong, by Lemma \ref{weakcopy}.
\endproof

Combining the last corollary with our work in the
next chapter, we obtain the proof of
the Erratic Orbits Theorem for all
$A \in \Delta_2$.  The reader who is
satisfied with this result can skip
most of Part IV.  The proof of
Lemma \ref{weakcopy} is really much
easier than the proof of
Lemma \ref{strongcopy}.   We close this discussion
with some observations on the size of the
sets $\Delta_k$.

\begin{lemma}
$\Delta_k$ has full measure in $(0,1)$ for any $k$.
\end{lemma}

\startproof
Any block of
$3$ consecutive odd terms $\geq k$ in
the continued fraction expansion of $A$ guarantees
a solution to Equation \ref{delta3}.  It follows
from the ergodicity of the Gauss map (or
the ergodicity of the geodesic flow on the
modular surface) that almost every $A$ has infinitely
many such blocks.  Hence $\Delta_k$ has full measure
in $(0,1)$.
\endproof

As Curt McMullen
pointed out to me, every irrational
in $(0,1)$ belongs to $\Delta_1$.
This result is similar in spirit to
Lagrange's famous theorem that every
irrational $A$ satisfies
$$\bigg|A-\frac{p}{q}\bigg|<\frac{1}{\sqrt 5 q^2}$$
infinitely often.  Lagrange's theorem doesn't
imply that every irrational lies in $\Delta_2$
because the conditions on $\Delta_2$ involve
a parity restriction.

For the interested reader, we sketch here McMullen's
argument that $\Delta_1=(0,1)-\Q$.
Consider the usual horodisk packing associated
to the modular group.  Remove all
horodisks except those based at
odd rationals.  Dilate each disk
(in the Euclidean sense) by a factor
of $2$ about its basepoint.  Observe
that the complement of these inflated
disks, in the hyperbolic plane
has infinitely many components.
Interpret this result in terms
of $\Delta_1$, using the usual
connection between the modular
horodisk packing and rational
approximation.
\newline

\subsection{The Decomposition Theorem}
\label{statedecomp}

Given an odd rational $A=p/q$, we construct the even rationals
$A_{\pm}=p_{\pm}/q_{\pm}$.  We let $A'$ be the inferior
predecessor of $A$ and we let $A^*$ be the superior predecessor.
For each rational, we use Equation \ref{boxvectors} to construct
the corresponding $V$ and $W$ vectors.  For instance,
$V_+=(q_+,-p_-)$ and $V_*=(q_*,-p_*)$.
Now we define the following lines.

\begin{itemize}
\item $L_0^-$ is the line parallel to $V$ and containing $W$.
\item $L_1^-$ is the line parallel to $V$ and containing $W^*$.
\item $L$ is the line parallel to $V$ through the $(0,0)$.
\item $L_0^+$ is the line parallel to $W$ through $(0,0)$.
\item If $q_+>q_-$ then $L_1^+$ is the line parallel to $W$ through $-V_-$.
\item If $q_+<q_-$ then $L_1^+$ is the line parallel to $W$ through $+V_+$.
\item If $q_+>q_-$ then $L_2^+$ is the line parallel to $W$ through $+V_+$.
\item If $q_+<q_-$ then $L_2^+$ is the line parallel to $W$ through $-V_-$.
\end{itemize}
Now we define the following parallelograms:
\begin{itemize}
\item $R_1$ is the parallelogram bounded by $L$ and $L_1^-$ and $L_0^+$ and $L_1^+$.
\item $R_2$ is the parallelogram bounded by $L$ and $L_0^-$ and $L_0^+$ and $L_2^+$.
\end{itemize}
The parallelogram $R_2$ is the bigger of the two parallelograms.
It is both wider and taller.
Note that translation by $V$ carries the leftmost edge of
$R_1 \cup R_2$ to the rightmost edge.
  
These might look like complicated definitions, but they are
exactly adapted to the structure of the arithmetic
graph. 
In Part IV we establish the following result.

\begin{theorem}[Decomposition]
$R_1 \cup R_2$ contains a
period of $\Gamma$.
\end{theorem}

The Decomposition Theorem is an improvement on the
containment result in the Room Lemma. It is our
main tool for Lemma \ref{strongcopy} and many of
the results we prove in Part VI.

\begin{center}
\resizebox{!}{5.1in}{\includegraphics{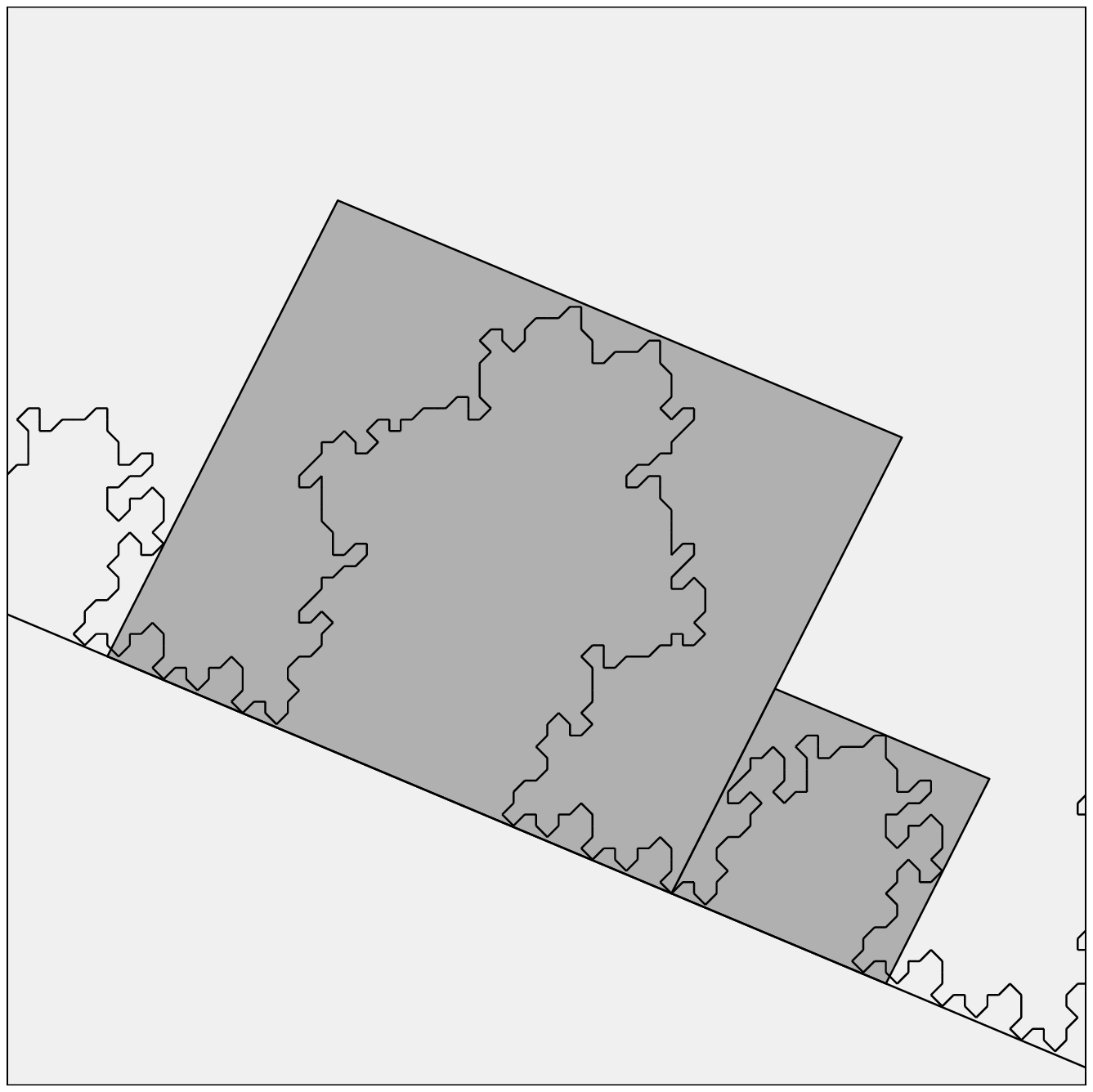}}
\newline
{\bf Figure 4.2:\/} $\Gamma(29/69)$ and
$R_1(29/69)$ and $R_2(29/69)$.
\end{center} 

Figure 4.2 shows the example $A=29/69$.
In this case,
$$A_-=21/50; \hskip 20 pt
A_+=8/19; \hskip 20 pt A'=A^*=13/31.$$
Since $q_+<q_-$, the smaller $R_1$ lies to the right of
the origin.  The ratio of heights of the two parallelograms
is $q^*/q=31/69$.  The ratio of widths is $q_+/q_-=19/50$.

Notice that the containment is
extremely efficient.  Notice also that each piece
$\Gamma \cap R_1$ and $\Gamma \cap R_2$ has approximate
bilateral symmetry.  This situation always happens.  We will
explain this symmetry in \S \ref{symm1proof}.

\newpage
\section{Proofs of the Main Results}
\label{bb}

\subsection{Proof of the Erratic Orbits Theorem}

the Erratic Orbits Theorem follows from
Lemma \ref{strongcopy},  Lemma \ref{dichotomy} and  Lemma \ref{stronglim} (stated below).
For the reader who wants to
take a shortcut, we remark again that we
prove the Erratic Orbits Theorem for all $A \in \Delta_2$, when
we use the much easier
Lemma \ref{weakcopy} in place of Lemma \ref{strongcopy}.

\begin{lemma}
\label{stronglim}
Suppose that $A$ is the limit of a strong sequence $\{p_n/q_n\}$.
Then the Erratic Orbits Theorem holds for $A$.
\end{lemma}

In our proof, we will 
consider the monotone increasing case. The other
case is essentially the same.
Note that our sequence remains strong if we pass to a subsequence.
Passing to a subsequence, we arrange that
\begin{equation}
\epsilon q_{n+1}>10q_n
\end{equation}

Let $V_n=(q_n,-p_n)$. 
Define
\begin{equation}
\Gamma_n^2=\Gamma_n^1+V_{n+1}, 
\end{equation}

\begin{lemma}
\label{erratic1}
\begin{equation}
\label{induct222}
\Gamma_n^1 \subset \Gamma_{n+1}^1; \hskip 30 pt
\Gamma_n^2 \subset \Gamma_m^1 \hskip 30 pt \forall m \geq n+2.
\end{equation}
\end{lemma}

\startproof
 We have 
$$\Gamma_n^{1+\epsilon} \subset \Gamma_{n+1}^1$$ by definition, and
$$\Gamma_{n}^1+V_{n+1} \subset \Gamma_{n+1}$$ because
$\Gamma_{n+1}$ is invariant under translation by
$V_{n+1}$.   Our choice of subsequence gives
\begin{equation}
\label{tighten3}
\Gamma_n^{1+\epsilon} \subset B_{10q_n}(0,0) \subset B_{\epsilon q_{n+1}}(0,0) \cap \Gamma_{n+1}.
\end{equation}
The first containment comes from the Room Lemma.
Translating by $V_{n+1}$, we have
\begin{equation}
\Gamma_n^1+V_{n+1}
\subset B_{\epsilon q_{n+1}}(V_{n+1}) \cap \Gamma_{n+1}^1 \subset \Gamma_{n+1}^{1+\epsilon} \subset \Gamma_{n+2}^1.
\end{equation}
Equation \ref{induct222} follows immediately.
\endproof

If follows from Equation \ref{induct222} and induction that
\begin{equation}
\label{vertexdef}
\omega_n=\omega(\sigma):=\sum_{k=1}^{n-1} \epsilon _k V_{2k+1}
\end{equation}
is a vertex of $\Gamma_{2n}^1$ 
for any binary sequence $\epsilon_1,...,\epsilon_{n-1}$.
Let $\Pi$ denote the set of not-eventually-constant sequences.
Given any $\sigma \in \Pi$, we form the sequence of
translated graphs
\begin{equation}
\Gamma_n'=\Gamma_{2n}^1 - \omega_n.
\end{equation}
Here $\omega_n$ is based on the first $n-1$ terms of
$\sigma$, as in Equation \ref{vertexdef}.

\begin{lemma}
$\{\Gamma_n'\}$ Hausdorff converges to $\Gamma$, an
open polygonal arc that rises unboundedly far,
in both directions from the line $L$ of slope
$(-A)$ through the origin.
\end{lemma}

\startproof
Figure 5.1 shows the sort of binary structure that we have established.
In this figure, the notation $ij$ stands for $\Gamma_i^j$.

\begin{center}
\resizebox{!}{1.8 in}{\includegraphics{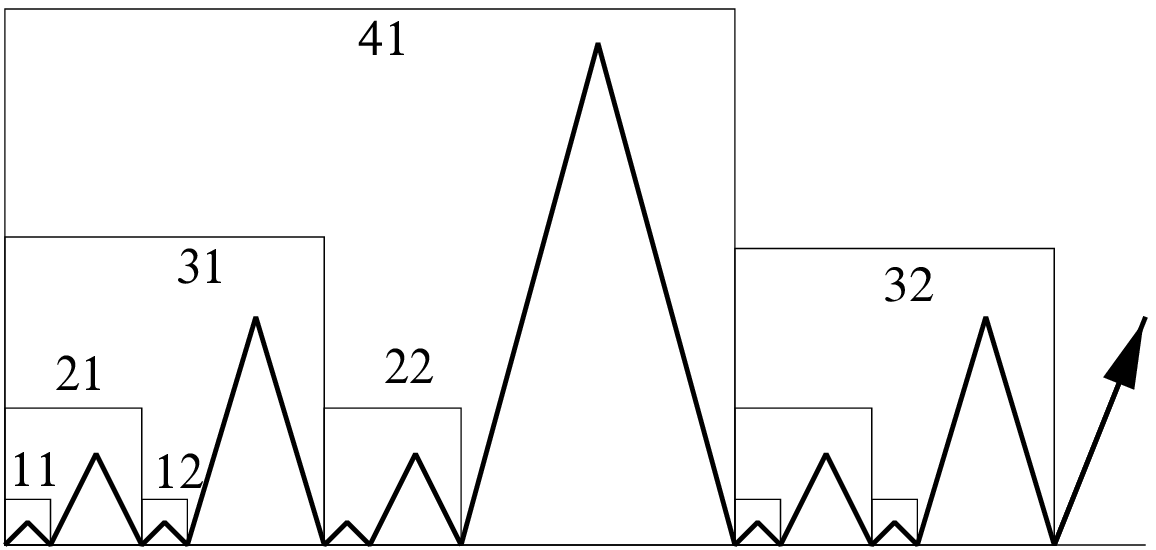}}
\newline
{\bf Figure 5.1:\/} large scale Cantor set structure
\end{center}

Figure 5.2 shows a simpler picture that retains the
structure of interest to us.

\begin{center}
\resizebox{!}{1.5in}{\includegraphics{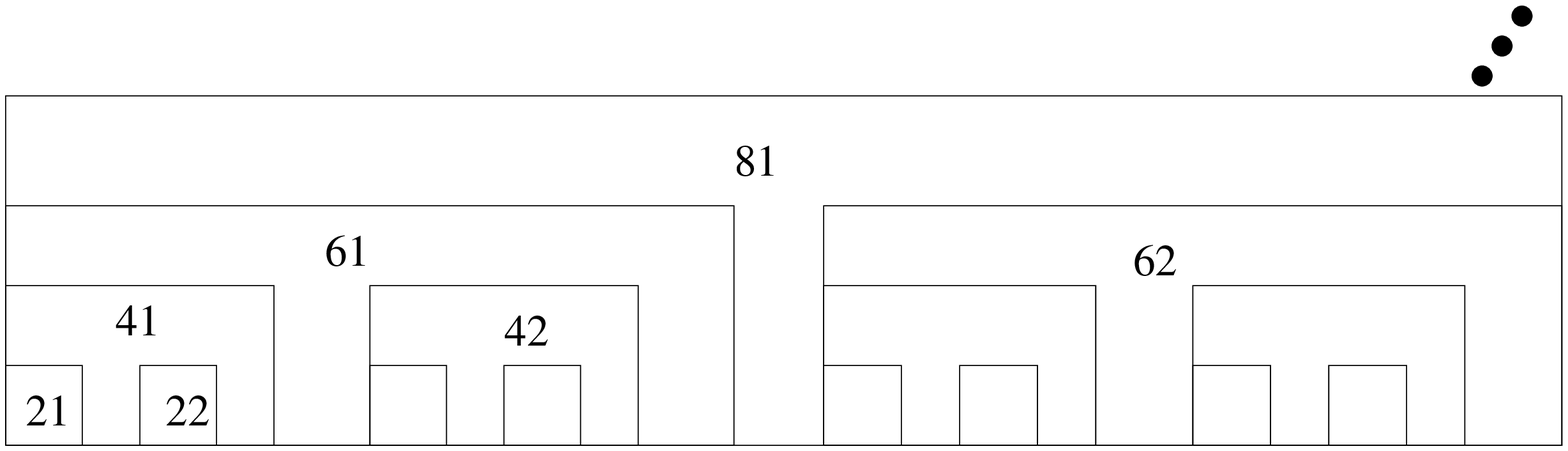}}
\newline
{\bf Figure 5.2:\/} large scale Cantor set structure
\end{center}

To make sense of Figures 5.1 and 5.2, we say that the
{\it box\/} containing $\Gamma_n^1$ is $R_n=R(A_n)$, the
box from the Room Lemma.
For instance, the $8$ smallest boxes in Figure 5.2 are
\begin{equation}
\label{binary}
R_2+\epsilon_1 V_3 + \epsilon_2 V_5 + \epsilon_3 V_7; \hskip 30 pt
\epsilon_j \in \{0,1\}.
\end{equation}
The larger boxes have a similar description.
The boxes are not quite nested, on account of the tiny mismatches
between the slopes of their boundaries, but they are very nearly nested.
See Property 4 below.
We rank each box according to the label of its leftmost translate.
The smallest boxes have rank $2$.  The next-smallest have rank $4$.  And so on.
The following structure emerges.

\begin{enumerate}
\item If two boxes have the same rank, then the
corresponding arcs are translates of each other.
\item The boxes of rank $n$ have diameter $O(q_n)$.
\item The arc inside a box of rank $n$, a translate
of $\Gamma_{2n}^1$, contains the
bottom corners of the box containing it
and rises up $O(q_n)$
units towards the top of this box.
This is a consequence of the Room Lemma.
\item The bottom edge of a box of rank $n$
lies within $O(1/q_n)$ of the bottom edge of
the box of rank $n+1$ that nearly contains it.
First we prove this for $R_n$ and $R_{n+1}$.
The bottoms of these boxes meet at the origin.
The difference in slopes of $O(1/q_n^2)$.  The
length of the bottom edge of $R(A_n)$ is $O(q_n)$.
The estimate follows immediately.  Once
we note that $V_{n+1}$ is $O(1/q_{n+1})$ units
away from the bottom of $R_{n+1}$, we get
the same result for $R_n+V_{n+1}$ and $R_{n+1}$.
The general case now follows from translation.
\end{enumerate}

By construction, the pattern of boxes surrounding
$\omega_n$ stabilizes when we view any fixed-radius
neighborhood of $\omega_n$.  More formally,
for any $R$, there is some $N$ such that $m,n>N$
implies that $\Gamma^1_{2m} \cap B_R(\omega_m)$
is a translate of $\Gamma_{2n}^1 \cap B_R(\omega_n)$.
Here we are crucially using the fact that
$\sigma \in \Pi$, so that our common pattern
of boxes grows both to the left and to the right
of the points of interest.
Hence, the sequence $\{\Gamma'_n\}$ Hausdorff
converges to a limit $\Gamma$.

From the $4$ properties listed above, $\Gamma$ is
an infinite open polygonal arc that rises unboundedly
far, in both directions, from $L$.
\endproof

It remains to recognize $\Gamma$.
Let $M$ be the map from Equation \ref{funm}, relative to
the limit parameter $A$.  Given
$\sigma=\{\epsilon_k\} \in \Pi$, the point
\begin{equation}
\alpha(\sigma)=\bigg(\sum_{k=1}^{\infty} 2\epsilon_k \Big( Aq_{2k+1}-p_{2k+1}\Big), - 1\bigg)
\end{equation}
is well defined because the $k$th term in the series has size $O(1/q_{2k+1})$,
and the sequence $\{q_{2k+1}\}$ grows exponentially. 
The union of such limits, taken over all of $\Pi$,
contains a pruned Cantor set.  Throwing out a countable
subset of $\Pi$, we can arrange that our pruned Cantor set
is disjoint from $2\Z[A]$.  
But then, and $\alpha=\alpha(\sigma)$ we consider has a well-defined
orbit, by 
Lemma \ref{irrat}.

\begin{lemma} $\Gamma$ is the arithmetic graph of $\alpha$.
\end{lemma}

\startproof
Define
\begin{equation}
\alpha_n=M_{2n}(\sigma_n)=
\bigg(\sum_{k=1}^{n-1} 2\epsilon_k \Big( A_{2n}q_{2k+1}-p_{2k+1}\Big), - 1\bigg)
\end{equation}
An easy argument shows that $\alpha_n \to \alpha$. By
construction, $\Gamma'_n$ is one period of the arithmetic
graph of $\alpha_n$ relative to $A_{2n}$.  The distance
that $\Gamma_n'$ extends from the origin, in either
direction, tends to $\infty$ with $n$.
By the Continuity Principle, $\{\Gamma_n'\}$ converges to
the arithmetic graph of $\alpha$.
\endproof

Given the structure of $\Gamma$, we know that
$\alpha$ has an unbounded orbit.  To finish our proof
of Lemma \ref{stronglim},
we just need to show that $\alpha$ has an erratic orbit.
Call an arc of $\Gamma_n'$ {\it stable\/} if this same
arc also belongs to $\Gamma_m'$ for $m>n$.
By construction, we get the following result.
For any $k$, there is some $n$ such that
$\Gamma_n'$ contains a stable arc of the form
$\beta-\omega_n$.  Here $\beta$ is a full period of
$\Gamma_k$, but contained in $\Gamma_{2n}^1$.
Some vertex $v$ of $\beta$ has the form
\begin{equation}
\sum_{j=k}^{n-1} \epsilon_j V_{2j+1}
\end{equation}
The distance from $v$ to the baseline of 
$\Gamma_{2n}$ is $O(1/q_{2k+1})$.
But then, the distance from $v-\omega_n$ to
the baseline of $\Gamma'_n$ is
$O(1/q_{2k+1})$.  But $v-\omega_n$ is
also a vertex of $\Gamma$ (by stability)
and its distance to the baseline of $\Gamma$
is also $O(1/q_{2k+1})$.   
 We can choose our arc
$\beta-\omega_n$ either to the left or
to the right of the origin.  Hence,
both sides of the limit $\Gamma$ come arbitrarily
close to the baseline of $\Gamma$.

\subsection{Proof of Theorem \ref{dichotomy0}}

The following result combines with the
Erratic Orbits Theorem to prove
Theorem \ref{dichotomy0}: Every
special orbit is either periodic
or else unbounded in both directions.
Note that the result does not quite
require the existence of erratic orbits,
but only the existence of orbits that
come fairly close to the kite vertex.

\begin{lemma}
\label{dichotomy}
Suppose that $A$ is a parameter, and
$p \in (0,2) \times \{1\}$ has an orbit
that is unbounded in both directions.
Then all special orbits relative to $A$ are either
periodic or else unbounded in both directions.
\end{lemma}

\startproof
We write $p=(2\zeta,1)$.  By hypothesis,
$\zeta \in (0,1)$.
Suppose that $\beta$ has an aperiodic orbit that is
forwards bounded. (The backwards case is similar.)
For ease of exposition, we suppose that
$\beta \not \in 2\Z[A]$, so that all components of
the arithmetic graph $\widehat \Gamma$ associated
to $\beta$ are well defined.  
In case $\beta \in 2\Z[A]$,
we simply apply our argument to a sequence
$\{\beta_n\}$ converging to $\beta$ and invoke
the Continuity Principle. Our robust geometric limit
argument works the same way with only notational
complications.

Let $\Gamma$ be the component of $\widehat \Gamma$
that tracks $\beta$.
The forwards direction
$\Gamma_+$ remains within a bounded
distance of the baseline $L$ of $\widehat \Gamma$ and
yet is not periodic.  Hence,
$\Gamma_+$ travels infinitely far either to the
left or to the right.
Since $L$ has irrational slope, 
we can find a sequence of vertices
$\{v_n\}$ of $\Gamma_+$ such that
the vertical distance from $v_n$ to $L$
converges to $\zeta+N$ for some integer $N$.
Let $w_n=v_n-(0,N)$.   Let $\gamma_n$ be the
component of $\widehat \Gamma_n$ containing
$w_n$.  Note that $M(w_n) \to p$.  Here
$M$ is as in Equation \ref{funm}.

Let $T_n$ be a translation so that $T_n(w_n)=(0,0)$.
By compactness, we can choose our sequence so that
$\{T_n(\Gamma_+)\}$ converges to an infinite polygonal 
arc $X$ that remains within a 
bounded distance of any line parallel to $L$.
By construction $X$ travels infinitely far both
to the left and to the right.  At the same time,
$\{T_n(\gamma_n)\}$ converges to 
the arithmetic graph $Y$ of $\zeta$. 
Here $Y$ starts at $(0,0)$, a point within $1$
unit of the baseline $L_{\infty}=\lim T_n(L)$
and rises unboundedly far from $L_{\infty}$.
Hence $Y$ starts out below $X$ and rises above
$X$, contradicting the Embedding Theorem.
\endproof

\subsection{The Rigidity Lemma}
\label{rigidity}

Here we prove a technical convergence result
that helps in the proofs of both
Theorem \ref{density} and Theorem \ref{cantor}.

\begin{lemma}[Rigidity]
Let $A_n$ be any sequence of parameters converging
to the irrational parameter $A$.
Let $\zeta_n \in [0,2] \times \{1\}$ be a sequence of points
converging to $(0,1)$. Let $\Gamma(\zeta_n,A)$ be the
arithmetic graph of $\zeta_n$ relative to $A$.
Then the sequence $\{\Gamma(\zeta_n,A)\}$ Hausdorff converges.
\end{lemma}
We think of our result as a rigidity result because
it implies that all possible limits we can take in the
above manner are the same.

Given $\epsilon>0$, let
$\Sigma_{\epsilon}(A) \subset (0,1)^2$ denote
those pairs $(s,A')$ where $s \in (0,\epsilon)$
and $|A'-A|<\epsilon$.  Let $O(s;A')$ denote the
outer billiards orbit of $(s,1)$ relative to
$K(A')$.

\begin{lemma}
\label{localstab}
For any $N$ there is some $\epsilon>0$ with the
following property.   The first $N$ iterates
of $O(s;A')$, forwards and backwards, are well defined
provided that $(s,A') \in \Sigma_{\epsilon}(A)$.
\end{lemma}

\startproof
Inspecting the proof of Lemma \ref{basic}, we
draw the following conclusion.
If $O(s;A')$ is not defined
after $N$ iterates, then 
$s=2A'm+2n$ for
integers $m,n \in (-N',N')$.
Here $N'$ depends only on $N$.
Rearranging this equation, we get
$$|A'-\frac{m}{n}|<\frac{s}{2m}.$$
For $s$ sufficiently small and
$A'$ sufficiently close to $A$, this
is impossible.   
\endproof

\begin{corollary}
\label{hausdorff1}
For any $N$ there is some $\epsilon>0$ with the
following property.   The combinatorics of the
first $N$ forward iterates 
of $O(s;A')$ is independent of the choice
of point $(s,A') \in \Sigma_{\epsilon}(A)$.
The same goes for the first $N$ backwards iterates.
\end{corollary}

\startproof
If all orbits in some interval are defined, then all
orbits in that interval have the same combinatorial
structure.
\endproof

The Rigidity Lemma is now a consequence of Corollary
\ref{hausdorff1} and the Return Lemma.  The Return
Lemma guarantees that as $N \to \infty$, the number
of returns to $\Xi$ tends to $\infty$ as well.

\subsection{Proof of Theorem \ref{density}}

First of all, since outer billiards is a piecewise
isometry, the set of periodic orbits is open
in $\R \times \Z_{\rm odd\/}$.  We just need
to prove density.

Let $A$ be an irrational parameter.  Let
$\widehat \Gamma$ be an arithmetic graph
associated to $A$, such that $\Gamma$
tracks an erratic orbit.   
Since $A$ is irrational, 
we can find a sequence of vertices $\{(m_k,n_k)\}$
of odd parity that converges to the baseline of
$A$.   Let $\gamma_k$ be the component of
$\widehat \Gamma$ that contains $(m_k,n_k)$.
Note that $\gamma_k \not = \Gamma$ because
$\Gamma$ only contains vertices of even parity.
By the Embedding Theorem, $\gamma_k$ is
trapped underneath $\Gamma$. 
Hence $\gamma_k$ is a polygon.
Let $|\gamma_k|$ denote the maximal distance
between a pair of low vertices on $\gamma_k$.

\begin{lemma}
\label{spread}
$|\gamma_k| \to \infty$ as $k \to \infty$.
\end{lemma}

\startproof
By the Rigidity Lemma, a very long
arc of $\gamma_k$, with one endpoint
$(m_k,n_k)$, 
agrees with the Hausdorff limit
$\lim_{n \to \infty} \Gamma(p_n/q_n)$. Here
$\{p_n/q_n\}$ is an approximating
strong sequence.  But this limit has vertices
within $\epsilon$ of the baseline and at least
$1/\epsilon$ apart for any $\epsilon>0$.  Our
result now follows from Hausdorff continuity.
\endproof

Let $S_k$ denote the set of components $\gamma'$ of $\widehat \Gamma$
such that $\gamma'$ is translation equivalent to $\gamma_k$
and the corresponding vertices are low.  The vertex
$(m,n)$ is low if the baseline of $\widehat \Gamma$
separates $(m,n)$ and $(m,n-1)$.

\begin{lemma}
There is some constant $N_k$ so that every point of $L$ is within
$N_k$ units of a member of $S_k$.
\end{lemma}

\startproof
Say that a lattice point $(m,n)$ is {\it very low\/} if
it has depth less than $1/100$ (but still positive.)
The polygon
$\gamma_k$ corresponds to a periodic orbit $\xi_k$.
Since $\xi_k$ is periodic, there is an open neighborhood
$U_k$ of $\xi_k$ such that all orbits in $U_k$ are
combinatorially identical to $\xi_k$. Let $M$ be fundamental map
associated
to $\widehat \Gamma$.   Then
$M^{-1}(U_k)$ is an open strip, parallel to $L$.
Since $L$ has irrational slope, there is some
constant $N_k$ so that every point of $L$ is within
$N_k$ of some point of $M^{-1}(U_k) \cap \Z^2$.
But the components of $\widehat \Gamma$ containing
these points are translation equivalent to $\gamma_k$.
Choosing $U_k$ small enough, we can guarantee that
the translations taking $\gamma_k$ to the
other components carry the very low vertices of
$\gamma_k$ to low vertices.
\endproof

Given two polygonal components $X$ and $Y$ of
$\widehat \Gamma$, we write   $X \bowtie Y$
if one low vertex of $Y$ lies to the left of
$X$ and one low vertex of $Y$ lies to the
right of $X$.  See Figure 5.3.  In this case,
$X$ is trapped underneath $Y$, by the
Embedding Theorem.

\begin{center}
\psfig{file=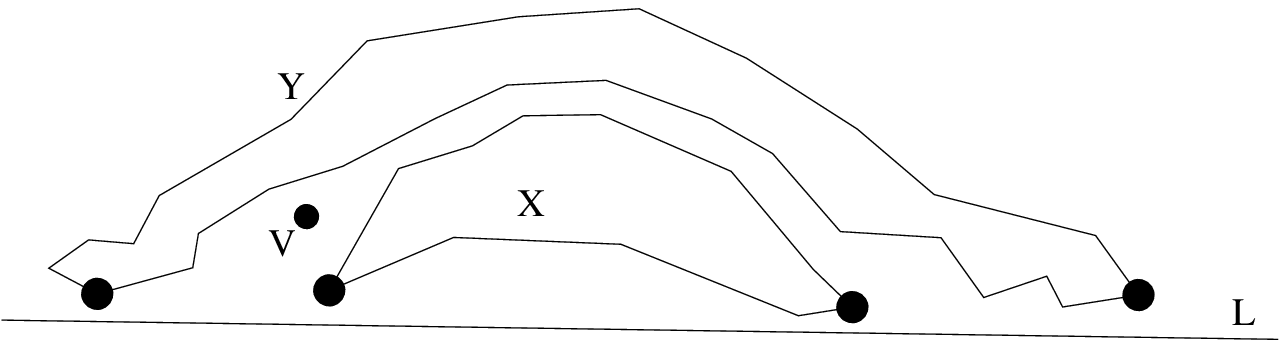}
\newline
{\bf Figure 5.3:\/} One polygon overlaying another.
\end{center} 

Now we pass to a subsequence so that
\begin{equation}
\label{bulge}
|\gamma_{k+1}|>10(N_k+|\gamma_k|).
\end{equation}
Equation \ref{bulge} has the following consequence.
For any integer $N$, we can find
components $\gamma_j$ of $S_j$, for $j=N,...,2N$ such
$\gamma_N \bowtie ... \bowtie \gamma_{2N}$.
Let $L_N$ denote the portion of $L$ between the two distinguished
low points of $\gamma_N$.  Let $\Lambda_N$ denote the set of
lattice points within $N$ units of $L_N$.  The set
$\Lambda_N$ is a parallelogram whose base is $L_N$,
a segment whose length tends to $\infty$ with $N$.
The height of $\Lambda_N$ tends to $\infty$ as well.

\begin{lemma}
The set $M(\Z^2 \cap \Lambda_N)$ consists entirely of
periodic orbits.
\end{lemma}

\startproof
Let $V$ be a vertical ray whose $x$-coordinate is an integer.
If $V$ starts out on $L_n$ then $V$
must travel upwards at least $N$ units
before escaping from underneath $\gamma_{2N}$. This is
an application of the pideonhole princple. The point is that $V$
must intersect each $\gamma_j$ for $j=N,...,2N$, in a different
lattice point.  Hence, any point of $\Lambda_N$ is
trapped beneath $\gamma_{2N}$.
\endproof

Given the fact that both base and height of
$\Lambda_N$ are growing unboundedly, and the
fact that $A$ is an irrational parameter, the
union
$\bigcup_{N=1}^{\infty} M(\Lambda_N \cap \Z^2)$
is dense in $\R_+$.
Hence, the set of periodic orbits starting in $\R_+ \times \{-1,1\}$ is
dense in the set of all special orbits. 
Our proof of the Pinwheel Lemma in Part II shows that every
special orbit eventually lands in 
$\R_+ \times \{-1,1\}$.   Hence, the set of periodic
special orbits is dense in $\R \times \Z_{\rm odd\/}$.

\newpage

\noindent
{\bf {\huge Part II\/} \/}
\newline

In this part of the monograph we will state and prove the
Master Picture Theorem.  All the auxilliary theorems left over from
Part I rely on this central result.
Here is an overview of the material.

\begin{itemize}
\item In \S \ref{masterpicture} we will state the Master Picture Theorem.  
Roughly, the Master Picture Theorem says that
 the structure of the return map $\Psi$
is determined by a pair of maps into a
flat $3$-torus, $\R^3/\Lambda$,
together with a partition of $\R^3/\Lambda$ into polyhedra.
Here $\Lambda$ is a certain $3$-dimensional lattice that
depends on the parameter.

\item In \S \ref{pinwheel}, we will prove
the Pinwheel Lemma, a key technical step along the
way to our proof of the Master Picture Theorem.  The Pinwheel Lemma states that we can
factor the return map $\Psi$ into a composition of
$8$ simpler maps, which we call {\it strip maps\/}.
A strip map is a very simple map from the
plane into an infinite strip.

\item  In \S \ref{torus} we prove
the Torus Lemma, another key result.
The Torus Lemma implies that there exists some
partition of our torus into open regions, such that
the regions determine the structure of the
arithmetic graph.  
The Torus Lemma reduces the Master Picture Theorem to
a rough determination of the singular set.
The singular set is the (closure of the) set of points
in the torus corresponding to points
where the return map is not defined.

\item In \S \ref{sing} we verify, with the
aid of symbolic manipulation, certain 
functional identities that arise in
connection with the Torus Lemma. These
function identities are the basis for
our analysis of the singular set.

\item In \S \ref{torus2} we combine the Torus
Lemma with the functional identities to prove the
Master Picture Theorem.

\item in \S \ref{computations} we will explain
how one actually makes computations with the
Master Picture Theorem.  \S \ref{master1} will be
very important for
Part IV of the monograph.

\end{itemize}
\newpage

\section{The Master Picture Theorem}
\label{masterpicture}
\label{returnproof}

\subsection{Coarse Formulation}

Recall that $\Xi=\R_+ \times \{-1,1\}$.  We distinguish two special
subsets of $\Xi$.
\begin{equation}
\label{XiPlus}
\Xi_+=\bigcup_{k=0}^{\infty} (2k,2k+2) \times \{(-1)^k\}; \hskip 30 pt
\Xi_-=\bigcup_{k=1}^{\infty} (2k,2k+2) \times \{(-1)^{k-1}\}.
\end{equation}
Each set is an infinite disconnected union of open intervals
of length $2$.  
Reflection in the $x$-axis interchanges $\Xi_+$ and $\Xi_-$.   The
union $\Xi_+ \cup \Xi_-$ partitions $(\R_+-2\Z) \times \{\pm 1\}$.

Define
\begin{equation}
R_A=[0,1+A] \times [0,1+A] \times [0,1]
\end{equation}
$R_A$ is a fundamental domain for the action of a
certain lattice $\Lambda_A$.  We have
\begin{equation}
\label{lattice}
\Lambda_A=\left[\matrix{1+A & 1-A & -1 \cr
0&1+A & -1 \cr
0&0&1}\right] \Z^3
\end{equation}
We mean to say that
$\Lambda_A$ is the $\Z$-span of the column vectors of the 
above matrix. 

We define $\mu_+: \Xi_+ \to R_A$ and $\mu_-: \Xi_- \to R_A$ by the
equations
\begin{equation}
\label{muplus}
\label{muminus}
\mu_{\pm}(t,*)=\bigg(\frac{t-1}{2},\frac{t+1}{2},\frac{t}{2}\bigg) \pm
\bigg(\frac{1}{2},\frac{1}{2},0\bigg) \hskip 10 pt
{\rm mod\/} \hskip 10 pt \Lambda.
\end{equation}
The maps only depend on the first coordinate.
In each case, we mean to map $t$ into $\R^3$ and then use the
action of $\Lambda_A$ to move the image into $R_A$.  It might
happen that there is not a unique representative in $R_A$.
(There is the problem with boundary points, as usual with
fundamental domains.)  However, if $t \not \in 2\Z[A]$, this
situation does not happen.  The maps
$\mu_+$ and $\mu_-$ are locally affine.

Here is a coarse formulation of the
Master Picture Theorem.   We will state the entire result
in terms of $(+)$, with the understanding that the same
statement holds with $(-)$ replacing $(+)$ everywhere.
Let $\Psi$ be the first return map.

\begin{theorem}
\label{weakmpt}
 For each parameter $A$ there is a partition
$({\cal P\/}_A)_+$ of $R_A$ into finitely many convex
polyhedra.  If $\Psi$ is defined on
$\xi_1,\xi_2 \in \Xi_+$ and
$\mu_+(\xi_1)$ and $\mu_{+}(\xi_2)$ lie in the
same open polyhedron
of $({\cal P\/}_A)_+$, then
$\Psi(\xi_1)-\xi_1=\Psi(\xi_2)-\xi_2$.
\end{theorem}

\subsection{The Walls of the Partitions}
\label{walls}

In order to make Theorem \ref{weakmpt} precise, we need to
describe the nature of the partitions $({\cal P\/}_A)_{\pm}$,
and also the rule by which the polygon in the partition
determines $\Psi(\xi)-\xi$.  We will make several passes
through the description, adding a bit more detail each time.

The polyhedra of $({\cal P\/}_A)_{\pm}$ are cut out by the
following $4$ families of planes.
\begin{itemize}
\item $\{x=t\}$ for $t=0,A,1,1+A$.
\item $\{y=t\}$ for $t=0,A,1,1+A$.
\item $\{z=t\}$ for $t=0,A,1-A,1$.
\item $\{x+y-z=t\}$ for $t=-1+A,A,1+A,2+A$.
\end{itemize}
The complements of the union of these planes are the open
polyhedra in the partitions.

\begin{center}
\psfig{file=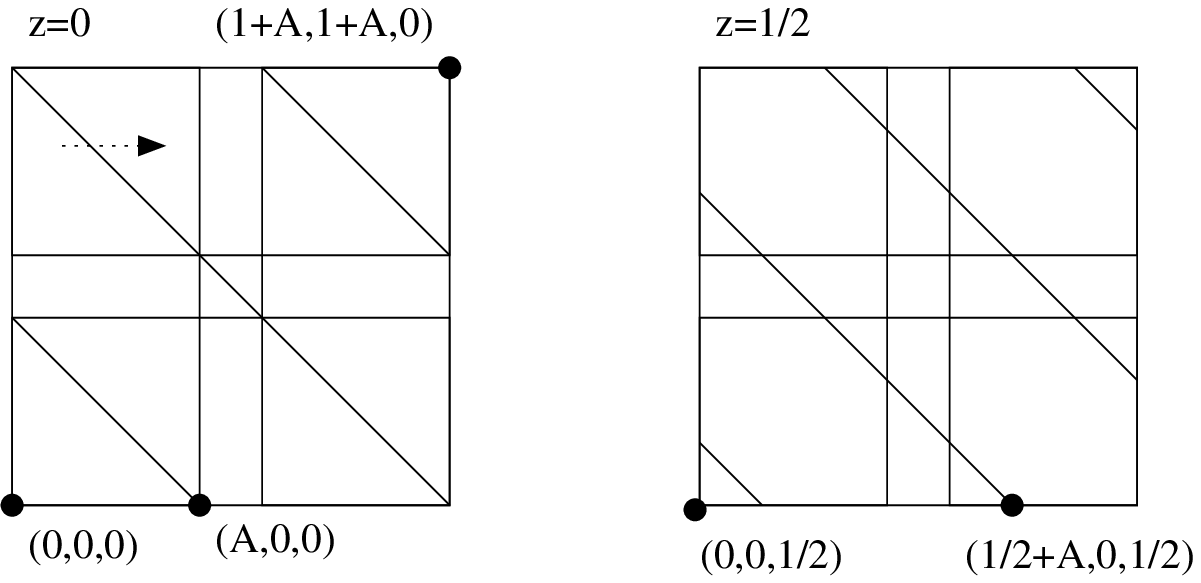}
\newline
{\bf Figure 6.1:\/} Two slices of the partition for $A=2/3$.
\end{center} 

Figure 6.1 shows a picture of two
slices of the partition for the parameter $A=2/3$.  We have sliced
the picture at $z=0$ and $z=1/2$. We have labelled several points
just to make the coordinate system more clear.
The little arrow in the picture
indicate the ``motion'' the diagonal lines would make were we to
increase the $z$-coordinate and show a kind of movie of the
partition.   The reader can see this partition for any parameter
and slice using Billiard King.

\subsection{The Partitions}
\label{assign}

For each parameter $A$ we get a solid body $R_A$ partitioned
into polyhedra.  We can put all these pieces together into
a single master picture.  We define
\begin{equation}
\label{FUND}
R=\bigcup_{A \in (0,1)} R_A \times \{A\} \subset \R^4.
\end{equation}
Each $2$-plane family discussed above
gives rise to a hyperplane family in $\R^4$.
 These hyperplane families are now all defined
over $\Z$, because the variable $A$ is just the $4$th
coordinate of $\R^4$ in our current scheme.
Given that we have two maps $\mu_+$ and $\mu_-$, it is
useful for us to consider two identical copies
$R_+$ and $R_-$ of $R$. 

We have a fibration $f: \R^4 \to \R^2$ given by
$f(x,y,z,A)= (z,A)$.   This fibration in turn gives
a fibration of $R$ over the unit square $B=(0,1)^2$. 
 Figure 6.1
draws the fiber $f^{-1}(3/2,1/2)$.  The base space $B$
has a partition into $4$ regions, as shown in
Figure 6.2.  

\begin{center}
\psfig{file=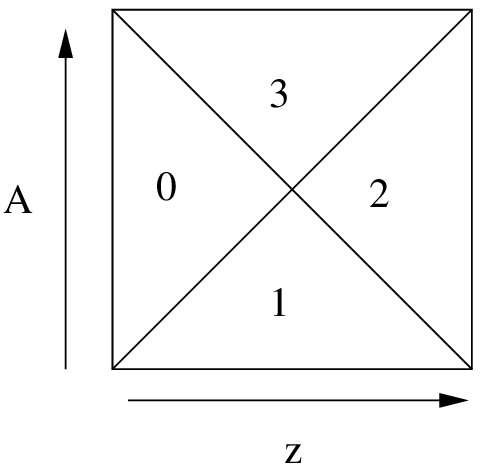}
\newline
{\bf Figure 6.2:\/} The Partition of the Base Space
\end{center} 

All the fibers above the same open region in the
base space have the same combinatorial structure.
Figure 6.3 explains precisely how the partition 
assigns the value of the return map.  Given a point
$\xi \in \Xi_+$, we have a pair of integers
$(\epsilon_1^+(\xi),\epsilon_2^+(\xi))$ such that
\begin{equation}
\Psi(\xi)-\xi=2(\epsilon_1^+,\epsilon_2^+,*).
\end{equation}
The second coordinate, $\pm 2$, is determined by the
parity relation in Equation \ref{shortreturn}.
Similarily, we have $(\epsilon_1^-,\epsilon_2^-)$
for $\xi \in \Xi_-$.

Figure 6.3 shows a schematic picture of $R$.
For each of the $4$ open triangles in the base,
we have drawn a cluster of 
$4$ copies of a representative fiber over that triangle.
The $j$th column of each cluster determines the
value of $\epsilon_j^{\pm}$.  The first row
of each cluster determines $\epsilon_j^+$ and
the second row determines $\epsilon_j^-$.
A light shading indicates a value of $+1$. A
dark shading indicates a value of $-1$.  No
shading indicates a value of $0$.

\begin{center}
\psfig{file=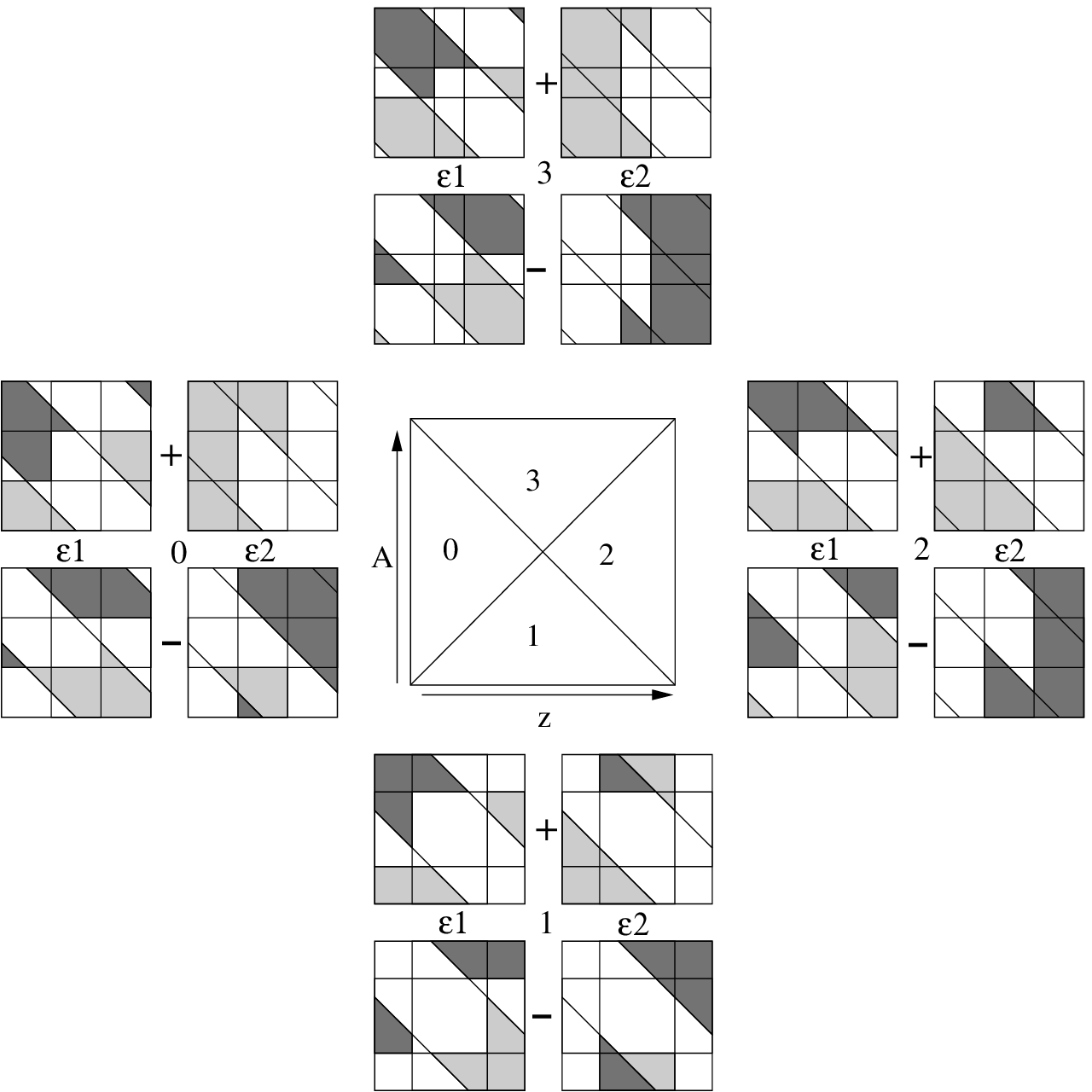}
\newline
{\bf Figure 6.3:\/} The decorated fibers
\end{center}

Given a generic point $\xi \in \Xi_{\pm}$, the image
$\mu_{\pm}(\xi)$ lies in some fiber.  We then use the
coloring scheme to determine $\epsilon_j^{\pm}(\xi)$
for $j=1,2$.  (See below for examples.)  Theorem \ref{weakmpt}, together with the
description in this section, constitutes the
Master Picture Theorem.  In \S \ref{computations}
we explain with more traditional formulas how
to compute these values.
The reader can get a vastly superior
understanding of the partition using Billiard King.

\subsection{A Typical Example}
\label{example}

Here we will explain how the Master Picture Theorem determines the
local structure of the arithmetic graph $\Gamma(3/5)$ at the point
$(4,2)$.  Letting $M$ be the fundamental map associated to
$A=3/5$ (and $\alpha=1/(2q)=1/10)$.
$$M(4,2)=\Big((8)(3/5)+(4)+(1/5),(-1)^{4+2+1}\Big)=(9,-1) \in \Xi_-.$$
So, $\mu_-(9,-1)$ determines the forwards direction and
$\mu_+(9,1)$ determines the backwards direction.
(Reflection in the $x$-axis conjugates $\Psi$ to its inverse.)

We compute
$$\mu_+(9,1) =(\frac{9}{2},\frac{11}{2},\frac{9}{2}) \equiv (\frac{1}{10},\frac{3}{2},\frac{1}{2})\ {\rm mod\/}\ \Lambda;$$
$$\mu_-(9,-1) =(\frac{7}{2},\frac{9}{2},\frac{9}{2}) \equiv  (\frac{7}{10},\frac{1}{2},\frac{1}{2})\ {\rm mod\/}\ \Lambda.$$
(In \S \ref{computations} we will explain algorithmically how to make these computations.)
We have $(z,A)=(1/2,3/5)$.  There we need to look at Cluster 3, the
cluster of fibers above region $3$ in the base.  Here is the
plot of the two points in the relevant fiber.
When we look up the regions in Figure 6.3, 
we find that $(\epsilon_1^+,\epsilon_2^+)=(-1,1)$ and $(\epsilon_1^-,\epsilon_2^-)=(1,0)$.
The bottom right of Figure 6 shows the corresponding local picture
for the arithmetic graph.

\begin{center}
\psfig{file=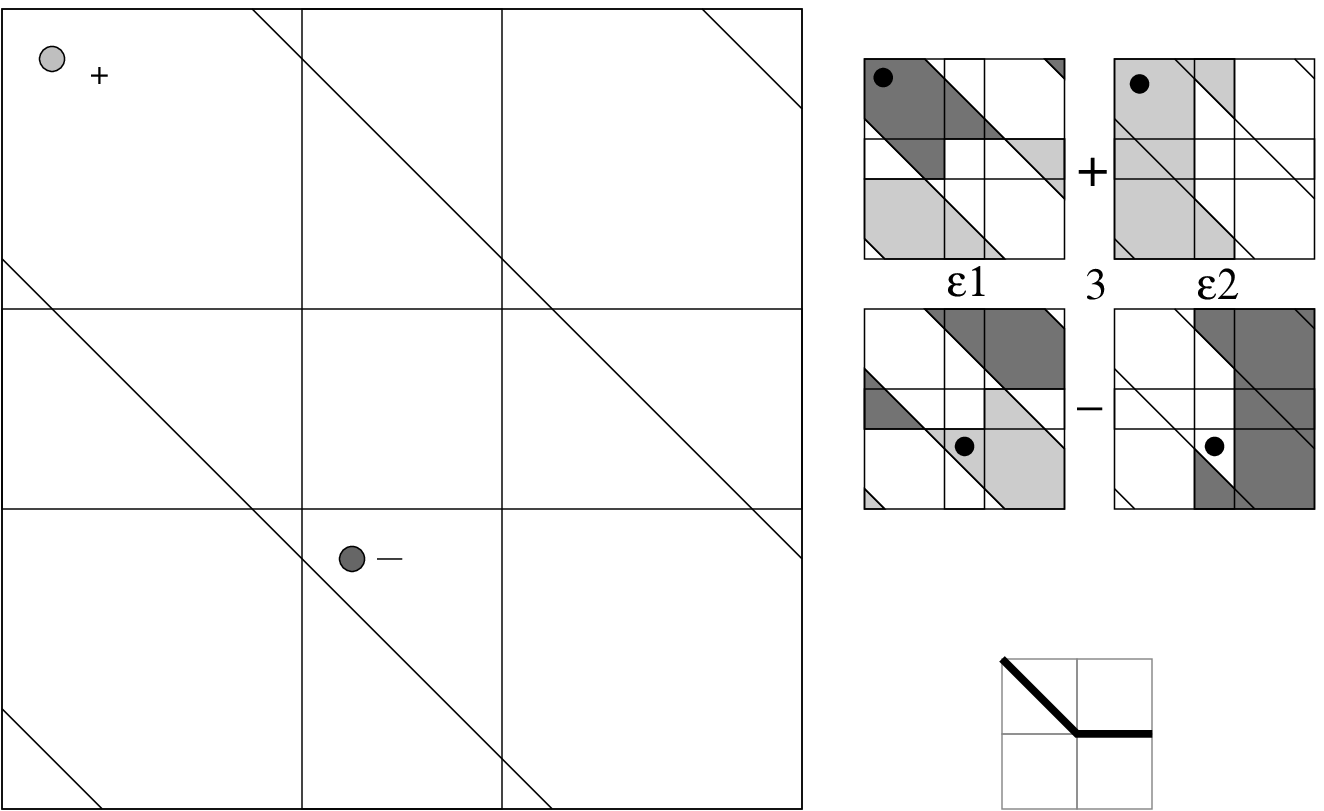}
\newline
{\bf Figure 6.4:\/} Points in the fiber.
\end{center} 

\subsection{A Singular Example}
\label{singular}

Sometimes it is an annoyance to deal with the tiny positive constant
$\alpha$ that arises in the definition of the
fundamental map.  In this section we will explain an alternate method
for applying the Master Picture Theorem.
One situation where this alternate approach proves useful is when
we need to deal with the fibers at $z=\alpha$.  We much prefer to
draw the fibers at $z=0$, because these do not contain any tiny
polygonal regions.  All the pieces of the partition can be
drawn cleanly.  However, in order to make sense of the
Master Picture Theorem, we need to slightly redefine how
the partition defines the return map.

Our method is to
redefine our polygonal regions to include their {\it lower\/} edges.
A lower edge is an edge first encountered by a line of slope $1$.  Figure
6.5 shows what we have in mind.

\begin{center}
\psfig{file=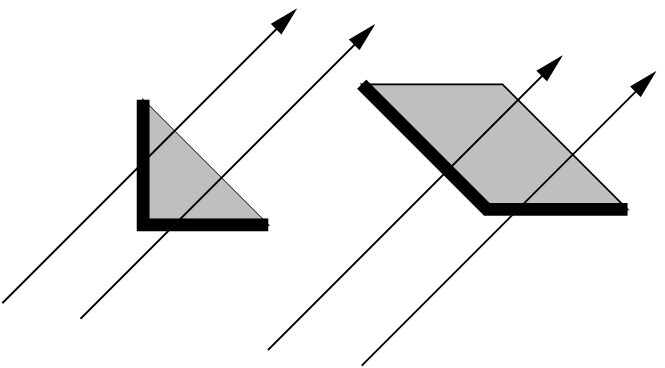}
\newline
{\bf Figure 6.5:\/} Polygons with their lower boundaries included.
\end{center} 

We then set $\alpha=0$ and determine the
relevent edges of the arithmetic graph by which {\it lower borded\/}
polygon contains our points.  if it happens that $z \in \{0,A,1-A\}$,
Then we think of the fiber at $z$ as being the geometric limit
of the fibers at $z+\epsilon$ for $\epsilon>0$.  That is, we
take a right-sided limit of the pictures.  When $z$ is not
one of these special values, there is no need to do this, for
the fiber is completely defined already.

We illustrate our approach with the example
$A=3/5$ and $(m,n)=(0,8)$.  We compute that
$t=8+\alpha$ in this case.  The relevant slices are the ones we get by
setting $z=\alpha$.  We deal with this by setting $\alpha=0$ and
computing
$$\mu_+(16,1)=(8,9,8) \equiv (\frac{4}{5},1,0) \hskip 5 pt {\rm mod\/} \hskip 5 pt \Lambda$$
$$\mu_-(16,-1)=(7,8,8) \equiv (0,\frac{7}{5},0) \hskip 5 pt {\rm mod\/} \hskip 5 pt \Lambda.$$
Figure 6.6 draws the relevant fibers.  The bottom right of Figure 6.6
shows the local structure of the arithmetic graph.  For instance,
$(\epsilon_1^+,\epsilon_2^+)=(0,1)$.

\begin{center}
\psfig{file=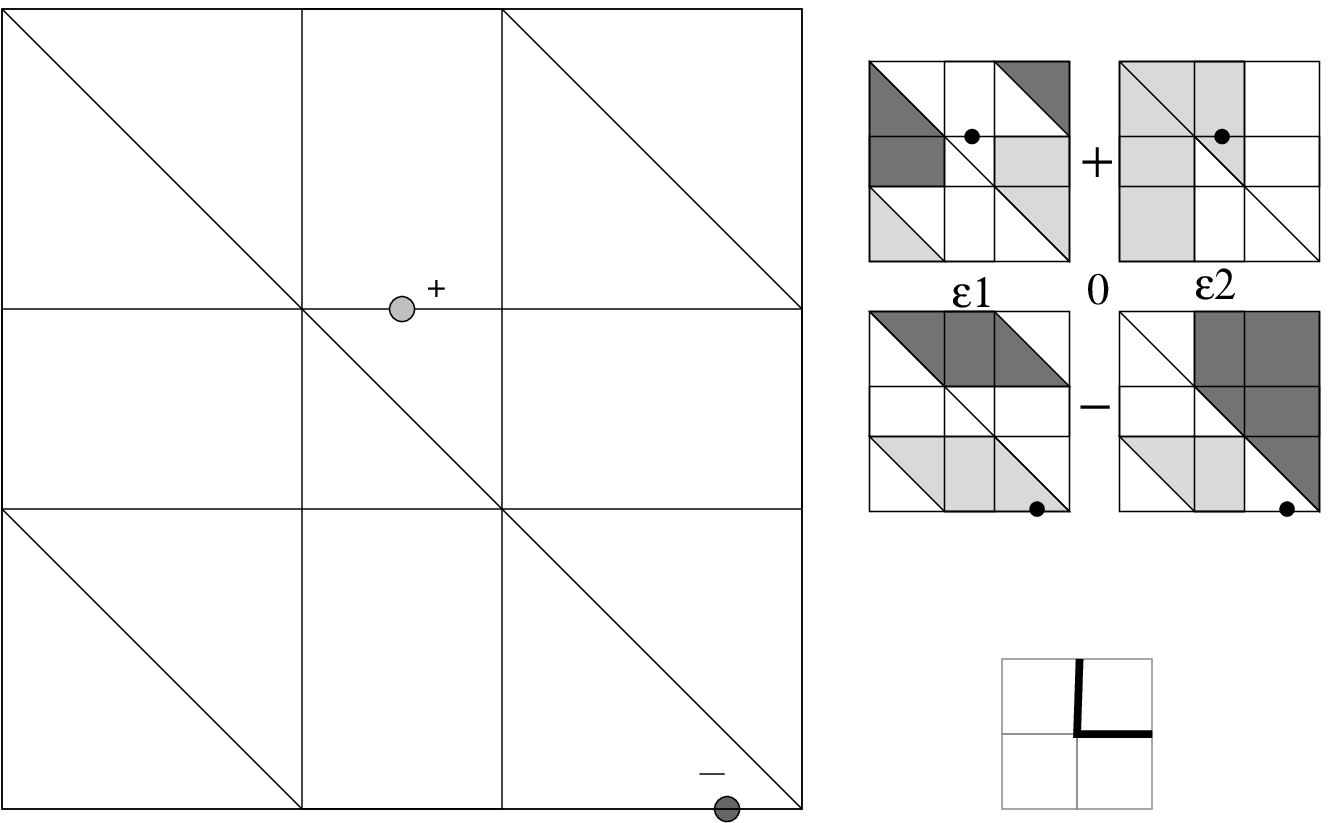}
\newline
{\bf Figure 6.6:\/} Points in the fiber.
\end{center} 

The only place where we need to use our special definition of a lower borded
polygon is for the point in the lower left fiber.  This fiber determines
the $x$ coordinate of the edge corresponding to $\mu_-$.  In this case, we
include our point in the lightly shaded parallelogram, because our
point lies in the lower border of this parallelogram.

There is one exception to our construction that requires an explanation.
Referring to the lower right fiber, suppose that the bottom point actually
was the bottom right vertex, as shown in Figure 6.7.  In this case,
the point is simultaneously the bottom left vertex, and we make the
definition using the bottom left vertex.  The underlying reason is that a tiny
push along the line of slope $1$ moves the point into the region on
the left.

\begin{center}
\psfig{file=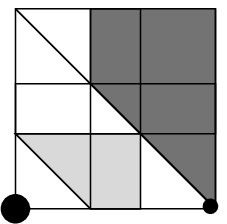}
\newline
{\bf Figure 6.7:\/} An exceptional case.
\end{center} 

\subsection{The Integral Structure}
\label{geometric2}

\subsubsection{An Affine Action}

We can describe Figure 6.3, and hence the Master Picture Theorem,
in a different way.
Let {\bf Aff\/} denote the $4$ dimensional affine group.
We define a discrete affine group action $\Lambda \subset {\bf Aff\/}$ on
the infinite slab
$\widetilde R=\R^3 \times (0,1)$. The group
$\Lambda$ is generated by the $3$ maps
$\gamma_1,\gamma_2,\gamma_3$.  Here
$\gamma_j$ acts on the first $3$ coordinates as translation
by the $j$th column of the matrix $\Lambda_A$, and on the
$4$th coordinate as the identity.  We think of the
$A$-variable as the $4$th coordinate.
Explicitly, we have
$$
\gamma_1\left[\matrix{x\cr y\cr z\cr A}\right]=\left[\matrix{x+1+A\cr y\cr z \cr A}\right]
$$
$$
\gamma_2 \left[\matrix{x\cr y\cr z\cr A}\right]=\left[\matrix{x+1-A\cr y+1+A\cr z\cr A}\right];
$$
\begin{equation}
\label{affine}
\gamma_3\left[\matrix{x\cr y\cr z \cr A}\right]=\left[\matrix{x-1\cr y-1\cr z+1\cr A}\right].
\end{equation}
These are all affine maps of $\R^4$.
The quotient $\widetilde R/\Lambda$ is
naturally a fiber bundle over $(0,1)$.  Each
fiber $(\R^3 \times \{A\})/\Lambda$ is isomorphic to
$\R^3/\Lambda_A$.  

The region $R$, from Equation \ref{FUND},
 is a fundamental domain for the action of
$\Lambda$. Note that $R$ is naturally an {\it integral polytope\/}.
That is, all the vertices of $R$ have integer coordinates.
$R$ has $16$ vertices, and they are as follows.
\begin{equation}
(\epsilon_1,\epsilon_2,\epsilon_3,0); \hskip 30 pt
(2\epsilon_1,2\epsilon_2,\epsilon_3,1); \hskip 30 pt
\epsilon_1,\epsilon_2, \epsilon_3 \in \{0,1\}.
\end{equation}

\subsubsection{Integral Polytope Partitions}

Inplicit in Figre 10.3 is the statement that the regions
$R_+$ and $R_-$ are partitioned into smaller convex polytopes.   
The partition is
defined by the $4$ families of hyperplanes discussed above.
An alternate point of view leads to a simpler partition.

For each pair $(\epsilon_1,\epsilon_2) \in \{-1,0,1\}$, we let
$R_+(\epsilon_1,\epsilon_2)$ denote the closure of the union of 
regions that assign $(\epsilon_1,\epsilon_2)$.
It turns out that
$R(\epsilon_1,\epsilon_2)$ if a finite union of convex integral
polytopes.  There are $14$ such polytopes, and they give an
integral partition
of $R_+$.  We list these polytopes in \S \ref{polytopelist}.

Let $\iota: R_+ \to R_-$ be given by the map
\begin{equation}
\label{iota}
\iota(x,y,z,A)=(1+A-x,1+A-y,1-z,A).
\end{equation}
Geometrically, $\iota$ is a reflection in the $1$-dimensional
line.   We have the general equation
\begin{equation}
R_-(-\epsilon_1,-\epsilon_2)=\iota(R_+(\epsilon_1,\epsilon_2)).
\end{equation}
Thus, the partition of $R_-$ is a mirror image of the partition
of $R_+$.   (See Example \ref{examplecalc} for an example
calculation.)

We use the action of $\Lambda$ to extend the partitions
of $R_+$ and $R_-$ to two integral polytope
tilings of $\widetilde R$.  (Again, see \S \ref{examplecalc} for
an example calculation.)   These $4$ dimensional tilings
determine the structure of the special orbits.

\subsubsection{Notation}
\label{mptnotation}

Suppose that $\widehat \Gamma$ is an arithmetic graph.
Let $M$ be the fundamental map associated to $\widehat \Gamma$.  We define
\begin{equation}
M_+=\mu_+ \circ M; \hskip 30 pt
M_-=\mu_- \circ \rho \circ M.
\end{equation}
Here $\rho$ is reflection in the $x$-axis.  Given a point $p \in \Z^2$, the
polytope of $R_+$ containing $M_+(p)$ determines the forward edge
of $\widehat \Gamma$ incident to $p$, and the polytope of
$R_-$ containing $M_-(p)$ determines the backward edge of
$\widehat \Gamma$ incident to $p$.  
Concretely, we have
$$
M_+(m,n)=(s,s+1,s) \hskip 10 pt {\rm mod\/} \Lambda;$$
$$
M_-(m,n)=(s-1,s,s) \hskip 10 pt {\rm mod\/} \Lambda;$$
\begin{equation}
s=Am+n+\alpha.
\end{equation}
As usual, $\alpha$ is the offset value.
Note that $\mu_+$ and $\mu_-$ only depend on the first
coordinate, and this first coordinate is not changed by $\rho$.
The map $\rho$ is present mainly for bookkeeping purposes,
because $\rho(\Xi_+)=\Xi_-$, and the domain of
$\mu_{\pm}$ is $\Xi_{\pm}$.

\newpage

\section{The Pinwheel Lemma}
\label{goodreturn}
\label{pinwheel}

\subsection{The Main Result}
\label{pinwheel2}

The Pinwheel Lemma gives a formula for the return
map $\Psi: \Xi \to \Xi$ in terms of maps we
call {\it strip maps\/}.
Similar objects are considered in [{\bf GS\/}] and
[{\bf S\/}].

Consider a pair $(\Sigma,L)$, where $\Sigma$ is an
infinite planar strip and $L$ is a line
transverse to $\Sigma$. 
The pair $(L,\Sigma)$ determines two vectors, $V_+$ and $V_-$,
each of which points from one boundary component of $\Sigma$ to
the other and is parallel to $L$.
Clearly $V_-=-V_+$.

For almost every point $p \in \R^2$, there is a
unique integer $n$ such that
\begin{equation}
E(p):=p+nV_+ \in \Sigma.
\end{equation}
We call $E$ the {\it strip map\/} defined relative to
$(\Sigma,L)$. 
The map $E$ is well-defined except on a countable collection of parallel
and evenly spaced lines.

\begin{center}
\psfig{file=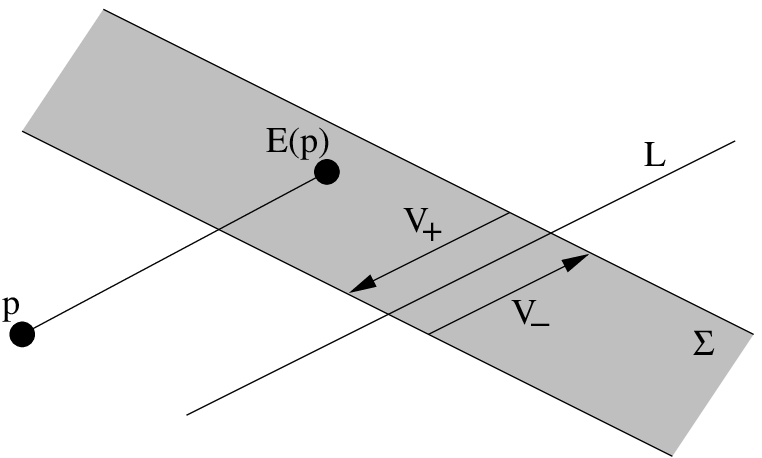}
\newline
{\bf Figure 7.1:\/} A strip map
\end{center} 

Figure 7.2 shows $4$ strips we associate to our kite.
To describe the strips in Figure 7.2 write
$(v_1,v_2,v_3)^t$ (a column vector) to
signify that $L=\overline{v_2v_3}$ and
$\partial \Sigma=\overline{v_1v_2} \cup I(\overline{v_1v_2})$,
where $I$ is the order $2$ rotation fixing $v_3$.
Here is the data for the
strip maps $E_1,E_2,E_3,E_4$. 

\begin{equation}
\label{strips}
\label{strip}
\left[\matrix{(-1,0) \cr (0,1) \cr (0,-1)}\right]; \hskip 15 pt
\left[\matrix{(A,0) \cr (0,-1) \cr (-1,0)}\right]; \hskip 15 pt
\left[\matrix{(0,1) \cr (A,0) \cr (-1,0)}\right]; \hskip 15 pt
\left[\matrix{(-1,0) \cr (0,-1) \cr (0,1)}\right].
\end{equation}
We set $\Sigma_{j+4}=\Sigma_j$ and $V_{j+4}=-V_j$.  Then
$\Sigma_{j+4}=\Sigma_j$. 
The reader can also reconstruct the strips from the information
given in Figure 7.2.  Figure 7.2 shows the parameter $A=1/3$, 
but the formulas in the picture are listed for general $A$.
In particular, the point $(3,0)$ is independent of $A$.
  Here is an explicit formula
for the vectors involved.
\begin{equation}
\label{vectors0}
V_1=(0,4); \hskip 15 pt V_2=(-2,2); \hskip 15 pt V_3=(-2-2A,0); \hskip 15 pt V_4=(-2,-2)
\end{equation}

\begin{center}
\resizebox{!}{4.6in}{\includegraphics{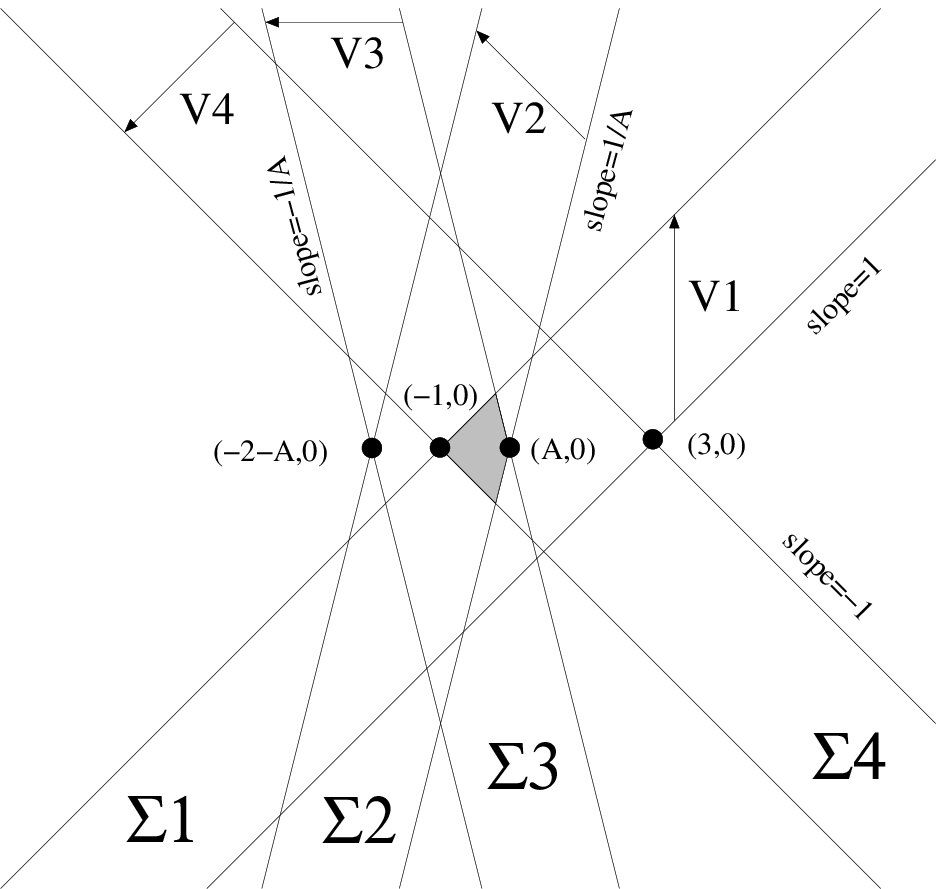}}
\newline
{\bf Figure 7.2:\/} The 4 strips for the parameter $A=1/3$.
\end{center} 

We also define a map $\chi: \R_+ \times \Z_{\rm odd\/} \to \Xi$ by the
formula
\begin{equation}
\chi(x,4n \pm 1)=(x,\pm 1)
\end{equation}

\begin{lemma}[Pinwheel]
$\Psi$ exists for any
point of $\Xi$ having a well-defined outer billiards orbit.
In all cases, $\Psi=\chi \circ (E_8...E_1)$.
\end{lemma}
We call the map in the Pinwheel Lemma the {\it pinwheel map\/}.
In \S \ref{pinwheelformulas} we give concrete formulas for this map.

\subsection{Some Corollaries}

Before we prove the Pinwheel Lemma, we list two corollaries.

\begin{corollary}
\label{proofparity}
The parity equation in Equation \ref{shortreturn} is true.
\end{corollary}

\startproof
The Pinwheel Lemma tells us that 
\begin{equation}
\Psi(x,1)-(x,1)=2(\epsilon_1 A+\epsilon_2,\epsilon_3); \hskip 30 pt
(\epsilon_1,\epsilon_2,\epsilon_3) \in \Z^2 \times \{-1,0,1\}.
\end{equation}
Given Equation \ref{vectors0}, we see that
the sum of the integer coefficients in each vector $V_j$ is
divisible by $4$. (For instance, $-2-2A$ yields $-2-2=-4$.)  Hence
$\epsilon_1+\epsilon_2+\epsilon_3$
is even.
\endproof

The Pinwheel Lemma gives a formula for the
quantities in Equation \ref{shortreturn}.

For $j=0,...,7$ we define points $p_{j+1}$ and integers $n_j$ by the
following equations.
\begin{equation}
p_{j+1}=E_{j+1}(p_{j})=p_{j}+n_j V_{j+1}.
\end{equation}
Given the equations
\begin{equation}
V_1=(0,4); \hskip 15 pt V_2=(-2,2); \hskip 15 pt V_3=(-2-2A,0); \hskip 15 pt V_4=(-2,-2)
\end{equation}
we find that
\begin{equation}
\label{spectrum2}
\epsilon_1=n_2-n_6; \hskip 30 pt
\epsilon_2=n_1+n_2+n_3-n_5-n_6-n_7; \hskip 40 pt
\end{equation}
We call $(n_1,...,n_7)$ the {\it length spectrum\/} of $p_0$.

The precise bound in Equation \ref{shortreturn} follows from
the Master Picture Theorem, but here we give a heuristic explanation.
If we define
\begin{equation}
m_1=n_7; \hskip 30 pt
m_2=n_6; \hskip 30 pt
m_3=n_5; \hskip 30 pt
\end{equation}
then we have
\begin{equation}
\epsilon_1(p)=n_2-m_2; \hskip 30 pt
\epsilon_2(p)=(n_1-m_1)+(n_2-m_2)+(n_3-m_3).
\end{equation}
The path with vertices $p_0,p_1,...,p_7,p_8,\chi(p_8)$
uniformly close to an octagon with dihedral symmetry.
See Figure 7.3 below.
For this reason, there is a universal bound to
$|n_i-m_i|$.  This is a heuristic explanation of the
bound in Equation \ref{shortreturn}.

\subsection{The Simplest Case}
\label{pinwheelsimple}

Here we prove the Pinwheel for points of $\Xi$ far from $K$.
Figure 7.3 shows a decomposition of
$\R^2-K'$ into $8$ regions, $S_0,...,S_7$.   Here $K'$
is a suitably large compact set.
Let $V_1,...,V_4$ we the vectors associated to our special
strip maps.  We set $V_{4+j}=-V_j$.  A calculation
shows that
\begin{equation}
\label{succession}
x \in S_j; \hskip 30 pt \Longrightarrow \hskip 30 pt
\psi(x)-x=V_j.
\end{equation}
One can easily see this using Billiard King or else our
interactive guide to the monograph.

\begin{center}
\psfig{file=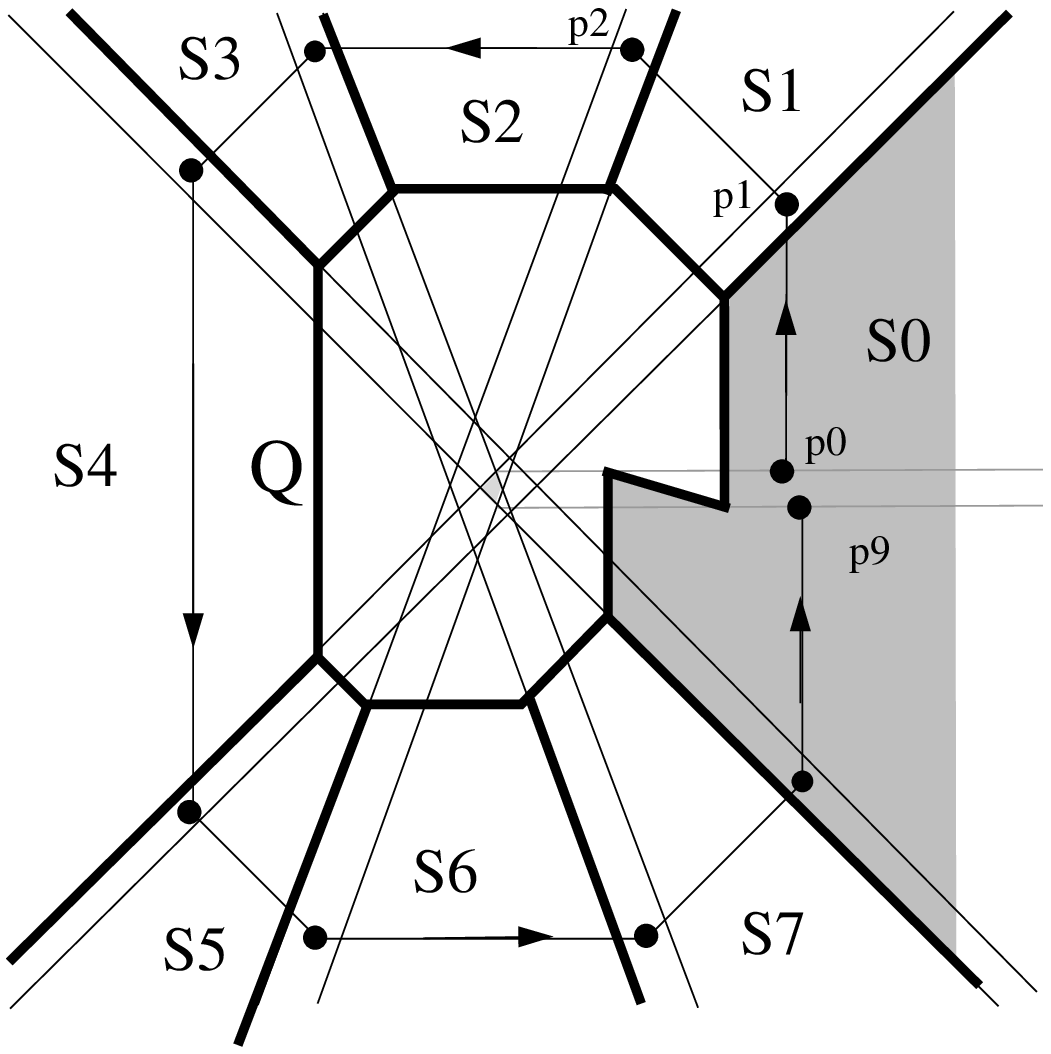}
\newline
{\bf Figure 7.3:\/} The Simplest Sequence of Regions
\end{center} 

Equation \ref{succession} tells the whole story for points
of $\Xi$ far away from $K$.  As above, let
$p_{j+1}=E_{j+1}(p_j)$ for $j=0,...7$.  let
$p_{9}=\chi(p_8)$. Here we have set $E_{j+4}=E_j$.
By induction and Equation \ref{succession},
$p_{j+1}$ lies in the forward orbit of $p_j$ for
each $j=0,...,8$.  
But then $p_{9}=\Psi(p_1)=\chi \circ E_8...E_1(p_0)$.

\subsection{Discussion of the General Case}

As we have just seen, the Pinwheel Lemma is a
fairly trivial result for points that are far from
the origin.  For points near the origin, the
Pinwheel Lemma is a surprising and nontrivial result.
In fact, it only seems to work because of a lucky
accident.  The fact that we consider the Pinwheel
Lemma to be an accident probably means that we
don't yet have a good understanding of what is
going on.  

Verifying the Pinwheel Lemma for any
given parameter $A$ is a finite calculation.
We just have to check, on a fine enough mesh of
points extending out sufficiently far away
from $K(A)$, that the equation in the
Pinwheel Lemma holds.  The point is that all the
maps involved are piecewise isometries for each parameter.
We took this approach in [{\bf S\/}] when we proved
the Pinwheel Lemma for $A=\phi^{-3}$.

Using Billiard King, we computed that the Pinwheel Lemma
holds true at the points $(x,\pm 1)$ relative
to the parameter $A$ for all
$$A=\frac{1}{256},...,\frac{255}{256}; \hskip 25 pt
x=\epsilon+\frac{1}{1024},...,\epsilon+\frac{16384}{1024};
\hskip 20 pt
\epsilon=10^{-6}.$$
The tiny number $\epsilon$ is included to
make sure that the outer billiards orbit
is actually defined for all the points we sample.
This calculation does not constitute a proof of
anything.  However, we think that it serves as
a powerful sanity check that the Pinwheel
Lemma is correct.
We have fairly well carpeted the region of doubt
about the Pinwheel Lemma with instances of its truth.

Our proof of the Pinwheel Lemma essentially
boils down to finding the replacement equation
for Equation \ref{succession}.  We will do this
in the section.  As the reader will
see, the situation in general is much more
complicated.  There is a lot of information packed
into the next section, but all this information
is easily seen visually on Billiard King.
We have programmed Billiard King so that the
reader can see pictures of all the regions
involved, as well as their interactions,
for essentially any desired parameter.

We think of the material in the next
section as something like a written description
of a photograph.
The written word is probably not the right medium for the proof
of the Pinwheel Lemma.  To put this in a different way,
Billiard King relates to the proof given here much
in the same way that an ordinary research paper
would relate to one that was written in crayon.

\subsection{A Partition of the Plane}
\label{partition}

Let $\psi=\psi_A$ be the square of the outer billiards map
relative to $K(A)$.
For each $x \in \R^2-K$ on which
$\psi$ is defined, there is a vector
$v_x$ such that $$\psi(x)-x=v_x.$$  This vector
is twice the difference between $2$ vertices of $K$,
and therefore can take on $12$ possible values. 
It turns out that $10$ of these values occur.  We call these
vectors $V_j$, with $j=1,2,3,4,4^{\sharp},5,6^{\flat},6,7,8$.
With this ordering, the argument of $V_j$ increases
monotonically with $j$.  Compare Figure 7.4.
For each of our vectors $V$, there is an open region
$R \subset \R^2-K$ such that $x \in R$ if and only
if $\psi(x)-x=V$.  The regions
$R_1,...,R_8$ are unbounded.
The two regions $R_4^{\sharp}$ and $R_6^{\flat}$ are bounded.

One can find the entire partition by extending the
sides of $K$ in one direction, in a pinwheel
fashion, and then pulling back these rays by
the outer billiards map.  To describe the regions, we use the notation
$\overrightarrow{q_1},p_1,...,p_k,\overrightarrow{q_2}$
to indicate that 
\begin{itemize}
\item The two unbounded edges are $\{p_1+tq_1|\ t \geq 0\}$ and
$\{p_k+tq_2|\ t \geq 0\}$.
\item $p_2,...,p_{k-1}$ are any additional intermediate vertices.
\end{itemize}
To improve the typesetting on our list, we set $\lambda=(A-1)^{-1}$.
Figure 6.3 shows the picture for $A=1/3$.  The
reader can see any parameter using Billiard King.
\newline
\newline
\begin{tabular}{ll}
 $V_1\!=\!(0,4)$.& $R_1: \overrightarrow{(1,-1)},(1,-2),\overrightarrow{(1,1)}$.
\\ $V_2\!=\!(-2,2)$.& $R_2: \overrightarrow{(1,1)},(1,-2),(0,-1),\overrightarrow{(A,1)}$.
\\ $V_3\!=\!(-2\!-\!2A,0)$ &$R_3: \overrightarrow{(A,1)},(2A,1),\lambda(2A^2,-1-A),\overrightarrow{(-A,1)}.$
\\ $V_4\!=\!(-2,-2)$&$R_4: \overrightarrow{(-A,1)},\lambda(2A,A-3),\overrightarrow{(-1,1)}.$
\\ $V_{4^{\sharp}}\!=\!(-2A,-2)$&$R_{4^{\sharp}}: (A,0),(2A,1),\lambda(2A^2,-1-A))$
\\ $V_5\!=\!(0,-4)$ &  $R_5: \overrightarrow{(-1,1)},\lambda (2A,A\!-\!3),(-A,2),\lambda(2A,3A\!-\!1),\overrightarrow{(-1,-1)}$
\\ $V_{6^{\flat}}\!=\!(2A,-2)$ &$R_{6^{\flat}}: (0,1),(-A,2),\lambda(2A,3A-1)$
\\ $V_6\!=\!(2,-2)$&$R_6: \overrightarrow{(-1,-1)},\lambda(2,A+1),\overrightarrow{(-A,-1)}$
\\ $V_7\!=\!(2+2A,0)$&$R_7: \overrightarrow{(-A,-1)},\lambda(2,A+1),(-2,-1),\overrightarrow{(A,-1)}$
\\ $V_8\!=\!(2,2)$&$R_8: \overrightarrow{(A,-1)},(-2,-1),(-1,0),\overrightarrow{(1,-1)}$.
\end{tabular}
\newline  

\begin{center}
\resizebox{!}{5.2in}{\includegraphics{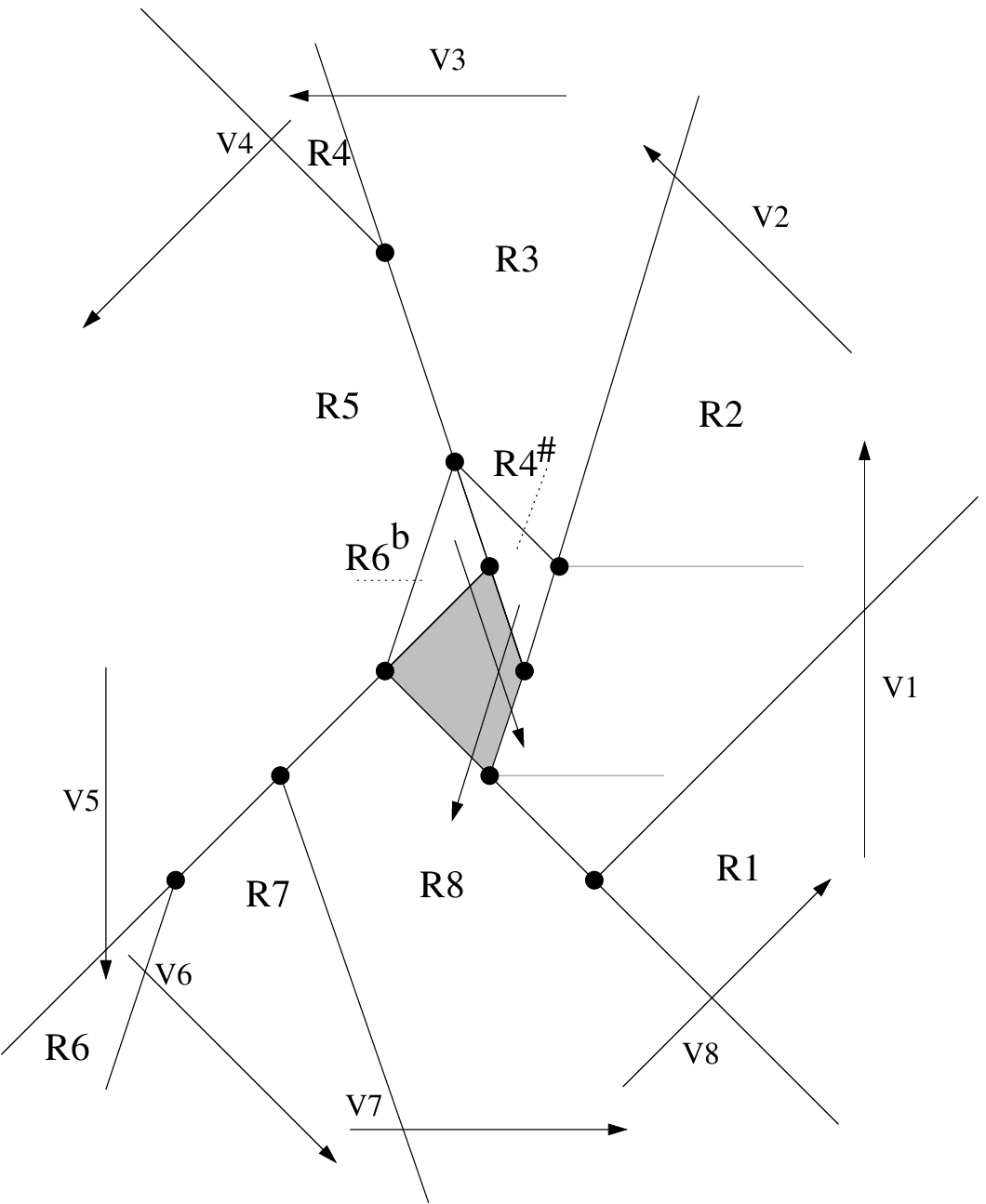}}
\newline
{\bf Figure 7.4:\/} The Partition for $A=1/3$.
\end{center} 

It is convenient to set
\begin{equation}
\widehat R_a=R_a+V_a=\{p+V_a|\ p \in R_a\}.
\end{equation}
One symmetry of the partition
is that reflection in the $x$-axis interchanges
$\widehat R_a$ with $R_{10-a}$, for all values of $a$.  (To make
this work, we set $R_9=R_1$,
and use the convention $4^{\sharp}+6^{\flat}=10$.)

We are interested in {\it transitions\/} between one region $R_a$, and another
region $R_b$.
If $\widehat R_a \cap R_b \not = \emptyset$ for
some parameter $A$ it means that there is some
$p \in R_a$ such that $\psi_A(p) \in R_b$.
(We think of our regions as being open.) 
We create a {\it transition matrix\/} using the following rules.
\begin{itemize}
\item A $0$ in the $(ab)$th spot indicates that 
$\widehat R_a \cap R_b=\emptyset$ for all $A \in (0,1)$.
\item A $1$ in the $(ab)$th spot indicates that 
$\widehat R_a \cap R_b\not =\emptyset$ for all $A \in (0,1)$.
\item A $t^+$ in the $(ab)$th spot indicates that
$R_a \cap R_b\not = \emptyset$ iff $A \in (t,1)$. 
\item A $t^-$ in the $(ab)$th spot indicates that
$R_a \cap R_b\not = \emptyset$ iff $A \in (0,t)$.
\end{itemize}

\begin{equation}
\matrix{
& R_1& R_2& R_3& R_4
& R_{4^{\sharp}}& R_5& R_{6^{\flat}}&R_6& R_7& R_8 \cr
\widehat R_1&1&1&(\frac{1}{3})^+&0&0&0&0&0&0&0 \cr
\widehat R_2&0&1&1&(\frac{1}{3})^-&1&1&1&0&0&0 \cr
\widehat R_3&0&0&1&1&(\frac{1}{2})^+&1&0&0&0&0 \cr
\widehat R_4&0&0&0&1&0&1&0&0&0&0 \cr
\widehat R_{4^{\sharp}}&0&0&0&0&0&(\frac{1}{3})^+&(\frac{1}{3})^+&0&0&1 \cr
\widehat R_5&0&0&0&0&0&1&(\frac{1}{3})^+&1&1&1 \cr
\widehat R_{6^{\flat}}&0&1&0&0&0&0&0&0&(\frac{1}{2})^+&1 \cr
\widehat R_6&0&0&0&0&0&0&0&0&1&(\frac{1}{3})^- \cr
\widehat R_7&(\frac{1}{3})^+&1&0&0&0&0&0&0&0&1 \cr
\widehat R_8&1&1&1&0&1&0&0&0&0&0}
\end{equation}
We have programmed Billiard King so that the interested reader can
see each of these relations at a single glance.  Alternatively,
they can easily be established using routine linear algebra.
For example, interpreting $\widehat R_3$ and $R_2$ as 
projectivizations of open convex cones $\widehat C_3$ and $C_2$
in $\R^3$, we easily
verifies that the vector $(-1,A,-2-A)$ has positive
dot product with all vectors in $\widehat C_3$ and
negative dot product with all vectors in $C_2$.
Hence $R_2 \cap \widehat R_3=\emptyset$.

We can relate all the nonempty intersections to our strips.  
As with the list of intersections,
everything can be seen at a glance using Billiard King, or else
proved using elementary linear algebra. 
First we list the intersections that comprise the
complements of the strips.

\begin{itemize}
\item $\widehat R_2 \cap R_2$ and $\widehat R_6 \cap R_6$ are the components
of $\R^2-(\Sigma_1 \cup \Sigma_2)$.
\item $\widehat R_4 \cap R_4$ and $\widehat R_8 \cap R_8$ are the components
of $\R^2-(\Sigma_3 \cup \Sigma_4)$.
\item $\widehat R_3 \cap (R_3 \cup R_{4^{\sharp}})$ and
$(\widehat R_{6^{\flat}} \cup \widehat R_7) \cap R_7$ are the components of
$\R^2-(\Sigma_2 \cup \Sigma_3)$.
\item $\widehat R_1 \cap R_1$ and
$(\widehat R_{4^{\sharp}} \cup \widehat R_5) \cap (R_5 \cup R_{6^{\flat}})$ are
the components of $\R^2-(\Sigma_1 \cup \Sigma_4)$.
\end{itemize}

Now we list the intersections that are contained in single strips.  To make our
typesetting nicer, we use the term $u$-{\it component\/} to denote an
unbounded connected component.  We use the term $b$-component to denote a
bounded connected component.

\begin{itemize}
\item $\widehat R_1 \cap R_2$ and $\widehat R_5 \cap R_6$ are the two $u$-components
of $\Sigma_1-(\Sigma_2 \cup \Sigma_4)$.
\item $\widehat R_8 \cap R_1$ and $\widehat R_4 \cap R_5$ are the two $u$-components
of $\Sigma_4-(\Sigma_1 \cup \Sigma_3)$.
\item $\widehat R_3 \cap R_4$ and $(\widehat R_{6^{\flat}} \cup \widehat R_7) \cap R_8$ are
the two $u$-components of $\Sigma_3-(\Sigma_2 \cup \Sigma_4)$.
\item $\widehat R_6 \cap R_7$ and $\widehat R_2 \cap (R_3 \cup R_{4^{\sharp}})$ are
the two $u$-components of $\Sigma_2-(\Sigma_1 \cup \Sigma_3)$.
\item $\widehat R_{6^{\flat}} \cap R_7$ is contained in the $b$-component of
$\Sigma_1-(\Sigma_2 \cup \Sigma_3)$.
\item $\widehat R_3 \cap R_{4^{\sharp}}$ is contained in the $b$-component of
$\Sigma_4-(\Sigma_2 \cup \Sigma_3)$.
\item $\widehat R_{4^{\sharp}} \cap (R_5 \cup R_{6^{\flat}})$ is contained in
the $b$-component of $\Sigma_3-(\Sigma_1 \cup \Sigma_4)$.
\item $(\widehat R_{4^{\sharp}} \cup \widehat R_5) \cap R_{6^{\flat}}$ is contained
in the $b$-component of $\Sigma_2-(\Sigma_1 \cup \Sigma_4)$.
\end{itemize}

Now we list the intersections of regions that are contained in double
intersections of strips.   In
this case, all the components are bounded:  Any two strips intersect
in a bounded region of the plane.

\begin{itemize}
\item $\widehat R_1 \cap R_3$ and $\widehat R_5 \cap R_7$ are the components of
$(\Sigma_1 \cap \Sigma_2)-(\Sigma_3 \cup \Sigma_4)$.
\item $\widehat R_7 \cap R_1$ and $\widehat R_3 \cap R_5$ are the components of
$(\Sigma_3 \cap \Sigma_4)-(\Sigma_1 \cup \Sigma_2)$.
\item $\widehat R_2 \cap R_4$ and $\widehat R_6 \cap R_8$ are
bounded components of $(\Sigma_2 \cap \Sigma_3)-(\Sigma_1 \cup \Sigma_4)$.
\item $\widehat R_8 \cap R_2=(\Sigma_1 \cap \Sigma_4)-(\Sigma_2 \cup \Sigma_3)$.
\end{itemize}

Now we list all the intersections of regions that are contained in
triple intersections of strips.
\begin{itemize}
\item $\widehat R_2 \cap (R_5 \cup R_{6^{\flat}})=\Sigma_2 \cap \Sigma_3 \cap \Sigma_4-\Sigma_1$.
\item $\widehat R_8 \cap (R_3 \cup R_{4^{\sharp}})=\Sigma_1 \cap \Sigma_2 \cap \Sigma_4-\Sigma_3$.
\item $(\widehat R_{4^{\sharp}} \cup \widehat R_5) \cap R_8=\Sigma_1 \cap \Sigma_2 \cap \Sigma_3-\Sigma_4$.
\item $(\widehat R_{6^{\flat}} \cup \widehat R_7) \cap R_2=\Sigma_1 \cap \Sigma_3 \cap \Sigma_4-\Sigma_2$.
\end{itemize}

Here we list a bit more information about the two
regions $R_{4^{\sharp}}$ and $R_{6^{\flat}}$ some of
the information is redundant, but it is useful to
have it all in one place.
\begin{itemize}
\item $R_{4^{\sharp}} \subset \Sigma_4-\Sigma_3$.
\item $R_{4^{\sharp}}+V_3=\Sigma_3-(\Sigma_2 \cup \Sigma_4)$.
\item $R_{6^{\flat}} \subset \Sigma_2-\Sigma_1$.
\item $R_{6^{\flat}}+V_5 \subset \Sigma_1-\Sigma_2$.
\item $R_{6^{\flat}}+V_5-V_6=\Sigma_2-(\Sigma_1 \cup \Sigma_3)$.
\end{itemize}

Finally, we mention two crucial relations between our various vectors:
\begin{itemize}
\item $V_3-V_4+V_5=V_{4^{\sharp}}$.
\item $V_5-V_6+V_7=V_{6^{\flat}}$.
\end{itemize}
These two relations are responsible for the lucky cancellation
that makes the Pinwheel Lemma hold near the kite.

We will change our notation slightly from the simplest case
considered above.
Given any point $z_1 \in \Xi$, we can
associate the {\it sequence of regions\/}
\begin{equation}
\label{general}
R_{a_1} \to ... \to R_{a_k}
\end{equation}
through which the forwards orbit of $z_1$
transitions until it returns as $\Psi(z_1)$.
The simplest possible sequence is the one where $a_j=j$ for $j=1,...,9$.
See Figure 7.2.  We already analyzed this case above.
We let $z_j$ denote the first point in the
forward orbit of $z_1$ that lies in $R_{a_j}$.

To prove the Pinwheel Lemma in general, we need to analyze
all allowable sequences and see that the equation in the
Pinwheel Lemma always holds. We will break the set of
all sequences into three types, and then analyze the types
one at a time.  Here are the types.
\begin{enumerate}
\item Sequences that do not involve the indices $4^{\sharp}$ or $6^{\flat}$.
\item Sequences that involve $4^{\sharp}$ but not $6^{\flat}$.
\item Sequences that involve $6^{\flat}$.
\end{enumerate}

\subsection{No Sharps or Flats}

\begin{lemma}
\label{intersect}
\label{intersect1}
If $j<k$ then
$\widehat R_j \cap R_{k} \subset \Sigma_j \cap ... \cap \Sigma_{k-1}.$
\end{lemma}

\startproof 
This is a corollary of the
the intersections listed above.
\endproof

Suppose by induction we have shown that
\begin{equation}
\label{induct}
z_{j}=E_{a_j-1}E_{a_j-2}...E_1(z_1).
\end{equation}
By construction and Lemma \ref{intersect},
$$z_{j+1}=E_{a_j}(z_j) \in \widehat R_{a_j} \cap R_{a_{j+1}}
\subset \Sigma_{a_j} \cap ... \cap \Sigma_{a_{j+1}-1}.$$
Therefore, 
$E_{a_j},...,E_{a_{j+1}-1}$ all act trivially on $z_{j+1}$, forcing
$$z_{j+1}=E_{a_{j+1}-1}E_{a_{j+1}-2}...E_1(z_1).$$
Hence, Equation \ref{induct} holds true for all
indices $j$. 

By the Intersection Lemma, we eventually reach
either a point $z_9$ or $z_{10}$.  (That is, we
wrap all the way around and return either to $R_9=R_1$ or
else to $R_{10}=R_2$.)  We will consider these two cases
one at a time.
\newline
\newline
{\bf Case 1:\/}  If we reach
$z_9=(x_9,y_9) \in R_9$ then we have
\begin{equation}
\label{gr1}
z_9=E_8...E_1(z_1); \hskip 30 pt
x_9>0;\hskip 30 pt y_9 \leq 1.
\end{equation}
From this we get that 
$\Psi(z_1)=\chi \circ (E_4...E_1)^2(z_1)$, as desired.
The last inequality in Equation \ref{gr1} requires
explanation.  By the Intersection
Lemma,  the point preceding $z_9$ on our
list must lie in $R_a$ for some
$a \in \{6^{\flat},6,7,8\}$.  However, the distance
between any point on
$\R_+ \times \{3,5,7...\}$ to any point in 
$R_a$ exceeds 
the length of vector $V_a$.
\newline
\newline
{\bf Case 2:\/}
If we arrive at $z_{10}=(x_{10},y_{10})$, then the
Intersection Lemma tells us that the point
preceding $z_{10}$ lies in $V_a$ for
$a=\{6^{\flat},6,7,8\}$ and $z_{10} \in \Sigma_9$.
Hence $E_9(z_{10})=z_{10}$.  That is
$$z_{10}=E_8...E_1(z_1); \hskip 20 pt
x_{10}>0; \hskip 30 pt y_{10}<3.$$
The last inequality works just as in Case 1.
All points in
$R_{10}$ have $y$-coordinate at least $-2$. Hence
$y_{10}=\pm 1$.  Hence
$\chi(z_{10})=z_{10}$.  Putting everything together
gives the same result as Case 1.

\subsection{Dealing with Four Sharp}
\label{4sharp}

In this section we will deal with orbits whose associated
sequence has a $4^{\sharp}$ in it, but not a
$6^{\flat}$.  The following result is an immediate
consequence the intersections discussed above.

\begin{lemma}
\label{intersect2} 
The following holds for all parameters.
$$
\widehat R_{4^{\sharp}} \cap R_{4^{\sharp}}=\emptyset;
\hskip 10pt R_{4^{\sharp}} \subset \Sigma_4\!-\!\Sigma_3;
\hskip 10pt R_{4^{\sharp}} + V_3 \in \Sigma_3\!-\!\Sigma_4;
\hskip 10pt \widehat R_{4^{\sharp}}\cap R_8 \subset \Sigma_1 \cap \Sigma_2 \cap \Sigma_3
$$
\end{lemma}

Let $z$ be the first point
in the forward orbit of $z_1$ such
that $z \in R_{4^{\sharp}}$.   Using Lemma
\ref{intersect} and the same analysis as in the
previous section, we get
\begin{equation}
\exists n \in \N \cup \{0\} \hskip 40 pt
z=E_2E_1(z_1)+nV_3,
\end{equation}
From Lemma
\ref{intersect} and Item 1 of Lemma \ref{intersect2}, the next point in
the orbit is
\begin{equation}
w=z+V_{4^{\sharp}} \in R_5 \cup R_8.
\end{equation}

Items 2 and 3 of Lemma \ref{intersect2} give
$$
E_3E_2E_1(z_1)=E_3(z)=z+V_3; \hskip 30 pt
E_4E_3(z)=z+V_3\!-\!V_4.$$
Figure 7.5 shows what is going on.
Since
$V_3-V_4+V_5=V_{4^{\sharp}}$,
\begin{equation}
w =z+V_{4^{\sharp}}=z+V_3-V_4+V_5 = E_4E_3(z)+V_5=E_4E_3E_2E_1(z_1)+V_5.
\end{equation}
The rest of the analysis
is as in the previous section.  We use Item 4 of
Lemma \ref{intersect2} as an {\it addendum\/} to
Lemma \ref{intersect1} in case $w \in R_8$.

\begin{center}
\psfig{file=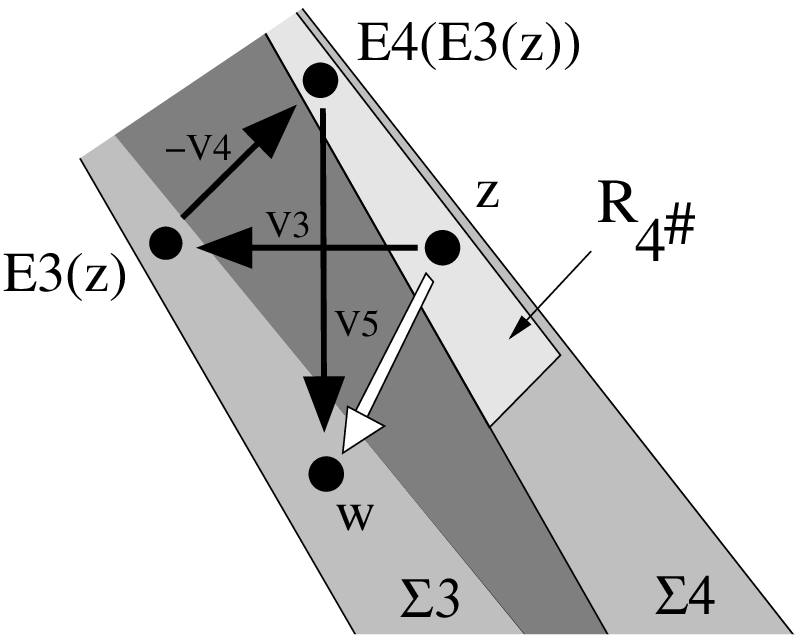}
\newline
{\bf Figure 7.5:\/} The orbit near $R_{4^{\sharp}}$.
\end{center} 

\subsection{Dealing with Six Flat}
\label{6flat}

Here is another immediate consequence of the intersections listed above.

\begin{lemma}
\label{intersect3}
The following is true for all parameters.
$$\hskip 10 ptV_{6^{\flat}} \subset \Sigma_2-\Sigma_1;
\hskip 10 ptV_{6^{\flat}}+R_5 \subset \Sigma_1-\Sigma_2;$$
$$
\widehat R_{6^{\flat}} \cap R_{2} \subset \Sigma_3 \cap \Sigma_4 \cap \Sigma_1;
\hskip 10 pt
\widehat R_2 \cap R_{6^{\flat}} \subset \Sigma_2 \cap \Sigma_3 \cap \Sigma_4;
$$
\end{lemma}

Let $z$ be the first point in the forwards orbit
of $z_1$ such that $z \in R_{6^{\flat}}$ and
let $w=\psi(z)$.  The same arguments
as in the previous section give
\begin{equation}
z=E_4E_3E_2E_1(z_1)+nV_5; \hskip 40 pt
w=z+V_{6^{\flat}} \in R_7 \cup R_2.
\end{equation}
Here $n \in \N \cup \{0\}$.
(The possibility of $w \in R_{6^{\flat}}$ is ruled out
by Item 1 of Lemma \ref{intersect2} and the
reflection symmetry.)
Items 2 and 3 of Lemma \ref{intersect3} give
$$E_5E_4E_3E_2E_1(z_1)=E_5(z)=z+V_5; \hskip 20 pt
E_6E_5E_4E_3E_2E_1(z)=z+V_5-V_6$$
Figure 7.6 shows what is going on.
Since
$V_5-V_6+V_7=V_{6^{\flat}}$,
\begin{equation}
w=E_6E_5E_4E_3E_2E_1(z)+V_7.
\end{equation}
The rest of the analysis is as in
the previous cases.  We use Item 3 of Lemma \ref{intersect3} as an
{\it addendum\/} to Lemma \ref{intersect1} in case $w \in R_{2}$.

\begin{center}
\psfig{file=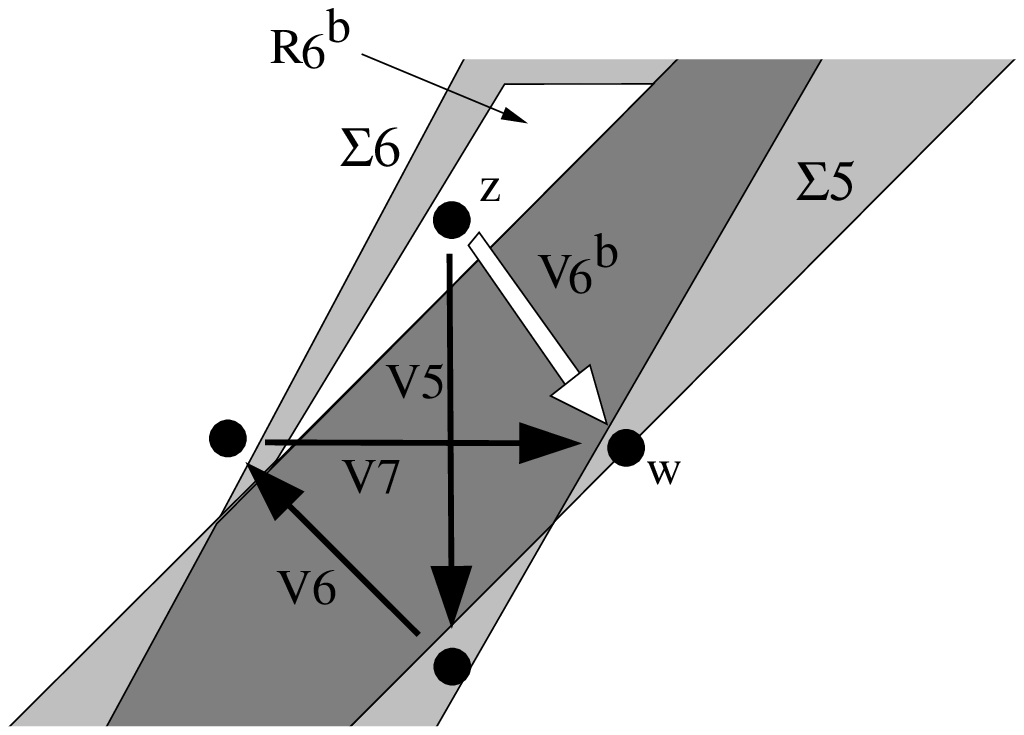}
\newline
{\bf Figure 7.6:\/}  The orbit near $R_{6^{\flat}}$.
\end{center} 

\newpage

\section{The Torus Lemma}
\label{torus}

\subsection{The Main Result}
\label{domainX}

For ease of exposition, we state and prove the
$(+)$ halves of our results. The $(-)$ halves
have the same formulation and proof. 

Let $T^4=\widetilde R/\Lambda$, the $4$ dimensional
quotient discussed in \S \ref{geometric2}.  Topologically,
$T^4$ is the product of a $3$-torus with $(0,1)$.
Let $(\mu_+)_A$ denote the map $\mu_+$ as defined for
the parameter $A$.  We now define
$\mu_+: \Xi_+ \times (0,1) \to T^4$ by
the obvious formula
$\mu_+(p,A)=((\mu_+)_A(p),A)$.
We are just stacking all these maps together.

The Pinwheel Lemma tells us that
$\Psi(p)=\chi \circ E_8...E_1(p)$
whenever both maps are defined.
This map involves the sequence
$\Sigma_1,...,\Sigma_8$ of strips.
We are taking indices mod $4$ so that 
$\Sigma_{j+4}=\Sigma_j$ and
$E_{4+j}=E_j$. 
Let $p \in \Xi_+$.
We set $p_0=p$ and indctively define
\begin{equation}
p_j=E_j(p_{j-1}) \in \Sigma_j.
\end{equation}
We also define
\begin{equation}
\theta(p)=\min \theta_j(p); \hskip 30 pt
\theta_j(p)={\rm distance\/}(p_j,\partial \Sigma_j).
\end{equation}
The quantity $\theta(p)$ depends on the parameter
$A$, so we will write $\theta(p,A)$ when we want to
be clear about this.

\begin{lemma}[Torus]
Let $(p,A), (q^*,A^*) \in \Xi_+ \times (0,1)$.
There is some $\eta>0$, depending only on
$\theta(p,A)$ and $\min(A,1-A)$, with 
the following property. Suppose that the pinwheel map
is defined at $(p,A)$.  Suppose also that
$\mu_+(p,A)$ and $ \mu_+(q^*,A^*)$ are within
$\eta$ of each other.  Then the pinwheel map
is defined at $(q^*,A^*)$ and
$(\epsilon_1(q^*),\epsilon_2(q^*))=(\epsilon_1(p),\epsilon_2(p))$.
\end{lemma}

\noindent
{\bf Remark:\/}  My proof
of the Torus Lemma owes a big intellectual debt to many sources.
I discovered the Torus Lemma experimentally, but I got some inspiration 
for its proof by reading [{\bf T2\/}], an account
of unpublished work by Chris Culter about the existence of periodic
orbits for polygonal outer billiards.  Culter's proof is closely related to
ideas in [{\bf K\/}].  
The paper [{\bf GS\/}] implicitly has some of these same ideas, though
they are treated from a different point of view.   
If all these
written sources aren't enough, I was also influenced by some
conversations with John Smillie.

\subsection{Input from the Torus Map}
\label{mapinput}

We first prove the Torus Lemma under the assumption that $A=A^*$.  We set
$q=q^*$.
In this section, we explain the significance of the map $\mu_+$.
We introduce the quantities
\begin{equation}
\widehat \lambda_j=\lambda_0 \times ... \times \lambda_j; \hskip 30 pt
\lambda_j=\frac{{\rm Area\/}(\Sigma_{j-1} \cap \Sigma_{j})}{{\rm Area\/}(\Sigma_{j} \cap \Sigma_{j+1})};
\hskip 20 pt j=1,...,7.
\end{equation}

Let $p=(x,\pm 1)$ and $q=(y,\pm 1)$.
We have
\begin{equation}
\mu_+(q)-\mu_+(p)=(t,t,t) \hskip 5 pt {\rm mod\/} \hskip 5 pt \Lambda;
\hskip 30 pt t=\frac{y-x}{2}.
\end{equation}

\begin{lemma}
\label{nearinteger}
If ${\rm dist\/}(\mu_+(x),\mu_+(y))<\delta$ in $T^3$,
then there is an integer $I_k$ such that
$t \widehat \lambda_k$ is within $\epsilon$ of $I_k$
for all $k$,
\end{lemma}

\startproof
We compute
$$
{\rm Area\/}(\Sigma_0 \cap \Sigma_1)=8; \hskip 30 pt
{\rm Area\/}(\Sigma_1 \cap \Sigma_2)=\frac{8+8A}{1-A}; \hskip 30 pt
$$
\begin{equation}
\label{areas} 
{\rm Area\/}(\Sigma_2 \cap \Sigma_3)=\frac{2(1+A)^2}{A}; \hskip 30 pt
{\rm Area\/}(\Sigma_3 \cap \Sigma_4)=\frac{8+8A}{1-A}.
\end{equation}
This leads to 
\begin{equation}
\label{quantities}
\widehat \lambda_0=\widehat \lambda_4=1; \hskip 10 pt
\widehat \lambda_1=\widehat \lambda_3=\widehat \lambda_5=\widehat \lambda_7=
\frac{1-A}{1+A}; \hskip 10 pt
\widehat \lambda_2=\widehat \lambda_6=\frac{4A}{(1+A)^2}.
\end{equation}
The matrix
\begin{equation}
\label{conjugate}
H=\left[\matrix{
\frac{1}{1+A} & \frac{A-1}{(1+A)^2} & \frac{2A}{(1+A)^2} \cr
0 & \frac{1}{1+A} & \frac{1}{1+A} \cr
0&0&1}\right]
\end{equation}
conjugates the columns of the matrix defining $\Lambda$ to
the standard basis.  Therefore, if $\mu_+(x)$ and $\mu_+(y)$ are
close in $T^3$ then $H(t,t,t)$ is close to a point
of $\Z^3$.  
We compute
\begin{equation}
\label{preperiod}
H(t,t,t)=\bigg(\frac{4A}{(1+A)^2},\frac{2}{1+A},1\bigg) t =
(\widehat \lambda_2,\widehat \lambda_1-1,1) t.
\end{equation}
Equations
\ref{quantities} and \ref{preperiod} now finish the proof.
\endproof

\subsection{Pairs of Strips}
\label{strippairs}

Suppose $(S_1,S_2,V_2)$ is triple, where
$V_2$ is a vector pointing from one corner
of $S_1 \cap S_2$ to an opposite corner.  Let
$p_1 \in S_1$ and $p_2=E_2(p_1) \in S_2$.
Here $E_2$ is the strip map
associated to $(S_2,V_2)$.
We define $n$ and $\alpha$ by the equations
\begin{equation}
p_2-p_1=nV_2; \hskip 30 pt
\alpha=\frac{{\rm area\/}(B)}{{\rm area\/}(S_1 \cap S_2)};
\hskip 30 pt
\sigma_j=\frac{\|p_j-p_j'\|}{\|V_2\|}
\end{equation}
All quantities are affine invariant functions of the quintuple
$(S_1,S_2,V_2,p_1,p_2)$.  

\begin{center}
\psfig{file=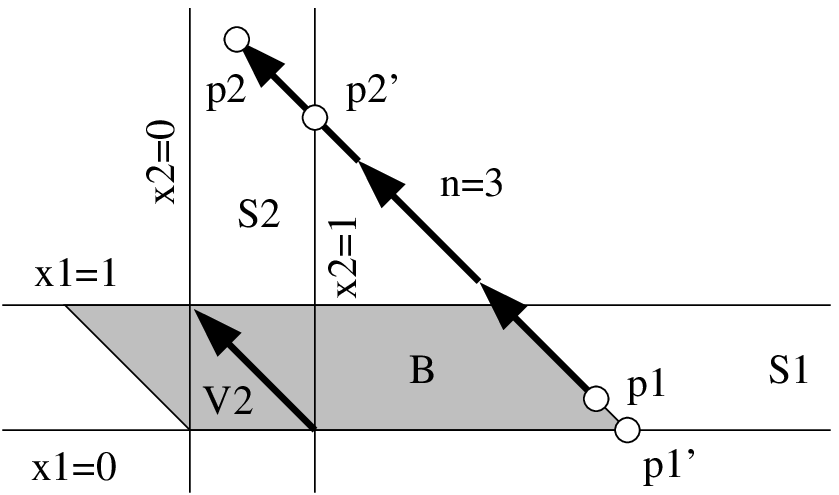}
\newline
{\bf Figure 8.1:\/} Strips and associated objects
\end{center} 

Figure 8.1 shows what we call the {\it standard pair\/} of
strips, where $\Sigma_j$ is the strip bounded by the lines
$x_j=0$ and $x_j=1$.
To get a better 
picture of the quantities we have defined, we consider
them on the standard pair.  We have
a\begin{equation}
\label{preaff}
\alpha=p_{11}+p_{12}=p_{21}+p_{22}; \hskip 15 pt
\sigma_1=p_{12}; \hskip 15 pt
\sigma_2=1-p_{22}; \hskip 15 pt
n={\rm floor\/}(p_{11}).
\end{equation}
Here $p_{ij}$ is the $j$th coordinate of $p_i$.
These equations lead to the following
affine invariant relations.
\begin{equation}
\label{aff}
n={\rm floor\/}(\alpha-\sigma_1); \hskip 30 pt
\sigma_2=1-[\alpha_1-\sigma_1]
\end{equation}
Here $[x]$ denotes the fractional part of $x$.
Again, the relations in Equation \ref{aff} hold
for any pair of strips. 

In our next result, we hold $(S_1,S_2,V_2)$ fixed but
compare all the quantities for $(p_1,p_2)$ and
another pair $(q_1,q_2)$.
Let $n(p)=n(S_1,S_2,V_2,p_1,p_2)$, etc.
Also, $N$ stands for an integer.   

\begin{lemma}
\label{strip1}
Let $\epsilon>0$.  There is some $\delta>0$ with the following
property.  If $|\sigma(p_1)-\sigma(q_1)|<\delta$ and
$|\alpha(q)-\alpha(p)-N|<\delta$ then
$|\sigma(p_2)-\sigma(q_2)|<\epsilon$ and $N=n(q)-n(p)$.
The number $\delta$ only depends on $\epsilon$ and
the distance from $\sigma(p_1)$ and $\sigma(p_2)$ to
$\{0,1\}$.
\end{lemma}

\startproof
If $\delta$ is small enough then
$[\alpha(p)-\sigma(p_1)]$ and
$[\alpha(q)-\sigma(q_1)]$ are very close,
and relatively far from $0$ or $1$.
Equation \ref{aff} now says that
$\sigma(p_2)$ and $\sigma(q_2)$ are close.
Also, the following two quantities are both
near $N$ while the individual 
summands are all relatively far from integers.
$$\alpha(q)-\alpha(p); \hskip 30 pt
\Big(\alpha(q)-\sigma(q_1)\Big)-
\Big(\alpha(p)-\sigma(p_1)\Big)$$
But the second quantity is near the integer
$n(q)-n(p)$, by Equation \ref{aff}.
\endproof

Suppose now that $S_1,S_2,S_3$ is a triple of strips, and $V_2,V_3$ is a
pair of vectors, such that $(S_1,S_2,V_2)$ and $(S_2,S_3,V_3)$ are as above.
Let $p_j \in S_j$ for $j=1,2,3$ be such that
$p_2=E_2(p_1)$ and $p_3=E_3(p_2)$.  Define,
\begin{equation}
\label{lambda}
\alpha_j=\alpha(S_j,S_{j+1},V_{j+1},p_j,p_{j+1}); \hskip 15 pt j=1,2; \hskip 30 pt
\lambda=\frac{{\rm Area\/}(S_1 \cap S_2)}{{\rm Area\/}(S_2 \cap S_3)}.
\end{equation}
It is convenient to set $\sigma_2=\sigma(p_2)$.

\begin{lemma}
\label{affinearea}
There are constants $C$ and $D$ such that
$\alpha_2=\lambda\alpha_1+C \sigma_2+D.$
The constants $C$ and $D$ depend on the strips.
\end{lemma}

\startproof
We normalize, as above, so that
Equation \ref{preaff} holds.  Then
\begin{equation}
\label{aa0}
p_2=(1-\sigma_2,\alpha_1+\sigma_2-1).
\end{equation}
There is a unique orientation preserving affine
transformation $T$ such that $T(S_{j+1})=S_j$ for
$j=1,2$, and $T$ the line $y=1$ to the line $x=0$.
Given that $S_1 \cap S_2$ has unit area, we have
$\det(T)=\lambda$.   Given the
description of $T$, we have
\begin{equation}
\label{aa2}
T(x,y)=\left(\matrix{a & \lambda \cr -1 & 0}\right)(x,y)+(b,1)=
(ax+b+\lambda y,1-x).
\end{equation}
Here $a$ and $b$ are constants depending on $S_2 \cap S_3$.
Setting $q=T(p_2)$, Equation \ref{preaff} gives $\alpha=q_1+q_2$.  Hence
\begin{equation}
\label{aa1}
\alpha_2=a(1-\sigma_2)+b+\lambda(\alpha_1+\sigma_2-1)+\sigma_2=\lambda\alpha_1+C\sigma_2+D.
\end{equation}
This completes the proof.
\endproof

\subsection{Single Parameter Proof}

We are still working under the assumption, in the Torus Lemma,
that $A=A^*$.
Our main argument relies on the Equation \ref{spectrum2},
which gives a formula for the return pairs in terms of 
the strip maps.
We define the point $q_j$ relative to $q$ just
as we defined $p_j$ relative to $p$.

\begin{center}
\resizebox{!}{3.7in}{\includegraphics{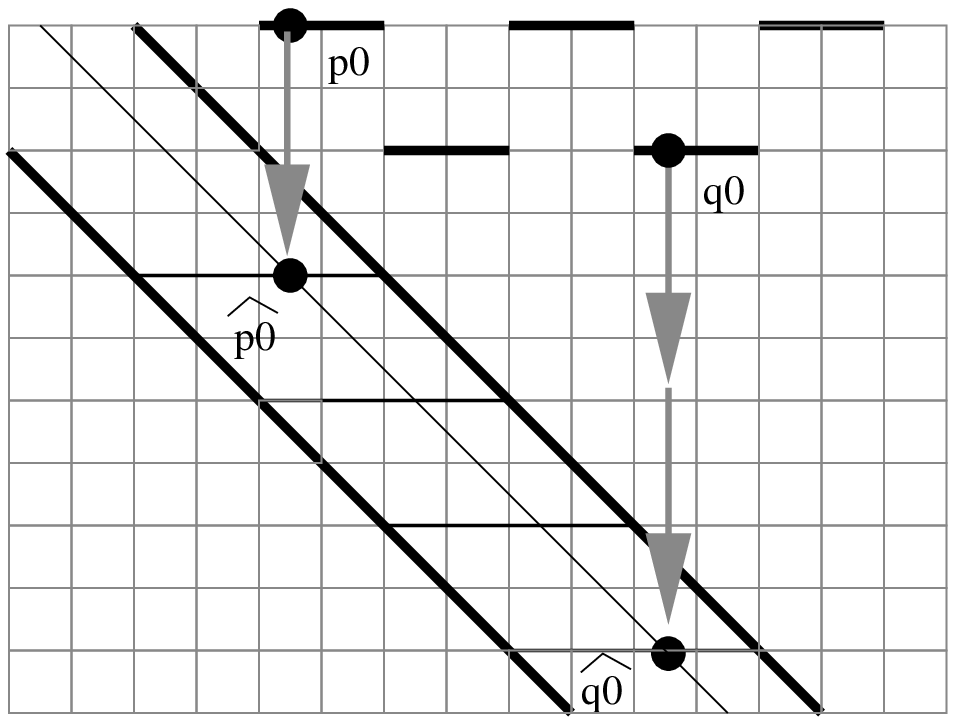}}
\newline
{\bf Figure 8.2:\/} The points $\widehat p_0$ and $\widehat q_0$.
\end{center} 

We would like to apply Lemmas \ref{nearinteger}, \ref{strip1}, and \ref{affinearea}
inductively. One inconvenience is that $p_0$ and $q_0$ do not
lie in any of our strips.  To remedy this situation we start with
the two points
\begin{equation}
\widehat p_0=E_0(p_0); \hskip 30 pt
\widehat q_0=E_0(q_0).
\end{equation}
We have $\widehat p_0, \widehat q_0 \in \Sigma_0$. 
Let $t$ be the near-integer from Lemma \ref{nearinteger}.
Looking at Figure 8.4, we see that
$|\sigma(\widehat q_0)-\sigma(\widehat p_0)|$ tends to $0$
as $\eta$ tends to $0$.

We define
\begin{equation}
\alpha_k(p)=\alpha(\Sigma_k,\Sigma_{k+1},V_{k+1},p_k,p_{k+1})
\end{equation}
It is also convenient to write
\begin{equation}
\sigma_k(p)=\sigma(p_k); \hskip 30 pt
\Delta \sigma_k=\sigma_k(q)-\sigma_k(p).
\end{equation}
For $k=0$, we use $\widehat p_0$ in place of $p_0$ and
$\widehat q_0$ in place of $q_0$ for these formulas. 

\begin{lemma}
\label{preaffine}
As $\eta \to 0$, the pairwise differences between
the $3$ quantities
$$\alpha_k(q)-\alpha_k(p); \hskip 30 pt
n_k(q)-n_k(p); \hskip 30 pt
t\widehat \lambda_k$$
converge to $0$ for all $k$.
\end{lemma}

\startproof
Referring to Figure 8.2, we have
$${\rm Area\/}(\Sigma_0 \cap \Sigma_1)=8; \hskip 20 pt
{\rm Area\/}(B(\widehat p_0))-
{\rm Area\/}(B(\widehat q_0))=4y-4x.
$$  
This gives us $\alpha_0(q)-\alpha_0(p)=t$.
Applying Lemma \ref{affinearea} inductively, we find that
\begin{equation}
\label{preAREA}
\alpha_k=\alpha_0 \widehat \lambda_k+
\sum_{i=1}^k \xi_i \sigma_i + C_k.
\end{equation}
for constants $\xi_1,...,\xi_k$ and $C_k$ that depend
analytically on $A$. 
Therefore
\begin{equation}
\label{AREA}
\alpha_k(q)-\alpha_k(p)=t \widehat \lambda_k+
\sum_{i=1}^k \xi_i \Delta \sigma_i;
\hskip 30 pt k=1,...,7
\end{equation}
By Lemma \ref{nearinteger}, the term $t\lambda_k$ is
near an integer
for all $k$.  By Lemma \ref{strip1} and induction, 
the remaining terms on the right hand side are near $0$.  
This lemma now follows from Lemma \ref{strip1}.
\endproof

Combining our last result with
Equation \ref{quantities}, we see that
$$
n_1(q)-n_1(p)=n_3(q)-n_3(p)=n_5(q)-n_5(p)=n_7(q)-n_7(p); 
$$
\begin{equation}
\label{cospectral}
n_2(q)-n_2(p)=n_6(q)-n_6(p).
\end{equation} 
once $\eta$ is small enough.
Given the dependence of constants in
Lemma \ref{strip1}, 
the necessary bound on $\eta$ only depends on
$\min(A,1-A)$ and $\theta(p)$.
Equation \ref{spectrum2} now tells us
that $\epsilon_j(p)=\epsilon_j(q)$ for $j=1,2$ once
$\eta$ is small enough.

\subsection{A Generalization of Lemma \ref{strip1}}

Now we turn to the proof of the Torus Lemma in the general case.
Our first result is the key step that allows
us to handle pairs of distinct parameters.  Once we set
up the notation, the proof is almost trivial.  Our
second result is a variant that will be useful in the
next chapter.

Suppose that $(S_1,S_2,V_2,p_1,p_2)$ and $(S_1^*,S_2^*,V_2^*,q_1^*,q_2^*)$
are two quintuples.  To fix the picture in our minds
we imagine that $(S_1,S_2,V_2)$ is near $(S_1^*,S_2^*,V_2^*)$, though
this is not necessary for the proof of the result to
follow.  We can define the quantities
$\alpha, \rho_j, n$ for each of these quintuples.  We
put a $(*)$ by each quantity associated to the second
triple.

\begin{lemma}
\label{strip22}
Let $\epsilon>0$.  There is some $\delta>0$ with the following
property.  If $|\sigma(p_1)-\sigma(q_1^*)|<\delta$ and
$|\alpha(q^*)-\alpha(p)-N|<\delta$ then
$|\sigma(p_2)-\sigma(q_2^*)|<\epsilon$ and $N=n(q^*)-n(p)$.
The number $\delta$ only depends on $\epsilon$ and
the distance from $\sigma(p_1)$ and $\sigma(p_2)$ to
$\{0,1\}$.
\end{lemma}

\startproof
There is an affine transformation
such that $T(X^*)=X$ for each object
$X=S_1,S_2,V_2$.  We set
$q_j=T(q_j^*)$.  Then 
$\alpha(q_1^*)=\alpha(q_1)$, by affine invariance.
Likewise for the other quantities.   Now we
apply Lemma \ref{strip1} to the
triple $(S_1,S_2,V_2)$ and the pairs
$(p_1,p_2)$ and $(q_1,q_2)$.  The conclusion
involves quantities with no $(*)$, but
returning the $(*)$ does not change any
of the quantities.
\endproof

For use in the next chapter, we state a variant of Lemma \ref{strip22}.
Let $[x]$ denote the image
of $x \in \R/\Z$.

\begin{lemma}
\label{strip222}
Let $\epsilon>0$.  There is some $\delta>0$ with the following
property.  If $|\sigma(p_1)-\sigma(q_1^*)|<\delta$ and
$|\alpha(q^*)-\alpha(p)-N|<\delta$ then the distance
from $[\sigma(p_2)]$ and $[\sigma(q_2)^*]$ in
$\R/\Z$ is less than $\epsilon$. 
$|\sigma(p_2)-\sigma(q_2^*)|<\epsilon$ and $N=n(q^*)-n(p)$.
The number $\delta$ only depends on $\epsilon$ and
the distance from $\sigma(p_1)$ to
$\{0,1\}$.
\end{lemma}

\startproof
Using the same trick as in Lemma \ref{strip1}, we reduce to the
single variable case.  In this case, we mainly repeat
the proof of Lemma \ref{strip1}.
If $\delta$ is small enough then
$[\alpha(p)-\sigma(p_1)]$ and
$[\alpha(q)-\sigma(q_1)]$ are very close,
and relatively far from $0$ or $1$.
Equation \ref{aff} now says that
$[\sigma(p_2)]$ and $[\sigma(q_2)]$ are close
in $\R/\Z$.   
\endproof

\subsection{Proof in the General Case}

We no longer suppose that $A=A^*$, and we return to the
original notation $(q^*,A^*)$ for the second point.
In our proof of this result, we 
attach a $(*)$ to any quantity that depends on
$(q^*,A^*)$.   
We first need to repeat the analysis from \S \ref{mapinput},
this time keeping track of the parameter.
Let $\eta$ be as in the Torus Lemma. 
We use the big O notation.

\begin{lemma}
\label{nearinteger2}
There is an integer $I_k$ such that
$|\alpha_0^* \widehat \lambda_k^*-
\alpha_0 \lambda_k-I_k| <O(\eta).$
\end{lemma}

\startproof
Let $[V]$ denote the distance from $V \in \R^3$ to the
nearest point in $\Z^3$.
Let $p=(x,\pm 1)$ and $q^*=(x^*,\pm 1)$.
Recalling the definition of
$\mu_+$, the hypotheses in the Torus Lemma imply that
\begin{equation}
\label{newaff}
\bigg[H^*\Big(\frac{x^*}{2},\frac{x^*}{2}+1,\frac{x^*}{2}\Big)-
H\Big(\frac{x}{2},\frac{x}{2}+1,\frac{x}{2})\bigg]<O(\eta)
\end{equation}
We compute that
$\alpha_0=x/2+1/2$, independent of parameter.
Therefore
$$H\Big(\frac{x}{2},\frac{x}{2}+1,\frac{x}{2}\Big)=
H(\alpha_0,\alpha_0,\alpha_0)+\frac{1}{2}H(-1,1,-1).$$
The same goes with the starred quantities.
Therefore,
$$
[(\widehat \lambda_2^*,\widehat \lambda_1^*-1,1) \alpha_0^*-
(\widehat \lambda_2,  \widehat \lambda_1-1,  1) \alpha_0]=$$
$$
[H^*(\alpha_0^*,\alpha_0^*,\alpha_0^*)-H(\alpha_0,\alpha_0,\alpha_0)]<
O(\eta)+\|(H^*-H)(-1,1,-1)\|<O(\eta).
$$
Our lemma now follows immediately from
Equation \ref{quantities}.
\endproof

The integer $I_k$ of course depends on 
$(p,A)$ and $(q^*,A^*)$, but in all cases
Equation \ref{quantities} gives us
\begin{equation}
\label{quantities2}
I_0=I_4; \hskip 30 pt I_1=I_3=I_5=I_7; \hskip 30 pt I_2=I_6,
\end{equation}

\begin{lemma}
As $\eta \to 0$, the pairwise differences between the $3$ quantities
$\alpha^*_k-\alpha_k$ and
$n_k^*-n_k$ and $I_k$
tends to $0$ for all $k$.
\end{lemma}

\startproof
Here $\alpha_k^*$ stands for $\alpha_k(q^*)$, etc.
Equation \ref{preAREA} works separately for each parameter. 
The replacement for Equation \ref{AREA} is
\begin{equation}
\label{analytic}
\alpha^*_k-\alpha_k=W+X+Y; \hskip 30 pt W=\alpha_0^*\widehat \lambda_k^*-\alpha_0 \widehat \lambda_k
\end{equation}
\begin{equation}
X=\sum_{i=1}^k \xi_i^* \sigma_i^*(q^*) - \sum_{i=1}^k \xi_i \sigma_i(p)
=\sum_{i=1}^k \xi_i \big(\sigma_i^*-\sigma_i\big)+O(|A-A^*|);
\end{equation}
\begin{equation}
Y=\sum_{i=1}^k C_i^* - \sum_{i=1}^k C_i=O(|A-A^*|).
\end{equation}
The estimates on $X$ and $Y$ comes from the fact 
$\xi_i$ and $C_i$ vary smoothly with $A$.
Putting everything together, we get the following.
\begin{equation}
\label{ccc}
\alpha^*_k-\alpha_k=\Big(\alpha_0^*\widehat \lambda_k^*-
\alpha_0\lambda_k\Big)+
\sum_{i=1}^k \xi_i \big(\sigma_i^*-\sigma_i\big)+O(|A-A^*|).
\end{equation}
In light of Lemma \ref{nearinteger2}, it suffices to show that
$\sigma_i^*-\sigma_i$ tends to $0$ as $\eta$ tends to $0$.
The same argument as in the single parameter case works here,
with Lemma \ref{strip22} used in place of Lemma \ref{strip1}.
\endproof

Similar to the single parameter case,
Equations \ref{spectrum2} and \ref{quantities2} now finish the proof.

\newpage

\section{The Strip Functions}
\label{sing}

\subsection{The Main Result}
\label{singular0}

The purpose of this chapter is to understand the
functions $\sigma_j$ that arose in the proof of
the Master Picture Theorem. 
We call these functions the {\it strip functions\/}.

Let $W_k \subset \Xi_+ \times (0,1)$ denote the set of
points where $E_k...E_1$ is defined but
$E_{k+1}E_k...E_1$ is not defined.
Let $S_k$ denote the closure of
$\mu_+(W_k)$ in $R$.   Finally,
let 
\begin{equation}
W_k'=\bigcup_{j=0}^{k-1} W_j;
\hskip 30 pt 
S_k'=\bigcup_{j=0}^{k-1}S_j; \hskip 30 pt k=1,...,7.
\end{equation}
The Torus Lemma applies to any point that
does not lie in the {\it singular set\/}
\begin{equation}
\label{lad}
S=S_0 \cup ... \cup S_7.
\end{equation}

If $p \in \Xi_+-W_k'$ then the points
$p=p_0,...,p_k$ are defined.  Here,
as in the previous chapter,
$p_j=E_j(p_{j-1})$. 
The functions $\sigma_1,...,\sigma_k$
and $\alpha_1,...,\alpha_k$ are defined
for such a choice of $p$.
Again, $\sigma_j$ measures the position
of $p_j$ in $\Sigma_j$, relative to
$\partial \Sigma_j$.  Even if
$E_{k+1}$ is not defined on $p_k$, the
equivalence class $[p_{k+1}]$ 
is well defined in the cylinder
$\R^2/\langle V_{k+1}\rangle$.
The corresponding function
$\sigma_{k+1}(q)=\sigma(q_{k+1})$ is well defined as
an element of $\R/\Z$.  

Let $\pi_j: \R^4 \to \R$ be the $j$th coordinate
projection.  Let $[x]$ denote the image of $x$ in $\R/\Z$.
The following identities refer to the
$(+)$ case.  We 
discuss the $(-)$ case at the end of the chapter.

\begin{equation}
\label{identity1}
\sigma_1=\bigg[\frac{2-\pi_3}{2}\bigg] \circ \mu_+ \hskip 30 pt {\rm on\/} \hskip 5 pt \Xi_+
\end{equation}

\begin{equation}
\label{identity2}
\sigma_2=\bigg[\frac{1+A-\pi_2}{1+A}\bigg] \circ \mu_+ \hskip 30 pt {\rm on\/}  \hskip 5 pt \Xi_+-W_1'
\end{equation}

\begin{equation}
\label{identity3}
\sigma_3=\bigg[\frac{1+A-\pi_1}{1+A}\bigg] \circ \mu_+ \hskip 30 pt  {\rm on\/} \hskip 5 pt \Xi_+-W_2'
\end{equation}

\begin{equation}
\label{identity4}
\sigma_4=\bigg[\frac{1+A-\pi_1-\pi_2+\pi_3}{2}\bigg] \circ \mu_+ \hskip 30 pt {\rm on\/}  \hskip 5 pt \Xi_+-W_3'
\end{equation}

In the next chapter we deduce the Master Picture Theorem from these identities
and the Torus Lemma.  In this chapter, we prove the identities.

\subsection{Continuous Extension}

Let $g=\sigma_j$ for $j=0,...,k+1$.
since the image
$\mu_+(\Xi \times (0,1))$ is dense in $R-S_k'$, we define
\begin{equation}
\widetilde g(\tau):=\lim_{n \to \infty} g(p_n,A_n);
\hskip 30 pt \tau \in R-S_k'.
\end{equation}
Here $(p_n,A_n)$ is chosen so that
all functions are defined and
$\mu_+(p_n,A_n) \to \tau$.  Note that the sequence
$\{p_n\}$ need not converge.

\begin{lemma}
The functions $\widetilde \sigma_1,...,\widetilde \sigma_{k+1}$,
considered as $\R/\Z$-valued functions, are
well defined and continuous on $R-S_k'$.
\end{lemma}

\startproof
For the sake of concreteness, we will give the proof in the
case $k=2$.  This representative case explains the idea.
First of all, the continuity follows from the well-definedness.
We just have to show that the limit above is always well defined.
$\widetilde \sigma_1$ is well defined and continuous
on all of $R$, by
Equation \ref{identity1}.

Since $S_1' \subset S_2'$, we see that 
$\tau \in R-S_1'$.
Hence $\tau$ does not lie in
the closure of $\mu_+(W_0)$.
Hence, there is some $\theta_1>0$ such 
that $\theta_1(p_n,A_n)>\theta_1$ for
all sufficiently large $n$.
Note also that there is a positive
and uniform lower bound to $\min(A_n,1-A_n)$.
Note that $[\alpha_1(p_n,A_n)]=[\pi_3(\mu_+(p_n,A_n)]$.
Hence
$\{[\alpha_1(p_n,A_n)]\}$ is a Cauchy sequence
in $\R/\Z$.

Lemma \ref{strip222} now applies uniformly to
$$(p,A)=(p_m,A_m); \hskip 30 pt
(q^*,A^*)=(p_n,A_n)$$ for all
sufficiently large pairs $(m,n)$.  Since
$\{\mu_+(p_n,A_n)\}$ forms a Cauchy sequence in
$R$,  Lemma \ref{strip222} implies that 
$\{\sigma_2(\tau_m,A_m)\}$ forms a Cauchy
sequence in $\R/\Z$.   Hence,
$\widetilde \sigma_2$ is well defined on
$R-S_1'$, and continuous.

Since $\tau \in R-S_2'$, we see that
$\tau$  does not lie in
the closure of $\mu_+(W_1)$.
Hence, there is some $\theta_2>0$ such 
that $\theta_j(p_n,A_n)>\theta_j$ for
$j=1,2$ and all sufficiently large $n$.
As in our proof of the General Torus Lemma,
Equation \ref{ccc} now says that
shows that $\{\alpha_2(p_n,A_n)\}$ forms a
Cauchy sequence in $\R/\Z$.
We now repeat the previous argument 
to see that
$\{\sigma_3(\tau_m,A_m)\}$ forms a Cauchy
sequence in $\R/\Z$.   Hence,
$\widetilde \sigma_3$ is well defined on
$R-S_2'$, and continuous. 
\endproof

Implicit in our proof above is the function
\begin{equation}
\beta_k=[\alpha_k] \in \R/\Z.
\end{equation}
This function will come in handy in our next result.

\subsection{Quality of the Extension}

Let $X=R-\partial R \subset \R^4$. Note that $X$ is open and convex

\begin{lemma}
\label{analy}
Suppose $X \subset R-S_k'$.  Then
$\widetilde \sigma_{k+1}$ is locally affine on $X_A$.
\end{lemma}

\startproof  
Since $\widetilde \sigma_{k+1}$ is continuous on $X$,
it suffices to prove this lemma for a dense set of $A$.
We can choose $A$ so that $\mu_+(\Xi_+)$ is dense in $X_A$.

We already know that
$\widetilde \sigma_1,...,\widetilde \sigma_{k+1}$
are all defined and continuous on $X$.
We already remarked that Equation \ref{identity1}
is true by direct inspection.   As we already remarked
in the previous proof, $\beta_0=\pi_3 \circ \mu_+$.
Thus, we define $\widetilde \beta_0=[\pi_3]$.
Let $\widetilde \beta_0=[\pi_3]$.
Both $\widetilde \sigma_0$ and $\widetilde \beta_0$ are locally
affine on $X_A$.

Let $m \leq k$.
The second half of Equation \ref{aff} tells us that
$\widetilde \sigma_m$ is a locally affine function of
$\widetilde \sigma_{m-1}$ and $\widetilde \beta_{m-1}$.
Below we will prove that
$\widetilde \beta_m$ is
defined on $X_A$, and locally affine, provided that
$\widetilde \sigma_1,...,\widetilde \sigma_m$ are locally affine.
Our lemma follows from this claim and induction.

Now we prove the claim.
All the addition below is done in
$\R/\Z$.
Since $\mu_+(\Xi_+)$ is dense in $X_A$, we can
at least define $\widetilde \beta_m$ on a dense
subset of $X_A$.   Define
\begin{equation}
p=(x,\pm 1); \hskip 10 pt
p'=(x',\pm 1); \hskip 10 pt
\tau=\mu_+(p); \hskip 10 pt \tau'=\mu_+(p');
\hskip 10 pt t=\frac{x'-x}{2}.
\end{equation}
We choose $p$ and $p'$ so that the pinwheel map is entirely defined.

From Equation \ref{AREA}, we have
\begin{equation}
\label{AREA2}
\widetilde \beta_m(\tau')-\widetilde \beta_m(\tau)=
[t\widehat \lambda_k] +\sum_{j=1}^m [\xi_j \times (\widetilde \sigma_j(\tau')-\widetilde \sigma_j(\tau))].
\end{equation}
Here $\xi_1,...,\xi_m$ are constants that depend on $A$.
Let $H$ be the matrix in Equation \ref{conjugate}.
We have
$H(t,t,t) \equiv H(\tau'-\tau)$ mod $\Z^3$
because $(t,t,t) \equiv \tau'-\tau$ mod $\Lambda$.
Our analysis in \S \ref{mapinput} shows that
\begin{equation}
[t\widehat \lambda_k]=[\pi \circ H(t,t,t)-\epsilon t]=
[(\pi-\epsilon \pi_3) \circ H(\tau'-\tau)].
\end{equation}
Here $\epsilon \in \{0,1\}$ and
$\pi$ is some coordinate projection.
The choice of $\epsilon$ and $\pi$ depends on $k$.
We now see that
\begin{equation}
\widetilde \beta_m(\tau')=
\widetilde \beta_m(\tau)+(\pi+ \epsilon_3 \pi) \circ H(\tau'-\tau)+\sum_{j=1}^m 
[\xi_j \times (\widetilde \sigma_j(\tau')-\widetilde \sigma_j(\tau))].
\end{equation}
The right hand side is everywhere defined and locally affine.
Hence, we define $\widetilde \beta_m$ on all of $X_A$ using
the right hand side of the last equation.  
\endproof

\begin{lemma}
\label{nicean}
Suppose $X \subset R-S_k'$.  Then
$\sigma_{k+1}$ is analytic on $X$.
\end{lemma}

\startproof
The constants $\xi_j$ in Equation \ref{AREA2}
vary analytically with $A$.
Our argument in Lemma \ref{analy}
therefore shows that the linear part of
$\sigma_{k+1}$ varies analytically with $A$.
We just have to check the linear term.
Since $X_A$ is connected
we can compute the linear term of $\sigma_{k+1}$ at $A$
from a single point.  We choose
$p=(\epsilon,1)$ where $\epsilon$ is very close to $0$.
The fact that
$A \to \sigma_k(p,A)$ varies analytically follows
from the fact that our strips
vary analytically. 
\endproof

\noindent
{\bf Remark:\/}
We have
$S_k \subset \widetilde \sigma_{k+1}^{-1}([0])$.
Given Equation \ref{identity1},
we see that $X \subset R-S_1'$.  Hence
$\sigma_2$ is defined on $X$.  Hence
$\sigma_2$ is analyic on $X$ and locally affine
on each $X_A$.  We use these two properties
to show that Equation \ref{identity2} is true.
But then $X \subset R-S_2'$. etc.
So, we will know at each stage of our
verification that Lemmas \ref{analy} and \ref{nicean}
apply to the function of interest.
\newline

Equations \ref{identity2}, 
\ref{identity3}, and \ref{identity4} are formulas
for $\widetilde \sigma_2$, $\widetilde \sigma_3$, and
$\widetilde \sigma_4$ respectively.
Let $f_{k+1}=\widetilde \sigma_{k+1}-\sigma_{k+1}'$, where
$k=2,3,4$.  Here $\sigma_{k+1}'$ is the
right hand side of the identity for $\widetilde \sigma_{k+1}$.
Our goal is to show that $f_{k+1}\equiv [0]$ for
$k=1,2,3$.     Call a parameter $A$ {\it good\/} if
$f_{k+1} \equiv [0]$ on $X_A$.
Call a subset $S \subset (0,1)$ {\it substantial\/} if
$S$ is dense in some open interval of $(0,1)$.  By 
analyticity, $f_{k+1} \equiv [0]$ provided that
a substantial set of parameters is good.

In the next section we explain how to verify that
a parameter is good.   If $f_{k+1}$ was
a locally affine map from $X_A$ into $\R$, we would
just need to check that $f_{k+1}=0$ on some tetrahedron
on $X_A$ to verify that $A$ is a good parameter.
Since the range of $f_{k+1}$ is $\R/\Z$, we have
to work a bit harder.

\subsection{Irrational Quintuples}

We will give a construction in $\R^3$.  When the time
comes to use the construction, we will identify
$X_A$ as an open subset of a copy of $\R^3$.

Let $\zeta_1,...,\zeta_5 \in \R^3$ be $5$ points.
By taking these points $4$ at a time, we
can compute $5$ volumes, $v_1,...,v_5$.  Here
$v_j$ is the volume of the tetrahedron obtained
by omitting the $j$th point.  We
say that $(\zeta_1,...,\zeta_5)$ is an
{\it irrational quintuple\/} if the
there is no rational relation
\begin{equation}
\sum_{i=1}^5 c_j \zeta_j=0; \hskip 30 pt
c_j \in \Q; \hskip 30 pt c_1c_2c_3c_4c_5=0.
\end{equation}
If we allow all the constants to be nonzero,
then there is always a relation.

\begin{lemma}
Let $C$ be an open convex subset of $\R^3$.
Let $f: C \to \R/\Z$ be a locally affine
function.  Suppose that there is an irrational
$(\zeta_1,...,\zeta_5)$ such that
$\zeta_j \in C$ and $f(\zeta_j)$ is the
same for all $j$.
Then $f$ is constant on $C$.
\end{lemma}

\startproof
Since $C$ is simply connected, we can lift
$f$ to a locally affine 
function $F: C \to \R$.   But then $F$ is
affinc on $C$, and we can extend $F$ to be an
affine map from $\R^3$ to $\R$.  By construction
$F(\zeta_i)-F(\zeta_j) \in \Z$ for all $i,j$.
Adding a constant to $F$, we can assume that $F$
is linear.
There are several cases.
\newline
\newline
{\bf Case 1:\/}
Suppose that $F(\zeta_j)$ is independent of $j$.
In this case, all the points lie in the same plane,
and all volumes are zero.  This violates
the irrationality condition.
\newline
\newline
{\bf Case 2:\/}
Suppose we are not in Case 1, and the following
is true.  For every index $j$ there is a second
index $k$ such that
$F(\zeta_k) = F(\zeta_j)$.
Since there are $5$ points total,
this means that the set
$\{F(\zeta_j)\}$ only has a total of $2$ values.
But this means that our $5$ points lie in a
pair of parallel planes, $\Pi_1 \cup \Pi_2$, with
$2$ points in $\Pi_1$ and $3$ points in $\Pi_2$.
Let's say that
that $\zeta_1,\zeta_2,\zeta_3 \in \Pi_1$ and
$\zeta_4,\zeta_5 \in \Pi_2$.  
But then $v_4=v_5$, and we violate the irrationality
condition.
\newline
\newline
{\bf Case 3:\/} If we are not in the above two cases,
then we can relabel so that
$F(\zeta_1) \not = F(\zeta_j)$ for $j=2,3,4,5$.
Let $$\zeta_j'=\zeta_j-\zeta_1.$$
Then $\zeta_1'=(0,0,0)$ and $F(\zeta_1')=0$.  But then
$F(\zeta_j') \in \Z-\{0\}$ for $j=2,3,4,5$.
Note that $v_j'=v_j$ for all $j$.
For $j=2,3,4,5$, let
$$\zeta_j''=\frac{\zeta_j'}{F(\zeta_j')}.$$
Then $v_j''/v_j' \in \Q$ for $j=2,3,4,5$.
Note that $F(\zeta_j'')=1$ for $j=2,3,4,5$.
Hence there is a plane $\Pi$ such that
$\zeta_j'' \in \Pi$ for $j=2,3,4,5$.

There is always a rational relation between
the areas of the $4$ triangles defined by $4$ points in the plane.
Hence, there is a rational relation between
$v_2'',v_3'',v_4'',v_5''$.  But then there
is a rational relation between
$v_2,v_3,v_4,v_5$.  This contradicts the
irrationality condition.
\endproof

\subsection{Verification in the Plus Case}

Proceeding somewhat at random, we define
\begin{equation}
\phi_j=\bigg(8jA+\frac{1}{2j},1\bigg); \hskip 30 pt j=1,2,3,4,5.
\end{equation}
We check that $\phi_j \in \Xi_+$ for $A$ near $1/2$.
Letting $\zeta_j=\mu_+(\phi_j)$, we check
that $f_{k+1}(\zeta_j)=[0]$ for $j=1,2,3,4,5$.
\newline
\newline
{\bf Example Calculation:\/}
Here is an example of what we do automatically in
Mathematica.
Consider the case $k=1$ and $j=1$.  When $A=1/2$, the
length spectrum for $\phi_1$ starts out $(1,1,2,1)$.
Hence, this remains true for nearby $A$.  Knowing the
length spectrum allows us to compute, for instance, that
$$E_2E_1(\phi_1)=\phi_1+V_1+V_2=\bigg(\frac{-3}{2}+8A,7\bigg) \in \Sigma_2$$
for $A$ near $1/2$.
The affine functional
\begin{equation}
(x,y) \to (x,y,1) \cdot \frac{(-1,A,A)}{2+2A}
\end{equation}
takes on the value $0$ on the line $x=Ay+A$ and $1$ on the line
$x=Ay-2-A$.    These are the two edges of $\Sigma_2$. (See \S \ref{pinwheelformulas}.)
 Therefore,
$$\sigma_2(\phi_1)=\bigg(\frac{-3}{2}+8A,7,1\bigg) \cdot \frac{(-1,A,A)}{2+2A}=
\frac{3}{4+4A}.$$
At the same time, we compute that
$$\mu_+(\phi_1)=\frac{1}{4}(-7+24A,1+4A,-7+16A),$$
at least for $A$ near $1/2$.  When $A$ is far from $1/2$
this point will not lie in $R_A$.   We then compute
$$\frac{1+A-\pi_2(\mu_+(\phi_1))}{1+A}=\frac{3}{4+4A}.$$
This shows that $f_2(\zeta_2)=[0]$ for all $A$ near $1/2$.
The verifications for the other pairs $(k,j)$ are similar.
\newline

\noindent
{\bf Checking Irrationality:\/}
It only remains to check that the points
$(\zeta_1,...,\zeta_5)$ form an irrational
quintuple for a dense set of parameters $A$.  In
fact this will true in the complement of a countable
set of parameters.

The $5$ volumes associated to our quintuple are as follows.
\begin{itemize}
\item $v_5=5/24 - 5A/12 + 5A^2/24$.
\item $v_4=71/40 + 19A/20 -787A^2/120 - 4A^3.$
\item $v_3=119/60 +7A/60-89A^2/15 - 4A^3$
\item $v_2=-451/240 - 13A/40 + 1349 A^2/240 + 4A^3$
\item $v_1=-167/80 - 13A/40 + 533 A^2/80+4A^3.$
\end{itemize}
If there is an open set of parameters for which
the first $4$ of these volumes has a rational
relation, then there is an infinite set on which
the same rational relation holds. Since every formula
in sight is algebraic, this means that there must be a single
rational relation that holds for all parameters.
But then the curve
$A \to (v_5,v_4,v_3,v_2)$ lies in a proper linear
subspace of $\R^4$.  

We evaluate this curve at
$A=1,2,3,4$ and see that the resulting points
are linearly independent in $\R^4$.  Hence,
there is no global rational relation. Hence,
on a dense set of parameters, there is no
rational relation between the first $4$
volumes listed.  A similar argument rules
out rational relations amongst any other
$4$-tuple of these volumes.

\subsection{The Minus Case}

In the $(-)$ case,
Equations \ref{identity2} and \ref{identity3} do not change,
except that $\mu_-$ replaces $\mu_+$ and all the sets are
defined relative to $\Xi_-$ and $\mu_-$.
Equations \ref{identity1} and \ref{identity4} become
\begin{equation}
\label{identity11}
\sigma_1=\bigg[\frac{1-\pi_3}{2}\bigg] \circ \mu_- \hskip 30 pt
{\rm on\/} \hskip 5 pt \Xi_-.
\end{equation}
\begin{equation}
\label{identity41}
\sigma_4=\bigg[\frac{A-\pi_1-\pi_2+\pi_3}{2}\bigg] \circ \mu_- \hskip 30 pt
{\rm on\/} \hskip 5 pt \Xi_--S_3'.
\end{equation}
Lemma \ref{analy} and Lemma \ref{nicean} have the same
proof in the $(-)$ case.  We use the same method
as above,  except that we use the points
\begin{equation}
\phi_j+(2,0); \hskip 50 pt j=1,2,3,4,5.
\end{equation}
These points all lie in $\Xi_-$ for $A$ near $1/2$.
The rest of the verification is essentially the same
as in the $(+)$ case.

\newpage
\section{Proof of the Master Picture Theorem}
\label{torus2}

\subsection{The Main Argument}

Let $S$ be the singular set defined in Equation \ref{lad}.
Let $\widetilde S$ denote the union of hyperplanes listed
in \S \ref{walls}.   let $d$ denote distance on
the polytope $R$. In this chapter we will prove

\begin{lemma}[Hyperplane]
$S \subset \widetilde S$ and
$\theta(p,A) \geq d(\mu_+(p,A),\widetilde S)$.
\end{lemma}

We finish the proof of the Master Picture Theorem assuming
the Hyperplane Lemma.

Say that a {\it ball of constancy\/} in $R-\widetilde S$ is an open ball $B$
with the following property.  If $(p_0,A_0)$ and $(p_1,A_1)$ are
two pairs and $\mu_+(p_j,A_k) \in B$ for $j=0,1$, then
$(p_0,A_0)$ and $(p_1,A_1)$ have the same return pair.
Here is a consequence of the Torus Lemma.

\begin{corollary}
Any point $\tau$ of $R-\widetilde S$ is contained in a ball of constancy.
\end{corollary}

\startproof
If $\tau$ is in the image of $\mu_+$, this result is an immediate
consequence of the Torus Lemma.  In general, the image 
$\mu_+(\Xi_+ \times (0,1))$ is dense in $R$.  Hence, we can
find a sequence $\{\tau_n\}$ such that $\tau_n \to \tau$ and
$\tau_n=\mu_+(p_n,A_n)$.   
Let $2\theta_0>0$ be the distance from $\tau$ to $S$.
From the triangle inequality and the
second statement of the Hyperplane Lemma,
$\theta(p_n,A_n) \geq \theta_0=\theta_1>0$ for large $n$.  By the Torus Lemma,
$\tau_n$ is the center of a ball $B_n$ of constancy whose
radius depends only on $\theta_0$.  In particular -- and
this is really all that matters in our proof --
the radius of $B_n$ does not tend to $0$.  Hence, for
$n$ large enough, $\tau$ itself is contained in $B_n$.
\endproof

\begin{lemma}
\label{continuee}
Let $(p_0,A_0)$ and $(p_1,A_1)$ be two
points of $\Xi_+ \times (0,1)$ such that
$\mu_+(p_0,A_0)$ and $\mu_+(p_1,A_1)$ lie in the
same path connected component of $R-\widetilde S$.
Then the return pair for $(p_0,A_0)$ equals
the return pair for $(p_1,A_1)$.
\end{lemma}

\startproof
Let $L \subset  R-\widetilde S$ be a path joining
points $\tau_0= \mu_+(p_0,A_0)$ and
$\tau_1= \mu_+(p_1,A_1)$. 
By compactness, we can cover $L$ by finitely many
overlapping balls of constancy.
\endproof

Now we just need to see that the Master Picture Theorem
holds for one component of the partition of $R-\widetilde S$.
Here is an
example calculation that does the job.
For each $\alpha=j/16$ for $j=1,...,15$, we
plot the image
\begin{equation}
\mu_A(2\alpha+2n); \hskip 30 pt n=1,...,2^{15};
\end{equation}
The image is contained in the slice $z=\alpha$.
We see that the Master Picture Theorem holds for
all these points.   The reader can use Billiard
King to plot and inspect millions of points for
any desired parameter.

We have really only proved the half of the Master Picture
Theorem that deals with $\Xi_+$ and $\mu_+$.  The half
that deals with $\Xi_-$ and $\mu_-$ is exactly the same.
In particular, both the Torus Lemma and the
Hyperplane Lemma hold {\it verbatim\/} in the $(-)$ case.
The proof of the Hyperplane Lemma in the $(-)$
case differs only in that the two identities
in Equation \ref{identity11} replace Equations
\ref{identity1} and \ref{identity4}.  We omit
the details in the $(-)$ case.

\subsection{The First Four Singular Sets}

Our strip function identites make short work of the
first four pieces of the singular set.
\begin{itemize}
\item
Given Equation \ref{identity1}, 
\begin{equation}
\label{consequence1}
S_0 \subset \{z=0\} \cup \{z=1\}.
\end{equation}

\item Given Equation \ref{identity2},
\begin{equation}
\label{consequence2}
S_1 \subset \{y=0\} \cup \{y=1+A\}.
\end{equation}

\item Given Equation \ref{identity3},
\begin{equation}
\label{consequence3}
S_2 \subset \{x=0\} \cup \{x=1+A\}.
\end{equation}

\item Give Equation \ref{identity4}, 
\begin{equation}
\label{consequence4}
S_3 \subset \{x+y-z=1+A\} \cup
\{x+y-z=-1+A\}.
\end{equation}

\end{itemize}
\subsection{Symmetry}
\label{symmetry00}

We use symmetry to deal with the remaining pieces.
Suppose we start with a point $p \in \Xi_+$.  We
define $p_0=p$ and $p_j=E_j(p)$.  As we go along
in our analysis, these points will be defined for
increasingly large values of $j$.  However, for
the purposes of illustration, we assume that
all points are defined.

Let $\rho$ denote reflection in the $x$-axis.
Then
\begin{equation}
\rho(\Sigma_9-j)=\Sigma_j; \hskip 30 pt
q_j=\rho(p_{9-j}); \hskip 30 pt j=1,2,3,4.
\end{equation}
Figure 10.1 shows a picture.  The disk in the center
is included for artistic purposes, to cover up some messy
intersections.
In the picture, we have included the coordinates for the
vectors $-V_1$ and $-V_2$ and $-V_2$ to remind the
reader of their values.  It is convenient to write
$-V_k$ rather than $V_k$ because there are
far fewer minus signs involved. 

\begin{center}
\psfig{file=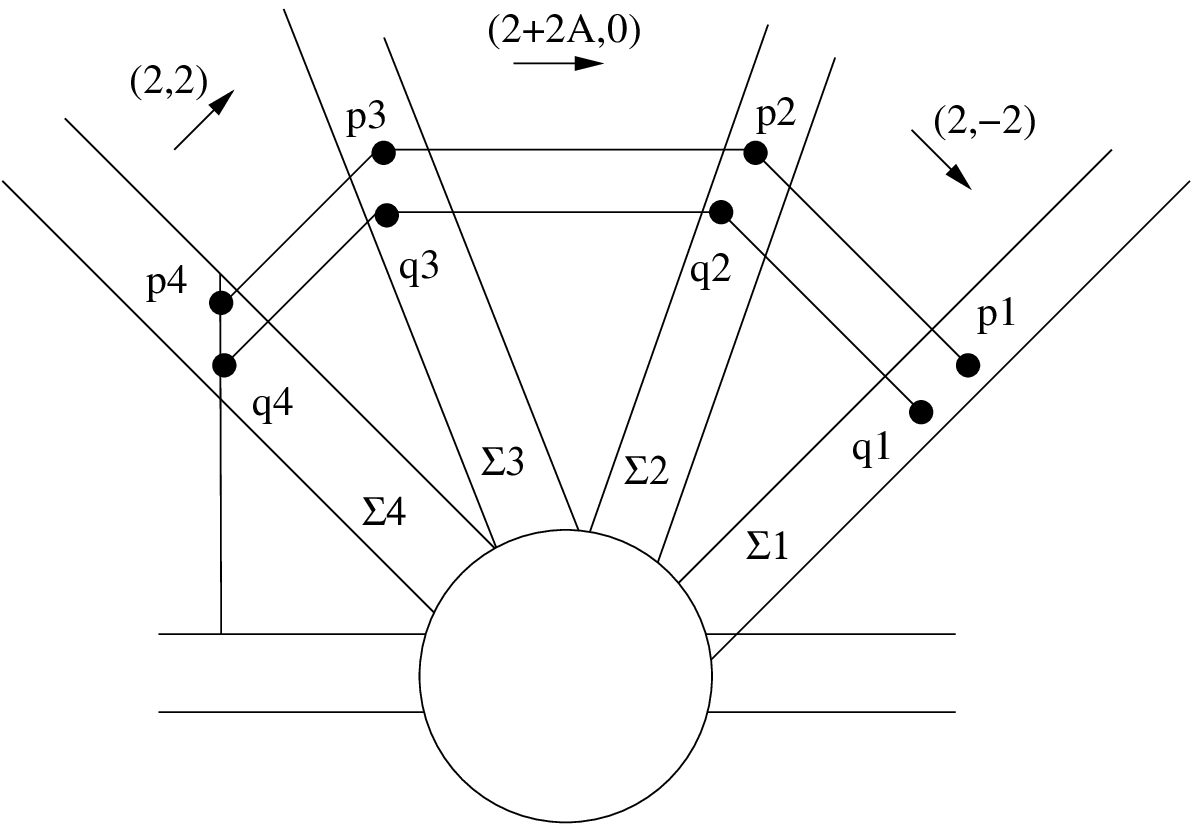}
\newline
{\bf Figure 10.1:\/} Reflected points
\end{center} 

Here is a notion we will use in our estimates.
Say that a strip $\Sigma$  {\it dominates\/} a vector $V$ if we can
translate $V$ so that it is contained in the interior of the
strip.  This is equivalent to the condition that we can translate
$V$ so that one endpoint of $V$ lies on $\partial \Sigma$ and
the other one lies in the interior.

\subsection{The Remaining Pieces}

\subsubsection{The set $S_4$}

Suppose $p \in W_4$.  Then $p_5$ and $q_4$ are
defined and $q_4 \in \partial \Sigma_4$.
Given that $V_5=(0,-4)$ and the $y$-coordinates of
all our points are odd integers, we have
$p_4-q_4=(0,2)+k(0,4)$ for some $k \in \Z$.
Given that $\Sigma_4$ dominates $p_4-q_4$ we have
$k \in \{-1,0\}$.  Hence
$p_4=q_4 \pm (0,2)$.  If
$p_5 \in \partial \Sigma_5$ then
$q_4 \in \partial \Sigma_4$.   
Any vertical line intersects $\Sigma_4$ in a
seqment of length $4$. From this we see that
$p_4$ lies on the centerline of $\Sigma_4$.
That is, $\sigma_4(p)=1/2$.   Given Equation
\ref{identity4}, we get
$$S_4 \subset \{x+y-z=A\} \cup
\{x+y-z=2+A\}.$$

\subsubsection{The Set $S_5$}

Suppose that $p \in W_5$.  Then
 $p_6$ and $q_3$ are defined, and
$q_3 \in \partial \Sigma_3$.
Given that $V_6=-V_4=(-2,2)$, we see that
$$p_3-q_3= \epsilon(0,2)+k(2,2); \hskip 30pt \epsilon \in \{-1,1\}; \hskip 20 pt 
k \in \Z.$$
The criterion that $\Sigma_3$ dominates a vector $(x,y)$ is that
$|x+Ay|<2+2A$.  

$\Sigma_3$ dominates the vector $q_3-p_3$.  
If $\epsilon=1$ then 
$|2k+2+2Ak|<2+2A$, forces $k \in \{-1,0\}$.
If $\epsilon=-1$, then the condition
$|2k-2+2Ak|<2+2A$ forces $k \in \{0,1\}$.
Hence $p_3-q_3$ is one of the vectors
$(\pm 2,0)$ or $(0,\pm 2)$.  Now we have a
case-by-case analysis.

Suppose that $q_3$ lies in the right boundary of $\Sigma_3$.
Then we have 
either $p_3=q_3-(2,0)$ or $p_3=q_3+(0,2)$.
Any horizontal line intersects $\Sigma_3$ in a strip
of width $2+2A$.  So,
$\sigma_3(p)$ equals either $1/(1+A)$ or $A/(1+A)$ depending on
whether or not $p_3=q_3-(2,0)$ or $p_3=q_3+(0,2)$.
A similar analysis reveals the same two values
when $q_3$ lies on the left boundary of $\Sigma_3$.
Given Equation \ref{identity3} we get
$$S_5 \subset \{x=A\} \cup \{x=1\}.$$

\subsubsection{The Set $S_6$}

Suppose that $p \in W_6$.  Then
 $p_7$ and $q_2$ are defined, and
$q_2 \in \partial \Sigma_2$.
We have
\begin{equation}
p_2-q_2=(p_3-q_3) + k(2+2A,0).
\end{equation}
The criterion that $\Sigma_3$ dominates a vector $(x,y)$ is that
$|x-Ay|<2+2A$.  

Let $X_1,...,X_4$ be the possible values for
$p_3-q_3$, as determined in the previous section.
Using the values of the vectors $X_j$, and the fact that
$\Sigma_2$ dominates $p_2-q_2$, we see
that 
\begin{equation}
\label{choices}
p_2-q_2=X_j + \epsilon (2A,2); \hskip 30 pt \epsilon \in \{-1,0,1\};
\hskip 30 pt j \in \{1,2,3,4\}.
\end{equation}

Note that the vector $(2A,2)$ is parallel to the
boundary of $\Sigma_2$.   Hence, for the purposes
of computing $\sigma_2(p)$, this vector plays no role.
Essentially the same calculation as in the previous
section now gives us the same choices for
$\sigma_2(p)$ as we got for $\sigma_3(p)$ in the
previous section.
Given Equation \ref{identity2} we get
$$S_6 \subset \{y=A\} \cup \{y=1\}.$$

\subsubsection{The Set $S_7$}

Suppose that $p \in W_7$.  Then
 $p_8$ and $q_1$ are defined, and
$q_1 \in \partial \Sigma_1$.
We have
\begin{equation}
p_1-q_1=(p_2-q_2)+k(-2,2).
\end{equation}
Note that the vector $(2,2)$ is parallel to $\Sigma_1$.
For the purposes of finding $\sigma_1(p)$, we can 
do our computation modulo $(2,2)$.  
For instance, $(2,-2) \equiv (0,4)$ mod $(2,2)$.
Given Equation \ref{choices}, we have
\begin{equation}
p_1-q_1=\epsilon_1 (0,2)+\epsilon_2 (2A,2)+k(0,4) \hskip 10 pt {\rm mod\/} \hskip 10 pt (2,2).
\end{equation}
Here $\epsilon_1,\epsilon_2 \in \{-1,0,1\}$.   Given that any vertical
line intersects $\Sigma_1$ in a segment of length $4$, we see that the only
choices for $\sigma_1(p)$ are
$$\bigg[\frac{k}{2}+2\epsilon A\bigg]; \hskip 30 pt \epsilon \in \{-1,0,1\}; \hskip 20 pt k \in \Z.$$
Given Equation \ref{identity1} we see that
$S_7 \subset \{z=A\} \cup  \{z=1-A\}.$

\subsection{Proof of The Second Statement}

Our analysis above establishes the first statement of the
Hyperplane Lemma.   For the second statement,
suppose that $d(\mu_+(p,A),\widetilde S)=\epsilon$.
Given Equations \ref{identity1}, \ref{identity2}, \ref{identity3},
and \ref{identity4}, we have 
$$\theta_j(p) \geq \epsilon; \hskip 30 pt j=1,2,3,4.$$
Given our analysis of the remaining points using symmetry,
the same bound holds for
$j=5,6,7,8$.  In these cases, $\theta_j(p,A)$ is a linear
function of the distance from $\mu_+(p,A)$ to
$S_{j-1}$, and the constant of proportionality 
is the same as it is for the index $9-j$.

\newpage

\section{Some Formulas}
\label{computations}

\subsection{Formulas for the Pinwheel Map}
\label{pinwheelformulas}

In this section we explain how to implement the pinwheel map.
We define
$$
V_1=(0,4); \hskip 20 pt
V_2=(-2,2);$$
\begin{equation}
V_3=(-2-2A,0); \hskip 20 pt
V_4=(-2,-2).
\end{equation}
Next, we define vectors
$$
W_1=\frac{1}{4}(-1, 1,3); \hskip 20 pt 
W_2=\frac{1}{2+2A}(-1,A,A);
$$
\begin{equation}
W_3=\frac{1}{2+2A}(-1,-A,A); \hskip 20 pt
W_4=\frac{1}{4}(-1,-1,3); \hskip 20 pt
\end{equation}
For a point $p \in \R^2$, we define
\begin{equation}
F_j(p)=W_{j} \cdot (p_1,p_2,1).
\end{equation}
$F(j,p)$ measures the position of $p$ relative to the strip
$\Sigma_j$.  This quantity lies in $(0,1)$ iff $p$ lies in the
interior of $\Sigma_j$.
\newline
\newline
{\bf Example:\/}
Let $p=(2A,1)$ and $q=(-2,1)$, we compute that
$$F_2(p)=\frac{1}{2+2A}(-1,A,A) \cdot (2A,1,1)=0.$$
$$F_2(q)=\frac{1}{2+2A}(-1,A,A) \cdot (-2,1,1)=1.$$
This checks out, because $p$ lies in one component of $\partial \Sigma_2$ and
$q$ lies in the other component of $\partial \Sigma_2$.
\newline

Here is a formula for our strip maps.
\begin{equation}
E_j(p)=p-{\rm floor\/}(F_j(p)) V_j.
\end{equation}

If we set $V_{4+j}=-V_j$ and $F_{j+4}=-F_j$ then we get the nice formulas
\begin{equation}
dF_j(V_j)=dF_j(V_{j+1})=1.
\end{equation}
with indices taken mod $8$.

\subsection{The Reduction Algorithm}
\label{master1}

Let $A \in (0,1)$ and $\alpha \in \R_+$ and
$(m,n) \in \Z^2$ be a point above the
baseline of $\Gamma_{\alpha}(A)$. In this section
we describe how we compute the points
$$\mu_{\pm}(M_{\alpha}(m,n)).$$
This algorithm will be important
when we prove the Copy Theorems in Part IV of the
monograph.

\begin{enumerate}
\item Let $z=Am+n+\alpha$.
\item Let $Z={\rm floor\/}(z)$.
\item Let $y=z+Z$.
\item Let $Y={\rm floor\/}(y/(1+A))$.
\item Let $x=y-Y(1-A)-1$.
\item Let $X={\rm floor\/}(x/(1+A))$.
\end{enumerate}
We then have
\begin{equation}
\mu_-(M_{\alpha}(m,n))=\left(\matrix{x-(1+A)X \cr y-(1+A)Y \cr z-Z}\right)
\end{equation}
The description of $\mu_+$ is identical, except that the
third step above is replaced by 
\begin{equation}
y=z+Z+1.
\end{equation}

\noindent
{\bf Example:\/}
Referring to \S \ref{example}, consider the case when
$A=3/5$ and $\alpha=1/10$ and $(m,n)=(4,2)$.  We get
$$z=\frac{9}{2}; \hskip 30 pt Z=4; \hskip 30 pt
y=\frac{17}{2}; \hskip 30 pt Y={\rm floor\/}\bigg(\frac{17/2}{8/5}\bigg)=5.$$
$$x=\frac{17}{2}-5(\frac{2}{5})-1=\frac{11}{2}; \hskip 30 pt
X={\rm floor\/}\bigg(\frac{11/2}{8/5}\bigg)=3.$$
$$\mu_-(M(4,2))=\Big(\frac{11}{2}-3(\frac{8}{5}),
\frac{17}{2}-5(\frac{8}{5}),
\frac{9}{2}-4\Big)= \Big(\frac{7}{10},\frac{1}{2},\frac{1}{2}\Big).$$

\subsection{Computing the Partition}
\label{master2}

Here we describe how Billiard King applies the
Master Picture Theorem.

\subsubsection{Step 1}

Suppose $(a,b,c) \in R_A$ lies in the range of
$\mu_+$ or $\mu_-$.
Now we describe how to attach a $5$-tuple
$(n_0,..,n_4)$ to $(a,b,c)$. 

\begin{itemize}
\item Determining $n_0$:
\begin{itemize}
\item If we are interested in $\mu_+$, then $n_0=0$.
\item If we are interested in $\mu_-$, then $n_0=1$.
\end{itemize}
\item Determining $n_1$:
\begin{itemize}
\item If $c<A$ and $c<1-A$ then $n_1=0$.
\item If $c>A$ and $c<1-A$ then $n_1=1$.
\item If $c>A$ and $c>1-A$ then $n_1=2$.
\item If $c<A$ and $c>1-A$ then $n_1=3$.
\end{itemize}
\item Determining $n_2$:
\begin{itemize}
\item If $a \in (0,A)$ then $n_2=0$.
\item If $a \in (A,1)$ then $n_2=1$.
\item If $a \in (1,1+A)$ then $n_2=2$.
\end{itemize}
\item Determining $n_3$.
\begin{itemize}
\item If $b \in (0,A)$ then $n_3=0$.
\item If $b \in (A,1)$ then $n_3=1$.
\item If $b \in (1,1+A)$ then $n_3=2$.
\end{itemize}
\item Determining $n_4$.
\begin{itemize}
\item Let $t=a+b-c$.
\item Let $n_4={\rm floor\/}(t-A)$.
\end{itemize}
\end{itemize}
Notice that each $5$-tuple $(n_0,...,n_4)$ corresponds
to a (possibly empty) convex polyhedron in $R_A$.  The polyhedron
doesnt depend on $n_0$.  It turns out that
this polyhedron is empty unless $n_4 \in \{-2,-1,0,1,2\}$.

\subsubsection{Step 2}

Let $n=(n_0,...,n_4)$.  We now describe two functions
$\epsilon_1(n) \in \{-1,0,1\}$ and $\epsilon_2(n) \in \{-1,0,1\}$.
\newline
\newline
\noindent
Here is the definition of $\epsilon_1(n)$.
\begin{itemize}
\item If $n_0+n_4$ is even then:
\begin{itemize}
\item If $n_2+n_3=4$ or $x_2<x_3$ set $\epsilon_1(n)=-1$.
\end{itemize}
\item If $n_0+n_4$ is odd then:
\begin{itemize}
\item If $n_2+n_3=0$ or $x_2>x_3$ set $\epsilon_1(n)=+1$.
\end{itemize}
\item Otherwise set $\epsilon_1(n)=0$.
\end{itemize}
\noindent
Here is the definition of $\epsilon_2(n)$.

\begin{itemize}
\item If $n_0=0$ and $n_1 \in \{3,0\}$.
\begin{itemize}
\item If $n_2=0$ let $\epsilon_2(n)=1$.
\item If $n_2=1$ and $n_4 \not = 0$ let $\epsilon_2(n)=1$.
\end{itemize}

\item If $n_0=1$ and $n_1 \in \{0,1\}$.
\begin{itemize}
\item if $n_2>0$ and $n_4 \not =0$ let $\epsilon_2(n)=-1$.
\item If $n_2<2$ and $n_3=0$ and $n_4=0$ let $\epsilon_2(n)=1$.
\end{itemize}

\item If $n_0=0$ and $n_1 \in \{1,2\}$.
\begin{itemize}
\item If $n_2<2$ and $n_4 \not = 0$ let $\epsilon_2(n)=1$.
\item If $n_2>0$ and $n_3=2$ and $n_4=0$ let $\epsilon_2(n)=-1$.
\end{itemize}

\item If $n_0=1$ and $n_1 \in \{2,3\}$.
\begin{itemize}
\item If $n_2=2$ let $\epsilon_2(n)=-1$.
\item If $n_2=1$ and $n_4 \not = 0$ let $\epsilon_2(n)=-1$.
\end{itemize}

\item Otherwise let $\epsilon_2(n)=0$.
\end{itemize}

\subsubsection{Step 3}
Let $A \in (0,1)$ be any parameter and let
$\alpha>0$ be some parameter such that $\alpha \not \in 2\Z[A]$.
Given any lattice point $(m,n)$ we perform the
following construction.
\begin{itemize}
\item Let $(a_{\pm},b_{\pm},c_{\pm})=\mu_{\pm}(A,m,n)$. See \S \ref{master1}.
\item Let $n_{\pm}$ be the
$5$-tuple associated to $(a_{\pm},b_{\pm},c_{\pm})$.  
\item Let $\epsilon_1^{\pm}=\epsilon_1(n_{\pm})$ and $\epsilon_2^{\pm}=\epsilon_2(n_{\pm})$.
\end{itemize}
The Master Picture Theorem says that the two edges
of $\Gamma_{\alpha}(m,n)$ incident to $(m,n)$ are
$(m,n)+(\epsilon_1^{\pm},\epsilon_2^{\pm})$.

\subsection{The List of Polytopes}
\label{polytopelist}

Referring to the simpler partition from \S \ref{geometric2},
we list the $14$ polytopes that partition $R_+$. 
In each case, we list some vectors, followed by the pair
$(\epsilon_1,\epsilon_2)$ that the polytope determines.
$$ $$

\lefteqn{\left[\matrix{0 \cr 0 \cr 0 \cr 0}\right]
\hskip 5 pt\left[\matrix{0 \cr 0 \cr 0 \cr 1}\right]
\hskip 5 pt\left[\matrix{0 \cr 0 \cr 1 \cr 0}\right]
\hskip 5 pt\left[\matrix{0 \cr 1 \cr 0 \cr 1}\right]
\hskip 5 pt\left[\matrix{0 \cr 1 \cr 1 \cr 1}\right]
\hskip 5 pt\left[\matrix{1 \cr 0 \cr 0 \cr 1}\right]
\hskip 5 pt\left[\matrix{1 \cr 0 \cr 1 \cr 0}\right]
\hskip 5 pt\left[\matrix{1 \cr 0 \cr 1 \cr 1}\right]
\hskip 5 pt\left[\matrix{1 \cr 1 \cr 1 \cr 1}\right]
\hskip 10 pt (1,1)}

\lefteqn{\left[\matrix{0 \cr 0 \cr 0 \cr 0}\right]
\hskip 5 pt\left[\matrix{0 \cr 1 \cr 0 \cr 0}\right]
\hskip 5 pt\left[\matrix{0 \cr 1 \cr 0 \cr 1}\right]
\hskip 5 pt\left[\matrix{0 \cr 2 \cr 0 \cr 1}\right]
\hskip 5 pt\left[\matrix{0 \cr 2 \cr 1 \cr 1}\right]
\hskip 5 pt\left[\matrix{1 \cr 1 \cr 0 \cr 1}\right]
\hskip 5 pt\left[\matrix{1 \cr 1 \cr 1 \cr 1}\right]
\hskip 5 pt\left[\matrix{1 \cr 2 \cr 1 \cr 1}\right]
\hskip 10 pt (-1,1)}

\lefteqn{\left[\matrix{0 \cr 1 \cr 0 \cr 0}\right]
\hskip 5 pt\left[\matrix{0 \cr 1 \cr 1 \cr 0}\right]
\hskip 5 pt\left[\matrix{1 \cr 1 \cr 1 \cr 0}\right]
\hskip 5 pt\left[\matrix{1 \cr 1 \cr 1 \cr 1}\right]
\hskip 5 pt\left[\matrix{1 \cr 2 \cr 1 \cr 1}\right]
\hskip 5 pt\left[\matrix{2 \cr 1 \cr 1 \cr 1}\right]
\hskip 10 pt (-1,-1)}
\lefteqn{\left[\matrix{0 \cr 1 \cr 0 \cr 0}\right]
\hskip 5 pt\left[\matrix{0 \cr 2 \cr 0 \cr 1}\right]
\hskip 5 pt\left[\matrix{1 \cr 0 \cr 0 \cr 0}\right]
\hskip 5 pt\left[\matrix{1 \cr 1 \cr 0 \cr 0}\right]
\hskip 5 pt\left[\matrix{1 \cr 1 \cr 0 \cr 1}\right]
\hskip 5 pt\left[\matrix{1 \cr 1 \cr 1 \cr 0}\right]
\hskip 5 pt\left[\matrix{1 \cr 2 \cr 0 \cr 1}\right]
\hskip 5 pt\left[\matrix{1 \cr 2 \cr 1 \cr 1}\right]
\hskip 10 pt (0,1)}
\lefteqn{\left[\matrix{0 \cr 0 \cr 0 \cr 0}\right]
\hskip 5 pt\left[\matrix{0 \cr 0 \cr 1 \cr 0}\right]
\hskip 5 pt\left[\matrix{0 \cr 1 \cr 0 \cr 1}\right]
\hskip 5 pt\left[\matrix{0 \cr 1 \cr 1 \cr 0}\right]
\hskip 5 pt\left[\matrix{0 \cr 1 \cr 1 \cr 1}\right]
\hskip 5 pt\left[\matrix{0 \cr 2 \cr 1 \cr 1}\right]
\hskip 5 pt\left[\matrix{1 \cr 0 \cr 1 \cr 0}\right]
\hskip 5 pt\left[\matrix{1 \cr 1 \cr 1 \cr 1}\right]
\hskip 10 pt (0,1)}
\lefteqn{\left[\matrix{0 \cr 0 \cr 0 \cr 0}\right]
\hskip 5 pt\left[\matrix{0 \cr 1 \cr 0 \cr 1}\right]
\hskip 5 pt\left[\matrix{1 \cr 0 \cr 0 \cr 1}\right]
\hskip 5 pt\left[\matrix{1 \cr 1 \cr 0 \cr 1}\right]
\hskip 5 pt\left[\matrix{1 \cr 1 \cr 1 \cr 1}\right]
\hskip 10 pt (0,1)}
\lefteqn{\left[\matrix{0 \cr 0 \cr 0 \cr 1}\right]
\hskip 5 pt\left[\matrix{0 \cr 0 \cr 1 \cr 0}\right]
\hskip 5 pt\left[\matrix{0 \cr 0 \cr 1 \cr 1}\right]
\hskip 5 pt\left[\matrix{0 \cr 1 \cr 1 \cr 1}\right]
\hskip 5 pt\left[\matrix{1 \cr 0 \cr 1 \cr 1}\right]
\hskip 10 pt (0,1)}
\lefteqn{\left[\matrix{0 \cr 0 \cr 0 \cr 0}\right]
\hskip 5 pt\left[\matrix{0 \cr 1 \cr 0 \cr 0}\right]
\hskip 5 pt\left[\matrix{0 \cr 1 \cr 1 \cr 0}\right]
\hskip 5 pt\left[\matrix{0 \cr 2 \cr 1 \cr 1}\right]
\hskip 5 pt\left[\matrix{1 \cr 1 \cr 1 \cr 1}\right]
\hskip 5 pt\left[\matrix{1 \cr 2 \cr 1 \cr 1}\right]
\hskip 10 pt (-1,0)}
\lefteqn{\left[\matrix{1 \cr 1 \cr 0 \cr 0}\right]
\hskip 5 pt\left[\matrix{1 \cr 2 \cr 0 \cr 1}\right]
\hskip 5 pt\left[\matrix{2 \cr 1 \cr 0 \cr 1}\right]
\hskip 5 pt\left[\matrix{2 \cr 2 \cr 0 \cr 1}\right]
\hskip 5 pt\left[\matrix{2 \cr 2 \cr 1 \cr 1}\right]
\hskip 10 pt (-1,0)}
\lefteqn{\left[\matrix{0 \cr 1 \cr 0 \cr 0}\right]
\hskip 5 pt\left[\matrix{1 \cr 1 \cr 0 \cr 1}\right]
\hskip 5 pt\left[\matrix{1 \cr 1 \cr 1 \cr 1}\right]
\hskip 5 pt\left[\matrix{1 \cr 2 \cr 1 \cr 1}\right]
\hskip 5 pt\left[\matrix{2 \cr 1 \cr 1 \cr 1}\right]
\hskip 10 pt (-1,0)}
\lefteqn{\left[\matrix{1 \cr 0 \cr 0 \cr 0}\right]
\hskip 5 pt\left[\matrix{1 \cr 1 \cr 0 \cr 0}\right]
\hskip 5 pt\left[\matrix{1 \cr 1 \cr 0 \cr 1}\right]
\hskip 5 pt\left[\matrix{1 \cr 1 \cr 1 \cr 0}\right]
\hskip 5 pt\left[\matrix{2 \cr 0 \cr 0 \cr 1}\right]
\hskip 5 pt\left[\matrix{2 \cr 1 \cr 0 \cr 1}\right]
\hskip 5 pt\left[\matrix{2 \cr 1 \cr 1 \cr 1}\right]
\hskip 10 pt (1,0)}
\lefteqn{\left[\matrix{1 \cr 0 \cr 0 \cr 1}\right]
\hskip 5 pt\left[\matrix{1 \cr 0 \cr 1 \cr 0}\right]
\hskip 5 pt\left[\matrix{1 \cr 0 \cr 1 \cr 1}\right]
\hskip 5 pt\left[\matrix{1 \cr 1 \cr 1 \cr 1}\right]
\hskip 5 pt\left[\matrix{2 \cr 0 \cr 1 \cr 1}\right]
\hskip 10 pt (1,0)}

\lefteqn{\left[\matrix{1 \cr 1 \cr 0 \cr 0}\right]
\hskip 5 pt\left[\matrix{1 \cr 1 \cr 0 \cr 1}\right]
\hskip 5 pt\left[\matrix{1 \cr 1 \cr 1 \cr 0}\right]
\hskip 5 pt\left[\matrix{1 \cr 2 \cr 0 \cr 1}\right]
\hskip 5 pt\left[\matrix{1 \cr 2 \cr 1 \cr 1}\right]
\hskip 5 pt\left[\matrix{2 \cr 1 \cr 0 \cr 1}\right]
\hskip 5 pt\left[\matrix{2 \cr 1 \cr 1 \cr 1}\right]
\hskip 5 pt\left[\matrix{2 \cr 2 \cr 1 \cr 1}\right]
\hskip 10 pt (0,0)}

\lefteqn{\left[\matrix{0 \cr 0 \cr 0 \cr 0}\right]
\hskip 5 pt\left[\matrix{0 \cr 1 \cr 0 \cr 0}\right]
\hskip 5 pt\left[\matrix{0 \cr 1 \cr 1 \cr 0}\right]
\hskip 5 pt\left[\matrix{1 \cr 0 \cr 0 \cr 0}\right]
\hskip 5 pt\left[\matrix{1 \cr 0 \cr 0 \cr 1}\right]
\hskip 5 pt\left[\matrix{1 \cr 0 \cr 1 \cr 0}\right]
\hskip 5 pt\left[\matrix{1 \cr 1 \cr 0 \cr 1}\right]
\hskip 5 pt\left[\matrix{1 \cr 1 \cr 1 \cr 0}\right]
\hskip 5 pt\left[\matrix{1 \cr 1 \cr 1 \cr 1}\right]
\hskip 5 pt\left[\matrix{2 \cr 0 \cr 0 \cr 1}\right]
\hskip 5 pt\left[\matrix{2 \cr 0 \cr 1 \cr 1}\right]
\hskip 5 pt\left[\matrix{2 \cr 1 \cr 1 \cr 1}\right]
\hskip 10 pt (0,0)}
$$ $$

\subsection{Calculating with the Polytopes}
\label{examplecalc}

We will illustrate a calculation with the polytopes we have listed.
Let $\iota$ and $\gamma_2$ be the maps from Equation \ref{geometric2}.
$R_+(0,0)$ consists of two polygons, $P_1$ and $P_2$.  These are the last
two listed above.  We will show that
$$\iota(P_2)+(1,1,0,0)=\gamma_2(P_2).$$
As above, the coordinates for $P_2$ are
$$
\left[\matrix{0 \cr 0 \cr 0 \cr 0}\right]
\hskip 5 pt\left[\matrix{0 \cr 1 \cr 0 \cr 0}\right]
\hskip 5 pt\left[\matrix{0 \cr 1 \cr 1 \cr 0}\right]
\hskip 5 pt\left[\matrix{1 \cr 0 \cr 0 \cr 0}\right]
\hskip 5 pt\left[\matrix{1 \cr 0 \cr 0 \cr 1}\right]
\hskip 5 pt\left[\matrix{1 \cr 0 \cr 1 \cr 0}\right]
\hskip 5 pt\left[\matrix{1 \cr 1 \cr 0 \cr 1}\right]
\hskip 5 pt\left[\matrix{1 \cr 1 \cr 1 \cr 0}\right]
\hskip 5 pt\left[\matrix{1 \cr 1 \cr 1 \cr 1}\right]
\hskip 5 pt\left[\matrix{2 \cr 0 \cr 0 \cr 1}\right]
\hskip 5 pt\left[\matrix{2 \cr 0 \cr 1 \cr 1}\right]
\hskip 5 pt\left[\matrix{2 \cr 1 \cr 1 \cr 1}\right]
$$
Recall that
$\iota(x,y,z,A)=(1+A-x,1+A-y,1-z,A).$
For example, $\iota(0,0,0,0)=(1,1,1,0)$.
The coordinates for $\iota(P_2)$ are
$$
\left[\matrix{1 \cr 1 \cr 1 \cr 0}\right]
\hskip 5 pt\left[\matrix{1 \cr 0 \cr 1 \cr 0}\right]
\hskip 5 pt\left[\matrix{1 \cr 0 \cr 0 \cr 0}\right]
\hskip 5 pt\left[\matrix{0 \cr 1 \cr 1 \cr 0}\right]
\hskip 5 pt\left[\matrix{1 \cr 2 \cr 1 \cr 1}\right]
\hskip 5 pt\left[\matrix{0 \cr 1 \cr 0 \cr 0}\right]
\hskip 5 pt\left[\matrix{1 \cr 1 \cr 1 \cr 1}\right]
\hskip 5 pt\left[\matrix{0 \cr 0 \cr 0 \cr 0}\right]
\hskip 5 pt\left[\matrix{1\cr 1 \cr 0 \cr 1}\right]
\hskip 5 pt\left[\matrix{0 \cr 2 \cr 1 \cr 1}\right]
\hskip 5 pt\left[\matrix{0 \cr 2 \cr 0 \cr 1}\right]
\hskip 5 pt\left[\matrix{0 \cr 1 \cr 0 \cr 1}\right]
$$
The coordinates for $\iota(P_2)+(1,1,0,0)$ are
$$
\left[\matrix{2 \cr 2 \cr 1 \cr 0}\right]
\hskip 5 pt\left[\matrix{2 \cr 1 \cr 1 \cr 0}\right]
\hskip 5 pt\left[\matrix{2 \cr 1 \cr 0 \cr 0}\right]
\hskip 5 pt\left[\matrix{1 \cr 2 \cr 1 \cr 0}\right]
\hskip 5 pt\left[\matrix{2 \cr 3 \cr 1 \cr 1}\right]
\hskip 5 pt\left[\matrix{1 \cr 2 \cr 0 \cr 0}\right]
\hskip 5 pt\left[\matrix{2 \cr 2 \cr 1 \cr 1}\right]
\hskip 5 pt\left[\matrix{1 \cr 1 \cr 0 \cr 0}\right]
\hskip 5 pt\left[\matrix{2 \cr 2 \cr 0 \cr 1}\right]
\hskip 5 pt\left[\matrix{1 \cr 3 \cr 1 \cr 1}\right]
\hskip 5 pt\left[\matrix{1 \cr 3 \cr 0 \cr 1}\right]
\hskip 5 pt\left[\matrix{1 \cr 2 \cr 0 \cr 1}\right]
$$
We have $\gamma_2(x,y,z,A)=(x+1-A,y+1+A,z,A)$.  For instance,
we compute that $\gamma_2(0,0,0,0)=(1,1,0,0)$.  
The coordinates for $\gamma(P_2)$ are
$$
\left[\matrix{1 \cr 1 \cr 0 \cr 0}\right]
\hskip 5 pt\left[\matrix{1 \cr 2 \cr 0 \cr 0}\right]
\hskip 5 pt\left[\matrix{1 \cr 2 \cr 1 \cr 0}\right]
\hskip 5 pt\left[\matrix{2 \cr 1 \cr 0 \cr 0}\right]
\hskip 5 pt\left[\matrix{1 \cr 2 \cr 0 \cr 1}\right]
\hskip 5 pt\left[\matrix{2 \cr 1 \cr 1 \cr 0}\right]
\hskip 5 pt\left[\matrix{1 \cr 3 \cr 0 \cr 1}\right]
\hskip 5 pt\left[\matrix{2 \cr 2 \cr 1 \cr 0}\right]
\hskip 5 pt\left[\matrix{1 \cr 3 \cr 1 \cr 1}\right]
\hskip 5 pt\left[\matrix{2 \cr 2 \cr 0 \cr 1}\right]
\hskip 5 pt\left[\matrix{2 \cr 2 \cr 1 \cr 1}\right]
\hskip 5 pt\left[\matrix{2 \cr 3 \cr 1 \cr 1}\right]
$$
These are the same vectors as listed for
$\iota(P_2)+(1,1,0,0)$ but in a different order.

\subsection{The Phase Portrait}

Here we explain how to derive the phase portrait
described in Figure 2.4.   Our discussion
refers to \S \ref{computations}.
Consider the two rectangles
$$Q_+=\{(t,t+1,t)|\ t \in (0,1)\} \times [0,1];$$
$$Q_-=\{(t-1,t,t)|\ t \in (0,1)\} \times [0,1].$$
Intersect $Q_{\pm}$ with the polytope $R_{\pm}$.
These intersections partition $Q_+$ and $Q_-$
into a small finite number of polygons.  
The partition of $Q_{\pm}$ tells the behavior of
$\Psi^{\pm}$ on points of $(0,2) \times \{1\}$.
Bu symmetry, the partition of $Q_{\mp}$ tells the
behavior of $\Psi^{\pm}$ on $(0,2) \times \{-1\}$.
The partition of $Q_{\pm}$ gives us the information
needed to build Figure 2.4.  Given the simplicity
of the partitions involved, we can determine the picture
just by plotting (say) $10000$ fairly dense points in
our rectangles.  This is what we do.

\newpage

\noindent
{\bf {\huge Part III\/} \/}
\newline

In this part of the monograph we use the Master Picture
Theorem to prove all the results quoted in 
Part I of the monograph. 

\begin{itemize}
\item In \S \ref{embeddingproof} we prove the Embedding Theorem.
\item In \S \ref{symm1proof} we prove some results about
the symmetries of the arithmetic graph and the hexagrid.
\item In \S \ref{odddoor} we establish some information
about the doors.  These special points were defined
in connection with the Hexagrid Theorem.
\item In \S \ref{hex1proof} we prove Statement 1 of the
Hexagrid theorem, namely that the arithmetic graph does not
cross any floor lines.
\item In \S \ref{hex2proof} we prove Statement 2 of the
Hexagrid theorem, namely that the arithmetic graph only
crosses the walls near the doors. 
The two statements of the Hexagrid Theorem have
similar proofs, though Statement 2 has a more
elaborate proof.
\item In \S \ref{barrier} we prove a variant of Statement 1 of
the Hexagrid Theorem.  We call the result that Barrier Theorem.
Though we don't need this result until Part VI, the proof
fits best right after the proof of the Hexagrid Theorem.
\end{itemize}

Many of the proofs in this part of the monograph require us to 
prove various disjointness results about some $4$ dimensional
polytopes.  We will give short computer-aided proofs of
these disjointness results.  The proofs only involve
a small amount of integer arithmetic.  An energetic
mathematician could do them all by hand in an afternoon.
To help make the proofs surveyable, we will include
extensive computer pictures of $2$ dimensional
slices of our polytopes.  These pictures, all reproducible
on Billiard King, serve as sanity checks for the
computer calculations.

\newpage

\section{Proof of the Embedding Theorem}
\label{embeddingproof}

Let $\widehat \Gamma=\widehat \Gamma_{\alpha}(A)$ be the
arithmetic graph for a parameter $A$ and some
number $\alpha \not \in 2\Z[A]$.  In this chapter we 
prove that $\widehat \Gamma$ is a disjoint union of
embedded polygons and infinite polygonal arcs.
This is the Embedding Theorem.

\subsection{Step 1}

We will first prove that every nontrivial vertex
of $\widehat \Gamma$ has valence $2$.
Each point $p  \in \widehat \Gamma$ is connected to two
points $q_+$ and $q_-$.   Hence, each non-trivial
vertex has valence either $1$ or $2$.
The following two cases are
the only cases that lead to valence $1$ vertices:
\begin{itemize}
\item $p=q_+$ and $q_+ \not = q_-$.
\item $q_+=q_-$ and $q_{\pm} \not = p$.
\end{itemize}
The following lemma rules out the first of these cases.

\begin{lemma}
\label{step1}
If $p=q_+$ or $p=q_-$ then $p=q_+=q_-$.
\end{lemma}

\startproof
Our proof refers to \S \ref{geometric2}.
Recall that $R_+(0,0)$ consists of $2$ convex integer polytopes.
Likewise $R_-(0,0)$ consists of $2$ convex integer polytopes.
It suffices to show that
\begin{equation}
\label{simul}
(t,t+1,t) \in R_+(0,0) \hskip 15 pt \Longleftrightarrow \hskip 15 pt
(t-1,t,t) \in R_-(0,0).
\end{equation}
This is equivalent to the statement that
$$R_-(0,0)+(1,1,0,0) \subset \Lambda R_+(0,0)$$
Here $\Lambda R_+(0,0)$ is the orbit of
$R_+(0,0)$ under the action of $\Lambda$.
Let $\iota$ be the involution from Equation \ref{iota}.
Recall that $R_-(0,0)=\iota(R_+(0,0))$.  Hence, 
Equation \ref{simul} equivalent to the statement that
\begin{equation}
\iota(R_+(0,0))+(1,1,0,0) \subset
\Lambda R_+(0,0).
\end{equation}
Let $P_1$ and $P_2$ denote the two polytopes comprising
$R_+(0,0)$, as listed at the end of \S \ref{masterpicture}.
Let $\gamma_2$ be the element of $\Lambda$ described in \S \ref{geometric2}.
We compute that
\begin{equation}
\iota(P_1)+(1,1,0,0)=P_1; \hskip 50 pt
\iota(P_2)+(1,1,0,0))=\gamma_2(P_2).
\end{equation}
We did the second calculation in \S \ref{examplecalc}, and the
first computation is similar.
This does it for us.
\endproof

\subsection{Step 2}

Our next goal is to rule out the possibility that
$p \not = q_{\pm}$, but $q_+=q_-$.    This situation happens
iff there is some $(\epsilon_1,\epsilon_2) \in \{-1,0,1\}$
such that
\begin{equation}
\label{bad1}
\Lambda R_+(\epsilon_1,\epsilon_2) \cap \big(R_-(\epsilon_1,\epsilon_2)+(1,1,0,0)\big)
\not = \emptyset.
\end{equation}
A visual inspection and/or a compute computer search -- we did both --
reveals that at least one of the two sets above
is empty unless $(\epsilon_1,\epsilon_2)$ is one of
\begin{equation}
(1,1); \hskip 30 pt (-1,-1); \hskip 30 pt (1,0); \hskip 30 pt (-1,0).
\end{equation}

To rule out Equation \ref{bad1} for each of these pairs,
we need to consider all possible pairs
$(P_1,P_2)$ of integral convex polytopes such that
\begin{equation}
P_1 \subset \Lambda R_+(\epsilon_1,\epsilon_2); \hskip 30 pt
P_2 \subset (R_-(\epsilon_1,\epsilon_2)+(1,1)\big)
\end{equation}
Recall that $\Lambda$ is generated by the three elements
$\gamma_1,\gamma_2,\gamma_3$.   Let
$\Lambda' \subset \Lambda$ denote the subgroup generated by
$\gamma_1$ and $\gamma_2$.   We also define
$\Lambda'_{10} \subset \Lambda'$ by the equation
\begin{equation}
\label{gammafinite}
\Lambda'_{10}=\{a_1 \gamma_1 + a_2 \gamma_2|\ |a_1|,|a_2| \leq 10\}.
\end{equation}

\begin{lemma}
\label{far1}
Let $\gamma \in \Lambda-\Lambda'$. 
Suppose that
$$P_1=\gamma(Q_1);  \hskip 20 pt
Q_1 \subset R_+(\epsilon_1,\epsilon_2); \hskip 20 pt
P_2 \subset R_-(\epsilon_1,\epsilon_2)+(1,1,0,0).$$
Then $P_1$ and $P_2$ have disjoint interiors.
\end{lemma}

\startproof
The third coordinates of points in
$P_1$ lies between $n$ and $n+1$ for
some $n \not = 0$ whereas the third coordinates of 
points in $P_2$ lie in $[0,1]$.
\endproof

\begin{lemma}
\label{far2}
Let $\gamma \in \Lambda'-\Lambda'_{10}$.
$$P_1=\gamma(Q_1);  \hskip 20 pt
Q_1 \subset R_+(\epsilon_1,\epsilon_2); \hskip 20 pt
P_2 \subset R_-(\epsilon_1,\epsilon_2)+(1,1,0,0).$$
Then $P_1$ and $P_2$ have disjoint interiors.
\end{lemma}

\startproof
$Q_1$ is contained in the ball of radius $4$ about
$P_2$, but $\gamma$ moves this ball entirely off itself.
\endproof

The last two results leave us with a finite problem.
Given a pair $(\epsilon_1,\epsilon_2)$ from our list above, and
$$\gamma \in \Lambda'_{10}; \hskip 10 pt
P_1=\gamma(Q_1);  \hskip 10 pt
Q_1 \subset R_+(\epsilon_1,\epsilon_2); \hskip 10 pt
P_2 \subset R_-(\epsilon_1,\epsilon_2)+(1,1,0,0),$$
we produce a vector 
\begin{equation}
w=w(P_1,P_2) \in \{-1,0,1\}^4
\end{equation}
such that
\begin{equation}
\max_{v \in {\rm vtx\/}(P_1)} v \cdot w  \leq 
\min_{v \in {\rm vtx\/}(P_2)} v \cdot w.
\end{equation}
This means that a hyperplane separates the interior
of $P_1$ from $P_2$.  In each case we find
$v(P_1,P_2)$ by a short computer search, and perform
the verification using integer arithmetic.  It is
a bit surprising to us that such a simple vector
works in all cases, but that is how it works out.

Using Billiard King, the interested reader can 
draw arbitrary $(z,A)$ slices of the sets
$\Lambda R_+(\epsilon_1,\epsilon_2)$ and
$\Lambda R_-(\epsilon_1,\epsilon_2)+(1,1,0,0)$,
and see that the interiors of the polygons from
the first set are disjoint from the interiors of
the polygons from the second set.
We will illustrate this with pictures in \S \ref{visualtour}.

\subsection{Step 3}

Given that every nontrivial vertex of $\widehat \Gamma$ has valence $2$,
and also that the edges of $\widehat \Gamma$ have length at most
$\sqrt 2$, the only way that $\widehat \Gamma$ can fail to be embedded
is if there is situation like the one shown in Figure 12.1.

\begin{center}
\psfig{file=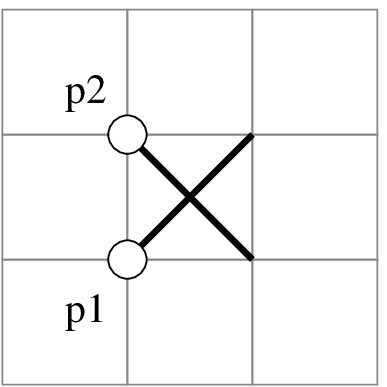}
\newline
{\bf Figure 12.1:\/} Embedding Failure
\end{center} 

Let $M_+$ and $M_-$ be the maps from \S
\ref{mptnotation}.
Given the Master Picture Theorem, this situation arises only in
the following $4$ cases:

\begin{itemize}
\item $M_+(p_1) \in \Lambda R_+(1,1)$ and $M_+(p_2) \in \Lambda R_+(1,-1)$.
\item $M_-(p_1) \in \Lambda R_-(1,1)$ and $M_-(p_2) \in \Lambda R_-(1,-1)$.
\item $M_-(p_1) \in \Lambda R_-(1,1)$ and $M_+(p_2) \in \Lambda R_+(1,-1)$.
\item $M_+(p_1) \in \Lambda R_+(1,1)$ and $M_-(p_2) \in \Lambda R_-(1,-1)$.
\end{itemize}

Note that $p_2=p_1+(0,1)$ and hence
\begin{equation}
\label{shift}
M_{\pm}(p_2)=M_{\pm}(p_1)+(1,1,1,0)\ {\rm mod\/}\ \Lambda.
\end{equation}

In particular, the two points $M_{\rm}(p_1)$ and $M_{\rm}(p_2)$ lie
in the same fiber of $R_{\rm}$ over the $(z,A)$ square.  We
inspect the picture and see that this situation never occurs
for the types $(1,1)$ and $(1,-1)$.  Hence,
Cases 1 and 2 do not occur.
More inspection shows that there are $R_+(1,-1)=\emptyset$.  Hence,
Case 3 does not occus.  This leaves Case 4, the only
nontrivial case.

Case 4 leads to the statement that
$$
(t,t,t,A)+(0,1,0,0) \in \Lambda R_+(1,1);$$
\begin{equation}
(t,t,t,A)-(1,0,0,0)+(1,1,1,0)=(t,t,t,A)+(0,1,1,0) \in \Lambda R_-(1,-1).
\end{equation}
Setting $p$ equal to the first of the two points above, we get
\begin{equation}
\label{basic2}
p \in \Lambda R_+(1,1); \hskip 30 pt
p+(0,0,1,0) \in \Lambda R_-(1,-1).
\end{equation}
Letting $\gamma_3 \in \Lambda$ be as in Equation
\ref{affine}, we have
\begin{equation}
\label{basic3}
p+(1,1,0,0)=\gamma_3^{-1}(p+(0,0,1,0)) \in \Lambda R_-(1,-1).
\end{equation}

For any subset $S \subset \widetilde R$, we have
\begin{equation}
\label{commute}
(\Lambda S) + (a,b,c,0)= \Lambda(S+(a,b,c,0)).
\end{equation}
The point here is that $\Lambda$ acts as a group of translations
on each set of the form $R^3 \times \{A\}$, and addition by
$(x,y,z,0)$ commutes with this action on every such set.
Equations \ref{basic3} and \ref{commute} combine to 
give
\begin{equation}
p \in \Lambda\Big(R_-(1,-1)-(1,1,0,0)\Big)
\end{equation}
Now we see that
$$
\Lambda R_+(1,1) \cap \Lambda(R_-(1,-1)-(1,1,0,0)) \not = \emptyset.
$$
Since the whole picture is $\Lambda$-equivariant, we have
\begin{equation}
\Lambda R_+(1,1) \cap (R_-(1,-1)-(1,1,0,0)) \not = \emptyset.
\end{equation}
We mean that there is a pair $(P_1,P_2)$ of polytopes, with
$P_1$ in the first set and $P_2$ in the second set,
such that $P_1$ and $P_2$ do not have disjoint interiors.

We rule out this intersection using exactly the same method as
in Step 2.  In \S \ref{visualtour} we illustrate this with
a convincing picture.

\subsection{A Visual Tour}
\label{visualtour}

The theoretical part of our proof
amounts to reducing the Embedding Theorem
to the statement that finitely many pairs of polytopes
have disjoint interiors.  The computer-aided
part of the proof amounts to verifying the disjointness
finitely many times.  Our verification used a very
fragile disjointness test.
We got a lucky, because
many of our polytope pairs share a $2$-dimensional face.
Thus, a separating hyperplane
has to be chosen very carefully.  Needless
to say, if our simple-minded approach did not
work, we would have used a more robust
disjointness test.

If we could write this monograph on
$4$-dimensional paper, we could simply replace
the computer-aided part of the proof with a direct appeal
to the visual sense.  Since we don't have
$4$-dimensional paper, we need to rely on the
computer to ``see'' for us.  In this case,
``seeing'' amounts to finding a hyperplane that
separates the interiors of the two polytopes.
In other words, we are getting the computer
to ``look'' at the pair of polytopes in such
a way that one polytope appears on one side
and the other polytope appears on the other side.

We do not have $4$ dimensional paper, but we can draw slices of
all the sets we discussed above.  The interested user of
Billiard King can see any desired slice.  We will just draw
typical slices.  In our pictures below, we will
draw the slices of $R_+$ with dark shading and the
slices of $R_-$ with light shading.  in our discussion,
the {\it base space\/} $B$ refers to the $(z,A)$ square
over which our picture fibers.
Let $B_j$ denote the $j$th component of $B$, as determined
by the characteristic $n_1$ discussed in \S \ref{master2}.

In reference to Step 2, our pictures for the pair $(-\epsilon_1,-\epsilon_2)$
look like rotated versions of the pictures for the pair $(\epsilon_1,\epsilon_2)$.
Accordingly, we will just draw pictures for $(1,1)$ and $(1,0)$.

Figure 12.2 shows a
slice of $\Lambda R_+(1,1)$ and $\Lambda(R_-(1,1)+(1,1,0,0))$
over $B_0$.  Both slices are nonempty over $B_1$ as well,
and the picture is similar.  

\begin{center}
\resizebox{!}{2.6in}{\includegraphics{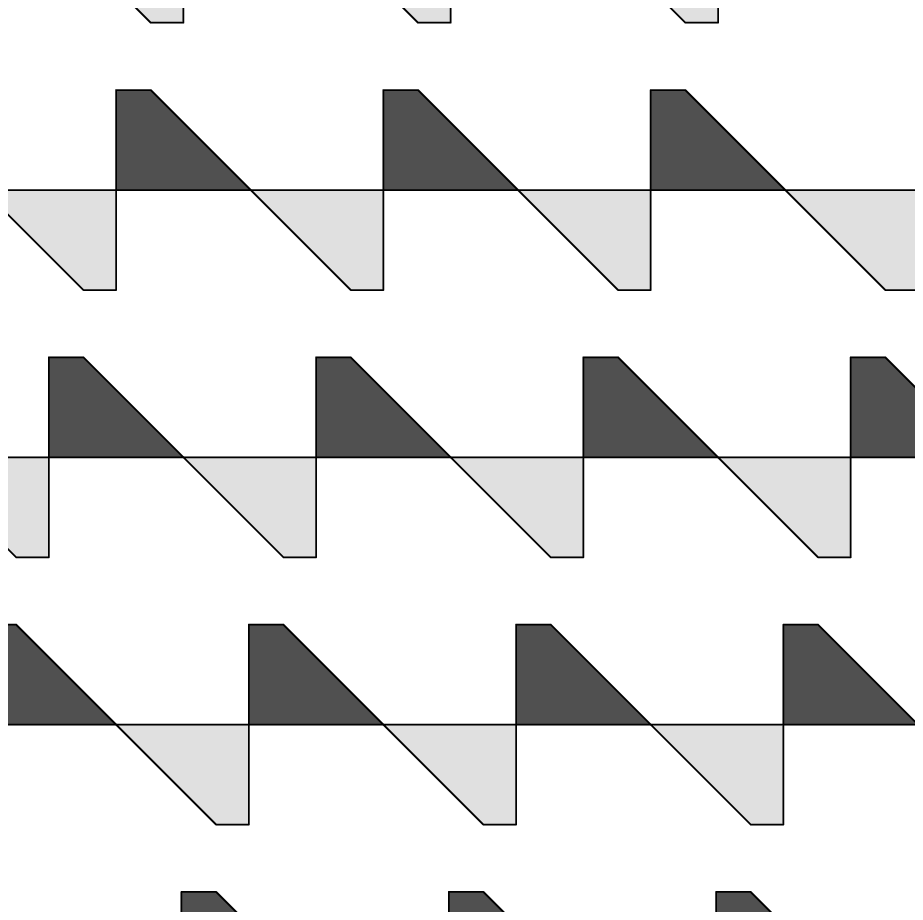}}
\newline
{\bf Figure 12.2:\/} 
A slice of $\Lambda R_+(1,1)$ and $\Lambda(R_-(1,1)+(1,1,0,0))$
\end{center} 

Figure 12.3 shows a
slice of $\Lambda R_+(1,0)$ and $\Lambda(R_-(1,0)+(1,1,0,0))$
over $B_0$.  The picture over $B_1$ is similar. 
Figure 12.4 shows a
slice of $\Lambda R_+(1,0)$ and $\Lambda(R_-(1,0)+(1,1,0,0))$
over $B_2$.  The picture over $B_3$ is similar.  

\begin{center}
\resizebox{!}{2.6in}{\includegraphics{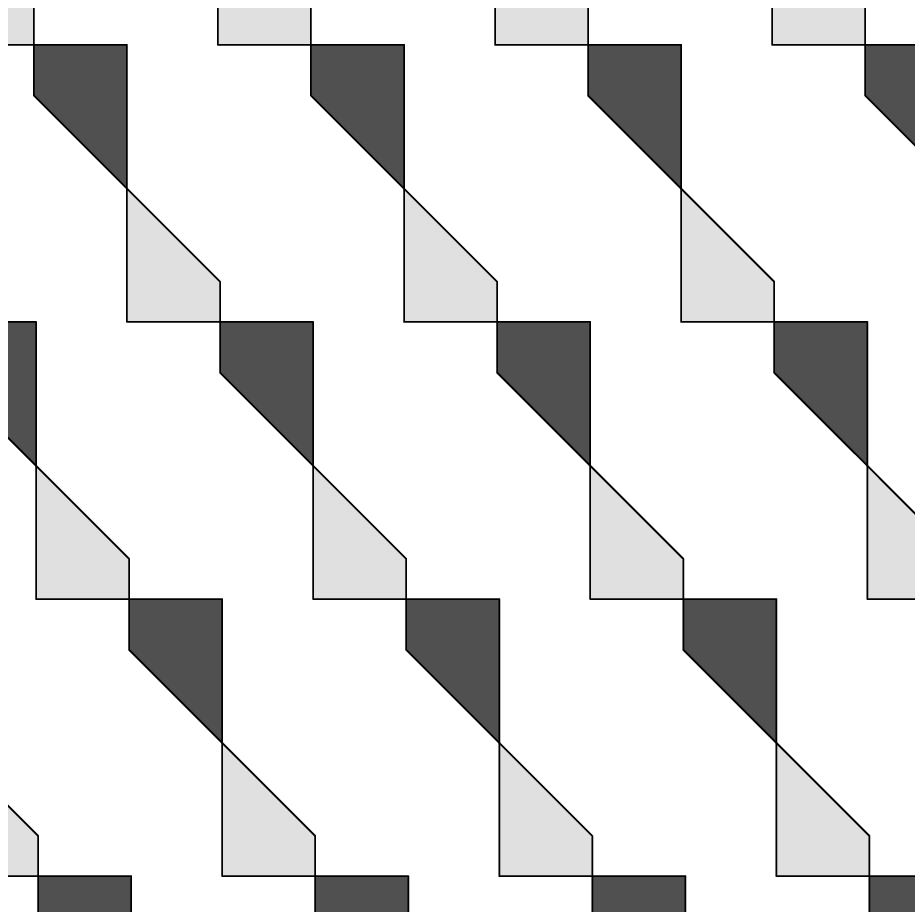}}
\newline
{\bf Figure 12.3:\/} 
A slice of $\Lambda R_+(1,1)$ and $\Lambda(R_-(1,-1)-(1,1,0,0))$.
\end{center} 

Figure 12.4 shows a
slice of $\Lambda R_+(1,0)$ and $\Lambda(R_-(1,0)+(1,1,0,0))$
over $B_2$.  The picture over $B_3$ is similar.  

\begin{center}
\resizebox{!}{2.6in}{\includegraphics{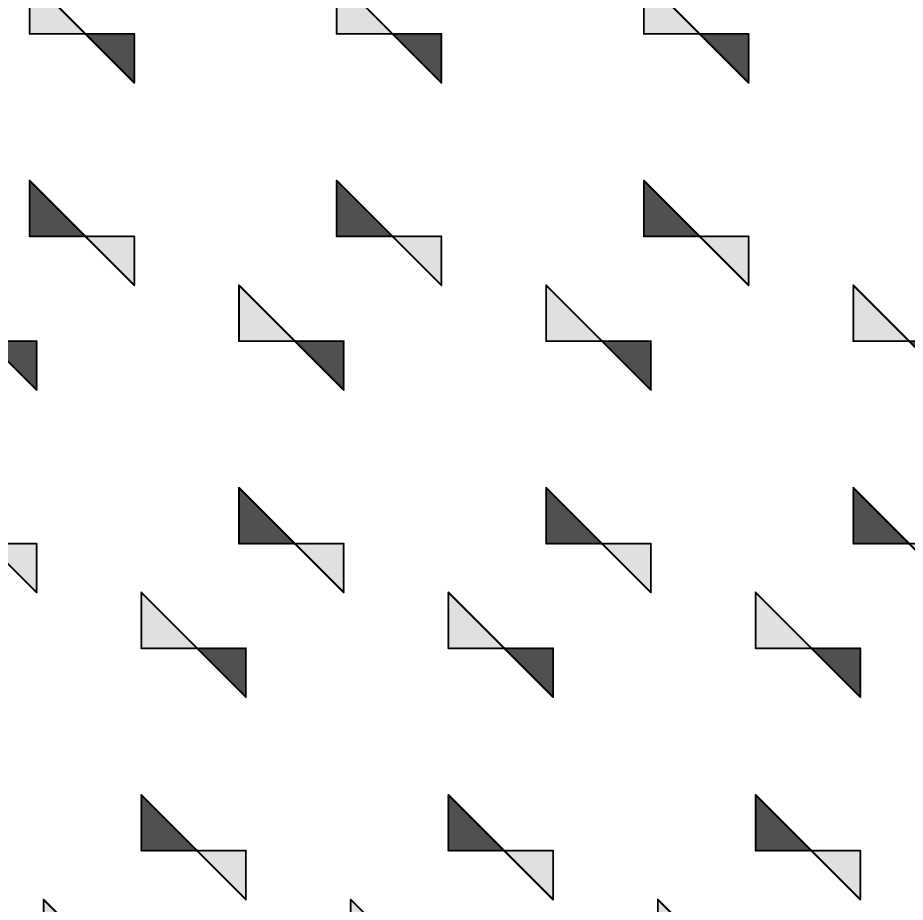}}
\newline
{\bf Figure 12.4:\/} 
A slice of $\Lambda R_+(1,1)$ and $\Lambda(R_-(1,-1)-(1,1,0,0))$.
\end{center} 

Figure 12.5 shows a slice of $\Lambda R_+(1,1)$ and $\Lambda(R_-(1,-1)-(1,1,0,0))$
over $B_2$.  The picture looks similar over $B_3$ and otherwise
at least one of the slices is empty.

\begin{center}
\resizebox{!}{2.6in}{\includegraphics{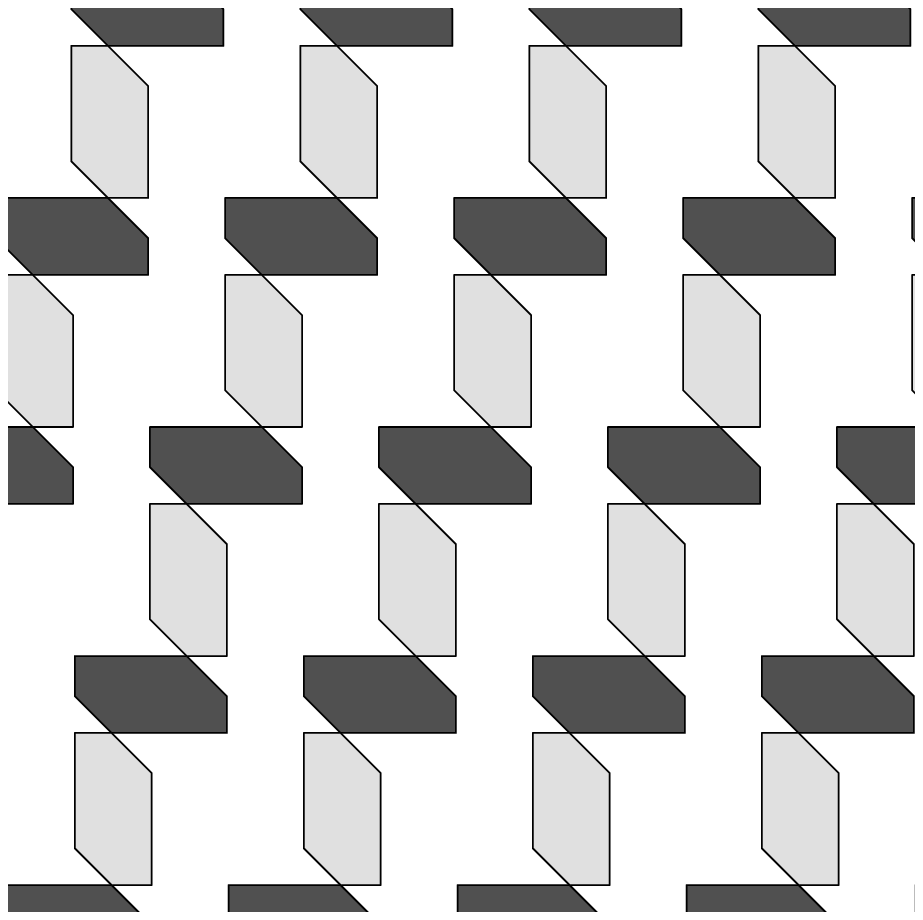}}
\newline
{\bf Figure 12.5:\/} 
A slice of $\Lambda R_+(1,1)$ and $\Lambda(R_-(1,-1)-(1,1,0,0))$.
\end{center} 

\newpage

\section{Extension and Symmetry}
\label{symm1proof}

\subsection{Translational Symmetry}
\label{periodicproof2}

Referring to \S \ref{mptnotation}, the maps $M_+$ and $M_-$ are
defined on all of $\Z^2$. This gives the extension of the
arithmetic graph to all of $\Z^2$. 

\begin{lemma}
The extended arithmetic graph does not cross the
baseline.
\end{lemma}

\startproof
By the Pinwheel Lemma, the arithmetic graph
describes the dynamics of the pinwheel map, $\Phi$.
Note that $\Phi$ is generically defined and invertible
on $\R_+ \times \{-1,1\}$.   Reflection in the
$x$-axis conjugates $\Phi$ to $\Phi^{-1}$.
By the Pinwheel Lemma, $\Phi$ maps
$\R_+ \times \{-1,1\}$ into itself.  By symmetry
the same goes for $\Phi^{-1}$.  Hence
$\Phi$ and $\Phi^{-1}$ also map
$\R_- \times \{-1,1\}$ into itself.
If some edge of $\widehat \Gamma$
crosses the baseline, then one of $\Phi$ or $\Phi^{-1}$
would map a point of $\R^+ \times \{-1,1\}$ into
$\R_- \times \{-1,1\}$.  This is a contradiction.
\endproof

Let $\lambda(p/q)=1$ if $p/q$ is odd,
and $\lambda(p/q)=2$ if $p/q$ is even.
Define

\begin{equation}
\label{generate}
\Theta=\Z V + \Z V';
\hskip 30 pt
V'=\lambda^2\bigg(0,\frac{(p+q)^2}{4}\bigg); \hskip 30 pt
\lambda=\lambda(p/q).
\end{equation}

\begin{lemma}
\label{invariant1}
The arithmetic graph $\widehat \Gamma(p/q)$ is invariant under $\Theta$.
\end{lemma}

\startproof
We will give the proof in the case when $p/q$ is odd.  The
even case is similar.
We have already seen that $\widehat \Gamma$ is invariant under $V$.
We just have to show invariance for $V'$.
By the Master Picture Theorem, it suffices
to prove that $(t,t,t) \in \Lambda$ when
$t$ is the second coordinate of $V'$.  Here
$\Lambda$ is as in Equation \ref{lattice}.

We have $(t,t,t) \equiv (2t,2t,0)$ mod $\Lambda$ 
because $t$ is an integer.  Setting
\begin{equation}
a=pq; \hskip 30 pt b=\frac{pq+q^2}{2},
\end{equation}
We compute that
\begin{equation}
a \left[\matrix{1+A\cr 0\cr 0}\right]+b  \left[\matrix{1-A\cr 1+A\cr 0}\right]=
 \left[\matrix{2t \cr 2t \cr 0}\right].
\end{equation}
This completes the proof.
\endproof

\begin{lemma}
the hexagrid is invariant under the action of $\Theta$.
\end{lemma}

\startproof
Again, we treat the odd case only.
Let $G=G(p/q)$ denote the hexagrid. 
As in the previous result, we just have to show
that $G$ is invariant under $V'$.
Let
$$W=\bigg(\frac{pq}{p+q},\frac{pq}{p+q}+\frac{q-p}{2}\bigg)$$
be the vector from the definition of the hexagrid $G$.
It suffices to prove that $6$ lines of $G$ contain
$V'$.   We compute that
\begin{equation}
V'=-\frac{p}{2}V + \frac{p+q}{2}W.
\end{equation}
The second coefficient is an integer.
Given that the room grid 
$RG$ is invariant under
the lattice $\Z[V/2,W]$, we see that
$RG$ is also invariant under
translation by $V'$.  This gives $2$ lines,
$L_1$ and $L_2$,
one from each family of $RG$.

Note that
$DG$ is only invariant under $\Z[V]$, so we have to work harder.
We need to produce $4$ lines of $DG$ that contain $V'$.
Here they are.
\begin{itemize}
\item The vertical line $L_3$ through the
origin certainly contains $V'$.  
This line extends the bottom left edge of $Q$
and hence belongs to $DG$.
\item Let $L_4$ be the line containing $V'$ and
point $-(p+q)V/2 \in \Z[V]$. 
We compute that the slope of $L_4$ coincides with the
slope of the top left edge of $Q$.  The origin
contains a line of $DG$ parallel to the top left edge of $Q$,
and hence every point in $\Z[V]$ contains such a line.
Hence $L_4$ belongs to $DG$.  To avoid a repetition
of words below, we call our argument here the {\it translation 
principle\/}.
\item Let $L_5$ be the line containing $V'$ and
point $-pV \in \Z[V]$. 
We compute that the slope of $L_5$ coincides with the
slope of the bottom right edge of $Q$.  The translation principle
shows that $L_5$ belongs to $DG$.
\item Let $L_6$ be the line containing $V'$ and
point $(q-p)V/2 \in \Z[V]$. 
We compute that the slope of $L_6$ coincides with the
slope of the top right edge of $Q$.  The translation
principle shows that
$L_6$ belongs to $DG$.
\end{itemize}
The reader can see these lines, for any desired parameter,
using Billiard King.
\endproof

\subsection{Rotational Symmetry}
\label{symmetry}

Let $p/q$ be an odd rational.
Let $p_+/q_+$ be as in Equation \ref{induct0}.
Let $\iota$ be the rotation
\begin{equation}
\label{iota5}
\iota(m,n)=V_+-(m,n).
\end{equation}
Here $V_+=(q_+,-p_+)$.
The fixed point of $\iota$ is $(1/2)V_+$.
This point lies very close to the baseline
of $\widehat \Gamma(p/q)$.
Figure 13.1 shows $\Gamma(7/17)$ centered on
this fixed point.

\begin{center}
\resizebox{!}{3.6in}{\includegraphics{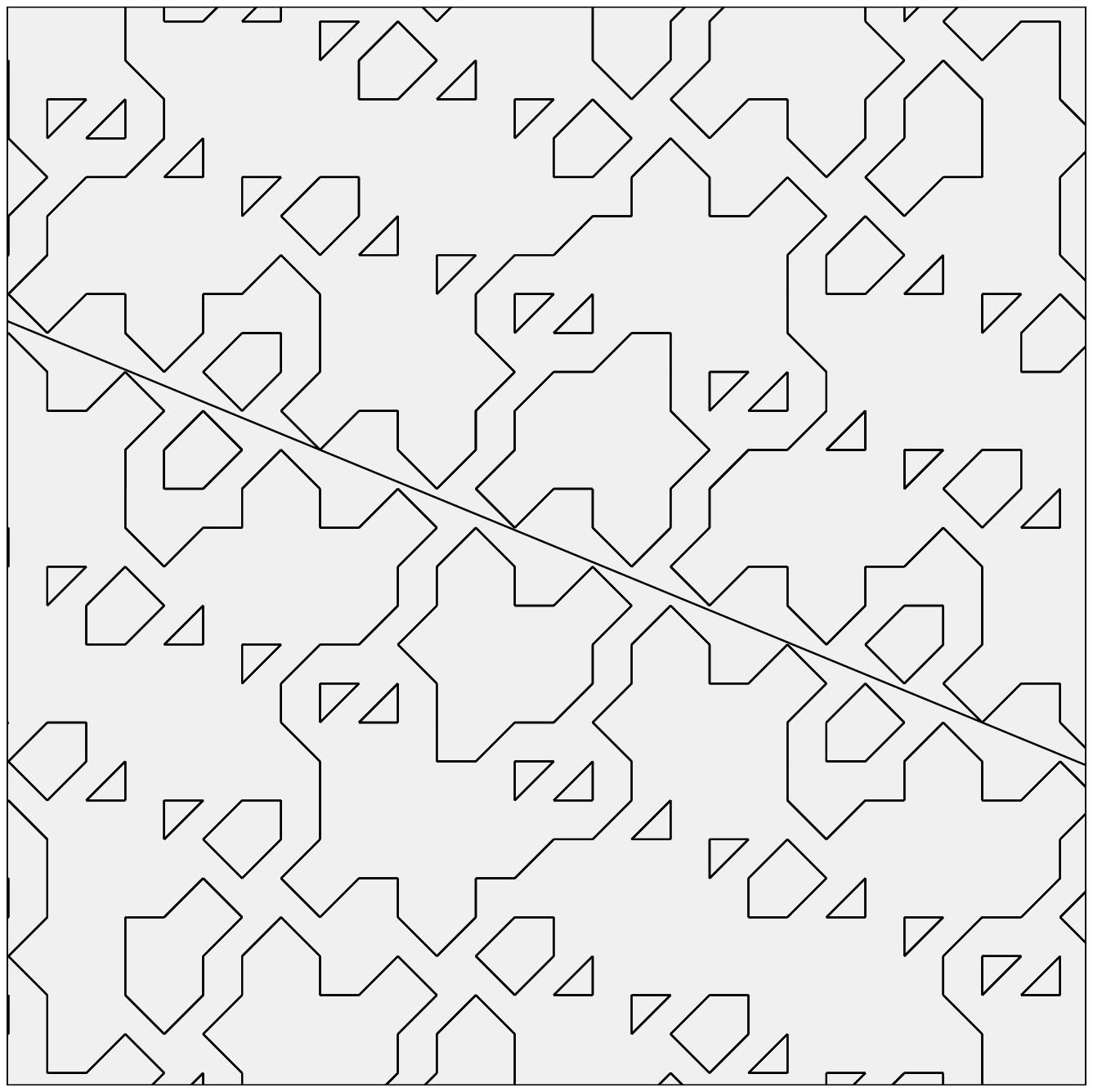}}
\newline
{\bf Figure 13.1:\/} $\widehat \Gamma(7/17)$ centered on 
the point $(12,-5)/2$.
point of symmetry.
\end{center} 

Below we prove that $\iota(\widetilde \Gamma)=\widehat \Gamma$,
as suggested by Figure 13.1. 
Combining this result with the translation symmetry above, we see
that rotation by $\pi$ about any of the points
\begin{equation}
\label{rotpoints}
\beta+\theta; \hskip 30 pt \beta=(1/2)V_+; \hskip 30 pt \theta \in \Theta
\end{equation}
is a symmetry of $\widehat \Gamma$.  
\newline
\newline
{\bf Remark:\/} In particular, there is an
involution swapping $(0,0)$ and $V_++dV$ for any $d \in \Z$.
\newline

\begin{lemma}
$\iota(\widetilde \Gamma)=\widetilde \Gamma$.
\end{lemma}

\startproof
Let $M_+$ and $M_-$ be as in \S \ref{mptnotation}.
As usual, we take $\alpha=1/(2q)$.
We will first compare $M_+(m,n)$ with
$M_-(\iota(m,n))$.
We have
\begin{equation}
M_+(m,n)=(t,t+1,t) \hskip 5 pt {\rm mod\/} \hskip 5 pt
\Lambda; \hskip 30 pt \frac{pm}{q}+n+\frac{1}{2q}
\end{equation}
Next, using the fact that $q_+p-p_+q=-1$, we have
$$
M_-(\iota(m,n))=(t'-1,t',t') \hskip 5 pt {\rm mod\/} \hskip 5 pt \Lambda;$$
$$
t'=\bigg(\frac{q_+p}{q}-p_+\bigg)- \bigg(\frac{pm}{q}+n\bigg)+\frac{1}{2q}=
-\bigg(\frac{pm}{q}+n\bigg)-\frac{1}{2q}=-t.$$
In short
\begin{equation}
M_-(\iota(m,n))=(-t-1,-t,-t) \hskip 5 pt {\rm mod\/} \hskip 5 pt \Lambda.
\end{equation}

Recall that $R_A$ is the fundamental domain for the action
of $\Lambda=\Lambda_A$.  We mean to equate $\Lambda$ with the
$\Z$ span of its columns.  There is some
$v \in \Lambda$ such that
\begin{equation}
(s_1,s_2,s_3)=(t,t+1,t)+(v_1,v_2,v_3) \in R_A
\end{equation}
Given Equation \ref{lattice}, we have
$(2+A,A,1) \in \Lambda$.  Hence
\begin{equation}
w=(-v_1+2+A,-v_2+A,-v_3+1) \in \Lambda.
\end{equation}
We compute that
\begin{equation}
(-t-1,-t,-t)+w=(1+A,1+A,1)-(s_1,s_2,s_3).
\end{equation}
So, we have
\begin{equation}
\label{reflect1}
M_+(m,n)=\rho \circ M_-(\iota(m,n)),
\end{equation}
where $\rho$ is reflection through the midpoint of the space $R_A$.
Similarly,
\begin{equation}
M_-(m,n)=\rho \circ M_+(\iota(m,n)),
\end{equation}

But now we just verify by inspection that our partition of
$R_A$ is symmetric under $\rho$, and has the labels
appropriate to force the type determined by
$$\rho \circ M_+(m,n), \hskip 5 pt
\rho \circ M_-(m,n)$$
to be the $180$ degree rotation of the type forced by
$$M_-(m,n), \hskip 5 pt
M_+(m,n)$$
Indeed, we can determine this with an experiment performed on
any rational large enough such that all regions are
sampled.
\endproof

\subsection{Near Bilateral Symmetry}
\label{nbs}

Our pictures of arithmetic graphs show that they
have an approximate bilateral symmetry.  For example,
in Figure 4.2 the two arcs
$\Gamma \cap R_1$ and $\Gamma \cap R_2$ both have
near bilateral symmetry.  In Figure 13.1 we see a
similar phenomenon.  Here we will explain this
near-symmetry. 

We say that a map $\J$ from $\widehat \Gamma$ to
$\widehat \Gamma$ is a {\it combinatorial
isomorphism\/} if $\J$ maps vertices to
vertices and edges to edges.  We say that
$\J$ is {\it pseudo-linear\/} if there
is a linear isomorphism $J: \R^2 \to \R^2$
such that $J$ is a bounded distance from
$\J$ (in the sup norm.)  In this case,
we call $J$ the {\it model\/} for $\J$.
Here is our main result.

\begin{lemma}
\label{pseudolinear}
For any irrational $A$, there exists an involution
$\J: \widehat \Gamma \to \widehat \Gamma$
with the following properties.
\begin{enumerate}
\item $\J$ is a combinatorial isomorphism that
swaps the components of $\widehat \Gamma$ above
the baseline with the ones below.
\item $\J$ is a translation when restricted to
low vertices.  More precisely,
if $v$ is a low vertex then $\J(v)=v-(0,1)$.
\item $\J$ is pseudo-linear, modelled on the
affine map $J$ such that $J(V)=V$ and
$J(W)=-W$.  Here $V$ and $W$ are as in
Equation \ref{boxvectors}.
\end{enumerate}
\end{lemma}

\noindent
{\bf Remarks:\/} \newline
(i) We think that $\J$ is within $2$ units of
$J$.  Probably an analysis similar to what we
did for the Pinwheel Lemma would prove this.
\newline
(ii) Let $\iota$ be the symmetry discussed in the
previous section.  Then $\iota \circ \J$ permutes
the components of $\widehat \Gamma$ above the
baseline. In particular, $\iota \circ \J$ 
preserves $\Gamma$ but reverses its direction.
This is the near-bilateral symmetry that we see
in the pictures.  \newline
(iii) One can probably see the action of $\J$
by looking at Figure 13.1 again.  Notice the
symmetry between components above the baseline
and components below it.
\newline

The construction of $\J$ is almost completely
soft. It only uses the easy case of the
Pinwheel Lemma, and basic symmetries of
the outer billiards map.  
Our construction uses the following definition.
Say that a {\it low component\/} is a component
of $\widehat \Gamma$ above the baseline
that contains a low vertex.
\newline
\newline
{\bf Constructing the Involution:\/}
We turn now to the construction of $\J$.
Recall that $\Xi=\R_+ \times \{-1,1\}$.  Let
$\Xi_{\pm}=\R_{\pm} \times \{-1,1\}$.  Then
$\Xi=\Xi_+$.   Recall that $\Psi: \Xi_+ \to \Xi_+$
is the first return map.  We can extend
$\Psi$ to that it is also the return map
from $\Xi_-$ to $\Xi_-$.   The points of the
arithmetic graph below the baseline correspond
to this extended notion of $\Psi$.

Let $\Psi^{1/2}$ denote the first return map to $\R \times \{-1,1\}$.
Then $\Psi$ is the square of $\Psi^{1/2}$.  In terms of
the Pinwheel Lemma, we start with (say) a point in
$\Xi_+$ and then watch it wind halfway around the kite
until it lands in $\Xi_-$.  This is the point $\Psi^{1/2}(\xi)$.
The correspondence $\xi \to \Psi^{1/2}(\xi)$ gives a bijection
between $\Psi$-orbits in $\Xi_+$ and $\Psi$-orbits in $\Xi_-$.

In terms of the arithmetic graph, there is a combinatorial
isomorphism $\J_+$ of $\widehat \Gamma$ that 
swaps the components above the
baseline with the ones below it.  Here
$\J_+(m,n)=(m',n')$, where $(m,n)$ corresponds to
$\xi$ and $(m',n')$ corresponds to $\Psi^{1/2}(\xi)$.
There
is a second involution that is equally good.
We used the forwards direction of $\Psi$ to define
$\J_+$, but we could have used the backwards
direction.  That is, we would match $\xi \in \Xi_+$
to the point $\Psi^{-1/2}(\xi) \in \Xi_-$.
Call this map $\J_-$.

Our map $\J$ is made from $\J_+$ and
$\J_-$ in a not-completely-canonical way.
We will define $\J$ on components above
the baseline.  We then define $\J$ for components
below the baseline so as
to make $\J$ an involution.

Recall that the {\it parity\/} of a low vertex $(m,n)$ to be
the parity of $m+n$.  By Lemma \ref{parity}, 
a low component only has vertices of one
kind of parity.  We call the low component
{\it even\/} or {\it odd\/} depending on the parity
of its low vertices.
If $\gamma$ is a component
of $\widehat \Gamma$ above the baseline that is
not low, we use (say) $\J=\J_+$.  (We don't care
about these components.)  For even low components
we use $\J=\J_+$. For odd low components, we use
$\J=\J_-$.  From the discussion above, we see that
$\J$ is a graph isomorphism of $\widehat \Gamma$.
\newline
\newline
{\bf Remark:\/} There might be a canonical choice of
$\J_-$ or $\J_+$ for components that are not low,
but we don't know this.
\newline
\newline
{\bf Action on Low Vertices:\/}
Let's see what happens to low vertices. 
Let $(m,n)$ be an even low vertex and
let $(x,-1)=M(m,n)$.  We compute easily
that 
\begin{equation}
\Psi^{1/2}(x,-1)=\psi^2(x,-1)=(x-2,1)=M(m,n-1).
\end{equation}
Hence, $\J(m,n)=(m,n-1)$.
Similarly, if $(m,n)$ has odd parity, then
\begin{equation}
\Psi^{-1/2}(x,1)=\psi^{-2}(x,1)=(x-2,-1)=M(m,n-1).
\end{equation}
Hence $\J(v)=v-(0,1)$ when $v$ is a low vertex.
\newline
\newline
{\bf Pseudo-Linearity:\/}
It remains to show that $\J$ is pseudo-linear, modelled
on $J$.  Since we don't need this final result for any
purpose, we will only sketch the argument. Let
$(x,1)$ be a point 
on $\Xi_+$ about $N$ units from the origin,
we roughly trace out the Pinwheel map.  First we
some integer multiple of the vector $(0,4)$, then
we add some integer multiple of the vector $(-2,2)$,
etc.  When we reach $\Xi_-$ we have a vector of
the form 
$$(x+2Am_1+2n_1,\pm 1),$$
where the pair $(n_1,m_1)$ depends linearly
on $N$, up to a uniformly bounded error.
But, for the corresponding point $v \in \widehat \Gamma$,
we have $\J(v)=v+(m_1,n_1)$.  This shows that
$\J$ is pseudo-linear, and a simple calculation
shows that $\J$ is modelled on $J$.

\newpage

\section{The Structure of the Doors}
\label{odddoor}

\subsection{The Odd Case}

We suppose that $p/q$ is an odd rational.
Say that a {\it wall line\/} is a line of positive slope in
the room grid.  The doors are the intersection points
of lines in the door grid with the wall lines.
Let $L_0$ be the wall line through $(0,0)$.  Let
$L_1$ be the wall line through $V/2$.   

\begin{lemma}
Any two wall lines are equivalent mod $\Theta$.
\end{lemma}

\startproof
We check explicitly that the vector
$$V'+\frac{p+1}{2}V \in \Theta \cap L_1$$
Hence $L_0$ and $L_1$ are equivalent mod $\Theta$.
But any other wall line is obtained from one of
$L_0$ or $L_1$ by adding a suitable integer
multiple of $V$.
\endproof

\begin{lemma}
\label{doorint}
The first coordinate of any door is an integer.
\end{lemma}

\startproof
Any wall line is equivalent mod $\Theta$ to $L_0$.
Since $\Theta$ acts by integer translations,
it suffices to 
door lies on $L_0$.  Such a door is an integer multiple of
the point $v_3$ in Figure 3.1.  That is, our door has
coordinates
\begin{equation}
\label{dooreq1}
\frac{k}{2q}(2pq,(p+q)^2-2p^2).
\end{equation}
The first coordinate here is certainly an integer.
\endproof

It could happen that the second coordinate of a door is
an integer.  Call such a door exceptional.

\begin{lemma}
Modulo the action of $\Theta$, there are only
two exceptional doors.
\end{lemma}

\startproof
The point $(0,0)$ gives rise to two exceptional doors with (the same)
integer coordinates.  One of these doors is associated
to the wall above $(0,0)$, and one of these doors is
associated to the door below.  Hence, it suffices to
show that any door with integer coordinates lies in
$\Theta$.

As in the preceding result, it suffices to consider
doors on $L_0$.   Given Equation \ref{dooreq1},
we see that
$$k\frac{(p+q)^2-2p^2}{2q} \in \Z$$
for an exceptional door.
Expanding this out, and observing that $q$ divides both $q^2$ and $pq$, we get that
$$k\frac{q^2-p^2}{2q} \in \Z.$$
But $q$ and $q^2-p^2$ are relatively prime.  Hence $k=jq$ for some $j \in \Z$.
But $$qv_3=2V'+V \in \Theta.$$
Hence $jqv_3=kv_3 \in \Theta$ as well.
\endproof

Here is a related result.

\begin{lemma}
\label{specialcross}
Any lattice point on a wall line is equivalent
to $(0,0)$ mod $\Theta$.
\end{lemma}

\startproof
By symmetry, it suffices
to consider the cases when $(m,n) \in L_0$.

Looking at Figure 3.1, we see that any point on $L_0$ has the form
\begin{equation}
\label{onL0}
sv_5=\frac{s}{2(p+q)}(2pq,(p+q)^2-2p^2).
\end{equation}
In order for this point to lie in $\Z^2$, the first coordinate
must be an integer.  Since $p$ and $q$ are relatively
prime, $pq$ and $p+q$ are relatively prime.  Hence,
the first coordinate is an integer only if $s=k(p+q)$ for
some $k \in \Z$.   Hence $(m,n)$ is an integer multiple of the
point
$$(p+q)v_5=\bigg(pq,\frac{(p+q)^2}{2}-p^2\bigg)=2V'+pV \in \Theta.$$
Here $V$ and $V'$ are the vectors generating $\Theta$, as in
Equation \ref{generate}. 
\endproof

The vertical lines in the door grid have the
form $x=kq$ for $k \in \Z$.  Say that a
{\it Type 1 door\/} is the intersection 
of such a line with a wall line.

\begin{lemma}
\label{door1}
Let $(kq,y)$ be a Type 1 door.  Then $py \in \Z$.
\end{lemma}

\startproof
The group $\Theta$ acts transitively, by integer translations, on
the vertical lines of the door grid.  Hence, is suffices to
prove this lemma for the case $k=0$.  In other words, we
need to show that $py \in \Z$ if $(0,y)$ lies on a wall line.

We order the wall lines according to the order in which they
intersect the line of slope $-A=-p/q$ through the origin.
Let $y_n$ be such that $(0,y_n)$ lies on the $k$th wall line.
The sequence $\{y_n\}$ is an arithmetic progression.  Hence,
it suffices to prove our result for two consecutive values of $n$.
Note that $(0,0)$ is a type A door.  We might as well normalize
so that $y_0=0$.  Then $(0,y_1)$ lies on the wall line $L_1$
through $(-q,p)$.   Referring to Equation \ref{boxvector},
two points on $L_1$ are $-V$ and $-V+W$.  These points
are given by
$$
-V=(-q,p); \hskip 30 pt 
-VW=(-q,p)+\bigg(\frac{pq}{p+q},\frac{pq}{p+q}+\frac{q-p}{2}\bigg).
$$
From this information, we compute that
$y_1=(p+q)^2/2p$.  Since $p+q$ is even,
$py_1=(p+q)^2/2 \in \Z$.
\endproof

Recall that $\underline y$ is the greatest integer
less than $y$.

\begin{corollary}
\label{nohalf}
Suppose that $(x,y)$ is a door of type $1$, then
$y-\underline y \not = 1/2$.
\end{corollary}

\startproof
$p(y-\underline y)=p/2$ is an integer, by the previous
result.  But $p/2$ is not an integer.  This is a
contradiction.
\endproof

Say that a {\it Type 2\/} door
is a door on $L_0$ that is not of Type 1.
One obtains a Type 2 door by intersecting
$L_0$ with a line of the door grid that is
parallel to the top left (or right) 
edge of the arithmetic kite. 

\begin{lemma}
\label{door2}
The Type 2 
doors are precisely the points on $L_0$ of the form
$(kp,y_k)$, where $k \in \Z$ and $y_k$ is a number that
depends on $k$.
\end{lemma}

\startproof
Referring to Figure 3.1,
two consecutive doors on $L_0$ are $(0,0)$ and $v_3=(p,y_1)$.
Our lemma now follows from the fact that the sequence of
doors on $L_0$ forms an arithmetic progression.
\endproof

\subsection{The Even Case}

Now we revisit all the results above in case
$p/q$ is even.

\begin{lemma}
Any two wall lines are equivalent mod $\Theta$.
\end{lemma}

\startproof
This is easy in the even case.
Translation by $V$ maps each
wall line to the adjacent one.
\endproof

\begin{lemma}
The first coordinate of any door is an integer.
\end{lemma}

\startproof
The first door on $L_0$ is the same in the
even case as in the odd case.  The rest
of the proof is the same as in the odd case.
\endproof

\begin{lemma}
Modulo the action of $\Theta$, there are only
two exceptional doors.
\end{lemma}

\startproof
As in the odd case, we just have to show that any
door with integer coordinates is equivalent to
$(0,0)$ mod $\Theta$.   As in the odd case,
the doors on $L_0$ have the form $kv_3$.  As in the
odd case, this
leads to the statement that
$$k\frac{q^2-p^2}{2q} \in \Z.$$
Now the proof is a bit different.
Here $2q$ and $q^2-p^2$ are relatively prime.  
Hence $k=2jq$ for some $j \in \Z$.
But $$2qv_3=V'+2V \in \Theta.$$
Hence $2jqv_3=kv_3 \in \Theta$ as well.
\endproof

\begin{lemma}
Any lattice point on a wall line is equivalent
to $(0,0)$ mod $\Theta$.
\end{lemma}

\startproof
As in the proof of Lemma \ref{specialcross}, we see that
\begin{equation}
sv_5=\frac{s}{2(p+q)}(2pq,(p+q)^2-2p^2) \in \Z.
\end{equation}
As in the odd case, we look at the first coordinate and
deduce the fact that
$s=k(p+q)$ for some $k \in \Z$.
This is not enough for us in the even case.
Looking now at the second coordinate, we see that
$$\frac{k(p+q)}{2} - kp^2 \in \Z.$$
Hence $k$ is even.
Hence $(m,n)$ is an integer multiple of the
point
$$2(p+q)v_5=(2pq,(p+q)^2-2p^2)=V'+2pV \in \Theta.$$
\endproof

We don't repeat the proof of Lemma \ref{door1} because
we don't need it. We only need the even
version of Corollary \ref{nohalf}.
In the even case, we have simply forced
Corollary \ref{nohalf} to be true by eliminating
the crossings for which it fails.

We say that a {\it Type 2\/} door is a door
on $L_0$ that is not of Type 1, and also
is not one of the crossings we have eliminated.
Once we make this redefinition, we have
the following result

\begin{lemma}
\label{door22}
The Type 2 
doors are precisely the points on $L_0$ of the form
$(kp,y_k)$, where $k \in \Z$ is not an odd multiple
of $q$, and $y_k$ is a number that
depends on $k$.
\end{lemma}

\startproof
The Type two doors are as in the odd case,
except that we eliminate the points $(kp,y_k)$
where $k$ is an odd multiple of $q$.
\endproof

\newpage

\section{Proof of the Hexagrid Theorem I}
\label{hex1proof}

\subsection{The Key Result}

We will assume that $p/q$ is an odd rational until
the end of the chapter.

Say that a {\it floor line\/} is a negatively sloped
line of the floor grid. 
Say that a {\it floor point\/} is a point on a floor
line.  Such a point need not have integer coordinates.
Let $M_+$ and $M_-$ denote the maps from Equation \ref{mptnotation}.

\begin{lemma}
\label{oddkey1}
If $(m,n)$ is a floor point, 
then $M_-(m,n)$ is equivalent mod
$\Lambda$ to a point of the form
$(\beta,0,0)$.
\end{lemma}

\startproof
The map $M_-$ is constant when restricted to
each floor line, because these lines have slope
$-A$.  Hence, it suffices to prove this result for
one point on each floor line. 
The points
\begin{equation}
\label{hex1integral}
(0,t); \hskip 30 pt t=
\frac{k(p+q)}{2}; \hskip 30 pt k \in \Z.
\end{equation}
form a sequence of floor points, one per floor line.
Note that $t$ is an integer, because $p+q$ is even.

To compute the image of the point $(0,t)$, we just
have to subject the point $t$ to our reduction algorithm
from \S \ref{master1}.  
The first $4$ steps of the algorithm lead to the following
result.
\begin{enumerate}
\item $z=t.$
\item $Z={\rm floor\/}(t)=t$, because $t$ is an integer.
\item $y=2t=k(p+q)=kq(1+A)$.
\item $Y={\rm floor\/}(y/(1+A))=kq$.
\end{enumerate}
Hence $z=Z$ and $y=(1+A)Y$.  Hence
\begin{equation}
M_-(0,t)=(x-(1+A)X,y-(1+A)Y,z-Z)=(\beta,0,0),
\end{equation}
for some number $\beta \in \R$ that depends on $A$ and $k$.
\endproof

\subsection{Two Special Planes}

Let $\Pi_- \subset \R^3$ denote the plane
given by $y=z$.   We can think of
$\Pi_-$ as the plane through the origin
generated by the vectors $(1,0,0)$ and
$(1,1,1)$.   In particular, the
vector $(1,1,1)$ is contained in $\Pi_-$.
Let $\Pi_-(0)$ denote the line through the
origin parallel to $(1,0,0)$.  Then
$\Pi_-(0)$ is a line in $\Pi_-$.
Define
\begin{equation}
\Pi_+=\Pi_-+(1,1,0); \hskip 30 pt
\Pi_+(0)=\Pi_-(0)+(1,1,0).
\end{equation}

\begin{lemma}
\label{oddkey2}
If $(m,n)$ is a floor point, then
$M_{\pm}(m,n)$ is equivalent mod
$\Lambda$ to a point in
$\Pi_{\pm}(0)$.
\end{lemma}

\startproof
The $(-)$ case of this result is just a restatement
of Lemma \ref{oddkey1}. The $(+)$ case follows
from the $(-)$ case and symmetry.  That is,
we just translate the $(-)$ case by
the vector $(1,1,0)$ to get the $(+)$ case.
\endproof

Define
\begin{equation}
\label{oddkey4}
\Pi_{\pm}(r)=\Pi_{\pm}(0) + (r,r,r).
\end{equation}
Let $\Pi_{\pm}(r,s)$ denote the open infinite strip
that is bounded by $\Pi_{\pm}(r)$ and $\Pi_{\pm}(s)$.
In the case of interest to us, we will have $r=0$ and $s=\lambda>0$.

For each pair $(\epsilon_1,\epsilon_2) \in \{-1,0,1\}^2$, let
$\Sigma(\epsilon_1,\epsilon_2)$ denote the set of
lattice points $(m,n)$ such that $(m,n)$ and
$(m,n)+(\epsilon_1,\epsilon_2)$ are separated by some
floor line.  The set $\Sigma(\epsilon_1,\epsilon_2)$ is
obtained by intersecting $\Z^2$ with an infinite
union of evenly spaced infinite strips, each of which
has a floor line as one boundary component.
For our purposes, it suffices to consider the pairs
\begin{equation}
(-1,0); \hskip 20 pt (-1,-1); \hskip 20 pt (0,-1); \hskip 20 pt (1,-1).
\end{equation}
For these pairs, the floor lines are the lower boundaries of
the strips.  We define
\begin{equation}
\label{oddkey3}
\lambda(\epsilon_1,\epsilon_2)=-(A\epsilon_1+\epsilon_2).
\end{equation}

\begin{lemma}
\label{oddkey5}
Let $\lambda=\lambda(\epsilon_1,\epsilon_2)$.
Suppose that 
$(m,n) \in \Sigma(\epsilon_1,\epsilon_2)$.
Then $M_{\pm}(m,n) \in \Pi_{\pm}(0,\lambda)$.
\end{lemma}

\startproof
We consider the case of $M_-$ and the pair $(-1,0)$.  The other
cases have essentially the same proof.  If $(m,n) \in \Sigma(-1,0)$,
Then there is some $x \in (m-1,m)$ such that $(x,n)$ is a floor point.
Then $M_+(x,n)$ is $\Lambda$-equivalent to a point $p$ in $\Pi(0)$. But
then $M_+(m,n)$ is $\Lambda$-equivalent to $p+(m-x)(A,A,A) \in \Pi(0,A)=\Pi(0,\lambda(-1,0))$.
\endproof

\subsection{Critical Points}

Say that a point $v \in \Sigma(\epsilon_1,\epsilon_2)$ is {\it critical for\/}
$(\epsilon_1,\epsilon_2)$ if the arithmetic graph contains the edge
joining $(m,n)$ to $(m+\epsilon_1,n+\epsilon_2)$.   Statement 1 of the
Hexagrid Theorem says, in particular, that there are no such points
like this.

\begin{lemma}
There are no critical points.
\end{lemma}

\startproof
Let ${\cal R\/}_+$ denote the tiling of $\R^3$ by polyhedra, according to
the Master Picture Theorem.  Let ${\cal P\/}_+$ denote the intersection
of ${\cal R\/}_+$ with the plane $\Pi$.  We make the same definitions
in the $(-)$ case.  If
$(m,n)$ is critical for $(\epsilon_1,\epsilon_2)$, then
one of two things is true.
\begin{enumerate}
\item $\Pi_+(0,\lambda)$ nontrivially intersects
a polygon of ${\cal P\/}_+$ labelled by $(\epsilon_1,\epsilon_2)$.
\item $\Pi_-(0,\lambda)$ nontrivially intersects
a polygon of ${\cal P\/}_-$ labelled by $(\epsilon_1,\epsilon_2)$.
\end{enumerate}
Here we have set $\lambda=\lambda(\epsilon_1,\epsilon_2)$.  Considering
the $4$ pairs of interest to us, and the $2$ possible signs, we
have $8$ conditions to rule out.  We check, in all cases,
that the relevant strip is disjoint from the relevant
polygons.

We can check the disjointness for all parameters at once.
The union 
$$S_{\pm}(\epsilon_1,\epsilon_2):=\bigcup_{A \in (0,1)}\bigg( \Pi_{\pm}(\epsilon_1,\epsilon_2;A) \times \{A\}\bigg)$$
is a polyhedral subset of $\R^4$.    To get an honest polyhedron,
we observe that $S$ is invariant the action of the lattice
element $\gamma_1$ from Equation \ref{affine}, and we take
a polyhedron whose union under translates by $\gamma_3$
tiles $S$.  In practice, we simply restrict the $x$-coordinate
to lie in $[0,2]$.

We check that $S_{\pm}(\epsilon_1,\epsilon_2)$, or rather
the compact polyhedron replacing it, is disjoint
from all $\Lambda$-translates of the polytope $P_{\pm}(\epsilon_1,\epsilon_2)$,
the polytope listed in \S \ref{master2}.  In practice,
most translates are very far away, and we only need to check
a small finite list.  This
is a purely algebraic calculation.
\endproof

Rather than
dwell on the disjointness calculation, which gives no insight
into what is going on, we will draw pictures
for the parameter $A=1/3$.
The combinatorial type changes with the parameter, but not the
basic features of interest to us.  The interested reader
can see the pictures for any parameter using Billiard
King.

To draw pictures, we identify the planes $\Pi_{\pm}$ with
$\R^2$ using the projection $(x,y,z) \to (x,(y+z)/2)$.  Under
this identification, all the polygons in question are
rectangles!  The coordinates of the rectangle vertices are
small rational combinations of $1$ and $A$, and can
easily be determined by inspection. The whole
picture is invariant under translation by $(1+A,0)$.
The thick line in the first picture 
corresponds to $\Pi_-(0)$.  In terms of
$\R^2$ coordinate, this is the $x$-axis.
The black dot is $(0,0)$.
dot is $(4/3,0)=(1+A,0)$.

We explain by example the notation on the right hand
side of the fiture.  The label $\lambda(-1,-1)$
denotes the line $\Pi(\lambda)$, where
$\lambda=\lambda(-1,-1)$.
 In each case, the relevant strip lies below
the relevant shaded piece.  While the combinatorics
of the picture changes as the parameter changes,
the basic disjointness stays the same. 

\begin{center}
\resizebox{!}{3.6in}{\includegraphics{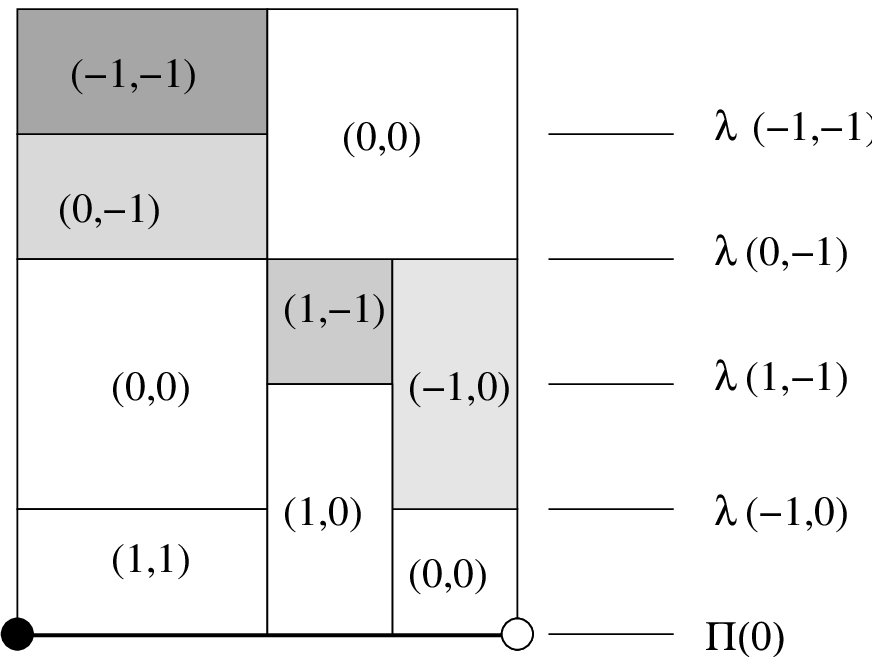}}
\newline
{\bf Figure 15.1:\/} The $(-)$ picture for $A=1/3$.
\end{center} 

Figure 15.2 shows the same thing for the $(+)$ case.
This time the black dot is $(1/2,1/2)$ and the white
dot is $(1/2,1/2)+(1+A,0)$.   The thick line represents
$\Pi_+(0)$. In $\R^2$ coordinates, this is the line $y=1/2$.
In the $(+)$ case is isn't even a close call.

\begin{center}
\resizebox{!}{3.8in}{\includegraphics{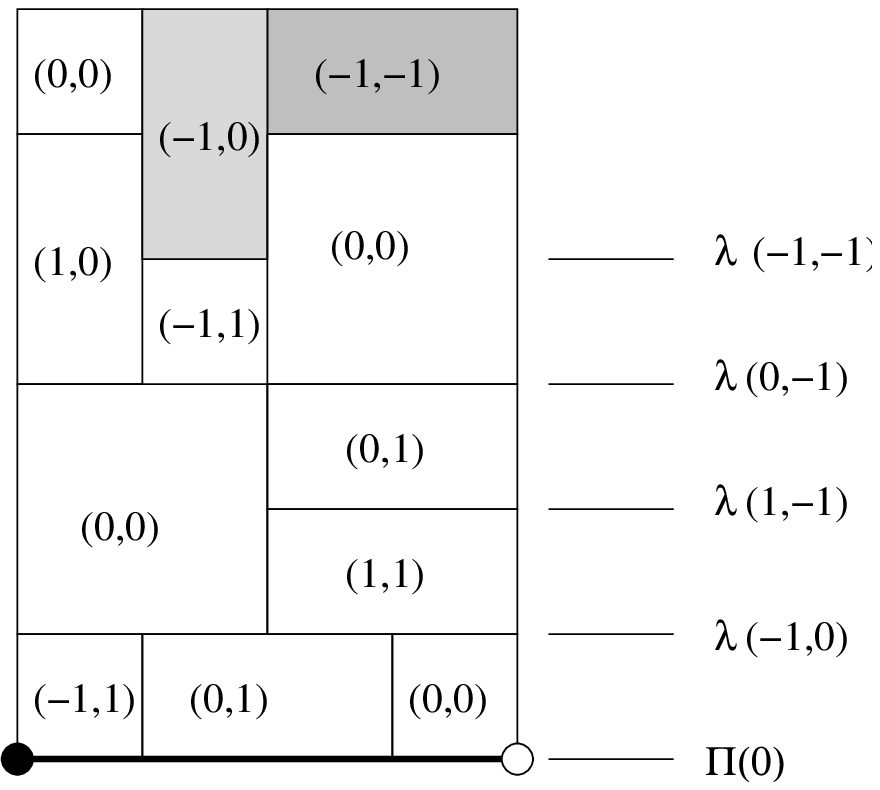}}
\newline
{\bf Figure 15.2:\/} The $(+)$ picture for $A=1/3$.
\end{center} 

\subsection{The End of the Proof}

Now we know that there are no critical points.  The only other
way that the arithmetic graph could cross a floor line would
be at a floor point that was also a lattice point.
It might happen that one edge emanating from
such a floor point lies above the floor line, and the
other lies below.  

Define
\begin{equation}
\label{floorpoint}
\zeta_k=\bigg(0,\frac{k(p+q)}{2}\bigg); \hskip 30 pt k \in \Z.
\end{equation}

\begin{lemma}
Modulo the symmetry group $\Theta$, the
only lattice floor points are the ones listed in Equation \ref{floorpoint}. 
\end{lemma}

\startproof
If $(m,n)$ is a lattice floor point, then 
$2Am+2n \in \Z$.   But his means that $q$ divides $m$.
Subtracting off a suitable multiple of $V=(q,-p) \in \Theta$,
we can arrange that the first coordinate of our
lattice floor point is $0$.  But, now we must have
one of the points in Equation \ref{floorpoint}.
\endproof

The slices as shown in Figure 6.3 determine the nature of the
edges of the arithmetic graph, although the slices currently
of interest to us are not shown there.  We are interested
in following the method discussed in \S \ref{singular},
where we set $\alpha=0$ and consider the singular situation.
The points $M_-(\zeta_k)$ and $M_+(\zeta_k)$ both lie in the $(0,A)$ slices
of our partitions.  Figure 15.1 does for these slices what Figure 6.3 does for
the generic slice.  The point $M_-(\zeta_k)$ always lies along the bottom
edge of the fiber, and the point $M_+(\zeta_k)$ just above the
edge contained in the line $y=1$.  The relevant edges are highlighted.

\begin{center}
\resizebox{!}{3.8in}{\includegraphics{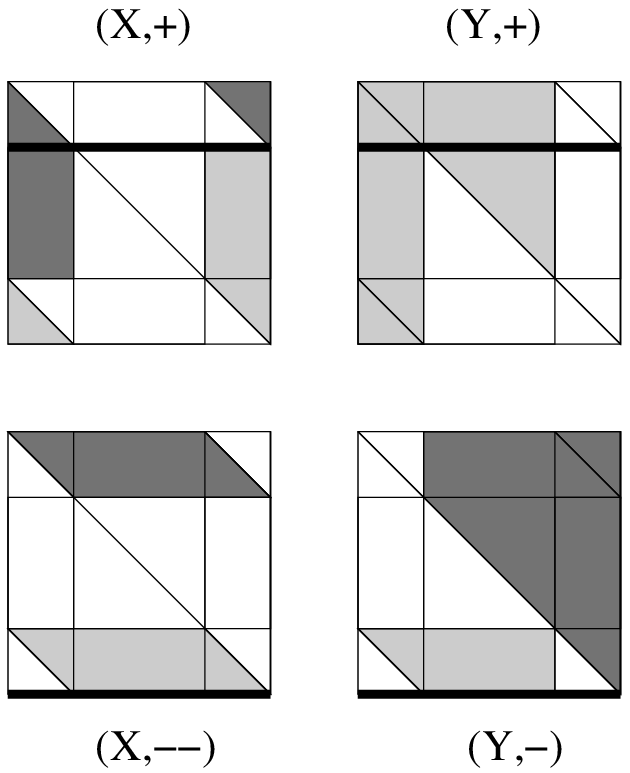}}
\newline
{\bf Figure 15.1:\/} The $(0,A)$ slices.
\end{center}

From this picture we can see that the only edges emanating from
$\zeta_k$ are those corresponding to the pairs
$$(0,1); \hskip 30 pt (1,0); \hskip 30 pt (1,1); \hskip 30 pt (-1,1).$$
All of these edges point into the halfplane above the relevant
floor line.  This
what we wanted to establish.

\subsection{The Even Case}

The only place where we used the fact that $p/q$ is odd was
in Lemma \ref{oddkey1}.  We needed to know that the
number $t$ in Equation \ref{oddkey1} was odd.  This
no longer works when $p+q$ is odd.
However, when $p/q$ is even, the floor grid has a different
definition:  Only the even floor lines are present in
the grid. That is, the number $k$ in Equation \ref{oddkey1}
is an even integer.  Hence, for the floor lines in the even
case, the number $t$ is an integer.  The rest of the
proof of Lemma \ref{oddkey1} works word for word.
The rest of the proof of Statement 1 goes through word for word.

\newpage

\section{Proof of the Hexagrid Theorem II}
\label{hex2proof}

\subsection{The Basic Definitions}

As in the previous chapter, we will take $p/q$ odd until the very end.
It turns out that the secret to proving Statement 2 of the
Hexagrid Theorem is to use variants of the maps
$M_+$ and $M_-$ from Equation \ref{mptnotation}.
Let $A \in (0,1)$ be any parameter.
Let $\Lambda$ the lattice
from the Master Picture Theorem.   Let $\Pi \subset \R^3$ be
the plane defined by the relation $x+y=A$.

For $(m,n) \in \R^2$ we define 
$\Delta_+(m,n)=(x,y,z)$, where
\begin{equation}
x=2A(1-m+n)-m; \hskip 30 pt y=A-x; \hskip 30 pt z=Am.
\end{equation}
We also define
\begin{equation}
\Delta_-(m,n)=\Delta_+(m,n)+(-A,A,0).
\end{equation}
Note that $\Delta_{\pm}(m,n) \in \Pi$.
Indeed, $\Delta$ is an affine isomorphism
from $\R^2$ onto $\Pi$.

\begin{lemma}
Suppose that $(m,n) \in \Z^2$.  Then
$\Delta_{\pm}(m,n)$ and $M_{\pm}(m,n)$ are equivalent
mod $\Lambda$.
\end{lemma}

\startproof
Let $v_1,v_2,v_3$ be the three columns of the matrix
defining $\Lambda$.  So, $v_1=(1+A,0,0)$ and
$v_1=(1-A,1+A,0)$ and $v_3=(-1,-1,1)$.
Let
$$c_1=-1+2m; \hskip 30 pt c_2=1-m+2n; \hskip 30 pt c_3=n.$$
We compute directly that
$$M_+(m,n)-\Delta_+(m,n)=c_1v_1+c_2v_2+c_3v_3.$$
$$M_-(m,n)-\Delta_-(m,n)=c_1v_1+(c_2-1)v_2+c_3v_3.$$
This completes the proof.
\endproof

We introduce the vector
\begin{equation}
\zeta=(-A,A,1) \in \Lambda.
\end{equation}
Referring to the proof of our last result, we have $\zeta=v_2+v_3$.  This
explains why $\zeta \in \Lambda$.  Note that $\Pi$ is invariant under
translation by $\zeta$.

\subsection{Interaction with the Hexagrid}

Now we will specialize to the case when $A=p/q$ is an odd rational.
The results above hold, and we can also define the hexagrid.
We will see how the maps $\Delta_+$ and $\Delta_-$ interact
with the Hexagrid.  Let $L_0$ denote the wall line
through the origin.

\begin{lemma}
\label{parallel1}
$\Delta_{\pm}(L_0)$ is parallel to $\zeta$ and contains
$(-2A,A,0)$.
\end{lemma}

\startproof
We refer to the points in Figure 3.1.  The points
$v_5$ and $v_1$ both lie on $L_0$.  We compute
$$\Delta_+(v_5)-\Delta_+(v_1)=\frac{p^2}{p+q}\zeta.$$
Hence $\Delta_+(L_0)$ is parallel to $\zeta$.
We compute that $\Delta_+(0,0)=(2A,-A,0)$.
\endproof

We introduce the notation 
$\Pi(x)$ to denote the line in $\Pi$ that is parallel to $\zeta$ and
contains the point $(x,A-x,0)$.  For instance,
\begin{equation}
\Delta_+(0,0) \subset \Pi(2A); \hskip 30 pt
\Delta_-(0,0) \subset \Pi(A).
\end{equation}
Let $\Pi(r,s)$ denote the infinite
strip bounded by the lines $\Pi(r)$ and $\Pi(s)$.

For each pair of indices $(\epsilon_1,\epsilon_2) \in \{-1,0,1\}^2$, we let
$\Sigma(\epsilon_1,\epsilon_2)$ denote the set of lattice points
$(m,n)$ such that $L_0$ separates
$(m,n)$ from $(m+\epsilon_1,n+\epsilon_2)$.   
Now we define constants
$$
\lambda(0,1)=2A \hskip 30 pt
\lambda(-1,-1)=1-A^2;$$
\begin{equation}
\lambda(-1,0)=1+2A-A^2 \hskip 30 pt
\lambda(-1,1)=1+4A-A^2
\end{equation}

\begin{lemma}
Let $(\epsilon_1,\epsilon_2)$ be any of the $4$ pairs listed above.  
Let $\lambda=\lambda(\epsilon_1,\epsilon_2)$.  The following
$3$ statements are equivalent.
\begin{enumerate}
\item $(m,n) \in \Sigma(\epsilon_1,\epsilon_2)$.
\item $\Delta_+(m,n)$ is congruent mod $\Lambda$ to a point in the
interior of $\Pi(2A-\lambda,2A)$. 
\item $\Delta_-(m,n)$ is congruent mod $\Lambda$ to a point
in the interior of $\Pi(A-\lambda,A)$.
\end{enumerate}
\end{lemma}

\startproof
The formula $\Delta_-=\Delta_++(-A,A,0)$ immediately
implies the equivalence of the second and third statements.
So, it suffices to prove the equivalence of the first two statements.
We will consider the pair $(-1,0)$.  The other cases have
the same treatment. The set
$\Sigma(-1,0)$ is the intersection of
$\Z^2$ with the interior of some infinite strip, one
of whose boundaries is $L_0$.  To find the
image of this strip under $\Delta_+$, we just have to see what
$\Delta_+$ does to two points, one per boundary component
of the strip.  We choose the points
$(0,0)$ and $(1,0)$.    We already know that
$\Delta_+(0,0) \subset \Pi(2A)$.  We just have to
compute $\Delta_+(1,0)$.  We compute
$$\Delta(1,0)=(1,A-1,A) \subset \Pi(1-A^2).$$
This gives us $\lambda(-1,0)=1+2A-A^2$.
Our lemma follows from this fact, and from
the fact that $\Delta_+$ is an affine
isomorphism from $\R^2$ to $\Pi$.
\endproof

\subsection{Determining the Local Picture}

A crossing 
cell can consist of either $1$ edge or $2$, depending on
whether or not a vertex of the cell lies on a wall line.
According to Lemma \ref{specialcross}, the only
crossing cells with one edge are equivalent
mod $\Theta$ to the one whose center vertex is $(0,0)$.
For these {\it special\/} crossing cells, 
Statement 2 of the Hexagrid Theorem is obvious.
The door is just the central vertex.

The remaining crossing cells are what we call {\it generic\/}.
Each generic crossing cell has one vertex in
one of our sets $\Sigma(\epsilon_1,\epsilon_2)$, for
one of the $4$ pairs considered above.
We call $v \in \Sigma(\epsilon_1,\epsilon_2)$ a
{\it critical for\/} $(\epsilon_1,\epsilon_2)$.
$v$ and $v+(\epsilon_1,\epsilon_2)$ are the two
vertices of a crossing cell. To prove Statement 2
of the Hexagrid Theorem, we need to 
understand the critical vertices.  This means
that we need to understand the local picture of the
arithmetic graph in terms of the maps
$\Delta_+$ and $\Delta_-$.

We want to draw pictures as in the previous chapter,
but here we need to be more careful.  In the previous chapter,
our plane $\Pi$ contained the vector $(1,1,1)$.  Thus, we could
determine the structure of the arithmetic graph just
by looking at the intersection $\Pi \cap \cal R$.  Here
$\cal R$ is the polyhedron partition for the given parameter.
The situation here is different.
The vector $(1,1,1)$ is transverse to the plane $\Pi$.  
What we really need to do is to understand the way that
the plane $\Pi_{\alpha}$ intersects the our partition.
Here $\Pi_{\alpha}$ is the plane satisfying the
equation $x+y=A+2\alpha$. We think of
$\alpha$ an infinitesimally small but positive number.
More formally, we  take the geometric limit of
the set $\Pi_{\alpha} \cap {\cal R\/}$ as $\alpha \searrow 0$.

We say that a subset $S \subset \Pi$ is {\it painted\/} $(\epsilon_1,\epsilon_2,+)$
if $\Delta_+(m,n) \in S$ implies that
$\Delta_+(m,n)$ determines the pair $(\epsilon_1,\epsilon_2)$.
This is to say that $S$ is contained in the
Hausdorff limit of $\Pi_{\alpha} \cap {\cal R\/}_+(\epsilon_1,\epsilon_2)$
as $\alpha \to 0$.
We make the same definition with $(+)$ in place of $(-)$.
We think of $(\epsilon_1,\epsilon_2,\pm)$ as a kind of
color, because these regions are assigned various
colors in Billiard King.  For instance $(0,1,\pm)$ is green.
There is essentially one {\it painting\/} of 
$\Pi$ for $(+)$ and one for $(-)$.

To visualize the painting, we identify $\Pi$ with $\R^2$
using the map $(x,y,z) \to (x,z)$.  We just drop the second
coordinate. The vector
$\zeta$ maps to the $(-A,1)$.  Thus, our whole
painting is invariant under translation by this vector.
Each wall
of $\cal R$ intersects $\Pi$ in a line segment
whose image in $\R^2$ is either horizontal or
vertical.   The endpoints of each such segment
have coordinates that are simple rational
combinations of $1$ and $A$.  For this reason,
we can determine the intersection we seek
just by inspecting the output from Billiard King.
In practice, we take $\alpha=10^{-5}$, examine
the resulting picture, and then adjust the
various vertices slightly so that their coordinates
are small rational combinations of $1$ and $A$.

\subsection{An Extended Example}

We consider the pair $(0,1)$ in detail.
We will draw pictures for the parameter $A=1/3$,
though the same argument works for any parameter.
There is no polyhedron ${\cal R\/}_-(0,1)$, so
no points are painted $(0,1,-)$.  The interesting
case is $(0,1,+)$.
First of all, we only care about
points in our strip $\Sigma(0,1)$.  So, we only
need to understand the portion of our painting
that lies in our strip $\Pi(0,2A)$.   In
$\R^2$ (considered as the $xz$ plane), our
strip is bounded by the lines $x=-zA$ and $x=-zA+2A$.

We will first study the picture when $z=0$.
Referring to Figure 16.1,
the shaded triangles correspond to
${\cal R\/}(0,1)$.  The thick line corresponds
to the intersection of $\Pi$ with our fiber.
The black dot is the point $(A,0,0,A)$.
Moving away from the black dot, the white dots are
$$(0,A,0,A); \hskip 15 pt
(-A,2A,0,A); \hskip 15 pt
(-1,1+A,0,A.$$

It we move the thick line an infinitesimal amount in the
direction of $(1,1)$, we see that it crosses
though a shaded region whose diagonal edge is
bounded by the points
$(0,A)$ and $(A,0)$. 
The only tricky part of the analysis is that the
point $(A,0)$ determines the pair $(0,0)$ and the
point $(0,A)$ determines the pair $(-1,1)$.

\begin{center}
\resizebox{!}{2.7in}{\includegraphics{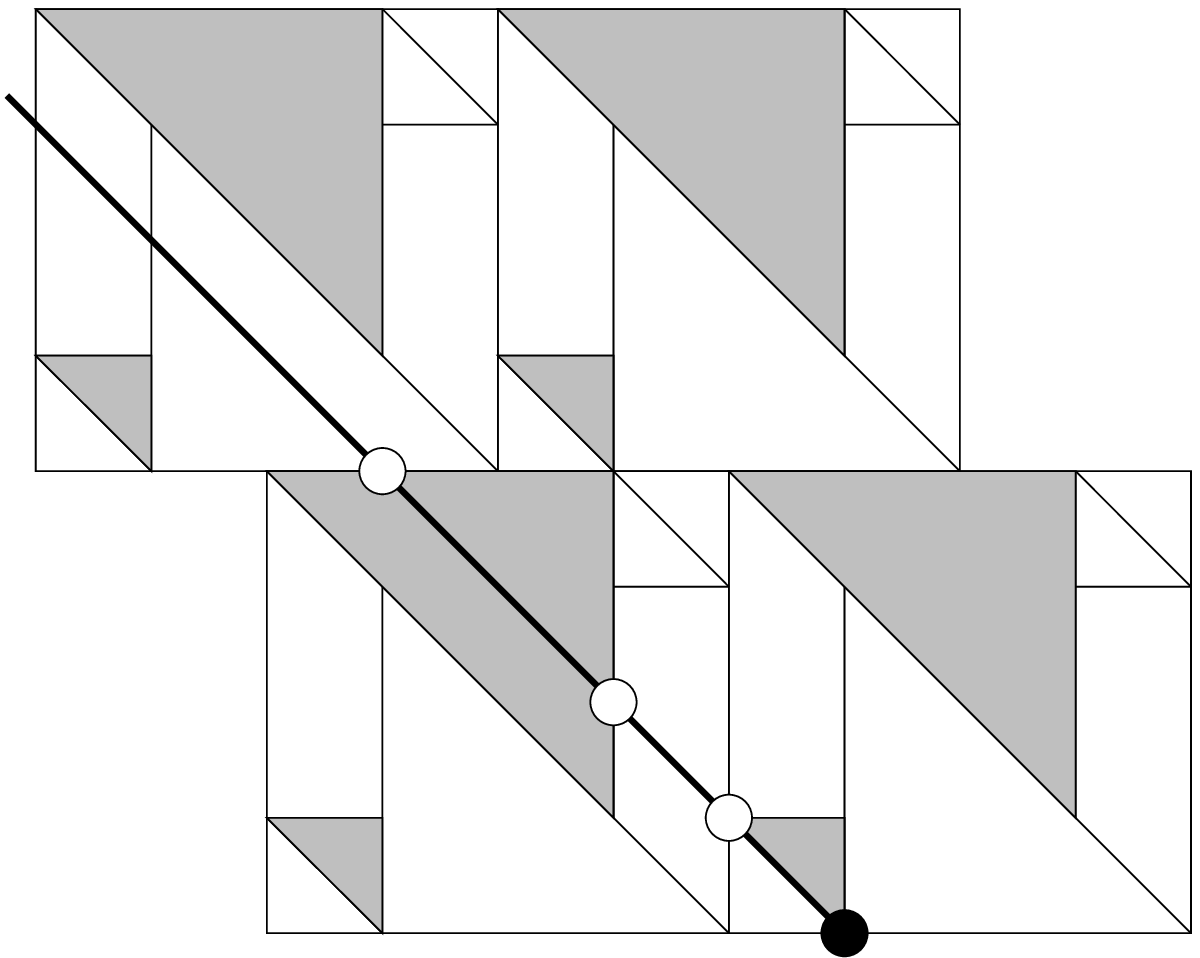}}
\newline
{\bf Figure 16.1:\/} Slicing the $0$ fiber.
\end{center} 

From this discussion, we conclude that 
$(0,A) \times \{0\}$ is painted $(0,1,+)$.
For later use, we remark that $(0,0)$ is painted
$(-1,1,+)$ and $(A,0)$ is painted $(0,0,+)$.
Looking at the picture, we also see that
$(-1,-A) \times \{0\}$
is painted $(0,1,+)$.  Notice, however, that
this set lies outside our strip.  It is irrelevant.

Figure 16.2 shows the picture for a typical parameter
$z \in (0,1)$.  We choose $z=1/6$, though the
features of interest are the same for any choice of
$z$. The interested reader can see essentially any
slice (and in color) using Billiard King.

The black dot and the white dots have the same
coordinates as in Figure 16.1.   Notice that the
point $(0,A,z,A)$ lies at the bottom corner
of a shaded region.  This remains true for all $z$.
We conclude that the open line 
segment $\{0\} \times (0,1)$
is painted $(0,1,+)$.   Similarly, the rectangle
$(-1,-A) \times [0,1]$ is painted $(0,1,+)$.  However,
this rectangle is disjoint from the interior of
our strip.  Again, it is irrelevant.

Recalling that our painting is invariant under translation
by $(0,1,+)$, we can now draw the portion of plane
painted $(0,1,+)$ that is relevant to our analysis.
To give the reader a sense of the geometry, we
also draw one copy of the irrelevant rectangle.
Again, we draw the picture for the parameter $A=1/3$.
The interested reader can see the picture for
any parameter using Billiard King.

\begin{center}
\resizebox{!}{2.7in}{\includegraphics{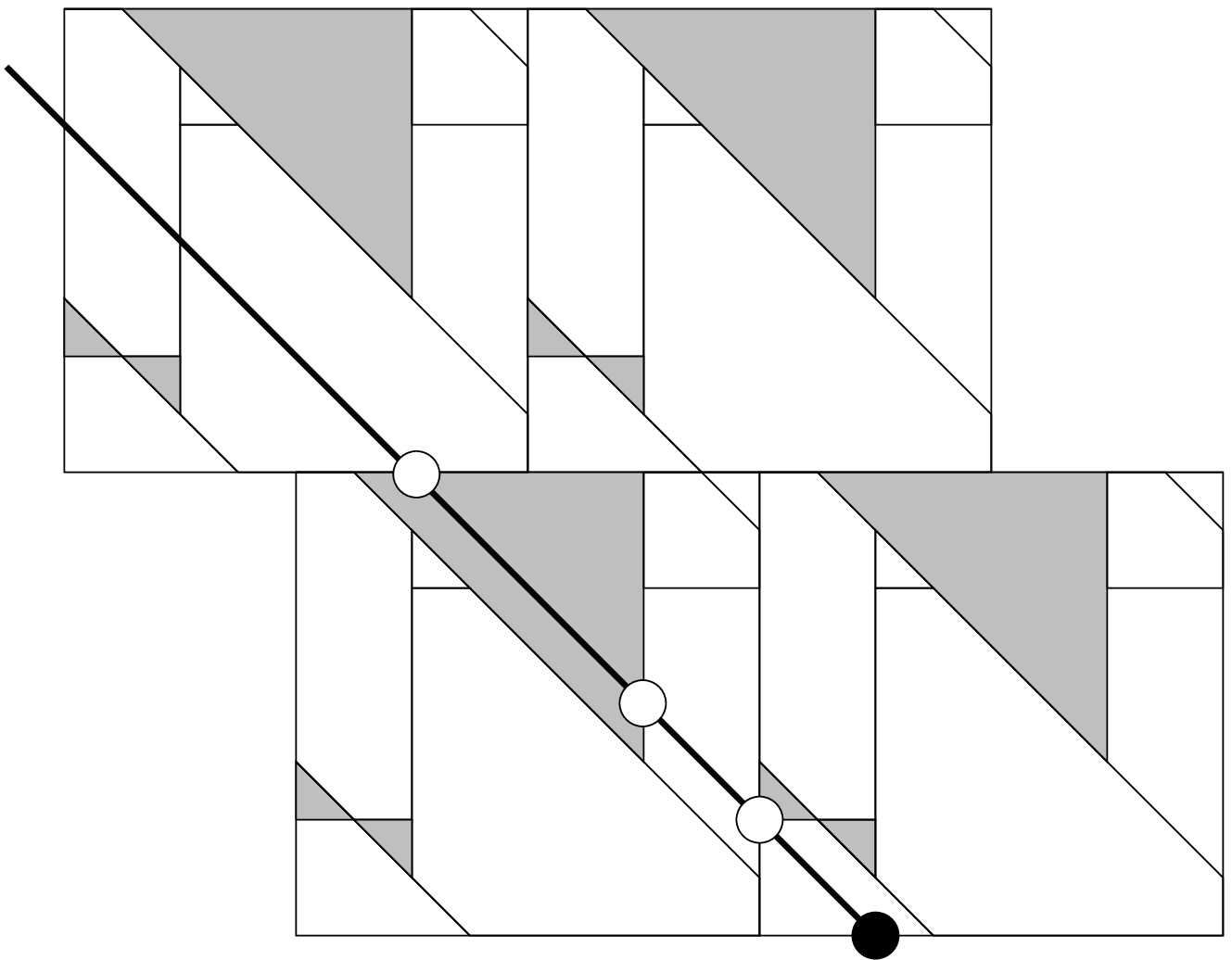}}
\newline
{\bf Figure 16.2:\/} Slicing a typical fiber.
\end{center} 

  In Figure 16.3, the arrow
represents the vector $(-A,1)$. 
The black dot is $(0,0)$ and the white dot
is $(A,0)$.  The thick zig-zag, which is
meant to go on forever in both directions,
is the relevant part of the painting.  The
lightly shaded region is the strip of
interest to us.

\begin{center}
\resizebox{!}{3in}{\includegraphics{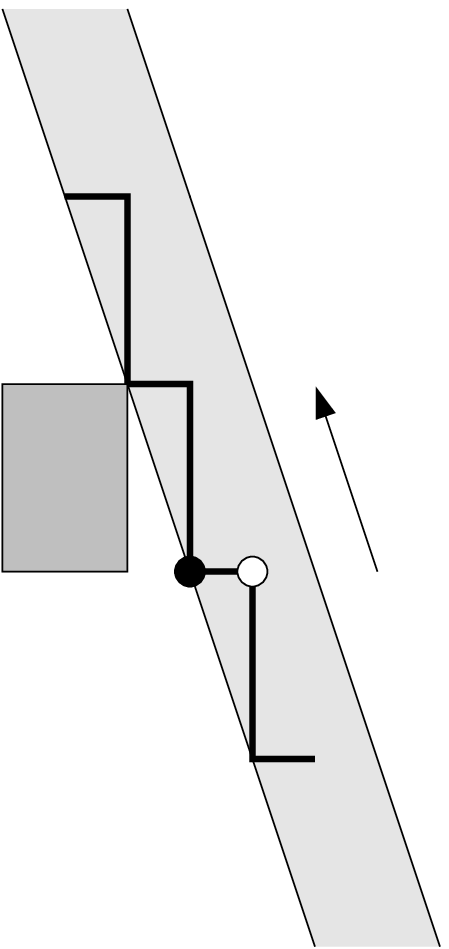}}
\newline
{\bf Figure 16.3:\/} The relevant part of the $(0,1,+)$ painting.
\end{center} 

\subsection{The Rest of the Painting}

We determine the rest of the painting using the
same techniques.  The interested reader can see
everything plotted on Billiard King.
The left side of
igure 19.4 shows the relevant part of
the $(+)$ painting. The right side shows
the relevant part of the $(-)$ painting.  The dots
are exceptional points in the painting.  The
two grey dots at the endpoints correspond
to the right endpoints of the special crossing cells.
We have shown a ``fundamental domain'' for the
paintings.  The whole painting is obtained taking
the orbit under the group $\langle \zeta \rangle$.
In our picture, $\zeta$ acts as translation by the vector
$(-A,1)$, because we are leaving off the $y$ coordinate.
In particular, the two endpoints of the $L$ are
identified when we translate by this group.

The small double-braced labels, such as $((0,1))$, indicate
the paint colors.  The large labels, such
as $(0,0)$, indicate the coordinates in the plane.
Note that the point $(x,z)$ in the plane actually
corresponds to $(x,A-z,y)$ in $\Pi$.  The grey
vertices on the left corresponds to $\Delta_+(0,0)$.
The grey vertices correspond to the various
images of points on the special crossing cell.
These vertices are not relevant to our analysis
of the points that are critical relative to
our $4$ pairs.

\begin{center}
\resizebox{!}{3.5in}{\includegraphics{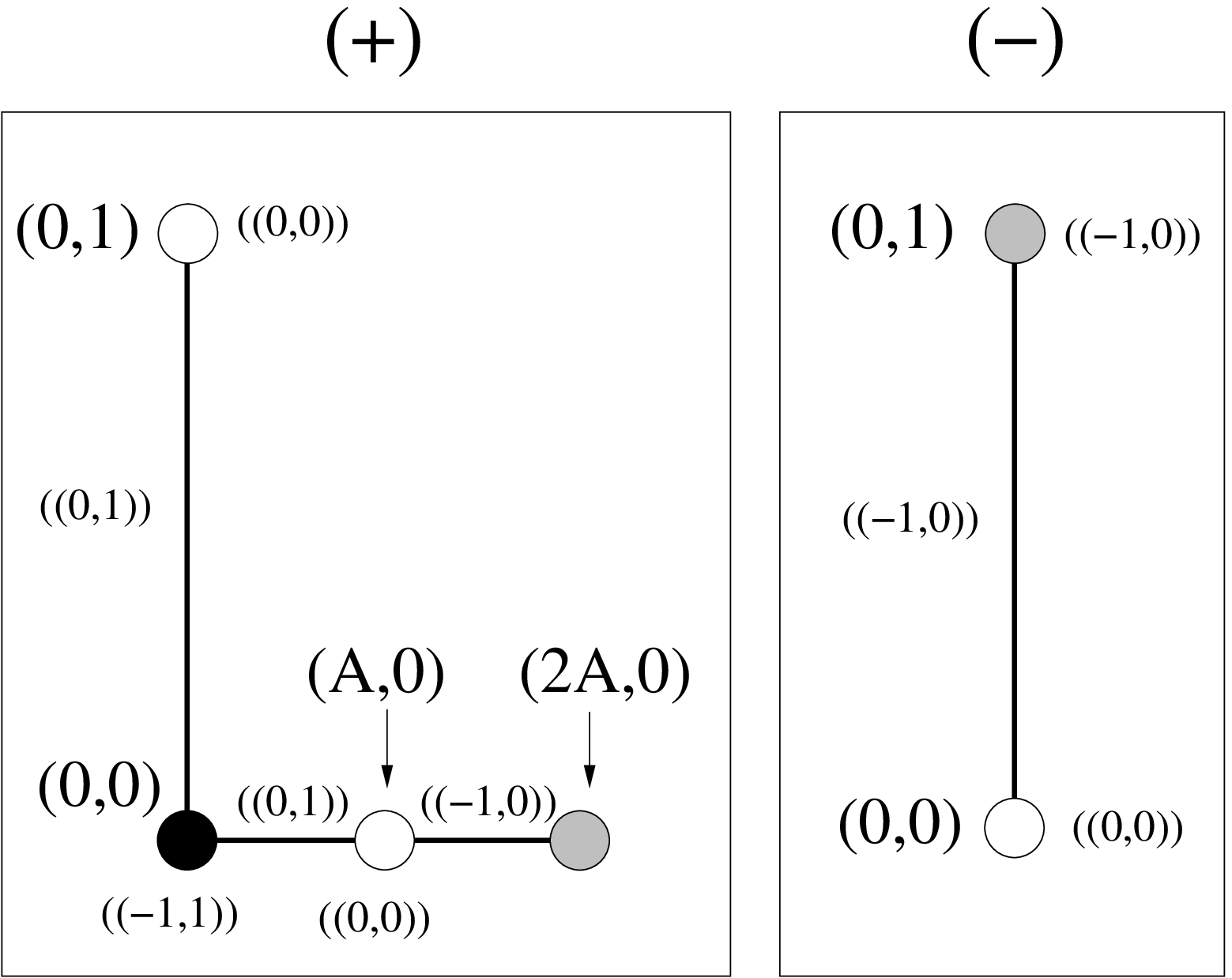}}
\newline
{\bf Figure 16.4:\/} The relevant part of the $(+)$ painting
\end{center} 

Say that a vertex $v=(m,n)$ is {\it critical\/} if it either
lies in $\Theta$ or else is critical for
one of our strips.  The point $(0,0)$ is the
center vertex of a special crossing cell. Hence,
By Lemma \ref{specialcross},
the critical vertices are in bijection with the
crossing cells.
Given our analysis above, we see that $v \in \Z^2$ is
{\it critical\/} if and only if it satisfies the following criterion.
{\it Modulo the action of $\Lambda$, the point
$\Delta_+(v)$ (respectively the point $\Delta_-(v)$) lies in one of the
colored parts of the painting on the left (respectively right) in Figure 16.4.\/}

Recalling that $\Delta_+=\Delta_-+(A,-A,0)$, we can eliminate
$\Delta_-$ from our discussion.  We translate the right hand side
of Figure 16.4 by $(A,0)$ and then
superimpose it over the left hand side.
(This translation does not reflect the way the two halves
of Figure 16.4 are related to each other on the page.)
See Figure 16.5.   The result above has the following
reformulation.

\begin{lemma}[Critical]
A vertex $v$ is critical if and only of
$\Delta_+(v)$ is equivalent mod $\Lambda$ to a point
colored portion of Figure 16.5.
\end{lemma}

\begin{center}
\resizebox{!}{4in}{\includegraphics{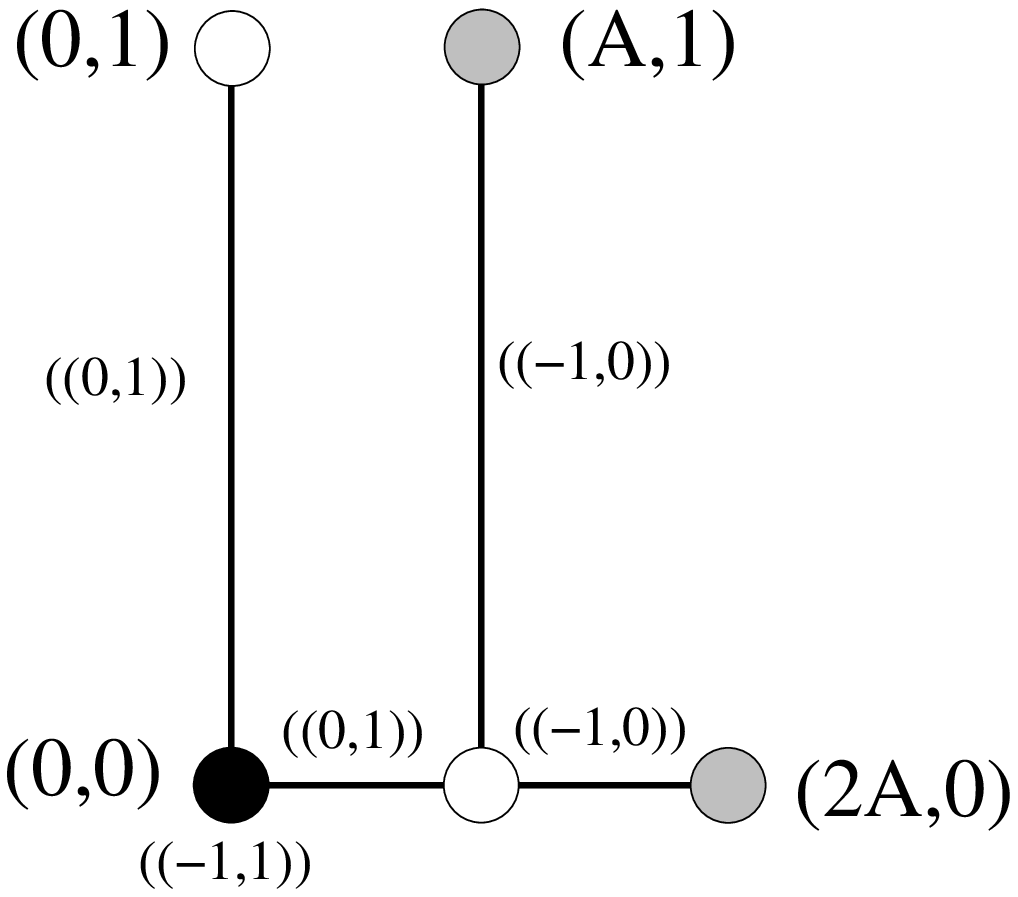}}
\newline
{\bf Figure 16.5:\/} Superimposed paintings
\end{center} 

Our drawing of Figure 16.5 somewhat hides the symmetry of
our picture.  In Figure 16.6, we show several translates
of this fundamental domain at the same time, without
the labels.  We also show the strip $\Pi(0,2A)$.
The pattern is meant to repeat endlessly in both 
directions.   The line on the left is $\Pi(0)$ and
the line on the right is $\Pi(2A)$.  Again, we
are drawing the picture for the parameter $A=1/3$.
The combinatorial pattern is the same for any $A$.

\begin{center}
\resizebox{!}{5.2in}{\includegraphics{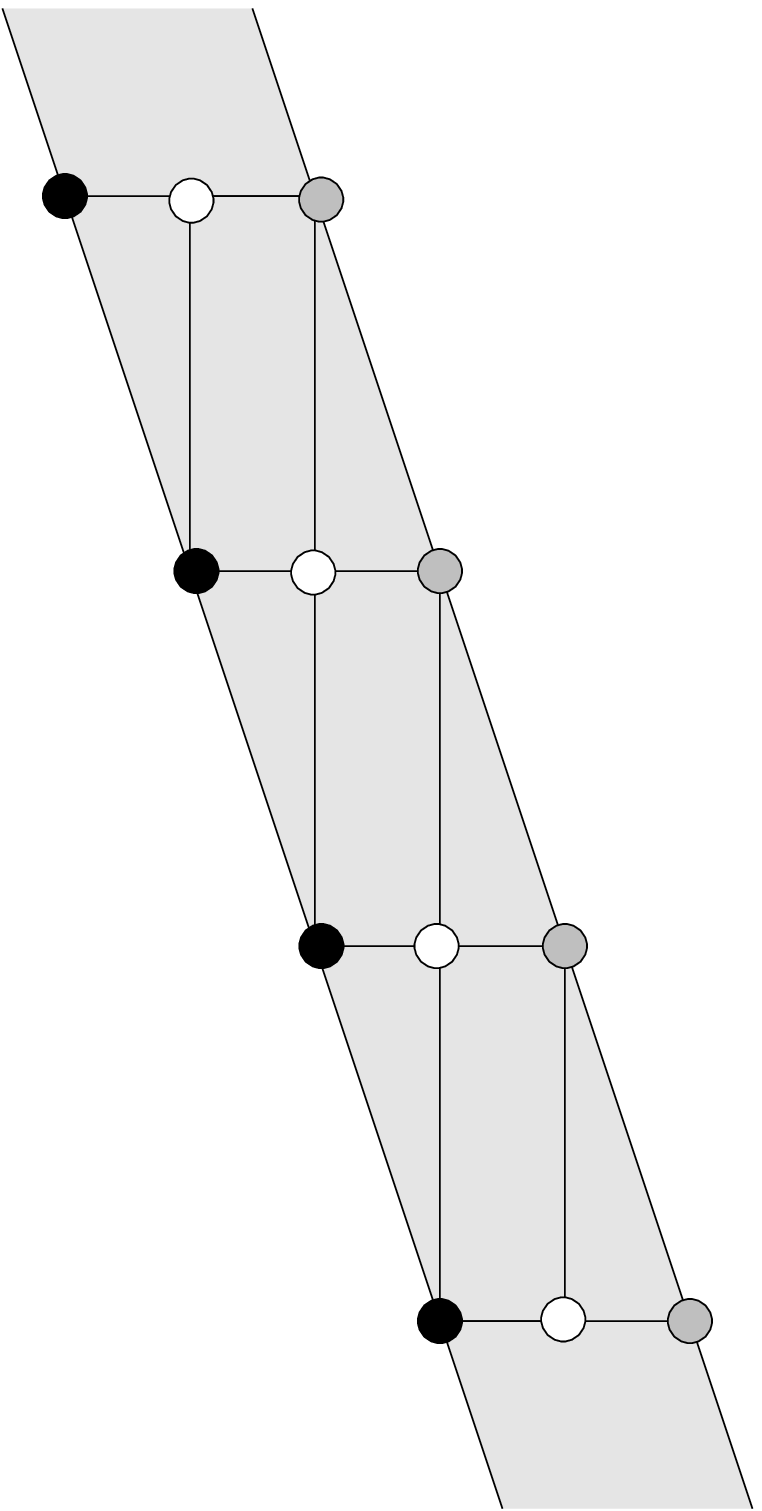}}
\newline
{\bf Figure 16.6:\/} Superimposed paintings
\end{center}

 To prove the hexagrid theorem, it only remains to
identify the lattice points in the Critical Lemma with the 
doors from the Hexagrid Theorem.

\subsection{The End of the Proof}

Now we interpret the Critical Lemma algebraically.
A vertex $v \in \Z^2$ is {\it critical\/} if and only
if $\Delta_+(v)$ is equivalent mod $\Lambda$ to one
of the following kinds of points.

\begin{enumerate}
\item $(2A,-A,0)$.
\item $(0,A,0)$.
\item $(x,A-x,0)$, where $x \in (0,2A)-\{A\}$.
\item $(0,A,z)$, where $z \in (0,1)$.
\end{enumerate}
As we point out in our 
subsection headings, each case corresponds to a
different feature of our painting in Figures 16.5 and 16.6

\subsubsection{Case 1:  The Grey Dots}

Note that $\Delta_+(0,0)=(2A,-A,0)$.  Moreover,
$\Delta_+(v) \equiv \Delta_+(v')$ mod $\Lambda$ iff
$M_+(v) \equiv M_-(v)$ mod $\Lambda$ iff
$v \equiv v'$ mod $\Theta$.  Hence,
Case 1 above corresponds precisely
to the special crossing cells.  
The door associated to $v$ is precisely $v$.
In this case, the door is associated to the
wall above it.

\subsubsection{Case 2: The Black Dots}

Note that $\Delta_+(0,-1)=(0,A,0)$.  Hence
the second case occurs iff $v$ is equivalent
mod $\Theta$ to $(0,-1)$. But $(0,-1)$ is the
vertex of a crossing cell whose other
vertex is $(-1,0)$.  The door associated to
$(0,-1)$ is $(0,0)$.   In this case, the door
is associated to the wall below it.

\subsubsection{Case 3: Horizontal Segments}
\label{vbias}

We are going to demonstrate the bijection between
the Type 1 doors not covered in Cases 1 and 2
and the critical points that arise from Case 3 above.

Let $v$ be a critical point.
Using the symmetry of $\Theta$, we can arrange
that our point $v$ is closer to $L_0$
than to any other wall line.  In this
case, $\Delta_+(v)$ lies in the strip $\Pi(0,2A)$.
Hence $v \in \Sigma(0,1)$.  Hence $L_0$
separates $v$ from $v+(0,1)$.  Let
$y \in (n,n+1)$ be such that $(m,y) \in L_0$.

The third coordinate of $\Delta_+(v)$ is
an integer.  Setting $v=(m,n)$, we see that
$Am \in \Z$.  Hence $q$ divides $m$.  Hence
$v=(kq,n)$ for some $k \in \Z$.  
Hence $(kq,y)$ is a Type 1 door.

For the converse, suppose that the point
$(kq,y)$ is a Type 1 door and that
$n=\underline y$.  Let $v=(kq,n)$. We
want to show that $v$ is critical.
By construction $(kq,n) \in \Sigma(0,1)$.
Hence $\Delta_+(kq,n) \in \Pi(0,2A)$.
But the third coordinate of
$\Delta_+(kq,n)$ is an integer.
Hence $\Delta_+(kq,n)$ is equivalent
mod $\Lambda$ to a point of the
form $(x,A-x,0)$.  Here $x \in (0,2A)$.

If $x=A$ then $v$ lies on the centerline
of the strip $\Sigma(0,1)$.  But then
$y-\underline y=1/2$.  This contradicts
Lemma \ref{nohalf}.  Hence $A \not = x$.

Now we know that
$\Delta_+(kq,n)$ satisfies Case 3 above.
Hence $(kq,n)$ is critical, either
for $(0,1)$ or for $(-1,0)$.  These are
the relevant labellings in Figure 16.5.
Note that $L_0$ has positive slope greater than $1$.
Hence $$(kq,n) \in \Sigma_+(0,1) \cap \Sigma_+(-1,0).$$
If $\Delta_+'(v)$ is colored
$(0,1)$, then $v$ is critical for $(0,1)$. If
 $\Delta_+'(v)$ is colored
$(-1,0)$, then $v$ is critical for $(-1,0)$. 
So, $v$ is always  vertex of a crossing cell.

\subsubsection{Case 4: Vertical Segments}
\label{hbias}

We are going to demonstrate the bijection between
Type 2 doors and the critical points that
arise from Case 4 above.

We use the symmetry of $\Theta$ to guarantee that
our critical point is closer to $L_0$ than to any
other wall line.
As in Case 3, the point
 $\Delta_+(v) \in \Pi(0,A)$.  Hence
$v \in \Sigma(0,1)$ and $L_0$ 
separates $v$ from $v+(0,1)$.  We define
$y$ as in Case 3.  We want to show that
$(m,y) \in L$ is a Type 2 door.

Since we are in Case 4, the
first coordinate of 
$\Delta_+(v)$ lies in $A\Z$.  The idea here is
that $\Delta_+(v)$ is equivalent mod $(-A,A,1)$ to
a point whose first coordinate is either $0$ or $A$. Hence
$$x=2A(1-m+n)-m \in A\Z.$$
Hence $x/A \in \Z$.  Hence $m/A \in \Z$.
Hence $m=kp$. 
By Lemma \ref{door2}, the
point $(x,y)$ is a door.

Conversely, suppose that the point $(kp,y)$ is door contained in $L_0$.
Let $n=\underline y$.  Then 
$(kp,n) \in \Sigma(0,1)$ and the
first coordinate of $\Delta_+(kp,n)$ lies
in the set $A\Z$.   Also, $\Delta_+(kp,n) \in \Pi(0,2A)$.
Hence, $(kp,n)$ satisfies Case 4 above.  Hence
$(kp,n)$ is critical for either $(0,1)$ of $(-1,0)$.
In either case, $(kp,n)$ is a vertex of a
crossing cell.

\subsection{The Pattern of Crossing Cells}

Our proof of the Hexagrid Theorem is done, but we can
say more about the nature of the crossing cells.
First of all, there are two crossing cells 
consisting of edges of slope $\pm 1$.  These
crossing cells correspond to the black and
grey corner dots in Figure 16.6.  

The remaining
crossing cells involve either vertical or
horizontal edges.  These crossing cells
correspond to the interiors of the segments
in Figures 16.5 and 16.6.  
Let $v=(m,n)$ be the critical vertex
associated to the door $(m,y)$.   
Then $v$ is critical either for
$(0,1)$ or $(-1,0)$.  In the former
case, the crossing cell associated to
$v$ is vertical, and in the latter
case it is horizontal.  Looking
at the way Figure 16.5 is labelled,
we see that 
\begin{itemize}
\item The crossing cell is vertical if
$y-\underline y>1/2$.
\item The crossing cell is horizontal if
$y-\underline y<1/2$.
\end{itemize}
The case $y-\underline y=1/2$ does not occur,
by lemma \ref{nohalf}.

There are exactly $p+q$ crossing
cells mod $\Theta$.  These cells
are indexed by the value of
$y-n$.  The possible numbers are
$$\{0,\frac{1}{p},...,\frac{p-1}{p},\frac{1}{q},...\frac{q-1}{q},1\}.$$
In all cases we have $y-n=y-\underline y$, except when
$n=y-1$.

Figure 16.7 shows, for the
case $p/q=3/5$, the images of the critical
vertices, on one fundamental domain for
Figure 16.6. (The fundamental domain here
is nicer than the one in Figure 16.5.)
We have labelled the image points by
the indices of the corresponding crossing
cells.  The lines inside the dots show the nature of
the crossing cell.
The dashed grid lines in the figure are
present to delineate the structure.
The lines inside the dots show the nature of
the crossing cell.

One can think of the index values in the following way.
Sweep across the plane from right to left
by moving a line of slope $-5/3$ parallel to itself.
(The diagonal line in Figure 16.7 is one such line.)
The indices are ordered according to how
the moving line encounters the vertices.
The lines we are using correspond to the lines in
$\Pi$ that are parallel to the vector $\zeta$.

\begin{center}
\resizebox{!}{3.6in}{\includegraphics{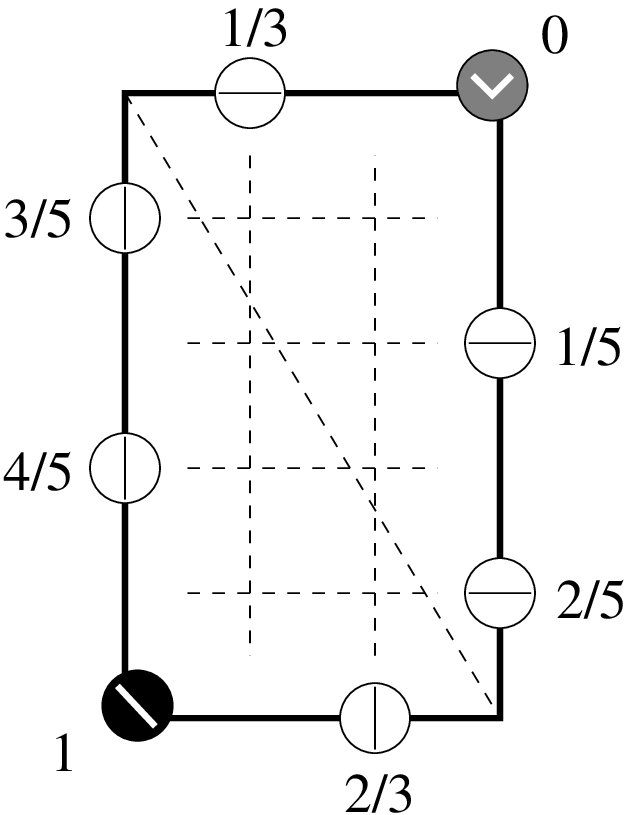}}
\newline
{\bf Figure 16.7:\/} Images of the critical points
\end{center}

Figure 16.7 is representative of the general
case.  It is meant to suggest the general pattern.
We hope that the pattern is clear.

\subsection{The Even Case}
When $p/q$ is an even rational, all the constructions
in this chapter go through word for word.

It appears that we used the fact that $p/q$ is
odd in Cases 3 and 4 in the last section,
but this isn't so.  We only used the
fact that Corollary \ref{nohalf} was true,
and that Lemma \ref{door2} was true. At the
time, we had only proved these results in the
odd rational case.  However, since these
results hold in the even case,
the arguments for Cases 3 and 4 go through
word for word.

A final remark on Case 3:
Case 3 required us to use Corollary \ref{nohalf} to
rule out the possibility that the point
$(m,n)$ is equivalent mod $\Lambda$ to 
the points $(A,0,0)$.  This can happen in
the even case, and indeed it happens 
when $(m,y)$ is a crossing of the kind
we are no longer calling a door.  In other
words, this does not happen for a {\it door\/}
because we have forced the situation.  

\newpage

\section{The Barrier Theorem}
\label{barrier}	

We remind the reader that we don't need the material in this
chapter until Part VI.

\subsection{The Result}
\label{minor}

Let $A=p/q$ be an even rational.
All the components of $\widehat \Gamma=\widehat \Gamma(A)$
are embedded polygons.  Say that a {\it low component\/}
is one that contains a low vertex.  The component
$\Gamma$ containing $(0,0)$ is a distinguished
low component.  The infinite set of components
$\Gamma+kV$, with $k \in \Z$ are translates
of $\Gamma$. Here $V=(q,-p)$ as usual.
We call these components {\it major\/} components.
We call the remaining low components {\it minor\/}
components.

\begin{center}
\resizebox{!}{4.7in}{\includegraphics{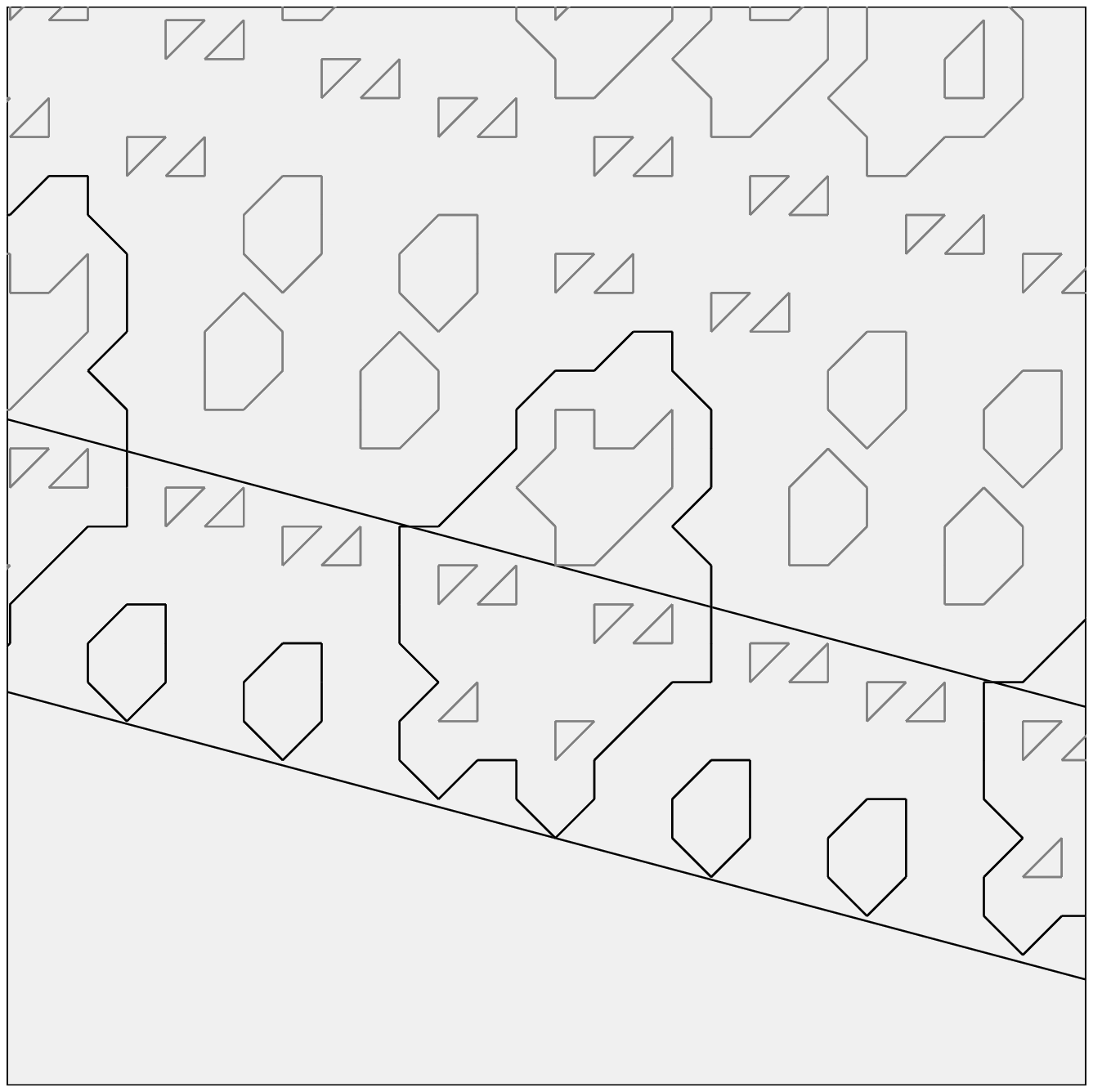}}
\newline
{\bf Figure 17.1:\/} Components of $\widehat \Gamma(4/15)$ and
a barrier.
\end{center}

Figure 17.1 shows some of $\widehat \Gamma(4/15)$.
The three biggest polygons are major components,
and the little polygons along the bottom
are minor components.  Figure 17.1 also shows a
barrier, parallel to the baseline, which is only
crossed by the major components.  The Barrier
Theorem describes this barrier and establishes
its basic properties.

Referring to Equation \ref{induct0}, one of the two
rationals $p_{\pm}/q_{\pm}$ is even and one is odd.
Let $p'/q'$ denote whichever of these
rationals is odd.
We call $p'/q'$ the {\it odd predecessor\/}
of $p/q$.  We say that the {\it barrier\/} is the line
parallel to $V$ that contains the point
\begin{equation}
\bigg(0,\frac{p'+q'}{2}\bigg)
\end{equation}

\begin{theorem}[Barrier]
Modulo translation by $V$,
only $2$ edges of $\widehat \Gamma$ cross
the barrier, and these lie on major components.
Hence, no minor component of $\widehat \Gamma$ crosses the barrier.
\end{theorem}

One could think of the Barrier Theorem as an
improvement of Statement 1 of the Hexagrid Theorem.
Statement 1 of the Hexagrid Theorem bounds the
distance that any low component can rise above
the baseline.  The Barrier Theorem gives a
bound that is at least twice as good for all
the minor components.

We have stated the precise version of the Barrier Theorem
that we need for our applications, but the Barrier
Theorem is really part of a more robust general
theorem.  If $A^*$ is a parameter that is close to
$A$ in the sense of Diophantine approximation, then
the line $\Lambda^*$ parallel to $V$ and containing the point
\begin{equation}
\bigg(0,\frac{p^*+q^*}{2}\bigg)
\end{equation}
is not frequently crossed by $\widehat \Gamma$.
The basic reason is that
$\Lambda^*$ serves as a kind of memory
of the Hexagrid Theorem for the parameter $A^*$.
The two graphs $\widehat \Gamma$ and $\widehat \Gamma^*$
mainly agree along $\Lambda^*$, and the only crossings
take places at the few mismatches in the graphs.

We will prove the Barrier Theorem using the same
ideas that we used to prove Statement 1 of the 
Hexagrid Theorem.  Mainly we shall be interested
in the differences between the Barrier Theorem and
Statement 1 of the Hexagrid Theorem.

\subsection{Review of the Hexagrid Proof}

Let us recall the idea behind the proof of Statement 1
of the Hexagrid Theorem, given in \S \ref{hex1proof}.
Our main idea was to analyze points just above the
floor line and observe that no such point contained
an edge of $\widehat \Gamma$ that crossed the
floor line.  Given the Master Picture Theorem,
this amounted to checking that the image of such
points never landed in a ``bad polygon'' of
the partition -- one that would assign to the
vertex a crossing edge.  Here were the $3$ main
ideas.

\begin{enumerate}
\item In Lemma \ref{oddkey1}, we
computed that $M_-$ (one of the two
classifying maps from the Master Picture
Theorem) maps
each floor line to a point of the form
$(\beta,0,0) \in \R^3/\Lambda$. 
\item In Lemma \ref{oddkey2} we deduced from
Lemma \ref{oddkey1}
that $M_{\pm}$ mapped lattice points just above the floor
line into a plane $\Pi_{\pm}$.   Both planes were translates
of the plane $\Pi$ containing $(1,0,0)$ and $(1,1,1)$.
The relevant lattice points were contained in
the strips $\Sigma(\epsilon_1,\epsilon_2)$ for
various choices of $\epsilon_1,\epsilon_2 \in \{-1,0,1\}$.
\item We sliced our partition ${\cal P\/}_{\pm}$ by the
planes $\Pi_{\pm}$ and simply checked that no relevant
vertex landed in a bad polygon.  Figures 15.1 and 15.2
showed the relevant picture for $A=1/3$.  The relevant 
vertices are contained in
the strips $\Sigma(\epsilon_1,\epsilon_2)$ for
various choices of $\epsilon_1,\epsilon_2 \in \{-1,0,1\}$
and the relevant domains were strips
$\Pi_{\pm}(0,\lambda(\epsilon_1,\epsilon_2)) \subset \Pi_{\pm}$.
See Lemma \ref{oddkey5}.
\end{enumerate}

Now we explain the change that occurs when we pass to the
present situation.
Let $\Lambda$ be the barrier line. 

\begin{lemma}
There is some real $\beta$ such that
\begin{equation}
\label{nearimage}
M_-(\Lambda)=\bigg(\beta,\pm \frac{1}{q},0\bigg)
\end{equation}
\end{lemma}

\startproof
We think of $M_-$ as acting on all of $\R^2$.  Then
$M_-$ is constant along $\Lambda$. 
If we use the map $M'$, relative
to $A'$, then we get a point $(\beta',0,0)$ by
the previous calculation. But $M(m,n)-M(m,n)=(\pm 2/q,0)$.
But $M_-=\mu_- \circ M$, where $\mu_-(2t)=(t-1,t,t)$
mod $\Lambda$.  Putting these two facts together
gives the proof.
\endproof

\subsection{Proof of the Barrier Theorem}

We will suppose that $A'<A$ until the end of the  section.
When $A'>A$, the $(+)$ option
in Equation \ref{nearimage} is taken.

The maps $M_{\pm}$ map $\Lambda$ into planes
$\Pi'_{\pm}$ that are obtained by translating
$\Pi_{\pm}$ in the $y$ direction by $1/q$.
The same argument as in Lemma \ref{oddkey2}
shows that $M_{\pm}$ maps the relevant lattice points -- namely
those in the strips $\Sigma(\epsilon_1,\epsilon_2)$ -- into
the regions $\Pi'_{\pm}(0,\lambda)$.
Here $\Pi'_{\pm}(0,\lambda) \subset \Pi'$
is a translate of the strip 
considered in Lemma \ref{oddkey5}.

The planes $\Pi_{\pm}$ are transverse
to the walls defining our partition.  When we translate
$\Pi_{\pm}$ off itself by $(0,1/q,0)$, the intersections we
see are practically the same.  Note that
$\Pi_{\pm}$ is not transverse to the partition itself,
just to the walls.  
When we translate, some new regions
pop into view.  Figure 17.2 shows one period of
the exact picture for
$\Pi_+$ relative to the parameter $4/15$.  The little
lines in the middle refer to similar lines drawn
in Figure 15.1.

\begin{center}
\resizebox{!}{4.2in}{\includegraphics{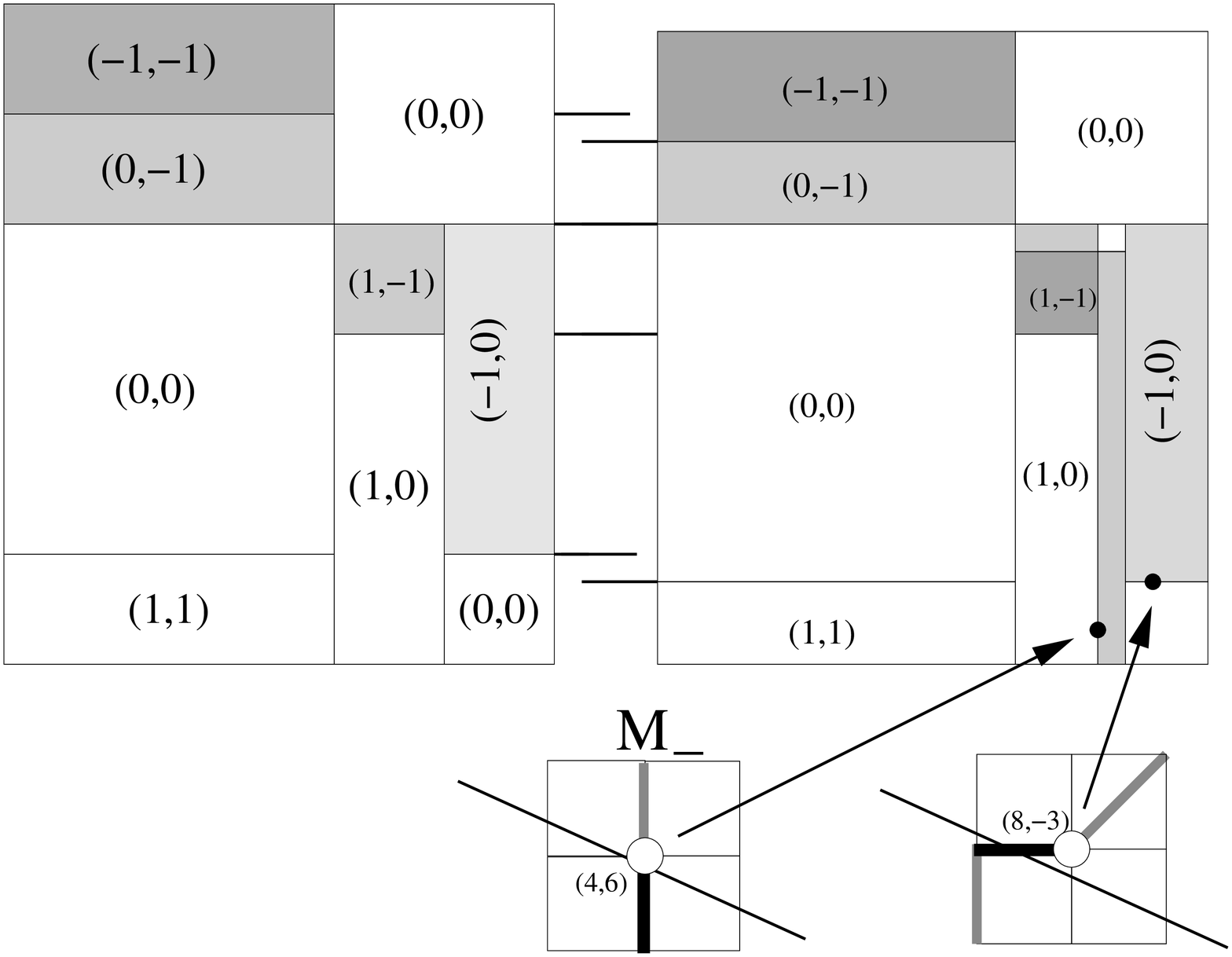}}
\newline
{\bf Figure 17.2:\/} The slices $\Pi_-$ and $\Pi'_-$.
\end{center}

The left hand side shows the slice $\Pi_-$.  The right hand
side shows the slice $\Pi_-'$.  The point on the right is
$M(-3,8)$, the image of a vertex above the barrier incident
to one of the crossing edges.  The lightly shaded region above the point
assigns the edge $(-1,0)$ to $M_-(-3,8)$, and this edge
crosses the barrier. The bottom of the figure shows this.
Likewise, the skinny rectangle assigns the edge
$(0,-1)$ to $M_-(4,6)$.  This is also shown at
the bottom of the figure.  
Were we to analyze the picture relative
to the parameter $A'=3/11$, these offending points would
get assigned non-crossing edges.  
 
In the new setting, our analysis for Statement 1 of the
Hexagrid Theorem does not {\it completely succeed\/} 
for two potential reasons.
\begin{enumerate}
\item The image of a relevant vertex might lie in one of the
newly appearing regions.  These regions all have width $1/q$.
\item Tee image of a relevant vertex might lie in a different
one of the old regions because the region has a slightly
different rectangle and/or location in the new slice.
This is what happens in our example.
In this case, the change in each edge is at most $1/q$.
\end{enumerate}
The bounds on the changes come from the equations for the
walls defining our partitions.

Now we make an analysis of how many crossings one
can get in Figure 17.2.  The images of the relevant
vertices all lie on a diagonal line of slope $1$.
This line starts on the bottom edge (on the
right hand figure.)
The difference in the coordinates between
successive points is $1/q$. 
Thus, each modified rectangle can give rise to
one new crossing.  Likewise, each new rectangle
can give rise to one new crossing.  This 
implies that there are at most $5$ crossings
entailed by the picture.

We can reduce this number by looking more closely.
The rectangle labelled $(0,-1)$ remains out of
range of the image $M_-(\Sigma(0,-1))$.  The
same argument as in the Hexagrid Theorem I applies here.
Thus, this rectangle entail no crossings.  This gets us
down to $4$.
Here is a trick to get down to $2$.  Figure 17.3 shows
$3$ of the relevant rectangles.

\begin{center}
\resizebox{!}{2in}{\includegraphics{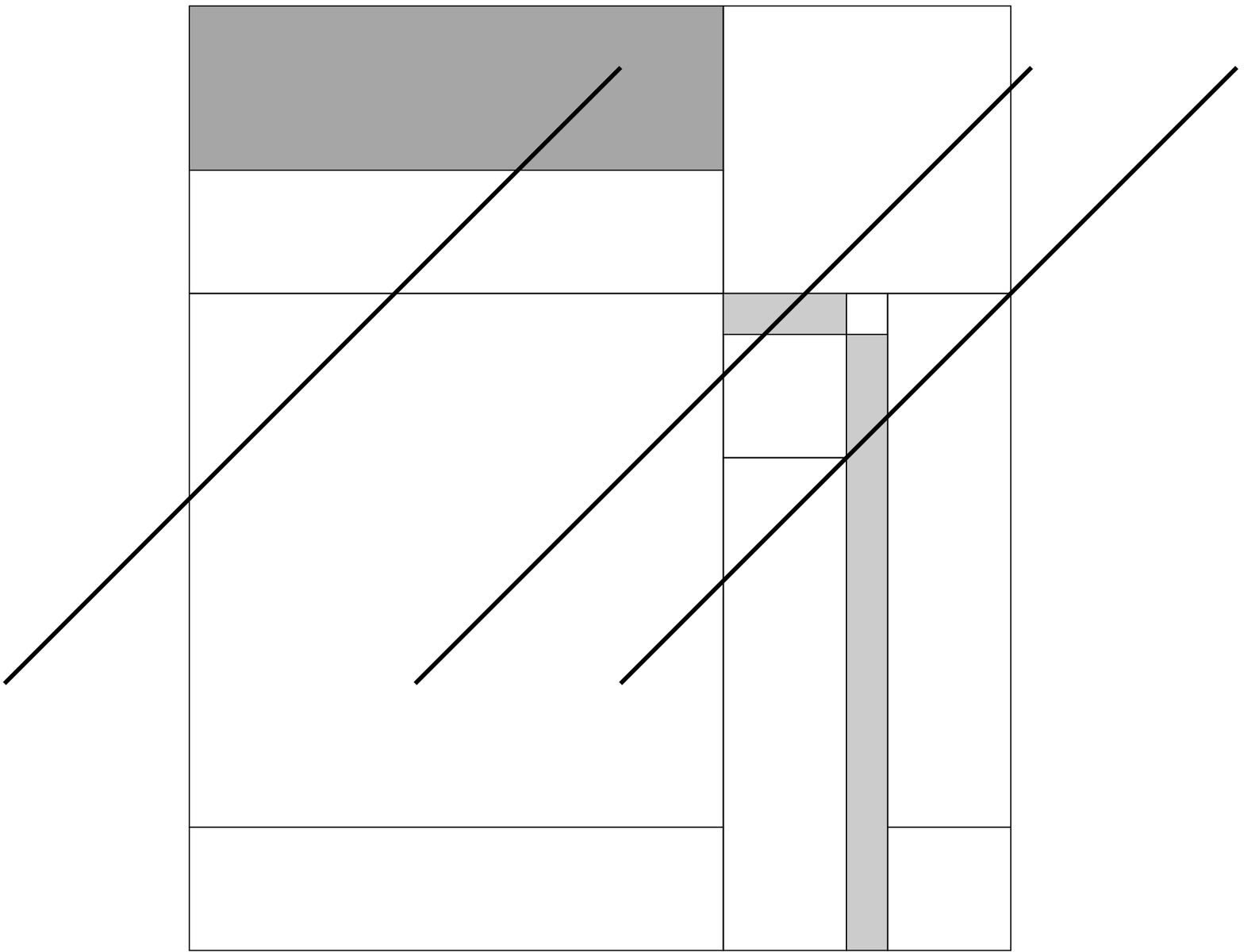}}
\newline
{\bf Figure 17.3:\/} The bottom row of $\Pi_+$ and $\Pi'_+$.
\end{center}

  Note that any diagonal line
intersects at most one of the $3$ relevant rectangles.  Therefore,
what seems like $3$ potential crossings is just $1$.
All in all, there are $2$ potential crossings created
by our perturbation.  This estimate is sharp. The $2$
crossings can happen.

The picture for $\Pi_+$ is easier to analyze.  Recall
from the proof of the Hexagrid Theorem that all
the relevant rectangles were well above the range
of the corresponding vertices.  See Figure 15.2.
Thus, we only have to worry about the emergennce of
new rectangles.  The only new rectangle to emerge within
range is a rectangle labelled $(-1,0)$ that emerges at
the very bottom.  Hence, there is at most $1$ crossing.

\begin{center}
\resizebox{!}{.6in}{\includegraphics{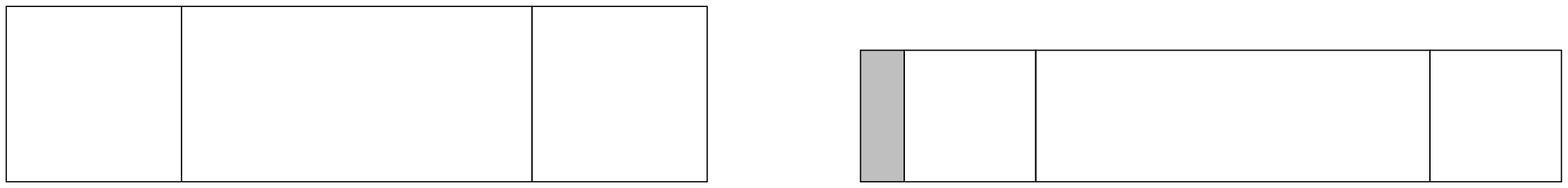}}
\newline
{\bf Figure 17.4:\/} The bottom row of $\Pi_+$ and $\Pi'_+$.
\end{center}

All in all, there are at most $3$ barrier crossings within
one period.  Also,
the number of barrier crossings is even because
every component is a polygon.  Hence there are
exactly $2$ barrier crossings. The major components
do cross the barrier, and hence this accounts
for the $2$ crossings.
\newline
\newline
{\bf Remark:\/}
As in the proof of the Hexagrid Theorem, the pictures
we drew above look somewhat different for other
parameters.  However, the main point of these
pictures is depict the relative changes between
the two pictures.  This works the same for any
parameter.  The reader can see the picture for
any parameter using Billiard King.
\newline

Now we consider the case when $A'<A$.
We can analyze this case just as above,
though the details are somewhat trickier.
The $(+)$ slice still entails only one
crossing.  However, just from looking
at the pictures, we cannot easily rule
out the possibility that the $(-)$ slice
entails $3$ crossings.  This gives us
a bound of $4$ crossings, which is not
quite good enough.  A more subtle analysis
of the picture could get us down to
$2$ crossings, but we prefer to take
a different approach.

We will use symmetry.  Let
$\Lambda_+$ denote the barrier line.
 There is nothing
special about the fact that $\Lambda_+$ lies
above the baseline.  We could consider the
corresponding line $\Lambda_-$ below the baseline.
Here $\Lambda_-$ is parallel to $V$ and contains
\begin{equation}
P_-=\bigg(0,\frac{-p'+q'}{2}\bigg)
\end{equation}
Actually, to get things exactly right, we
think of $\Lambda_+$ and $\Lambda_-$ lying
infinitesimally near, but below, the lines
we have defined.  This, in particular,
$P_-$ lies above $\Lambda_-$.

We compute that
$$
M(P_-)=\bigg(\beta+\frac{1}{q},0,0\bigg)
$$ 
for some $\beta \in \R$.
Thus, by considering $\Lambda_-$ in place of
$\Lambda_+$, we are back in the case we have
already analyzed.  But now we can apply
the rotational symmetry $\iota$ considered
in \S \ref{symmetry}.   Assuming that
$\iota(\Lambda_-)=\Lambda_+$, the result
for $\Lambda_+$ follows from the result
for $\Lambda_-$.

It is not quite true that $\iota(\Lambda_-)=\Lambda_+$.
In fact, $\iota(\Lambda_-)$ is parallel to
$\Lambda_+$ and exactly
$1/q$ vertical units beneath $\Lambda_-$.  Thus,
we have actually proved the Barrier Theorem
for a barrier that is lower by a tiny bit.
This result suffices for all our
purposes.  

To get the stated result right on the nose, we note that
$P_-$ is the only point adversely effected:
$\iota(P_-)$ lies beneath $\Lambda_+$ whereas
$P_-$ lies on $\Lambda_-$.  However, recall that
we consider our lines to be infinitesimally
beneath the lines through integer points.
Thus, as we mentioned above,
$P_-$ lies above $\Lambda_-$.  So, even
though $\iota(\Lambda_-) \not = \Lambda_+$, all
the relevant lattice points lie on the correct
sides.

This completes the proof of the Barrier Theorem.

\newpage

\noindent
{\bf {\huge Part IV\/} \/}
\newline

Here is an overview of this part of the monograph.

\begin{itemize}

\item In \S \ref{ssp} we prove
the Superior Sequence Lemma from \S \ref{pc}.
The analysis here, especially
Lemma \ref{indX}, is central to all our
arguments in this part.  In \S \ref{omegadef}
we introduce a function $\Omega$, closely related
to our sequences, that plays an important
role in subsequent chapters.  We
call $\Omega$ the {\it Diophantine constant\/}.
The reader interested mainly in 
Lemma \ref{weakcopy} can skip everything
in this chapter except \S \ref{omegadef}.

\item In \S \ref{copyproof} we prove
the Diophantine Lemma.  This result
is the source of most of our period copying
results.

\item In \S \ref{strongexist} we prove
Lemma \ref{weakcopy} and Lemma \ref{strongcopy}.
Lemma \ref{strongcopy} is the final ingredient
in the proof of the Erratic Orbits Theorem.
Lemma \ref{weakcopy} is an easier result
that is the final ingredient in our proof of
the Erratic Orbits Theorem for almost all
parameters.  The reader who is satisfied with
the Erratic Orbits Theorem for almost
all parameters can stop reading the monograph
after \S \ref{weakcopyproof}.

\item In \S \ref{decompproof} we prove
the Decomposition Theorem, stated
in \S \ref{statedecomp}.  Our proof
Lemma \ref{strongcopy} relies on one
case of the Decomposition Theorem, namely
Lemma \ref{strong4}, which we prove
in \S \ref{most}. The reader who
is satisfied with the Erratic Orbits Theorem
can stop reading the monograph after
\S \ref{most}.

\end{itemize}

\newpage

\section{Proof of the Superior Sequence Lemma}
\label{ssp}

\subsection{Existence of the Inferior Sequence}
\label{farey1}

We will give a hyperbolic geometry construction
of the inferior sequence.
Our proof is similar to what one does for ordinary
continued fractions.  Our model for the
hyperbolic plane is the upper halfplane
$\H^2 \subset \C$.   The group $SL_2(\R)$ of real
$2 \times 2$ matrices acts isometrically by
linear fractional transformations. The geodesics are 
vertical rays or semicircles centered on $\R$. 
See [{\bf B\/}].

The {\it Farey graph\/} is a tiling of $\H^2$ by ideal
triangles.  
We join $p_1/q_1$ and $p_2/q_2$ by a geodesic
iff $|p_1q_2-p_2q_1|=1$.
The resulting graph divides the hyperbolic plane into
an infinite symmetric union of ideal geodesic triangles.
The Farey graph and the associated triangulation
is one of the most beautiful
pictures in all of mathematics.

We modify the Farey graph by erasing all the lines that connect
even fractions to each other.  The remaining edges
partition $\H^2$ into an infinite union of ideal
squares.  The subgroup $\Gamma_2 \subset SL_2(\R)$, consisting of
matrices congruent to the identity mod $2$, acts in such a
way as to preserve the tiling by idea squares.

We say that a {\it basic square\/} is one of these
squares that has all vertices in the interval $(0,1)$.
Each basic square 
has two opposing vertices that are labelled by positive
odd rationals,
$p_1/q_1$ and $p_2/q_2$.  These odd rationals 
satisfy $|p_1q_2-p_2q_1|=2$.  Ordering so that $q_1<q_2$,
we call $p_1/q_1$ the {\it head\/} of the square
and $p_2/q_2$ the {\it tail\/} of the square.
We draw an arrow
in each odd square that points from the tail to the head.
That is $p_1/q_1 \leftarrow p_2/q_2$.
We call the odd square {\it right biased\/} if the rightmost
vertex is an odd rational, and {\it left biased\/} if
the leftmost vertex is an odd rational.

The general form of a
left biased square is
\begin{equation}
\label{quad1}
\frac{a_1}{b_1}; \hskip 30 pt
\frac{a_1+a_2}{b_1+b_2}; \hskip 30 pt
\frac{a_1+2a_2}{b_1+2b_2}; \hskip 30 pt
\frac{a_2}{b_2}.
\end{equation}
The leftmost vertex in a left-biased square is the head,
and the rightmost vertex in a right-biased square is the head.
One gets the equation for a right-biased square just by
reversing Equation \ref{quad1}.

For an irrational parameter $A$, we simply
drop the vertical line down from $\infty$ to $A$, and record
the sequence of basic squares we encounter.  To form the
inferior sequence, we list the heads of the encountered squares 
and weed out repeaters.  The nesting properties of the squares
guarantees convergence.

\subsection{Structure of the Inferior Sequence}
\label{inferiorstruct}

Now suppose that $\{p_n/q_n\}$ is the inferior
sequence approximating $A$.   Referring to
Equation \ref{induct0}, we write
$A_n=p_n/q_n$ and $(A_n)_{\pm}=(p_n)_{\pm}/(q_n)_{\pm}$. 
We have $(A_n)_-<A_n<(A_n)_+$, and these numbers form
$3$ vertices of an ideal square.  $A_n$ is the tail of
the square.

\begin{lemma}
\label{squeeze}
\label{squeeze2}
The following is true for all indices $m$.
\begin{enumerate}
\item Let $N >m$.  Then $A_{m-1}<A_m$ iff $A_{m-1}<A_N$.
\item If $A_{m-1}<A_m$ then $(q_m)_-=q_{m-1}+(q_m)_+$.
\item If $A_{m-1}>A_m$ then $(q_m)_+=q_{m-1}+(q_m)_-$.
\item Either $A_m<A<(A_m)_+$ or $(A_m)_-<A<A_m$.
\end{enumerate}
\end{lemma}

\startproof
Statement 1 follows from the nesting properties
properties of the ideal squares encountered by the vertical geodesic
$\gamma$ as it converges to $A$. 

 For Statement 2,
note that $A_{m-1}<A_m$ iff these two rationals participate in a
left-biased basic square, which happens iff $(q_m)_+<(q_m)_-$.
By definition $q_{m-1}=|(q_m)_--(q_m)_+|$.  When $(q_m)_+<(q_m)_-$,
we can simply remove the absolute value symbol and solve for
$(q_m)_-$.  Statement 3 is similar.
  
For Statement 4,
we will consider the case when $A_m<A_{m-1}$.  The
other case is similar.  At some point
$\gamma$ encounters the basic square with vertices
$$(A_m)_-<A_m<(A_m)_+<A_{m-1}.$$
If $A_{m+1}<A_m$, then $\gamma$ exits $S$ between
$(A_m)_-$ and $A_m$. So, $(A_m)_-<A<A_m$.  If
$A_{n+1}>A_m$, then $\gamma$ exits $S$ to
the right of $A_m$.  If $\gamma$ exits $S$ to
the right of $(A_m)_+$, then $\gamma$ next encounters a
basic square $S'$ with vertices
$$(A_m)_+<O<E<A_{m-1},$$
where $O$ and $E$ are odd and even rationals.  But then
$A_m$ would not be the term in our sequence after
$A_{m-1}$.  The term after $A_{m-1}$ would lie in
the interval $[O,A_{m-1})$.  This is a contradiction.
\endproof

Let $[x]$ denote the floor of $x$.   Let
$d_n$ be as in Equation \ref{renorm}.  Relatedly,
define 
\begin{equation}
\label{enhanced}
\delta_0=q_1-1; \hskip 30 pt
\delta_n=\bigg[\frac{q_{n+1}}{q_n}\bigg]; 
\hskip 30 pt n=1,2,3...
\end{equation}

Now we come to our main structural result about
the inferior sequence.

\begin{lemma}
\label{indX}
The following is true for any index $m \geq 1$.
\begin{enumerate}
\item If  $A_{m-1}<A_m<A_{m+1}$ then 
\begin{itemize}
\item $\delta_m$ is odd;
\item  $(q_m)_+<(q_m)_-$;
\item  $(q_{m+1})_+=d_mq_m+(q_m)_+$;
\item  $(q_{m+1})_-=(d_m+1)q_m+(q_m)_+.$
\end{itemize}
\item If  $A_{m-1}>A_m<A_{m+1}$ then
\begin{itemize}
\item $\delta_m$ is even;
\item $(q_m)_-<(q_m)_+$;
\item $(q_{m+1})_+=d_mq_m-(q_m)_-$;
\item $(q_{m+1})_-=d_mq_m+(q_m)_+.$
\end{itemize}
\item If $A_{m-1}>A_m>A_{m+1}$  then
\begin{itemize}
\item $\delta_m$ is odd;
\item $(q_m)_-<(q_m)_+$;
\item $(q_{m+1})_+=(d_m+1)q_m+(q_m)_-$;
\item $(q_{m+1})_-=d_mq_m+(q_m)_-.$
\end{itemize}
\item If $A_{m-1}<A_m>A_{m+1}$ then \begin{itemize}
\item $\delta_m$ is even;
\item $(q_m)_+<(q_m)_-$;
\item $(q_{m+1})_+=d_mq_m+(q_m)_-$;
\item $(q_{m+1})_-=d_mq_m-(q_m)_+.$
\end{itemize}
\end{enumerate}
\end{lemma}

\noindent
{\bf Remarks:\/} \newline
(i) 
Here $A_{m-1}<A_m>A_{m+1}$ means that
$A_{m-1}<A_m$ and $A_m>A_{m+1}$. \newline
(ii) There is a basic symmetry
in this result.  We we swap all inequalities,
then the signs $(+)$ and $(-)$ all switch.
This symmetry swaps Cases 1 and 3, and likewise
swaps Cases 2 and 4. \newline
(iii) The same result holds for $p$ in place of $q$.
We used $q$ just for notational convenience.

\startproof
Cases 3 and 4 follow from Cases 1 and 2 by symmetry:
Whatever argument we would give for Cases 1 and 2,
we would just switch all the $(+)$ signs to $(-)$
signs and reverse all the inequalities to prove
the corresponding statement for Cases 3 and 4.
Thus, it suffices to consider Cases 1 and 2.
We will consider Case 1 in detail, and only treat
Case 2 briefly at the end.

In Case 1, the
vertical geodesic $\gamma$
to $A$ passes through the basic square $S$ with vertices
$$
A_{m-1}<(A_m)_-<A_m<(A_m)_+.$$
Since $A_n<A_{m+1}$, the geodesic $\gamma$ next crosses
through the geodesic $\alpha_m$ connecting $A_m$ to $(A_m)_+$.
Following this, $\gamma$ encounters
the basic squares $S'_k$ for $k=0,1,2...$ until it
crosses a geodesic that does not have $A_m$ as
a left endpoint.  By Equation \ref{quad1} and induction,
we get the following list of vertices for the
equare $S'_k$:
\begin{equation}
\label{trap}
\frac{p_m}{q_m}<
\frac{(k+1)p_m+(p_m)_+}{(k+1)q_m+(q_m)_+}< 
\frac{(2k+1)p_m+2(p_m)_+}{(2k+1)q_m+2(q_m)_+}< 
\frac{kp_m+(p_m)_+}{kq_m+(q_m)_+}.
\end{equation}
Here $S'_k$ is a left-biased square.  But then
there is some $k$ such that
\begin{equation}
\label{kappa}
\frac{p_{m+1}}{q_{m+1}}=\frac{(2k+1)p_m+2(p_m)_+}{(2k+1)q_m+2(q_m)_+}; \hskip 30 pt
\frac{(p_{m+1})_+}{(q_{m+1})_+}=\frac{kp_m+(p_m)_+}{kq_m+(q_m)_+}
\end{equation}

Since $(q_m)_+<(q_m)_-$, we have
$2(q_m)_+<q_m$.
Since $2(q_m)_+<q_m$, we have
\begin{equation}
\frac{p_{m+1}}{q_{m+1}}-\frac{p_m}{q_m}=\frac{2}{(2k+1)q_m^2+2q_m(q_m)_+}
\in \bigg( \frac{2}{(2k+2)q_m^2},\frac{2}{(2k+1)q_m}\bigg).
\end{equation}
Hence $\delta_m=(2k+1) \equiv 1$ mod $2$.  
Here $k=d_m$.
This takes care of the second implication.
Equation \ref{kappa} gives the formula for
$(q_{m+1})_+$.  Lemma \ref{squeeze2} now gives
the formula for $(q_{m+1})_-$.

In Case 2, the vertical geodesic $\gamma$ again encounters the
basic square $S$.  This time, $\gamma$ exits $S$ through the
geodesic joining $(A_m)_-$ to $A_m$.  This fact follows
from the inequality $A_m>A_{m-1}>(A_m)_-$, a result of
Lemma \ref{squeeze}.  Following this, $\gamma$ encounters
the basic squares $S_k''$ for $k=0,1,2$ until
it crosses a geodesic that does not have $A_m$ as a right
endpoint.  The coordinates for the vertices of $S''_k$ are
just like those in Equation \ref{kappa}, except that
all terms have been reversed and each $(\cdot)_+$ is
switched to $(\cdot)_-$.  The rest of the proof is
similar.
\endproof

\subsection{Existence of the Superior Sequence}
\label{existsup}

The following result completes the proof of the Superior Sequence Lemma.

\begin{lemma}
\label{bigkappa}
$d_m \geq 1$ infinitely often.
\end{lemma}

\startproof
We can sort the indices of our sequence into $4$ types, depending
on which case holds in Lemma \ref{indX}.
If this lemma is false, then $n$ eventually has odd type.
But, it is impossible for $n$ to have Type $1$ and for
$n+1$ to have Type $3$.  Hence, eventually $n$ has
constant type, say Type 1.  (The Type 3 case has a similar treatment.)
Looking at the formula in Case 1 of Lemma \ref{indX}, we see that
the sequence $\{(q_n)_+\}$ is eventually constant.  But then
$$r=\lim_{n \to \infty} \frac{(q_n)_+ p_n}{q_n}$$
exists.   Since $(q_n)_+p_n \equiv -1$ mod $q_n$ and
$q_n \to \infty$, we must have $r \in \Z$.  But then
$\lim p_n/q_n \in \Q$, and we have a contradiction.
\endproof

\begin{lemma}
\label{superior}
If $d_m \geq 1$ then
$$\bigg|\frac{p_N}{q_N}-\frac{p_m}{q_m}\bigg|<\frac{2}{d_mq_m^2} \hskip 10 pt
\forall N>m; \hskip 40 pt
\bigg|A-\frac{p_m}{q_m}\bigg| \leq \frac{2}{d_mq_m^2}$$
\end{lemma}

\startproof
The first conclusion implies the second.
We will consider the case when $A_m<A_{m+1}$.
By Lemma \ref{squeeze}, we have
\begin{equation}
\label{appx1}
|A_N-A_m| \leq |(A_{m+1})_+-A_m|=
\frac{1}{q_m(q_{m+1})_+}.
\end{equation}
If $m$ is an index of Type 1, then
\begin{equation}
\label{appx2}
(q_{m+1})_+=d_mq_m+(q_m)_+>d_mq_m.
\end{equation}
if $m$ is an index of Type 2, then
Lemma \ref{indX} 
tell us that
\begin{equation}
\label{appx3}
(q_{m+1})_+=(q_{m+1})_--q_m=
d_mq_m+(q_m)_+-q_m>\bigg(d_m-\frac{1}{2}\bigg)q_m \geq \frac{1}{2}d_mq_m.
\end{equation}
Combining Equations \ref{appx1},\ref{appx2}, and
\ref{appx3} we
get our result.
\endproof

\subsection{The Diophantine Constant}
\label{omegadef}

We have two odd rationals
$A_1=p_1/q_1$ and $A_2=p_2/q_2$.
 We define the real number $a=a(A_1,A_2)$ by the formula
\begin{equation}
\label{dio1}
\bigg|\frac{p_1}{q_1}-\frac{p_2}{q_2}\bigg|=\frac{2}{a q_1^2}.
\end{equation}
We call $(A_1,A_2)$ {\it admissible\/} if $a(A_1,A_2)> 1$.

Define
\begin{equation}
\lambda_1=\frac{(q_1)_+}{q_1} \in (0,1).
\end{equation}
If $A_1<A_2$ we define
\begin{equation}
\label{determineA}
\Omega={\rm floor\/}\bigg(\frac{a}{2}-\lambda_1\bigg)+1+\lambda_1.
\end{equation}
if $A_1>A_2$ we define
\begin{equation}
\label{determineB}
\Omega={\rm floor\/}\bigg(\frac{a}{2}+\lambda_1\bigg)+1-\lambda_1.
\end{equation}

Here we explain the geometric meaning of the
Diophantine Constant.  For ease of
exposition, assume that $A_1<A_2$. Assume
also that $(A_1,A_2)$ is admissible.
Let $\epsilon$ denote
an infinitesimally small negative number.
Consider two infinite
rays $R_1$ and $R_2$ starting at $(0,\epsilon)$.
Let $R_j$ have slope $-A_j$.  Thus, $R_j$ is
contained in the baseline of the arithmetic
graph associated to $A_j$.  Then there is no
lattice point between $R_1$ and $R_2$ whose
first coordinate lies in $(0,\Omega q_1)$.
Compare Lemma \ref{good} from the
next chapter.   
Typically there is
such a lattice point with first coordinate
exactly $\Omega q_1$. (We think that this
is always the case, but we did not try to prove it.)
\newline
\newline
{\bf Remarks:\/} \newline
(i) We have formulated our description in terms
of an infinitesimal number because
e.g. the lattice point $(q_1,-p_1)$
is closer to the origin and lies on the line
of slope $-A_1$ through the origin.  However,
this lattice point lies above both rays
$R_1$ and $R_2$, on account of the infinitesimal
downward push we have given these rays.
This is the same bit of silliness we dealt
with when defining the baseline of the
arithmetic graph.
\newline
(ii) The only fact relevant for Lemma \ref{weakcopy} is that
$a>4$ implies that $\Omega>2$.  The reader who cares mainly
about Lemma \ref{weakcopy} can skip the rest of this chapter.
\newline

\subsection{Structure of the Diophantine Constant}

Let $A=p/q$ be an odd rational.  
We say that $A'$ is a {\it near predecessor\/} of
$A$ if $A'$ precedes $A$ in the inferior sequence,
but does not precede the superior predecessor of
$A$.  The inferior and superior predecessors of $A$
are the two extreme examples of near predecessors of $A$.
Here is nice characterization 
of the Diophantine constant for these pairs of rationals.

\begin{lemma}
\label{keycomp1} 
If $A'$ is a near predecessor of $A$ then the following is true.
\begin{enumerate}
\item If $A'<A$ then $\Omega q'=q'+q_+$.
\item If $A'>A$ then $\Omega q'=q'+q_-$.
\end{enumerate}
\end{lemma}

\startproof
There is a finite chain
\begin{equation}
\label{chain}
A'=A_1 ... \leftarrow ... \leftarrow A_m=A.
\end{equation}
Referring to Equation \ref{renorm},
we have $d_1 \geq 0$ and $d_2=...=d_{m-1}=0$.
 By Lemma \ref{squeeze},
$A_1<A_2$ iff $A'<A$.
We will consider the case when $A_1<A_2$.
The other case is similar. 
Recall that
\begin{equation}
\label{recalldef}
A-A'=\frac{2}{a(q')^2}; \hskip 30 pt
\Omega={\rm floor\/}\bigg(\frac{a}{2}-\lambda\bigg)+1+\lambda;
\hskip 30 pt \lambda=\frac{q'_+}{q'}.
\end{equation}
Hence
\begin{equation}
\Omega q'=q'(N+1)+q_+'; \hskip 30 pt N={\rm floor\/}\bigg(\frac{a}{2}-\lambda\bigg)
\end{equation}

There are two cases to consider, depending on whether
$\delta_1$ is odd or even.  Here $\delta_1$ is as in
Equation \ref{enhanced}.  If $\delta_n$ is odd then we
have Case 1 of Lemma \ref{indX}.   In this case, we will
show below that
$d_1=N$. By
Case 1 of Lemma \ref{indX}, we get
\begin{equation}
(q_2)_+=d_1q_1+(q_1)_+=Nq_1+(q_1)_+.
\end{equation}
If $\delta_n$ is even then we show
below that $d_1=N+1$.  By Case 2 of
Lemma \ref{indX}, we have
\begin{equation}
(q_2)_+=d_1q_1-(q_1)_-=(d_1-1)q_1+q_+=Nq_1+(q_1)_+.
\end{equation}
We get the same result in both cases.

Repeated applications of Lemma \ref{indX},
Case 1, give us
\begin{equation}
q_+=(q_m)_+=...=(q_2)_+=Nq'+q'_+=(N+1)q'-q'+q'_+=\Omega q'-q'.
\end{equation}
Rearranging this gives Statement 1.
\endproof

\begin{lemma}
If $A_1<A_2$ and $\delta_1$ is odd, then $d_1=N$.
\end{lemma}

\startproof
Rearranging the basic definition of $a(A',A)$, and
using $A'=A_1$ and $A=A_m$ in equation \ref{chain}, we have
$$\frac{a}{2}=\frac{1}{q_1^2|A_1-A_m|}.$$
By Lemma \ref{squeeze} and monotonicity, we have
\begin{equation}
\frac{1}{q_1^2|A_1-(A_2)_+|}<\frac{a}{2}<\frac{1}{q_1^2|A_1-A_2|}.
\end{equation}
After some basic algebra, we get
\begin{equation}
\label{change1}
d_1+\lambda_1=^*\frac{(q_2)_+}{q_1}<\frac{a}{2}<\frac{q_2}{2q_1}.
\end{equation}
The starred inequality is Case 1 of Lemma \ref{indX}.
The lower bound
gives us
\begin{equation}
d_1<\frac{a}{2}-\lambda_1
\end{equation}
Here $\lambda_1$ is the same as $\lambda$ in 
Equation \ref{recalldef}.  Since $d_1 \in \Z$, we get
$d_1 \leq N$.
On the other hand, the upper bound gives us
\begin{equation}
\label{change2}
N={\rm floor\/}\bigg(\frac{a}{2}-\lambda_1\bigg) \leq
{\rm floor\/}\bigg(\frac{q_2}{2q_1}-\lambda_1\bigg) \leq d_1.
\end{equation}
In short, $N \leq d_1$.  Combining the two halves gives $N=d_1$.
\endproof

\begin{lemma}
If $A_1<A_2$ and $\delta_1$ is even, then $d_1=N+1$.
\end{lemma}

\startproof
  The proof is very similar to what we did in the other case.
Here we mention the $2$ changes. The
first change is that $(d_1-1)+\lambda_1$ occurs on the left
hand side of Equation \ref{change1}, by Case 2 of Lemma \ref{indX}.
This gives us $d_1 \leq N+1$.  The second change occurs on
the right hand side of Equation \ref{change2}.
By Case 2 of Lemma \ref{indX}, we know that
${\rm floor\/}(q_2/q_1)$ is even.  Hence
$q_2/(2q_1)$ has fractional part less than $1/2$.   But,
also by Case 2 of Lemma \ref{indX}, $\lambda_1$ has
fractional part greater than $1/2$.  Hence
$$
{\rm floor\/}\bigg(\frac{q_1}{2q_1}-\lambda_1\bigg)=
{\rm floor\/}\bigg(\frac{q_1}{2q_1}\bigg)-1 \leq d_1-1.$$
This gives us the bound $N \leq d_1-1$, or
$N+1 \leq d_1$.  Putting the two halves together, we
get $d_1=N+1$.
\endproof

\newpage

\section{The Diophantine Lemma}
\label{copyproof}

\subsection{Three Linear Functionals}
\label{3lf}

Let $p/q$ be an odd rational. 

Consider the following linear functionals.

\begin{equation}
F(m,n)=\bigg(\frac{p}{q},1\bigg) \cdot (m,n)
\end{equation}

\begin{equation}
\label{defG}
G(m,n)=\bigg(\frac{q-p}{p+q},\frac{-2q}{p+q}\bigg) \cdot (m,n).
\end{equation}

\begin{equation}
\label{defH}
H(m,n)=\bigg(\frac{-p^2+4pq\!+\!q^2}{(p+q)^2},\frac{2q(q-p)}{(p+q)^2}\bigg) \cdot (m,n).
\end{equation}

We have $F=(1/2)M$, where $M$ is the fundamental map from
Equation \ref{funm}.  We can understand $G$ and $H$ by evaluating
them on a basis:
\begin{equation}
\label{affine1}
H(V)=G(V)=q; \hskip 30 pt
H(W)=-G(W)=\frac{q^2}{p+q}.
\end{equation}
Here $V=(q,-p)$ and $W$ are the vectors from
Equation \ref{boxvectors}.  We can also understand $G$ by evaluating
on a simpler basis.
\begin{equation}
G(q,-p)=q; \hskip 30 pt G(-1,-1)=1.
\end{equation}

We can also (further) relate $G$ and $H$ to the hexagrid from
\S 3.
A direct calculation establishes the following result.

\begin{lemma}
\label{Hcalc}
The fibers of $G$ are parallel to the top left
edge of the arithmetic kite.   The fibers of $H$ are parallel to the top right
edge of the arithmetic kite.
 Also $\|\nabla G\| \leq 3$ and $\|\nabla H\| \leq 3$.
\end{lemma}
Here $\nabla$ is the gradient.

Given any interval $I$, define
\begin{equation}
\label{deltadefined}
\Delta(I)=\{(m,n)|\ G(m,n), H(m,n) \in I\} \cap \{(m,n)|\ F(m,n) \geq 0\}.
\end{equation}
This set is a triangle whose bottom edge is the baseline
of $\Gamma(p/q)$.  

\subsection{The Main Result}
\label{diotheorem}

\begin{lemma}[Diophantine]
Let $(A_1,A_2)$ be an admissible pair of odd rationals.
\begin{enumerate}
\item If $A_1<A_2$ let  $I=[-q_1+2,\Omega q_1-2]$.
\item If $A_1>A_2$ let  $I=[-\Omega q_1+2,q_1-2]$.
\end{enumerate}
Then $\widehat \Gamma_1$ and $\widehat \Gamma_2$ agree on $\Delta_1(I) \cup \Delta_2(I)$.
\end{lemma}

\begin{center}
\resizebox{!}{4.6in}{\includegraphics{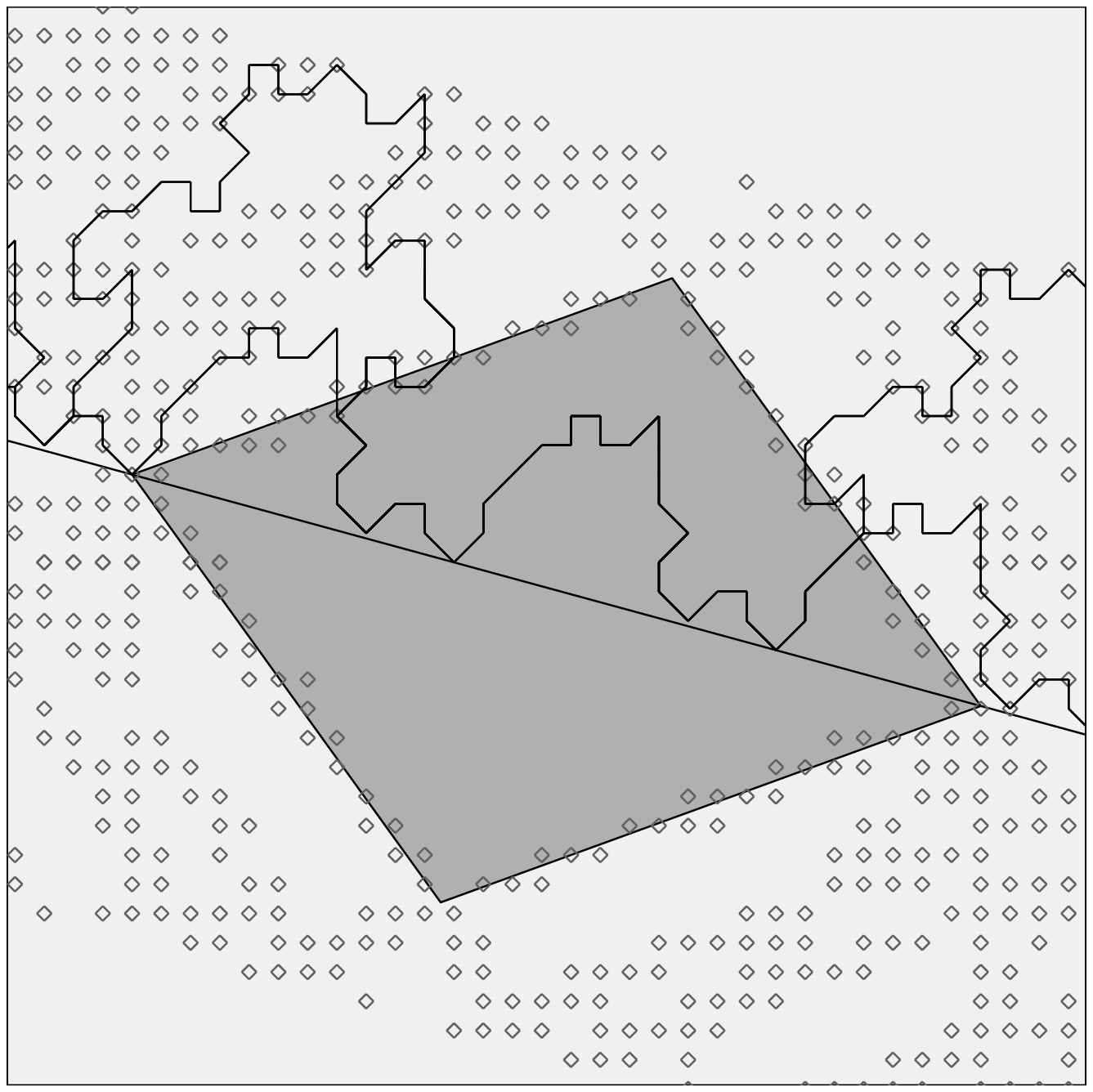}}
\newline
{\bf Figure 19.1:\/} The Diophantine Lemma in action.
\end{center}

Figure 19.1 illustrates our result
for $A_1=3/11$ and $A_2=7/25$.  The portion of the
shaded parallelogram above the baseline is
 $\Delta_1(-q_1,\Omega q_1)$, a set slightly larger than
$\Delta_1(I)$.  The sets $\Delta_1(I)$ and
$\Delta_2(I)$ are almost identical.
\newline
\newline
{\bf Remarks:\/} \newline
(i)  The Diophantine Lemma also works for
points below the baseline, but for technical reasons
we ignore these points.
We plot the points near $\Delta$ where
$\widehat \Gamma_1$ and $\widehat \Gamma_2$
disagree.   
Starting from $(0,0)$ and tracing $\Gamma_1$ and $\Gamma_2$
in either direction, we get 
agreement until we nearly hit the edges of $\Delta$.
\newline
(ii) The Diophantine Lemma is quite nearly sharp.  We think that
the sharp version runs as follows.  The two arithmetic
graphs agree at any point in $\Delta_1(-q_1,\Omega q_1)$ that
is not adjacent to a point that lies outside of $\Delta_1(-q_1,\Omega q_1)$.
One can see this structure plotting pictures on
Billiard king.
\newline
(iii) The Diophantine Lemma is defined in terms of somewhat
complicated formulas, but the domains involved have
simple geometric
descriptions.   Consider the two triangles
$\Delta_1(-q_1,\Omega q_1)$ and
$\Delta_2(-q_1,\Omega q_1)$.  The bases of these two
nearly identical triangles have no lattice points between
them.  The triangles are then constructed from the bases
by extending lines parallel to the top edges of the
arithmetic kites.  Since $A_1$ and $A_2$ are nearby
rational parameters, the two kites have about the
same shape, and so do the two triangles. \newline
\newline

Here we outline the proof of the Diophantine Lemma.   
We will establish the case when $A_1<A_2$.  The
other case has a nearly identical proof.

We say that an integer $\mu$ is {\it good\/} if 
$\mu A_1$ and $\mu A_2$ have the same floor.  
We call $\mu$ $1$-{\it good\/} if
$\mu+\epsilon$ is good for all $\epsilon \in \{-1,0,1\}$.
We can subject a lattice point $(m,n)$ to the
reduction algorithm from \S \ref{master1}.
For $\theta \in \{1,2\}$, we perform the algorithm
relative to the parameter $A_{\theta}$. This
produces integers $X_{\theta}$ and $Y_{\theta}$ and $Z_{\theta}$. 
Below we prove the following result.

\begin{lemma}[Agreement]
\label{agree}
Suppose, for at least one choice of
$\theta \in \{1,2\}$, that
$m$ and $m-X_\theta$ and $m-Y_\theta$ and
$m+Y_\theta-X_\theta$ are all $1$-good.  Then
$\widehat \Gamma_1$ and $\widehat \Gamma_2$ agree at $(m,n)$.
\end{lemma}

Next, we give a criterion for an integer to be good.

\begin{lemma}[Goodness]
\label{good}
If $\mu \in (-q_1,\Omega q_1) \cap \Z$ then $\mu$ is a good integer.
\end{lemma}

Finally, we show that $(m,n) \in \Delta_1(I) \cup \Delta_2(I)$
implies that the integers in the Agreement Lemma satisfy the
criterion in the Goodness Lemma.  This completes the proof.

\subsection{Proof of the Agreement Lemma}

In our technical lemmas, we will use integers
$\mu$ and $\nu$ roughly in place of $m$ and $n$.
Sometimes $\mu$ and $\nu$ will take on values
other than $m$ and $n$, however.

\begin{lemma}
\label{floor2}
Let $\mu,\nu,N_j \in \Z$ and
$$N_j={\rm floor\/}\bigg(\frac{\mu A_j + \nu}{1+A_j}\bigg).$$
Suppose, for at least one choice of $\theta \in \{1,2\}$
that both $\mu-N_{\theta}$ and $\mu-N_{\theta}+1$ are good.
Then $N_1=N_2$.
\end{lemma}

\startproof
For the sake of contradiction, assume w.l.o.g. that $N_1<N_2$. Then
$$\mu A_1+\nu<N_2(A_1+1); \hskip 30 pt
(\mu-N_2)A_1<N_2-\nu$$
$$N_2(A_2+1) \leq \mu A_2+\nu; \hskip 30 pt
N_2-\nu \leq (\mu-N_2)A_2.$$
The first equation implies the second in each case.
The second items imply that $\mu-N_2$ is not good.
On the other hand, we have
$$\mu A_1 +\nu <(N_1+1)(A_1+1) ; \hskip 30 pt
A_1(m-N_1+1)<N_1+1-n.$$
$$(N_1+1)(1+A_2) \leq \mu A_2+\nu; \hskip 30 pt
A_2(m-N_2+1) \geq N_1-1+n.$$
The first equation implies the second in each case.
The second items imply that $\mu-N_1+1$ is not good.
Now we have a contradiction.
\endproof

\begin{corollary}
Let $(X_{\theta},Y_{\theta},Z_{\theta})$ be
as in the Agreement Lemma.
Under the hypotheses of the Agreement Lemma, we have
$(X_1,Y_1,Z_1)=(X_2,Y_2,Z_2)$.
\end{corollary}

\startproof
We go through the reduction algorithm.
First we deal with $(-)$ case.

\begin{enumerate}

\item Let $z_j=A_jm+n$.
\item   Let $Z_j={\rm floor\/}(z_j)$.  Since $m$ is good,
we have $Z_1=Z_2$. Call this common integer $Z$.

\item $y_j=z_j+Z_j=z_j+Z$.  Hence
$y_j=mA_j+n'$ for some $n' \in \Z$.

\item
Recall that $Y_j={\rm floor\/}(y_j/(1+A))$.
To see that $Y_1=Y_2$ we apply 
Lemma \ref{floor2} to
$(\mu,\nu,N_j)=(m,n',Y_j).$
Here we use the fact that
$m-Y_{\theta}$ and $m-Y_{\theta}+1$ are good.
We set $Y=Y_1=Y_2$.

\item Let $x_j=y_j-Y(1-A_j)-1$.  Hence $x_j=(m+Y)A_j+n''$.

\item Recall that
$X_j={\rm floor\/}(x_j/(1+A))$.  To see
that $X_1=X_2$ we apply 
Lemma \ref{floor2} to 
$(\mu,\nu,N_j)=(m+Y,n'',X_j).$
Here we use the fact that
$m+Y-X_{\theta}$ and $m+Y-X_{\theta}+1$ are
good integers.  We set $X=X_1=X_2$.
\end{enumerate}

Now we deal with the $(+)$ case.   The only difference is
that
$$y_j=z_j+Z+1.$$   
This time we have
$y_j=mA_j+n'+1$, and the argument works exactly
the same way.  We apply Lemma \ref{floor2} to
$(\mu,\nu,N_j)$ to $(m,n'+1,Y_j)$.
\endproof

In the next result, all quantities except $A_1$ and $A_2$
are integers.

\begin{lemma}
\label{floor3}
If $\mu-dN-\epsilon_1$ is good, then the statement
$$(\mu A_j + \nu)-N(dA_j+1)<\epsilon_1 A_j+\epsilon_2$$
is true or false independent of $j=1,2$.
\end{lemma}

\startproof
Assume w.l.o.g. that
the statement is true for $j=1$ and false for $j=2$. Then
$$(\mu-dN-\epsilon_1)A_1<\epsilon_2+N-\nu \leq (\mu-dN-\epsilon_1)A_2,$$
a contradiction.
\endproof

Now we finish the proof of the Agreement Lemma.
Let $M_+$ and $M_-$ be as in \S \ref{mptnotation}.
By the Master Picture Theorem, it suffices to
show that the two images $M_+(m,n)$ and $M_-(m,n)$ land in
the same polyhedra for both $A_1$ and $A_2$.  We have already seen
that the basic integers $(X,Y,Z)$ are the same relative
to both parameters. It remains to locate the relevant
points inside our tori $\R^3/\Lambda_1$ and
$\R^3/\Lambda_2$.
The polyhedra of interest to us are
cut out by the following partitions.

\begin{itemize}
\item
$\cal Z$, the union 
$\{z=0\} \cup \{z=A\} \cup \{z=1-A\} \cup \{z=1\}$.
\item $\cal Y$, the union
$\{y=0\} \cup \{y=A\} \cup \{y=1\} \cup \{y=1+A\}$.
\item $\cal X$, the union
$\{x=0\} \cup \{x=A\} \cup \{x=1\} \cup \{x=1+A\}$.
\item $\cal T$, the union
$\{x+y-z=A + j\}$ for $j=-2,1,0,2,1$.
\end{itemize}
Letting $\cal S$ stand for one of these partitions,
we say that $\cal S$ is {\it good\/} if 
the points $M_+(m,n)$ land in the same component
of $R_+-\cal S$ for both parameters $A_1$ and $A_2$,
and likewise 
the points $M_-(m,n)$ land in the same component
of $R_--\cal S$ for both parameters $A_1$ and $A_2$.
Here $R_{\pm}=\R^3/\Lambda$, the domain of the
maps $M_{\pm}$.  This domain depends on the
parameter.  By the Master Picture Theorem,
$\Gamma_1$ and $\Gamma_2$ agree at $(m,n)$
provided all the partitions are good.
The proof works the same
for the $(+)$ and the $(-)$ case.
\begin{itemize}
\item
For $\cal Z$,
we apply Lemma \ref{floor3} to
$(\mu,\nu,d,N)=(m,n,0,Z)$ to show that the statement
$z_j-Z<\epsilon_1A_j+\epsilon_2$ is true independent
of $j$, for $\epsilon_1 \in \{-1,0,1\}$ and
$\epsilon_2 \in \{0,1\}$.  The relevant
good integers are $m-1$ and $m$ and $m+1$.
\item
For $\cal Y$, 
we apply Lemma \ref{floor3} to
$(\mu,\nu,d,N)=(m,n',1,Y)$ to show that the statement
$z_j-Z<\epsilon_1A_j+\epsilon_2$ is true independent
of $j$, for $\epsilon_1 \in \{0,1\}$ and
$\epsilon_2 \in \{0,1\}$.  The relevant good
integers are $m-Y$ and $m-Y-1$.
\item
For $\cal X$, 
we apply Lemma \ref{floor3} to
$(\mu,\nu,d,N)=(m+Y,n'',1,X).$  The
relevant good integers are $m+Y-X$ and $m+Y-X-1$.
\item
For $\cal T$, we define
$$\sigma_j=(x_j-X(1+A_j))+(y_j-Y(1+A_j))-(z_j-Z).$$
We have $\sigma_j=(m-X)A_j+n'''$ for some $n''' \in \Z$.
Let $h \in \Z$ be arbitrary.  To see that the
statement $\sigma_j<A_j+h$ is true independent of $j$ we
apply  Lemma \ref{floor3} to 
$(\mu,\nu,d,N)=(m-X,n''',1,0).$
The relevant good integer is $m-X-1$.
\end{itemize}

\noindent
{\bf Remark:\/} Our proof does not use
the fact that $m-X+1$ is a good integer.
This technical detail is relevant for
Lemma \ref{YX}.

\subsection{Proof of the Goodness Lemma}

We prove the Goodness Lemma in two steps.
The first step takes care of the lower bound
and the second step takes care of the 
upper bound.  Before we start our proof,
we note that the Goodness Lemma is the result
that justifies out claims, made in \S \ref{omegadef},
 about the geometric
meaning of the Diophantine constant $\Omega$.

\begin{lemma}
\label{floor1}
If $\mu \in (-q_1,0) $, then
$\mu$ is good.
\end{lemma}

\startproof
Since $q_1$ is odd, we have unique
integers $j$ and $M$ such that
\begin{equation}
\label{quantum}
\mu A_1= M + \frac{j}{q_1}; \hskip 30 pt
|j|<\frac{q_1}{2}.
\end{equation}
By hypotheses, $a>1$.  Hence
\begin{equation}
|A_2-A_1|<2/q_1^2
\end{equation}
 in all cases.
If this result is false, then there is some integer $N$ such that
\begin{equation}
\label{hypo1}
\mu A_2<N \leq \mu A_1.
\end{equation}
Referring to Equation \ref{quantum}, we have
\begin{equation}
\label{longchain2}
\frac{|j|}{q_1}<\mu A_1-N \leq \mu A_1-\mu A_2<\frac{2|\mu|}{q_1^2}<\frac{2}{q_1}.
\end{equation}
If $j=0$ then $q_1$ divides $\mu$, which is impossible.  Hence
$|j|=1$.  If $j=-1$ then $\mu A_1$ is $1/q_1$ less than an integer.
Hence $\mu A_1-N \geq (q_1-1)/q_1$. This is false, so we must have $j=1$.

From the definition of $\lambda_1$, we have the following implication.
\begin{equation}
\label{congruence}
\mu \in (-q_1,0) \hskip 15 pt
{\rm and\/} \hskip 15 pt \mu p_1 \equiv 1 \hskip 5 pt
{\rm mod\/} \hskip 5 pt q_1 \hskip 30 pt
\Longrightarrow \hskip 30 pt \mu=-\lambda_1q_1.
\end{equation}

Equation \ref{quantum} implies
$$\frac{\mu p_1}{q_1}-\frac{1}{q_1} \in \Z.$$
But then
$\mu p_1 \equiv 1$ mod $q_1$. Equation
\ref{congruence} now tells us that $\mu=-\lambda_1q_1$.  Hence
$|\mu|<q_1/2$.  But now Equation \ref{longchain2} is twice as
strong and gives $|j|=0$.  
This is a contradiction.
\endproof

\begin{lemma}
\label{floor0}
If $\mu \in (0,\Omega q_1)$ then
$\mu$ is good.
\end{lemma}

\startproof
We observe that $\Omega<a$, by Equation \ref{determineA}.
If this result is false, then there is some integer $N$ such that
$\mu A_1<N \leq \mu A_2$.
If $\mu A_2=N$. Then $q_2$ divides $\mu$.   But then,
$$\mu \geq q_2 \geq aq_1>\Omega q_1.$$  This
is a contradiction.  Hence
\begin{equation}
\label{gap}
\mu A_1<N< \mu A_2.
\end{equation}
Referring to Equation \ref{quantum}, we have
\begin{equation}
\label{longchain0}
\frac{|j|}{q_1}\leq N-\mu A_1< \mu(A_2-A_1)=\frac{2\mu}{a q_1^2} <
\frac{2}{q_1}.
\end{equation}
Suppose that $j \in \{0,1\}$ in Equation \ref{quantum}.  Then
$$1-\frac{1}{q_1} \leq N-\mu A_1 \leq \mu A_2-\mu A_1< \frac{1}{q_1},$$
a contradiction. Hence $j=-1$.  Hence $\mu> aq_1/2$.

Since $j=-1$,  Equation \ref{quantum} now tells us
that $\mu p_1 +1 \equiv 0$ mod $q_1$.  But then
\begin{equation}
\mu=kq_1+(q_1)_+,
\end{equation}
for some $k \in \Z$.
On the other hand, from Equation \ref{determineA} and the fact that $\mu<\Omega q_1$,
we have
\begin{equation}
\mu<k' q_1+(q_1)_+; \hskip 20 pt
k'=\Big({\rm floor\/}(a/2-\lambda_1)+1\Big).
\end{equation}
Comparing the last two equations, we have $k \leq k'-1$.  Hence
\begin{equation}
k \leq \Big({\rm floor\/}(a/2-\lambda_1)\Big).
\end{equation}
Therefore
$$\mu \leq  \Big({\rm floor\/}(a/2-\lambda_1)\Big) q_1 + \lambda_1 q_1 \leq
a q_1/2.$$
But we have already shown that $\mu>aq_1/2$.  This is a contradiction.
\endproof

\subsection{The End of the Proof}

We will assume that $(m,n) \in \Delta_{\theta}(I)$, for
one of the two choices $\theta \in \{1,2\}$.
Here $I$ is as in the Diophantine Lemma.  Our proof
works the same for $\theta=1$ and $\theta=2$.
We set $p=p_{\theta}$ and $q=q_{\theta}$, etc.

We will show that all the integers that arise in
our proof of Lemma \ref{agree} lie in
$(-q_1,\Omega q_1)$.  These integers have the form
$N+\epsilon$ for $\epsilon \in \{-1,0,1\}$.  We
will show, for all relevant integers (except one),
that $N \in J:=(-q_1+1,\Omega q_1-1)$.
For the exceptional case, see the remark
after Lemma \ref{YX}.

\begin{lemma}
$m \in J$.
\end{lemma}

\startproof
We have $z=Am+n \geq 0$.  We compute
\begin{equation}
\label{algG}
-q_1+2 \leq G(m,n)=m - \frac{2z}{1+A} \leq m.
\end{equation}
\begin{equation}
\label{algH}
\Omega q_1-2 \geq H(m,n)=m+  \frac{2z(1-A)}{(1+A)^2} \geq m.
\end{equation}
These inequalities establish that $m\in J$.  
\endproof

\begin{lemma}
\label{efftrip1}
$m-Y \in J$.
\end{lemma}

\startproof
We have $Y \geq 0$.  Hence $m-Y \leq m \leq \Omega q_1-2$.
We just need the lower bound.
worry about the lower bound on $m-Y$.
We first deal with the algorithm in \S \ref{master1} for
the $(-)$ case.
Let $G=G(m,n)$.
We have $y=z+Z \leq 2z$.
By definition of $Y$, we have
\begin{equation}
\label{Yest}
Y \leq \frac{y}{1+A} \leq \frac{2z}{1+A}; \hskip 30 pt
Y<\frac{2z}{1+A}.
\end{equation}
At least one of the first two inequalities is sharp. This
gives us the second inequality.
Now we know that
\begin{equation}
m-Y> m-\frac{2z}{1+A}=G \geq -q_1+2.
\end{equation}
The last equality comes from Equation \ref{algG}.
In the $(+)$ case we add $1$ to $Y$. giving
us $m-Y> -q_1+1$.
\endproof

\begin{lemma}
\label{YX}
$m-X \in J \cup \{\Omega q_1 - 1\}$.
\end{lemma}

\startproof
The condition that $F(m,n) \geq 0$ implies that $y \geq Y \geq 0$.
Hence
\begin{equation}
\label{xy}
x=y-Y(1-A)-1 \in [-1,y-1].
\end{equation}
 Hence $X \in [-1,Y-1]$.
Hence $$m-X \in [m-Y+1,m+1] \subset J \cup \{\Omega q_1-1\},$$
by the two previous results.
\endproof

\noindent
{\bf Remark:\/}  As we remarked at the end of the proof
of Lemma \ref{agree}, 
the integer $m-X+1$ does not arise in our
proof of Lemma \ref{agree}.  The relevant
integers $m-X$ and $m-X-1$ are good, by
the result above.

\begin{lemma}
$m+Y-X\in J$.
\end{lemma}

\startproof
Our proof works the same in the $(+)$ and $(-)$ cases.
Lemma \ref{YX} gives us $Y-X \geq 0$.
Hence 
$$m+Y-X \geq m>-q_1+1.$$
This takes care of the lower bound.  Now we treat the upper bound.
We have
$$Y = {\rm floor\/}\bigg(\frac{y}{1+A}\bigg) \leq \frac{y}{1+A}; \hskip 30 pt
1+X ={\rm floor\/} \bigg(1+\frac{x}{1+A}\bigg) \geq \frac{x}{1+A}.$$
Hence
$$
Y-X-1 \leq \frac{y-x}{1+A} =^1
Y \frac{1-A}{1+A} + \frac{1}{1+A}<^* $$
$$
2z\frac{1-A}{(1+A)^2}+\frac{1}{1+A}=^2
H - m + \frac{1}{1+A}<H-m+1.$$
The first equality comes from Equation \ref{xy}.
The second equality comes from Equation \ref{algH}.
The starred inequality comes from the upper bound
in Equation \ref{Yest}.
Adding $m$ to both sides, we get
$$m+Y-X<H+1 \leq \Omega q_1-1.$$
This completes the proof.
\endproof

\newpage

\section{Existence of Strong Sequences}
\label{strongexist}

\subsection{Proof of Lemma \ref{weakcopy}}
\label{weakcopyproof}

We will prove the result when $A_1<A_2$.  The other case is similar.
By hypotheses, we have $a(A_1,A_2)>4$.  From Equation \ref{determineA},
we get $\Omega>2$.  Let $R_1=R(A_1)$ be the parallelogram from the
room lemma.  Let
\begin{equation}
u=W_1; \hskip 30 pt
w=V_1+W_1
\end{equation}
denote the top left and right vertices of $R_1$.
We compute

\begin{equation}
\label{strong11}
G_1(u)=-\frac{q_1^2}{p_1+q_1}>-q_1+2;
\hskip 30 pt
H_1(w)=\frac{q_1^2}{p_1+q_1}+q_1<\Omega q_1-2.
\end{equation}
The inequalities hold once $p_1$ is sufficiently large.
Given the description of the fibers of $G$, we have
\begin{equation}
G(u) \leq G(v) \leq H(v) \leq H(w); \hskip 30 pt
\forall v \in R_1.
\end{equation}
The middle inequality uses the fact that $F(v)\geq 0$.
 In short, we have made the
extremal calculations.  The extremal calculation shows that
$v \in \Delta_1(I)$ for all $v \in R_1$.
The Diophantine Lemma now shows that
$\Gamma_1$ and $\Gamma_2$ agree in $R_1$.

When $v$ lies in the bottom edge of $R_1$ we have
\begin{equation}
G_1(v), H_1(v) \in [0,q_1].
\end{equation}
Given our gradient bounds
$\|\nabla G_1\| \leq 3$ and
$\|\nabla H_1 \| \leq 3$, we see that
\begin{equation}
G_1(v), H_1(v) \in [-q_1+2,\Omega q_1+2]
\end{equation}
provided that $v$ is within $q_1/4$ from 
the bottom edge of $R_1$. 
Hence $\Gamma_1$ and $\Gamma_2$ agree in the
$q_1/4$ neighborhood of the bottom edge of
$R_1$.  

By the Room Lemma, $\Gamma_1^1 \subset R_1$. Hence
$\Gamma_1^1 \subset \Gamma_2$.  Our calculation
involving the bottom edge of $R_1$ shows that
$\Gamma_1^{1+\epsilon} \subset \Gamma_2$ for
$\epsilon=1/4$.
Since the right endpoint of $\Gamma_2^1$ is far to
the right of any point on $\Gamma_1^{1+\epsilon}$, we
have $\Gamma_1^{1+\epsilon} \subset \Gamma_2^1$, as
desired.
\newline
\newline
\noindent
{\bf Remark:\/} We proved Lemma \ref{weakcopy} for
$\epsilon=1/4$ rather than $\epsilon=1/8$, which
is what we originally claimed.  We don't care about
the value of $\epsilon$, as long as it is positive.

\subsection{Proof of Lemma \ref{strongcopy}}

The Decomposition Theorem is stated in \S \ref{statedecomp}.
Our proof requires a limited version of this result.
Define the complexity of an odd rational to be the number of
terms preceding it in the superior sequence.

\begin{lemma}
\label{limiteddecomp}
The Decomposition Theorem holds for all $A$ having sufficiently
large complexity.
\end{lemma}

\startproof
This is an immediate consequence of Lemma \ref{strong4},
proved in \S \ref{most}.
\endproof

Let $A$ be any irrational parameter.  Let
$\{p_n/q_n\}$ denote the superior sequence
associated to $A$.  Let $\cal S$ be a monotone
subsequence of the superior sequence.
We will treat the case when $\cal S$ is monotone
increasing.  If necessary, we cut off the
first few terms of $\cal S$ so that
Lemma \ref{limiteddecomp} holds for all terms.

For any odd rational $p/q$, let
$R^*(p/q)$ denote the rectangle with
vertices
\begin{equation}
-\frac{V}{2}; \hskip 30 pt -\frac{V+W}{2}; \hskip 30 pt \frac{V+W}{2}; \hskip 30 pt \frac{V}{2}.
\end{equation}
Here $V$ and $W$ are as in Equation \ref{boxvectors}.
The parallelogram $R^*$ is just as wide as $R$ but half
as tall.  Also, the bottom edge of $R^*$ is centered
on the origin.

\begin{lemma}
\label{strong1}
If $A_1 \leftarrow A_2$ and 
$p_1$ is sufficienly large, then
$\Gamma_1$ and $\Gamma_2$ agree in $A_1$.  Moreover,
$\Gamma_1$ and $\Gamma_2$ agree 
in the $q_1/8$ neighborhood of the bottom
edge of $R_1^*$.
\end{lemma}

\startproof
The proof works the same way regardless of the sign
of $A_1-A_2$.  The main point is that $\Omega>1$.
Note that $(A_1,A_2)$ is admissible.
We use the linear functionals $G_1$ and $H_1$ associated to
$A_1$.
Let 
$$u=\frac{-V+W}{2}; \hskip 30 pt
w=\frac{V+W}{2}
$$
denote the top left and right vertices of
$R_1^*$ respectively.
We compute
\begin{equation}
-G_1(u)=H_1(w)=\frac{q_1(p_1+2q_1)}{2p_1+2q_1)}<q_1-2<\Omega q_1-2.
\end{equation}
The same argument as in Lemma \ref{weakcopy} now finishes
the proof.
\endproof

\begin{lemma}
\label{strong2}
Suppose that $A_1<A_2$ and $A_1$ is the superior predecessor of $A_2$.
If $A_1$ has sufficiently large complexity, then
$\Gamma_1^{1+\epsilon} \subset \Gamma_2^1$.
\end{lemma}

\startproof
If $\Omega>2$, we have the same proof as in Lemma \ref{weakcopy}.
Equation \ref{determineA} does not allow $\Omega=2$. We just
need to consider the case $\Omega<2$.  By
Equation \ref{determineA}, we must have
${\rm floor\/}(a/2-\lambda_1)=0$.  Since $a>1$, we must have
$\lambda_1>1/2$.  Since $\lambda_1=(q_1)_+/q_1$ and
$q=q_++q_-$, we must have
have $(q_1)_-<(q_1)_+$.   This seemingly minor fact is crucial
to our argument.

Now we really need to use Lemma \ref{limiteddecomp}.
Let $R(A_1)$ denote the parallelogram from the Room Lemma.
In contrast, let $R_1(A_1)$ and $R_2(A_2)$ denote the
smaller parallelograms from the Decomposition Theorem.
Since $(q_1)_-<(q_-)_+$, we see that $R_2(A_1)$ lies to the
left of $R_1(A_1)$.   By the Decomposition Theorem,
\begin{equation}
\label{contain3}
\Gamma_1 \cap R(A_1) \subset R_2(A_1) \cup (R_1(A_1)+V_1)
\end{equation}
Figure 20.1 shows a schematic picture.

\begin{center}
\resizebox{!}{1.5in}{\includegraphics{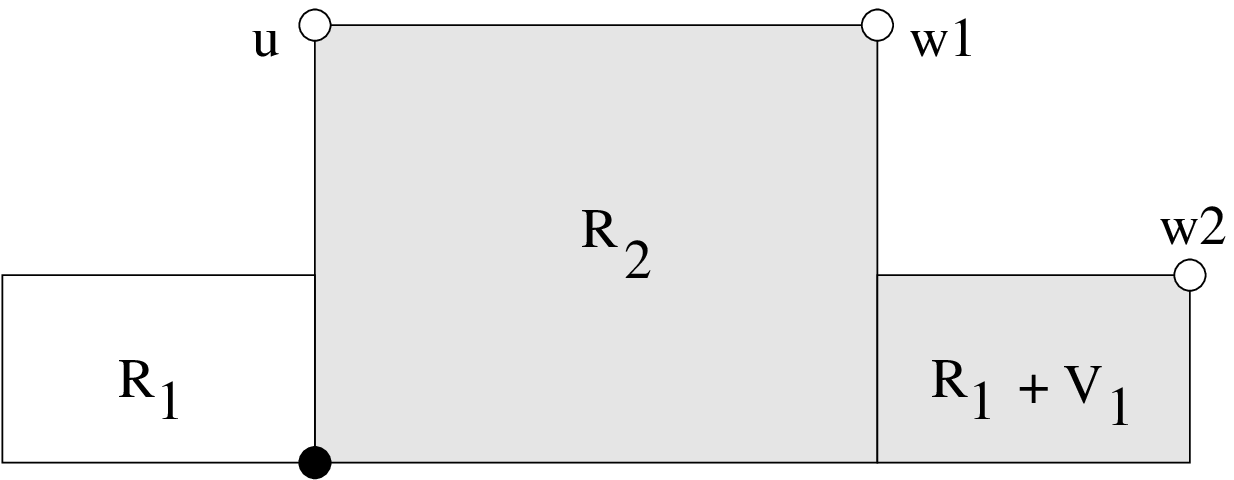}}
\newline
{\bf Figure 20.1:\/} $R_2(A_1)$ and $R_1(A_1)+V_1$.
\end{center}

The vertices shown in Figure 20.1 are
\begin{equation}
u=W_1; \hskip 30 pt
w_1 \approx W_1+\lambda_1 V_1; \hskip 30 pt
w_2 \approx V_1+\mu W_1.
\end{equation}
Here $\mu=q_0/q_1<1/2$, where
$A_0=p_0/q_0$ is the superior predecessor of $A_1$.
Also, $\lambda_1=(q_1)_+/(q_1)$, as in Equation
\ref{determineA}.

The approximation sign means that the distance
between the two points is at most $1$ unit.
For instance, $w_1$ is the intersection of
the line parallel to $W_1$ and containing
$V_+$, and the line parallel to $V_1$ and
containing $W_1$.   The point $V_+$ is
$O(q_1^{-2})$ of the point $\lambda_1 V_1$.
Hence $w_1$ is within $O(q_1^{-2})$ of
$W_1+\lambda_1 V_1$.  The argument for
$w_2$ is similar.

As in the proof of Lemma \ref{weakcopy}, we have
$G_1(u)>-q_1+2$ once $p_1$ is large.  The computations
for $H_1(w_1)$ and $H_1(w_2)$ are the interesting ones.
Case 1 of Lemma \ref{indX} gives $(q_2)_+ \geq (q_1)_+$.
Hence, for
$p_1$ is sufficiently large, we get the following inequalities.
$$
2+H_1(w_1) \leq \Big(2+\|\nabla H \|\Big)+H_1(W_1)+\lambda_1 H_1(V_1) \leq $$
\begin{equation}
 5+\frac{q_1^2}{p_1+q_1}+(q_1)_+< q_1+(q_1)_+ \leq 
q_1+(q_2)_+=\Omega q_1.
\end{equation}
Here we use the bound
$\|\nabla H\| \leq 3$. 
We already remarked that $(q_2)_+ \geq (q_1)_+$. We
also know that $(q_1)_+>q_1/2$.  Hence
\begin{equation}
\Omega q_1=q_1+(q_2)_+>(3/2)q_1.
\end{equation}
For $p_1$ large, we get
$$
2+H_1(w_2) \leq \Big(2+\|\nabla H\|\Big) + H_1(V_1)+ \mu H_1(W_1) <$$
\begin{equation}
5+H_1(V_1)+(1/2)H_1(W_1)=
5+q_1+\frac{q_1^2}{2(p_1+q_1)}<
(3/2)q_1<\Omega q_1.
\end{equation}
These arguments show that
$v \in \Delta_1(I)$ for all $v \in \Gamma_1^1$.
The rest of the proof is just like the proof of
Lemma \ref{weakcopy}.
\endproof

Suppose $A_1'<A_2'$ are two consecutive terms in $\cal S$, when
we have a finite chain
\begin{equation}
A_1'=A_1 \leftarrow A_2 \leftarrow ... \leftarrow A_n=A_2'; \hskip 30 pt
A_1<A_n; \hskip 30 pt q_2>2q_1.
\end{equation}
The following result finishes the proof of
Lemma \ref{strongcopy}.

\begin{lemma}
$\Gamma_1^{n+\epsilon} \subset \Gamma_n^1$.
\end{lemma}

\startproof
We will change our notation slightly from the
previous result.   We let $R_1=R(A_1)$ denote the
parallelogram from the Room Lemma.  Likewise,
let $R_k^*=R^*(A_k)$, the parallelogram from 
Lemma \ref{strong1}.
For any parallelogram $R_k$, let $XR_k$ denote the
union of $R$ with the points within $q_k/8$ units
from the bottom edge of $R_k$. Likewise define
$XR_k^*$.

Since $A_1<A_n$, we have $A_1<A_2$ by Lemma \ref{squeeze}.
We now have
\begin{equation}
\Gamma_1^{1+\epsilon} \subset \Gamma_1 \cap XR_1  \subset \Gamma_2.
\end{equation}
The first containment comes from the Room Lemma
and the definition of $\Gamma_1^{1+\epsilon}$.
The second containment is Lemma \ref{strong2}.
Lemma \ref{strong1} gives us
\begin{equation}
\Gamma_k \cap XR_k^* \subset \Gamma_{k+1}.
\hskip 30 pt k=2,...,n-1.
\end{equation}
Let us compare $R_1$ and $R_k^*$ for $k \geq 2$.
\begin{enumerate}
\item The sides of $R_1$ have length $O(q_1)$.
\item The slope of each side of $R_1$ is within
$O(q_1^{-2})$ of the slope of the corresponding
side of $R_k^*$.  This comes from Lemma \ref{superior}.
\item Each side of $R_1$ is less than half as
long as the corresponding side of $R_k^*$.  This
follows from the first two facts, and from the
fact that $2q_1<q_2 \leq q_k$.  Indeed, the quantity
$q_2-2q_1$ tends to $\infty$ with the complexity
of $A_1$. 
\end{enumerate}
These properties give us
\begin{equation}
XR_1 \subset XR_k^*; \hskip 30 pt k=2,...,n-1.
\end{equation}
Figure 20.2 shows a schematic picture.

\begin{center}
\resizebox{!}{1in}{\includegraphics{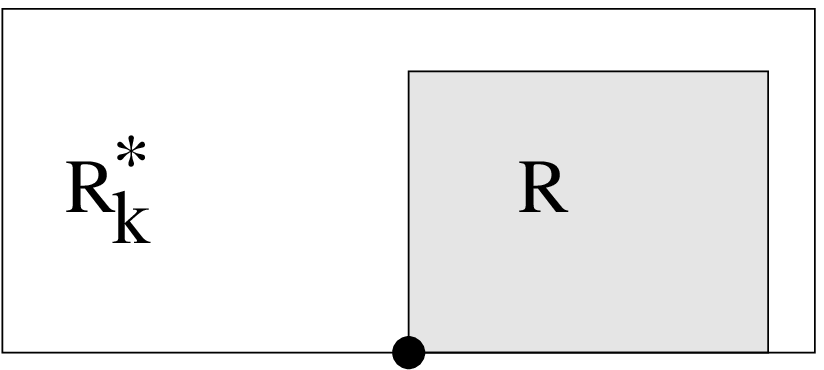}}
\newline
{\bf Figure 20.2:\/} $R_1$ and $R_k^*$ for any $k \geq 2$.
\end{center}

We already know that $\Gamma_1 \cap XR_1 \subset \Gamma_2$.
Suppose $\Gamma_1 \cap XR_1 \subset \Gamma_k$
for some $k \geq 2$.  Then
\begin{equation}
\Gamma_1 \cap XR_1 \subset \Gamma_k \cap XR_1 \subset
\Gamma_k \cap XR_k^*\subset \Gamma_{k+1}.
\end{equation}
Hence, by induction, $\Gamma_1^{1+\epsilon} \subset \Gamma_n$.
The right endpoint of $\Gamma_n^1$ lies far to the right of
any point on $\Gamma_1^{1+\epsilon}$.  Hence
$\Gamma_1^{1+\epsilon} \subset \Gamma_n^1$.
\endproof

\newpage

\section{Proof of the Decomposition Theorem}
\label{decompproof}

\subsection{Decomposition into Arcs}
\label{divreg}

Let $p/q$ be an odd rational in $(0,1)$.
There are $2$ cases of the Decomposition
Theorem, depending on whether $q_-<q_+$ or
$q_+<q_-$.  We will give our proofs mainly in the
case when $q_+<q_-$.  This case relies
on Statement 1 of the Diophantine Lemma.
The other case relies on Statement 2.

Referring to \S \ref{statedecomp}, the
lines $L^+_j$ for $j=0,1,2$ are all parallel to
the vector $W$.  When $q_+<q_-$, the line
$L_1^+$ contains $V_+$ and the line
$L_2^+$ contains $-V_-$.
By the Hexagrid Theorem, $\Gamma$ only crosses
$L_0^+$ once, at the point $(0,0)$. 

\begin{lemma}
\label{symm2}
$\Gamma$ crosses each of $L_1^+$ and $L_2^+$ only
once, and this crossing occurs within $1$ unit of
the baseline.
\end{lemma}

\startproof
This result follows from symmetry and the
Hexagrid Theorem. Our proof refers to Figure 21.1.

\begin{center}
\psfig{file=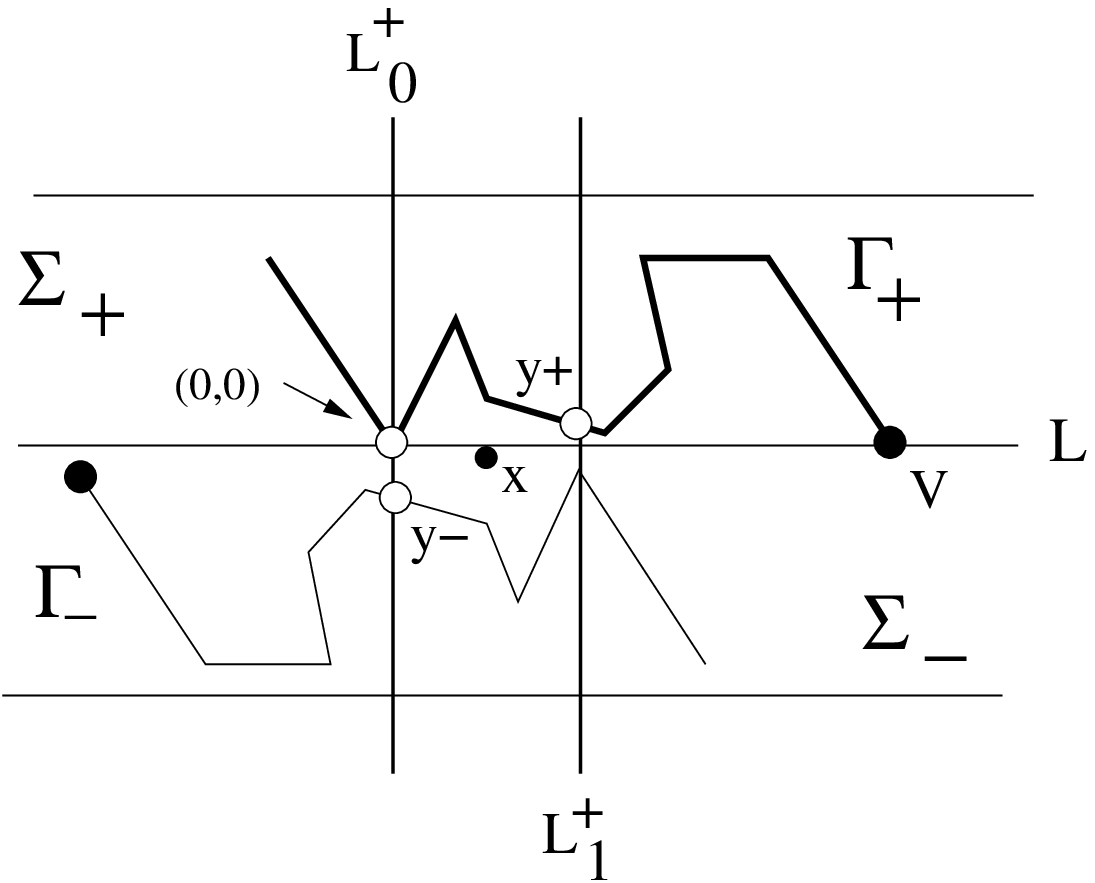}
\newline
{\bf Figure 21.1:\/} Rotating the graph
\end{center}

Let $L$ denote the line of slope $-A$ through the origin.
$\Sigma_+$ (respectively $\Sigma_-$)
 is the infinite strip
bounded by $L$ and the first ceiling line
above (respectively below) $L$.
By Theorem \ref{orbit structure2}, there is one
infinite component of $\widehat \Gamma$ in 
$\Sigma_{\pm}$. We call this component
$\Gamma_{\pm}$.  Here $\Gamma_+=\Gamma$ is the
component of interest to us.

The point $x=(1/2)V_+$ is the fixed point of
$\iota$, the rotation from Equation \ref{iota5}.
We have
\begin{equation}
\iota(L_0^+)=L_1^+; \hskip 30 pt \iota(\Gamma_-)=\Gamma_+;
\hskip 30 pt \iota(L) \downarrow L.
\end{equation}
Our last piece of notation means that $\iota(L)$ lies
(very slightly) beneath $L$.

By the Hexagrid Theorem, $(0,0)$ is the
door corresponding to the point where
$\Gamma_+$ crosses $L_0^+$ and {\it also\/}
to the point $y_-$ where $\Gamma_-$ crosses
$L_0^+$.  This point is the intersection of
$L_0^+$ with the edge connecting $(0,-1)$ to $(-1,0)$.
The image $y_+=\iota(y_-) \in \iota(L_0^+)=L^+$ is the
only point where $\iota(\Gamma_-)=\Gamma_+$
crosses $L_+$.   This point is less than $1$ unit
from $L$ because
$\iota(L)$ lies beneath $L$.  
This shows that $\Gamma=\Gamma_+$ only crosses
$L^+$ once, within $1$ units of $L$.
Since  $L_1^+=L_2^+ \pm V$, and $\Gamma$
is invariant under translation by $V$, it
suffices to prove the result for one
of the lines, as we have done.
\endproof

From this result we see that we can divide a period of
$\Gamma$ into the union of two connected arcs.
One of the arcs lies in what we call $R_0$ and
the other arc lies in $R_2$.  Each arc
connects points near the bottoms of the boxes
and otherwise does not cross the boundaries.
Figure 21.1 shows a schematic picture. Here
$R_0$ is the union of the two shaded regions. Our main goal
is to show that $\Gamma \cap R_0 \subset R_1$.

\begin{center}
\resizebox{!}{2.2in}{\includegraphics{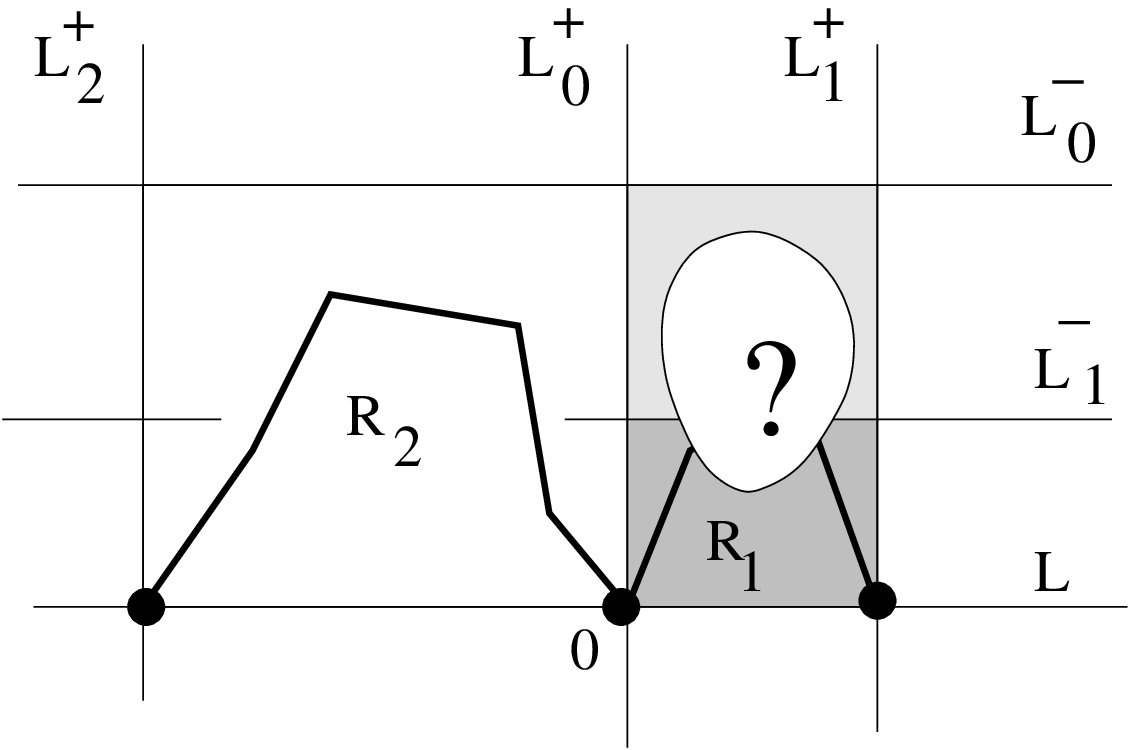}}
\newline
{\bf Figure 21.2:\/} Dividing $\Gamma^1$ into two arcs.
\end{center}

\subsection{The Superior Predecessor}

Let $A'=p'/q'$ denote the superior predecessor of $A$.
Let $\Omega=\Omega(A',A)$.  We consider the case
when $A'<A$.

\begin{lemma}
The second coordinate of any point in $R_1$ lies in
$(0,\Omega q_1'-1)$.
\end{lemma}

\startproof
By convexity, it suffices to consider the vertices of
$R_1$.
The bottom vertices of $R_1$ have first coordinates
$0$ and $q_+$, whereas $\Omega q'=q_++q'$. This
takes care of the bottom vertices.
Let $u=(u_1,u_2)$ be the top left vertex of $R_1$.
Since $R_1$ is a parallelogram, we can finish the proof
by showing that $u_1 \in (0,q'-1)$.
Let $y=(p'+q')/2 \leq q'-1$.  Note that $u$ lies
on a line of slope in $(1,\infty)$ through the origin.
Since the top edge of $R_1$ has negative slope and
contains $(0,y)$, we get $u_2<y$.  Hence $u_1<y$ as well.
\endproof

\begin{lemma}
\label{sup0}
Let $A'$ denote the superior predecessor of $A$.
Suppose that $A' \not = 1/1$.  Then
$\Gamma' \cap R_0 \subset R_{1}$.
\end{lemma}

\startproof
Let $\gamma=\Gamma' \cap R_0$. Since $\gamma$ starts
out in $R_1$ (at the origin), we just need
to see that $\gamma$ never cross the top edge
of $R_{1}$.   The top 
edge of $R_{1}$ contained in the line $\lambda=L_1^-$ of slope $-A$
though the point $X=(0,(p'+q')/2)$.
By the Room Lemma, $\gamma$ does not cross the (nearly identical) line
$\lambda'=(L_0^-)'$ of slope $-A'$ through $X$. 

If $\gamma$ crosses the top edge of $R_{1}$,
then there is a lattice point $(m,n)$ between
$\lambda$ and $\lambda'$, and within $1$ unit of $R_1$.
But then
\begin{equation}
{\rm floor\/}(Am) \not = {\rm floor\/}(A'm); \hskip 30 pt
m \in (-1,q'+q_+)=(-q',\Omega q').
\end{equation}
The second equation comes from our previous result.
Our last equations contradict Lemma \ref{good}.
\endproof

\begin{corollary}
\label{symm3}
Suppose that $\Gamma$ and $\Gamma'$ agree in $R_1$.  Then
The Decomposition Theorem holds for $A$.
\end{corollary}

\startproof
Let's trace $\Gamma \cap R_0$ from left to right,
starting at $(0,0)$.  By hypothesis, this arc does
not cross the top of $R_1$ until it leaves $R_0$.
Once $\Gamma \cap R_0$ leaves $R_0$ from the
right, it never re-enters.  This is a consequence
of Lemma \ref{symm2}.  
\endproof

\subsection{Most of the Parameters}
\label{most}

Here we prove the Decomposition Theorem for most
odd rationals.  We deal with the exceptional
cases in subsequent sections.  Here is the
result we prove.

\begin{lemma}
\label{strong4}
Let $A'=p'/q'$ be the superior predecessor of
$A$.  Then the Decomposition Theorem
holds for $A$ as long as $p' \geq 3$ and $q' \geq 7$.
\end{lemma}

By Corollary \ref{symm3}, it suffices to prove that
$\Gamma'$ and $\Gamma$ agree in $R_1$.

\begin{lemma}
\label{symm4}
$\Gamma' \cap R_1$ and $\Gamma \cap R_1$ have the same outermost
edges.
\end{lemma}

\startproof
The leftmost edge of both arcs is the edge connecting $(0,0)$ to $(1,1)$.
Looking at the proof of Lemma \ref{symm2}, we see that the rightmost
edge $e$ of $\Gamma \cap R_0$ connects $V_++(0,1)$ to $V_++(1,0)$.
Here $V_+=(q_+,-p_+)$.    Applying Lemma \ref{symm2} to
$\Gamma'$, we see that some edge $e'$ of $\Gamma'$ connects
$V_+'+(0,1)$ to $V_+'+(1,0)$. 
 By repeated applications of Case 1 or Case 2 of 
Lemma \ref{indX} tell us that
$V_+=V_+'+d V'$ for some $d \in Z$.
Since $\Gamma'$ is invariant under translation by $V'$, we
see that $e$ is also an edge of $\Gamma'$.
\endproof

\noindent
{\bf Adjacent Mismatch Principle:\/}
Lemma \ref{symm4} has the following corollary. If
$\Gamma'$ and $\Gamma$ fail to agree in $R_1$, then
there are two adjacent vertices of
$\Gamma' \cap R_1$ where our two arithmetic graphs
$\widehat \Gamma$ and
$\widehat \Gamma'$ do not agree.  One can see this by tracing the
two curves from left to right, starting at the
origin.  Once we get the first mismatch on $\Gamma'$ our
arc $\Gamma$ has veered off, and the next
vertex on $\Gamma'$ is also a mismatch.
\newline

In our analysis below,
we will treat the case when $A'<A$.  The other case is
similar.  The bottom right vertex of $R_1$ lies on 
a line of slope in $(1,\infty)$ that contains the point
$V_+$.  The point $V_+$ has the same first coordinate
as the very nearby point
\begin{equation}
\widetilde V_+=\frac{q_+}{q}V.
\end{equation}
Indeed, the two points differ by exactly $1/q$.
Let $\widetilde R_1$ denote the slightly smaller
parallelogram whose vertices are
\begin{equation}
(0,0); \hskip 30 pt u; \hskip 30 pt
\widetilde V_+; \hskip 30 pt
\widetilde w=: u+\widetilde V_+.
\end{equation}
If the Decomposition Theorem fails for $A$, then
at least one of the adjacent vertices of mismatch
will lie in $\widetilde R_1$.  (There are not
two adjacent vertices between the nearly identical
right edges of $R_1$ and $\widetilde R_1$.)

As in the previous chapter, it suffices to make
the extremal calcualation
\begin{equation}
\label{conditions}
G(u) \geq -q'+2; \hskip 30 pt H(\widetilde w) \leq \Omega q'-2=q'+q_+-2.
\end{equation}
The Diophantine Lemma then finishes the proof.

We first need to locate $u$.
There is some $r$ such that $v_1=rW$.   Letting $M$ be the
map from Equation \ref{funm}, relative to the parameter $A$,
we have
$$M(v_1)=M(rW)=p'+q'.$$
Solving for $r$ gives
\begin{equation}
v_1=\bigg(\frac{p'+q'}{p+q}\bigg)W.
\end{equation}
We compute
\begin{equation}
\label{affinecoincidence1}
G(u)=\frac{p'+q'}{p+q}G(W)=-\frac{p'+q'}{p+q} \times \frac{q^2}{p+q}=
\frac{-(1+A')q'}{(1+A)^2}>\frac{-q'}{1+A'}.
\end{equation}
\begin{equation}
\label{affinecoincidence2}
H(\widetilde w)=H(u)+(q_+/q)H(V)=\frac{(1+A')q'}{(1+A)^2}+q_+<
\frac{q'}{1+A'}+q_+.
\end{equation}
Our last inequality in each case uses the fact that $0<A'<A$.
Notice the great similarity in these two calculations.
One can ultimately trace this symmetry back to the
affine symmetry of the arithetic kite ${\cal K\/}(A)$
defined in \S 3.

The conditions in Equation \ref{conditions} are
simultaneously met provided
\begin{equation}
\label{angles}
\frac{-q'}{1+A'} \geq -q'+2; \hskip 50 pt
\bigg(\Longleftrightarrow \hskip 10 pt \frac{1}{p'}+\frac{1}{q'} \leq \frac{1}{2}\bigg)
\end{equation}
The equation on the right is equivalent to the one on the left.
We see easily that it holds as long as $p' \geq 3$ and
$q' \geq 7$.

In the next two sections we will make a more detailed study
of the few exceptions to Lemma \ref{strong4}.  The reader
mainly interested in the Erratic Orbits Theorem can stop
reading here.

\subsection{Some Tricks}
\label{tricks}

Here we take care of some more cases of the Decomposition Theorem.
We use the notation from the previous section.  We assume that
$A' \not = 1/1$ is one of the
rationals not covered by Lemma \ref{strong4}.
 In the previous section we used the
linear functionals $G$ and $H$ defined relative to $A$.
Given the statement of the Diophantine Lemma, we can try to
use the linear functionals $G'$ and $H'$ for the same
purpose.  Here $G'$ and $H'$ are associated to $A'$.
Before we begin our argument, we warn the reader that $G'$ is not
the derivative of $G$. We will
denote the partial derivatives of $G'$
by $\partial_x G'$ and $\partial_y G'$.

\begin{lemma}
\label{dc5}
$G'(v) \geq -q'+2$ for all
$v \in R_{1}$.
\end{lemma}

\startproof
We only have to worry about points near the top left corner
of $R_1$. Such points lie on the first period of $\Gamma'$
to the right of the origin.  Call this period $\beta'$.
When $A' \in \{3/5,3/7,5/7\}$ we check this result explicitly
for every point of $\beta'$.
When $A'=1/q'$ we note
that $\partial_x G'>0$ and $\partial_y G'<0$.
We also note that all points
in $R_{1}$ have positive first coordinate and second
coordinate at most $(q'-1)/2$.  Thus, the point that
minimizes $G'$ is $v=(1,(q'-1)/2)$.  We compute
$$G'(v)+q-2=\frac{q'-3}{q'+1} \geq 0.$$
The extreme case occurs when $q'=3$. 
\endproof

$H'$ is tougher to analyze because the points of
interest to us are near the top right corner of
$R_1$, and this corner can vary drastically with
the choice of $A$.  We will use rotational symmetry
to bring the points of interest back into view,
so to speak.
Let $\iota$ be the isometric involution that swaps
$(0,0)$ and $V_+$.   Repeated applications of
Lemma \ref{indX} show that $V_+=V_+'+dV'$
for some $d \in \Z$.  Hence $\iota$ is a
symmetry of $\widehat \Gamma'$.
See the remark after Equation \ref{rotpoints}.

The infinite arc $\iota(\Gamma')$ is the open component of
$\widetilde \Gamma'$ that lies just beneath the baseline.
One period of $\iota(\Gamma')$ connects the point
$(0,-1)$ to the point $(q',-p'-1)$.  Let's denote this
period by $\beta'$.  Compare the proof of
Lemma \ref{symm2}.  The points of $R_1$ near the top
right corner correspond to points on $\beta'$.  To
evaluate $H'$ on the points near the top right corner
of $R_1$, we evaluate $H'$ on points of $\beta'$ and
then relate the results. 

\begin{lemma}
For any $v \in \R^2$ we have
$$|H'(v)+H'(\iota(v)) - q_+| <\frac{2}{q'}.$$
\end{lemma}

\startproof 
Since $H'$ is a linear functional, it suffices to prove
our result for $v=(0,0)$.  In this case, we must
show that $|H'(V_+)-q_+|<2/(q')$.  We have already
remarked that
$V_+=V_+'+dV'$.  Hence $q_+=q_+'+dq'$. From Lemma
\ref{Hcalc}, we get
$H'(dV')=dq'$.  Hence, our equality is equivalent to
\begin{equation}
\label{target1}
|H'(V_+')-q_+'|<\frac{2}{q'}.
\end{equation}   The point
$V_+'$ lies on the same vertical line as the
point $u'=(q_+'/q')V'$, and exactly $1/q'$ units
away. Equation \ref{target1} now follows from the next $3$ facts.
\begin{equation}
\label{target2}
H'(u')=q_+'; \hskip 30 pt |\partial_y H'|<2; \hskip 20 pt
\|u'-V_+'\|=\frac{1}{q'}.
\end{equation}
The first fact comes from Lemma \ref{Hcalc}.  The second
fact is an easy calculus exercise.  The third fact, already
mentioned, is an easy exercise in algebra that uses $|q'p_+'-p'q_+'|=1.$
\endproof

The bound
$$H'(v) \leq \Omega q'-2=q'+q_+-2$$
only fails for points very near the top right vertex of $R_1$.
Any such point has the form $\iota(v)$ for some $v \in \beta'$.
Thus, to establish the above bound, it suffices to prove that
\begin{equation}
H'(v) \geq -q'+2+\frac{2}{q'}.
\end{equation}
This inequality can fail for very small choices of $q'$.
However, from the Adjacent Mismatch Principle, the
inequality must fail for at least $2$ vertices on
$\beta'$, and this does not happen.

We check all cases with $q' \leq 7$ by hand.  This
leaves only $A'=1/q'$ for $q' \geq 9$. Reasoning as we
did in Lemma \ref{dc5}, we see that
the extreme point is $v=(0,(1-q)/2)$. 
We compute
\begin{equation}
H'(v)-\bigg(-q'+2+\frac{1}{q'}\bigg)=\frac{2(q'^2-2q'-1)}{(1+q')^2} -\frac{2}{q'}>0.
\end{equation}
The last equation is an easy exercise in calculus.
This completes our proof of the Decomposition Theorem for
all parameters $A$ such that $A' \not = 1/1$.

\subsection{The End of the Proof}
\label{decompend}

Now we deal with the case 
when $A'=1/1$ is the superior predecessor of $A$.
We have the following structure
\begin{equation}
\frac{1}{1} \leftarrow A_1= \frac{2k-1}{2k+1} \leftarrow ... \leftarrow A_m=1.
\end{equation}
Here $k \geq 1 $.  
For instance, when $A=17/21$, we have $1/1 \leftarrow 9/11 \leftarrow 17/21$.
Figure 21.3 shows $\Gamma(17/21)$.  In this case
$\Gamma \cap R_{1}$ is the line segment connecting
$(0,0)$ to $(-5,5)=(-k,k)$.  We will establish
this structure in general.

\begin{center}
\resizebox{!}{4in}{\includegraphics{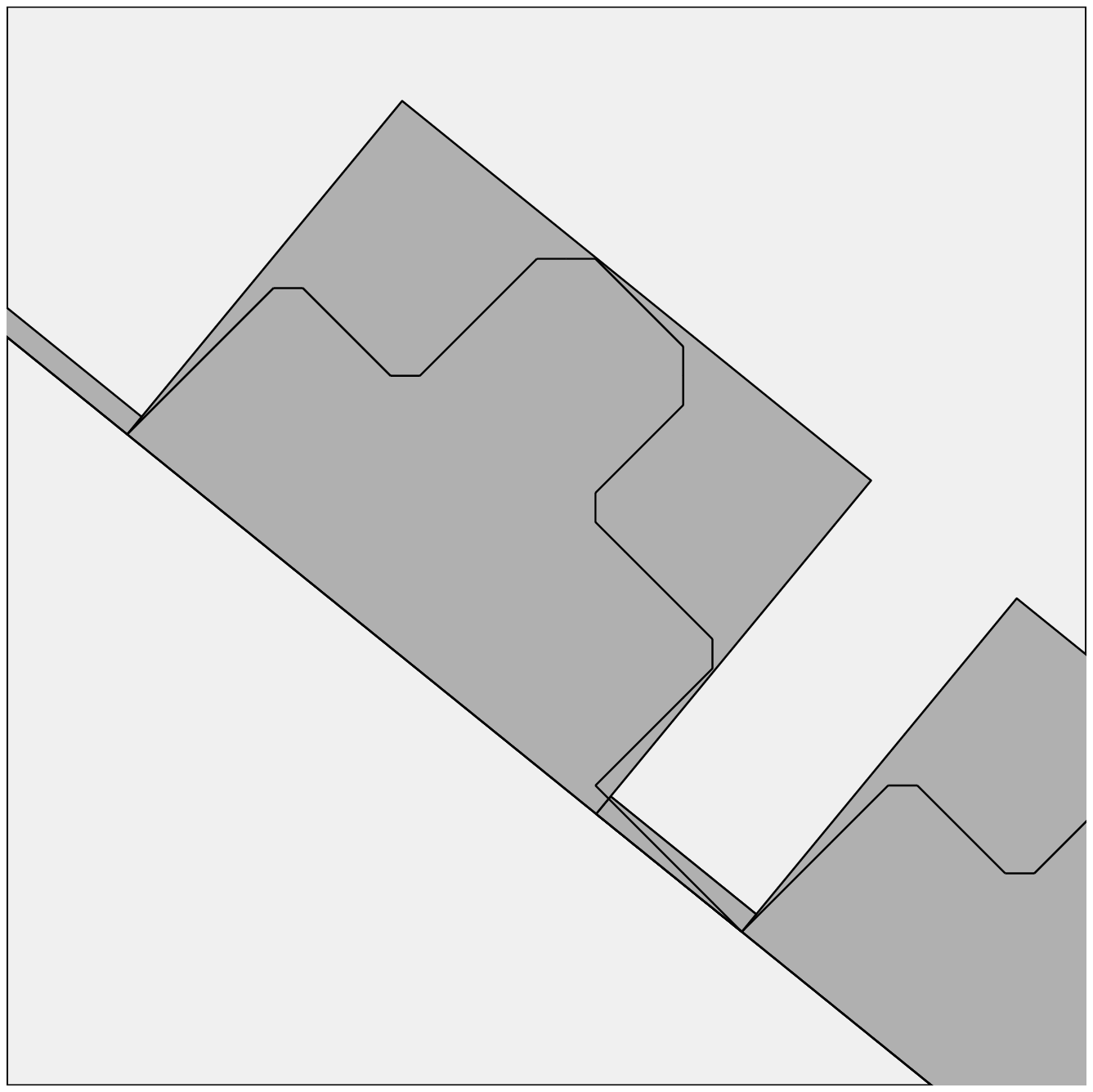}}
\newline
{\bf Figure 21.3:\/} $\Gamma(17/21)$
\end{center}

$R_1$ is the very short and squat parallelogram 
near the bottom right corner of
Figure 21.3.  This time $R_1$ lies to the left of
the origin. The left side of $R_1$ lies in $L_1^+$.  Repeated
applications of Lemma \ref{indX} show that
$(-k,k-1) \in L_1^+$. 
The right side of $R_0$ lines in $L_0^+$, the
parallel line through the origin.
The top of $R_{1}$ contains $(0,1)$ and is parallel
to the baseline.

Let $\gamma= \Gamma \cap R_0$.  The rightmost
vertex of $\gamma$ is $(0,0)$, and the rightmost
edge of $\gamma$ connects
$(0,0)$ to $(-1,1)$.  Compare the proof of the
Room Lemma.

\begin{lemma}
The leftmost edge of $\gamma$ connects $(-k,k)$ to $(-k+1,k-1)$.
\end{lemma}

\startproof
By Lemma \ref{symm2}, there is a
unique edge $e$ of $\Gamma$ that crosses $L_1^+$.
Looking at the proof of Lemma \ref{symm2} we see
$e=\iota(e')$, where $e'$ connects $(0,-1)$ to
$(-1,0)$ and $\iota$ is order $2$ 
rotation about the point
\begin{equation}
\bigg(\frac{-q_-}{2},\frac{p_-}{2}\bigg)=\bigg(\frac{-k}{2},\frac{k-1}{2}\bigg).
\end{equation}
From this, we conclude that $e$ connects $(-k,k)$ to $(-k+1,k-1)$.
The leftmost edge of $\gamma$ crosses $L_1^+$.  This edge
must be $e$.
\endproof

\begin{lemma}
The line segment $\gamma'$ connecting $(0,0)$ to $(-k,k)$ lies beneath $L_0^-$.
Hence $\gamma' \cap R_0 \subset R_{1}$.
\end{lemma}

\startproof
Letting $F(m,n)=Am+n$, we have $F(0,1)=1$.
Hence $F(x)=1$ for all $x \in L_0^-$.
On the other hand, we compute that
$F(0,0)=0$ and $F(-k,k)=2k/(2k+1)<1$.  By
convexity, $F(y)<1$ for all $y \in \gamma'$.
\endproof

To finish our proof, we just have to show that $\gamma'=\gamma$.
The first and last edges of $\gamma$ and $\gamma'$ agree,
and these edges are $\pm (1,-1)$, with the sign depending
on which way we orient our curves.
Let $p_j=(-j,j)$, for $j=2,...,k-1$.
By Lemma \ref{squeeze}, we have 
\begin{equation}
\label{boundsk}
(A_1)_-=\bigg(\frac{k-1}{k}\bigg)<A<\bigg(\frac{k}{k+1}\bigg)=(A_1)_+;
\hskip 20 pt
\frac{1}{k+1}<1-A<\frac{1}{k}
\end{equation}
The first equation implies the second.
We compute
\begin{equation}
M_+(p_j)=(x_j,y_j,z_j)=j(1-A,1-A,1-A)+(0,1,0); \hskip 20 pt
{\rm mod\/} \hskip 10 pt \Lambda
\end{equation}
Equation \ref{boundsk} combines with the fact that
$j \in \{1,...,k-1\}$ to give
\begin{equation}
x_{j}=z_j \in [1-A,A); \hskip 30 pt
x_{j}+y_{j}-z_{j} = y_j \in (1,1+A) \subset (A,1+A).
\end{equation}
We check that these inequalities always specify the edge $(-1,1)$.
Hence $\gamma'$ and $\gamma$ are both line segments.
Hence $\gamma=\gamma'$.

\newpage

{\bf {\Huge Part V\/}\/}
\newline

\begin{itemize}

\item In \S \ref{oddappx} we prove some further results
about the inferior and superior sequences.
We list the basic results in the first
section and then spend the rest of
the chapter  proving these results.

\item In \S \ref{fundamental}, we prove
Theorem \ref{discrete}. 
We also build a rough model for the
way the orbit $O_2(1/q_n,-1)$ returns to
the interval $I=[0,2] \times \{-1\}$.
Our work here depends on two technical
results, the Copy Theorem and the
Pivot Theorem, which we
establish in Part VI.

\item In \S \ref{comet2} we prove 
Statements 2,3,4 of the Comet Theorem,
modulo some technical details that we
handle in Part VI.
We defer the proof of Statement 1 of
the Comet Theorem until
Part VI.  

\item in \S \ref{comet2} we deduce
a number of dynamical consequences of
the Comet Theorem, including minimality
of the set of unbounded orbits.  We also
define the cusped solenoids and explain
how the time-one map of their geodesic
flow models the outer billiards dynamics.

\item in \S \ref{lengthdim} we analyze the
structure of the Cantor set $C_A$. This
chapter has a number of geometric results,
such as a formula for $\dim(C_A)$ when
$A$ is a quadratic irrational.

\end{itemize}

\newpage

\section{Odd Approximation Results}
\label{oddappx}

\subsection{The Results}

Let $\{p_n/q_n\}$ be the inferior sequence associated to
an irrational parameter $A \in (0,1)$,
and let $\{d_n\}$ be sequence obtained from
Equation \ref{renorm}.
We call $\{d_n\}$ the {\it inferior renormalization sequence\/}.
We call the subsequence of $\{d_n\}$ corresponding to the
superior terms the {\it superior renormalization sequence\/}
or just the {\it renormalization sequence\/}.   
Referring to the inferior sequences, we have
$d_n=0$ if and only if $n$ is not a superior term.
In this case, we call $n$ an {\it inferior term\/}.
So, the renormalization sequence is created from
the inferior renormalization sequence simply by
deleting all the $0$s.

For any odd rational $p/q \in (0,1)$, define
\begin{equation}
\label{star}
p^*=\min(p_-,p_+); \hskip 30pt
q^*=\min(q_-,q_+).
\end{equation}
Here $p^*/q^*$ is one of the rationals $p_{\pm}/q_{\pm}$.
It is convenient to define
\begin{equation}
\frac{p_0^*}{q_0^*}=\frac{1}{0}.
\end{equation}

Given the superior sequence $\{p_n/q_n\}$ we define
\begin{equation}
\lambda_n=|Aq_n-p_n|; \hskip 30 pt
\lambda^*_n=|Aq_n^*-p_n^*|; \hskip 30pt
\end{equation}
Note that 
\begin{equation}
\lambda_0^*=1.
\end{equation}
For the purposes of making a clean statement,
we define $\lambda_{-1}=+\infty$.
All our results are meant to apply
to the superior sequence, for indices $n \geq 0$.

\begin{equation}
\label{DIO1}
d_n\lambda_n<2q_n^{-1}; 
\end{equation}

\begin{equation}
\label{DIO4}
q_{2n}>(5/4)^n D_{2n}.
\end{equation}

\begin{equation}
\label{DIO2}
\label{DIO3}
\sum_{k=n}^{\infty}d_k\lambda_k=\lambda_{n}^*<\lambda_{n-1}
\end{equation}
Note that Equation \ref{DIO1} is an immediate
consequence of Lemma \ref{superior}.
The rest of the chapter is devoted
to proving Equations \ref{DIO4}
and \ref{DIO2}.

\subsection{The Growth of Denominators}
\label{groden}

Let $\delta_n$ be as in Equation \ref{enhanced}.
By Lemma \ref{indX},
the sequence $\{\delta_n\}$ determines the sequence $\{p_n/q_n\}$.
At each step, $(q_{n+1})_{\pm}$ is a non-negative integer linear
combination of $(q_n)_{\pm}$, and the precise linear
combination is determined by $\{\delta_n\}$.
Call this the {\it positivity property\/}.
Call the sequence $\{\delta_n\}$ the {\it inferior enhanced 
renormalization sequence\/}, or {\it IERS\/} for short.
Call the subsequence corresponding
to the superior indices the 
{\it enhanced renormalization sequence\/}.
The reason for the terminology is that we can
determine the inferior renormalization sequence
from the IERS,
but not {\it vice versa\/}.

Say that a parameter $A$ is {\it superior to\/} the parameter $A'$,
if the IERS for $A'$ is obtained
by inserting some $1$s into the IERS for $A$.
For instance, $\sqrt 5-2$ has IERS $2,1,2,1...$
and $\sqrt 2-1$ has IERS
sequence $2,2,2...$.  Hence $\sqrt 2-1$ is superior
to $\sqrt 5-2$.

\begin{lemma}
\label{growth}
Suppose that $A$ is superior to $A'$.  Then
$q_n \leq q_n'$ for all $n$.
\end{lemma}

\startproof
Consider the operation of inserting a $1$ into the $m$th position
in IERS for $A$ and recomputing $\{A_n\}$.
Call this new sequence the $A^*$-sequence.
We have
$(q_{m+1}^*)_{\pm} \geq (q_m)_{\pm}.$
By induction the positivity property, we get
$(q_{n+1}^*)_{\pm} \geq (q_n)_{\pm}.$
Now let's delete the $(m+1)$st term from the
$A^*$-sequence.  Call the new sequence the $A'$-sequence.
We have $q_n' \geq q_n$ for all $n$.   Our result
now follows from induction.
\endproof

Call $A$ {\it superior\/} if the corresponding 
inferior sequence has no inferior terms.  That
is, the IERS has no $1$s in it.  For instance
$\sqrt 2 -1$ is a superior parameter.
If we want to get a lower bound on the growth of
denominators, it suffices to consider only the
superior parameters. 
Equation \ref{DIO4} follows from induction and
our next lemma.

\begin{lemma}
\label{fivefour}
Suppose that $A_1,A_2,A_3$ are $3$ consecutive terms in
the superior sequence.  Let $d_1,d_2,d_3$ be the
corresponding terms of the renormalization sequence.
Then $q_3>(5/4)(d_1+1)(d_2+1)q_1$.
\end{lemma}

\startproof
It suffices to assume that $A$ is a superior parameter,
so that $A_1,A_2,A_3$ are (also) $3$ consecutive terms
in the inferior sequence. 

First of all, the estimates 
\begin{equation}
q_{n+1} > 2d_n q_n;  \hskip 40 pt
q_{n+1} > \delta_n q_n; 
\end{equation}
follow directly from the definitions.
Our notation is as in Lemma \ref{indX}.

Suppose first that $\min(d_1,d_2) \geq 2$.
Then
\begin{equation}
q_{3} > 4d_1d_{2}q_1>\frac{4}{3}(d_1+1)(d_{2}+1)q_1.
\end{equation}
Now suppose that $d_1=d_2=1$ and
$\min(\delta_1,\delta_2) \geq 3$.
Then
\begin{equation}
q_{3}>6q_{1}=\frac{3}{2}(d_1+1)(d_{2}+1)q_1.
\end{equation}

Suppose finally that
$d_1=d_2=1$ and $\delta_1=\delta_2=2$.
We will deal with the case that
$A_1<A_2$.  The other case is similar.
In this case, we must have
\begin{equation}
A_0>A_1<A_2>A_3
\end{equation}
 by Lemma \ref{indX}. 

 By Case 2 of Lemma \ref{indX},
\begin{equation}
\label{appx4}
(q_2)_-=q_1+(q_1)_+; \hskip 30 pt
(q_2)_+=(q_1)_+.
\end{equation}
By Case 4 of
By Lemma \ref{indX},
\begin{equation}
(q_3)_+=q_2+(q_2)_-; \hskip 30 pt
(q_3)_-=(q_2)_-.
\end{equation}
Hence,
\begin{equation}
\label{appx5}
q_3=(q_3)_++(q_3)_-=q_2+2(q_2)_-=^^*q_2+2q_1+2(q_2)_+.
\end{equation}
The starred equality comes from Lemma \ref{squeeze2},
since $A_1<A_2$.

 Since $A_0>A_1$, Lemma \ref{indX} says that
\begin{equation}
\label{appx6}
2(q_1)_+>(q_1)_++(q_1)_-=q_1.
\end{equation} 
Combining Equations \ref{appx4}, \ref{appx5}, and
\ref{appx6}, we have
\begin{equation}
q_3=q_2+2q_1+2(q_1)_+>q_2+3q_1>5q_1.
\end{equation}
Hence
\begin{equation}
q_3>\frac{5}{4}(d_1+1)(d_2+1)q_1.
\end{equation}
This completes our proof.
\endproof

\subsection{The Identities}

We first verify the identity in Equation \ref{DIO3}.
In this identity, we sum over the superior indices.
However, notice that
we get the same answer if we sum over all indices. The
point is that $d_n=0$ when $n$ is an inferior index.  So,
for our derivation, we work with the inferior sequence.
Let $\{p_n/q_n\}$ be the inferior sequence associated to $A$.
Define
\begin{equation}
\Delta(n,N)=|p_Nq_n-q_Np_n|; \hskip 30 pt
\Delta^*(n,N)=|p_Nq_n^*-q_Np_n^*|; \hskip 30 pt N \geq n.
\end{equation}

\begin{lemma}
\label{beauty}
$\Delta^*(n,N)-\Delta^*(n+1,N)=d_n \Delta(n,N)$.
\end{lemma}

\startproof 
The quantities relevant to the case $n=0$ are
$$A_0=\frac{1}{1}; \hskip 30 pt
A_0^*=\frac{1}{0}; \hskip 30 pt
A_1^*=\frac{d_0-1}{d_0}<A_1=\frac{2d_0-1}{2d_0+1}.$$
In this case, a simple calculation checks the
formula directly.  

Now suppose $n \geq 1$.
We suppose
that $A_{n-1}<A_n$.  The other case has a similar treatment.
Let $r$ stand for either $p$ or $q$.   There are two cases,
depending on whether the index $n$ has type 1 or type 4.
When $n$ has type 1, Lemma \ref{indX} gives
\begin{equation}
\label{CASE1}
r_n^*=(r_n)_+; \hskip 30 pt
r^*_{n+1}=(r_{n+1})_+; \hskip 30 pt
r_n^*=d_nr_n-r_{n+1}^*.
\end{equation}
We have $\Delta^*(n,N)=|a_1-a_2|$, where
$$
a_1=d_np_{N}q_n-d_nq_Np_n=d_n\Delta(n,N);$$
\begin{equation}
a_2= p_Nq_{n+1}^*-q_Np_{n+1}^*=-\Delta^*(n+1,N).
\end{equation}
The sign for $a_1$ is correct because
$A_N>A_n$.  The sign for $a_2$ is correct
because, by Lemma \ref{squeeze}, we have
$A_N<(A_{n+1})_+=A_{n+1}^*$.
The identity in this lemma follows immediately.

When $n$ has type 4, Lemma \ref{indX} gives
$$r_n^*=(r_n)_+; \hskip 30 pt r_{n+1}^*=(r_{n+1})_-; \hskip 30 pt
r_{n}^*=d_nq_n-r_{n+1}^*.$$
Hence $\Delta^*(n,N)=|a_1+a_2'|$, where $a_2'=-a_2$.  
 The sign changes for $a_2'$ because
$A_N>(A_{n+1})_-=A_{n+1}^*$.  In this case,
we get the same identity.
\endproof

Dividing the Equation in Lemma \ref{beauty} by $q_N$, we get
\begin{equation}
|A_Np^*_{n}-q^*_{n}|-|A_Np^*_{n+1}-q^*_{n+1}|=d_n|A_Np_n-q_n|.
\end{equation}
Taking the limit as $N \to\infty$, we get
\begin{equation}
\lambda^*_{n}-\lambda^*_{n+1}=d_n\lambda_n.
\end{equation}
Summing this equation from $n+1$ to $\infty$ gives the
equality in Equation \ref{DIO3}.

Now we verify the inequality in Equation \ref{DIO3}.

\begin{lemma}
$\lambda_{n+1}^*<\lambda_n$.
\end{lemma}

\startproof
There are two cases to consider, depending on whether
$A_n<A$ or $A_n>A$.   We will consider the case when
$A_n<A$. The other case has a similar treatment.
By Lemma \ref{squeeze}, we have $A_n<A_{n+1}$. Therefore,
by Lemma \ref{indX} (applied to $m=n+1$), we have
$(q_{n+1})_+<(q_{n+1})_-.$
But this means that
$A_{n+1}^*=(A_{n+1})_+$.
By Lemma \ref{squeeze}, we have
\begin{equation}
A_n<A<A_{n+1}^*.
\end{equation}
Given the above ordering, we have
$$\lambda_n=|Aq_n-p_n|=Aq_n-p_n$$
and
$$\lambda^*_{n+1}=|Aq^*_{n+1}-p^*_{n+1}|=p^*_{n+1}-Aq^*_{n+1}.$$
Hence
\begin{equation}
\label{goodgap}
\lambda_n-\lambda^*_{n+1}=A(q_n+q_{n+1}^*)-(p_n+p^*_{n+1})
\end{equation}
But
$$q_n+q_{n+1}^*=q_n+(q_{n+1})_+=\Big(q_{n+1})_--(q_{n+1})_+\Big)-(q_{n+1})_+=(q_{n+1})_-.$$
Likewise
$$p_n+p_{n+1}^*=(p_{n+1})_-.$$
Combining these identities with Equation \ref{goodgap}, we get
$$\lambda_n-\lambda^*_{n+1}=A(q_{n+1})_--(p_{n+1})_-=
(q_{n+1})_-(A-(A_{n+1})_-)>0.$$
This completes the proof.
\endproof
\newpage
\section{The Fundamental Orbit}
\label{fundamental}

\subsection{Main Results}
\label{maintwirl}

We will assume that $p/q=p_n/q_n$, the $n$th term
in a superior sequence.  We call $O_2(1/q_n,-1)$ the
{\it fundamental orbit\/}.
 Let $C_n$ denote the set from Theorem \ref{discrete}.
Let
\begin{equation}
C_n'=O_2(1/q_n,-1) \cap I; \hskip 30 pt
I=[0,2] \times \{-1\}.
\end{equation}
Theorem \ref{discrete} says that
$C_n=C_n'$.   In this chapter we will prove
Theorem \ref{discrete}, and establish a
some geometric results about how the
orbits return to $C_n$.

After we prove Theorem \ref{discrete},
we establish a coarse model for how
the points of $O_2(1/q_n)$ return to
$C_n$.  Statement 2 of the Comet Theorem
is the ``geometric limit'' of the
Discrete Theorem, and Statement 3 of
the Comet Theorem is the ``geometric
limit'' of the coarse model we build here.

Let $\Pi_n$ denote the truncation of the space
defined in Equation \ref{productspace}.  
Let $\chi: \Pi_n \to C_n$ denote the mapping
given in Theorem \ref{discrete}.  We will
describe the ordering on $\Pi_n$ sch that
$\chi(\kappa)$ returns to $\chi(\kappa_+)$, where
$\kappa_+$ is the successor of $\kappa$ in the ordering.

Here we will define two natural orderings on the
sequence space $\Pi_n$ associated to $p_n/q_n$.  
Let $\{d_n\}$ be the renormalization 
sequence.
\newline
\newline
{\bf Reverse Lexicographic Ordering:\/}
Given two finite  sequences $\{a_i\}$ and $\{b_i\}$ of the same length,
let $k$ be the largest index where $a_k \not = b_k$.  We define
$\{a_i\} \prec' \{b_i\}$ if $a_k<b_k$, and
$\{b_i\} \prec' \{a_i\}$ if $a_k>b_k$.  This ordering is known as
the {\it reverse lexicographic\/} ordering. 
\newline
\newline
{\bf Twist Automorphism:\/}
 Given a sequence
$\kappa=\{k_i\} \in \Pi_n$, 
we define $\widetilde k_i=k_i$ if $A_i<A_n$, and $\widetilde k_i=d_i-k_i$
if $A_i>A_n$.  We define $\widetilde \kappa=\{\widetilde k_i\}$.
The map $\kappa \to \widetilde \kappa$ is an involution
on $\Pi_n$.   We call this involution the {\it twist involution\/}.
\newline
\newline
{\bf Twirl Ordering:\/}
Any ordering on $\Pi_n$ gives an ordering on $C_n$, via
the formula in Theorem \ref{discrete}.  Now we describe
the ordering that comes from the first return map.
Given two sequences,
$\kappa_1, \kappa_2 \in \Pi_n$, we define
$\kappa_1 \prec \kappa_2$ if and only if
$\widetilde \kappa_1 \prec' \widetilde \kappa_2$.  We call the
ordering determined by $\prec$ the {\it twirl ordering\/}.
We think of the word ``twirl'' as a kind of
acronym for {\it twisted reverse lexicographic\/}.
We will give an example below.

\begin{lemma}
\label{returncomb}
 When $C_n$ is equipped with the twirl order,
each element of $C_n$ except the last returns to its immediate successor,
and the last element of $C_n$ returns to the first.
\end{lemma}

Our third goal is to understand 
$O_2(1/q_n,-1)$ far away from $I$.
Let $h_1(\kappa)$ denote the maximum distance the forward
$\Psi$-orbit of $\chi(\kappa)$ gets from the kite
vertex $(0,1)$ before returning as $\chi(\kappa_+)$.
Let $h_2(\kappa)$ denote the number of iterates
it takes before the forward $\Psi$-orbit of
$\chi(\kappa)$ returns as $\chi(\kappa_+)$.

Let $\sigma(\kappa)$ be the largest index $k$ such
that the sequences corresponding to $\kappa$ and
$\kappa_+$ differ in the $k$th position.  Here
$\sigma(\kappa) \in \{0,...,n-1\}$.  Finally, we
define $\sigma(\kappa)=n$ if $\kappa$ is the
last element of $\Pi_n$. 

\begin{lemma}
\label{coarsefar}
Let $m=\sigma(\kappa)$.  Then
$$q_m/2-4<h_1(\kappa)<2q_m+4; \hskip 50 pt
h_2(\kappa)<5q_m^2.$$
\end{lemma}

The table below encodes the example from the introduction:
$$\frac{p_0}{q_0}=\frac{1}{1}> \frac{1}{3}< \frac{5}{13}> \frac{19}{49}=\frac{p_3}{q_3}.$$
The first $3$ columns indicate the sequences.
The next column indicates the first coordinate of
$49 \chi(\kappa)$.  The first point of $C_3$
is $(65/49,-1)$.   The next column shows
$(m)=\sigma(\kappa)$.  The last column shows $q_m$. 
$$	
\matrix{
1&0&1& \longrightarrow & 65 & (0)&1 \cr
0&0&1&\longrightarrow &  5&   (1)&3 \cr
1&1&1&\longrightarrow & 81&   (0)&1 \cr
0&1&1&\longrightarrow & 21&(1)&3 \cr
1&2&1&\longrightarrow & 97& (0)&1 \cr
0&2&1&\longrightarrow &37& (2)&13 \cr
1&0&0&\longrightarrow & 61&(0)&1 \cr
0&0&0&\longrightarrow &1& (1)&3 \cr
1&1&0&\longrightarrow &77& (0)&1 \cr
0&1&0&\longrightarrow &17& (1)&3 \cr
1&2&0&\longrightarrow &93& (0)&1 \cr
0&2&0&\longrightarrow &33&  (3)& 49}
$$
For instance, the
the $\Psi$ orbit of $37/49$ to wanders between $13/2-4=5/2$ and $2*13+4=30$
units away before returning to $61/49$ in
less than $5 \times (13^2)$ steps.
This is not such an inspiring result.  We have included
this small example just to show how the chart works.
Larger examples would yield much more dramatic results.

\subsection{The Copy and Pivot Theorems}
\label{pivots}
\label{statepivot}

Here we describe the technical results that we
will establish in Part VI.

Relative to the parameter $A$, we associate a sequence
of {\it pairs of points\/} in $\Z^2$.  We call these
points the {\it pivot points\/}.  We make the construction
relative to the inferior sequence.

Define $E_0^{\pm}=(0,0)$ and $V_n=(q_n,-p_n)$.
Define
\begin{equation}
\label{pivotplus}
A_n<A_{n+1} \hskip 10 pt
\Longrightarrow \hskip 10 pt
E_{n+1}^{-}=E_n^{-}; \hskip 30 pt
E_{n+1}^{+}=E_n^{+}+d_nV_n
\end{equation}
\begin{equation}
A_n>A_{n+1}\hskip 10 pt
\Longrightarrow \hskip 10 pt
E_{n+1}^{-}=E_n^{-}-d_nV_n; \hskip 30 pt
E_{n+1}^{+}=E_n^{+}.
\end{equation}
We have set $A_n=p_n/q_n$.
Here is an example. 
$$\frac{1}{1} \stackrel {>}{\leftarrow} \frac{3}{5} \stackrel {>}{\leftarrow}
\frac{17}{29} \stackrel {<}{\leftarrow} \frac{37}{63} \stackrel {<}{\leftarrow} \frac{57}{97}
\stackrel{>}{\leftarrow}\frac{379}{645}$$
The inferior renormalization sequence is $2,2,1,0,3$.
We compute 

\begin{itemize}
\item $E^+_1=E^+_0=(0,0)$
\item $E^+_2=E^+_1=(0,0)$ 
\item $E^+_3=E^+_2+1(29,-17)$
\item $E^+_4=E^+_3+0(97,-57)=(29,-17)$
\item $E^+(379/645)=E^+_5=E^+_4$.
\end{itemize}

\begin{itemize}
 \item $E_1^-=E_0^--2(1,-1)=(-2,2)$
\item $E_2^-=E_1^--2(5,-3)=(-12,8)$
 \item $E_3^-=E_2^-=(-12,8)$ 
 \item $E_4^-=E_3^-=(-12,8)$ 
\item $E^-(379/645)=E_5^-=E_4^--3(97,-57)=(-303,197)$.
\end{itemize}

This procedure gives
inductive way to define the pivot points to a
pair of odd rationals.
We define the {\it pivot arc\/} $P\Gamma$ of
$\Gamma$ to be the arc whose
endpoints are $E^+$ and $E^-$.
It turns out that the pivot arc is well-defined -- this
is something we will prove simultaneously with our
Copy Theorem below.
This is to say that
$E^+$ and $E^-$ are both
vertices of $\Gamma$. 
In Part VI we prove the following result.

\begin{theorem}[Copy]
\label{copy1}
If $A_1 \leftarrow A_2$ then
$P\Gamma_2 \subset \Gamma_1$.  
\end{theorem}

Figures 22.1 and 22.2 together illustrate Theorem \ref{copy1}
for $17/29 \leftarrow 57/97$.

\begin{center}
\resizebox{!}{3.8in}{\includegraphics{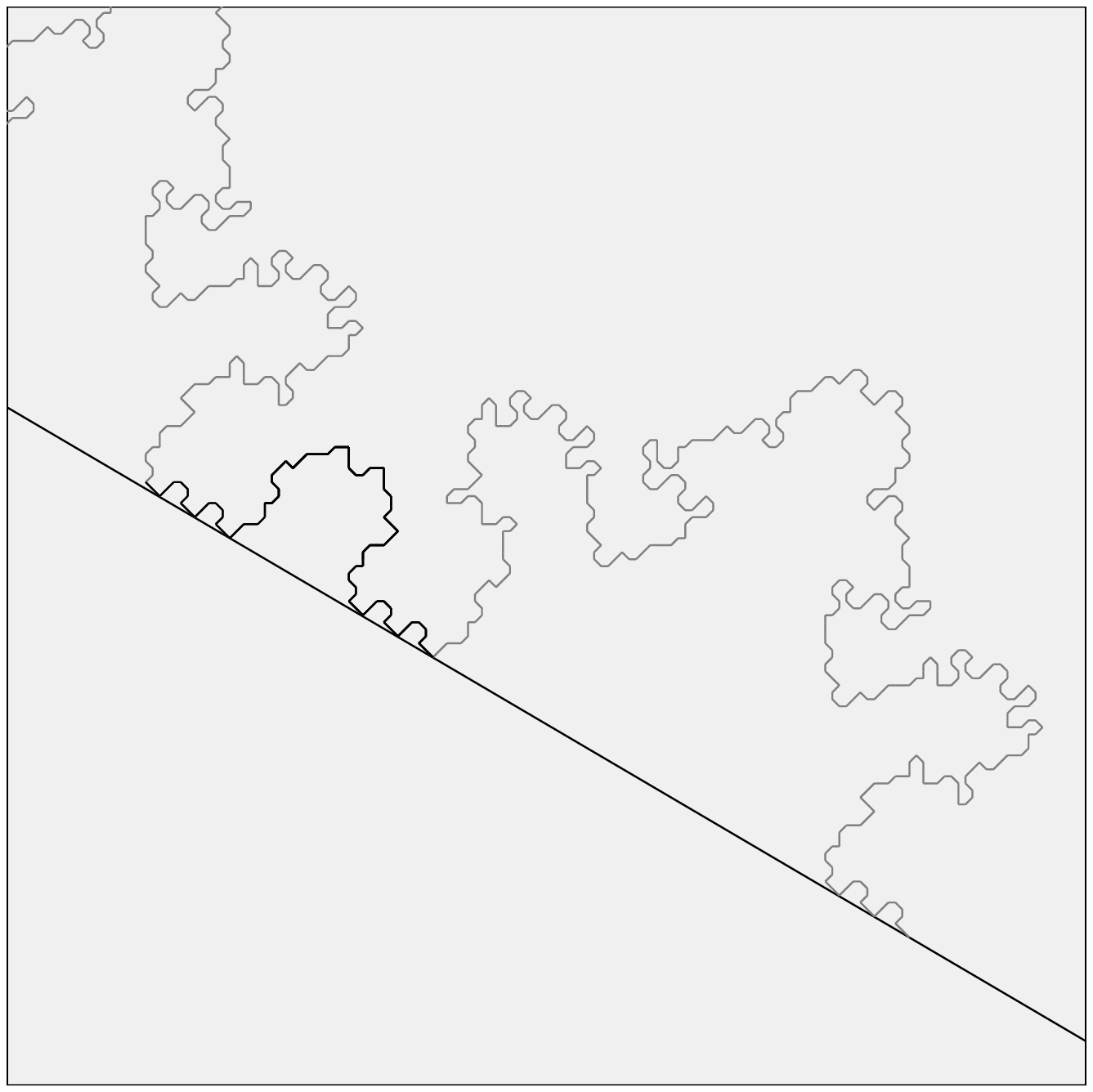}}
\newline
{\bf Figure 23.1:\/} $\Gamma(57/97)$ in grey and $P\Gamma(57/97)$ in black.
\end{center} 

\begin{center}
\resizebox{!}{2.2in}{\includegraphics{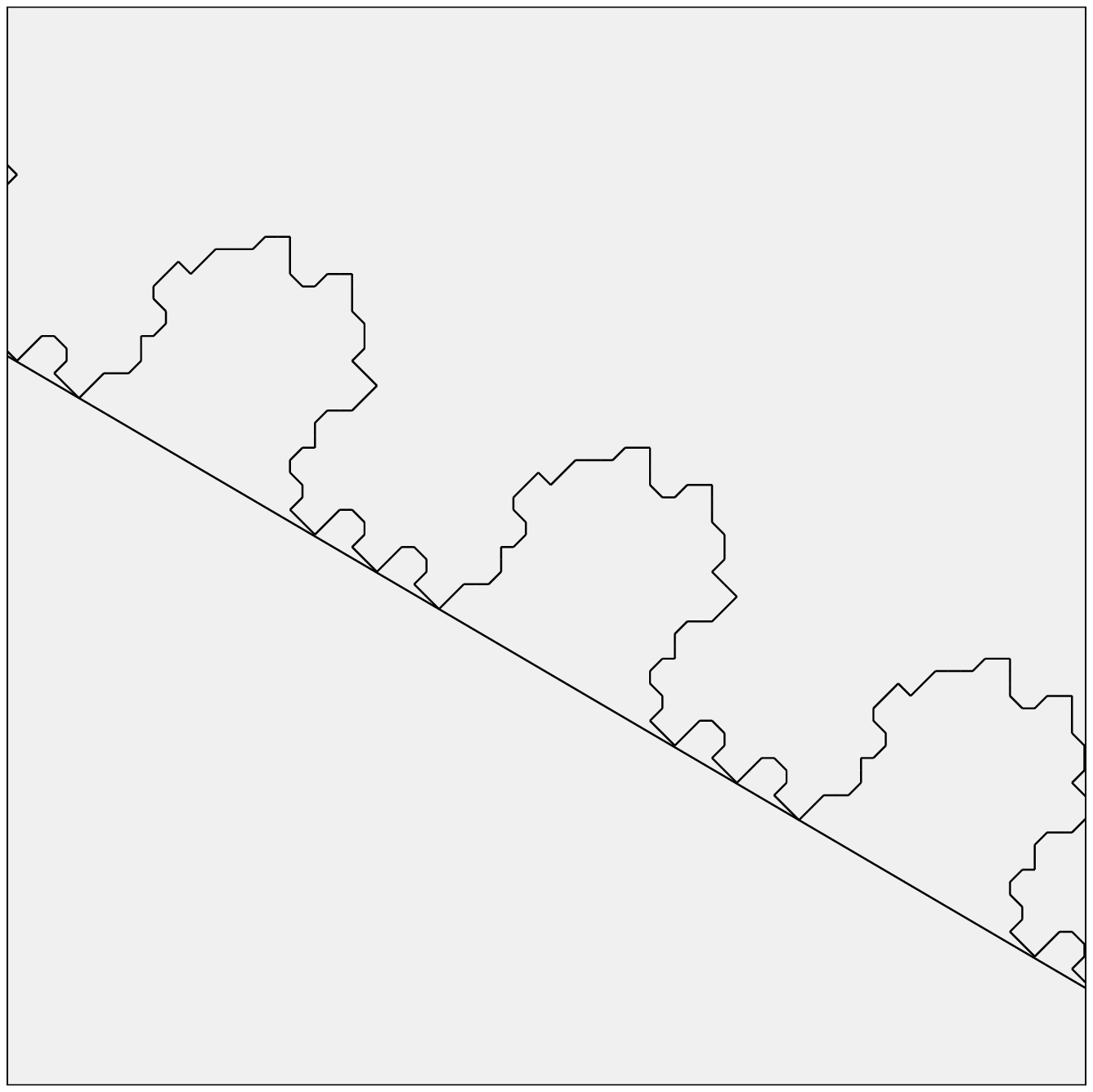}}
\newline
{\bf Figure 23.2:\/} $\Gamma(17/29)$
\end{center}

Now we turn to the statement of the Pivot Theorem.
Given an odd rational parameter $A=p/q$,
let $V$ be the vector from Equation
\ref{boxvectors}.  Let $\Z V$ denote
the group of integer multiples of $V=(q,-p)$.
In Part VI we prove the following result.

\begin{theorem}[Pivot]
Every low vertex of $\Gamma$ is equivalent mod
$\Z V$ to a vertex of $P\Gamma$. That is,
$P\Gamma$ contains one period's worth of low
vertices on $\Gamma$.
\end{theorem}

The Pivot Theorem makes a dramatic
statement. Another way to state the Pivot Theorem is that
there are no low vertices on the complementary
arc $\gamma-P\Gamma$. Here $\gamma$ is the arc
just to the right of $P\Gamma$ such that
$P\Gamma \cup \gamma$ is one full period of
$\Gamma$.  A glance at Figure 23.1 will make
this clear.  
We will prove the Pivot Theorem in Part VI.  We will
also prove the following easy estimate.

\begin{lemma}
\label{pivotbound}
$$
\label{pivotbind}
-\frac{q}{2}<\pi_1(E^-)<\pi_1(E^+)<\frac{q}{2}.
$$
\end{lemma}

\subsection{Half of Theorem \ref{discrete}}
\label{cantorformula}

We will prove that $C_n \subset C_n'$.
This almost an immediate consequence of
the Copy Theorem.   When $1/1 \to A$, the pivot arc $P\Gamma$
contains the points
\begin{equation}
k \widetilde V_1; \hskip 30 pt k=0,...,d_1; \hskip 30 pt
\widetilde V_1=(-1,1).
\end{equation}
This is a consequence of the argument in \S \ref{decompend}.

In general,
suppose $A_1 \leftarrow A_2$ are two parameters.
Then, by construction,
the pivot arc $P\Gamma_2$ contains all
points 
\begin{equation}
\label{copycorr}
v  + k \widetilde V_1 \hskip 30 pt k \in \{0,...,d\}; \hskip 30 pt
d={\rm floor\/}(q_2/2q_1).
\end{equation}
Here $v$ is any vertex of $P\Gamma_1$.
It now follows from induction that $P\Gamma_n$ contains all points
of the form
\begin{equation}
\label{total}
\sum_{j=0}^{n-1} k_j \widetilde V_j; \hskip 30 pt
k_j \in \{0,...,d_j\}.
\end{equation}
Let $M$ denote the map from Equation \ref{funm}. Usually
we take $M$ so that $M(0,0)=0$, but for our proof here,
we adjust so that $M(0,0)=(1/q_n,-1)$.  (This makes no
difference; see the discussion surrounding the definition of
$M$ in \S \ref{ag}.)  Call a lattice
point even if the sum of its coordinates is even.  Note that 
$\widetilde V_j$ is even for all $j$.  Hence, all points in
Equation \ref{total} are even.  The images of these points
under $M$ have second coordinate $-1$.  We just have to
worry about the first coordinate.     We have
\begin{equation}
M(\widetilde V_j)=\frac{1}{q_n}+2|A q_j-p_j|=\frac{1}{q_n}+\frac{1}{q_n}2|p_nq_j-q_np_j|.
\end{equation}
The absolute value in our equation comes from the fact that
$\widetilde V_j=(q_j,-p_j)$ iff $p_j/q_j<A_j$, and
$\widetilde V_j=(-q_j,p_j)$ iff $p_j/q_j>A$.

For convenience, we recall the definition of $C_n$.
Let $\mu_i=|p_nq_i-q_np_i|$.  
\begin{equation}
C_n=\bigcup_{\kappa \in \Pi_n} \Big(X_n(\kappa),-1\Big); \hskip 20 pt
X_n(\kappa)=\frac{1}{q_n} \bigg(1+\sum_{i=0}^{n-1} 2k_i\mu_i\bigg).
\end{equation}
  It now
follows from the affine nature of $M$ and from the definition of
$C_n$ that
\begin{equation}
C_n \subset O_2(1/q_n,-1).
\end{equation}
It follows from the case $n=0$ of Equation \ref{DIO2}, that
$C_n \subset [0,2] \times \{-1\}.$

\subsection{The Inheritance of Low Vertices}

The rest of Theorem \ref{discrete} follows
from the Pivot Theorem and from what we have
done by applying the information contained
in the Pivot Theorem to what we have
already done in the previous section.
To make the argument work, we first need to
deal with a tedious technical detail.
We take care of the detail in this section.

Let $A_1 \leftarrow A_2$ be two odd rationals.
  Let $\widetilde V_1=V_1$ if
$A_1<A_2$ and $\widetilde V_1=-V_1$ if $A_1>A_2$.
Let 
\begin{equation}
d_1={\rm floor\/}(q_2/2q_1).
\end{equation}

Let $v_1$ be a vertex on the pivot arc 
$P\Gamma_1$. Define
\begin{equation}
v_2=v_1+k \widetilde V_1; \hskip 30 pt
k \in \{0,...,d_1\}.
\end{equation}
Notice that these were precisely the vertices that
we considered in \S \ref{cantorformula}.   Now
we want to take a close look at these vertices.
Here is the main result of this section.

\begin{lemma}
\label{low}
$v_1$ is low with respect to $A_1$ iff $v_2$ is low
with respect to $A_2$.
\end{lemma}

\startproof
There are two cases to consider, depending on whether
$A_1<A_2$ or $A_2<A_1$.  We will consider the former
case.  The latter case has essentially the same
treatment.  In our case, we have $\widetilde V_1=V_1$.
Let $E_j^{\pm}$ be the pivot points for $\Gamma_j$.
Say that a vertex is {\it high\/} if it is not low.

We will first suppose that $v_1$ is low with respect
to $A_1$ and that $v_2$ is high with respect
to $A_2$.  This will lead to a contradiction.
We write $v_j=(m_j,n_j)$.  Let $M_j$ be
the fundamental map from Equation \ref{funm}.
Since $v_1$ is
low and $v_2$ is high, we have
$$
2A_1m_1+2n_1+\frac{1}{q_1}=M_1(v_1) \leq 2-\frac{1}{q_1};$$
$$
2A_2m_2+2n_2+\frac{1}{q_2}=M_2(v_2) \geq 2+\frac{1}{q_2}.
$$
Rearranging terms, 
\begin{equation}
\label{case11}
2\bigg(\frac{p_2}{q_2}m_2+n_2\bigg)-2\bigg(\frac{p_1}{q_1}m_1+n_1\bigg) \geq \frac{2}{q_1}.
\end{equation}
Plugging in the relations $m_2=m_1+kq_1$ and $n_2=n_1-kp_1$ and simplifying, we get
\begin{equation}
\label{case22}
\frac{(m_1+kq_1)(p_2q_1-p_1q_2)}{q_1q_2} \geq \frac{1}{q_1}.
\end{equation}
Since $A_1 \leftarrow A_2$ and $A_1<A_2$, we have
\begin{equation}
\label{case33}
p_2q_1-p_1q_2=2.
\end{equation}
Hence
\begin{equation}
\label{precon}
m_1+kq_1 \geq \frac{q_2}{2}.
\end{equation}
Combining Equation \ref{pivotplus} and Equation
\ref{precon}, we get
$$
E_1^+(A_2)=E_1^+(A_1)+d_1q_1 \geq^* m_1+kq_1>\frac{q_2}{2}>^*E_1^+(A_2).$$
This is a contradiction.
the first starred inequality comes from the Pivot 
Theorem and the fact that $k \leq d_1$.  The second
starred inequality comes from the Corollary \ref{pivotbound}.

Now we will suppose that $v_1$ is high with respect to
$A_1$ and $v_2$ is low with respect to $A_2$.  This
will also lead to a contradiction.  Let $M_1$ denote
the first coordinate of the fundamental map relative
to the parameter $A_1$, adjusted so that $M_1(0,0)=1/q_1$. That is
\begin{equation}
M_1(m,n)=2A_1m+2n+\frac{1}{q_1}.
\end{equation}
Since $v_1$ is high, we have the following
dichotomy.
\begin{equation}
M_1(v) \geq 2+\frac{1}{q_1}; \hskip 30 pt
M_1(v)>2+\frac{1}{q_1} \hskip 12 pt
\Longrightarrow \hskip 12 pt
M_1(v) \geq 2+\frac{3}{q_1}.
\end{equation}
We will consider these two cases in turn.
\newline
\newline
{\bf Case 1:\/}
If $M_1(v_1)=2+1/q_1$, then 
$$M_1(m_1,n_1-1)=\frac{1}{q_1}=M_1(0,0).$$
But than $(m_1,n_1-1)=jV_1$ for some integer $j$.
But then $|m_1| \geq q_1$. Since
$v_1 \in P\Gamma_1$, this contradicts 
Corollary \ref{pivotbound}.
Hence $$v_1=(0,1); \hskip 30pt v_2=kV_1+(0,1).$$
If $v_2$ is low then
$$0=2k(A_1 q_1-p_1)<2k(A_2 q_1-p_1)=M_2(v_2)-M_2(0,1) \leq 0.$$
This is a contradiction.  The first inequality comes from
$A_1<A_2$.
\newline
\newline
{\bf Case 2:\/}
If $M_1(v_1) \geq 2+3/q$, then the same reasoning as
in Equations \ref{case11}, \ref{case22}, and \ref{case33} (but with
signs reversed) leads to 
\begin{equation}
m_1+kq_1<-3q_2.
\end{equation}
But then
$$
\frac{-q_2}{2}<\frac{-q_1}{2}<^*m_1 \leq m_1+kq_1<-3q_2
$$
The starred inequality comes from Corollary \ref{pivotbound}.
Again we have a contradiction, this time by a wide margin.
\endproof

\subsection{Proof of Theorem \ref{discrete}}

Now we revisit the construction in \S \ref{cantorformula}
and show that actually $C_n=C_n'$.
 Let $\Lambda_n$ denote the
set of low vertices of $P\Gamma_n$.
By the Pivot Theorem, every low vertex on $P\Gamma_n$ is
equivalent to a point of $\Lambda_n$ modulo $\Z V_n$.

\begin{lemma}
\label{inductlow}
For any $n \geq 0$, we have
$$\Lambda_{n+1}=\bigcup_{k=0}^{d_n}(\Lambda_n+k \widetilde V_n).$$
\end{lemma}

\startproof
Induction.
For $n=0$ we have $E_1^-=(-d_0,d_0)$ and
$E_1^+=(0,0)$.  In this case, the right
hand side of our equation precisely describes
the set of points on 
the line segment joining the pivot points.
The case $n=0$ therefore follows directly
from the Pivot Theorem.

Let $\Lambda_{n+1}'$ denote the right hand side of
our main equation. Since
$\Gamma_n$ is invariant under translation by
$V_n$, every vertex of
$\Lambda_{n+1}'$ is low with respect to
$A_n$. Hence, by Lemma \ref{low}, every
vertex of $\Lambda_{n+1}'$ is low with
respect to $A_{n+1}$. 
 Combining this
fact with Equation \ref{copycorr}, we
see that 
$\Lambda_{n+1}' \subset \Lambda_{n+1}$.

By Lemma \ref{low} again, every $v \in \Lambda_{n+1}$
is also low with respect to $A_n$.  Hence
\begin{equation}
v=v'+k \widetilde V_n; \hskip 30 pt k \in \Z.
\end{equation}
for some $v' \in \Lambda_n$.  If
$k \not \in \{0,...,d_n\}$ then
$v$ either lies to the left of the left
pivot point of $\Gamma_{n+1}$ or to
the right of the right pivot point
of $\Gamma_{n+1}$.  Hence
$k \in \{0,...,d_n\}$.  This proves
that $\Lambda_{n+1} \subset \Lambda'_{n+1}$.
Combining the two facts completes our induction
step.
\endproof

We proved Lemma \ref{inductlow} with respect to the inferior
sequence.  However, notice that if $d_n=1$ then
$\Lambda_{n+1}=\Lambda_n$.  Thus, we get precisely the
same result for consecutive terms in the superior sequence.
We have shown that $v \in \Gamma_n$ is low if and only 
if $v \in \Lambda_n$ mod $\Z V$.  But then
\begin{equation}
O_2(1/q_n,-1) \cap I=M(\Lambda_n); \hskip 30 pt
I=[0,2] \times \{-1\}.
\end{equation}
Here $M$ is the fundamental map. 
Recognizing $\Lambda_n$ as the set from Equation \ref{total},
we get precisely the equality in Theorem \ref{discrete}.
There is one last detail.  One might worry that $M$ maps
some points of $\Lambda_n$ to points on $[0,2] \times \{1\}$,
but all points in $\Lambda_n$ have even parity. Hence,
this does not happen.

This completes the proof of Theorem \ref{discrete}.

\subsection{Proof of Lemmas \ref{returncomb} and \ref{coarsefar}}

Let $\Sigma_n$ denote the union of all points in
Equation \ref{total}.  Here $M(\Sigma_n)=C_n$.
The ordering on $\Sigma_n$ determines the ordering
of the return dynamics to $C_n$.
We set $\Sigma_0=\{(0,0)\}$, for convenience.
We can determine the ordering on $\Sigma_{n+1}$
from the ordering on $\Sigma_n$ and the sign of
$A_{n+1}-A_n$.   When $A_n<A_{n+1}$, we can write the
relation
\begin{equation}
\Sigma_n + k V_n \prec \Sigma_n + (k+1) V_n; \hskip 30 pt
k=0,...,(d_n-1).
\end{equation}
to denote that each point in the left hand set precedes
each point on the right hand set.  Within each set, the
ordering does not change.   
When $A_n>A_{n+1}$, we can write the relation
\begin{equation}
\Sigma_n - (k+1) V_n \prec \Sigma_n - k V_n; \hskip 30 pt
k=0,...,(d_n-1).
\end{equation}
Lemma \ref{returncomb} follows from these facts, and
induction.

Let $\beta_n$ denote the arc of $P\Gamma_n$, chosen
so that $P\Gamma_n \cup \beta_n$ is one period
of $P\Gamma_n$.
Let $L_n$ be the line of slope $-A_n$ through the origin.

\begin{lemma}
No point of $\beta_m$ lies more
than $q_m$ vertical units away from $L_m$ and
some point of $\beta_m$ lies at least $q_m/4$ vertical 
units away
from $L_m$. 
\end{lemma}

\startproof
By the Room Lemma,
$\beta_m \subset R(A_m)$.
The upper bound follows immediately from this containment.
For the lower bound,
recall from the Room Lemma that $P\Gamma_m$ crosses
the centerline $L$ of $R(A_m)$ once, and
this crossing point lies at least $(p_m+q_m)/4>q_m/4$ vertical units from $L_m$.
By Lemma \ref{pivotbound} and symmetry,
the left endpoint of $\beta_m$ lies to the left
of $L$ and the right endpoint of $\beta_m$ lies to the
right of $L$.  Hence, $\beta_m$ contains the crossing point we have
mentioned.   For an alternative argument, we note that
no point on the pivot arc crosses the line parallel to
the floor and ceiling of $R(A_m)$ and halfway between
them, whereas the crossing point lies above this midline.
\endproof

Notice that the line $L_n$ replaces the line $L_m$ in
our next lemma.

\begin{lemma}
\label{riseee}
Let $m \leq n$ and $q_m>10$. Then
some point of $\beta_m$ lies at least
$q_m/4-1$ vertical units from $L_n$.
Moreover, no point of $\beta_m$ lies more than
than $q_m+1$ vertical units away from $L_n$.
\end{lemma}

\startproof
Some point $v$ of $\beta_m$ at least $q_m$ vertical units from
$L_m$ by the previous result.
From Lemma \ref{superior}, we have
\begin{equation}
\label{dioX}
|A_m-A_n|<\frac{2}{q_m^2}.
\end{equation}
On the other hand, by the Room Lemma and by construction,
$P\Gamma_m$ is contained in two
consecutive translates of $R(A_m)$, one of which
is $R(A_m)$ itself.   Hence, $P\Gamma_m$ lies entirely
inside the ball $B$ of radius $4q_m$ about the origin.
By Equation \ref{dioX}, the the Hausdorff distance
between the segments
seqments $L_m \cap B$ and $L_n \cap B$ is less than $1$ once
$m>10$.
By construction, the vertical line segment starting at $v$ and
dropping down $q_m-1$ units is disjoint from $L_n \cap B$.
But this segment is disjoint from $L_n-B$ as well.  
Hence $v$ is at least $q_m/2-1$ vertical units from $L_n$.
The upper bound has a similar proof.
\endproof

\begin{lemma}
\label{length}
$\beta_m$ has length at most $5q_m^2$.
\end{lemma}

\startproof
$\beta_m$ is contained in one period of $P\Gamma_m$.
Hence, it suffices to bound the length of any one
period of $P\Gamma_m$.  By the Room Lemma,
one such period is contained in $R(A_m)$.  We
compute easily that the area of $R(A_m)$ is much
less than $5q_m^2$.  Hence, there are less than
$5q_m^2$ vertices in $R(A_m)$.  Hence, the
length of one period of $P\Gamma_m$ is less than
$5q_m^2$.
\endproof

Suppose now that $\kappa$ and $\kappa_+$ are two consecutive
points on $\Sigma_n$.  We want to understand the
arc of $P\Gamma_n$ that joins these points.  Suppose that $\sigma(\kappa)=m$.
It follows from induction and from the Copy Theorem that
there is some translation $T$ such that $T(\kappa)$ and
$T(\kappa_+)$ are the endpoints of the arc $\beta_m$.
The arc joining $\kappa$ to $\kappa_+$ has the
same length as $\beta_m$, and this length is less than
$5q_m^2$.  This gives us the estimate for $h_2$.

Now we deal with $h_1$.  We check the result by
hand for $q_n<10$.  So, suppose that $q_n>10$.
All the vertices $\kappa$, $\kappa_+$, $T(\kappa)$, and
$T(\kappa_+)$ lie within $1$ vertical unit of
the baseline $L_n$.
We know that the vertical
distance from some point of $\beta_m$ to $L_n$ is
at least $q_m/2-1$.  Hence, the vertical distance
from some point on $T(\beta_m)$ to $L_n$ is at
least $q_m/2-2$.   Similarly, the vertical
distance from any point of $\beta_m$ to
$L_n$ is at most $q_m+2$. 
If two points in $\Z^2$ have
vertical distance $d$ then the images of these
points under the fundamental map $M_n$ have
horizontal distance $2d$.  In short, the
fundamental map doubles the relevant distances.
This fact gives us our estimate on $h_1$.

This completes the proof of Lemma \ref{coarsefar}.

\subsection{Theorem \ref{discrete} in the Even Case}
\label{extra}

Here we discuss Theorem \ref{discrete} in the even case.
For each even rational $A_1 \in (0,1)$ there is a unique
odd rational $A_2$ such that (in the language of
Equation \ref{induct0}) $A_1=(A_2)_{\pm}$ and
$q_2<2q_1$.  In Lemma \ref{welldefined} we will
show that $\Gamma_1$ (a closed polygon) contains a copy of
$P\Gamma_2$, and all low vertices of
$\Gamma_1$ lie on this arc.  
From this fact, we see that

\begin{equation}
O(1/q_1,-1)=M_1(\Sigma_1),
\end{equation}
just as in the odd case.  Here $M_1$ is the fundamental
map defined relative to the parameter $A_1$ and
$\Sigma_1$ is the set of low vertices on $P\Gamma_1$.

Note that $\Sigma_1=\Sigma_2$, where $\Sigma_2$ is the
set of low vertices on $P\Gamma_2$.  The only difference
between the two sets $M_1(\Sigma_1)$ and
$M_2(\Sigma_2)$ is the difference in the maps $M_1$
and $M_2$.  Now we explain the precise form of 
Theorem \ref{discrete} that this structure entails.

Switching notation, let
$A$ be an even rational.  One of the two rationals
$A_{\pm}$ from Equation \ref{induct0} is odd, and
we call this rational $A'$.  
We can find the initial part of a superior sequence
$\{A_k\}$ such that $A'=A_{n-1}$.  We set
$A=A_n$ even though $A$ does not belong to this sequence.
Referring to Theorem \ref{discrete}, we define $\Pi_n$
exactly in the odd case, but for one detail.  In
case $2q'>q$, we simply ignore the $n$th factor of
$\Pi_n$. That is, we treat $q'$ as an inferior term.
With these changes, Theorem \ref{discrete} goes through
word for word.

Here we give an example. Let $A_1=12/31$.  Then
$A_2=19/49$, exactly is in the introduction.  We have
$n=3$ and our sequence is
$$\frac{p_0}{q_0}=\frac{1}{1}, \frac{1}{3}, \frac{5}{13}, \frac{12}{31}=\frac{p_3}{q_3}.$$
All terms are superior, so this is also the superior sequence.
$n=3$ in our example, and the renormalization sequence is $1,2,1$.
The $\mu$ sequence is $19,5,1$.   The first coordinates of the
 $12$ points of 
$O_2(1/49) \cap I$ are given by 
$$\bigcup_{k_0=0}^1\  \bigcup_{k_1=0}^2\  \bigcup_{k_2=0}^1 \frac{2(19k_0+5k_1+1k_2)+1}{31}.$$
Writing these numbers in a suggestive way, the union above works out to
$$\frac{1}{31} \times \big(1\hskip 8pt 3\hskip 20 pt    11\hskip 8pt 13\hskip 20pt   21\hskip 8pt 23 \hskip 60 pt   
39\hskip 8pt 41\hskip 20 pt    49\hskip 8pt 51\hskip 20 pt    59\hskip 8pt 61\big).$$

\subsection{A Conjectural Extension}
\label{gapmap}

Let $C_n$ be the set from Theorem \ref{discrete}.  
Each $\xi \in C_n$ is the midpoint of a
special interval, in the sense of \S \ref{definebasic}. Call this interval $J(\xi)$.
Define
\begin{equation}
\widehat C_n=\bigcup_{\xi \in C_n} J(\xi).
\end{equation}

Figure 23.1 shows three examples.
In the picture, we have thickened
the intervals to get a better picture.  We have also added
in the white bars to clarify the spacing.  

\begin{center}
\resizebox{!}{2in}{\includegraphics{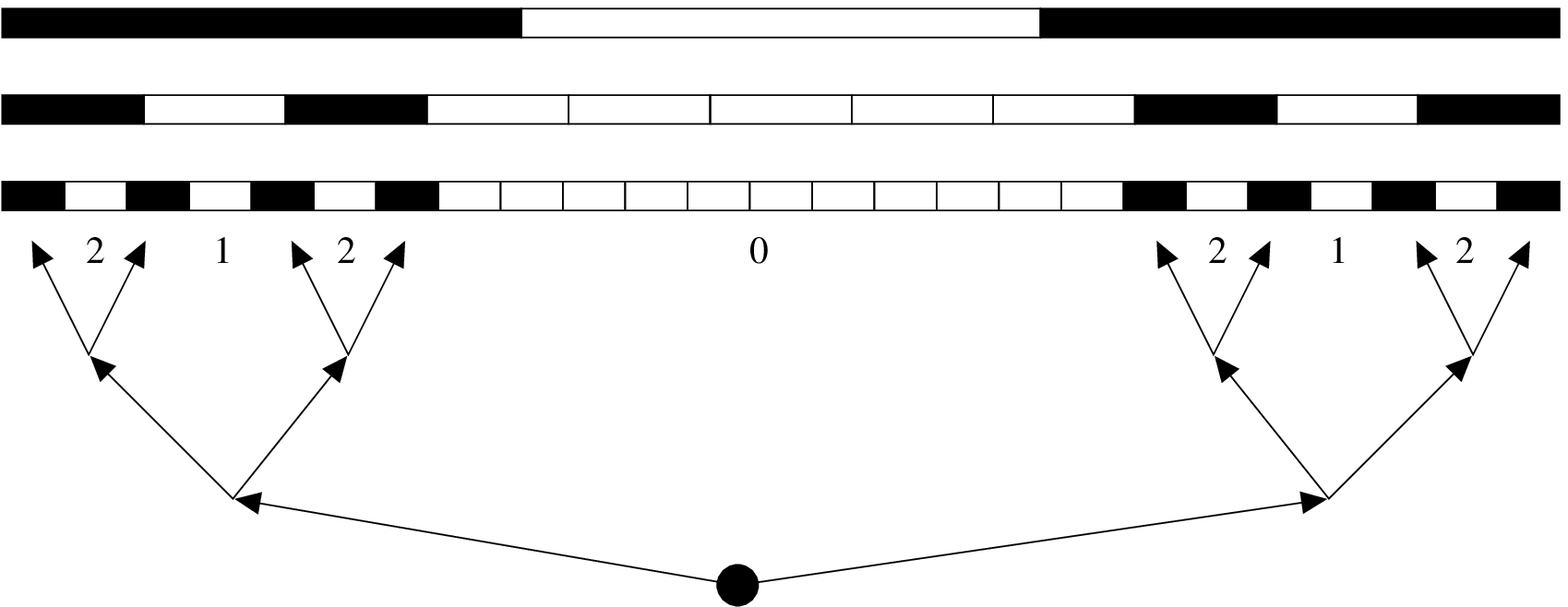}}
\newline
{\bf Figure 23.1:\/} $\widehat C(A)$ for $A=1/3$ and $3/11$ and $7/25$.
\end{center}

The three rationals in Figure 30.1 are part of a superior sequence,
one can see that each picture sort of refines the one above it.
It is a consequence of Lemma \ref{parity} that, in the odd case,
there is a gap between every pair of intervals in $\widehat C(A)$.
In the even case, this need not be true.  One can compute
the positions of the intervals using the formula in
Theorem \ref{discrete}. 

Say that a {\it gap\/} is an maximal interval of
$I-\widehat C$.  For $\widehat C(7/25)$ there are $7$ gaps.
Each gap has a {\it level\/}, as indicated in the figure.
The levels go from $0$ to $n-1$ in $\widehat C_n$.
Informally, the gaps of level $k \leq n-2$ are
inherited from simpler rationals, and the gaps of
level $n-1$ are newly created with the new parameter.
Generally speaking, the higher-level gaps are smaller,
but this need not be the case. For $\widehat C(7/25)$,
the gaps of level $1$ and $2$ have the same size.

Given this notion of levels there is a natural identification
of $C_n$ with the ends of a directed finite binary tree.  The return
map $\Theta: C_n \to C_n$ comes from an automorphism
of this tree.  The union of all the gaps is bijective
with the forward cones of the tree.  The automorphism
of the tree induces a bijection on its forward cones.

\begin{conjecture}
\label{gapmap2}
The outer billiards map is entirely defined on the interior
of any gap, and the return map to the interval $I$ is
naturally conjugate to the map on the forward cones
induced by the tree automorphism.
\end{conjecture}

Some reflection will convince the reader that this is the
simplest possible answer to the question of what happens
in the gaps.  Our Inheritance Lemma from \S \ref{period}
makes some progress in proving this conjecture, but
it doesn't have quite enough juice in it.

\newpage

\section{Most of The Comet Theorem}
\label{comet1}

\subsection{Preliminaries}

In this chapter we prove Statements 2,3,4 of
the Comet Theorem.  
We defer the proof
of Statement 1 until Part VI.  Statement 4
assumes the truth of Statement 1, but our proof of
Statement 1, given in Part VI, does not
depend on Statement 4.  (That is, our
argument isn't circular.)

Suppose that $\{A_n\}$ is the superior
sequence approximating some irrational $A$.
Let $C_n$ be the set from Theorem \ref{discrete}.
Let $U_A$ and $I$ be as in the Comet Theorem.
We also prove the following result in Part VI.

\begin{theorem}[Period]
For any $\epsilon>0$ there is an $N>0$ with the following
property.  If $\zeta \in I$ is more than
$\epsilon$ units from $C_n$, then the period of $\zeta$ is
at most $N$.  The constant $N$ only depends on $\epsilon$.
\end{theorem}

\begin{corollary}
\label{trim2}
$U_A \cap I \subset C_A$.
\end{corollary}
\startproof
We will suppose that $U_A$ contains a point
$\zeta \not \in C_A$ and derive a contradiction.
By compactness,
there is some $\epsilon>0$ such that
$\zeta$ is at least $3\epsilon$ from any point
of $C_A$.  Since $C_A$ is the geometric
limit of $C_n$, we see that there is some
$N_1$ such that $n>N_1$ implies that $\zeta$
is at least $2\epsilon$ from $C_n$.

Let $\{\zeta_n\} \in I$ be a sequence of points
converging to $\zeta$.  We can choose these
points so that the orbit of $\zeta_n$ relative to
$A_n$ is well defined.  There is a constant
$N_2$ such that $n>N_2$ implies that
$\zeta_n$ is at least $\epsilon$ from $C_n$.
But then, by the Period Theorem, there is
some $N_3$ such that the period of
$\zeta_n$ is at most $N_3$.  

On the other hand, by the Continuity Principle,
the arithmetic graph $\Gamma(\zeta_n,A_n)$
converges to the arithmetic graph
$\Gamma(\zeta,A)$. In particular, the period
of $\Gamma(\zeta_n,A_n)$ tends to $\infty$.
This is a contradiction.   Hence, 
$\zeta$ cannot exist.  
\endproof

Let $\Pi_A$ be the sequence space from
\S \ref{productspace}.   Say that
two sequences in $\Pi_A$ are
equivalent if they have the
same infinite tail ends.  Given the
nature of the odometer map, we have the following
useful principle.
\newline
\newline
{\bf Odometer Principle:\/}
Any two equivalent sequences are in the
same orbit of the odometer map.  
Call this the {\it odometer principle\/}.
We will use this principle several times in our proofs.

\subsection{Overview of the Proof}

We first prove a preliminary version of the Comet Theorem.
Let 
\begin{equation}
C_A'=C_A-(2\Z[A] \times \{-1\}).
\end{equation}

\begin{theorem}
\label{precomet}
Let $U_A$ denote the set of unbounded special orbits
relative to an irrational $A \in (0,1)$.
\begin{enumerate}
\item $C_A' \subset U_A$.
\item 
The first return map $\rho_A: C_A' \to C_A'$ is defined precisely
on $C_A'-\phi(-1)$. The map $\phi^{-1}$ conjugates $\rho_A$
to the restriction of the odometer on ${\cal Z\/}_A$.
\item  For any $\zeta \in C_A'-\phi(-1)$, the
orbit-portion
between $\zeta$ and $\rho_A(\zeta)$ has excursion
distance in $$\bigg[\frac{d^{-1}}{2}-4,2d^{-1}+20\bigg]$$ and length in
$$\bigg[\frac{d^{-2}}{32}-\frac{d^{-1}}{4},100d^{-3}+100d^{-2}\bigg].$$
 Here $d=d(-1,\phi^{-1}(\zeta))$.
\end{enumerate}
\end{theorem}

\noindent
{\bf Remarks:\/}
\newline
(i)
Our constants in Item 3 are not optimal; some tedious
elementary arguments would improve them.
\newline
(ii)  Since $d^{-1} \geq 1$, the estimates
in Item 3 above imply the less precise
estimates in the Comet Theorem -- once we
establish that $C_A^{\#}=C_A'$. \newline
(iii) As we remarked after the Comet Theorem, the only
non-sharp bound in Item 3 is the length upper-bound.
For instance, our proof in [{\bf S1\/}], which
establises a kind of coarse self-similarity structure,
would give a better bound for $A=\sqrt 5-2$ if carefully
examined.
We conjecture that $-3$ is the best bound that
works for all parameters at once.
\newline

Next, we prove a double identity.
\begin{lemma}
\label{doubleid}
$U_A \cap I = C_A^{\#}=C_A-(2\Z[A] \times \{-1\})$.
\end{lemma}
Statements 2 and 3 of the Comet Theorem
follow from this result
and Lemma \ref{precomet}. 
Lemma \ref{doubleid}  
also contains the first claim
in Statement 4 of the Comet Theorem.

At the end of the chapter, we will prove
the second claim made in Statement 4 of
the Comet Theorem.

\subsection{The Cantor Set}
\label{pf1}

We first need to resolve the technical point that our set
$C_A$ is actually well defined.  For convenience,
we repeat the definition.
\begin{equation}
\label{cantor2}
C_A=\bigcup_{\kappa \in \Pi} \Big(X(\kappa),-1\Big); \hskip 30 pt
X(\kappa)=\sum_{i=0}^{\infty} 2k_i|Aq_i-p_i|.
\end{equation}

\begin{lemma}
The infinite sums in Equation \ref{cantor2} converge.
Hence $C_A$ is well defined.
\end{lemma}

\startproof
Combining Equation \ref{DIO1} with the bound $0 \leq k_n<d_n$, we
see that the $n$th term in the sum defining $X(\kappa)$ is at most
$2q_n^{-1}$.  Given that $2q_k<q_{k+1}$ for all $k$, we get
$2q_n^{-1}<2^{-n+1}$.  The sequence defining $X(\kappa)$ 
decays exponentially and hence
converges. 
\endproof

For the purposes of this section we equip the product
space $\Pi$ with the lexicographic ordering and
the product topology.

\begin{lemma}
The map $X: \Pi \to C_A$ is a
homeomorphism that maps the lexicographic
order to the linear order. Hence
$C_A$ is a Cantor set.
\end{lemma}

\startproof
We first show that the map $X$ is injective.
In fact, we will show that $X$ is order preserving.
If $\kappa=\{k_i\} \prec \kappa'=\{k_i'\}$ in the lexicographic
ordering, then there is some smallest index $m$ such that
$k_i=k_i'$ for all indices $i=0,...,(m-1)$ and $k_m<k_m'$.
Let $\lambda_m=|Aq_m-p_m|$, as in Equation \ref{DIO1}. Then
\begin{equation}
X(\kappa')-X(\kappa) \geq 2\lambda_m- \sum_{k=m+1}^{\infty}2d_k\lambda_k=\lambda_m-\lambda'_{m+1}>0
\end{equation}
by Equation \ref{DIO3}.

The map $X: \Pi \to [0,2]$ is continuous with
respect to the topology on $\Pi$, because
the $n$th term in the sum defining $X$ is
always less than $2^{-n+1}$.  We also know
that $X$ is injective.  Hence, $X$ is bijective
onto its image.
Any continuous bijection from
a compact space to a Hausdorff topological space
is a homeomorphism.
\endproof

\subsection{Convergence of the Fundamental Orbit}

Let $\{p_n/q_n\}$ denote the superior sequence associated to $A$.
We use the notation from the previous chapter.  Here
$\Gamma_n$ denotes the corresponding arithmetic graph and
\begin{equation}
C_n=\bigcup_{\kappa \in \Pi_n} \Big(X_n(\kappa),-1\Big); \hskip 30 pt
X_n(\kappa)=\frac{1}{q_n}+\sum_{i=0}^{n-1} 2k_i|A_nq_i-p_i|.
\end{equation}
We have already proved that $C_n \subset O_2(1/q_n,-1)$.

Let $\kappa \in \Pi$ be some infinite sequence.
Let $\kappa_n \in \Pi_n$ be the truncated sequence.
Let 
\begin{equation}
\sigma_n=(X_n(\kappa_n),-1); \hskip 30 pt \sigma=(X(\kappa),-1)
\end{equation}
Here is our basic convergence result.

\begin{lemma}
$\sigma_n \to \sigma$ as $n \to \infty$.
\end{lemma}

\startproof
For $i<n$, let $\tau_{i,n}$ denote the $i$th term in the
sum for $X_n(\kappa_n)$.  Let $\tau_n$ be the corresponding
term in the sum for $X(\kappa)$.  When we construct
the superior sequence, we will see that
the sign of $A-A_i$ is the same as the sign of
$A_n-A_i$.  Therefore
\begin{equation}
|\tau_n-\tau_{i,n}|=2k|A-A_n|q_n<2q_n^{-1}<2^{-n+1}.
\end{equation}
Therefore
$$
|X(\kappa)-X(\kappa_n)|=\sum_{i=0}^{n-1}|\tau_n-\tau_{i,n}|+
\sum_{i=n}^{\infty} \tau_i<
$$
\begin{equation}
2\sum_{i=0}^{n-1}2^{-n}+2\sum_{i=n}^{\infty} 2^{-i}<2^{n-3}.
\end{equation}
This completes the proof.
\endproof

\subsection{All but the Last Sequence}

  We call the sequence
$\{k_i\}$ {\it first\/} if $\widetilde k_i=0$ for all $i$
and {\it last\/} of $\widetilde k_i=d_i$ for all $i$.
The map
$\phi_2: \Pi_A \to C_A$ is a homeomorphism.
Using $\phi_2$, we transfer the notions of
{\it first\/} and {\it last\/} to points of $C_A$.

Let $\zeta \in C_A'$ denote a point that is not last.
Let $\kappa$ denote the corresponding sequence
in $\Pi_A$.
Say that two sequences in $\Pi$ are {\it equivalent\/}
if they have the same infinite tail end.  We can
define the reverse lexicographic order on any equivalence.
Likewise we can extend the twirl order to any 
equivalence class.
In particular, we extend the twirl order to the
equivalence class of $\kappa$, the sequence currently
of interest to us.

Since $\kappa$ is not last, we can find some
smallest index $m=m(\zeta)$ such that
where $\widetilde k_m<d_i$.  In other words,
$m$ is the smallest index such that $\kappa$
differs from the last sequence in the $m$th spot.

The successor $\kappa_+$ of $\kappa$ is obtained by
incrementing $\widetilde k_m$ by $1$ and setting
$\widetilde k_i=0$ for all $i<m$.  This notion of
successor is compatible with the twirl ordering
on the finite truncations $\Pi_n$. Define
\begin{equation}
\zeta_+=(X(\kappa_+),-1); \hskip 30 pt
(\zeta_n)_+=(X(\kappa_n)_+,-1).
\end{equation}

\begin{lemma}
\label{freturn}
Let $\zeta \in C_A'$ be a point that is not last.
Let $m=m(\zeta)$.
The forward $\Psi$ orbit of $\zeta$ returns to
$C_A$ as $\zeta_+$ in at most
$5q_m^2$ steps.  Along the way, this portion
of the orbit wanders between $q_m/2-2$ units and
$2q_m+2$ units away from $(0,-1)$.
\end{lemma}

\startproof
By Lemma \ref{irrat}, the orbit of $\zeta$ is well-defined.
Referring to the notation in Lemma \ref{coarsefar}, we
get $\sigma(\kappa_n)=m$ for $n$ large enough.
Hence the forward $\Psi_n$ orbit of $\zeta_n$ returns
to $(\zeta_n)_+$ after at most $5q_m^2$ steps, moving
away from $(0,-1)$ by at least $q_m/2-2$ units and
at most $2q_m+2$ steps. Here $m$ is
independent of $n$.
Since $X$ is continuous, we have
$(\zeta_n)_+ \to \zeta_+$ as $n \to \infty$.
The Contintuity Principle implies that
the forward $\Psi$ orbit of $\zeta$ returns as
$\zeta_+$ after at most $5q_m^2$ steps, moving
away from $(0,-1)$ at least $q_m/2-2$ units and
at most $2q_m+2$ steps.
\endproof

There is an entirely analogous result for the
backwards return map.  This analogous result holds
for all but the first point.

\subsection{Statement 1 of Lemma \ref{precomet}}
\label{st2proof}

We call a sequence of $\Pi_A$ {\it equivalent-to-first\/} if it
differs from the first sequence in only a finite number
of positions.  We call a sequence {\it equivalent-to-last\/}
if it differs from the last sequence in a finite
number of positions.  As in the previous section,
we transfer these notions to $C_A$.

\begin{lemma}
No sequence in $\Pi_A$ is both equivalent-to-first
and equivalent-to-last.
\end{lemma}

\startproof
This is immediate from the definitions.
\endproof

Let $\zeta$ be a point in $C_A'$ that is not equivalent-to-last.
We will show that the forwards orbit of
$\zeta$ is unbounded.
Let $m=m(\kappa)$ be as in the proof of Lemma \ref{freturn}.
Lemma \ref{coarsefar} says that the portion
of the orbit between $\zeta$ and $\zeta_+$ 
wanders at least $q_m/5$ from the origin.  Since
we can achieve any initial sequence we like with
iterated successors of $\kappa$, we can find
iterated successors $\kappa'$ of $\kappa$
such that $m(\kappa')$ is as large as we like.
But this shows that the forwards orbit of
$\zeta$ is unbounded.  Here we are using the
fact that $\lim_{m \to \infty} q_m=\infty$.
This shows that $\zeta$ has an unbounded
forwards orbit.  

  Essentially the same
argument works for the backwards orbit of points
that are not equivalent-to-first.  This
establishes Statement 1.

\subsection{Statement 2 of Lemma \ref{precomet}}

The successor map on $\Pi_A$ is defined
except on the last sequence $\kappa$ of $\Pi_A$.
Referring to the homeomorphism
$\phi_1$ given in Equation \ref{phi1}, we
have $$\phi_1(-1)=\kappa.$$
Thus, the point $\phi_2(\kappa) \in C_A$
corresponding to $\kappa$ is precisely
$\phi(-1)$.
By Lemma \ref{freturn}, the return map
$\rho_A: C_A' \to C_A'$ is defined
on $C_A'-\phi(-1)$.

The map $\phi_1$ conjugates the odometer
map on ${\cal Z\/}_A$ to the successor map on $\Pi_A$.
Combining this fact with Lemma \ref{freturn},
we see that $\phi^{-1}$ conjugates $\rho_A$
to the restriction of the odometer map
on ${\cal Z\/}_A$.

It remains to understand what happens to
the forward orbit of $x=\phi(-1)$, in case
$x \in C_A'$.  The following result
completes the proof of Statement 2.

\begin{lemma}
\label{noreturn}
If $x \in C_A'$ then the forward orbit of $x$
does not return to $C_A'$.
\end{lemma}

\startproof
Suppose that the forward orbit of $x$ returns
to $C_A'$ after $N$ steps.
Since outer billiards
is a piecewise isometry, there is
some open neighborhood $U$ of
$x$ such that every point
of $C_A' \cap U$ returns to 
$C_A'$ in at most $N$ steps.  But
there is some uniformly small $m$ such that
every point $\zeta \in C_A'-U$ differs from
the last sequence $\kappa$ at or before
the $m$th spot.  Lemma \ref{freturn}
says that such points return to
$C_A'$ in a uniformly bounded number of
steps. In short, all points of
$C_A'$ return to $C_A'$ in a uniformly
bounded number of steps. But then
all orbits in $C_A'$ are bounded.
This is a contradiction.
\endproof

\subsection{Statement 3 of Lemma \ref{precomet}}

Let $\zeta \in C_A'$.
Let $O_{\zeta}$ denote the portion of
the forward outer billiards orbit of $\zeta$ between $\zeta$ and
$\rho_A(\zeta)$.   We mean to use the original
outer billiards map $\psi'$ here.
 Let $m$ be such that
\begin{equation}
d(\phi^{-1}(\zeta),-1)=q_m^{-1}.
\end{equation}
By definition 
$\phi^{-1}(\zeta)$ and $-1$ disagree
by $\Z/D_{m+1}$, but agree in
$\Z/D_k$ for $k=1,...,m$.
In case $m=0$, the points
$\phi^{-1}(\zeta)$ and $-1$ already disagree in
$\Z/D_1$.
  Let $\kappa \in \Pi_A$
denote the sequence corresponding to $\zeta$.

\begin{lemma}
$\sigma(\kappa)=m$,
\end{lemma}

\startproof
Let $\lambda$ be the sequence corresponding
to $\phi(-1)$.  Then $\lambda$ is the last
sequence in the twirl order.  The sequences
$\kappa$ and $\lambda$ agree in positions
$k=0,...,m-1$ but then disagree in position
$m$.  When $m=0$, the sequences already
disagree in position $0$.   This is to
say that $m$ is the first index where
$\kappa$ disagrees with the last sequence
in the twirl order.  But then,
$\kappa$ and $\kappa_+$ disagree
in positions $0,...,m$ and agree
in position $k$ for $k>m$.
\endproof

\begin{lemma}
$O_{\zeta}$ has excursion distance between
$q_m/2-4$ and $2q_m+20$,
\end{lemma}

\startproof
Lemma \ref{freturn} tells us that the
$\Psi$-orbit of $\zeta$ between $\zeta$
and $\rho_A(\zeta)$ wanders between
$q_m/2-4$ and $2q_m+4$ units from the origin.
Here are interested in the full outer billiards
$O_{\zeta}$.
Since the $\Psi$ orbit of $\zeta$ between
$\zeta$ and $\rho_{\zeta}$ is a subset
of $O_{\zeta}$, the lower bound on the
excursion distance is an immediate
corollary of the lower bound from
Lemma \ref{freturn}.

The upper bound follows from a simple geometric
analysis of the Pinwheel Lemma.
Looking at the proof of 
the Pinwheel Lemma, we see the following geometry.
Starting at a point on $\Xi$ that is $R$ units from
the origin, the $\psi$-orbit remains within 
$2R+8$ units of the origin before returning
to $\Xi$. 
Essentially, the $\psi$-orbit follows
an octagon once around the kite before returning,
as shown in Figure 7.3.   The constant
of $10$ takes care of the small deviations from
the path in Figure 7.3, as discussed in \S \ref{4sharp}
and \S \ref{6flat}.  

Recall that $\psi$ is the square of the outer billiards
map $\psi'$.   Since $\psi'$ is always reflection in a vertex
that is within $1$ unit of the origin, we see that
the entire $\psi'$ orbit of interest to us is at
most $2R+12$ units from the origin.  Hence,
the portion of the outer billiards orbit of
interest to us wanders at most $2(q_m+4)+12=2q_m+20$
units from the origin.
\endproof

\begin{lemma}
$O_{\zeta}$ has length at most $100q_m^3+100q_m^2$.
\end{lemma}

\startproof
We know that the $\Psi$ orbit of $\zeta$ between
$\zeta$ and $\rho_A(\zeta)$ has length at most
$5q_m^2$.  Examining the proof of the Pinwheel
Lemma, we see that the a point on $\Xi$
that is $R$ units from the origin returns to
$\Xi$ in less than $10R$ iterates.    Given our bound
of $R=2q_m+2$, we see that
the orbit $O_{\zeta}$ is at most $20q_m+20$ times
as long as the corresponding $\Psi$-orbit.
This gives us a length bound of $100q_m^3+100q_m^2$.
\endproof

\begin{lemma}
$O_{\zeta}$ has length at least 
$q_m^2/32-q_m/4$.
\end{lemma}

\startproof
Some point in the $\Psi$-orbit of $\zeta$ between
$\zeta$ and $\rho_A(\zeta)$ lies at least
$q_m/2-4$ vertical units from the origin.
Consecutive iterates in the $\Psi$-orbit
have vertical distance at most $4$ units
apart.  Hence, there are at least
$q_m/8-1$ points in the $\Psi$-orbit
that are at least $q_m/4$ horizontal
units from the origin.  Inspecting the
Pinwheel Lemma, we see that the length of
the $\psi'$-orbit between two such points
is at least $q_m/4$.   Hence, 
$O_{\zeta}$ has length at least 
$q_m^2/32-q_m/4$.
\endproof

This completes the proof of Statement 3.

\subsection{Proof of Lemma \ref{doubleid}}

\begin{lemma}
\label{trim0}
No point of $C_A-C_A^{\#}$ has a well-defined orbit.
\end{lemma}

\startproof
Call a sequence in $\Pi_A$ {\it equivalent-to-trivial\/} if
either differs from the $0$ sequence by a finite number
of terms, or it differs from the sequence
$\{d_i\}$ by a finite number of terms. The
homeomorphism $\phi_2$ bijects the equivalent-to-trivial
points in $\Pi_A$ to 
$C_A-C_A^{\#}$.

Suppose first that the superior sequence for $A$ is
not eventually monotone.
In this case, an equivalent-to-trivial sequence is
neither equivalent-to-first nor equivalent-to-last.
See \S \ref{st2proof} for definitions of these terms.

Suppose $\sigma \in C_A-C_A^{\#}$ has a well-defined orbit.
Let $\kappa$ be the equivalent-to-trivial sequence corresponding
to $\sigma$.  By Lemma \ref{freturn} and the analogue
for the backwards orbit,
both directions of the orbit of $\sigma$ return infinitely
often to $C_A-C_A^{\#}$.
If $\kappa$ is eventually $0$, then by the
Odometer Principle $\kappa$
is in the same sequence orbit as the $0$ sequence $\kappa_0$.
But the point in $C_A$ corresponding to $\kappa_0$
is exactly the vertex $(0,-1)$.  This vertex does
not have a well defined orbit.  This is a contradiction.
If $\kappa$ is such that $k_i=d_i$ for large $i$, then
by the Odometer Principle, $\kappa$ is in the same
orbit as the sequence $\{d_i\}$.  By Equation \ref{DIO1},
the corresponding point in $C_A$ is $(2,-1)$.
One checks easily that the orbit of $(2,-1)$ is not defined
after the second iterate.  Again we have a contradiction.

Now suppose the superior sequence is eventually monotone.
We will treat the case when $A-A_n$
is eventually positive. In this case,
$\{A_n\}$ is eventually monotone increasing.  Suppose
that $\kappa$ is equivalent to the $0$-sequence.  We can
iterate backwards a finite number of times until
$\sigma$ returns as the first point of $C_A$.
Hence, without loss of generality, we can assume
that $\kappa$ is the first sequence in $\Pi_A$.
But now we can iterate forwards
indefinitely, and we will reach every equivalent-to-zero
sequence by the Odometer Principle.
Eventually we reach the $0$ sequence and get the
same contradiction as above.
If $\kappa$ is such that $k_i=d_i$ for large $i$,
we run the same argument abckwards.
\endproof

\begin{lemma}
\label{trim1}
No point of
$C_A^{\#}$ has first coordinate in $2\Z[A]$.
\end{lemma}

\startproof
Let $\{A_n\}$ be the superior sequence approximating $A$. 
We assume that $A_n<A$ infinitely often.  The other case
has the same treatment.
Suppose that 
\begin{equation}
\alpha=(2MA+2N,-1) \in C^{\#}.
\end{equation}
By Equation \ref{DIO2}, the set
$C^{\#}_A$ is invariant under the map
$(x,-1) \to (2-x,-1).$  Indeed, the twist automorphism
of $\Pi$ induces this map on $C_A$.
From this symmetry, we can assume that $M>0$.

Let $P\Gamma_k$ denote the pivot arc.
We claim that $(M,N)$ is not a vertex of $P\Gamma_k$ for any $k$.
Here is the proof.
Suppose that $(M,N) \subset P\Gamma_k$ for some $k$.
Then $2AM+2N$ is a finite
sum of terms $\lambda_j=|2Aq_j-p_j|$, by Theorem \ref{discrete}.
But such points
all lie in $C_A-C_A^{\#}$.
To avoid a contradiction, $(M,N) \not \in P\Gamma_k$ for any $k$.
This completes the proof of the claim.

Let $P\Gamma_k^+$ denote the forwards portion
of $P\Gamma_k$.  From the definition of the pivot points,
the length of $P\Gamma_k^+$ tends to $\infty$ with $k$.
Hence $\{\Gamma_k^+\}$ and $\{P\Gamma_k^+\}$ have
the same Hausdorff limit.
We can choose $k$ so large enough so
that $P\Gamma_k^+$ contains a low vertex $(M',N')$ to the
right of $(M,N)$.  So, $P\Gamma_k^+$ connects $(0,0)$ to
$(M',N')$ and skips right over $(M,N)$.

Since $\alpha \in C_A^{\#}$, we can find a sequence of
points $\{\alpha_n\} \in C_A^{\#}-\Z[A]$ such that the
first coordinate of $\alpha_n-\alpha$ is positive.
Let $\zeta_n=\alpha_n-\alpha$.   Note that
$\zeta_n \not \in 2\Z[A]$.
Let $\widehat \Gamma(\zeta_n,A)$ be the whole arithmetic
graph corresponding to $\zeta_n$.  Let
$\gamma_n=\Gamma(\zeta_n,A)$ be the
component containing $(0,0)$.
By the Rigidity Lemma, the sequences
$\{\Gamma(\zeta_n,A_n)\}$ and
$\{\Gamma_n\}$ have the same Hausdorff limit.
Hence $P\Gamma^+_k \subset \gamma_n$ once $n$ is large.
In particular, some arc of $\gamma_n$ connects
$(0,0)$ to $(M',N')$ and skips over $(M,N)$.
Call this the {\it barrier arc\/}. 

Since $\alpha_n-\zeta_n=\alpha \in 2\Z[A]$, there is another
component $\beta_n \subset \widehat \Gamma(\zeta_n)$ that
tracks the orbit of $\alpha_n$.   One of the vertices of
$\beta_n$ is exactly $(M,N)$.
The component $\beta_n$ is unbounded in both
directions, because all defined orbits in $C_A^{\#}$
are unbounded.  On the other hand
$\beta_n$ is trapped beneath the barrier arc.
It cannot escape out either end, and it cannot
intersect the barrier arc, by the Embedding Theorem.
But then $\beta_n$ cannot be unbounded in either
direction.  This is a contradiction.   
\endproof

Corollary \ref{trim2} and Lemma \ref{trim0} show
that $U_A \cap I \subset C_A^{\#}$.   Lemma
\ref{trim1} shows that
$C_A^{\#} \subset C_A'$.  Lemma \ref{precomet}
shows that $C_A' \subset U_A \cap I$. Putting
all this together gives Equation \ref{doubleid}.

\subsection{Statement 4 of the Comet Theorem}
\label{2distinct}

We have already established the first part of
Statement 4.  Now we prove the second part.

By Statements 1 and 2 of the Comet Theorem,
it suffices to consider
pairs of points in $C_A^{\#}$.  (This is
where we use the truth of Statement 1.)
 It follows  immediately
from Equation \ref{trans} that two points
of $C_A^{\#}$ lie on the same orbit only if
their first coordinates differ by an element of
$2\Z[A]$.  Our goal is to prove the converse.

\begin{lemma}
All but at most $2$ orbits in $C_A^{\#}$ are
erratic.
\end{lemma}

\startproof
By Lemma \ref{precomet} and Lemma
\ref{freturn}, and the backwards analogue
of Lemma \ref{freturn}, all orbits in
$C_A^{\#}$ are erratic except for those
corresponding to the eqivalent-to-first
sequences and the equivalent-to-last
sequences.    By the Odometer Principle,
all the points in $C_A^{\#}$ corresponding
to equivalent-to-first sequences lie on
the same orbit.  Likewise, all the
points in $C_A^{\#}$ corresponding
to equivalent-to-last sequences lie in
the same orbit.
These two orbits are the only ones which
can fail to be erratic.
\endproof

\begin{lemma} Suppose that two points in 
$C_A^{\#}$ have first coordinates that differ
by $2\Z[A]$.  Suppose also that at least one
of the points has an erratic orbit.  Then
the two points lie on the same orbit.
\end{lemma}

\startproof
One direction follows immediately from Equation \ref{trans}.
For the converse, suppose that our two points have first
coordinates that differ by $2\Z[A]$.  The first coordinates
of our points do not lie in $2\Z[A]$, by Lemma
\ref{trim1}.  Hence, one
and the same arithmetic graph $\widehat \Gamma$ 
contains components $\gamma_1$ and $\gamma_2$
that respectively track our two orbits. 

Since both orbits are dense in $C_A^{\#}$, we know that
both orbits are erratic in at least one direction.
Suppose first that $\gamma_1$ is erratic in both
directions.  Since $\gamma_2$ is erratic in one
direction, we can find a low vertex $v$ of $\gamma_1$
that is not a vertex of $\gamma_2$.  Since
$\gamma_2$ is erratic in both directions, we can
find vertices $w_1$ and $w_2$ of $\gamma_1$,
lying to the left and to the right of $v$.  But
then the arc of $\gamma_1$ starting at $v$
is trapped beneath the arc of $\gamma_2$
connecting $w_1$ to $w_2$.  This contradicts
the Embedding Theorem.  In short, $\widehat \Gamma$
is not big enough to contain both components.
\endproof

It only remains to deal with the case when both
points lie on orbits that are only erratic
in one direction..

\begin{lemma} Suppose that two points in 
$C_A^{\#}$ have first coordinates that differ
by $2\Z[A]$.  Suppose also that neither point
lies on an erratic orbit. Then
the two points lie on the same orbit.
\end{lemma}

\startproof
Let $\alpha \in C_A^{\#}$ (respectively $\beta$) be the unique point such that
the forwards (respectively backwards) first return map to $C_A^{\#}$ at $\alpha$
(respectively $\beta$) does not
exist. There are exactly $2$ one-sided erratic orbits.
$\alpha$ is one orbit and $\beta$ is on the other. 
It suffices to prove that $\alpha-\beta \not \in 2\Z[A] \times \{0\}$.
We will suppose the contrary, and derive a contradiction.
Suppose that $\alpha-\beta = (2Am+2n,0)$ for some $(m,n)\in \Z^2$.

$\alpha$ is the last point in the twirl order and $\beta$ is
the first point.  In terms of sequences, $\alpha$ corresponds
to the sequence $\{\widetilde d_i\}$ and $\beta$ corresponds
to the sequence $\{\widetilde 0_i\}$.   Let $\{\alpha_j\}$
be a sequence of points in $C_A^{\#}$ converging to
$\alpha$, chosen so that the corresponding orbit is erratic.
Define $\beta_j=\alpha_j+(\beta-\alpha)$.  Then
$\alpha_j-\beta_j=(2Am+2n,0)$.  By the case we have
already considered, $\beta_j$ lies in the same orbit
as $\alpha_j$.   

For $j$ large, the sequence corresponding to $\alpha_j$ matches
the terms of the sequence for $\alpha$ for many terms.
Likewise, the sequence corresponding to $\beta_j$ matches
the terms of the sequence for $\beta$ for many terms. 
Hence, these two sequences disagree for many terms.
Given that our return dynamics to $C_A^{\#}$ is conjugate
to the odometer map on the sequence space, we have
\begin{equation}
2Am+2n=\pi_1(\alpha_j-\beta_j)=\sum_{i=0}^{N_j} a_{ji} \lambda_i; \hskip 30 pt
|a_{ji}| \leq d_i.
\end{equation}
Here $N_j \to \infty$ as $j \to \infty$, and $\pi_1$ denotes
projection onto the first coordinate. 

Let $M$ be the map from Equation \ref{funm}.  We have
\begin{equation}
M(m,n)=\sum_{i=0}^{N_j} b_{ji}M(V_i); \hskip 30 pt |b_{ji}| \leq d_i.
\end{equation}
Here $b_{ji} = \pm a_{ji}$, depending on the sign of $A_i-A$.  Since
$A$ is irrational, $M$ is injective.  Therefore, setting $N=N_j$
for ease of notation, we have
\begin{equation}
(m,n)=\sum_{i=0}^{N} b_{ji} V_{i}=b_{Ni}V_N + \sum_{i=0}^{N-1} b_{ji}V_i.
\end{equation}
Looking at the second coordinates, we see that
\begin{equation}
q_N-\sum_{i=0}^{N-1}d_i q_i \leq \bigg|b_{Ni} q_N -\sum_{i=0}^{N-1} b_{ji}q_i\bigg|=|n|.
\end{equation}
However, it follows fairly easily from Equation \ref{DIO4} that the left hand
side tends to $\infty$ as $N_j \to \infty$.  This contradiction
finishes the proof.
\endproof

\newpage

\newpage

\section{Dynamical Consequences}
\label{comet2}

In this chapter we discuss some dynamical 
consequences of the Comet Theorem.

\subsection{Minimality and Homogeneity}
\label{comets}

Now we deduce some consequences of the Comet Theorem.
Let $U_A$ denote the set of unbounded special orbits.
Since every orbit in $U_A$ intersects
$C_A^{\#}$, it suffices to prove that
every point of $C_A^{\#}$ lies on an orbit
that is either forwards dense in $U_A$ or
backwards dense or both.

Let $\zeta \in C_A^{\#}$ be our point.  By the Comet Theorem, the
orbit of $\zeta$ is either forwards dense in $C_A^{\#}$, or backwards
dense in $C_A^{\#}$, or both.  Assume that $\zeta$ lies on an orbit
that is forwards dense in $C_A^{\#}$.  The case of backwards dense
orbits has a similar treatment.

Let $\beta \in U_A$ be some other point. Some point 
$\alpha \in C_A^{\#}$ lies in the orbit of $\beta$.
Hence, $(\psi')^k(\alpha)=\beta$ for some $k$.
Here $\psi'$ is the outer billiards map.
But $(\psi')^k$ is a piecewise
isometry.  Hence, $(\psi')^k$ maps small intervals
centered at $\alpha$ isometrically to small intervals
centered at $\beta$.
The forwards orbit
of $\zeta$ enters any interval about $\alpha$
infinitely often.  Hence, the
forwards orbit of $\zeta$ enters every interval about
$\beta$ infinitely often.  

Say that a subset $S \subset \R^2$ is
{\it locally homogeneous\/} if every two points
of $S$ have arbitrarily small neighborhoods
that are translation equivalent.  Note that the
points themselves need not sit in the same positions
within these sets. 

\begin{lemma}
For any irrational $A$, the set $U_A$ is locally
homogeneous.
\end{lemma}

\startproof
For any $p \in U_A$, there is some integer $k$
such that $(\psi')^k(p) \in C_A^{\#}$.
Here $\psi'$ is the outer billiards map.  But
$\psi^k$ is a local isometry.   Hence, there
are arbitrarily small neighborhoods of $p$
that are isometric to neighborhoods
of points in $C_A^{\#}$.  

Hence, it suffices to prove that $C_A^{\#}$ is
locally homogeneous. Let $\{d_k\}$ denote
the renormalization sequence.  The set $C_A$ breaks
into $d_0+1$ isometric copies of a smaller Cantor
set Each of these breaks into $d_1+1$ isometric
copies of still smaller Cantor sets.  And so on.
From this we see that both $C_A$ and 
$C_A^{\#}$ are locally homogeneous.
\endproof

\subsection{Tree Interpretation of the Dynamics}
\label{treedyn}

Let $A$ be an irrational kite parameter. We can 
illustrate the return dynamics to $C_A^{\#}$
using infinite trees.   The main point here
is that the dynamics is conjugate to an
odometer.   The conjugacy is given by the
map $\phi: {\cal Z\/}_A \to C_A$ from
the Comet Theorem.   Our pictures
encode the structure of $\phi$ graphically.

We think of $C_A$ as the ends of a tree $T_A$.
We label $T_A$ according to the sequence of signs 
$\{A-A_n\}$.  Since $A-A_0$ is negative, we label the
level $1$ vertices $0,...,d_0$ from right to left.
Each level $1$ vertex has $d_1$ downward vertices.
We label all these vertices from left to right
if $A-A_1>0$ and from right to left if $A-A_1$ is
negative.  And so on.  This business of switching
left and right according to the sign of $A-A_k$
corresponds precisely to our method of
identification in Equations \ref{phi1} and
\ref{tildek}.
Figure 25.1 shows the example for the renormalization
sequence $\{1,3,2\}$ and the sign sequence $-,+,-$.

\begin{center}
\resizebox{!}{1.9in}{\includegraphics{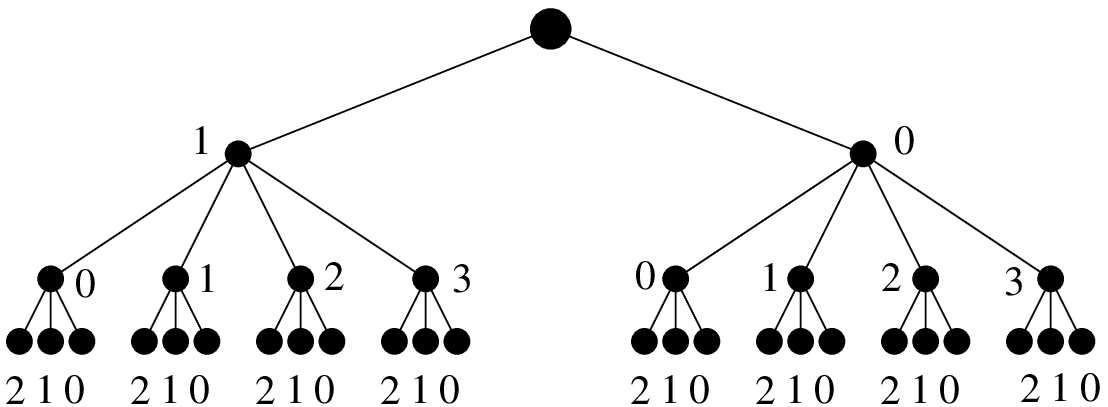}}
\newline
{\bf Figure 25.1:\/} Tree Labelling
\end{center}

We have the return map 
$$\rho_A: C_A^{\#}-\phi(-1) \to C_A^{\#}-\phi(-1),$$
and this map is conjugate to the restriction of the
odometer on ${\cal Z\/}_A$.  Accordingly, we can
extend $\rho_A$ to all of $C_A$, even though the
extension no longer describes outer billiards dynamics
on the extra points.   Nonetheless, it is convenient
to have this extension.

To see what $\rho_A$ does, we write the code for
a given end. Then we add $1$, carrying to the right.
Referring to our example above, we have
$000... \to 100...$ and $130... \to 001...$.  This
map is exactly what is called an odometer.

\subsection{Periodic Orbits}

One might wonder about the other orbits in the interval $I$.
First of all, we have the following result.

\begin{theorem}
Any defined orbit in $I-C_A$ is periodic.  There is a uniform bound
on the period, depending only on the distance from the point to 
$C_A$.
\end{theorem}

\startproof
The Comet Theorem combines with the Dichotomy Theorem to
prove any defined orbit in $I-C_A$ is periodic.  
The period bound comes from taking a limit of the
Period Theorem as $n \to \infty$ in our rational
approximating sequence.  In other words, if this
result was false, then we could contradict the
Period Theorem using our Continuity Principle.
\endproof

In the next chapter we will prove that
$C_A$ has length $0$. See Lemma \ref{length0}.
  By the local
homogeneity, $U_A$ also has length $0$.
Hence, by the Dichotomy Theorem, almost
all special orbits are periodic.
\newline
\newline
{\bf A Conjectural Picture:\/}
It we knew Conjecture \ref{gapmap2},
we could give a very nice account of what happens.
We now describe this conjectural picture.

We can naturally identify $C_A$ with the ends of an
infinite directed tree $T_A$.  Using the homeomorphism
$\phi: {\cal Z\/}_A \to C_A$, we can formally
extend the return map on $C_A^{\#}-\phi(-1)$ to
all of $C_A$, even though the extended return map
does not correspond to the outer billiards
dynamics on the extra points.  This is
exactly what we did in \S \ref{treedyn} above.

The extended return map to
$C_A$ induced  by an automorphism 
\begin{equation}
\Theta_A: T_A \to T_A
\end{equation}
as discussed in \S \ref{treedyn}.  The complementary
open intervals in $I-C_A$ -- the {\it gaps\/}-- are naturally in bijection
with the forward cones of $T_A$.

\begin{conjecture}
\label{gapmap3}
The outer billiards map is entirely defined on a gap.
The return map to $I-C_A$ permutes the gaps
according to the action of $\Theta_A$ on the
forward cones of $T_A$.
\end{conjecture}

\subsection{Proper Return Models and Cusped Solenoids}
\label{model0}

Here we will describe the sense in which the Comet Theorem
allows us to combinatorially model the dynamics on
$U_A$.  The results in this section are
really just a repackaging of some of the statements
of the Comet Theorem.

Let $X$ be an unbounded metric space and let
$f: X \to X$ be a bijection.  We assume that
$f^2$ moves points by a small amount.  That is,
there is a universal constant $C$ such that
\begin{equation}
d(x,f^2(x))<C \hskip 30 pt \forall x \in X.
\end{equation}
The example we have in mind, of course, is the
outer billiards map 
\begin{equation}
\label{basicmap}
\psi': U_A \to U_A.
\end{equation}
The square map $\psi$ moves points by at
most $4$ units.

We say that a compact
subset $X_0 \subset X$ is a {\it proper section\/}
for $f$ if for every $N$ there is some $N'$
such that $d(x,X_0)<N$ implies that
$f^k(x) \in X_0$ for some $|k|<N'$. 
In particular, every orbit of $f$ intersects $X_0$.
This condition is just the abstract version of
Statement 1 of the Comet Theorem. Informally,
all the orbits either head directly to $X_0$ or directly
away from $X_0$.  

Let $f_0: X_0 \to X_0$ be the first return map.
This is a slight abuse of notation, because
$f_0$ might not be defined on all points of
$X_0$.  Some points might exit $X_0$ and
never return.  
We define
two functions $e_1,e_2: X_0 \to (0,\infty]$.
The function $e_1(x)$ is the maximum
distance the forward orbit of $x$
gets away from $X_0$ before returning
as $f_0(x)$. The function $e_2(x)$ is
the length of this same portion of
the orbit.  If $f_0$ is not defined
on $x$ then obviously $e_2(x)=\infty$.
The proper section condition guarantees
that $e_1(x)=\infty$ as well.

The condition that $X_0$ is a proper
section guarantees that $e_1$ and $e_2$ are
proper functions of each other.
That is, if $\{x_n\}$
is a sequence of points in $X_0$, then
$e_1(x) \to \infty$ if and only if
$e_2(x) \to \infty$.  This observation
includes the statement that
$e_1(x)=\infty$ iff $e_2(x)=\infty$ iff
$f_0$ is not defined on $x_0$.
For the purposes of getting a rough
qualitative picture of the orbits,
we just consider the function $e_1$.
We set $e=e_1$, and call $e$ the
{\it excursion function\/}.

Suppose now that $f':X' \to X'$ is another
bijection, and $X_0'$ is a proper section.
Let $e': X_0' \to (0,\infty]$ denote the
excursion function for this system.  We
say that $(X,X_0,f)$ is {\it properly
equivalent\/} to $(X',X_0',f')$ if there
is a homeomorphism
$\phi: X \to X'$ such that
\begin{itemize}
\item $\phi$ conjugates $f_0$ to $f_0'$. 
\item $e' \circ \phi$ and $e$ are proper functions
of each other on $X_0$. 
\end{itemize}
These conditions guarantee that $\phi$ carries the
points where $f_0$ is not defined to the points
where $f_0'$ is not defined.

The notion of proper equivalence turns out to be 
a tiny bit too strong for
our purposes.  We saye that $(X,X_0,f)$ and $(X',X_0',f')$
are {\it essentially properly equivalent\/} if
$\phi$ has all the above properties but is only
defined on the complement of a finite number of
orbits of $X_0$.   In this case, the inverse map
will have the same property:  It will be well
defined on all but a finite number of orbits
of $X_0'$.   In other words, an essential proper
equivalence is a proper equivalence provided that
we first delete a finite number of orbits from
our spaces.   We call $(X,X_0,f)$ an
{\it essentially proper model\/} for
$(X',X_0',f')$.

Statement 1 of the Comet Theorem says that
$C_A^{\#}$ is a proper section for the map in
Equation \ref{basicmap}.
Now we can describe our proper models for
the triple $(U_A,C_A^{\#},\psi')$.  Statements
2 and 3  in particular describe the excursion
function up to a b-lipschitz constant.
Here we convert this information into a
concrete essentially proper model for
this dymamics.

Let ${\cal Z\/}_A$ denote the metric abelian
group from the Comet Theorem.
For convenience, we recall the definition
of the metric $d$ here.
$d(x,y)=q_{n-1}^{-1}$,
where $n$ is the smallest index
such that $[x]$ and $[y]$ disagree in $\Z/D_n$.
Here $\{p_n/q_n\}$ is the superior sequence
approximating $A$.

We denote the odometer map on ${\cal Z\/}_A$ by
$f_0$.   That is, $f_0(x)=x+1$.   Topologically,
the {\it solenoid\/} based on ${\cal Z\/}_A$ is
defined as the mapping cylinder
\begin{equation}
{\cal S\/}_A={\cal Z\/}_A \times [0,1]/\sim; \hskip 30 pt
(x,1) \sim (x+1,0).
\end{equation}
This is a compact metric space.

We now modify this space a bit.  First of all, we
remove the point 
$$(-1,1/2)$$
from ${\cal S\/}_A$.  This deleted point, our cusp,
lies halfway between $(-1,0)$ and $(0,0)$.  We now
change the metric on our space by declaring the
length of the segment between $(x,0)$ and $(x,1)$ to
be $$
\frac{1}{d(x,-1)}$$
Metrically, we simply rescale the
length element on each interval by the appropriate
amounts.   We call the resulting space
${\cal C\/}_A$.  We call ${\cal C\/}_A$ the
{\it cusped solenoid\/} based on $A$.

We define $f: {\cal C\/}_A \to {\cal C\/}_A$ to be
the map such that
\begin{equation}
f(x,t)=\bigg(x,\frac{t}{d(x,-1)}\bigg)
\end{equation}
From the way we have scaled the distances,
$f$ maps each point by $1$ unit. Indeed,
some readers will recognize $f$ as the
time-one map of the geodesic slow on
${\cal C\/}_A$.  The original set
${\cal Z\/}_A$ is a proper section
for the map, and the return map is
precisely $f_0$.   Put another way,
$f$ is a suspension flow over $f$.
Note that $f$ also depends on $A$, but
we suppress this from our notation.

\begin{theorem} 
\label{model1}
The triple $({\cal C\/}_A,{\cal Z\/}_A,f)$ is an
essentially proper model for $(U_A,C_A^{\#},\psi')$.
\end{theorem}

\startproof
This is just a repackaging (and weakening) 
of Statements 2 and 3 of
the Comet Theorem.
\endproof

\noindent
{\bf Remarks:\/} \newline
(i)  Our model forgets the
linear ordering on $C_A^{\#}$ that
comes from its inclusion in $I$, but one
can recover this
from the discussion in \S \ref{treedyn}.
\newline
(ii)
In a certain
sense, the triple $({\cal C\/}_A,{\cal Z\/}_A,f)$
provides a {\it bi-lipschitz model\/} for
the nature of the unboundedness of the orbits
in $U_A$.
However, it would be misleading
to call our model an actual bi-lipschitz model
for the dynamics on $U_A$ 
because we are not saying much about what happens
to the orbits in the two systems after they
leave their proper section. For instance,
the excursion times could be wildly different
from each other, even though they are
proper functions of each other.
\newline

Here is a universality result.
\begin{theorem}
\label{universal}
The time-one map of the
geodesic flow on any cusped solenoid serves as
an essentially proper model for the dynamics of
the special unbounded orbits relative to uncountably
many different parameters.
\end{theorem}

\startproof
Up to a proper change of the excursion function, our model
only depends on the renormalization sequence, and there are
uncountbly many parameters realizing any renormalization
sequence.
\endproof

\subsection{Equivalence and Universal Behavior}
\label{UNIVERSAL}

To each parameter $A$, we associate the
renormalization sequence $\{d_n\}$.   We then
associate the sequence $\{D_n\}$, where
\begin{equation}
D_n=\prod_{i=0}^{n-1}(d_i+1).
\end{equation}

We call $A$ and $A'$ {\it broadly 
equivalent\/} iff
for each $m$ there is some $n$ such that $D_m$ divides
$D_n'$ and $D_m'$ divides $D_n$.  Each broad
equivalence class has uncountably many members.

\begin{lemma}
  If $A$ and $A'$
are broadly equivalent then there is a homeomorphism
from ${\cal Z\/}_A$ to ${\cal Z\/}_{A'}$ that
conjugates the one odometer to the other.
\end{lemma}

\startproof 
Each element of
${\cal Z\/}_A$ is a compatible sequence
$\{a_m\}$ with $a_m \in \Z/D_m$.  Using the
divisbility relation, this element determines
a corresponding sequence $\{a_m'\}$.  Here
$a_m'$ is the image of $a_n$ under the factor
map $\Z/D_n \to \Z/D_m'$, where $n$ is
such that $D_m'$ divides $D_n$.  One checks easily
that this map is well defined and determines
the desired homeomorphism.
\endproof

\begin{theorem}
If $A$ and $B$ are broadly equivalent then
there is an essentially proper
equivalence between $(U_A,C_A^{\#},\psi'_A)$ and
$(U_{B},C_{B}^{\#},\psi'_{B})$.
In particular, the return maps to
$C_A$ and $C_B$ are topologically
conjugate modulo countable sets.
\end{theorem}

\startproof
The homeomorphism from
${\cal Z\/}_A$ to ${\cal Z\/}_B$ maps
$-1$ to $-1$.  By construction, this
homeomorphism sets up a proper equivalence
between  $({\cal C\/}_A,{\cal Z\/}_A,f_A)$ and
$({\cal C\/}_B,{\cal Z\/}_B,f_B)$.
This result now follows from Theorem
\ref{model1}.
\endproof

One might wonder about the nature of the topological
equivalence between the return maps to
$C_A^{\#}$ and $C_B^{\#}$.  One can reconstruct
the conjugacy from the tree labellings given
in \S \ref{treedyn}.  The conjugacy is well defined
for all points of $C_A$ and $C_B$, but we typically
have to ignore the countable sets of points on which
the relevant return maps are not defined.  This
acconts for the precise statement of our theorem
above.

Let ${\cal Z\/}$ denote the inverse limit over all
finite cyclic groups.  The map $x \to x+1$ is
defined on ${\cal Z\/}$. This dynamical
system is called the {\it universal odometer\/}.
Sometimes ${\cal Z\/}$ is called the
{\it profinite completion\/} of $\Z$.

We call $A$ {\it universal\/} if every $k \in \N$ divides
some $D_n$ in the sequence.  If $A$ is universal,
then there is a group isomorphism from
${\cal Z\/}$ to ${\cal Z\/}_A$ that respects
the odometer maps.  In short, when $A$ is
universal, ${\cal Z\/}_A$ is the universal odometer.
See [{\bf H\/}, \S 5] for a proof of this fact -- stated
in slightly different terms -- and
for a detailed discussion of the
universal odometer.

\begin{lemma}
Almost every parameter is universal.
\end{lemma}

\startproof
A sufficient condition for a parameter to be universal
is that every integer appears in the renormalization
sequence. We can express the fact that
a certain number appears in the renormalization sequence
as a statement that a certain combination appears in
the continued fraction expansion of $A$.  Geometrically,
as one drops a geodesic down from $\infty$ to $A$, the
appearance of a certain pattern of geodesics in the
Farey graph forces a certain number in the renormalization
sequence.  As is well known, the continued fraction
expansion for almost every number in $(0,1)$ contains
every finite string of digits.  
\endproof

\begin{theorem}
For almost every $A \in (0,1)$, the triple
$(U_A,C_A^{\#},\psi')$ is properly modelled
by the time-one map of the geodesic flow
on the universal cusped solenoid.  In particular,
the return map to $C_A^{\#}$ is topologically
conjugate to the universal odometer, modulo a
countable set. 
\end{theorem}

\startproof
This is an immediate consequence of the
previous result and Theorem  \ref{model1}.
\endproof

One might wonder if there is a concrete parameter that
exhibits this universal behavior.  Here we give an
example.  Let $A$ be the parameter whose inferior sequence satisfies
$$
\frac{1}{1} \leftarrow \frac{5}{7} \leftarrow \frac{51}{71} \leftarrow \frac{719}{1001} \ldots
;\hskip 30 pt
r_{n+1}=(4n+2)r_n+r_{n-1}
$$
$r$ stands for either $p$ or $q$.  All terms are superior.  The
renormalization sequence is $3,5,7,9...$.   Hence
$$D_0=4; \hskip 15 pt D_1=4 \times 6; \hskip 15 pt
D_2=4 \times 6 \times 8 \ldots$$
One can see easily in this example
that $\Sigma(A)=\N$.  Hence $A$ is universal.
 Let $E=A+2$. 
It seems that $E=e$, the base of the natural log.  We didn't work out
a proof, but this should follow from the 
famous continued fraction expansion $e=[2;1,4,1,1,6,1,1,8,...]$.
In short, $e-2$ is universal.

\subsection{Some other Equivalence Relations}

Call $A$ and $B$ {\it narrowly equivalent\/} if they have
the same renormalization sequence and if
the sign of $A-A_j$ is the same as the sign of
$B-B_j$ for all $j$.  Here $\{A_j\}$ and
$\{B_j\}$ are the superior sequences approximating
$A$ and $B$ respectively.  Referring
to Equation \ref{tildek}, 
the definition of $\widetilde k_j$ relative
to the narrowly equivalent parameters is the same for every index.
Each narrow equivalence class again has uncountably many
members.   

\begin{theorem}
If $A$ and $B$ are narrowly equivalent then there
is an order-preserving homeomorphism from
$I$ to $I$ that conjugates the return map
on $C_A^{\#}$ to the return map on
$C_{B}^{\#}$.  This map is a proper equivalence
from $(U_A,C_A^{\#},\psi'_A)$ to
from $(U_{B},C_{B}^{\#},\psi'_{B})$ to
\end{theorem}

\startproof
The two spaces $\Pi_A$ and $\Pi_{B}$ are exactly the same,
and the extended twirl orders on the (equivalence classes)
of these spaces are the same.  Thus, the successor maps
on the two spaces are identical.  The map
$h=\phi_2' \circ \phi_2^{-1}$ is a homeomorphism from
$C_A$ to $C_{B}$ that carries $C_A^{\#}$ and
$C_{B}^{\#}$ and conjugates the one return 
dynamics to the other.   By construction,
$h$ preserves the linear ordering on $I$, and we
can extend $h$ to the gaps of $I-C_A$ in the
obvious way.  By construction, this map
carries $\phi_A(-1)$ to $\phi_B(-1)$ and is
continuous.  Hence, it is a proper equivalence
in the sense discussed above.
\endproof

The {\it first renormalization\/} of the odometer map $x \to x+1$ on
the inverse system 
\begin{equation}
\ldots \to \Z/D_3 \to \Z/D_2 \to \Z/D_1
\end{equation}
is the $D_1$st power of the map.  This corresponds to the map
$x \to x+1$ on the inverse system
\begin{equation}
\ldots \to \Z/D_3' \to \Z/D_2' \to \Z/D_1'; \hskip 30 pt
D_n'=D_{n+1}/D_1.
\end{equation}
As in the Comet Theorem, each $D_n$ divides $D_{n+1}$ for all $n$,
so the construction makes sense.  In terms of the symbolic dynamics
on the sequence space $\Pi$, the renormalization consists of the
first return map to the subspace
\begin{equation}
\Pi'=\{\kappa \in \Pi|\ k_0=0.\}
\end{equation}
In terms of the dynamics on $C_A$, the first renormalization
is the first return map to the Cantor subset corresponding
to $\Pi'$.   The {\it second renormalization\/} if the
first renormalization of the first renormalization. And so on.

Let $\Gamma_2 \subset SL_2(\Z)$ denote the subgroup of matrices
congruent to the identity mod $2$.  Then $\Gamma_2$ acts on
$\Q \cup \infty$ by linear fractional transformations.
The action preserves the parity of the rationals.   Even
though $\Gamma_2$ does not preserve the parameter interval
$(0,1)$, it still makes sense to say that $A \sim B$ mod
$\Gamma_2$.  This is to say that
\begin{equation}
\label{gamma2}
B=\frac{a A+b}{c A+d}; \hskip 30 pt \left[\matrix{a&b \cr c&d}\right] \in \Gamma_2.
\end{equation}

Here we recall our construction of the inferior
sequence for $A$.
Our construction is based on the graph in the hyperbolic
plane obtained from the Farey graph by deleting
the edges connecting even rationals to each other.
The result is the $1$-skeleton of a tiling by ideal squares.
$\Gamma_2$ preserves this tiling.  We construct the
inferior sequence by dropping a vertical geodesic down
to $A$ and recording the sequence of ideal squares the
geodesic enters as it limits to $A$.  From this
description, we see that the renormalization and
sign sequences
for $A$ and $B$ are eventually the same.   This
gives us the following result.

\begin{corollary}
Suppose that $A$ and $B$ are equivalent under $\Gamma_2$.
Then the return maps to $C_A^{\#}$ and $C_{B}^{\#}$ have a common
renormalization.  The conjugacy between the one
renormalization to the other is
implemented by a homeomorphism that
preserves the order on the interval $I$.
\end{corollary}

\newpage

\section{Geometric Consequences}
\label{lengthdim}

\subsection{Hausdorff Dimension}
\label{hausreview}

In this chapter we study the structure of $C_A$.

We first review
some basic properties of the
the Hausdorff dimension, including its definition.
\newline
\newline
{\bf Basic Definition:\/}
Given an interval $J$, let $|J|$ denote its length.
Let $I=[0,2] \times \{-1\}$ be our usual interval.  
Given a subset $S \subset I$, and $s \in [0,1]$, and
some $\delta>0$, we define
\begin{equation}
\label{cover}
\mu(S,s,\delta)=\inf \sum |J_n|^s
\end{equation}
The infimum is taken over all countable covers of
$S$ by intervals $\{J_n\}$ such that
${\rm diam\/}(J_n)<\delta$.
Next, we define
\begin{equation}
\mu(S,s)=\lim_{\delta \to 0} \mu(S,s,\delta) \in [0,\infty].
\end{equation}
This limit exists because $\mu(S,s,\delta)$ is a
monotone function of $\delta$.
Note that $\mu(S,1)<\infty$ because $I$ has finite total length.
Finally, 
\begin{equation}
\dim(S)=\inf \{s|\ \mu(S,s)<\infty\}.
\end{equation}
The number $\dim(S)$ is called the {\it Hausdorff dimension\/} of $S$.
\newline
\newline
{\bf Bi-Lipschitz Invariance:\/}
Let $f: \R \to \R$ be a map.  $f$ is called $K$-{\it bi-lipschitz\/} if
\begin{equation}
K^{-1}\|x-y|<\|\phi(x)-\phi(y)\|<K\|x-y\|
\end{equation}
$f$ is called {\it bi-lipschitz\/} if it is $K$-bi-lipschitz
for some $K$.   It follows easily from the definitions
that $\dim(S)=\dim(S')$ if $f(S)=S'$ for some bi-lipschitz
function $f$.
\newline
\newline
{\bf Borel Slicing Property:\/}
Let $S \subset [0,1]^2$ be a Borel subset.  Let
$S_A$ denote the intersection of $S$ with the line
$\{y=A\}$.   Let $f(A)=\dim(S_A)$.  It is
known that $f$ is a Borel measurable function.
See [{\bf MM\/}].  In our application, we shall
apply this criterion to the set $C$ from
Equation \ref{butterfly}.  This is a very
explicit example of a Borel measurable set.

\subsection{Ubiquity of Periodic Orbits}

\begin{lemma}
\label{length0}
$C_A$ has length $0$.
\end{lemma}

\startproof
Let $\lambda_n=|Aq_n-p_n|$, as in Equation \ref{DIO1}.
We define
\begin{equation}
G_n=\sum_{k={n+1}}^{\infty}2\lambda_k d_k.
\end{equation}
Then
\begin{equation}
C_A \subset \sum_{\kappa \in \Pi_n} \Big(I_n+X(\kappa)\Big).
\end{equation}
Here $I_n$ is the interval with endpoints $(0,1)$
and $(G_n,1)$.
In other words, $C_A$ is contained in  $D_n$ translates
of an interval of length $G_n$.  We just need to
prove that $D_n G_n \to 0$.   It suffices to
prove this when $n$ is even.  By Equation \ref{DIO4}, 
\begin{equation}
\label{measure0}
D_n<\epsilon^{-n} q_n; \hskip 30 pt
\epsilon=\sqrt{5/4}.
\end{equation}
By Equation \ref{DIO1} we have
\begin{equation}
\label{measure1}
G_n<2\sum_{k=n+1}^{\infty} q_k^{-1}<2q_{n}^{-1} \sum_{k=1}^{\infty}2^{-k}<2q_n^{-1}.
\end{equation}
Here we have used the trivial bound that
$q_{m}/q_n<2^{n-m}$ when $m>n$.
Therefore
\begin{equation}
\label{measure2}
D_nG_n<2\epsilon^{-n}.
\end{equation}
This completes the proof.
\endproof

\begin{theorem}
Relative to any irrational parameter, almost every point on
$\R \times \Z_{\rm odd\/}$ has a periodic outer billiards
orbit.
\end{theorem}

\startproof
Since $U_A$ is locally homogeneous and
$C_A^{\#}$ has length $0$, the set
$U_A$ has length $0$.  The point here is that
$U_A$ cannot have any points of Lebesgue density.
There are only countably many points in $\R \times \Z_{\rm odd\/}$
with undefined orbits, and the rest are periodic by
the Dichotomy Theorem.
\endproof

\subsection{A Dimension Formula}

Now we prove the following result.

\begin{theorem}
\label{dimformula}
Let $A$ be an irrational parameter.  Let
$\{p_n/q_n\}$ be the superior sequence associated to $A$.
Suppose that $q_{n+1}<Cq_n$ for some constant $C$ that is independent of $n$.
Then
$$u(A)=\frac{D}{Q}; \hskip 30 pt
D=\lim_{n \to \infty}\frac{\log(D_n)}{n}; \hskip 20 pt
Q=\lim_{n \to \infty}\frac{\log(q_n)}{n},$$
provided that these limit exist.
Limits are taken with respect to the superior terms.
\end{theorem}

We call $A$ {\it tame\/} if $A$ satisfies the hypotheses of
Theorem \ref{dimformula}.  We leave it as an exercise
to the interested reader to show that all quadratic
irrational parameters are tame.

\begin{lemma}
\label{tamebound}
Suppose $A$ is a tame parameter.
Let $\{p_n/q_n\}$ be the associated superior
sequence.  Then $\lambda_n \in [C_1,C_2]q_n^{-1}$
for positive constants $C_1, C_2$.
\end{lemma}

\startproof
For tame parameters, the renormalization
sequence $\{d_n\}$ is bounded.
We have
$$\lambda_n=q_n|A-A_n|<2d_n^{-1}q_n^{-1}<C_2q_n^{-1}$$
by Lemma \ref{superior}.   For the lower bound
note first that
$\lambda_{n+1}<\lambda_{n+1}'<\lambda_n$,
by Equation \ref{DIO3}.   By the triangle
inequality
$$|A-A_n|+|A-A_{n+1}| \geq |A_n-A_{n+1}| \geq
\frac{2}{q_nq_{n+1}}.$$
Hence
$$
2\lambda_n>\lambda_n+\lambda_{n+1}=
q_n|A-A_n|+q_{n+1}|A-A_{n+1}|>$$
$$
q_n\Big(|A-A_n|+|A-A_{n+1}|\Big) \geq 2q_{n+1}^{-1} \geq 2C_1q_n^{-1}.
$$
This gives the lower bound.
\endproof

Now we derive our dimension formula for tame
parameters $A$.  
The constants
$C_1,C_2,...$ denote positive constants that
depend on $A$.   
Let ${\cal C\/}_n$ be the covering constructed in the
proof of Lemma \ref{length0}.
The intervals in ${\cal C\/}_n$ are pairwise disjoint and
all have the same length.   
Each interval of ${\cal C\/}_n$ contains $(d_n+1)$
evenly and maximally spaced intervals of ${\cal C\/}_{n+1}$.
From these properties, it suffices to use
the covers ${\cal C\/}_n$ to compute $u(A)$. 

There are $D_n$ intervals in ${\cal C\/}_n$, all having
length $G_n$.
Choose any $\epsilon>0$. 
For $n$ large, we have
\begin{equation}
D_n \in \bigg(\exp\Big(n(D-\epsilon)\Big),\exp\Big(n(D+\epsilon)\Big)\bigg)
\end{equation}
We have
$$
G_n=2\lambda^*_{n+1} \in [2\lambda_{n+1},\lambda_n] \in [C_1q_{n+1}^{-1},C_2q_n^{-1}]
\in [C_3,C_2]q_n^{-1},$$
by the preceding lemma.
Hence
\begin{equation}
G_n \in \bigg(\exp\Big(-n(Q+\epsilon)\Big),\exp\Big(-n(Q-\epsilon)\Big)\bigg).
\end{equation}
From these estimates, we get
$u(A) \in [(D/Q)-\epsilon,(D/Q)+\epsilon]$.
But $\epsilon$ is arbitrary.
This establishes our dimension formula.

\subsection{Modularity}
\label{further}

The level $2$ congruence subgroup $\Gamma_2$ acts on
$\R \cup \infty$ by linear fractional transformations.
To say that $A$ and $A'$ are in the same
$\Gamma_2$ orbit is to say that
\begin{equation}
\label{GAMMA2}
A'=\frac{a A+b}{c A+d}; \hskip 30 pt \left[\matrix{a&b \cr c&d}\right] \in \Gamma_2.
\end{equation}

\begin{lemma}
If $A$ and $A'$ belong to the same $\Gamma_2$-orbit, then
$C_A$ and $C_A$ are asymptotically similar.  In
particular, $\dim(C_A)=\dim(C_A^{\#})$.
\end{lemma}

\startproof
Recall that $C_A$ is defined by the formula
\begin{equation}
C_A=\bigcup_{\kappa \in \Pi} (X(\kappa),1); \hskip 30 pt
X(\kappa)=\sum_{i=0}^{\infty}2k_i \lambda_i; \hskip 30 pt
\lambda_i=|Aq_i-p_i|.
\end{equation}

We say that two
 sequences $\{\lambda_i\}$ and $\{\lambda_i'\}$ are
{\it asymptotically proportional\/} if there is an
integer $m$ and a constant $C$ such that
\begin{equation}
\lim_{k \to \infty} \frac{\lambda'_k}{\lambda_{k+m}}=C \not =0 .
\end{equation}
The integer $m$ just serves to shift the terms
appropriately.

\begin{lemma}
\label{apro}
Suppose that $A$ and $A'$ are $\Gamma_2$-equivalent.  Then
the corresponding sequences $\{\lambda_k\}$ and $\{\lambda_k'\}$ are
asymptotically proportional.
\end{lemma}

\startproof
Let $T \in \Gamma$ be such that $T(A)=A'$, as in
Equation \ref{GAMMA2}.
Let $T'(A)$ be the derivative of $T$ at $A$.
  Our construction of the
inferior sequence is such that $T(A_{m+k})=A'_k$ for some $m$ and
all sufficiently large $k$.  Therefore
\begin{equation}
\label{pro1}
\lim_{k \to \infty}\frac{q'_k \lambda'_k}{q_{k+m}\lambda_{k+m}}=
\lim_{k \to \infty} \frac{\|T(A)-T(A_k)\|}{\|A_{k+m}-A\|}=\|T'(A)\|.
\end{equation}
We compute
$$T(p/q)=\frac{ap+bq}{cp+dq}.$$
It is an exercise in modular arithmetic to show that the fraction on
the right is already in lowest terms.  Therefore
\begin{equation}
\label{pro2}
\lim_{k \to \infty} \frac{q'_k}{q_{k+m}}=\lim_{k \to \infty}
\frac{c p_{m+k}+d q_{m+k}}{q_{k+m}}=cA+d.
\end{equation}
Combining Equations \ref{pro1} and \ref{pro2}, we get
\begin{equation}
\lim_{k \to \infty}\frac{\lambda'_k}{\lambda_{k+m}}=\frac{\|T'(A)\|}{cA+d}
\end{equation}
This completes the proof.
\endproof

\begin{corollary}
\label{ergodic2}
If $A$ and $A'$ are $\Gamma_2$-equivalent, then
$C_A$ and $C_A'$ are asymptotically similar.
\end{corollary}

\startproof
If $A$ and $A'$ are $\Gamma_2$-equivalent, then we have an
obvious map
\begin{equation}
\sum_{i=k_0}^{\infty} k_{m+i}\lambda_{m+i} \to
\sum_{i=k_0}^{\infty} k_i \lambda_i',
\end{equation}
which is defined if $k_0$ is taken large enough.  This
map makes sense because the corresponding sequence spaces
are the same.  Given the asymptotic proportionality
of the sequences, the map above is $f(k_0)$-bi-lipschitz.
Here $f(k_0)$ is a function that converges to $1$ as
$k_0$ to $\infty$. 
\endproof

\begin{lemma}
The function $A \to \dim(C_A)$ is Borel measurable.
\end{lemma}

\startproof
When $A=p/q$, we define $C_A=O_2(J) \cap I$.
Here $J$ is the interval of length $2/q$ in $I$ whose left
endpoint is $(0,1)$.  Thus, $C_A$ is just a thickened version
of part of the fundamental orbit.  Having made this definition,
we define $C$ as in Equation \ref{butterfly}.
In the proof of Lemma \ref{length0} we produced a covering
${\cal C\/}_n$ of $C_A$ by intervals all having the same
length.  One can extend this definition to the rational
case in a fairly obvious way.
 Let $C_A^{(n)}$ denote the union of these intervals.
Let $C^{(n)}$ be the corresponding union, with 
$C_A^{(n)}$ replacing $C_A$ in Equation \ref{butterfly}.
The sizes and positions of the intervals in $C_A^{(n)}$
vary with $A$ in a piecewise continuous way.  Hence $C^{(n)}$
is a Borel set.  Hence $C=\cap C^{(n)}$ is a Borel set.
Then $C_A$ is obtained by intersecting a Borel subset
of $[0,1]^2$ with the line $\{y=A\}$.  By the
Borel Slicing Property, $u$ is Borel measurable.
\endproof

\noindent
{\bf Remark:\/} My website has a picture of
the beautiful set $C$.

\subsection{The Dimension Function}

In this section we study the function $u(A)=\dim(C_A)=\dim(U_A)$.

\begin{lemma}
$u$ is almost everywhere constant.
\end{lemma}

\startproof
We have already seen that $u$ is constant on
$\Gamma_2$ orbits.   Since $\Gamma_2$ acts
ergodically on $\R \cup \infty$, we see that
$u$ is almost everywhere constant.
\endproof

\noindent
{\bf Remarks:\/} 
We would guess that $\dim(C_A)=1$ for almost
all $A$.  We don't know.
\newline

Now we derive some corollaries of the dimension
formula.  Say that $A$ is {\it superior\/} if all
the terms in the inferior sequence are superior.

\begin{theorem}
\label{supbound}
Let $A$ be a tame parameter.  Let $D=D_A$.  Then
$$u(A) \leq \frac{D}{D+\log(\sqrt 5/2)}.$$
If $A$ is also superior, then
$$u(A) \geq \frac{D}{D+\log(2)}.$$
\end{theorem}

\startproof
The upper bound follows from Equation \ref{DIO4} and
our dimension formula.  Now we prove the
lower bound.
Referring to the inferior sequence $\{p_n/q_n\}$ and
the inferior renormalization sequence $\{d_n\}$, we
always have 
$q_{n+1} < 2(d_n+1)q_n$.
This bound directly applies to the superior
sequence when $A$ is superior.
By induction, $q_n \leq 2D_n$.
Hence $Q \leq D+\log 2$.
Our bound follows immediately.
\endproof

\begin{lemma}
\label{subexp}
Let $A$ be a superior parameter whose renormalization
sequence $\{d_n\}$ diverges to $\infty$.  If
$d_{n+1}/d_n$ grows sub-exponentially, $u(A)=1$.
\end{lemma}

\startproof
The same argument as in Lemma \ref{tamebound} shows that
\begin{equation}
\lambda_n>(h_n q_n)^{-1}.
\end{equation}
Here $h_n$ grows sub-exponentially.
From Equation \ref{DIO3}, we get
\begin{equation}
G_n=2\lambda'_n>2\lambda_n>2(h_{n}q_{n})^{-1}.
\end{equation}
Therefore
$$
\lim_{n \to \infty} \frac{\log(D_n)}{\log(G_n^{-1})} \geq
\lim_{n \to \infty} \frac{\log(D_n)}{\log(h_nq_n)}=^*$$
\begin{equation}
\lim_{n \to \infty} \frac{\log(D_n)}{\log(q_n)} \geq
\lim_{n \to \infty} \frac{\log(D_n)}{\log(D_n)+\log(2)}=1.
\end{equation}
The starred equality comes from the sub-exponential
growth of $h_n$.
Essentially the same derivation as above now shows that
$u(A) \geq 1$.  But, of course $u(A) \leq 1$ as well.
Hence $u(A)=1$.
\endproof

\begin{lemma}
$u$ maps $(0,1)-\Q$ onto $[0,1]$.
\end{lemma}

\startproof
By the previous result, we can get $u(A)=1$.
It is easy to get $u(A)=0$ by taking
$A$ so that the IERS is $2,0,...,0,2,0,...,0...$ where
the number of $0$s grows rapidly enough.
(See \S \ref{groden} for a definition of the IERS.)

Suppose we want find $A$ such that $u(A)=x \in (0,1)$.
Let $A$ be the parameter whose IERS is $N,N,N...$.
Here $N$ must be even because the $0$th term is
even.  Choosing $N$ large enough, we can arrange
that $u(A)>x$, by Lemma \ref{supbound}. 
Let $B_m$ denote the parameter whose IERS is
$N,0_m,N,0_m,...$.  Here $0_m$ represents $m$ zeros in a row.
As $m \to \infty$, we have $u(B_m) \to 0$.  Thus, we can
find an integer $m$ such that $u(B_{m+1})<\delta_0< u(B_m)$.
(We're already done if we have equality on either side.)

Given a binary sequence $\eta=\{\epsilon_k\}$ we define
$A(\eta)=N,0_{m+\epsilon_1},N,0_{m+\epsilon_2},...$
The parameter $A(\eta)$ is tame, and $D(\eta)=\log N$,
independent of $\eta$.
Letting $\eta_0$ and
$\eta_1$ denote the $0$ sequence and the $1$ sequence,
we have $Q(\eta_0)<Dx$ and $Q(\eta_1)>Dx$.  Essentially
by the intermediate value theorem, we can adjust
$\eta$, so that $Q(\eta)=Dx$.   Then $A(\eta)$ is
the desired parameter.
\endproof

Since $u$ is $\Gamma_2$-invariant, and
$\Gamma_2$ orbits are dense in $(0,1)$, the
function $u$ maps any open subset of $(0,1)-\Q$
onto $[0,1]$.

\subsection{Example Calculations}

{\bf Example 1:\/}
Let $A=\sqrt 5-2$, the Penrose kite parameter.
The inferior sequence is
$$ {\bf \frac{1}{1}\/} \leftarrow \frac{1}{3} {\bf \leftarrow \frac{1}{5}\/} 
\leftarrow \frac{3}{13} \leftarrow {\bf \frac{5}{21}\/} \leftarrow \frac{13}{55}
 {\bf \leftarrow \frac{21}{89}\/} \leftarrow  \frac{55}{233}\leftarrow {\bf \frac{89}{377}\/} \ldots$$
The superior sequence is
$$ \frac{1}{1}, \frac{1}{5},
\frac{5}{21}, \frac{21}{89} \ldots \hskip 30 pt
p_{n+1}=q_n; \hskip 30 pt q_{n+1}=4q_n+p_n.$$
The inferior renormalization sequence
is $2,0,2,0...$.  The renormalization sequence is $2,2,2...$.  Hence $D=\log(3)$.  The
superior sequence satisfies the recurrence relation
$$
p_{n+1}=q_n; \hskip 30 pt q_{n+1}=4q_n+p_n.
$$
This gives $Q=\log(\phi^3)$.  Theorems \ref{dimformula} combines with
the modularity to show that
$$
\left[\matrix{a&b \cr c&d}\right] \in \Gamma_2; \hskip 15 pt
A=\frac{a \sqrt 5 +b}{c \sqrt 5 + d} \in (0,1)\hskip 20 pt \Longrightarrow \hskip 5 pt
u(A)=\frac{\log(3)}{\log(\phi^3)}
$$
\noindent
{\bf Example 2:\/}
The renormalization sequence for the parameter $A=E-2$
considered at the end of \S \ref{UNIVERSAL} is
$3,5,7,9...$.  
The example here satisfies Lemma \ref{subexp}.  Therefore
\begin{equation}
\label{sqe}
\left[\matrix{a&b \cr c&d}\right] \in \Gamma_2; \hskip 15 pt
A=\frac{a E +b}{c E + d} \in (0,1)\hskip 20 pt \Longrightarrow \hskip 5 pt
u(A)=1.
\end{equation}
Again, it seems that $E=e$.

\newpage

{\bf {\Huge Part VI\/}\/}
\newline

\begin{itemize}

\item In \S \ref{copyproof2} we prove the Copy Theorem from
\S \ref{pivots}.

\item In \S \ref{evenpivot} we define what we mean by
the pivot arc relative to an even rational kite parameter.
Along the way we will prove another version of the
Diophantine Lemma from \S \ref{diotheorem}. The
Diophantine Lemma works for pairs of odd rationals,
and the result here works for pairs of Farey-related
rationals, even or odd.

\item In \S \ref{pivotproof} we prove the Pivot Theorem
from \S \ref{statepivot}.  
The Pivot Theorem works in both
the even and the odd case, and is proved in an inductive way
that requires both cases.

\item In \S \ref{period} we prove the Period Theorem.

\item In \S \ref{gold} we prove Statement 1 of the Comet Theorem.

\end{itemize}

\newpage

\section{Proof of the Copy Theorem}
\label{copyproof2}

\subsection{A Formula for the Pivot Points}
\label{formulapivot}

Let $A$ be an odd rational.  Let $A_-$ be as
in Equation \ref{induct0}.  Let
$V_-=(q_-,-p_-)$.   Here we give a formula for
the pivot points $E^{\pm}$ associated to $A$.

\begin{lemma}
\label{swap}
The following is true.
\begin{itemize}
\item If $q_-<q_+$ then
$E^++E^-=-V_-+(0,1)$.
\item If $q_+<q_-$ then
$E^++E^-=V_++(0,1)$.
\end{itemize}
\end{lemma}

\startproof
We will establish this result inductively.  
Suppose first that $1/1 \leftarrow A$. Then
$$A=\frac{2k-1}{2k+1}; \hskip 30 pt
E^-=(-k,k); \hskip 30 pt E^+=(0,0); \hskip 30 pt
V_-=(k,-k+1).
$$
$$A_-=\frac{k-1}{k}; \hskip 30 pt
q_-=k-1<k=q_+.
$$
The result works in this case.

In general, we have $A=A_2$ and
 $A_0 \leftarrow A_1 \leftarrow A_2.$
There are $4$ cases, depending on Lemma \ref{indX}.
Here the index is $m=1$.
We will consider Case 1.  The other cases are similar.
By Case 1, we have $(q_1)_+<(q_1)_-$.   Hence, by
induction
$$E_1^++E_1^-=(V_1)_++(0,1).$$
Since $A_1<A_2$ we have
$$E_2^-=E_1^-; \hskip 30 pt
E_2^+=E_1^++d_1V_1.$$
Therefore
$$E_2^++E_2^-=(V_1)_++d_1V_1+(0,1)=(V_2)_++(0,1).$$
The last equality comes from
Case 1 of Lemma \ref{indX}.
As we remarked after stating Lemma \ref{indX}, this
result works for both numerators and denominators.)
In Case 1, we have $(q_2)_+<(q_2)_-$, so the result holds.
\endproof

\begin{lemma}
\label{stickout}
$E^-$ lies to the left of $R_1$ and
$E^+$ lies to the right of $R_1$.
\end{lemma}

\startproof 
Let $\pi_1$ denote the projection to the first coordinate.
One or the other bottom vertices of $R_1$ is $(0,0)$.
We will consider the case when the left bottom
vertex is $(0,0)$.  In all cases one checks
easily from our definitions that $\pi_1(E^-) \leq -1$.
Hence $E^-$ lies to the left of $R_1$.

Consider the right side.
We have $q_+<q_-$ in our case. By
Case 2 of Lemma \ref{swap}, and the
result for the left hand side, we have
$$\pi_1(E^+) \geq \pi_1(V_+)+1.$$
But $V_+$ lies on the line extending the
bottom right edge of $R_1$, exactly
$1/q$ vertical units beneath the bottom
edge of $R_1$. This right edge
has slope greater than $1$.  Finally, the
line connecting $V_+$ to $\pi_1(E^+)$ has
nonpositive slope because $E^+$ is a
low vertex lying to the right of
$V_+$. From all this geometry,
we see that $E_1^+$ lies to the right of
$R_1$.
\endproof

While we are in the neighborhood, we clear up
a detail from Part V.
\newline
\newline
{\bf Proof of Lemma \ref{statepivot}\/}:
We will prove this result inductively.  Suppose
that $A_1 \leftarrow A_2$, and the result is true
for $A_1$.   We consider the case when
$A_1<A_2$.  The case when $A_1>A_2$ has the
same treatment.  When $A_1<A_2$, we have
$E_1^-=E_2^-$, so certainly the bound holds for
$E_1^-$.   On the other hand, we have
\begin{equation}
\pi_1(E_2^+)=\pi_1(E_1^+)+d_1 q_1 \hskip 30 pt
d_1={\rm floor\/}\bigg(\frac{q_2}{2q_1}\bigg)
\end{equation}
There are two cases to consider.  Suppose
first that
$\delta_1={\rm floor\/}(q_2/q_1)$ is odd.
In this case 
$$(2d_1+1)q_1<q_2; \hskip 30 pt \Longrightarrow
\hskip 30 pt
d_1q_1<\frac{q_2}{2}-\frac{q_1}{2}.$$
The first equation implies the second.
Hence, by induction
$$\pi_1(E_2^+)<\frac{q_1}{2}+\frac{q_2-q_2}{2}<\frac{q_2}{2}.$$

Suppose that $\delta_1$ is even.  Then we have Case 2 of
Lemma \ref{indX}, applied to the index $m=1$.
This is to say that
$(q_1)_-<(q_1)_+$.
From our formula above,
the first coordinate of $E_2^-+E_2^+$ is negative.
Hence
$$|\pi_1(E^-)|>|\pi_1(E^+)|.$$
This fact finishes the proof.
\endproof

\subsection{Good Parameters}

Our pivot points are well defined vertices, but
so far, we don't know that the pivot arc is
well defined.  That is, we don't know that 
$E^-$ and $E^+$ are actually vertices of
$\Gamma$.  These points might be vertices
of some other component of $\widetilde \Gamma$.
To start things off right, we deal with the
base case.

\begin{lemma}
If $1/1 \leftarrow A$, then the pivot arc is
well-defined relative to $A$.
\end{lemma}

\startproof
Here $A=(2k-1)/(2k+1)$ for $k \in \N$.
In \S \ref{decompend}, 
we showed that the line segment connecting
$(0,0)$ to $(-k,k)$ is contained in the arithmetic
graph.  So, the pivot arc is well defined.
\endproof

Now we consider the general case.
Let $A_1$ be an odd rational.
For each integer $\delta_1 \geq 1$, there is a unique
odd rational $A_2=A_2(\delta_1)$ 
such that $A_1 \leftarrow A_2$ and
$$\delta_1={\rm floor\/}\bigg(\frac{q_2}{q_1}\bigg).$$
Lemma \ref{indX} gives the recipe for how to construct
$A_2$.  As in Equation \ref{renorm}, we define
$$d_1={\rm floor\/}\bigg(\frac{q_2}{2q_1}\bigg).$$

Recall that $E_j^{\pm}$ are the pivot points associated to
$\Gamma_j$.   Let $P\Gamma_1(\delta_1)$ denote the arc of
$\Gamma_1$ whose endpoints are $E_2^-$ and $E_2^+$.

\begin{lemma}
$P\Gamma_1(\delta_1)$ is a well defined arc of
$\Gamma_1$.
\end{lemma}

\startproof
Suppose that $A_1<A_2$.  When $A_2<A_1$ the proof is similar.
Then, by Equation \ref{pivotplus}, we have
$$E_2^-=E_1^-; \hskip 30 pt
E_2^+=E_1^++d_1 V_1; \hskip 30 pt V_1=(q_1,-p_1).$$
But $\Gamma_1$ is invariant under translation by $V_1$.
Hence $E_2^{\pm}$ is a vertex of $\Gamma_1$.
\endproof

Call $A_1$ a {\it good parameter\/} if
\begin{equation}
P\Gamma_1 \subset \Delta_1(I).
\end{equation}
Here $\Delta_1(I)$ is the region from the
Diophantine Lemma, defined relative to
the pair $(A_1,A_2(1))$.   We call $I$ the
{\it base interval\/}.  We will give a
formula below.

\begin{lemma}
If $A_1$ is good, then the Copy Theorem holds for
$A_1$ and $A_2(1)$.
\end{lemma}

\startproof
Note that $P\Gamma_1(1)=P\Gamma_1$, the pivot arc of
$\Gamma_1$. The pivot points
do not change in this case: $E_1^{\pm}=E_2^{\pm}$. 
So, if $A_1$ is good then the Diophantine
Lemma immediately implies that $P\Gamma_1(1)=P\Gamma_1 \subset \Gamma_2$.
But then, there is an arc of $\Gamma_2$ that connects
$E_2^-$ to $E_2^+$, the two endpoints of $P\Gamma_1(1)$.
This shows that the Pivot arc for $A_2$ is well defined,
and that this pivot arc is a subarc of $\Gamma_1$. 
\endproof

\begin{lemma}
\label{oddodd}
If $A_1$ is good, then the Copy Theorem holds for
$A_1$ and $A_2(3)$.
\end{lemma}

\startproof
Let $A_0$ be such that
the sequence $A_0 \leftarrow A_1 \leftarrow A_2(1)$ is a
fragment of the inferior sequence.
We will consider the case when
$A_0<A_1$.  In this case $A_1<A_2(1)$ by Lemma \ref{indX}.
The base interval is given by
\begin{equation}
\label{oddfuture}
I=[-q_1+2,q_1+(q_2)_+-2]=[-q_2+2,q_1+(q_1)_+-2].
\end{equation}
The first equality is Lemma \ref{keycomp1}. The
second equality is Case 1 of Lemma \ref{indX},
with $d_1=0$.

Let $R_j=R_j(A_1)$, as in the Decomposition Theorem for $A_1$.
As in Lemma \ref{stickout}, we know that
$R_1$ lies to the right of the origin
and $R_2$ to the left. (This is because $(q_1)_+<(q_1)_-$ in the case
we are considering.)
The arc $P\Gamma_1(3)$ is obtained from
$P\Gamma_1(1)$ by concatenating one period of
$\Gamma_1$ to the right.

We claim that
\begin{equation}
\label{decomp10}
P\Gamma_1(3)=P\Gamma_1 \cup \gamma \cup \Big(P\Gamma_1+V_1\Big);
\hskip 30 pt \gamma \in (R_2+V_1).
\end{equation}

\begin{center}
\resizebox{!}{1.7in}{\includegraphics{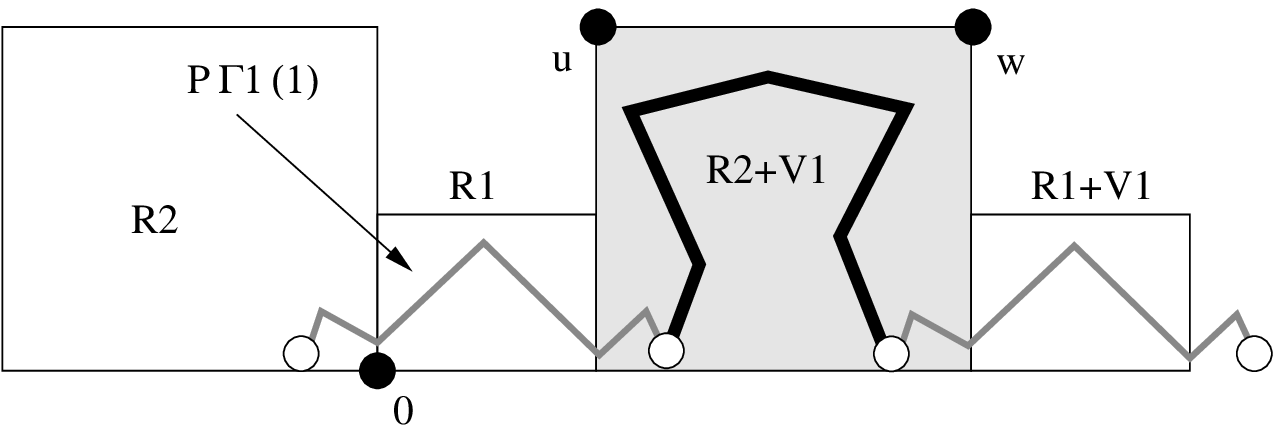}}
\newline
{\bf Figure 27.1:\/} Decomposition of $P\Gamma_1(3)$.
\end{center}

Here is the proof.
By Lemma \ref{stickout}, the arc $P\Gamma_1$ completely
crosses $R_1$.  The left endpoint lies in $R_2$ and
the right endpoint lies in $R_2+V_1$, the translate of
$R_2$ that lies on the other side of $R_1$.   By
symmetry, one endpoint of $P\Gamma_1(1)$ enters
$R_2+V_1$ from the left and one endpoint of
$P\Gamma_1(1)+V_1$ enters $R_2+V_1$ from the right.
The arc $\gamma$
joins two points already in $R_2+V_1$.  This 
arc cannot cross out of $R_2+V_1$, by Lemma \ref{symm2}.

Now we know that Equation \ref{decomp10} is true.
Let $A_2=A_2(1)$ and $A_2^*=A_2(3)$.
We attach a $(*)$ to objects associated to $A_2^*$.
Let $I$ be the base interval. Let $I^*$ denote
the interval corresponding to the pair $(A_1,A_2^*)$.
By Lemma \ref{indX}, we have
$(q_2^*)_+=q_1+(q_1)_+.$  Hence,
by Lemma \ref{keycomp1} and by definition,
\begin{equation}
\label{oddfuture0}
I^*=[-q_1+2,2q_1+(q_1)_+-2]=[I_{\rm left\/},I_{\rm right\/}+q_1].
\end{equation}
We have $P\Gamma_1 \subset \Delta_1(I) \subset \Delta_2(I)$.
The first containment is the definition of goodness.
Any $v^* \in P\Gamma_1+V_1$ has the form
$v+V_1$, where $v \in P\Gamma_1$.  By Lemma
\ref{Hcalc} we have
$$G_1(v^*)=G_1(v)+q_1; \hskip 30 pt
H_1(v^*)=H_1(v)+q_1.$$
Hence $v \in \Delta_1(I)$ implies $v^* \in \Delta(I^*)$.
It remains to deal with the arc $\gamma$.

We will use the same argument that we used in \S \ref{most}.
Let $u$ and $w$ respectively be the upper left and upper
right vertices of $R_2+V_1$.  We have
\begin{equation}
u\approx W_1+\frac{(q_1)_+}{q_1}V_1; \hskip 30 pt
w=W_1+V_1.
\end{equation}
Here the vectors are as in Equation \ref{boxvectors}, as usual. 
The approximation is good to within $1/q_1$. To avoid approximations,
we consider the very slightly altered parallelogram
$\widetilde R_2+V_1$.  The vertices are
\begin{equation}
(V_1)^+; \hskip 30 pt \widetilde u=W_1+\frac{(q_1)_+}{q_1}V_1; \hskip 30 pt
V_1; \hskip 30 pt w=V_1+W_1.
\end{equation}
Each vertex of the new parallelogram is within $1/q_1$
of the corresponding old parallelogram.  Using
the Adjacent Mismatch Principle, it suffices to do the
calculation in $\widetilde R_2+V_1$. 
The following calculation combines with
the Diophantine Lemma to show that
$\gamma \subset \Gamma_2(A_2^*)$.
\begin{equation}
\label{pivotbounds}
G_1(\widetilde u)-(-q_1)=(2q_1+q_+)-H_1(w)=
q_1+(q_1)_+-\frac{q_1^2}{p_1+q_1} \geq 2.
\end{equation}
This completes the proof.
\endproof

\begin{lemma}
If $A_1$ is good and $\delta_1$ is odd, then the Copy Theorem holds for
$A_1$ and $A_2(\delta_1)$.
\end{lemma}

\startproof
Now consider the case when $\delta_1=5$.  In this case,
$P\Gamma_1(5)$ is obtained by concatenating $2$ periods
of $\Gamma_1$ to the right of $P\Gamma_1(1)$.  We have 
decomposition of the form
\begin{equation}
P\Gamma_1(5)=P\Gamma_1(1) \cup \gamma \cup \Big(P\Gamma_1(1)+2V_1\Big); \hskip 30 pt
\gamma \subset (R_2+V_1) \cup (R_2+2V_1).
\end{equation}
Here $\gamma$ is contained in a parallelogram that is twice as long as
in the case $\delta=3$.  The calculations are exactly the same in this
case.  The key point is that $I^*=[a,b+2q_1]$.
The cases $\delta=7,9,11...$ have the same treatment.
\endproof

\begin{lemma}
If $A_1$ is good, then the Copy Theorem holds for
$A_1$ and $A_2(2)$.
\end{lemma}

\startproof
In this case,
$P\Gamma_1(2)$ is obtained from $P\Gamma_1$
by concatenating one period of $\Gamma_1$
to the left.  See Figure 27.2 below.
We have the decomposition
\begin{equation}
P\Gamma_1(2)=\Big(P\Gamma_1(1)-V_1\Big) \cup \gamma \cup P\Gamma_1(1);
\hskip 30 pt
\gamma \subset R_2.
\end{equation}
The proof is the same as in Lemma \ref{oddodd}.

\begin{center}
\resizebox{!}{1.7in}{\includegraphics{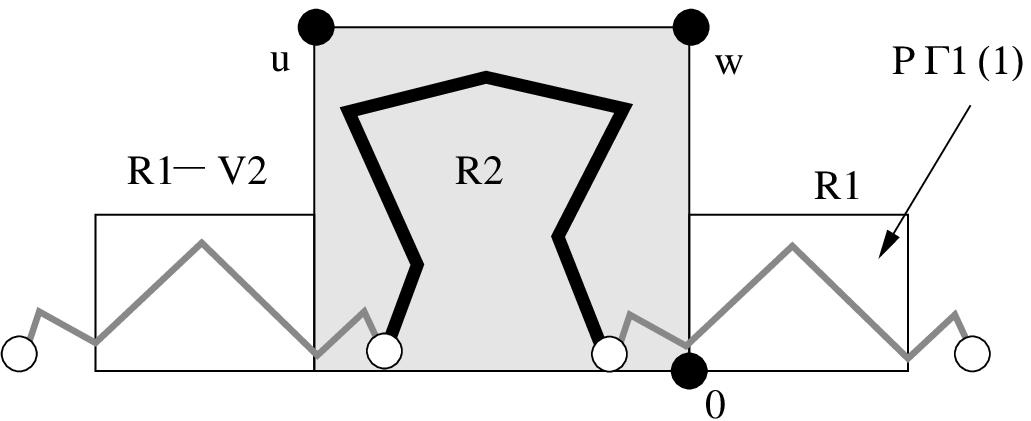}}
\newline
{\bf Figure 27.2:\/} Decomposition of $P\Gamma_1(2)$.
\end{center}

We use the same notational conventions as in the odd case.
The same argument as above works here, provided that
we can get the right estimates on the top vertices
$u$ and $w$ of $R_2=R_2(A_1)$.  Case 4 of Lemma \ref{indX}
tells us that
$$(q_2^*)_-=q_1^-.$$
Combining this fact with Lemma \ref{keycomp1}, we get
\begin{equation}
I^*=[-q_1-(q_2^*)_-+2,q_1-2]=[-q_1-(q_1)_-+2,q_1-2].
\end{equation}
We have
\begin{equation}
u \approx \frac{-(q_1)_-}{q_1}V_1 +W_1; \hskip 30 pt
w=W_1.
\end{equation}
Again, the approximation holds up to $1/q_1$.
To avoid approximations, we use the modified parallelogram
$\widetilde R_2$ with vertices
\begin{equation}
\frac{-(q_1)_-}{q_1} V_1; \hskip 30 pt
\widetilde u=\frac{-(q_1)_-}{q_1}V_1 +W_1;
 \hskip 20 pt
(0,0); \hskip 20 pt w=W_1.
\end{equation}
Again, this is justified by our Adjacent Mismatch Principle.
The following estimate combines with the Diophantine
Lemma to show that $\gamma \subset \Gamma_2(A_2^*)$.
\begin{equation}
G_1(\widetilde u)-(-q_1-(q_1)_-)=q_1-H(w)=
q_1-\frac{q_1}{p_1+q_1} \geq 2.
\end{equation}
As in \S \ref{most}, this estimate holds as long
as $p_1 \geq 3$ and $q_1 \geq 7$.  We handle the
few exceptional cases as we did in \S \ref{tricks}.
\endproof

\begin{lemma}
If $A_1$ is good and $\delta_2$ is even, 
then the Copy Theorem holds for
$A_1$ and $A_2(\delta_1)$ when $\delta_1$.
\end{lemma}

\startproof
The cases $\delta_1=4,6,8...$ relate to the case
$\delta_1=2$ exactly as the cases $\delta_1=5,7,9...$
relate to the cases $\delta_1=3$.
\endproof

\subsection{The End of the Proof}

It remains to show that any odd rational is good.
We will give an inductive argument.

\begin{lemma}
If $1/1 \leftarrow A$, then $A$ is good.
\end{lemma}

\startproof
In this case, Lemma \ref{indX} tells us that
$1/1 >A>\widehat A$.  (The first inequality is obvious.)
We have
$$A=\frac{2k-1}{2k+1}; \hskip 30 pt
\widehat A=\frac{4k-3}{4k+1};  \hskip 30 pt
\widehat q_-=k.$$
By Lemma \ref{keycomp1}, we have
$$I=[-q-q_-+2,q-2]=[-3k+1,2k-1]$$
The left vertex of $P\Gamma_1$ is $u=(-k,k)$
and the right vertex is $v=(0,0)$.
We compute
$$G(u)=-k-1 \geq -3k+1; \hskip 30 pt
H(w)=0 \leq 2k-1.$$
The extreme case happens when $k=1$.
\endproof

\begin{lemma}
$A=p/q$ is good if $q<20$ or if $p=1$.
\end{lemma}

\startproof
We check the case $q<20$ by hand.
If $p=1$, the pivot arc is just the
edge connecting $(-1,1)$ to $(0,0)$ whereas
the interval $I$ contains $[-q,q]$, a huge
interval.  This case is obvious.
\endproof

Now we establish the inductive step.  Suppose
that $A_1 \leftarrow A_2$ and that $A_1$ is good.
Having eliminated the few exceptional cases by
our result above, our argument in the previous section shows
that $P\Gamma_2 \subset \Delta_1(I_1)$.  Here
$I_1$ is the interval based on the constant
$\Omega(A_1,A_2)$. This is the Diophantine constant
defined in \S \ref{omegadef} relative to the
pair $(A_1,A_2)$. 
To finish of the proof of the Copy Theorem, we just
have to establish the following equation.
\begin{equation}
\label{copygoal}
P\Gamma_2 \subset \Delta_2(I_2),
\end{equation}
where $I_2$
is the different interval based on the pair
$A_2 \leftarrow A_3$, with $\delta(A_2,A_3)=1$.
Here we establish two basic facts.

\begin{lemma}
$I_1 \subset I_2$, and either endpoint of $I_1$ is more than
$1$ unit from the corresponding endpoint of $I_2$. 
\end{lemma}

\startproof
By Lemma \ref{keycomp1}, applied to both parameters, we have
$$I_1 \subset [-q_2+3,q_2-3] \subset [-q_2-2,q_2-2] \subset I_2.$$
This completes the proof.
\endproof

\begin{lemma}
\label{lipschitz}
$|G_1(v)-G_2(v)|<1$ and $|H_1(v)-H_2(v)|<1$ for $v \in \Delta_1(I_1)$.
\end{lemma}

\startproof
From Lemma \ref{keycomp1} and a bit of geometry, we get the bound
\begin{equation}
(m,n) \in \Delta_1(I_1) \hskip 30 pt
\Longrightarrow \hskip 30 pt
\max(|m|,|n|) \leq q_2.
\end{equation}
Looking at Equation \ref{defG}, we see that
$$
G(m,n)=\bigg(\frac{1-A}{1+A},\frac{-2}{1+A}\bigg) \cdot (m,n)=(G_1,G_2) \cdot(m,n)
$$
\begin{equation}
H(m,n)=\bigg(\frac{1+4A-A^2}{(1+A)^2},\frac{2-2A}{(1+A)^2}\bigg) \cdot (m,n)=(H_1,H_2) \cdot (m,n)
\end{equation}
A bit of calculus shows that
\begin{equation}
|\partial_A G_j| \leq 2; \hskip 30 pt
|\partial_A H_1| \leq 6; \hskip 30 pt
|\partial_A H_2| \leq 2.
\end{equation}
Since $A_1 \leftarrow A_2$, we have
\begin{equation}
|A_1-A_2|=\frac{2}{q_1q_2}.
\end{equation}
Putting everything together, and using basic calculus, we arrive
at the bound
\begin{equation}
|G_1(v)-G_2(v)|, |H_1(v)-H_2(v)|<\frac{16}{q_1}<1,
\end{equation}
at least for $q_1>16$.
\endproof

We have already remarked, during the proof of the Decomposition Theorem,
that no lattice point lies between the bottom of
$\Delta_2(I_2)$ and the bottom of $\Delta_1(I_2)$. 
Hence $F_1(v)>0$ iff $F_2(v)>0$.   Our two lemmas now show that
$\Delta_1(I_1) \subset \Delta_2(I_2)$.  This was our
final goal, from Equation \ref{copygoal}.

This completes the proof of the Copy Theorem.

\newpage

\newpage

\section{Pivot Arcs in the Even Case}
\label{evenpivot}

\subsection{Main Results}
\label{predecessor}
Given two rationals $A_1=p_1/q_1$ and
$A_2=p_2/q_2$, we introduce the notation
\begin{equation}
A_1 \bowtie A_2 \hskip 30 pt
\Longleftrightarrow \hskip 30 pt
|p_1q_2-q_1p_2|=1; \hskip 15 pt
q_1<q_2.
\end{equation}
In this case, we say that
$A_1$ and $A_2$ are {\it Farey related\/}.
We sometimes call $(A_1,A_2)$ a
{\it Farey pair\/}.

We have the notions of
{\it Farey addition\/} and {\it Farey subtraction\/}:
\begin{equation}
A_1 \oplus A_2=\frac{p_1+p_2}{q_1+q_2}; \hskip 40 pt
A_2 \ominus A_1 = \frac{p_2-p_1}{q_2-q_1}.
\end{equation}
Note that $A_1 \bowtie A_2$ implies that
$A_1 \bowtie (A_1 \oplus A_2)$ and that
$A_1$ is Farey related to $A_2 \ominus A_1$.

\begin{lemma}
Let $A_1$ be an even rational.  Then
there is a unique odd rational $A_2$ such
that $A_1 \bowtie A_2$ and $2q_1>q_2$.
\end{lemma}

\startproof
Equation \ref{induct0} works for both even and
odd rationals.  When $A_1$ is even, exactly one
of the rationals $(A_1)_{\pm}$ is also even.
Call this rational $A_1'$.  Then
$A_1' \bowtie A_1$.  We define
$A_2=A_1 \oplus A_1'$.  If $B_2$ was another
candidate, then $B_2 \ominus A'$ would be the
relevant choice of $(A_1)_{\pm}$.  Hence
$B_2=A_2$.
\endproof

We will write $A_1 \bowtie !\ A_2$ to denote
the relationship between $A_1$ and $A_2$ discussed
in the previous result. 
We can think about this relation in a different way.
Let $A$ be an odd rational.  Then 
either $A_- \bowtie !\ A$ or $A_+ \bowtie !\ A$ when
$A$ is an odd rational.  If $A_- \bowtie !\ A$ then
we write $A_+ \Leftarrow A$.  The relationship
implies that $2q_+<q$.  Likewise we write
$A_- \Leftarrow A$ when $2q_-<q$.  
Here is an example:  Let $A=3/7$.
Then  
$$A_+=1/2 \Leftarrow 3/7; \hskip 30 pt
A_-=2/5 \bowtie !\ 3/7.$$

So far, we have defined
pivot points and arcs for odd parameters.
Now we define them for even parameters.  We define

\begin{equation}
E^{\pm}(A_1):=E^{\pm}(A_2); \hskip 30 pt
A_1 \bowtie !\ A_2.
\end{equation}
This makes sense because we have already
defined the pivot points in the odd case.
We still need to prove that these
vertices lie on $\Gamma_1$.
We will do this below.

Assuming that the pivot points $E_1^{\pm}$ are vertices of
$\Gamma_1$, we define
$P\Gamma_1$ to be the lower arc of $\Gamma_1$ that
connects $E_1^-$ to $E_1^+$.  Since $\Gamma_1$
is a polygon in the even case, it makes sense to
speak of the lower arc.
Figure 28.1 shows an example.  Here
$P\Gamma_1=P\Gamma_2$.  We will show that
this always happens. 

\begin{center}
\resizebox{!}{5.2in}{\includegraphics{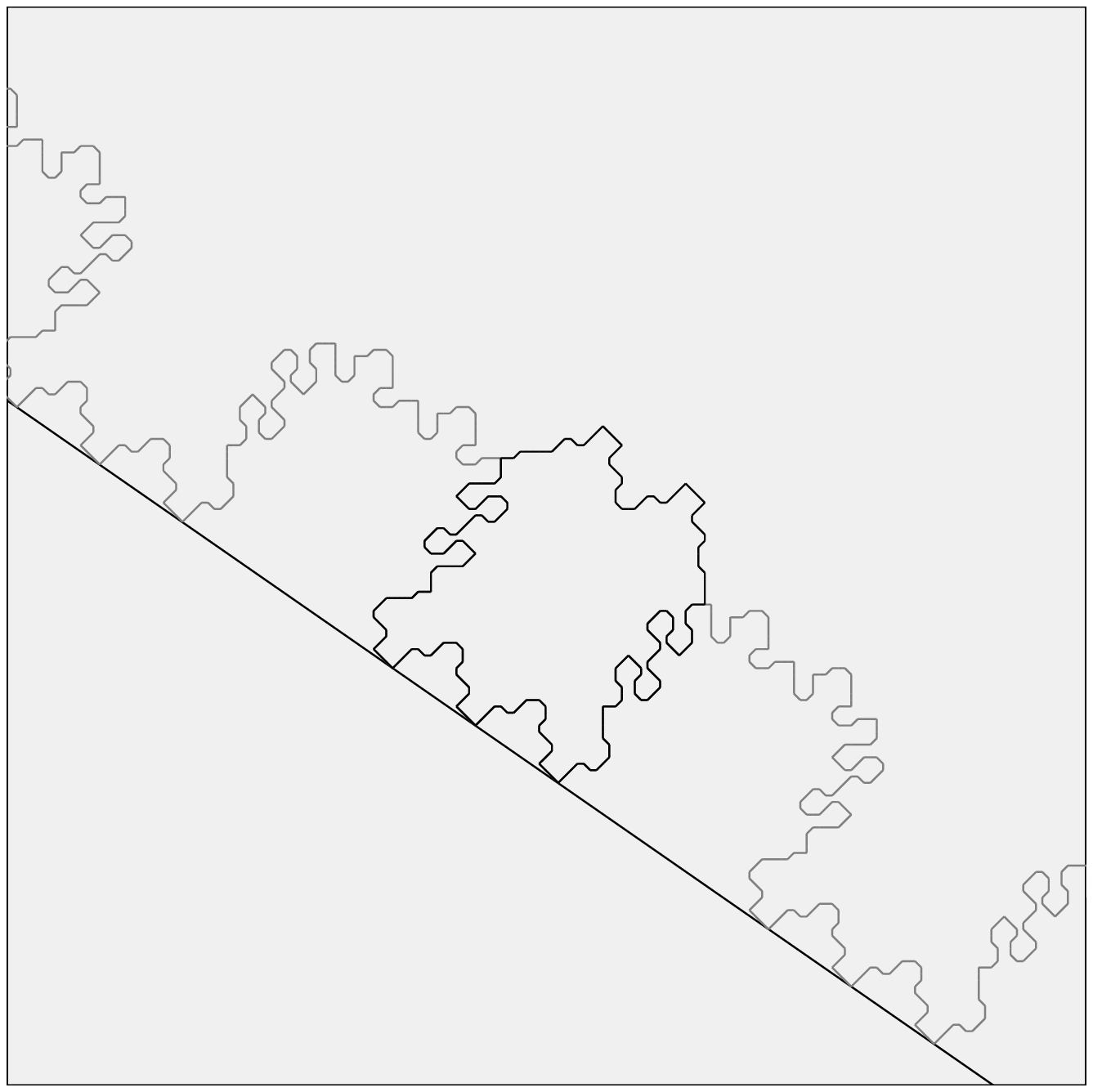}}
\newline
{\bf Figure 28.1:\/} $\Gamma(41/59)$ in grey and $\Gamma(25/36)$ in black
\end{center} 

In this chapter we prove the following results.

\begin{lemma}
\label{welldefined}
Let $A_1 \bowtie !\ A_2$.  Then $P\Gamma_1$ is well
defined and $P\Gamma_1=P\Gamma_2$.
\end{lemma}

\begin{lemma}[Structure]
\label{str1}
The following is true.
\begin{enumerate}
\item If $A_- \Leftarrow A$ then $E^+(A)=E^+(A_-)$.
\item If $A_+ \Leftarrow A$ then $E^-(A)=E^-(A_+)$.
\item If $A_- \Leftarrow A$ then $E^-(A)+V=E^-(A_-)+kV_-$ for some $k \in \Z$.
\item If $A_+ \Leftarrow A$ then $E^+(A)-V=E^+(A_+)+kV_+$ for some $k \in \Z$.
\end{enumerate}
\end{lemma}
The Structure Lemma is of crucial importance in our proof of
the Pivot Theorem and the Period Theorem. Here we illustrate its meaning
and describe a bit of the connection to the Pivot Theorem.

\begin{center}
\resizebox{!}{5.2in}{\includegraphics{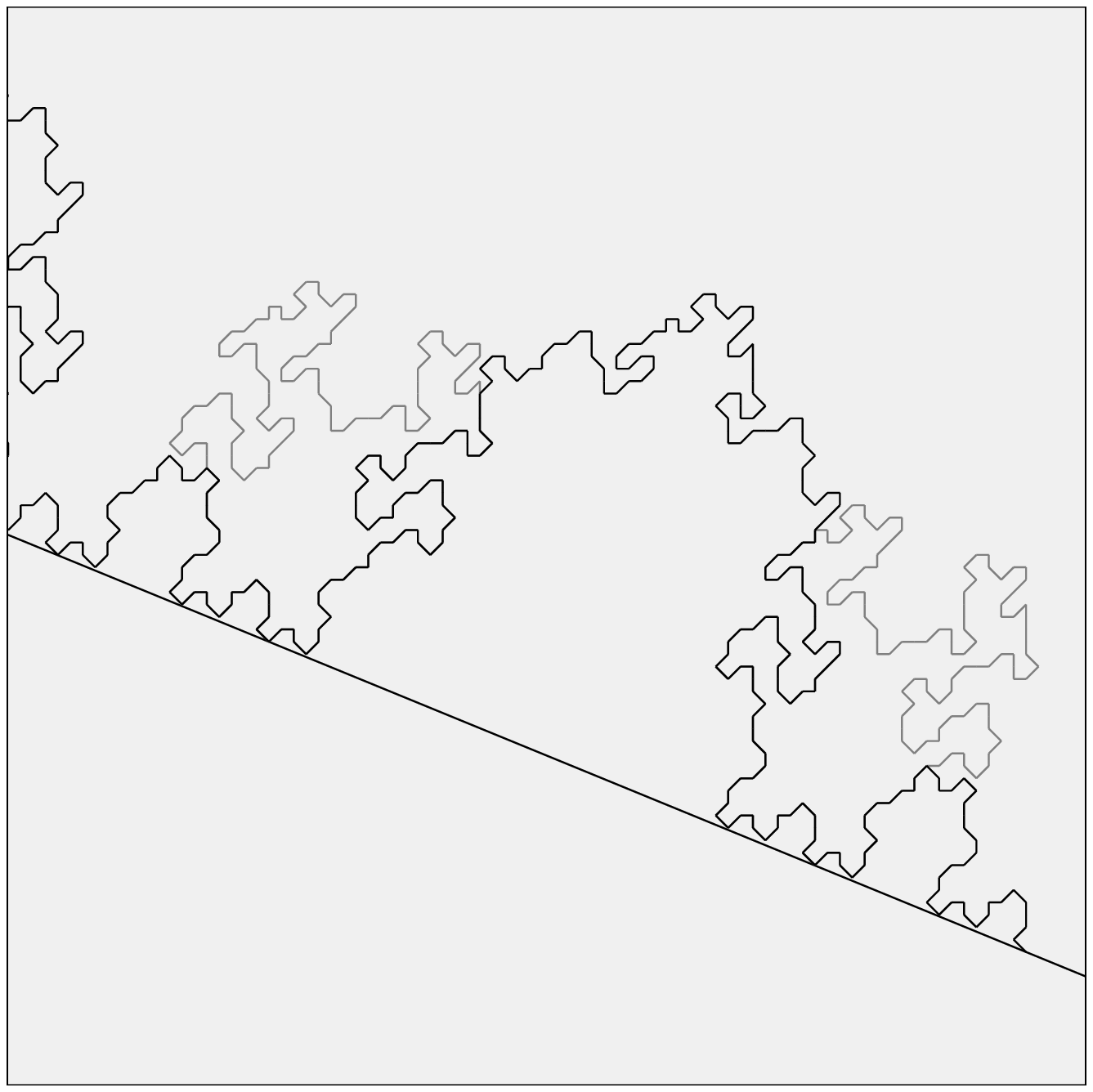}}
\newline
{\bf Figure 28.2:\/} $\Gamma(25/61)$ overlays several components
of $\widehat \Gamma(9/22)$.
\end{center} 

Figure 28.2 shows slightly more than one period of $\Gamma(25/61)$ in
black.  This black arc overlays $\Gamma(9/22)$ on the left and
$$\Gamma(9/22)+2V(9/22)$$ on the right.  Call these two grey
components {\it the eggs\/}.  Here $$9/22 \bowtie !\ 25/61.$$
The points 
$$E^+(25/61); \hskip 30 pt E^-(25/61)+V(25/61)$$ are the left
and right endpoints respectively 
of the big central hump of $\Gamma(25/61)$.
Call this black arc {\it the hump\/}.   The content of the
Structure Lemma (in this case) is that the endpoints of
the hump are simultaneously pivot points on the eggs.
The reader can draw many pictures like this on Billiard King.

The content of the Pivot Theorem for $25/61$ is that
the hump has no low vertices except its endpoints.
Note that the ends of the hump copy pieces of the eggs.
If we already understand the behavior of the eggs -- meaning
how they rise away from the baseline -- then we understand
the behavior of the ends of the hump.  The eggs are based
on a simpler rational.  In this way, the behavior of the
arithmetic graph for a simpler rational gives us information
about what happens for a more complicated rational. This
is (some of) the strategy for our proof of the Pivot Theorem.
In the first section of the next chapter we will
a long and somewhat informal discussion about the
remainder of the strategy.
\newline
\newline
{\bf Remarks:\/} \newline
(i)
In \S \ref{dpa} below we will give the precise
relationship between the two pivot arcs in the
cases of interest to us.  
\newline
(ii) Notice in Figure 28.2 that the grey curves
lie completely above the black one, except for
the edges where they coincide.  There is nothing
in our theory that explains such a clean kind of
relationship, but it always seems to hold.
\newline
(iii) The Structure Lemma has a crisp result,
easy to check computationally for individual
cases.  However, as the reader will see, our
proof is rather tedious. We wish we had a
better proof.

\subsection{Another Diophantine Lemma}

Here we prove a copying lemma that helps
with Lemma \ref{welldefined}.  Our
result works for Farey pairs.  Let
$\Delta_1(I)$ and $\Delta_2(I)$ be the
sets defined exactly as in the Diophantine Lemma.
See \S \ref{diotheorem}.  The result we
prove here is actually more natural than our 
original result.  However, the original
result better suited our more elementary
purposes.

\begin{lemma}
\label{dio2}
Suppose that $A_1 \bowtie A_2$.
\begin{enumerate}
\item If $A_1<A_2$ let  $I=[-q_1+2,q_2-2]$.
\item If $A_1>A_2$ let  $I=[-q_2+2,q_1-2]$.
\end{enumerate}
Then $\widehat \Gamma_1$ and $\widehat \Gamma_2$ agree on $\Delta_1(I) \cup \Delta_2(I)$.
\end{lemma}

\startproof
We will consider the case when $A_1<A_2$.  The other case
has a very similar treatment.
In our proof of the Diophantine Lemma we only used the
oddness of our rationals in Lemma \ref{good}.  Once
we prove the analogue of this result in the even setting,
the rest of the proof works {\it verbatim\/}.

Recall that an integer $\mu$ is {\it good\/} if 
${\rm floor\/}(A_1 \mu)={\rm floor\/}(A_2 \mu)$.   
The analogue of Lemma \ref{good} is the statement
that an integer $\mu$ is good provided that
$\mu \in (-q_1,q_2)$.
We will give a geometric proof.  
Let $L_1$ (respectively $L_2$)  denote the line 
segment of slope $-A_1$ (respectively $-A_2$) joining the two
points whose first coordinates are $-q_1$ and $q_2$.  
If we have a counterexample to our claim then there
is a lattice point $(m,n)$ lying between $L_1$ and $L_2$.

If $m<0$, we consider the triangle $T$ with vertices
$(0,0)$ and $-V_1$ and $(m,n)$.  Here
$V_1=(q_1,-p_1)$.
The vertical distance between the left endpoints
of $L_1$ and $L_2$ is $1/q_2$.   By the base-times-height
formula for triangles, ${\rm area\/}(T)<q_1/(2q_2)<1/2$.
 But this contradicts the fact that
$1/2$ is a lower bound for the area of a lattice triangle.
If $m>0$ we consider the triangle $T$ with vertices
$(0,0)$ and $V_1$ and $(m,n)$.
The lattice point $(m,n)$ is closer to the line containing
$L_1$ than is the right endpoint of $L_2$, namely
$(q_2,-p_2)$. Hence, ${\rm area\/}(T)<{\rm area\/}(T')$,
where $T'$ is the triangle with vertices
$(0,0)$ and $V_1$ and $V_2$.
But ${\rm area\/}(T')=1/2$ because $A_1$ and $A_2$ are
Farey related.  We get the same contradiction as in the
first case.
\endproof

\subsection{Proof of Lemma \ref{welldefined}}
\label{construct}

Suppose that $A_1 \bowtie !\ A_2$.  To show that
$P\Gamma_1$ is well defined, we just have to show
$P\Gamma_2 \subset \Gamma_1$. This simultaneously
shows that $P\Gamma_1=P\Gamma_2$, because
the endpoints of these two arcs are the same by definition.
We will consider the
case when $A_1<A_2$.  The other case is similar.
In this case, we have
$A_1=(A_2)_-$.  To simplify our
notation, we write $A=A_2$.  Then $A_1=A_-$.

  By Lemma \ref{dio2},
it suffices to prove that 
\begin{equation}
P\Gamma \subset \Delta(J); \hskip 30 pt
J=[-q_-+2,q-2].
\end{equation}
We have actually already proved this, but it takes
some effort to recognize the fact.

Let $A' \leftarrow A$ denote the inferior
predecessor of $A$.  Since $q_->q_+$, we have
\begin{equation}
A'=A_- \ominus A_+.
\end{equation}
In the previous chapter, when we proved the Copy Theorem,
we established (except for a few special cases)
\begin{equation}
\label{bigint}
P\Gamma \subset \Delta'(J'); \hskip 30 pt J'=[-q'+2,q'+q_+-2].
\end{equation}
Here $\Delta'$ is defined relative to the
linear functionals $G'$ and $H'$, which are
defined relative to $A'$.
The right endpoint in Equation \ref{bigint}
comes from Lemma \ref{keycomp1}.
Now observe that
\begin{equation}
q'=q_--q_+<q_-; \hskip 30 pt
q'+q_+<(q_--q_+)+q_+=q_-<q.
\end{equation}
These calculations show that 
$J \subset J'$.
Usually $J'$ is much larger.

The region $\Delta(J)$ is computed relative to
the parameter $A$ whereas the region
$\Delta'(J')$ is computed relative to the parameter
$A'$.
The same argument as in Lemma \ref{lipschitz} shows that
\begin{equation}
\Delta(J) \subset \Delta'(J')
\end{equation}
except when $q_2<20$. The point is that the much
larger size of $J'$ compensates for any tiny
difference between the pairs $(G,G')$ 
and $(H,H')$ defining the sets.

We check the remaining few
cases by hand.  This completes the proof.
\newline
\newline
{\bf Remark:\/} By taking $A_1<A_2$ we omitted the case
when $B_1=1/(2k)$ and $A_2=1/(2k+1)$.  In this nearly
trivial case,
$E^-=(-1,1)$ and $E^+=(0,0)$.

\subsection{Proof of the Structure Lemma}
\label{struct1}

We will consider the case when $A_- \Leftarrow A$.  The other
case is similar.  Let $B$ be the odd rational such that
$A_- \bowtie !\ B$. Then $P\Gamma(A_-)=P\Gamma(B)$ by
definition.

\begin{lemma}  The Structure Lemma holds when
$1/1 \leftarrow A$. 
\end{lemma}

\startproof
 In this case
\begin{equation}
A=\frac{2k-1}{2k+1}; \hskip 30 pt
A_-=\frac{k-1}{k}; \hskip 30 pt B=\frac{2k-3}{2k-1}.
\end{equation}
Then $P\Gamma(A)$ is the line segment connecting
$(0,0)$ to $(-k,k)$ and $P\Gamma(B)$ is the line
segment connecting $(0,0)$ to $(-k+1,k-1)$.
\endproof

In all other cases, we have $A' \leftarrow A$, where
$A' \not = 1/1$.    As in Lemma \ref{indX}, let
$$\delta=\delta(A',A)={\rm floor\/}\bigg(\frac{q'}{q}\bigg).$$

\begin{lemma} 
If $\delta=1$ then
the structure Theorem holds by induction.
\end{lemma}

\startproof
If $\delta(A',A)=1$ then $d(A',A)=0$.
If $d(A',A)=0$ then $P\Gamma=P\Gamma'$ by the
Copy Theorem and the definition of pivot arcs.
At the same time, we can apply Lemma \ref{indX} to
the pair $A_m=A'$ and $A_{m+1}=A$.  Since
$\delta(A',A)=1$, we must have Case 1 or Case 3.
But we also have $A_-<A_+$. Hence, we have Case 3.
But then $A'_-=A_-.$  Hence, we can replace the
pair $(A_-,A)$ by the pair $(A'_-,A')$, and the
result follows by induction on the size of the
denominator of $A$.
\endproof

\begin{lemma}
\label{tedious0}
Suppose that $\delta=2$. Then $A'=B$.
\end{lemma}

\startproof
$B$ is characterized by the property that
$A_-$ and $B$ are Farey related, and
$$2q_- > {\rm denominator\/}(B)>q_-.$$
We will show that $A'$ has this same
property.  Note that $A'$ and $A_-$ are
Farey related.  The equations
$$2q'<q; \hskip 30 pt q=q_++q_-; \hskip 30 pt
q'=q_+-q_-$$
lead to 
$$3q_->q_+ \hskip 30 pt \Longrightarrow \hskip 30 pt
2q_->(q_+-q_-)=q'.$$
This establishes the first property for $A'$.
The fact that $\delta=2$ gives $3q'>q$.  This
leads to
$$q_+>2q_-; \hskip 30 pt \Longrightarrow \hskip 30 pt  q'=q_+-q_->q_-.$$
This is the second property for $A'$.
\endproof

\begin{lemma}
\label{tedious}
Suppose $\delta \geq 3$. Then $A' \leftarrow B$.
\end{lemma}

\startproof
There is some even rational $C$ such that 
\begin{equation}
B=A_- \oplus C
\end{equation}
The denominator of $C$ is smaller than the
denominator of $A_-$, because of the fact that
$A_- \bowtie !\ B$.
The inferior predecessor of $B$ is $A_- \ominus C$.
At the same time, 
\begin{equation}
A'=A_+ \ominus A_-
\end{equation}
 So, we are
trying to show that
$A_+ \ominus A_-=A_- \ominus C.$
This is the same as showing that
\begin{equation}
C=D:=A_- \oplus A_- \ominus A_+.
\end{equation}
Since $A_+$ and $A_-$ are Farey-related,
$D$ and $A_-$ are Farey related.
We claim that
\begin{equation}
\label{denom}
2q_--q_+={\rm denominator\/} \in (0,q_-).
\end{equation}
The upper bound comes from the fact that $q_+>q_-$.
The lower bound comes from the fact that
$q_+<2q_-$.  To see this last equation, note
$$q=q_++q_-; \hskip 30 pt
q'=q_+-q_-; \hskip 30 pt 3q'<q.$$
But $C$ is the only even rational that is Farey
related to $A_-$ and satisfies equation
\ref{denom}.  Hence $C=D$.
\endproof

As we already proved, the case
$\delta=1$ is handled by induction on the
denominator of $A$.  The case $\delta=2$
gives $$P\Gamma_-=P\Gamma'.$$  In this case,
the Structure Lemma follows from the definition
of the pivot points.

When $\delta \geq 3$, the rational $A'$ is a common
inferior predecessor of $A$ and $B$. 
Since $A_+=A' \oplus A_-$ and $A_-<A_+$, we have
$A'>A_+$.  Hence $A'>A$.   

\begin{lemma}
$A'>B$.
\end{lemma}

\startproof
Lemma \ref{tedious} gives
\begin{equation}
A'= A_- \ominus C; \hskip 10 pt
A_+=A' \oplus A_-; \hskip 10 pt
A=  A_+ \oplus A_-; \hskip 10 pt
B=  A_- \oplus C.
\end{equation}
This gives us
$$
A \ominus B=A_+ \ominus C = A' \oplus A_- \ominus C = A' \oplus A'.
$$
Hence
\begin{equation}
\label{tedious2}
A=B \oplus A' \oplus A'.
\end{equation}
Since $A_+=A' \oplus A_-$ and $A_-<A_+$, we have
$A'>A_+$.  Hence $A'>A$.
By Equation \ref{tedious2}, $A$ lies between $A'$ and $B$.
Hence $B<A<A'$.  Hence $A'>B$. In short, $A'>A$ and $A'>B$.
\endproof

Finally,
from the definition of Pivot Points, we have $E^+(A)=E^+(B)$.
This establishes Statement 1.
Statement 2 has a similar proof.

Now for Statement 3.
By Lemma \ref{swap},
$$
E^+(A)+E^-(A)=-A_-+(0,1); \hskip 30 pt
E^+(B)+E^-(B)=-B_++(0,1).
$$
Since $E^+(A)=E^+(B)$, we have
\begin{equation}
\label{depleted}
E^-(B)-E^-(A)=A_--B_+=V(C)
\end{equation}
Here $V(C)$ is as in Equation \ref{boxvectors} defined
relative to $C$.
We now have
$$
E^-(A)+V-E^-(A_-)=E^-(A)-E^-(B)+V=-V(C)+V(A)=
$$
\begin{equation}
V(A_+ \oplus A_-)-V(A_+ \ominus A_-)=
2V(A_-) \in \Z(V_-).
\end{equation}
This completes the proof of Statement 3.  Statement 4 is similar.

\subsection{The Decrement of a Pivot Arc}
\label{dpa}

Here we work out the precise relationship
between the pivot arcs in the Structure Lemma.

Let $A$ be an odd rational, and let
$A'$ be the superior predecessor of $A$.
By the Copy Theorem, $P\Gamma$ contains
at least one period of $\Gamma'$, starting
from either end.  Let $\gamma'$ be one
period of $\Gamma'$ starting from
the right endpoint of $P\Gamma$.
We define $DP\Gamma$ by the following
formula.
\begin{equation}
\label{decr}
P\Gamma = DP\Gamma * \gamma',
\end{equation}
The operation on the right
hand side of the equation is the
concatenation of arcs.
 We call $DP\Gamma$ the
{\it decrement\/} of $P\Gamma$.

The arc $DP\Gamma$ is a pivot arc
relative to a different parameter.
(See the next lemma.)
$DP\Gamma$ is obtained
from $P\Gamma$ by deleting one
period of $\Gamma'$. 
Now we give an {\it addendum\/} to the
Structure Lemma.

\begin{lemma}
\label{tocones}
If $B \Leftarrow A$ then $P\Gamma(B)=DP\Gamma(A)$, up
to translation.
\end{lemma}

\startproof
We will consider the case when $A_- \Leftarrow A$.
The other case, when $A_+ \Leftarrow A$, has
essentially the same proof.
We re-examine Lemmas \ref{tedious0}
and \ref{tedious}.  In Lemma \ref{tedious0}, we
have 
$$P\Gamma_-=P\Gamma'.$$
However, in this case, 
$\delta(A,A')=2$, and
from the definition of pivot points we see
that $P\Gamma$ is obtained from $P\Gamma'$
by concatenating a single period of $\Gamma'$.
This gives us what we want.

In Lemma \ref{tedious}, we have Equation \ref{tedious2}, which
implies
\begin{equation}
{\rm denominator\/}(A)={\rm denominator\/}(B)+2q'.
\end{equation}
But this implies that $d(A',A)=d(A',B)+1$.
Applying the Copy Theorem to both pairs, we see
that $P\Gamma$ is obtained from
$P\Gamma'$ by concatenating $d(A,B)+1$ periods
of $\Gamma'$ where $P\Gamma_-$ is obtained
from $P\Gamma'$ by contatenating $d(A',B)$
periods of $\Gamma'$.  This gives us the
desired relationship.
\endproof

\subsection{A Corollary of the Structure Lemma}
\label{even2even}

For each even rational $A_2 \in (0,1)$ that is not of the
form $1/q_2$, there is another even rational $A_1=p_1/q_1 \in (0,1)$
such that $q_1<q_2$ and $A_1 \bowtie A_2$.   In this section
we prove that the Structure Lemma above implies the
same result for $A_1$ and $A_2$.

Consider Statement 1.  Let $A_3=A_1 \oplus A_2$.  Then
$A_1 \Leftarrow A_3$ and $A_2 \bowtie !\ A_3$.
Note that $E_2^+=E_3^+$ by definition.
Also, $E_1^+=E_3^+$ by the Structure Lemma.
Hence $E_1^+=E_2^+$.  This proves Statement 1
for the pair $(A_1,A_2)$.  Statement 2 has
the same kind of proof.

Consider Statement 3.  We have
$E_2^-=E_3^-$ and
\begin{equation}
E_3^--E_1^-+V_3 \in \Z V_1.
\end{equation}
On the other hand
\begin{equation}
V_3=V_2+V_1; \hskip 30 pt \Longrightarrow \hskip 30 pt
E_3^--E_1^-+V_2 \in \Z V_1.
\end{equation}
The first equation implies the second.
But $E_3^-=E_2^-$.  This finishes the proof
of Statement 3.  Statement 4 has the same 
kind of proof.

\begin{center}
\resizebox{!}{3.5in}{\includegraphics{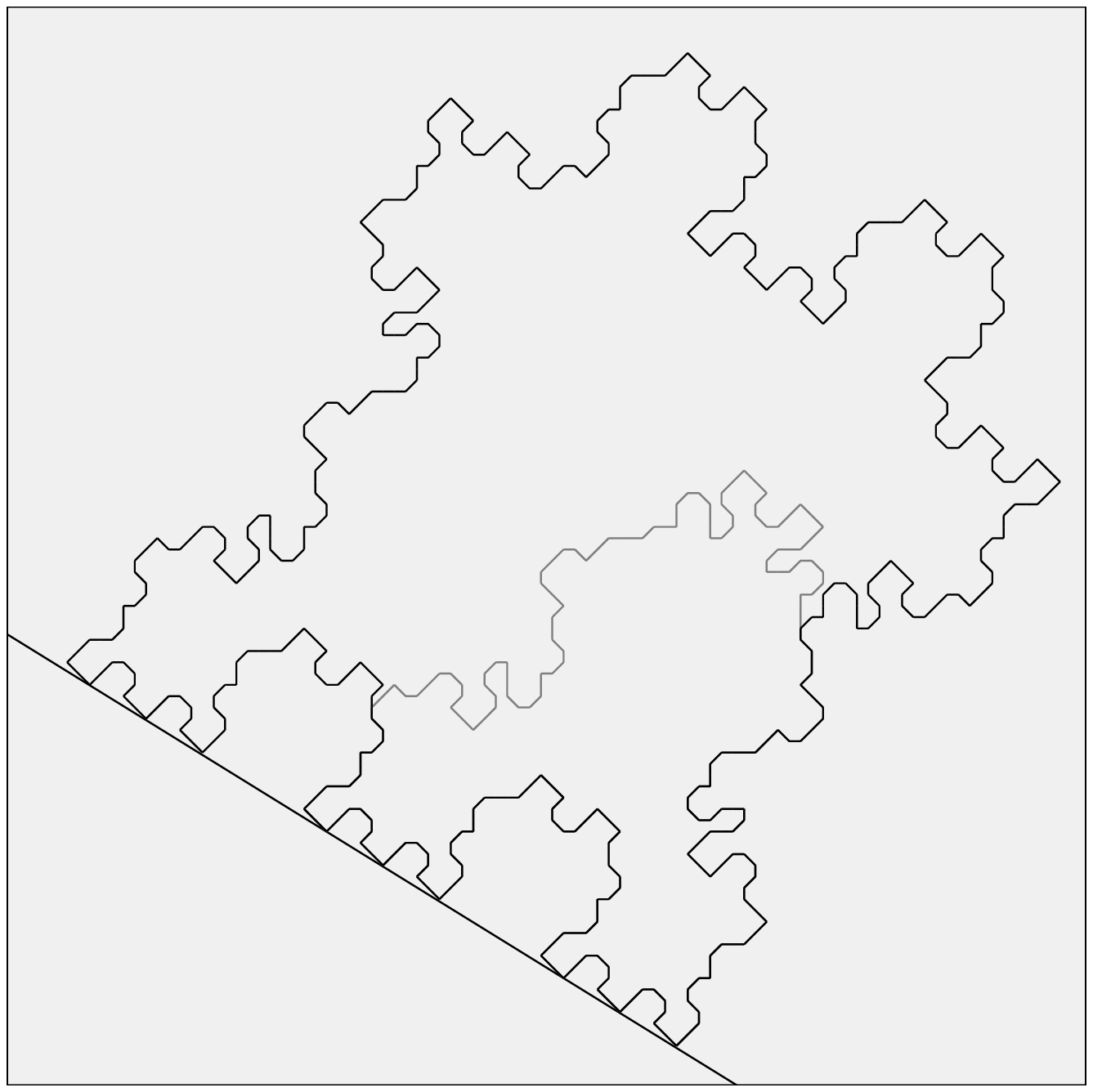}}
\newline
{\bf Figure 28.3:\/} $\Gamma(34/55)$ in black and $\Gamma(21/34)$ in grey.
\end{center} 

\subsection{An Even Version of the Copy Theorem}

Let $A_2$ be an even rational. We write $A_2=A_0 \oplus A_1$
where $A_0$ is odd and $A_1$ is even.

\begin{lemma}
\label{copyth}
$P\Gamma_2 \subset \Gamma_0$.
\end{lemma}

\startproof
We have $P\Gamma_2=P\Gamma(A_3)$, where
$A_3$ is the odd rational such that
$A_2 \bowtie !\ A_3$.  Since $A_1 \bowtie A_2$
and both $A_1$ and $A_2$ are even, we have
$A_3=A_1 \oplus A_2$.  At the same time, we have
$A_0=A_2 \ominus A_1$.  Hence $A_0 \leftarrow A_3$.
But now we can apply the Copy Theorem to
the pair $(A_0,A_3)$ to conclude that
$P\Gamma_3 \subset \Gamma_0$.   But
$P\Gamma_3=P\Gamma_2$.
\endproof

\newpage

\section{Proof of the Pivot Theorem}
\label{pivotproof}

\subsection{An Exceptional Case}
\label{easycase}

We first prove the Pivot Theorem for the parameter
$A=1/q$, with $q \geq 2$ being either even or odd.
This case does not fit the general pattern of proof.

Let $\Gamma$ be the arithmetic graph associated to
$A=1/q$, and
let $P\Gamma$ denote the pivot arc.  In
all cases, $P\Gamma$ contains the vertices
$(0,0)$ and $(-1,1)$.  These vertices
correspond to the two points
\begin{equation}
\bigg(\frac{1}{q},-1\bigg); \hskip 30 pt
\bigg(2-\frac{1}{q},-1\bigg)
\end{equation}
These two points are the midpoints of the
special intervals
\begin{equation}
I_1=\bigg(0,\frac{2}{q}\bigg) \times \{-1\} \hskip 40 pt
I_2=\bigg(2-\frac{2}{q},2\bigg) \times \{-1\}.
\end{equation}
These intervals appear at either end of 
\begin{equation}
I=[0,2] \times \{-1\}.
\end{equation}
  When we
say {\it special interval\/}, we refer to the
discussion in \S \ref{definebasic}.  These
special intervals are permuted by the
outer billiards dynamics.

For any $A<1/2$, our phase portrait
in Figure 2.4 shows that 
the interval 
\begin{equation}
I_3=(2A,2-2A) \times \{-1\}
\end{equation}  returns to itself
under one iterate of $\Psi$.
When $A=1/q$, we have
\begin{equation}
I-I_3=I_1 \cup I_2.
\end{equation}
But then the orbit of $I_1$ only intersects
$I$ in $I_1 \cup I_2$.
In terms 
of the arithmetic graph, this is to say that
the onlyl low vertices on $\Gamma$ are equivalent to
$(0,0)$ and $(-1,1)$ modulo translation by $V=(-q,1)$.
This establishes the Pivot Theorem for $A=1/q$.

\subsection{Discussion of the Proof}
\label{discussion}
\label{discuss}

Now we consider the general case of the Pivot Theorem.
We will consider the odd case until the last section
of the chapter.  At the end, we will explain the
minor differences in the even case.
For any odd rational $A_2 \not = 1/q_2$, we have
$A_1 \Leftarrow A_2$, where $A_1 \in (0,1)$ is an
even rational.  See \S \ref{predecessor}.  By induction, we can assume
that the Pivot Theorem is true for $A_1$.  
Our discussion refers to Figure 29.1.

Lemma \ref{dio2} gives a large region $\Delta$
where $\widehat \Gamma_1$ and $\widehat \Gamma_2$ agree.
$\Delta$ is white
in Figure 29.1.
The arc $\Gamma_2$ is drawn in black and the
relevant components of $\widehat \Gamma_1$ are drawn in grey.
The black dots are the endpoints of the black arc.
This black arc is ``the hump'' that we discussed
in connection with the Structure
Lemma in the previous chapter.  Indeed,
Figure 29.1 is a cartoon of Figure 28.2, with
other relevant details added. 

\begin{center}
\resizebox{!}{2.1in}{\includegraphics{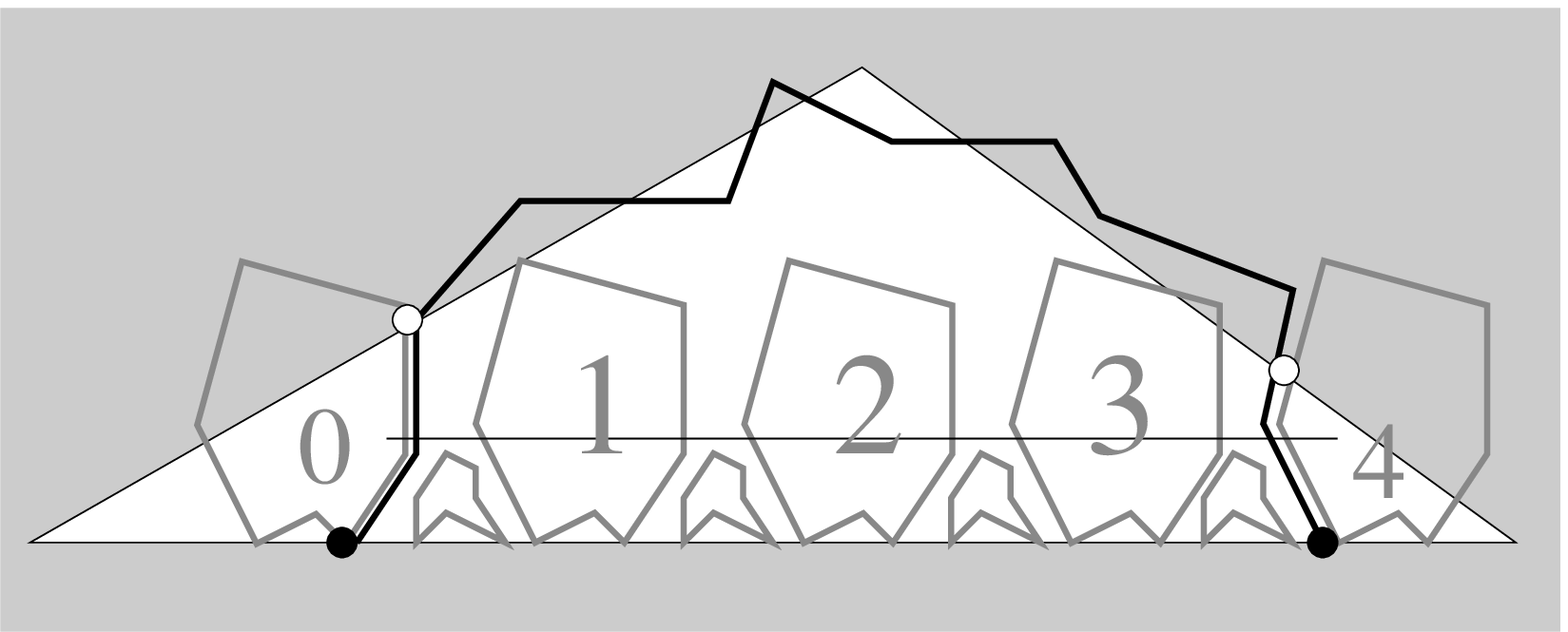}}
\newline
{\bf Figure 29.1:\/} Cartoon view of the proof
\end{center} 

We want to see that the black arc has no low
vertices except for its endpoints.   By the
structure Lemma, the endpoints of the black arc are
also endpoints of the pivot arcs of $C_0$ and $C_4$.
By induction, the only low vertices of $C_0$ and $C_4$
are contained on the pivot arcs.  These pivot
arcs are on the other sides of the endpoints
we are considering. 
 Hence there are
no low vertices on the black arc as long as it
coincides with either $C_0$ or $C_4$.

There is one subtle point to our argument.
When we refer to {\it low vertices\/} of
the black arc, the vertices are low with
respect to the parameter $A_2$.  However,
when we refer to low vertices of $C_0$ and
$C_4$, the vertices are low with respect to
$A_1$.  Will discuss this subtle point in
the next section.  What saves us is that
the two notions of {\it low\/} coincide,
due to the way in which $A_1$ approximates
$A_2$.

So, either end of our black arc starts out well:
It rises away from the baseline.  What could
go wrong?  One of the ends could dip back down
into $\Delta$ and (at the boundary) merge with
a component of $\widehat \Gamma_1$.  
In other words, some component of $\widehat \Gamma_1$
would have to stick out of $\Delta$.

There are two kinds of components of $\widehat \Gamma_1$
we need to consider.  First, there are the
major components. Recall from \S \ref{minor}
that these are the translates
of $\Gamma_1$ by vectors in $\Z V_1$.  We have labelled
these components $C_1,C_2,C_3$.  Second, there are
the minor components -- the low components
that are not major. 

The components that seem to give us the most trouble
are $C_1$ and $C_3$. These come the closest to sticking
out of $\Delta$.  In fact, we will not be able to
show that these components are contained in
$\Delta$, even though experimentally it is always
the case.  However, Lemma \ref{parity} comes to
the rescue.  The low vertices on these components
have odd parity, and the low vertices on the
black arc (a subset of $\widehat \Gamma_2$) have
even parity.   Hence, the black arc cannot
merge with $C_1$ and $C_3$.  The parity argument
steps in where our geometry fails.

The remaining major components are much farther
inside $\Delta$, and do not pose a threat.
We will give an explicit estimate to show that the
other major components are contained entirely
inside $\Delta$.  In this case, we are referring
to $C_2$, though in general there could be many
such components.

This leaves the minor components.  The
Barrier Theorem from \S \ref{barrier} 
handles these.  
The black horizontal
line in Figure 29.1 represents the barrier, which
no minor component can cross.  Equipped with the
Barrier theorem, we will be able to show that
all minor components lie in $\Delta$.

This takes care of all the potential problems.
Since the black arc can't merge with any of the
grey components, it just skips over everything
and has no low vertices, except for its
endpoints.  The rest of
the chapter is devoted to making this
cartoon description precise.

As with the proof of the Decomposition Theorem,
the estimates we make are true by a wide margin
when $A_1$ is large.  However, when $A_1$ is small,
the estimates are close and we need to deal with
the situation in a case-by-case way.  We hope that
this fooling around with small cases doesn't obscure
the basic ideas in the proof.

We close this section by remarking on a phenomenon that
we cannot establish.  Experimentally, we see that
$\widehat \Gamma_2$ copies all the low components of
$\widehat \Gamma_1$ beneath ``the hump''.  The interested
reader can see this in action using Billiard King.

\subsection{Confining the Arc}
\label{subtle}

We continue with the notation from the previous section.
For ease of exposition, we assume that $A_1<A_2$.  The
other case is similar. For ease of notation, we
set $A=A_2$.  Until the end of this section,
we only consider $A$.  We write one period of
$\Gamma$ as $P\Gamma \cup \gamma$.  Here $P\Gamma$ is
the pivot arc, and $\gamma$ is the black arc
considered in the previous section.

Let $W$ be the vector from Equation \ref{boxvectors}.
Let $S$ be the infinite strip whose left edge is
the line through $(0,0)$ parallel to $W$ and whose
right edge is the line through $V_+$ and parallel to
$W$.  Here $V_+=(q_+,-p_+)$, and $p_+/q_+$ is as
in Equation \ref{induct0}.   

\begin{lemma}
\label{nocross}
$\gamma$ does not cross the lines bounding $S$.
\end{lemma}

\startproof
The lines of $S$ are precisely the extensions of the
sides of $R_2$, the larger of the two parallelograms
from the Decomposition Theorem.  We know that
$\Gamma$ crosses these lines only once.  The left
crossing point is $(0,0) \in P\Gamma$. Hence,
the left crossing point is not a vertex of $\gamma$.

The right crossing point is $x=V_++(0,1)$.  
Let $\J$ be the symmetry from Lemma \ref{pseudolinear}.
Let $\iota(v)=V_+-v$.  Consider the map
$\phi=\iota \circ \J$.  On low vertices $v$,
we have 
\begin{equation}
\phi(v)=V_+-v+(0,1).
\end{equation}
Hence $\phi(0,0)=x$.  By Lemma \ref{swap},
$\phi$ swaps the endpoint of $\gamma$.
Moreover, $\phi$ permutes the set of low
vertices of $\gamma$. 

Since $V_+$ lies beneath the baseline, $x$ is
a low vertex.  If $x$ is a low vertex of $\gamma$
then $\phi(x)=(0,0)$ is a low vertex of $\gamma$.
This is a contradiction.  Hence $x$ is not a
low vertex of $\gamma$.
\endproof

Now we can clear up the subtlety discussed in the previous
section.  We set $S_2=S$, the strip defined relative to
the odd rational $A_2$.

\begin{lemma}
\label{lowlow}
A vertex in $S_2$ is low with respect to $A_1$ iff
it is low with respect to $A_2$.  Hence, a vertex
of $\gamma$ is low with respect to $A_1$ iff it
is low with respect to $A_2$.
\end{lemma}

\startproof
  Let $L_j$ denote the baseline with respect to $A_j$.
The conclusion of this lemma is equivalent to the
statement that
there is no lattice point between $L_1 \cap S$
and $L_2 \cap S$.  This is a consequence of our proof of
Lemma \ref{dio2}. 
\endproof

\subsection{A Topological Property of Pivot Arcs}
\label{toparc}

Let $A$ be a rational kite parameter, either even or odd.
Let $P\Gamma$ denote the pivot arc of $\Gamma=\Gamma(A)$.
The two endpoints of $P\Gamma$ are low vertices.  Here
we prove a basic structural result about $P\Gamma$.

\begin{lemma}
\label{nospoil1}
$P\Gamma$ contains no low vertex to the right of its right endpoint.
Likewise $P\Gamma$ contains no low vertex to the left of its left
endpoint.
\end{lemma}

\startproof
We will prove the first statement.  The second statement
has the same proof.  We give an argument like the one in
the proof of Lemma \ref{parity}.
Note that $\Gamma$ right-travels at $(0,0)$.  Hence
$P\Gamma$ right-travels at its right endpoint $\rho$.
Suppose that $P\Gamma$ contains a low vertex
$\sigma$ to the right of $\rho$.  Then some arc $\beta$
of $P\Gamma$ connects $\rho$ to $\sigma$.  Since
$\Gamma$ right travels at $\rho$, some arc $\gamma$
of $\Gamma-P\Gamma$ enters into the region
between $\rho$ and $\sigma$ and beneath
$\beta$.  But $\gamma$ cannot escape from this region,
by the Embedding Theorem.   The point here is that
$\gamma$ cannot squeeze
beneath a low vertex, because the only vertices below
a low vertex are also below the baseline.
Figure 29.2 shows the situation.

\begin{center}
\resizebox{!}{2in}{\includegraphics{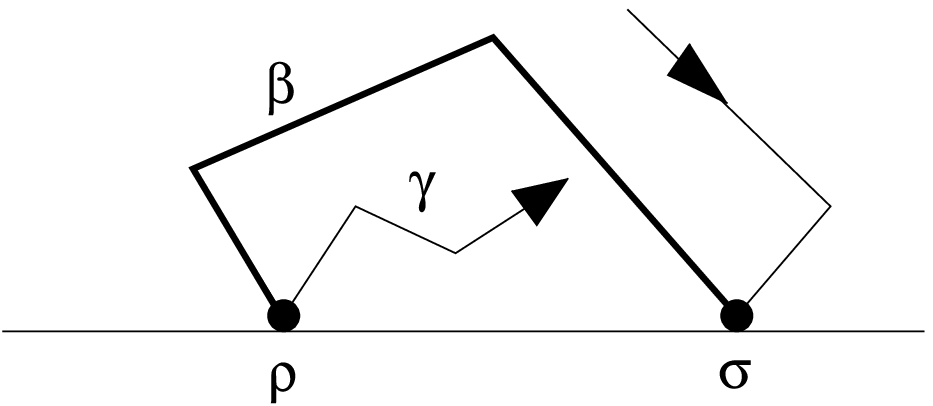}}
\newline
{\bf Figure 29.2:\/} $P\Gamma$ creates a pocket.
\end{center}

In the odd case we have an immediate contradiction.
In the even case, we see that there must be a loop
containing both $\rho$ and $\sigma$.  This loop
must be a closed polygon, and a subset of
$P\Gamma$.  Since $P\Gamma$ is also a closed (and
embedded) polygon, we our loop must
equal $P\Gamma$.  But by definition, $P\Gamma$ lies
below $\Gamma-P\Gamma$.  From Figure 29.3, we see
that $P\Gamma$ (which contains $\beta$)
 in fact lies above $\Gamma-P\Gamma$ (which
contains $\gamma$.)
This is a contradiction.
\endproof

\subsection{Corollaries of the Barrier Theorem}

Here we derive a few corollaries of the Barrier
Theorem.  See \S \ref{barrier} for the statement.

\begin{corollary}
\label{leftedge}
A minor component of $\widehat \Gamma$ cannot
cross the line through $(0,0)$ that is parallel to $W$.
\end{corollary}

\startproof
Our line is one of the lines in the Hexagrid Theorem.
By the Hexagrid Theorem, only $\Gamma$ crosses
this line beneath the barrier, and the crossing takes place
at $(0,0)$.
\endproof

We are trying to construct a parallelogram that bounds the
minor components.  The baseline contains the bottom edge.
The barrier contains the top edeg.  The line in 
Corollary \ref{leftedge} contains the left edge.  Now we
supply the right edge. Actually, there are many choices
for this right edge.

\begin{lemma}
\label{rightedge}
Let $V_+=(q_+,-p_+)$. 
A minor component of $\widehat \Gamma$ cannot cross
the line through $(0,0)$ that is parallel to $V_++kV$
for any $k \in \Z$.
\end{lemma}

\startproof
Since $\widehat \Gamma$ is invariant under translation by $V$, it
suffices to prove this result for $k=0$.  Let $L$ be the
line through $V_+$ parallel to $W$.  Our result really
follows from the bilateral symmetry discussed in
\S \ref{nbs}.  Here we work out the details, using the
rotational symmetry instead.  (We made more precise statements
about the rotational symmetry.)

Let $\Lambda$ be the barrier. 
Consider the symmetry $\iota$ defined in \S \ref{symmetry}.
The two lines $\Lambda$ and $\iota(\Lambda)$ are
equally spaced above and below the baseline up to an
error of at most $1/q$.
Suppose that
some minor component $\beta$ crosses our line $L$.  Then the
component $\iota(\beta)$ crosses the line $\iota(L)$.  But
$\iota(L)$ is the line from Lemma \ref{leftedge}.   Inspecting
the hexagrid, we see that $\iota(L)$ contains the door $(0,0)$,
but no other door between the baseline and $\iota(\Lambda)$.
Indeed, the doors above and
below the baseline are just about evenly spaced away
from $(0,0)$ going in either direction.  See Figure 3.2,
a representative figure.  (In this figure, we are talking
about the long axis of the kite, and $(0,0)$ is the
bottom tip of the kite.)

The component $\gamma'$ of
$\widehat \Gamma$ that crosses $\iota(L)$ near
$(0,0)$ has the same size as $\Gamma$. Hence,
this component crosses through $\iota(\Lambda)$.
Hence $\iota(\gamma')$ is a major component.
Hence $\beta \not = \iota(\gamma')$.  Hence
$\iota(\beta) \not = \gamma$.  Hence
$\iota(\beta)$ does not cross $\iota(L)$.
Hence $\beta$ does not cross $L$.
\endproof

\subsection{Juggling Two Parameters}
\label{juggle}

In our proof of the Pivot Theorem, we have two parameters
$A_1 \Leftarrow A_2$.  As above, we focus our attention
on the case when $A_1<A_2$.  The other case has a completely
parallel discussion.  See \S \ref{othercase} for a
brief discussion of the other case.

Lemma \ref{rightedge} applies to vectors
defined in terms of $A_1$, but we would like to
apply it to a special line defined partly 
in terms of $A_2$.  Let $(V_j)_+$ be as
in \S \ref{subtle}.
Then Lemma \ref{rightedge} applies to the vectors
of the form $(V_1)_+ + kV_1$.  However, we are
also interested in the vector $(V_2)_+$.

\begin{lemma}
\label{affine3}
Suppose that $A_1<A_2$.  Then,
there is some integer $k$ such that
$(V_2)_+=(V_1)_++kV_1$.
\end{lemma}

\startproof
We set $A=A_2$.  Then $A_-=A_1$.
Let $A_{-+}$ denote the parameter that relates to $A_-$
in the same way that $A_+$ relates to $A$.  That is
$A_{-+}>A_-$ are Farey related and $A_{-+}$
has smaller denominator than $A_-$.
We want to prove that
$V_+=V_{-+}+kV_-$ for some $k$.
The rationals $A_{-+}$ and $A_-$ are Farey-related.
Therefore, so are the parameters
\begin{equation}
\label{iterateFarey}
A_-; \hskip 40 pt
A_{-+} \oplus A_- \oplus ... \oplus A_-.
\end{equation}
Here we are doing Farey addition.   Conversely, if any rational
$A'$ is Farey related to $A_-$, and has bigger denominator,
then the Farey difference $A' \ominus A_-$ is also Farey
related to $A_-$.  Thus, the rationals in
Equation \ref{iterateFarey} account for all the rationals
$A'$ with the properties just mentioned.  But $A_+$ is one
such rational.  Hence $A_+$ has the form given
in Equation \ref{iterateFarey}.  This does it.
\endproof

Let $R$ denote the parallelogram defined by the following
lines
\begin{itemize}
\item The baseline relative to $A_1$.
\item The barrier for $A_1$.
\item The line parallel to $W_1$ through $(0,0)$.
\item The line parallel to $W_1$ through $(V_2)_+$.
\end{itemize}
Then any minor component with one vertex in $R$ stays
completely in $R$.   This is a consequence of the
Barrier Theorem, its corollaries, and the lemma
in this section.  Modulo a tiny adjustment in the
slopes, the left and right edges of $R$ are
contained in the left and right edges of the strip $S$
considered in \S \ref{subtle}.

\subsection{A Bound for Minor Components}
\label{boundminor}

Let $A_1 \Leftarrow A_2$ be as above.  Again, we
assume that $A_1<A_2$ for ease of exposition.
Define
\begin{equation}
\label{deltadefined1}
\Delta=\Delta_1(I)\cup \Delta_2(I); \hskip 30 pt
I=[-q_1+2,q_2-2].
\end{equation}
Here $\Delta$ is as in Lemma \ref{dio2}.
Let $R$ be the parallelogram discussed in the previous section.

\begin{lemma}
Let $\beta \subset \widehat \Gamma_1$ be any component
that is contained in $R$.  Then $\beta \subset \widehat \Gamma_2$.
\end{lemma}

\startproof
Our proof follows the same strategy as in the Decomposition
Theorem. We will work with the functionals
$G_1$ and $H_1$ defined relative to $A_1$.

Essentially, we want to show that $R \subset \Delta$
and then apply Lemma \ref{dio2}.  However, to avoid a messy
calculation, we invoke the Adjacent Mismatch Principle, and
replace $R$ by the extremely nearby parallelogram
$\widetilde R$ with vertices
\begin{equation}
(0,0); \hskip 30 pt \lambda W_1; \hskip 30 pt
(V_2)_+; \hskip 30 pt (V_2)_++\lambda W_1.
\end{equation}
The constant $\lambda$ has the following definition.
The top left vertex of $R$ lies on the line through
$(0,0)$ and parallel
to $W_1$, as we discussed above.
Hence this vertex has the form $\lambda W_1$.
We have
\begin{equation}
M_1(W_1)=M_1\bigg(0,\frac{p_1+q_1}{2}\bigg)=p_1+q_1; \hskip 30pt
M_1(\lambda W_1)=p_1'+q_1'<p_1+q_1.
\end{equation}
Here $A_1'$ is the rational that appears in the Barrier Theorem.
The point here is that the barrier contains the point
$(0,(p_1'+q_1')/2)$.  In particular, $\lambda<1$.

Let $u$ and $w$ be the top left and top right vertices of
$R$.  As usual, it suffices to show that the quantities
\begin{equation}
\label{minorgoal}
G_1(u)-(-q_1+2); \hskip 30 pt
(q_2-2)-H_1(w)
\end{equation}
are both positive.  In fact these quantities are equal.
As in Equation \ref{affine1}, we compute
\begin{equation}
G_1(u)-(-q_1+2)=q_1-\lambda \frac{q_1^2}{p_1+q_1} -2
\end{equation}
We will do the second calculation by symmetry.
By Lemma \ref{affine3}, we have
$$
(V_2)_++V_1=(V_2)_++(V_2)_-=V_2.
$$ 
Hence
$$V_2-w=V_1-\lambda W_1.$$  
Hence
$$
(q_2-2)-H_1(w)=
-2+H_1(V_2-w)=$$
\begin{equation}
-2+H_1(V_1-\lambda W_1)=
q_1-\lambda \frac{q_1^2}{p_1+q_1}-2
\end{equation} 
We get exactly the same answer in both cases.
This is a reflection of an 
underlying affine
symmetry, as we remarked after
Equations \ref{affinecoincidence1} and
\ref{affinecoincidence2}.

Since $\lambda \leq 1$, the quantities in
Equation \ref{minorgoal} are non-negative as long as
$p_1 \geq 3$ and $q_1 \geq 7$.  This is exactly
the same estimate as in Lemma \ref{strong4}.
When $p_1=2$ we see that
$$p_1'=1; \hskip 30 pt q_1'=\frac{q_1-1}{2}.$$
Thus $\lambda \approx 1/2$, and we get a massive
savings.   When $p_1 \geq 2$ and $q_1 \leq 7$ we
check the cases by hand, using the same trick
as in \S \ref{tricks}.

It remains to consider the case $p_1=1$.  In this case
$\widehat \Gamma_1$ has no minor components, as we
saw in \S \ref{easycase}.
\endproof

\subsection{A bound for Major Components}
\label{boundmajor}

We keep the parameters $A_1 \Leftarrow A_2$ as
above, with $A_1<A_2$.
We have already defined the pivot points
of $\Gamma_1$.  We define the pivot
points of the translates $C_k=\Gamma_1+kV_1$
in the obvious way, by translation.  

By the Structure Lemma, there is some
component $C_k$ whose left pivot point
is $E_2^-+V_2$, the right endpoint of
the ``hump'' discussed in \S \ref{discussion}.
The components $C_0,...,C_k$ are
exactly as in \S \ref{discussion}. 
By Lemma \ref{parity}, the index $k$ is
even.  More generally, $C_j$ contains
low vertices of even parity if and only
if $j$ is even.   

As in \S \ref{discussion}
we are interested in bounding the
components $C_2,...,C_{k-2}$.  Actually,
we only care about the even components,
but our bound works equally well for
the odd components between $C_2$ and
$C_{k-2}$.  If $k=2$ one can just ignore
the construction in this section.

By the Hexagrid Theorem, $C_0$
is contained in the parallelogram
$R_0$ with vertices
\begin{equation}
-V_1; \hskip 30 pt -V_1+2W_1; \hskip 30 pt
V_1+2W_1; \hskip 30 pt V_1.
\end{equation}
This means that $C_j$ is contained in 
translated parallelogram 
\begin{equation}
R_j=R_0+jV_1
\end{equation}
We choose $j \in \{2,...,k-2\}$.

Here we describe some features of $R_j$, as
well as a recipe for symmerizing it.
\begin{enumerate}
\item
The bottom edge of $R_j$ is contained in the line
through $(0,0)$ and parallel to $V_1$--i.e. the
baseline, as usual.
\item
The top edge of $R_j$ is contained
in the line through $2W_1$ and parallel to $V_1$.
These lines are independent of $j$.
\item
The left edge of $R_j$ is parallel to, and to the
right of, the line $\Lambda$ parallel to $W_1$ and containing
$V_1$.  When $j=2$ the left edge of $R_j$ is
contained in $\Lambda$.
\item The same symmetry argument as in Lemma \ref{rightedge}
shows that $C_2$ lies to the left of the line through
the point $(V_2)_+-V_1$ and parallel to $W_1$.  Referring
to the symmetry $\iota$ in Lemma \ref{rightedge}, this
is the line $\iota(\Lambda)$.  
In brief, if $C_j$ crosses $\iota(\Lambda)$, then
$\iota(C_j)$ crosses $\Lambda$, and this contradicts
the Hexagrid Theorem, applied below the baseline.
\end{enumerate}

Let $R$ be the parallelogram defined by the $4$ lines
above.  By construction $C_j \subset R$
for $j \in \{2,....,k-2\}$.   

\begin{lemma}
Let $\beta \subset \Gamma_1$ be any component of
$\widehat \Gamma_1$ that is contained in $R$. Then
$\beta \subset \widehat \Gamma_2$.
\end{lemma}

\startproof
The proof is exactly the same.  Let $u$ and $w$ denote the
top left and top right vertices of $R$.  We get the same
symmetry as in the previous bound, and so we just
have to compute $G_1(u) \geq -q_1+2$.
We compute
\begin{equation}
G_1(u)-(-q_2+2)=2q_1-\frac{2q_1^2}{p_1+q_1}-2.
\end{equation}
This time we always get a positive number, though
in small cases it is pretty close.
\endproof

\subsection{Even implies Odd}
\label{evenimpliesodd}

Let $P(A)$ be the statement that the Pivot Theorem is
true for $A$.

\begin{lemma}
\label{evenodd}
Let $A_1 \Leftarrow A_2$.  Then $P(A_1)$ implies $P(A_2)$.
\end{lemma}

Our proof follows the format of the discussion in
\S \ref{discuss}.  As in \S \ref{subtle}, we define the
{\it complementary arc\/} $\gamma_2 \subset \Gamma$ to 
be the arc to the right of $P\Gamma_2$ such that
$P\Gamma_2 \cup \gamma_2$ is one period of $\Gamma_2$.
The endpoints of $\gamma_2$ are
\begin{equation}
E^+_2; \hskip 30 pt E_2^- +V_2.
\end{equation}
This is the ``hump'' we discussed in \S \ref{discuss}.

We say that a
{\it spoiler\/} is a low vertex of $\gamma_2$
that is not an endpoint of $\gamma_2$.
The Pivot Theorem is equivalent to the statement
that there are no spoilers.

Let $L(\gamma_2)$ denote the left endpoint
of $\gamma_2$. Likewise, let $R(\gamma_2)$
denote the right endpoint of $\gamma_2$.

\begin{lemma}
Any spoiler lies between $L(\gamma_2)$ and
$R(\gamma_2)$.
\end{lemma}

\startproof
We will show that any spoiler lies to the right of
$L(\gamma_2)$.  The statement that any spoiler
lies to the left of $R(\gamma_2)$ is similar.
By Lemma \ref{nocross}, all spoilers lie in the
strip $S_2$.
But $P\Gamma_2$ crosses the left boundary of $S_2$.
Any low vertices in $S_2$ to the left of $L(\gamma_2)$
either lie on $P\Gamma_2$ or beneath it.  By the
Embedding Theorem, $\gamma$ cannot contain these
vertices.
\endproof

\begin{lemma}
\label{nospoil4}
$\Delta$ contains all the spoilers.
\end{lemma}

\startproof
We will work with the linear functionals
$G_2$ and $H_2$ defined relative to $A_2$.
Thus, we are really showing that the smaller
set $\Delta_2(I)$ contains all the spoilers.

Let $v=(m,n)$ be a spoiler.
If suffices to prove that
$G_2(v) \geq -q_0+2$ and
$H_2(v) \leq q_2-2.$
We have $m \geq 1$.  Since $v$ is a low vertex,
we have $n \leq 0$.  
We compute that $\partial_y G_2<0$.  Hence
$$G_2(v) \geq G_2(m,0)=m \frac{1-A_2}{1+A_2}>0 \geq -q_1+2.$$
This takes care of $G_2$.

Let $w=v-V_2=(r,s)$.
By Lemma \ref{Hcalc}, it suffices to show
$H_2(w) \leq -2$.  We compute $\partial_y H>0$.
Since $w$ lies at most one vertical unit above
the line of slope $-A_2$ through the origin, we have
\begin{equation}
H_2(w) \leq H_2(w'); \hskip 30 pt w'=(r,-A_2r+1).
\end{equation}
We compute
\begin{equation}
\label{elim}
H_2(w')=r+\frac{2(1-A_2)}{(1+A_2)^2}<r+2.
\end{equation}
This shows that $H(w)<-2$ as long as $r \leq -4$.
By Lemma \ref{parity}, we have $r+s$ even.
We just have to rule out $(-2,2)$ and $(-3,1)$ 
as spoilers.

If $A_2<1/2$ then $(-2,2)$ is not a low vertex.
If $A_2>1/2$ then 
$$\frac{2k-1}{2k+1} \leftarrow ... \leftarrow A_2$$
for some $k \geq 2$.  In this case,
$E_2^-$ has first coordinate $\leq -2$.  But then
$r \leq -3$.   This
rules out $(-2,2)$.  
We compute that $H_2(-3,1)<-2$ when $A \geq 1/9$.
When $A<1/9$, we use the phase portrait in
\S \ref{low vertex} to check that $\widehat \Gamma_2$
is trivial at $(-3,1)$. This rules out $(-3,1)$.
\endproof

Let $v$ be a spoiler.  By the previous result, there is
some component $\beta$ of $\widehat \Gamma_1$ that
has $v$ as a vertex.

\begin{lemma}
\label{nospoil3}
$\beta$ is not a subset of $\widehat \Gamma_2$.
\end{lemma}

\startproof
Let's start at $v$ and trace $\gamma_2$ in some direction.
If the conclusion of this lemma is false, we 
remain simultanously on $\gamma_2$ and $\beta$ until we loop
around and  return
to $v_2$ -- because $\beta$ is a closed polygon. This
contradicts the fact that $\gamma_2$ never visits the
same vertex twice.
\endproof

Here is the end of the argument.
$\beta$ cannot be a minor component,
given the bound in \S \ref{boundminor}.
Next, $\beta \not \in \{C_2,...,C_{k-2}\}$
given the bounds in \S \ref{boundmajor}.
Next, $\beta \not \in \{C_1,C_{k-1}\}$
by Lemma \ref{parity}.
Next, $\beta \not = C_0$:  By induction,
all the low vertices
of $C_0$ lie on $PC_0$.  By Lemma \ref{nospoil1}
these low vertices all lie to the left of the spoiler.
Likewise $\beta \not = C_k$.
We have exhausted all the possibilities.
$\beta$ cannot exist.  Hence there is
no spoiler.  Hence $P(A_2)$ holds.

\subsection{A Decomposition in the Even Case}

In this section we revisit the construction
in \S \ref{subtle}, but for even parameters.
Now $A_1$ and $A_2$ are both even parameters,
with $A_1 \bowtie A_2$. 
We set $A=A_2$ and just consider objects relative
to $A$.  We define the strip $S$ exactly as in
\S \ref{subtle}.  This time we define
\begin{equation}
\gamma = \Big(\Gamma \cup (\Gamma+V)\Big) \cap S.
\end{equation}
\begin{center}
\resizebox{!}{2.6in}{\includegraphics{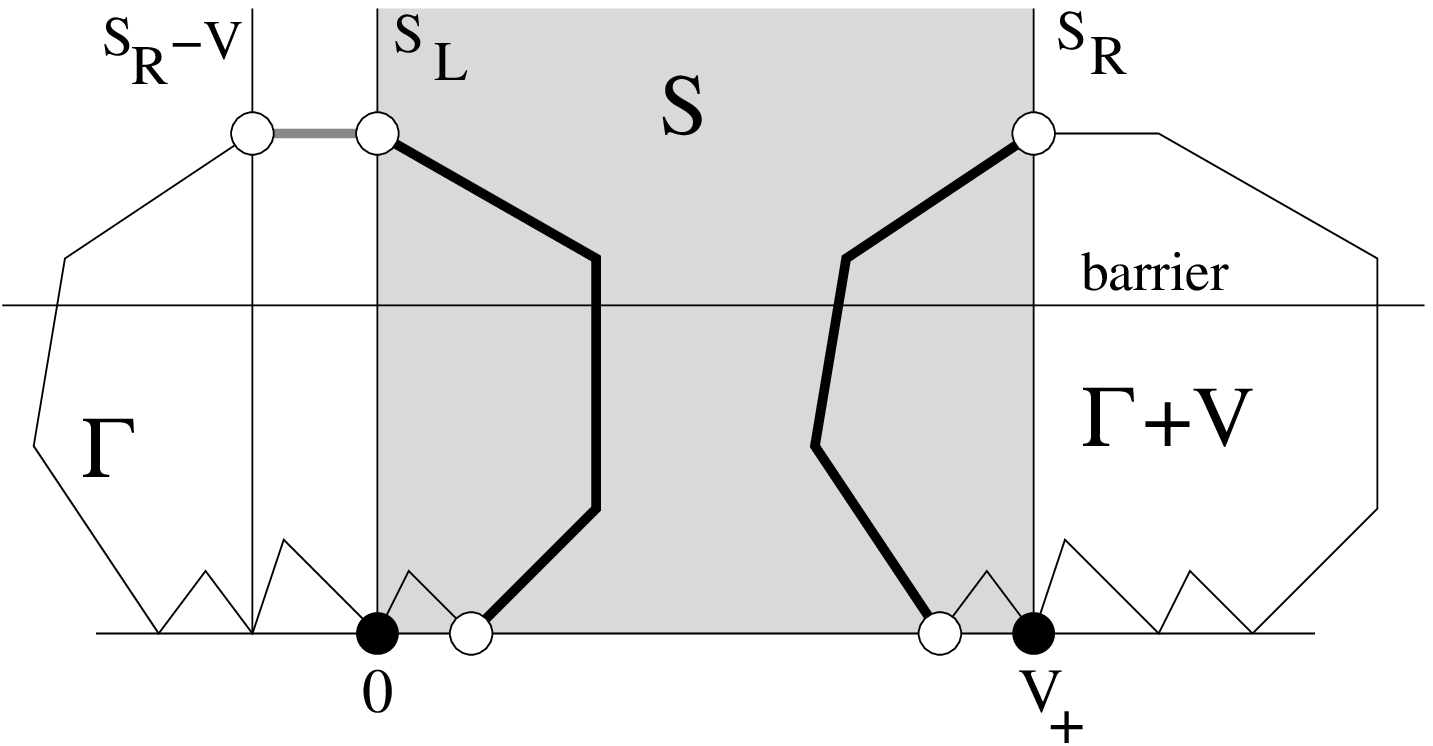}}
\newline
{\bf Figure 29.3:\/} The even version of $\gamma$.
\end{center} 

\begin{lemma}
\label{even nocross}
$\gamma$ consists of two connected arcs.  Any low vertex
of $\Gamma-P\Gamma$ is translation equivalent to a low vertex of
$\gamma$.
\end{lemma}

\startproof
By the Hexagrid Theorem $\Gamma$ only crosses $S_L$
once.  The door on $S_L$ lies above the
barrier line.  Hence, the crossing occurs above
the barrier line.
  Likewise, $\iota(\Gamma+V)$
only crosses $S_L$ once.  The relevant door lies
below the image of the barrier line under $\iota$.
  Here $\iota$ is as in
the proof of Lemma \ref{nocross}.  But then
$\Gamma+V$ only crosses $S_R$ once, and the crossing
occurs above the barrier line.  Hence $\gamma$
consists of $2$ connected arcs.

The line $S_R-V$ is parallel to $S_L$ and lies to the
left of $S_L$.
By symmetry, $\Gamma$ only crosses $S_R-V$ once, and
the crossing takes place above the barrier line.
By the Barrier Theorem, the grey arc of
$\Gamma$ between $S_L$ and $S_R-V$ lies above
the barrier line and hence has no low vertices.
Finally, any vertex of $\Gamma-P\Gamma$ not translation equivalent to
a vertex of $\gamma$ lies on the grey arc of
$\Gamma$ between $S_L$ and $S_R-V$. 
\endproof

\subsection{Even implies Even}
\label{evenimplieseven}

Let $A_1 \bowtie A_2$ be a pair of even rationals
as in \S \ref{even2even}.
This pair exists as long as $A_2 \not = 1/q_2$.
Referring to the terminology in Lemma \ref{evenodd}, we prove
the following result
in this section.

\begin{lemma}
\label{eveneven}
Let $A_1 \Leftarrow A_2$.  Then $P(A_1)$ implies $P(A_2)$.
\end{lemma}

We have already taken care of the base case of
our induction, the case $A=1/q$.  Lemma \ref{eveneven}
and Lemma \ref{evenodd} then imply the Pivot
Theorem by induction.   The proof is essentially the
same as in the odd case, once we see that the
basic structural results hold.  The result in
\S \ref{even2even} gives us the even/even
version of the structure lemma.

We consider the case when $A_1<A_2$.  The other
case is similar.  We define spoilers just in
the odd case.  We just need to show that the arc
$\gamma_2$ defined in the previous section
has no spoilers.
The same argument as in the odd case shows that
a spoiler must lie between $L(\gamma_2)$ and
$R(\gamma_2)$, the left and right endpoints.

Let $\Delta$ be the region
of agreement between $\widehat \Gamma_1$ and
$\widehat \Gamma_2$ as above.  The formulas
are exactly the same.  Here is the even version
of Lemma \ref{nospoil4}.

\begin{lemma}
\label{nospoil5}
$\Delta$ contains all the spoilers.
\end{lemma}

\startproof
The general argument in Lemma \ref{nospoil4} works
exactly the same here.  It is only at the end,
when we consider the vertices $(2,-2)$ and
$(3,-1)$ that we use the fact that $A_2$ is
odd.  Here we consider these special cases
again.  The argument for $(-3,1)$ does not
use the parity of $A_2$.  We just
have to consider $(-2,2)$.

If $A_2<1/2$ then $(-2,2)$ is not a low vertex.
We don't need to treat the extremely trivial
case when $A_2=1/2$.  When $A_2>1/2$ we have
$A_1>1/2$ as well. The point is that no edge
of the Farey graph crosses from $(0,1/2)$
to $(1/2,1)$.  Hence $A_3=A_1 \oplus A_2>1/2$
as well.  But, by definition, the pivot points
relative to $A_2$ are the same as for 
$A_3$.  This is as in \S \ref{even2even}.  Hence,
the same argument as in Lemma \ref{nospoil4}
now rules out $(2,-2)$.
\endproof

Essentially the same argument as in the odd
case now shows that $\gamma_2$ contains
no spoilers.

\newpage

\section{Proof of the Period Theorem}
\label{period}

\subsection{Inheritance of Pivot Arcs}
\label{inheritX}

Let $A$ be some rational parameter.  For each polygonal low component
$\beta$ of $\Gamma(A)$, we define the pivot arc
$P\beta$ to be the lower arc of $\beta$ that joins
the two low vertices that are farthest apart.
We say {\it lower arc\/} because all the components
are closed polygons, and hence two arcs join
the pivot points in all cases. 
When $A$ is an even rational and $\beta=\Gamma$,
this definition coincides with the definition
of $P\Gamma$, by the Pivot Theorem.  In general,
we say that a pivot arc of $\Gamma$ is a pivot
arc of some low component of $\widehat \Gamma$.
We call a pivot arc of $\widehat \Gamma$
{\it minor\/} if it is not a translate of
$P\Gamma$.  

Here we recall the definitions of the odd and even
predecessors of rationals in $(0,1)$. Aside from 
a few trivial cases, the predecessors exist and
are rationals in $(0,1)$.

\begin{enumerate}
\item When $A$ is odd, $A'$ is as in the inferior sequence.
\item When $A$ is odd, $A''$ is as in the Structure Lemma and
Lemma \ref{evenodd}.
\item When $A$ is even, $A'$ is as in the Barrier Theorem.
\item When $A$ is even, $A''$ is as in Lemma \ref{eveneven}.
\end{enumerate}

It is worth mentioning another characterization of these
numbers. 
\begin{equation}
\label{evenfarey}
A \hskip 5 pt {\rm even\/} \hskip 30 pt
\Longrightarrow \hskip 20 pt A=A' \oplus A''.
\end{equation}
\begin{equation}
\label{oddfarey1}
A \hskip 5 pt {\rm odd\/} \hskip 30 pt
\Longrightarrow \hskip 20 pt A=A' \oplus A'' \oplus A''.
\end{equation}

\begin{lemma}[Inheritance]
\label{inherit}
Let $A$ be any rational.  Suppose that
$$A' \leftarrow A; \hskip 30 pt A'' \Leftarrow A.$$  
Then, every minor pivot arc $\beta$ of $\widehat \Gamma$ is either a
minor pivot arc of 
$\widehat \Gamma'$ or a pivot arc of $\widehat \Gamma''$.
The set of low vertices of $\beta$ is the same when
considered either in $A$ or in the relevant predecessor.
\end{lemma}

We first prove the odd case and then we prove the even case.
The proof is almost the same in both case.

{\bf Proof in the Odd Case:\/}
Recall that $P\Gamma \cup \gamma$ is
one period of $\Gamma$.  There are two
kinds of minor components of
$\widehat \Gamma$.

\begin{enumerate}
\item Those pivot arcs that lie underneath $P\Gamma$.
\item Those pivot arcs that lie underneath $\gamma$.
\end{enumerate}

We can push harder on Lemma \ref{welldefined}.
Since $P\Gamma$ lies in the set $\Delta$
from Lemma \ref{dio2}, so does every low
component of $\widehat \Gamma$ underneath
$P\Gamma$.   To
see this, recall that our proof involved
showing that $P\Gamma \subset \Delta$.
But, if a point of $P\Gamma$ lies
in $\Delta$, then so does the entire
line segment connecting this point to
the baseline.  Hence, all components
of $\widehat \Gamma$ beneath $P\Gamma$
also belong to $\Delta$.  
Hence, the low components of
$\widehat \Gamma$ lying underneath
$P\Gamma$ coincide with the
low components of $\widehat \Gamma'$
lying underneath $P\Gamma'$.
This takes care of the first case.

In the second case, our
 proof of Lemma \ref{evenodd} shows that every
minor component of $\widehat \Gamma''$ lying inside
$\Delta(A'',A)$ is contained in $\widehat \Gamma$.
We showed the same result for every major component
except the ones we labelled $C_1$ and $C_{k-1}$.
Note that the pivot arcs are subject to the Barrier Theorem.
That is, the two crossings from the Barrier theorem
occur on the upper arcs rather than on the pivot arcs.
Hence, the pivot arcs behave exactly as the minor
components.  Hence, the pivot arcs of $C_1$
and $C_{k-1}$ are copied by $\widehat \Gamma$,
even though the upper arcs might not be.
By Lemma \ref{nospoil4}, every low vertex of
$\widehat \Gamma$ lying underneath
$\gamma$ lies on the pivot arcs of the
components we have just considered.  This
takes care of the second case.

There is only one detail we need to take care of.
A vertex of the kind we are considering is low
relative to $A'$ or $A''$ is low if and
only if it is low with respect to $A$.
This follows from the basic property of
$\Delta$.  See our geometric proof
of Lemma \ref{dio2}.   Thus, every
low component of $\widehat \Gamma$
of the kind we have considered is also
low relative to $\widehat \Gamma'$
or $\widehat \Gamma''$, whichever
is relevant.  Likewise,
the converse holds.
\endproof

\noindent
{\bf Proof in the Even Case:\/} 
The minor
pivot arcs of $\widehat \Gamma$ come in two kinds, those that
lie underneath $P\Gamma$ and those that do not.
By the same argument as in the odd case,
the pivot arcs of the first kind are all minor
pivot arcs of $\Gamma(A^*)$ where $A^*$ is such that
$A \bowtie !\ A^*$.  But then $A^*=A \oplus A''$.
Hence $A'' \Leftarrow A^*$.  At the same time,
$A'=A \ominus A''$.  Hence $A' \leftarrow A^*$.
Applying the odd case of the Inheritance Lemma to the triple
$(A^*,A',A'')$, we see that every pivot
arc of $\widehat \Gamma$ beneath $P\Gamma$
is a pivot arc of either $\widehat \Gamma'$
or $\widehat \Gamma''$.  This takes care
of the first case.  The second case is just like
the odd case.
\endproof

\subsection{Freezing Numbers}
\label{freezenumber}

Every rational parameter has an odd and an even predecessor.
Starting with (say) an odd rational $A$, we can iterate the
construction and produce a tree of simpler rationals.
If $B$ lies on this tree we write $B \prec A$.  Here
is an immediate corollary of the Inheritance Lemma.

\begin{corollary}
Every minor pivot arc of $\widehat \Gamma(A)$ is a
pivot arc of $\widehat \Gamma(B)$ for some even $B$
such that $B \prec A$.
\end{corollary}

Let $A$ be an odd rational.
Let $\beta$ be a minor component of
$\widehat \Gamma(A)$.   We define
$F(\beta,A)$ to be the smallest denominator
of a rational $B \prec A$ such that
$P\beta$ is a pivot arc of $\widehat \Gamma(B)$.
We call $F(\beta,A)$ the
{\it freezing number\/} of $\beta$.

\begin{lemma}
\label{freeze}
The $\Psi$-period of a minor component $\beta$ is at most
$20s^2$, where $s=F(\beta,A)$.
\end{lemma}

\startproof
This is an immediate consequence of the Hexagrid Theorem,
applied to the rational $B=r/s$ such that
$\beta$ is a component of $\widehat \Gamma(B)$.  The
Hexagrid Theorem confines $\beta$ to a parallelogram
of area less than $20s^2$.
\endproof

Let $x \in I$ correspond to a point not on $C(A_n)$.
We let $$F(x,n)=F(\beta_x,A_n),$$
where $\beta_x$ is the component of $\widehat \Gamma_n$
corresponding to $x$.  We say that a {\it growing sequence\/}
is a sequence $\{x_n\}$ such that
$F(x_n,n) \to \infty$.   Recall that $C_A$ is the Cantor
set from the Comet Theorem.

\begin{lemma}
Suppose every growing sequence has $(0,-1)$ as a limit
point.  Then the Period Theorem is true for $A$.
\end{lemma}

\startproof
If the Period Theorem is
false, then we can find a sequence of points $\{x_n\}$ in
$G_n$ such that the distance from $x_n$ to $C_n$
is uniformly bounded away from $0$, and yet the period of
$x$ tends to $\infty$.  But then Lemma \ref{freeze}
shows that $\{x_n\}$ is a growing sequence. By construction
$\{x_n\}$ does not have a limit point on $C_A$.  In particular,
$(0,-1)$ is not a limit point.
\endproof

\subsection{A Weak Approximation Result}

Let $\{A_n\}$ be the odd sequence of rationals above.
For each $n$ we can form the tree of predecessors,
as above.  Suppose we choose some proper function
$m(n)$ such that $B_m \prec A_n$ is some even
rational in the tree for $A_n$.   

\begin{lemma}
\label{treeapprox}
$\lim_{n \to \infty} B_m=A$.
\end{lemma}

\startproof
We consider the picture in the
hyperbolic plane, relative to the Farey triangulation.
See \S \ref{farey1} for definitions.
We consider the portion $G$ of the Farey graph consisting
of edges having both endpoints in $[0,1]$.  We
direct each edge in $G$ so that it points
from the endpoint of smaller denominator to the
endpoint of larger denominator.  The two endpoints
never have the same denominator, so our
definition makes sense. Say that the {\it displacement\/}
of a directed path in $G$ is the maximum distance between
a vertex of the path and its initial vertex.
  
Given and $\epsilon>0$ there are only finitely many
vertices in $G$ that are the initial points of directed
paths having displacement greater than $\epsilon$.   This
follows from the nesting properties of the
half-disks bounded by the edges in $G$, and from the
fact that there are only finitely many edges in $G$
having diameter greater than $\epsilon$.   

Given the nature of the tree of predecessors, there
is a directed path in $G$ connecting $B_m$ to $A_n$.
The displacement of this path tends to $0$ as $n \to \infty$
because $\{B_m\}$ is an infinite list of rationals
with only finitely many repeaters.  Also, the distance
from $A_n$ to $A$ tends to $0$.  Hence, the distance
from $B_m$ to $A$ tends to $0$ by the triangle inequality.
\endproof

\subsection{The End of the Proof}

To finish our proof, we must show that every growing sequence
has an accumulation point on $C_A$.  We will prove this
indirectly, using the Rigidity Lemma from
\S \ref{rigidity}.  
 Let us first explain the input
from the Rigidity Lemma.

Let $\{B_m\}$ be any sequence of even rationals converging
to the irrational parameter $A$.  Then the Rigidity Lemma
implies that the limits
\begin{equation}
\label{hconverge}
\lim_{m \to \infty} \Gamma(A_m); \hskip 50 pt
\lim_{m \to \infty} \Gamma(B_m)
\end{equation}
agree.  In other words, longer and longer
portions of $\Gamma(A_m)$ look like longer
and longer pictures of $\Gamma(B_m)$.
This is all we need to know from the
Rigidity Lemma. 
\newline

Now, let $M_{m,A}$ be the fundamental map associated to
$A_m$. This map is defined in
Equation \ref{funm}.
 In our proof of Theorem \ref{cantor}, we
showed that 
\begin{equation}
C_A=\lim_{m \to \infty} M_{m,A}(\Sigma(A_m)).
\end{equation}
The limit takes place in the Hausdorff topology.
Here $\Sigma(A_m)$ is the set of low vertices on $\Gamma_m$.
Given Equation \ref{hconverge}, we get the analogous
result
\begin{equation}
C_A=\lim_{n \to \infty} M_{m,B}(\Sigma(B_m)).
\end{equation}

Let's generalize this result.  For each $m$ suppose
there is some $n \geq m$. We also have 
\begin{equation}
C_A=\lim_{m \to \infty} M_{n,A}(\Sigma(B_m)).
\end{equation}
The reason is that the maps $M_{m,A}$ and $M_{n,B}$
converge to each other on any compact subset of 
$\R^2$, and compact pieces of our limit
in Equation \ref{hconverge} determine increasingly
dense subsets of $C_A$.

\begin{lemma}
\label{low2low}
Suppose that $\Sigma_n \subset \widehat \Gamma(A_n)$
is a translate of $\Sigma_m$, consisting entirely
of low vertices.   Then
$$
C_A=\lim_{m \to \infty} M_{n,A}(\Sigma_n).
$$
\end{lemma}

\startproof
We have some vector $U_m$ such that
\begin{equation}
\Sigma_n=\Sigma(A_m)+U_m.
\end{equation}
Since $M_{n,A}$ is affine, we have
\begin{equation}
 M_{n,A}(\Sigma_n)=M_{n,A}\Sigma(A_m)+\lambda_m
\end{equation}
Now we get to the moment of truth.
Since $\Sigma(B_m)$ consists entirely of
low vertices, we have
$M_{A,n}(x) \in [0,2]$ for all $x \in \Sigma(B_m)$.
Since $\Sigma_n$ consists entirely of
low vertices, we have
$M_{A,n}(x)+\lambda_n \in [0,2]$ as well.
Putting $t=M_{A,n}(x)$, we have
\begin{equation}
t; \hskip 30 pt t+\lambda_m \in [0,2].
\end{equation}
This last equation puts constraints on
$\lambda_m$.

By the case $n=0$ of Equation \ref{DIO2}, the
set $C_A$ contains both $0$ and $2$.  Therefore,
once $m$ is large, we can choose $x \in \Sigma(B_m)$
such that $t=M_{A,n}(x)$ is very close to $0$.  
But this forces $$\lim \inf \lambda_m \geq 0.$$
At the same time, we can choose $x$ such that
$M_{A,m}(x)$ is very close to $2$.  This
shows that $$\lim \sup \lambda_m \leq 0.$$
in short $\lambda_m \to 0$. 
\endproof

We just have to tie the discussion above together with our
notion of a growing sequence.
  Suppose that $\{x_n\}$ is a
growing sequence.

Let $\beta_n$ denote the component of $\widehat \Gamma_n$ corresponding
to $x_n$.  There is a proper function $m=m_n$ such that
the pivot arc $P\beta_n$ is a translate
of the major pivot arc $P\Gamma(B_m)$.  Here
$\{B_m\}$ is a sequence of even rationals that satisfies the
hypotheses of Lemma \ref{treeapprox}.  Hence
$\{B_m\} \to A$.  Hence, the application of the
Rigidity Lemma above applies.

Every low vertex on $P\beta_n$ is a translate of a low
vertex on $P\Gamma(B_m)$.  By the Inheritance Lemma,
every low vertex on $P\beta_n$ relative to
$B_m$ is also low with respect to $A_n$.  Thus,
we have exactly the situation described in
Lemma \ref{low2low}.

Let $\Sigma_n$ denote the set of low vertices of
$P\beta_n$.  Then $\Sigma_n$ is a translate of the set 
$\Sigma(B_m)$ of low vertices on $P\Gamma(B_m)$, as in our
lemma above.  Since
\begin{equation}
x_n \in M_{A,n}(\Sigma_n)
\end{equation}
we see that the Hausdorff distance from $\{x\}$ to
$C_A$ tends to $0$ as $n$ (and $m$) tends to $\infty$.

This completes the proof of the
Period Theorem.

\newpage

\section{The End of the Comet Theorem}
\label{gold}

\subsection{The Main Argument}

In this chapter we finish the proof of the
Comet Theorem by proving Statement 1.
Our proof does not use any of the other statements
of the Comet Theorem, so our proof is not circular.
In this first section we give the main argument
modulo several details.  Following this section,
we clear up the several details.

For each rational $B$, we form the
depth-$2$ tree by considering the
two predecessors of $B$, and their two predecessors.
We define the complexity of $B$ to be the
minimum value of all the numerators of the
rationals involved in this list of $7$ rationals.
The point is this definition is that these
are the only rationals that arise in the
geometric constructions we make in this
chapter.

Let $A \in (0,1)$ be irrational.
In this section we consider a sequence $\{B_n\}$ of
ratonals that converges to $A$.  In our applications,
this sequence is the superior sequence, but our
results hold more generally.
Recall that any rational parameter $B$ has its tree $T(B)$
of predecessors. We can consider $T(B_n)$
for each parameter $B_n$ in our sequence.

\begin{lemma}
\label{finitelowcomplexity}
Let $N$ be any integer.  Then there are only finitely many
rationals in the union
$$
\bigcup_{n=1}^{\infty} T(B_n)
$$
having complexity less than $N$.
\end{lemma}

\startproof
We will argue as in the proof of Lemma \ref{treeapprox}.
Suppose $C=r/s$ is a rational in the tree $T(B_n)$ such
that $r$ is small and $s$ and $n$ are large.  Then
the directed Farey path connecting $C$ to $B_n$ has
tiny displacement, and $|B_n-A|$ is small.
Hence $|C-A|$ is small.  Also, $C$ is near $0$.
Hence $A$ is near $0$.  This is a contradiction
once $s$ and $n$ are large enough.
Hence, there is some function $f$,
depending on the sequence, such that $s<f(r)$.
Hence, our union only contains finitely many rationals
having numerator less than $N$.  Our result follows
from this fact.
\endproof

Let $\beta$ be a component of the arithmetic graph.
We call $\beta$ a {\it hovering component\/} if it 
has no low vertices.  More specifically,
we call $\beta$ a {\it $D$-hovering component\/}
of $\widehat \Gamma(A)$ if $\beta$ has no
low vertices, and if $\beta$ contains a vertex
with is within $D$ vertical units of the baseline.
Here $D \geq 2$ is an integer.  We use this name because we think
of $\beta$ as hovering somewhere above the
baseline without coming really close.

Below we prove the following result.

\begin{lemma}
\label{inherit2}
Let $A_2$ be any rational, having the 
predecessors
$A_0\leftarrow A_2$ and $A_1 \Leftarrow A_2$.
Let $\beta$ be a $D$-hovering 
component of $\widehat \Gamma(A)$.  Assuming
that $A_2$ has sufficiently high complexity,
$\beta$ is either a translate of a
$D$-hovering component of $\widehat \Gamma_0$ or a
translate of a $D$-hovering component of $\widehat \Gamma_1$.
\end{lemma}

\begin{corollary}
\label{nohover}
Let $\{A_n\}$ be the superior sequence approximating $A$.
Let $D$ be fixed.  
Then there is a constant $D'$ with
the following property.  If $n$ is sufficiently large,
then $\widehat \Gamma_n$ has no $D$-hovering
components having diameter greater than $D'$.
Here $D'$ is independent of $n$.
\end{corollary}

\startproof
Applying Lemma \ref{inherit} recursively, we see that
$\beta$ is the translate of a $D$-hovering component
of $\widehat \Gamma(C_n)$, where $C_n$ belongs
to the tree of predecessors of $A_n$ and has uniformly
bounded complexity.  But then, by Lemma \ref{finitelowcomplexity},
the sequence $\{C_n\}$ has only finitely many different terms.
Hence $\beta$ is the translate of one of finitely
many different polygons.
\endproof

Below we prove the following result.

\begin{lemma}
\label{pivot2}
Let $\beta$ be a low component of $\widehat \Gamma(B_n)$.
There is some constant $D'$ such that every $D$-low  vertex
of $\beta$
can be connected to a low vertex of $\beta$ in less
than $D'$ steps.  Here $D'$ depends on $D$ and on $A$,
but not on $n$.
\end{lemma}

\begin{lemma}
\label{quick0}
Let $\zeta \in U_A$ such that $\|\zeta\|<N$.
Then there is some $N'$ such that
$\psi^k(\zeta) \in \Xi$ for some $|k|<N'$.
\end{lemma}

\startproof
This is a fairly immediate consequence of our proof
of the Return Lemma in \S 2 and the Pinwheel
Lemma in \S \ref{pinwheel}.  The Return Lemma
takes care of the case when $N$ is small, and the
Pinwheel Lemma takes care of the case when $N$
is large.
\endproof

Recall that $J$ is the interval from Equation \ref{atmosphere}.

\begin{lemma}
\label{quick1}
Let $\zeta \in \Xi$ be such that $\|\zeta\|<N$.
Then there exists $N'$ such that
$\Psi^k(\zeta) \in J$.
\end{lemma}

\startproof
We choose a special interval relative to $A_n$
whose closure contains $\zeta$.  The term
{\it special interval\/} refers to \S \ref{definebasic}.
 Typically the
choice is unique, but when $\zeta$ lies in the
boundary of a special interval there are
two choices and we pick one arbitrarily.
Let $\beta_n$ be the component of $\Gamma(A_n)$ that
tracks this special interval.  There is some uniform
$D$ such that the vertex of $\beta_n$ corresponding
to our special interval is $N$-low.

By the Continuity Principle, ${\rm diam\/}(\beta_n) \to \infty$.
By Corollary \ref{nohover}, we see that
$\beta_n$ is a low component for $n$ sufficiently large.
By Lemma \ref{pivot2}, the vertex
corresponding to $\zeta$ can be connected
to a low vertex within $N'$ steps.  But then
there is a sequence $\{k_n\}$ such
that
\begin{equation}
\Psi_n^{k_n}(\zeta) \in [0,2] \times \{-1,1\}; \hskip 30 pt
|k_n|<N'.
\end{equation}
Here $\Psi_n$ is the first return map defined relative to $A_n$.
The important point here is that $N'$ is independent of $n$.
This lemma now follows from the Continuity Principle
from \S 2.
\endproof

\begin{lemma}
\label{quick2}
If $x \in \Xi$ is such that $|x|<R$, then there
is some $R'$ such that the portion of the
$\psi$-orbit of $x$ between $x$ and $\Psi(x)$
has cardinality at most $R'$.
\end{lemma}

\startproof
same proof as Lemma \ref{quick0}.
\endproof

Statement 1 of the Comet Theorem now follows
from Equation \ref{atmosphere}, Lemma
\ref{quick0}, Lemma \ref{quick1}, and
Lemma \ref{quick2}.

Our work is almost done.
The two remaining details are Lemma
\ref{inherit2} and Lemma \ref{pivot2}.
We establish these results in the sections
below.

\subsection{Proof of Lemma \ref{inherit2}}

\subsubsection{Traps}

For $j=0,1$, let $\Delta_j$ denote the region of agreement between
$\widehat \Gamma_j$ and $\widehat \Gamma_2$, as in Lemma
\ref{dio2}.   Call a parallelogram $X_j$ a {\it trap\/} if
$X_j \subset \Delta_j$, and no hovering component relative to
$A_j$ crosses $\partial X_j$.  
Say that a pair $(X_1,X_2)$ is a $D$-{\it trap\/}
if $X_j$ is a trap relative to
$(A_j,A_2)$ and if every $D$-low vertex, relative to
$A_2$, is translation equivalent to a vertex in one of the
traps. As usual, the translation takes place in $\Z V_2$.

Below we will prove the following result.
\begin{lemma}
\label{trapp}
If $A_2$ has high complexity, then there
exists a $D$-trap.
\end{lemma}

Lemma \ref{inherit2} follows immediately from Lemma \ref{trapp}.   Let
$\beta_2$ be a $D$-hovering component of $\widehat \Gamma_2$.
Let $v$ be a $D$-low vertex of $\beta$.  We can translate so that
$v$ lies in either $X_0$ or $X_1$.   Suppose that $v \in X_0$.
Let $\beta_0$ be the component of $\widehat \Gamma_0$ that contains
$v$.  Since $X_0$ is a trap, $\beta_0 \subset X_0$.  But
$X_0$ is a region of agreement between $\widehat \Gamma_0$
and $\widehat \Gamma_2$.   Hence $\beta_0=\beta_2$.
The same argument works if $v \in X_1$.  

Now we define the traps.
First we make some general comments.  The rational
$A_0$ is odd and the rational $A_1$ is even.  The parallelogram
$X_0$ is always bounded by the lines used in the Decomposition
Theorem.  The parallelogram $X_1$ is the one we used in the
proof of the Pivot Theorem.   The reader can see the traps
drawn, in all cases, using Billiard King.

Now we get down to specifics.  There are $4$ cases:
\begin{enumerate}
\item $A_2$ is odd and $A_1<A_2$.
\item $A_2$ is odd and $A_1>A_2$.
\item $A_2$ is even and $A_1<A_2$.
\item $A_2$ is even and $A_1>A_2$.
\end{enumerate}
In our proof of the Pivot Theorem, we considered
Cases 1 and 3.  Here we will consider Cases 1 and 3
in detail, and just remark briefly on Cases 2 and 4.
The reader will see that Case 2 is essentially
identical to Case 1 and Case 4 is essentially identical
to Case 2.

The reader can see the traps drawn in Billiard King,
for any desired smallish parameter.

\subsubsection{Case 1}

We first reconcile some bits of notation.  In this case we have
\begin{equation}
A_0=(A_2)_-; \hskip 30 pt
A_1+A_0=(A_2)_+.
\end{equation}
Both $A_0$ and $A_2$ are odd rationals.

We define $X_0=R_1(A_2)$, the parallelogram
used in the decomposition Theorem for $A_2$.   We define
$X_1=R$, the parallelogram defined at the end
of \S \ref{juggle}.  

In \S \ref{evenimpliesodd} we showed that $R_1 \subset \Delta_1$
when $A_2$ has high complexity.  Lemma \ref{symm4} shows
that $R_0 \subset \Delta_0$ when $A_2$ has high complexity.
However, since our notation has changed slightly, and
since we are considering the opposite case from the one
in Lemma \ref{symm4} -- namely, $R_0$ here lies to
the left of the origin -- we will re-work the proof.

\begin{lemma}
$X_0 \subset \Delta_0$.
\end{lemma}

\startproof
We will apply the Diophantine Lemma.
We work with the linear functionals $G_2$ and $H_2$ associated
to the parameter $A_2$. 
Let $u$ and $w$ denote the top left and right vertices
of $X_0$ respectively.  The interval in the
Diophantine Lemma is
\begin{equation}
[-(q_2)_--q_0+2,q_0-2].
\end{equation}
The lower bound comes from Case 2 of Lemma \ref{keycomp1}.

Hence, it suffices to show that
to show that
\begin{equation}
G_2(u)>>-(q_2)_--q_0 \hskip 50 pt
H_2(w)<<q_0
\end{equation}
The $(>>)$ symbol means an inequality in which the difference
between the two sides tends to $\infty$ with the
complexity of $A_2$.

 We have the estimates
\begin{equation}
u \approx -(V_2)_-+ \lambda W_2; \hskip 30 pt
w \approx \lambda W_2; \hskip 30 pt
\lambda = \frac{q_2^*}{q_2} \leq \frac{q_0}{q_2}.
\end{equation}
Here $A_2^*=p_2^*/q_2^*$ is the superior predecessor of $A_2$.
The approximation becomes arbitrarily good as the complexity
of $A_2$ tends to $\infty$.  In particular, the approximation
is good to within $1$ unit once $A_2$ has sufficiently
high complexity.

We compute
$$
G_2(u) \approx  -(q_2)_- - \lambda \frac{q_2^2}{p_2+q_2}>>
-(q_2)_- - \lambda (q_2) \geq
-(q_2)_- - q_0.$$
This takes care of the vertex $u$.  Now we compute
$$H_2(w) \approx \lambda \frac{q_2^2}{p_2+q_2}<<
\lambda q_2=q_0.$$
This takes care of the vertex $w$.
\endproof

Now we know that $X_j \subset \Delta_j$ for $j=0,1$.
Indeed, our proof shows that a large neighborhood of
$X_j$ is contained in $\Delta_j$ once the complexity
of $A_2$ is large.

\begin{lemma}
$X_0$ is a trap.
\end{lemma}

\startproof
By the Decomposition Theorem, the only component of
$\widehat \Gamma_2$ that crosses $X_0$ is $\Gamma_2$,
a low component.  The crossings occur within $1$ unit
of the bottom vertices of $X_0$.  Since
$\widehat \Gamma_0$ and $\widehat \Gamma_2$ agree
in a neighborhood of $X_0$, the same structure holds
for $\widehat \Gamma_0$.  The only places where a
component of $\widehat \Gamma_0$ crosses $X_0$ are
at low vertices.  In particular, no hovering
component of $\widehat \Gamma_0$ crosses $X_0$.
\endproof

\begin{lemma}
$X_1$ is a trap.
\end{lemma}

\startproof
We already saw in \S \ref{boundminor} that the
only components of $\widehat \Gamma_1$ that cross
$X_1$ are the ones on major components of
$\widehat \Gamma_1$.  But, such components are
not hovering components.  Hence $X_1$ is a trap.
\endproof

Here is the last remaining step.
\begin{lemma}
The pair $(X_0,X_1)$ is a $D$-trap
once $A_2$ has high complexity.
\end{lemma}

\startproof
  The left bottom vertex of
$X_0$ is $-(V_2)_-$ whereas the bottom right vertex of
$X_1$ is $(V_2)_+$.  These two vertices differ by $V_2$.
The bottom right vertex of $X_0$ is $(0,0)$, the same
as the bottom left vertex of $X_1$.
Figure 31.1 shows the picture.

\begin{center}
\resizebox{!}{1.8in}{\includegraphics{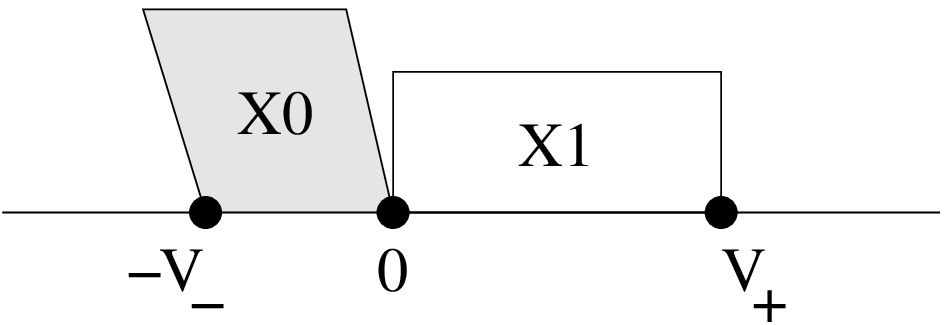}}
\newline
{\bf Figure 31.1:\/} The trap
\end{center} 

Suppose for the moment that the sides of $X_0$ have the same
slope as the sides of $X_1$.  Then, once $A_2$ has high complexity,
the tops of both parallelograms are more than $D$ units from the
baseline.  But then the union of translations
\begin{equation}
\label{boxorbit}
\bigcup_{k \in \Z} \Big(X_0+X_1+kV_2)
\end{equation}
contains all $D$-low vertices, as desired.

The slight complication is that the sides of $X_0$ are parallel to
$W_2$ whereas the sides of $X_1$ are parallel to $W_1$. These
are the vectors from Equation \ref{boxvectors}, relative to
$A_2$ and $A_1$. As the complexity of $A_2$ tends to $\infty$,
the slopes converge, and no $D$-low lattice point lies between
the two lines emanating from the same point.  Thus, our union
in Equaton \ref{boxorbit}
still
contains all $D$-low vertices once $A_2$ has high complexity.
\endproof

\subsubsection{Case 2}
\label{othercase}

In this case we have $A_1>A_2$.  We take $X_0=R_1(A_2)$, the
smaller of the two parallelograms in the Decomposition
Theorem.  This time $X_0$ lies to the right of the origin.
We take $X_1$ to be just like the parallelogram
defined at the end of \S \ref{juggle}, except that
$-(V_2)_-$ replaces $(V_2)_+$.
Here $X_1$ lies to the
left of the origin.  The picture looks exactly like
Figure 31.1, except that the roles of left and right are
reversed, and the subscripts $(+)$ and $(-)$ are switched 
in the labels.
Aside from switching the roles placed by left and right, and
$(+)$ and $(-)$, the proofs for Case 2 are exactly the
same as the proofs for Case 1.

\subsubsection{Case 3}

We first reconcile some bits of notation.  In this case we have
\begin{equation}
A_0=(A_2)_+; \hskip 30 pt
A_1=(A_2)_-.
\end{equation}
Here $A_1$ and $A_2$ are even and $A_0$ is odd.
We define $X_1$ exactly as in Case 1, using the
rectangle $R$ described at the end of \S \ref{juggle}.
The same argument as in Case 1 shows that $X_1$ is a
trap.  

We define $X_0$ to be the parallelogram bounded by the
following lines.

\begin{enumerate}
\item The baseline relative to $A_0$.
\item The line parallel to $V_0$ and containing $W_0$.  This
is the top of the room $R(A_0)$ from the Room Lemma.
\item The line parallel to $W_0$ and containing $(0,0)$.
\item The line parallel to $W_0$ and containing $-(V_2)_-$.
\end{enumerate}

\begin{lemma}
$X_0 \subset \Delta_0$ once $A_2$ has
sufficiently high complexity.
\end{lemma}

\startproof
We will apply Lemma \ref{dio2}.  This time,
we work with the linear functionals $G_0$ and $H_0$ associated
to the parameter $A_0$. 
Let $u$ and $w$ denote the top left and right vertices
of $X_0$ respectively.  The interval in the
Diophantine Lemma is
\begin{equation}
[-q_2+2,q_0-2].
\end{equation}

Hence, it suffices to show that
to show that
\begin{equation}
G_2(u)>>-q_2 \hskip 50 pt
H_2(w)<<q_0
\end{equation}

We have
\begin{equation}
u = -(V_2)_-+W_0; \hskip 30 pt
w = W_0.
\end{equation}

We compute
$$
G_0(u) \approx -(q_2)_- - \frac{q_0^2}{p_0+q_0} >> -(q_2)_--q_0=-(q_2)_--(q_2)_+=-q_2.$$
This takes care of the vertex $u$.  Now we compute
$$H_2(w) =\frac{q_0^2}{p_0+q_0}<<q_0.$$
This takes care of the vertex $w$.
\endproof

\begin{lemma}
$X_0$ is a trap.
\end{lemma}

\startproof
The same argument as in Lemma \ref{affine3} shows that
\begin{equation}
-(V_2)_-=-(V_0)_-+kV_0.
\end{equation}
for some $k \in \Z$.
Geometrically, this says that 
\begin{equation}
X_0=Y_1 \cup ... \cup Y_k \cup Z; \hskip 30 pt
Y_j=R(A_0)-j V_0.
\end{equation}
Here $R(A_0)$ is the parallelogram from the Room Lemma.  The
parallelogram $Z$ is $\Z V_0$ translate of the parallelogram
$Z'$ bounded by the following lines.
\begin{enumerate}
\item The baseline relative to $A_0$.
\item The line parallel to $V_0$ and containing $W_0$.  This
is the top of the room $R(A_0)$ from the Room Lemma.
\item The line parallel to $W_0$ and containing $(0,0)$.
\item The line parallel to $W_0$ and containing $(V_2)_+$.
\end{enumerate}
But each $Y_j$ separately is a trap by the Room Lemma.
The parallelogram $Z'$ is also a trap, by the same
argument we gave in the proof of the Decomposition Lemma.
This argument is repeated in the proof of Lemma
\ref{nocross}.  Hence, by symmetry $Z$ is also a trap.

Hence, $X_0$ is a finite union of traps, all meeting
edge to edge.   Hence, $X_0$ is also a trap.
\endproof

It only remains to show that
the pair $(X_0,X_1)$ is a $D$-trap.  The bottom
vertices in this case have the same description as
in Case 1, and the argument there works here
word for word.

\subsubsection{Case 4}

We define $X_1$ just as in Case 2.  We define $X_0$ as in Case 3,
except that we replace the vector $-(V_2)_-$ by the
vector $(V_2)_+$.  The rest of the proof is the same
as in Case 3, modulo the same switching of ``left'' and
``right''.

\subsection{Proof of Lemma \ref{pivot2}}

\subsubsection{Major Components}

We keep the notation from the previous section.

\begin{lemma}
\label{arcs2}
When $A_2$ has sufficiently high complexity,
the set
$$(\Gamma_2-P\Gamma_2) \cap X_1$$
consists of $2$ connected arcs, each joining an
endpoint of $\gamma_2$ to the top of $X_1$.
\end{lemma}

\startproof
In the even case, this is a restatement of
Lemma \ref{even nocross}.

\begin{center}
\resizebox{!}{1.8in}{\includegraphics{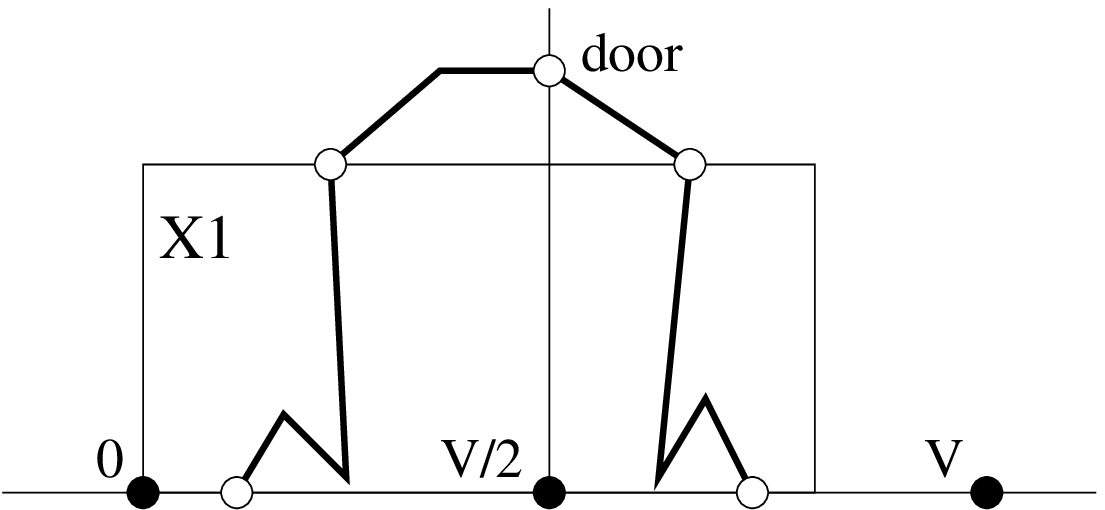}}
\newline
{\bf Figure 31.2:\/} The arc $\gamma_2$ in the odd case.
\end{center} 

We consider the odd case.   
We describe the case when $A_1<A_2$.  The other
case is entirely similiar. The two endpoints
of $\gamma_2$ are $E_2^+$ and $E_2^-+V_2$.
Both these points belong to $X_1$.
The line parallel to $W_2$ through $V_2/2$
divides $X_1$ into two pieces.  By the
Hexagrid Theorem, $\gamma_2$ crosses a door
on this line.  This door lies above the top
of $X_1$.  At the same time, $\gamma_2$ can
only cross the top of $X_1$ twice.  This
follows from the Barrier Theorem, as applied
to $A_1$, and from the fact that
$\widehat \Gamma_1$ and $\widehat \Gamma_2$
agree in a neighborhood of $X_1$.  This
structure forces the following structure.
Starting from the left endpoint of
$\gamma_2$, some initial arc of $\gamma_2$
rise up to the top of $X_1$.  Following this,
the next arc of $\gamma_2$ crosses through a
door and returns to the top of $X_1$.  The
final arc of $\gamma_2$ connects the top of
$X_1$ to the right endpoint of
$\gamma_2$.   
\endproof

Now we derive some corollaries from our structure result.
We say that a $D$-arc of $\Gamma_2$ is a connected arc $\alpha$ that
joins a low vertex to a $D$-low vertex.
Let $|\alpha|$ denote the smallest integer $N$ such that
$\alpha$ contains no vertices that are $N$ vertical
units above the baseline.  So $\alpha$ remains within
$N$ vertical units of the baseline.
Given a $D$-low vertex $v \in \alpha$, let
$f(v;A_2))=\min |\alpha|$, where the minimum is taken over
all $D$-arcs having $v$ as an endpoint.
Let $F(A_2)=\max f(v;A_2)$, where the maximum is taken
over all $D$-low vertices $v$ of $\Gamma_2$.
These functions depend implicitly on $D$, which
is fixed throughout the discussion.

\begin{lemma}
If $A_2$ has sufficiently high complexity, then
$$F(A_2) \leq \max\Big(F(A_0),F(A_1)\Big).$$
\end{lemma}

\startproof
Suppose first that $A_2$ is odd.
Choose a vertex $v \in \Gamma_2$ such that $F(A_2)=f(v)$.
By symmetry, we can choose $v \in P\Gamma_2 \cup \gamma_2$.
Suppose that $v \in P\Gamma_2$.  By the Copy Theorem,
$P\Gamma_2 \subset \Gamma_0$.  Hence $v \in \Gamma_0$.
The argument in Lemma \ref{lowlow} shows that a vertex
on $P\Gamma_2$ is $E$-low relative to $A_0$ if and only
if it is $E$-low relative to $A_2$.  Here $E \in \{1,2,3...\}$.
 Call this the
{\it low principle\/}.  Let $\alpha$ be a $D$-arc of
$\Gamma_0$ such that $f(v;A_0)=|\alpha|$.
Since both endpoints of $P\Gamma_2$ are low relative
to both parameters, we can take
$\alpha \subset P\Gamma_2$.  Hence
$$F(A_0) \geq f(v;A_0)= |\alpha| \geq f(v;A_2)= F(A_2).$$
Note that $|\alpha|$ is the same
relative to both parameters, by the low principle.

Suppose that $v \in \gamma_2$. Then $v$ is in one of the
two arcs from Lemma \ref{arcs2}.  Let's say that
$v$ is on the left arc, $\lambda$.  The left endpoint
of $\lambda$ is common to $\Gamma_1$ and $\Gamma_2$,
and $\lambda \subset X_1$, a region of agreement for
the two arithmetic graphs.  Hence $\lambda \subset \Gamma_1$.
The low principle applies
to any vertex in $X_1$.
 Let $\alpha$ be a $D$-arc of
$\Gamma_1$ such that $f(v;A_1)=|\alpha|$.
The left endpoint of $\lambda$ is low, and the right
endpoint lies on the top of $X_1$.  When $A_2$ has
high complexity, $\alpha \subset \lambda$.
The idea here is that the $D$-arc connecting
$v$ to the left endpoint of $\lambda$ remains
in $X_1$ whereas any $D$-arc exiting $\lambda$
must pass through the top of $X_1$.
Since $\alpha \subset \lambda$, we get
$F(A_1) \geq F(A_2)$ by the same argument as in the
previous case.

When $A_2$ is even, the proof is the same except for two small
changes.  First, we need to invoke Lemma \ref{even nocross} (rather
than just symmetry) to get $v \in P\Gamma_2 \cup \gamma_2$.
Second, when $v \in P\Gamma_2$, we use Lemma
\ref{copyth} in place of the Copy Theorem.
\endproof

Now let $\{B_n\}$ be the sequence in Lemma \ref{pivot2}.

\begin{corollary}
\label{uniform short}
$F(B_n)$ is uniformly bounded, independent of $n$.
\end{corollary}

\startproof
Applying the previous result recursively, we see that
there is some parameter $C_n \in T(B_n)$, of uniformly
bounded complexity, such that $$F(B_n) \leq F(C_n).$$
But the sequence $\{C_n\}$ has only finitely many
distinct members, by Lemma \ref{finitelowcomplexity}.
\endproof

\begin{corollary}
Let $v_n$ be a $D$-low vertex on $\Gamma(B_n)$.  Then $v_n$
can be connected to a low vertex of $\Gamma_n$ by an arc
of length less than $D'$ for some $D'$ that is independent
of $n$.   Hence, Lemma \ref{pivot2} is true for
points that lie on major components.
\end{corollary}

\startproof
By Corollary \ref{uniform short} we can find a $D$-arc
$\alpha_n$ connecting $v_n$ to a low vertex of
$\Gamma(B_n)$ such that $|\alpha_n|<N$ and $N$ is
independent of $n$.  But the same argument as in
the proof of Lemma \ref{dichotomy} shows that
the diameter of $\alpha_n$ is uniformly bounded.
The idea here is that $\alpha_n$ cannot grow
a long way in a thin neighborhood of the baseline.
This takes care of $D$-low vertices on
$\Gamma(B_n)$.  
Any other major
component of $\widehat \Gamma(B_n)$ is translation
equivalent to $\Gamma(B_n)$.
\endproof

\subsubsection{A Quick but Unjustified Finish}
\label{unjust}

It remains to prove Lemma \ref{pivot2} for points that
lie on minor components.  First we will give a short
proof that we cannot quite justify.  Then we will
patch up the argument.
Experimentally, we observe the following strengthening
of the Inheritance Lemma.

\begin{conjecture}
\label{extreme}
Let $A_2$ be any rational, having the 
predecessors
$A_0\leftarrow A_2$ and $A_1 \Leftarrow A_2$.
Then every minor low component of $\widehat \Gamma_2$
is either the translate of a low component
of $\widehat \Gamma_0$ or the translate of a
low component of $\widehat \Gamma_1$.
\end{conjecture}

Assuming this conjecture, we can quickly finish
the proof of Lemma \ref{pivot2}.   
Suppose that Lemma \ref{pivot2} is false.
Then we can find a sequence of
triples $\{(v_n,\beta_n,B_n)\}$ with the
following properties.
\begin{enumerate}
\item $\beta_n$ is a minor component of $\widehat \Gamma(B_n)$;
\item $v_n$ is a vertex of $\beta_n$ that lies within
$D$ units of the baseline;
\item The $n$-neighborhood of $v_n$ in $\beta_n$ contains no low
vertices.
\end{enumerate}

By Conjecture \ref{extreme}, the component
$\beta_n$ is the translate of 
$\Gamma(B_n')$ for some $B_n' \in T(B_n)$.
Since the diameter of $\beta_n$ tends to $\infty$ with $n$,
the complexity of $B_n'$ tends to $\infty$.  Hence,
by Lemma \ref{treeapprox}, $B_n' \to A$.  Thus,
Thus, a counterexample to Lemma \ref{pivot2} involving minor components
leads to a counterexample involving major components.

We can't quite prove Conjecture \ref{extreme}.  Our
approach to Lemma \ref{pivot2} is to prove a slightly
weaker version of Conjecture \ref{extreme}, and
then scramble to finish the proof.

\subsubsection{The End of the Proof}

Let $A$ be an even
rational.   Previously, we had divided
the polygon $\Gamma(A)$ into two arcs,
the pivot arc $P\Gamma(A)$ and the upper
arc.  These two arcs join together at the
pivot points.  

\begin{center}
\resizebox{!}{1.1in}{\includegraphics{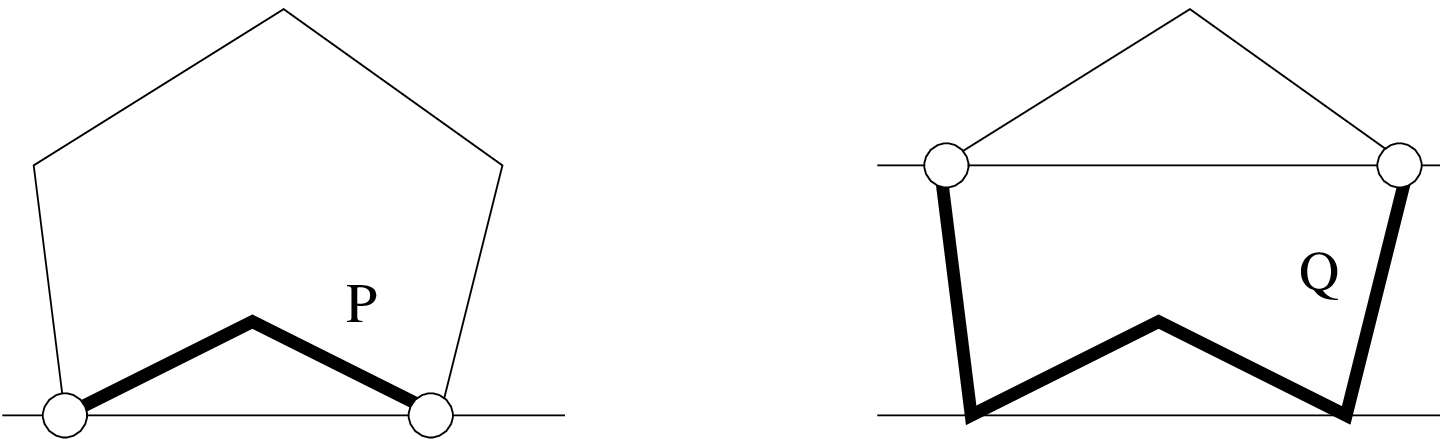}}
\newline
{\bf Figure 31.3:\/} $P\Gamma$ and $Q\Gamma$.
\end{center} 

Now we consider a new
decomposition of $\Gamma$.  Referring
to the Barrier Theorem, recall that
$\Gamma(A)$ passes through the barrier
at $2$ points.  One arc of $\Gamma$ lies
below the barrier and one above.  
Let $Q\Gamma$ denote the component that lies
below.  Then $P \Gamma \subset Q \Gamma$.
We call $Q\Gamma$ an {\it extended pivot
arc\/}.   
We think of $Q\Gamma$ as a kind
of compromise between the whole component
$\Gamma$ and the pivot arc $P\Gamma$.
If $A$ has sufficiently high complexity, then
$Q\Gamma$ contains all the vertices within
$D$ of the baseline.  This is a consequence of 
the Barrier Theorem.

So far we have only defined $Q\beta$ when
$\beta=\Gamma(A)$ and $A$ is an even rational.
Our strengthening of the Inheritance
Lemma extends this definition to all polygonal
low components of $\widehat \Gamma(A)$,
when $A$ is any rational parameter.

\begin{lemma}
\label{extreme2}
Let $A_2$ be any rational, having the 
predecessors
$A_0\leftarrow A_2$ and $A_1 \Leftarrow A_2$.
If $A_2$ has sufficiently high complexoty, then
every low component of $\widehat \Gamma_2$
has a well-defined extended pivot arc, and
this pivot arc is the translate of an extended
pivot arc of $\widehat \Gamma_j$ for one
of $j=0,1$.
\end{lemma}

\startproof
We will use the existence of the traps
$X_0$ and $X_1$.  Let $\beta$ be a low
component of $\widehat \Gamma_2$.  Let
$v$ be a low vertex of $\beta$.  If
$v \in X_0$ then $\beta$ is copied whole
from $\widehat \Gamma_0$.  If
$\beta \in X_1$, then $\beta$ is copied
whole from $\widehat \Gamma_1$ unless
(in the language of \S \ref{evenimpliesodd})
$\beta=C_1$ or $\beta=C_{k-1}$. In these cases,
$QC_1$ and $QC_{k-1}$, both subsets of $X_1$,
are copied whole by $\widehat \Gamma_2$.
Hence, the portion of $\beta$ lying
in $X_2$ is copied from $\widehat \Gamma_1$.

Let $\widetilde \beta$ denote the component
of $\widehat \Gamma_2$ that contains
$\beta \cap X_1$.  We define
$Q\widetilde \beta=Q\beta$.  Then
$Q\widetilde \beta$ is copied from
$\widehat \Gamma_1$ by construction.
\endproof

The following result is an {\it addendum\/} to the
proof of Lemma \ref{extreme2}.  We would like to
say that the components $\widetilde C_0$ and
$\widetilde C_{k-1}$, though perhaps
imperfect copies of $C_0$ and $C_{k-1}$, still
retain a basic property of the original components.

\begin{lemma}
\label{nobadarc}
Let $N$ be fixed.  If $A_2$ has sufficiently high complexity,
then $\widetilde C_1-Q\widetilde C_1$ does not contain
any vertices within $N$ units of the baseline.  The
same goes for $C_{k-1}$.
\end{lemma}

\startproof
As in our proof of the Pivot Theorem, we consider the
case when $A_1<A_2$.  The other case is entirely similar.

Let $\gamma=\widetilde C_1-Q\widetilde C_1$.
Here $\gamma$ is an arc of $\widehat \Gamma_2$.
Let $R$ be the parallelogram we used in the proof of the
Pivot Theorem.  We
defined $R$ at the end of \S \ref{juggle}.
Recall that $\widehat \Gamma_1$
and $\widehat \Gamma_2$ agree in $R$. 
The component $\widetilde C_1$ has a low vertex in
$R$.   The arc $\gamma$ has both its endpoints on
the top edge of $R$.

Let $S$ denote the infinite strip
obtained by extending the left and right sides of $R$.
We claim that $\widetilde C_1$
does not cross either side of $S$.
To prove this claim, let
$S_L$ and $S_R$ denote the left and right
boundaries of $S$.  Then $\widetilde C_1$
does not cross $S_L$ by the
Hexagrid Theorem applied to $A_2$. 
Likewise, $\iota(\widetilde C_1)$ does not
cross $S_L$ by the Hexagrid Theorem.
Here $\iota$ is the same symmetry as in
Lemma \ref{nocross}.  By construction,
$\iota$ swaps $S_L$ and $S_R$.  Hence,
$\widetilde C_1$ does not cross $S_R$.
This establishes our claim.

Now we know that $\gamma$ does not
cross the sides of $S$.  Hence, if $\gamma$
contains a vertex within $N$ units of the
baseline, this vertex must lie in $R$.  But
then $\widetilde C_1$ crosses the
top edge of $R$ at least $4$ times.  But these
$4$ crossing points are then copied from
$\widehat \Gamma_1$.  This contradicts
the Barrier Theorem, because the top edge of
$R$ is contained in the barrier line for
$\widehat \Gamma_1$.
\endproof

Let $\{B_n\}$ be the sequence from Lemma \ref{pivot2}.

\begin{corollary}
\label{nobadarc2}
Let $\{\beta_n\}$ be a sequence of components,
with $\beta_n$ a low component of $\widehat \Gamma(B_n)$.
Suppose that the diameter of $\beta_n$ tends to $\infty$.
Then the distance from any point on $\beta_n-Q\beta_n$
to the baseline of $\widehat \Gamma(B_n)$ tends to
$\infty$ as well.
\end{corollary}

\startproof
This is a fairly immediate consequence of the
previous result.  Each $\beta_n$ is a translate
of a component of the form $\widetilde C$,
$C=\Gamma(B_n')$.  Here $B_n'$ is on the tree
of predecessors of $B_n$.  Since the diameter
of $\widetilde C$ tends to $\infty$ with $n$,
we see than the complexity of $B'_n$ tends to
$\infty$ with $n$ by Lemma \ref{finitelowcomplexity}.
Hence, the distance from
$\widetilde C-Q\widetilde C$ to the relevant baseline
tends to $\infty$ with $n$.
\endproof

Now we redo the argument in \S \ref{unjust}
equipped with our weaker but sufficient results.
From Lemma \ref{extreme2}, we conclude that $Q\beta_n$ is
the translate of $Q\Gamma(B'_n)$ for some other
sequence $B'_n \to A$.  This is just as in
the proof of Lemma \ref{nobadarc2}.  Thus, a
counterexample to Lemma \ref{pivot2} involving minor components
leads to a counterexample involving major components.
Since we have already taken care of the major components,
our proof is done.

\newpage

\section{References}

[{\bf B\/}] P. Boyland, {\it Dual Billiards, twist maps, and impact oscillators\/},
Nonlinearity {\bf 9\/} (1996) 1411-1438
\newline
\newline
[{\bf De\/}] N .E. J. De Bruijn, {\it Algebraic Theory of Penrose's Nonperiodic Tilings\/},
Nederl. Akad. Wentensch. Proc. {\bf 84\/} (1981) pp 39-66
\newline
\newline
[{\bf Da\/}] Davenport, {\it The Higher Arithmetic: An Introduction to the Theory of Numbers\/},
Hutchinson and Company, 1952
\newline
\newline
[{\bf D\/}], R. Douady, {\it These de 3-eme cycle\/}, Universite de Paris 7, 1982
\newline
\newline
[{\bf DF\/}] D. Dolyopyat and B. Fayad, {\it Unbounded orbits for semicircular
outer billiards\/}, preprint (2008)
\newline
\newline
[{\bf DT\/}] F. Dogru and S. Tabachnikov, {\it Dual Billiards\/},
Math Intelligencer vol. 27 No. 4 (2005) 18--25
\newline
\newline
[{\bf G\/}] D. Genin, {\it Regular and Chaotic Dynamics of
Outer Billiards\/}, Penn State Ph.D. thesis (2005)
 \newline \newline
[{\bf GS\/}] E. Gutkin and N. Simanyi, {\it Dual polygonal
billiard and necklace dynamics\/}, Comm. Math. Phys.
{\bf 143\/} (1991) 431--450
\newline
\newline
[{\bf H\/}] M. Hochman, {\it Genericity in Topological Dynamics\/},
Ergodic Theory and Dynamical Systems {\bf 28\/} (2008) pp 125-165
\newline
\newline
[{\bf Ke\/}] R. Kenyon, {\it Inflationary tilings with a similarity
structure\/}, Comment. Math. Helv. {\bf 69\/} (1994) 169--198
\newline
\newline
[{\bf Ko\/}] Kolodziej, {\it The antibilliard outside a polygon\/},
Bull. Polish Acad Sci. Math.
{\bf 37\/} (1989) 163--168
\newline
\newline
[{\bf M1\/}] J. Moser, {\it Is the Solar System Stable?\/},
Mathematical Intelligencer, 1978
\newline
\newline
[{\bf M2\/}] J. Moser, {\it Stable and Random Motions in Dynamical Systems, with
Special Emphasis on Celestial Mechanics\/},
Annals of Math Studies {\bf 77\/}, Princeton University Press (1973)
\newline
\newline
[{\bf MM\/}] P. Mattila and D. Mauldin, {\it Measure and dimension functions:
measurability and densities\/} Math. Proc. Cambridge Philos. Soc. {\it 121\/} No. 1 (1997)
\newline
\newline
[{\bf N\/}] B.H. Neumann, {\it Sharing Ham and Eggs\/},
\newline
summary of a Manchester Mathematics Colloquium, 25 Jan 1959 
\newline
published in Iota, the Manchester University Mathematics students' journal
\newline
\newline
[{\bf S\/}] R. E. Schwartz, {\it Unbounded Orbits for Outer Billiards\/},
Journal of Modern Dynamics {\bf 3\/} (2007) 
\newline
\newline
[{\bf T1\/}] S. Tabachnikov, {\it Geometry and Billiards\/},
A.M.S. Math. Advanced Study Semesters (2005)
\newline
\newline
[{\bf T2\/}] S. Tabachnikov, {\it A proof of Culter's theorem on the existence of periodic
orbits in polygonal outer billiards\/}, preprint (2007)
\newline
\newline
[{\bf VL\/}] F. Vivaldi and J.H. Lowenstein, Arithmetical properties of a family
of irrational piecewise rotations, {\it Nonlinearity\/}, (2006) 19:1069-1097
\newline
\newline
[{\bf VS\/}] F. Vivaldi, A. Shaidenko, {\it Global stability of a class of discontinuous
dual billiards\/}, Comm. Math. Phys. {\bf 110\/} (1987) 625--640 
\newline
\newline

\end{document}